\documentclass[preprint,authoryear,10pt,3p]{elsarticleFTP}

\usepackage{amssymb}
\usepackage{tikz}
\usepackage{graphicx}
\usepackage{caption}
\usepackage{mathtools}
\usepackage{amsmath}
\usepackage{multirow}
\usepackage{colortbl}
\usepackage{color}
\usepackage{tabulary}
\usepackage{bm}
\usepackage{amsmath}
\usepackage{eqparbox}
\usepackage{tabularx}
\usepackage{tabulary}
\usepackage[utf8]{inputenc}
\usepackage{booktabs}
\usepackage{kpfonts}
\usepackage{rotating}
\usepackage{enumitem}
\usepackage{hyperref}
\usepackage{breakurl}
\usepackage{setspace}
\usepackage{epigraph}
\usepackage{longtable}
\usepackage{pdflscape}

\journal{Journal of the Operational Research Society}

\begin{document}

\begin{frontmatter}
\title{Operational Research: Methods and Applications\footnote{This is an original manuscript of an article published by Taylor \& Francis in the \textit{Journal of the Operational Research Society} on 27th December 2023, available at: \url{https://doi.org/10.1080/01605682.2023.2253852}}}

\author[label1,label1b]{Fotios~Petropoulos\corref{cor1}}
\cortext[cor1]{Corresponding author: f.petropoulos@bath.ac.uk; fotios@bath.edu}
\author[label7a,label1,label401]{Gilbert~Laporte}
\author[label100]{Emel~Aktas}
\author[label2]{Sibel~A.~Alumur}
\author[label3]{Claudia~Archetti}
\author[label4]{Hayriye~Ayhan}
\author[label1]{Maria~Battarra}
\author[label5]{Julia~A.~Bennell}
\author[label1001]{Jean-Marie~Bourjolly}
\author[label6]{John~E.~Boylan}
\author[label7a]{Michèle~Breton}
\author[label9]{David~Canca}
\author[label7a,label7b]{Laurent~Charlin}
\author[label10]{Bo~Chen}
\author[label301]{Cihan~Tugrul~Cicek}
\author[label11a,label11b]{Louis~Anthony~Cox,~Jr}
\author[label12]{Christine~S.M.~Currie}
\author[label13]{Erik~Demeulemeester}
\author[label14]{Li~Ding}
\author[label15]{Stephen~M.~Disney}
\author[label16]{Matthias~Ehrgott}
\author[label17]{Martin~J.~Eppler}
\author[label1]{Güneş~Erdoğan}
\author[label18c,label18a,label18b]{Bernard~Fortz}
\author[label19a,label19b]{L.~Alberto~Franco}
\author[label20]{Jens~Frische}
\author[label21a,label21b]{Salvatore~Greco}
\author[label1002]{Amanda~J.~Gregory}
\author[label22]{Raimo~P.~Hämäläinen}
\author[label13]{Willy~Herroelen}
\author[label23]{Mike~Hewitt}
\author[label24]{Jan~Holmström}
\author[label25]{John~N.~Hooker}
\author[label26]{Tuğçe~Işık}
\author[label27]{Jill~Johnes}
\author[label28]{Bahar~Y.~Kara}
\author[label28]{Özlem~Karsu}
\author[label29]{Katherine~Kent}
\author[label31]{Charlotte~Köhler}
\author[label32]{Martin~Kunc}
\author[label33a,label33b]{Yong-Hong~Kuo}
\author[label16]{Adam~N.~Letchford}
\author[label34]{Janny~Leung}
\author[label19a]{Dong~Li}
\author[label35]{Haitao~Li}
\author[label201]{Judit~Lienert}
\author[label36]{Ivana~Ljubić}
\author[label37]{Andrea~Lodi}
\author[label38]{Sebastián~Lozano}
\author[label40]{Virginie~Lurkin}
\author[label41]{Silvano~Martello}
\author[label42]{Ian~G.~McHale}
\author[label43a,label43ab,label43b,label43c,label43f]{Gerald~Midgley}
\author[label44]{John~D.W.~Morecroft}
\author[label45]{Akshay~Mutha}
\author[label46]{Ceyda~Oğuz}
\author[label47]{Sanja~Petrovic}
\author[label61]{Ulrich~Pferschy}
\author[label48]{Harilaos~N.~Psaraftis}
\author[label49]{Sam~Rose}
\author[label50]{Lauri~Saarinen}
\author[label102]{Said~Salhi}
\author[label51]{Jing-Sheng~Song}
\author[label56]{Dimitrios~Sotiros}
\author[label52]{Kathryn~E.~Stecke}
\author[label20]{Arne~K.~Strauss}
\author[label46]{İstenç~Tarhan}
\author[label62]{Clemens~Thielen}
\author[label41]{Paolo~Toth}
\author[label58]{Tom~Van~Woensel}
\author[label8]{Greet~Vanden~Berghe}
\author[label1]{Christos~Vasilakis}
\author[label53]{Vikrant~Vaze}
\author[label101]{Daniele~Vigo}
\author[label54a,label54b]{Kai~Virtanen}
\author[label55]{Xun~Wang}
\author[label56]{Rafał~Weron}
\author[label57]{Leroy~White}
\author[label15]{Mike~Yearworth}
\author[label59]{E.~Alper~Yıldırım}
\author[label7a]{Georges~Zaccour}
\author[label60]{Xuying~Zhao}


\address[label1]{School of Management, University of Bath, Bath, UK}
\address[label1b]{Makridakis Open Forecasting Center, University of Nicosia, Nicosia, Cyprus}
\address[label7a]{Department of Decision Sciences, HEC Montréal, Montreal, Canada}
\address[label401]{Molde University College, Molde, Norway}
\address[label100]{Cranfield School of Management, Cranfield University, Cranfield, UK}
\address[label2]{Department of Management Sciences, University of Waterloo, Waterloo, Canada}
\address[label3]{ESSEC Business School in Paris, Department of Information Systems, Decision Sciences and Statistics, Cergy, France}
\address[label4]{H. Milton Stewart School of Industrial and Systems Engineering, Georgia Institute of Technology, Atlanta, USA}
\address[label5]{Centre for Decision Research, Leeds University Business School, University of Leeds, Leeds, United Kingdom}
\address[label1001]{Université du Québec à Montréal, Montreal, Canada}
\address[label6]{Centre for Marketing Analytics and Forecasting, Lancaster University Management School, Lancaster University, Lancaster, UK}
\address[label9]{School of Engineering, Department of Industrial Engineering and Management Science I, University of Seville, Seville, Spain}
\address[label7b]{Mila-Quebec AI Institute, Montreal, Canada}
\address[label10]{Warwick Business School, University of Warwick, Coventry, UK}
\address[label301]{Department of Industrial Engineering, Atilim University, Ankara, Turkey}
\address[label11a]{Department of Business Analytics, University of Colorado School of Business, Denver, USA}
\address[label11b]{Cox Associates, Denver, USA}
\address[label12]{Mathematical Sciences, University of Southampton, Southampton, UK}
\address[label13]{Faculty of Economics and Business, Research Center for Operations Management, KU Leuven, Leuven, Belgium}
\address[label14]{Durham University Business School, Durham University, Durham, UK}
\address[label15]{Center for Simulation, Analytics, and Modelling, University of Exeter Business School, Exeter, UK}
\address[label16]{Department of Management Science, Lancaster University Management School, Lancaster University, Lancaster, UK}
\address[label17]{University of St. Gallen, St. Gallen, Switzerland}
\address[label18c]{HEC - Management School of the University of Liège, Liège, Belgium}
\address[label18a]{Département d'Informatique, Université Libre de Bruxelles, Brussels, Belgium}
\address[label18b]{Inria Lille-Nord Europe, Villeneuve d'Ascq, France}
\address[label19a]{School of Business and Economics, Loughborough University, Loughborough, UK}
\address[label19b]{Universidad del Pacífico, Lima, Perú}
\address[label20]{WHU – Otto Beisheim School of Management, Vallendar, Germany}
\address[label21a]{Department of Economics and Business, University of Catania, Catania, Italy}
\address[label21b]{Portsmouth Business School, Centre of Operations Research and Logistics (CORL), University of Portsmouth, Portsmouth, UK}
\address[label1002]{Centre for Systems Studies, Faculty of Business, Law and Politics, University of Hull, Hull, UK}
\address[label22]{Systems Analysis Laboratory, Aalto University, Finland}
\address[label23]{Department of Information Systems and Supply Chain Management, Loyola University, Chicago, USA}
\address[label24]{Department of Industrial Engineering and Management, Aalto University, Aalto, Finland}
\address[label25]{Tepper School of Business, Carnegie Mellon University, Pittsburgh, USA}
\address[label26]{Department of Industrial Engineering, Clemson University, Clemson, USA}
\address[label27]{Department of Accounting, Finance and Economics, University of Huddersfield, Huddersfield, UK}
\address[label28]{Department of Industrial Engineering, Bilkent University, Ankara, Turkey}
\address[label29]{Office for National Statistics, Newport, UK}
\address[label31]{Europa Universität Viadrina, Frankfurt (Oder), Germany}
\address[label32]{Southampton Business School, University of Southampton, UK}
\address[label33a]{Department of Industrial and Manufacturing Systems Engineering, The University of Hong Kong, Pokfulam, Hong Kong, China}
\address[label33b]{HKU Musketeers Foundation Institute of Data Science, The University of Hong Kong, Pokfulam, Hong Kong, China}
\address[label34]{State Key Laboratory of Internet of Things for Smart City, Choi Kai Yau College, University of Macau, Taipa, Macau, China}
\address[label35]{Supply Chain \& Analytics Department, University of Missouri-St. Louis, St. Louis, USA}
\address[label201]{Eawag: Swiss Federal Institute of Aquatic Science and Technology, Switzerland}
\address[label36]{ESSEC Business School, Cergy-Pontoise, France}
\address[label37]{Jacobs Technion-Cornell Institute, Cornell Tech and Technion - IIT, USA}
\address[label38]{Department of Industrial Management, Escuela Superior de Ingenieros, University of Seville, Sevilla, Spain}
\address[label40]{HEC Lausanne Faculty of Business and Economics, University of Lausanne, Lausanne, Switzerland}
\address[label41]{DEI ``Guglielmo Marconi'', Alma Mater Studiorum Università di Bologna, Bologna, Italy}
\address[label42]{Centre for Sports Business, University of Liverpool Management School, UK}
\address[label43a]{Centre for Systems Studies, Faculty of Business, Law and Politics, University of Hull, Hull, UK}
\address[label43ab]{Birmingham Leadership Institute, University of Birmingham, Birmingham, UK}
\address[label43b]{Department of Informatics, Faculty of Technology, Linnaeus University, Växjö, Sweden}
\address[label43c]{School of Innovation, Design and Engineering, Mälardalen University, Eskilstuna, Sweden}
\address[label43f]{School of Agriculture and Food Sciences, University of Queensland, Brisbane, Queensland, Australia}
\address[label44]{London Business School, London, UK}
\address[label45]{Grossman School of Business, University of Vermont, Burlington, Vermont, USA}
\address[label46]{Department of Industrial Engineering, Koç University, İstanbul, Turkey}
\address[label47]{Nottingham University Business School, University of Nottingham, Nottingham, UK}
\address[label61]{Department of Operations and Information Systems, University of Graz, Graz, Austria}
\address[label48]{Department of Technology, Management and Economics, Technical University of Denmark, Lyngby, Denmark}
\address[label49]{Department for Transport, London, UK}
\address[label50]{Department of Industrial Engineering and Management, Aalto University, Aalto, Finland}
\address[label102]{Centre for Logistics and Heuristic optimisation, Kent Business School, University of Kent, Kent, UK}
\address[label51]{The Fuqua School of Business, Duke University, Durham, North Carolina, USA}
\address[label56]{Department of Operations Research and Business Intelligence, Wrocław University of Science and Technology, Wrocław, Poland}
\address[label52]{Naveen Jindal School of Management, The University of Texas at Dallas, Richardson, Texas, USA}
\address[label62]{TUM Campus Straubing for Biotechnology and Sustainability, Technical University of Munich, Straubing, Germany}
\address[label58]{School of Industrial Engineering, Eindhoven University of Technology, Eindhoven, The Netherlands}
\address[label8]{Department of Computer Science, KU Leuven, Gent, Belgium}
\address[label53]{Thayer School of Engineering, Dartmouth College, Hanover, USA}
\address[label101]{Department of Electrical, Electronic, and Information Engineering ``G. Marconi'' and CIRI-ICT, University of Bologna, Bologna, Italy}
\address[label54a]{Department of Mathematics and Systems Analysis, Aalto University, Helsinki, Finland}
\address[label54b]{Department of Military Technology, National Defence University, Helsinki, Finland}
\address[label55]{Logistics Systems Dynamics Group, Cardiff Business School, Cardiff University, UK}
\address[label57]{Management Department, University of Exeter Business School, Exeter, UK}

\address[label59]{School of Mathematics, The University of Edinburgh, Edinburgh, UK}
\address[label60]{Mendoza College of Business, University of Notre Dame, Notre Dame, Indiana, USA}


\end{frontmatter}

\epigraph{Operations research is neither a method nor a technique; it is or is becoming a science and as such is defined by a combination of the phenomena it studies.}{\cite{Ackoff1956-qh}}

\section*{Abstract}
\label{sec:Abstract}

\noindent Throughout its history, Operational Research has evolved to include a variety of methods, models and algorithms that have been applied to a diverse and wide range of contexts. This encyclopedic  article consists of two main sections: methods and applications. The first aims to summarise the up-to-date knowledge and provide an overview of the state-of-the-art methods and key developments in the various subdomains of the field. The second offers a wide-ranging list of areas where Operational Research has been applied. The article is meant to be read in a nonlinear fashion. It should be used as a point of reference or first-port-of-call for a diverse pool of readers: academics, researchers, students, and practitioners. The entries within the methods and applications sections are presented in alphabetical order.

\noindent The authors dedicate this paper to the 2023 Turkey/Syria earthquake victims. We sincerely hope that advances in OR will play a role towards minimising the pain and suffering caused by this and future catastrophes.

\noindent \textbf{Keywords:} review; encyclopedia; theory; practice; principles; optimisation; programming; systems; simulation; decision making; models.

\clearpage
\setcounter{tocdepth}{3}
\small
\tableofcontents
\normalsize

\clearpage

\section[Introduction (Gilbert Laporte)]{Introduction\protect\footnote{This subsection was written by Gilbert Laporte.}}
\label{sec:introduction}

The year 2024 marks the 75\textsuperscript{th} anniversary of the \textit{Journal of the Operational Research Society}, formerly known as \textit{Operational Research Quarterly}. It is the oldest Operational Research (OR) journal worldwide. On this occasion, my colleague Fotios Petropoulos from University of Bath proposed to the editors of the journal to edit an encyclopedic article on the state of the art in OR. Together, we identified the main methodological and application areas to be covered, based on topics included in the major OR journals and conferences. We also identified potential authors who responded enthusiastically and whom we thank wholeheartedly for their contributions.

Modern OR originated in the United Kingdom during World War II as a need to support the operations of early radar-detecting systems and was later applied to other operations \citep{McCloskey1987-gw_GL}. However, one could argue that it precedes this period in history since it is partly rooted in several mathematical fields such as probability theory and statistics, calculus, and linear algebra, developed much earlier. For example, the Fourier-Motzkin elimination method \citep{Fourier1826b-hx_GL,Fourier1826a-hx_GL} constitutes the main basis of linear programming.  Queueing theory, which plays a central role in telecommunications and computing, already existed as a distinct field of study since the early 20\textsuperscript{th} century \citep{Erlang1090-er_GL}, and other concepts, such as the economic order quantity \citep{harris1913many_JSS} were developed more than one century ago. Interestingly, while many recent advances in OR are rooted in theoretical or algorithmic concepts, we are now witnessing a return to the practical roots of OR through the development of new disciplines such as business analytics.

After the war ended, several industrial applications of OR arose, particularly in the manufacturing and mining sectors which were then going through a renaissance. The transportation sector is without doubt the field that has most benefited from OR, mostly since the 1960s. The aviation, rail, and e-commerce industries could simply not operate at their current scale without the support of massive data analysis and sophisticated optimisation techniques. The application of OR to maritime transportation is more recent, but it is fast gaining in importance. Other areas that are less visible, such as telecommunications, also deeply depend on OR. The success of OR in these fields is partly explained by their network structures which make them amenable to systematic analysis and treatment through mathematical optimisation techniques. In the same vein, OR also plays a major role in various branches of logistics and project management, such as facility location, forecasting, inventory planning, scheduling, and supply chain management. 

The public sector and service industries also benefit greatly from OR. Healthcare is the first area that comes to mind because of its very large scale and complexity. Decision making in healthcare is more decentralised than in transportation and manufacturing, for example, and the human issues involved in this sector add a layer of complexity. OR methodologies have also been applied to diverse areas such as education, sports management, natural resources, environment and sustainability, political districting, safety and security, energy, finance and insurance, revenue management, auctions and bidding, and disaster relief, most of which are covered in this article.

Among OR methodologies, mathematical programming occupies a central place. The simplex method for linear programming, conceived by Dantzig in 1947 but apparently first published later \citep{Dantzig1951-qp_GL}, is arguably the single most significant development in this area. Over time, linear programming has branched out into several fields such as nonlinear programming, mixed integer programming, network optimisation, combinatorial optimisation, and stochastic programming. The techniques most frequently employed for the exact solution of mathematical programs are based on branch-and-bound, branch-and-cut, branch-and-price (column generation), and dynamic programming. Game theory and data envelopment analysis are firmly rooted in mathematical programming. Control theory is also part of continuous mathematical optimisation and relies heavily on differential equations.

Complexity theory is fundamental in optimisation. Most problems arising in combinatorial optimisation are $\mathcal{NP}$-hard and typically require the application of heuristics for their solution. Much progress has been made in the past 40 years or so in the development of metaheuristics based on local search, genetic search, and various hybridisation schemes. Many problems in fields such as vehicle routing, location analysis, cutting and packing, set covering, and set partitioning can now be solved to near optimality for realistic sizes by means of modern heuristics. A recent trend is the use of open-source software which not only helps disseminate research results, but also contributes to ensuring their accuracy, reproducibility and adoption.

Several modelling paradigms such as systems thinking and systems dynamics approach problems from a high-level perspective, examining the inter-relationships between multiple elements. Complex systems can often be analysed through simulation, which is also commonly used to assess the performance of heuristics. Decision analysis provides a useful framework for structuring and solving complex problems involving soft and hard criteria, behavioural OR, stochasticity, and dynamism. Recently, issues related to ethics and fairness have come to play an increasing role in decision making.

Because the various topics of this review paper are listed in alphabetical order, the subsection on “Artificial intelligence, machine learning and data science” comes first, but this topic constitutes one of the latest developments in the field. It holds great potential for the future and is likely to reshape parts of the OR discipline. Already, machine learning-based heuristics are competitive for the solution of some hard problems.

This paper begins with a quote from Russell L. Ackoff who has been a pioneer of OR. In 1979, he published in this journal two articles \citep{Ackoff1979-zx_GL,Ackoff1979-vt_GL} that presented a rather pessimistic view of our discipline. The author complained about the lack of communications between academics and practitioners, and about the fact that some OR curricula in universities did not sufficiently prepare students for practice, which is still true to some extent. One of his two articles is entitled “The Future of Operational Research is Past”, which may be perceived as an overreaction to this diagnosis. In my view, the present article provides clear evidence to the contrary. Soon after the publication of the two Ackoff papers, we witnessed the development of micro-computing, the Internet and the World Wide Web. It has become much easier for researchers in our community to access information, software and computing facilities, and for practitioners to access and use our research results. We are now fortunate to have access to sophisticated open-source software, data bases, bibliographic sources, editing and visualisation tools, and communication facilities. Our field is richer than it has ever been, both in terms of theory and applications. It is constantly evolving in interaction with other disciplines, and it is clearly alive and well and has a promising future. 

\clearpage

\section{Methods}
\label{sec:methods}

\subsection[Artificial intelligence, machine learning, and data science (Laurent~Charlin \& Andrea~Lodi)]{Artificial intelligence, machine learning, and data science\protect\footnote{This subsection was written by Laurent~Charlin and Andrea~Lodi.}}
\label{sec:Artificial_intelligence_machine_learning_data_science}
Machine learning (ML) comprises  techniques for modelling predictive tasks, i.e.\ tasks that involve the prediction of an unknown quantity from other observed quantities. Ideas of learning in an artificial system and the term machine learning were first discussed in the 1950s \citep{Samuel_IBM_5392560_LCAL} and their development and popularity have seen enormous growth over the last two decades in part due to the availability of large-scale datasets and increased computational resources to model them.

\cite{mitchell1997machine_LCAL} provides this concrete definition of machine learning ``A computer program is said to learn from experience E with respect to some class of tasks T, and performance measure P, if its performance at tasks in T, as measured by P, improves with experience E''. The program is a model or a function and its experience E is the type of data it has access to. There are three types of experiences supervised, unsupervised, and reinforcement learning. The performance measure (P) allows for model evaluation and comparison including model selection.

Supervised learning is an experience where a model aims at predicting one or more unobserved target (dependent) variables given observed ibackprout (independent) variables. In other words, a supervised model is a function that map inputs to outputs. The process of solving a supervised problem involves first learning a model, that is adjusting its parameters using a training dataset with both input and target variables. The training set is drawn IID (independently and identically distributed) from an underlying distribution over inputs and targets. Once trained, the model can provide target predictions for new unseen samples from the same distribution. The most common tasks in supervised learning are regression (real dependent variable) and classification (categorical dependent variable). Evaluating a supervised system is usually performed using held-out data referred to as the test data while held-out validation data is used for model development and selection using procedures such as $k$-fold cross-validation.

Supervised models can be dichotomised into linear and nonlinear models. Linear models perform a linear mapping from inputs to outputs (e.g., linear regression). Machine learning mostly investigates nonlinear supervised models including deep neural network (DNN) models \citep{Goodfellow-et-al-2016_LCAL}. DNNs are composed of a succession of parametrised nonlinear transformations called layers and each layer contains a set of transformations called neurons. Layers successively transform an input datum into a target. The parameters of the layers are adjusted to iteratively obtain better predictions using a procedure called backpropagation, a form of gradient descent \cite[\S 6.5]{Goodfellow-et-al-2016_LCAL}. DNNs are state-of-the-art methods for many large-scale non-structured datasets across domains (see also \S\ref{sec:Power_markets_and_systems}). DNNs can be adapted to different sizes of inputs and targets as well as variable types. They can also be specialised for specific types of data. Recurrent neural networks (RNNs) are auto-regressive models for sequential data \citep{Rumelhart86_LCAL}. The sequential data are tokenised and an RNN transforms each token sequentially along with a transformation of the previous tokens. Convolutional neural networks (CNNs) are specialised networks for modelling data that is arranged on a grid \citep[e.g., an image][]{LeCun89_LCAL}. Their layers contain a convolution operation between an input and a parameterised filter followed by a nonlinear transformation, and a pooling operation. Each layer processes data locally and so requires fewer parameters compared to vanilla DNNs. As a result, CNNs can model higher-dimensional data. Graphical neural networks (GNNs) are specialised architectures for modelling graph data \citep[e.g., a social network;][]{Scarselli09_LCAL}. In GNNs, the data are transformed by following the topology of the graph. 
Last, attention layers dynamically combine their inputs (tokens) based on their values. Transformer models use successions of attention and feed-forward layers to model sequential input and output data \citep{NIPS2017_3f5ee243_LCAL}. Transformers are more efficient to train than RNNs and can be trained on internet-scale data given enormous computational power. The availability of such broad datasets especially in the text and image domains has given rise to a class of very-large-scale models (also referred to as foundation models) that display an ability to adapt to and obtain high performance across a diversity of downstream supervised tasks \citep{Bommasani2021FoundationModels}

Last, attention is a mechanism that considers data to be unordered and uses transformations dynamically. Transformers are models based on attention. They provide more efficient training than RNNs for very large-scale datasets \citep{NIPS2017_3f5ee243_LCAL}.  

Neural networks currently outperform other methods when learning from unstructured data (e.g., images and text). For tabular data, data that is naturally encoded in a table and that has heterogeneous features \citep{NEURIPS2022_Varoquaux_LCAL}, best-performing methods use ideas first proposed in tree-based classifiers, bagging, and boosting. They include random forests \citep{Breiman01_LCAL}, XGBoost \citep{Chen16_LCAL} which both scale to large-scale datasets as well as kernel methods including support vector machines \citep[SVMs see, e.g.,][]{Schlkopf18_LCAL} and probabilistic Gaussian Processes  \citep[GPs see, e.g.,][]{Rasmussen05_LCAL}. These methods are used across regression and classification tasks.

In unsupervised learning, the second type of experience, the data consist of independent variables (features or covariates) alone. The aim of unsupervised learning is to model the structure of the data to better understand their properties. As a result, evaluating an unsupervised model is often task and application-dependant \cite[\S 1.3.4]{pml1Book_LCAL}. The prototypical unsupervised-learning task is clustering. It involves learning a function that groups similar data together according to a similarity measure and desiderata often expressed as an objective function. Several standard algorithms divided into hierarchical and non-hierarchical methods exist. The former uses the similarity between all pairs of data and finds a hierarchy of clustering solutions with a different number of clusters using either a bottom-up or top-down approach. Agglomerative clustering is a standard hierarchical approach. Non-hierarchical methods tend to be more computationally efficient in terms of dataset size. 
For example, $K$-means clustering is a well-known non-hierarchical method that finds a single solution using $K$ clusters \citep{MacQueen1967_LCAL}. Other unsupervised learning tasks include dimensionality reduction for example for visualisation or to prepare data for further analysis. Density modelling is another unsupervised task where a probabilistic model learns to assign a probability to each datum \cite[\S 1.3]{pml1Book_LCAL}. Probabilistic models can be used to learn the hidden structure in large quantities of data \citep[e.g.,][]{hoffman13_LCAL}. Further, probabilistic models are also used to generate high-dimensional data (e.g., images of human faces or English text) with high fidelity \citep{Karas19_LCAL} and often referred in this context as generative models. Large Language Models are examples of such generative models \citep{Bommasani2021FoundationModels}.

Reinforcement learning (RL) is the third type of experience. RL models collect their own data by executing actions in their environment to maximise their reward. RL is a sequential decision-making task and is formalised using Markov decision processes (MDPs) \cite[\S 3.8]{sutton2018reinforcement_LCAL}. An MDP encodes a set of states, available actions, distribution over next states given current states and action, a reward function, and a discount factor. Partially observable MDPs (or POMDPs) extend the formalism to environments where the exact current state is unknown \citep{KAELBLING199899_LCAL}. In RL, an agent's objective is to learn a policy, a distribution over actions for each state in an environment. Tasks are defined by rewards attached to different states. Exact and approximate methods exist for solving RL problems. Whereas exact solutions are appropriate for smaller tabular problems only, deep neural networks are widely used for solving larger-scale problems that require approximate solutions yielding a set of techniques known as deep reinforcement learning \citep{Mnih_DRL_2015_LCAL}. An RL agent can also learn to imitate an expert either by learning a mapping from states/observations to actions as in supervised learning \citep[a technique known as imitation learning; for a survey, see][]{Hussein17_LCAL} or by trying to learn the expert's reward function \citep[inverse reinforcement learning][]{Russell98_LCAL}.  

In addition to learning models for solving prediction tasks using one of the three experiences above, machine learning also studies methods for enabling the reuse of information learned from one or multiple datasets and environments to other similar ones. Representation learning studies how to learn such reusable information and it can use both supervised and unsupervised experiences \cite[\S32]{pml2Book_LCAL}. When using a deep learning model, a representation is obtained after one or more layer transformations of the data. Representation learning is used in a variety of situations including for transfer learning tasks, where a trained model is reused to solve a different supervised task \cite[for a survey, see][]{Zhuang_TransferLearning21_LCAL}. 

In the last decade, machine learning models have achieved high performance on a variety of tasks including perceptual ones (e.g.\ recognising objects in images and words from speech) as well as natural language processing ones thereby becoming a core component of artificial intelligence (AI) methods. The goal of AI methods is to develop intelligent systems. Some of these advances shine a bright light on the ethical aspects of machine learning techniques and are active areas of study \citep[see, e.g.,][]{Dignum19_LCAL,barocas19_LCAL}. Another area of active study is explainability \citep{XAI_NIST_21}. Some of the most effective ML tools make predictions and recommendations that are hard to explain to users (for example when neural networks are employed). Clearly, lack of explainability slows down ML use in those contexts where decisions made due to those predictions and recommendations are life changing and involve a human in the loop, healthcare (applying a treatment), finance (refusing a mortgage), or justice (granting parole) to mention a few. So, explainability is currently one of the most crucial challenges for ML and AI and, at the same time, a tremendous opportunity for their wider applicability. 

Further, advances in machine learning alongside statistics, data management, and data processing, as well as the wider availability of datasets from a variety of domains have led to the popularisation and development of data science (DS), a  discipline whose goal is to extract insights and knowledge from these data. DS uses statistics and machine-learning techniques for inference and prediction, but it also aims at enabling and systematising the analysis of large quantities of data. As such, it includes components of data management, visualisation, as well as the design of (efficient) data processing algorithms \citep{grus2019data_LCAL}.
 
\subsubsection*{Resources} 
\citet{pml1Book_LCAL} provides a thorough introduction to the field following a probabilistic approach and its sequel \citep{pml2Book_LCAL} introduces advanced topics. \cite{Goodfellow-et-al-2016_LCAL} provide a self-contained introduction to the field of deep learning \citep[the field evolves rapidly and more advanced topics are covered through recent papers and in][]{pml2Book_LCAL}. 
Open-source software packages in Python and other languages are essential. They include data-wrangling libraries such as pandas \citep{mckinney-proc-scipy-2010_LCAL} and plotting ones such as matplotlib \citep{Hunter:2007_LCAL}. The library scikit-learn \citep{scikit-learn_LCAL} in Python offers an extensive API that includes data processing, a toolbox of standard supervised and unsupervised models, and evaluation routines. For deep learning, PyToch \citep{PyToch_NEURIPS2019_9015_LCAL} and TensorFlow \citep{tensorflow2015-whitepaper_LCAL} are the standard. 

\subsubsection*{Learning for combinatorial optimisation} 
The impressive success of machine learning in the last decade made it natural to explore its use in many scientific disciplines, such as drug discovery and material sciences. Combinatorial optimisation (CO; \S\ref{sec:Combinatorial_optimisation}) is no exception to this trend and we have witnessed an intense exploration (or, better, revival) of the use of machine learning for CO. Two lines of work have strongly emerged. On the one side, ML has been used to learn crucial decisions within CO algorithms and solvers. This includes imitating an algorithmic expert that is computationally expensive like in the case of strong branching for branch and bound, the single application that has attracted the largest amount of interest \citep{lodi2017learning_LCAL,gasse2019exact_LCAL}. The interested reader is referred to two recent surveys \citep{bengio2021machine_JLYHK,cappart2021combinatorial_LCAL}, the latter highlighting the relevance of GNNs for effective CO representation. On the other side, ML has been used end to end, i.e., for solving CO problems directly or leveraging ML to devise hybrid methods for CO. The area is surveyed in \cite{fioretto2021_LCAL}.

\subsection[Behavioural OR (L.~Alberto~Franco \& Raimo~P.~Hämäläinen)]{Behavioural OR\protect\footnote{This subsection was written by L.~Alberto~Franco and Raimo~P.~Hämäläinen.}}
\label{sec:Behavioural_OR}

Behavioural OR (BOR) is concerned with the study of human behaviour in OR-supported settings. Specifically, BOR examines how the behaviour of individuals affects, or is affected by, an OR-supported intervention\footnote{We use the term ‘intervention’ to describe a structured process comprised of designed OR-related activities such as, for example, modelling, model use, data collection, interviews, meetings, workshops, and presentations.}. The individuals of interest are those who, acting in isolation or as part of a team, design, implement and engage with OR in practice. These individuals include OR practitioners playing specific intervention roles (e.g., modellers, facilitators, consultants), and other individuals with varying interests and stakes in the intervention (e.g., users, clients, domain experts, sponsors). 

A concern with the behavioural aspects of the OR profession can be traced back to past debates in the 1960s, 1970s and 1980s \citep{Churchman1970-hv_GM,Dutton1964-wf_AFRH,Jackson1989-fk_AFRH}. Although these debates dwindled down in subsequent years, the emergence of BOR as a field of study represents a return to these earlier concerns \citep{Franco2016-qk_MYLW,Hamalainen2013-zr_AFRH}. What motivates this resurgence is the recognition that the successful deployment of OR in practice relies heavily on our understanding of human behaviour. For example, overconfidence, competing interests, and the willingness to expend effort in searching, sharing, and processing information are three behavioural issues that can negatively affect the success of OR activities. Attention to behavioural issues has been central in disciplines such as economics, psychology and sociology for decades, and BOR studies draw heavily from these reference disciplines \citep{Franco2021-sm_JL}.

It is important to distinguish between the specific focus of BOR and the broader focus of behavioural modelling. The creation of models that capture human behaviour has a long tradition within OR, but it is not necessarily concerned with the study of human behaviour in OR-supported settings. For example, in the last 20 years operational researchers have produced an increasing number of robust analytical models that describe behaviour in, and predict its impact on, operations management settings \citep{Cui2018-xl_AFRH,Donohue2020-mw_AFRH,Loch2007-ij_AFRH}. Operational researchers also have produced simulation models that capture human behaviour within a system with different levels of complexity. For example, systems dynamics models incorporate high-level variables representing average behaviour \citep[][\S\ref{sec:Systems_dynamics}]{Morecroft2015-vr_MCJM,Sterman2000_SMD}, and discrete event simulation models capture human processes controlled by simple behavioural rules \citep[][\S\ref{sec:Simulation}]{Brailsford2003-rs_AFRH,Robinson2014-bd_AFRH}. More complex agent-based simulation models represent behaviour as emergent from the interactions of agents with particular behavioural attributes \citep[][\S\ref{sec:Simulation}]{Sonnessa2017-la_AFRH,Utomo2018-uy_AFRH}. Overall, behavioural modelling within the OR field is concerned with examining human behaviour in a system of interest in order to improve that system\footnote{To date, the practice of studying OR-supported intervention as a system of interest has been somewhat overlooked by behavioural modellers, with the notable exception of behavioural forecasting researchers; see \S\ref{sec:forecasting}, and also \cite{Arvan2019-xy_AFRH} and \cite{Petropoulos2016-uk_FP}.}. In contrast, BOR takes an OR-supported intervention as the core system of interest where human behaviour is examined. The ultimate goal of BOR is to generate an improved understanding of the behavioural dimension of OR practice, and use this understanding to design and implement better OR-supported interventions. 

Another important distinction worth stating is that between BOR and Soft OR. At first glance, this distinction may seem unnecessary as BOR is a field of study within OR, while Soft OR refers to a specific family of problem structuring approaches (\S\ref{sec:Soft_OR_and_problem_structuring_methods}). Soft OR approaches have been developed to help groups reach agreements on problem structure and, often, appropriate responses to a problem of concern \citep{Franco2022-mn_AFRH,Rosenhead2001-ai_JL}. However, while Soft OR intervention design and implementation typically require the consideration of behavioural issues, this is not the same as choosing human behaviour in a Soft OR intervention context as the unit of analysis. Of course, a study with such a focus would certainly fall within the BOR remit \citep[e.g.,][]{Tavella2021-mj_AFRH}. But note that BOR is also concerned with the study of human behaviour in other OR-supported settings, such as those involving the use of ‘hard’ and ‘mixed-method’ OR approaches.

Studies of behaviour in OR-supported settings assume implicitly or explicitly that human behaviour is either influenced by cognitive and external factors, or is in itself an influencing factor \citep{Franco2021-sm_JL}.  In the first case, observed individual and collective action is taken to be guided by cognitive structures (e.g., personality traits, cognitive styles) manifested during OR-related activity -- behaviour is \textit{influenced}. In contrast, the second case assumes that individuals and collectives are responsible for determining how OR-related activity will unfold -- behaviour is \textit{influencing}. This raises the practical possibility that the same OR methodology, technique, or model could be used in distinctive ways by various individuals or groups according to their cognitive orientations, goals and interests \citep{Alberto_Franco2013-yw_MYLW}. Whilst behaviour in practice is likely to lie somewhere between the influenced and influencing assumptions, BOR studies tend to foreground one of the extremes as the focus, while backgrounding the other.

BOR studies can adopt three different research methodologies to examine behaviour: variance, process, and modelling.  A \textit{variance methodology} uses variables that represent the important aspects or attributes of the OR-supported activity being examined. Variance explanations of behavioural-related phenomena take the form of causal statements captured in a theoretically-informed research model that incorporates these variables (e.g., A causes B, which causes C). The research model is then tested with data generated by the activity, and the research findings are assessed in terms of their generality \citep{Poole2004-cq_AFRH}. Adopting a variance research methodology typically requires the implementation of experimental, quasi-experimental, or survey research designs\footnote{It should be noted that a variance approach could also be implemented through field research designs where pre and post intervention measures of key variables are used to assess changes in behaviour or surrogates of behaviour. Studies adopting this approach are common in the System Dynamics field \citep[see, for example,][]{Scott2013-fg_AFRH}.}. This involves careful selection of independent variables, which might be either manipulated or left untreated, and of dependent variables that act as surrogates for specific behaviours. Once information about all variables is collected, data is quantitatively analysed using a wide range of variance-based methods (e.g., analysis of variance, regression, structural equation modelling).

Behavioural studies that use a variance research methodology can produce a good picture of the generative mechanisms underpinning behavioural processes if they test hypotheses about those mechanisms. For example, variance studies in BOR have examined the impact of individual differences in cognitive motivation and cognitive style on the conduct of OR-supported activity \citep{Fasolo2014-mp_AFRH,Franco2016-yn_AFRH,Lu2001-ow_AFRH}. There is also a long tradition of testing the behavioural effects of reconfiguring different aspects of OR-supported activity such as varying model or information displays \citep{Bell1995-lw_AFRH,Gettinger2013-mi_AFRH}, and preference elicitation procedures \citep{Cavallo2019-sj_AFRH,Hamalainen2016-gr_AFRH,Poyhonen2001-dw_AFRH,Von_Nitzsch1993-hh_AFRH}.

A \textit{process methodology} is used to examine OR-supported activity as a series of events that bring about or lead to some behaviour-related outcome. Specifically, it considers as the unit of analysis an evolving individual or group whose behaviour is led by, or leading, the occurrence of events \citep{Poole2004-cq_AFRH}. Process explanations take the form of theoretical narratives that account for how event dynamics lead to a final outcome \citep{Poole2007-cg_AFRH}. These narratives are often derived from observation, but it is also possible to use an established narrative (e.g., a theory) to guide observation that further specifies the narrative. 

Diverse and eclectic research designs are used to implement a process research methodology. Central to these designs is the task of identifying or reconstructing the process through the analysis of events taking place over time. For example, there is an important stream of BOR studies that examines the process of building models by experts and novices \citep{Tako2015-dk_AFRH,Tako2010-ua_AFRH,Waisel2008-uw_AFRH,Willemain1995-qz_AFRH,Willemain2007-hj_AFRH}. There is also an increasing interest to use process methodologies to take a closer look at actual behaviour in OR-supported settings both, before, during and after OR-related activity is undertaken \citep{Franco2018-nk_AFRH,Kaki2019-bo_AFRH,Velez-Castiblanco2016-kh_AFRH,White2016-sv_MYLW}.

The variance and process approaches may seem opposite to each other, but instead they should be seen as complementary \citep{Franco2021-sm_JL,Van_de_Ven2005-ug_AFRH}. BOR studies using a variance research methodology can explore and test the mechanisms that drive process explanations of behaviour, while BOR studies adopting a process research methodology can explore and test the narratives that ground variance explanations of behaviour. One way of combining a variance and process approach within a single BOR study is by adopting \textit{modelling} as a research methodology. A modelling approach would create models that capture the mechanisms that generate a process of interest such as, for example, trust on an OR-derived solution, and the model can be run to generate the characteristics of that process. Model parameters and structure can then be varied systematically to enable variance-based comparisons of trust levels. Furthermore, the trajectory of trust levels over time can be used to gain insights into the nature of the trust development process.  As already mentioned, there is a long behavioural modelling tradition within OR but, as far as we know, its potential as a research methodology tool to specifically examine behaviour in OR-supported settings is yet to be realised. 

In sum, the variance, process and modelling methodologies offer rich possibilities for the study of human behaviour in OR-supported settings. Which is best for a particular study will depend on the types of question being addressed by BOR researchers, their assumptions about human behaviour, and the data they have access to. Ultimately, a thorough understanding of behaviour in OR-supported settings is likely to require all three research methodologies. 

For  a detailed review of BOR studies the reader is referred to \cite{Franco2021-sm_JL}. A review of behavioural studies in the context of OR in health has been written by \cite{Kunc2020-fn_MCJM}. There are also two collections edited by \cite{Kunc2016-gy_AFRH} and \cite{White2020-fw_AFRH}. The \textit{European Journal of Operational Research} published a feature cluster on BOR edited by \cite{Franco2016-ns_AFRH}.  Finally, BOR-related news and events can be found on the sites of the European Working Group on Behavioural OR\footnote{\url{https://www.euro-online.org/websites/bor/}}, and the UK BOR Special Interest Group\footnote{\url{https://www.theorsociety.com/get-involved/society-groups/special-interest-groups-and-networks/behavioural-or/}}. 

\subsection[Business analytics (John~E.~Boylan)]{Business analytics\protect\footnote{This subsection was written by John~E.~Boylan.}}
\label{sec:Business_analytics}

Business Analytics has its origins in practice, rather than theory, as illustrated by some of the earliest publications on the subject \citep[e.g.,][]{Kohavi2002-rl_JEB}. Senior executives began to realise the importance of analytics in the first decade of the new millennium because of the ready availability of large amounts of data, the maturity of business performance management, the emergence of self-service analytics and business intelligence, and the declining cost of computing power, data storage and bandwidth \citep{Acito2014-my_JEB}.

\cite{Davenport2007-jn_JEB} gave examples of companies becoming ‘analytical competitors’ by using analytics to support distinctive organisational capabilities. To achieve this level of maturity, it was argued that analytics needs to become a strategic competency. In the 1990s, \cite{Fildes1997-bo_JEB} reported on the closure or dispersal of Operational Research groups. \cite{Davenport2010-vu_JEB} reflected a reversal of that trend, by focusing on how analytical talent can be organised as an internal resource. They suggested that there are four categories of people to be considered when finding, developing and managing analysts: champions, professionals, semi-professionals and amateurs. In 2012/13, the Institute for Management Science and Operations Research (INFORMS) introduced the Certified Analytics Professional program and examination. This covers the broad spectrum of skills required of analytics professionals, including business problem framing, analytics problem framing, data (handling), methodology selection, model building, deployment and lifecycle management \citep{INFORMS2022}.

The development of talent is just one of the prerequisites for Business Analytics to create value. \cite{Vidgen2017-vp_JEB} recommended ‘coevolutionary change’, aligning their analytics strategy with their strategies for Information and Communications Technology, human resources and the whole business. This helps to ensure that the necessary data assets are available, the right culture is developed to build data and analytics skills, and that there is alignment with the business strategy for value creation. \cite{Hindle2018-zo_JEB} proposed a Business Analytics Methodology based on four activities, namely problem situation structuring, business model mapping, business analytics leverage and analytics implementation. They advocated a soft OR approach, Soft Systems Methodology \citep{Checkland2006-tv_JEB}, to support structuring and mapping activities.

Many definitions of Business Analytics have been proposed; for a review of early definitions, see \cite{Holsapple2014-gf_JEB}. According to \cite{Davenport2013-xh_JEB}, “\textit{By analytics we mean the extensive use of data, statistical and quantitative analysis, explanatory and predictive models, and fact-based management to drive decisions and actions}” (p. 7). \cite{Mortenson2015-ak_JEB} suggested that analytics is at the intersection of quantitative methods, technologies and decision making. \cite{Rose2016-lq_JEB} considered analytics as the union of Data Science (which is data centric) and Operational Research (which is problem centric). \cite{Power2018-kx_JEB} proposed the following definition: “\textit{Business Analytics is a systematic thinking process that applies qualitative, quantitative and statistical computational tools and methods to analyse data, gain insights, inform and support decision-making}". \cite{Delen2018-hf_JEB} pointed out that, although analytics includes analysis, it also involves synthesis and subsequent implementation. These broad perspectives, emphasising synthesis as well as analysis, and qualitative as well as quantitative approaches, are consistent with earlier writings on the use of a broad range of methods in Management Science \citep[e.g.,][]{Mingers1997-nd_MYLW,Pidd2009-ny_JEB}.

Business Analytics can be viewed from different orientations. From a methodological viewpoint, the subject covers descriptive, predictive and prescriptive methods \citep{Lustig2010-ty_JEB}. These three  categories are sometimes extended to four, with a distinction being drawn between ‘descriptive’ and ’diagnostic’ analytics, following the Gartner analytics ascendancy model \citep{Maoz2013-zz_JEB}. \cite{Lepenioti2020-mq_JEB} argue that it is preferable to maintain the threefold categorisation to ensure consistency, with each category addressing both  ‘What?’ and ‘Why’ questions. (Descriptive: ‘What happened?’, ‘Why did it happen?’; Predictive: ‘What will happen?’, ‘Why will it happen?’; Prescriptive: ‘What should I do to make it happen?’, ‘Why should I make it happen?’).  For detailed literature reviews on descriptive, predictive and prescriptive analytics, the reader is directed to \cite{Duan2015-bh_JEB}, \cite{Lu2017-oe_JEB}, and \cite{Lepenioti2020-mq_JEB}, respectively.

From a technological viewpoint, Business Analytics is facilitated by the integration of transactional data with big data streaming from social media platforms and the Internet of Things into a unified analytics system \citep{Shi2018-eb_JEB}. These authors suggest that this integration can be achieved in two stages, starting with integration of traditional Enterprise Resource Planning (ERP) and big data, and proceeding to integration of big-data ERP with Business Analytics. \cite{Ruivo2020-kq_JEB} reported that analytics ranked second in extended ERP capabilities (behind collaboration) according to the views of 20 experts engaged in a Delphi study. \cite{Romero2022-nw_JEB} suggested that cloud-based big data analytics software will not provide competitive advantage to firms that have not installed a large ERP system, although it will ensure that they do not lag further behind their sector-leading competitors.

From an ethical viewpoint, Business Analytics faces a number of challenges. \cite{Davenport2010-vu_JEB} recognised that issues of data privacy can be difficult to address, especially if an organisation operates in a wide range of territories or industries. \cite{Ram_Mohan_Rao2018-jb_JEB} summarised major privacy threats in data analytics, namely surveillance, disclosure, discrimination, and personal embarrassment and abuse, and reviewed privacy preservation methods, including randomisation and cryptographic techniques. A further ethical issue is that AI algorithms are likely to replicate and reinforce existing social biases \citep{ONeil2016-rs_JEB}. Such algorithmic bias is said to occur when the outputs of an algorithm benefit or disadvantage certain individuals or groups more than others without a justified reason. \cite{Kordzadeh2022-yp_JEB} reviewed the literature on algorithmic bias and showed that most studies had examined the issue from a conceptual standpoint, with only a limited number of empirical studies. Similarly, \cite{Vidgen2020-vy_JEB} reviewed papers on ethics in Business Analytics and found that most were at the level of guiding principles and frameworks, with little of direct applicability for the practitioner. Their case study demonstrated how ethical principles (utility, rights, justice, common good and virtue) can be embedded in analytics development. For further discussions on ethics and OR, the reader is referred to \cite{Ormerod2013-km_JH}, \cite{Le_Menestrel2004-wy}, and \cite{Mingers2011-id} but also \S\ref{sec:Ethics_and_fairness}.

Analytics maturity models have been developed to describe, explain and evaluate the development of analytics in an organisation. \cite{Krol2020-yp_JEB} reviewed 11 maturity models and assessed them in terms of the number of assessment dimensions, scoring mechanism, number of maturity levels, and the public availability of the methodology. They found that the most common assessment dimensions were technical infrastructure, analytics culture and human resources, including staff’s analytics competencies. \cite{Lismont2017-qz_JEB} undertook a survey of companies, based on the DELTA maturity model \citep{Davenport2010-vu_JEB} of data, enterprise, leadership, targets and analysts. They identified four analytics maturity levels from their survey. The most advanced companies tended to use a wider variety of analytics techniques and applications, to organise analytics more holistically, and to have a more mature data governance policy.

A crucial empirical question is whether Business Analytics adds value to an organisation. An early study on the effect of Business Analytics on supply chain performance was conducted by \cite{Trkman2010-nx_JEB}. They examined over 300 companies, showing a statistically significant relationship between self-assessed analytical capabilities and performance. \cite{Oesterreich2022-ez_JEB} conducted a meta-analysis of 125 firm-level studies, spanning ten years of research in 26 countries.  They found evidence of Business Analytics having a positive impact on operational, financial and market performance. They also found that human resources, management capabilities and organisational culture were major determinants of value creation, whereas technological factors were less important.   

\subsection[Combinatorial optimisation (Silvano~Martello \& Paolo~Toth)]{Combinatorial optimisation\protect\footnote{This subsection was written by Silvano~Martello and Paolo~Toth.}}
\label{sec:Combinatorial_optimisation}

A {\em Combinatorial Optimisation} (CO) problem consists of searching for the optimal element in a finite collection of elements. More formally, given a set of elements and a family of its subsets, each defining a feasible solution and having an associated value, a CO problem is to find  a subset having the minimum (or, alternatively, the maximum) value. The subsets may be proper, like, e.g., in the knapsack problem, or represented by permutations, like, e.g., in the assignment problem (see below). Typically, the  feasible solutions are not explicitly listed, but are described in a concise manner (like a set of equalities and inequalities, or a  graph structure) and their number is huge, so scanning all feasible solutions to select the optimal one is intractable. A CO problem can usually be modelled as an {\em Integer Program} (IP, see also \S\ref{sec:Mixed_integer_programming}) with linear or nonlinear objective function and constraints, in which the variables can take a finite number of integer values.

Consider for example the problem of assigning $n$ tasks to $n$ agents, by knowing the time that each agent needs to complete each task, with the objective of finding  a solution that minimises the overall time needed to complete all tasks ({\em Assignment Problem}, AP). The solution could obviously be found by enumerating all permutations of the integers $1, 2\dots, n$ and selecting the best one. However, this number is so huge that such approach is ruled out even for small-size problem instances: for $n = 30$, we have $n! \cong 2.6\cdot10^{30}$, and the fastest supercomputer on hearth would need millions of years to scan all solutions. The challenge is thus to find more efficient methods. For example, one of the most famous CO algorithms (the {\em Hungarian algorithm}) can solve assignment problem instances with millions of variables in few seconds on a standard PC.

The algorithm mentioned above can be implemented so as to solve any AP instance in a time of order $n^3$, i.e., in a time bounded by a polynomial function of the input size. Unfortunately, only for relatively few CO problems we know algorithms with such property, while for most of them (${\cal NP}$-hard problems) the best known algorithms can take, in the worst case, a time that grows exponentially in the size of the instance. In addition, Complexity theory (see also \S\ref{sec:Computational_complexity}) suggests that the existence of polynomial-time algorithms for such problems is unlikely. On the other hand, CO problems arise in many industrial sectors (manufacturing, crew scheduling, telecommunication, distribution, to mention a few) and hence there is the prominent and practical need to obtain good quality solutions, especially to large-size instances, in reasonable times.

\subsubsection*{Origins}
Many problems arising on graphs and networks (see \S\ref{sec:Graphs_and_networks}) belong to CO (the AP discussed above can be described as that of finding a minimum weight perfect matching in a bipartite graph), and hence the origins of CO date back to the eighteen century. In the following, we narrow our focus to modern CO \citep[see][]{S05_SMPT}. Its roots can be found in the first decades of the past century, when Central European mathematicians developed seminal studies on \textit{matching problems} \citep{K16_SMPT}, \textit{paths} \cite{M27_SMPT}, and {\em Shortest Spanning Trees} (SST) (\citeauthor{J30_SMPT}, \citeyear{J30_SMPT}; \citeauthor{B26_SMPT}, \citeyear{B26_SMPT}, results independently rediscovered by \citeauthor{P57_SMPT}, \citeyear{P57_SMPT} and \citeauthor{K57_SMPT}, \citeyear{K57_SMPT}).
The Fifties produced major results on the AP (\citeauthor{K55_SMPT}, \citeyear{K55_SMPT,K56_SMPT}, on the basis of the results by \citeauthor{K16_SMPT}, \citeyear{K16_SMPT} and \citeauthor{E31_SMPT}, \citeyear{E31_SMPT}, also see \citeauthor{M10_SMPT}, \citeyear{M10_SMPT}), the {\em Travelling Salesman Problem} \citep{DFJ:1954_IL}, and {\em Network Flows} \citep{Ford-Fulkerson:1962_IL},
as well as fundamental studies on basic methodologies: {\em dynamic programming} \citep[DP;][see \S \ref{sec:Dynamic_programming}]{Bellman1957-ue_HL}, {\em cutting planes} \citep[][see \S \ref{sec:Mixed_integer_programming}]{Go58_ALAL}, and {\em branch-and-bound} \citep{LD60_ALAL}.

\subsubsection*{Problems and complexity}
The most important CO problems, for which we know there are polynomial algorithms, are the basic graph-theory problems mentioned in the previous section. Other important problems, which are relevant both from the theoretical point of view and from that of real-world applications, are instead ${\cal NP}$-hard. The main ${\cal NP}$-hard CO problems arise in the following areas.

{\em Scheduling}. Given a set of tasks which must be processed on a set of processors, a scheduling problem asks to find a processing schedule that satisfies prescribed conditions and minimises (or maximises) an objective function, frequently related to the time needed to complete all tasks. This huge area, that includes literally hundreds of problems and variants (mostly ${\cal NP}$-hard), is also discussed in \S \ref{sec:Timetabling}.

{\em Travelling Salesman Problem} (TSP). Given a weighted (directed or undirected) graph, the problem is to find a circuit that visits each vertex exactly once (Hamiltonian tour) and has minimum total weight. This is one of the most intensively studied problems of CO, and is treated in detail in \S \ref{sec:Graphs_and_networks}.

{\em Vehicle Routing Problems} (VRP). A VRP is a generalisation of the TSP which consists of finding a set of routes for a fleet of vehicles, based at one or more depots, to deliver goods to a given set of customers by satisfying a set of conditions and minimising the overall transportation cost.

{\em Facility Location}. These problems require to find the best placement of facilities on the basis  of geographical demands, installation costs, and transportation costs, so as to satisfy a set of conditions and to minimise the total cost (see \S \ref{sec:Location} for a detailed treatment).

{\em Steiner Trees}. Given a weighted graph and a subset $S$ of vertices, it is requested to find an SST connecting all vertices in $S$ (possibly containing additional vertices). These problems, which generalise both the shortest path problem and the SST, are treated in detail in \S \ref{sec:Graphs_and_networks}.

{\em Set Covering}. Given a set of elements and a collection of its subsets, each having a cost, we want to find the least cost sub-collection whose union includes (covers) all the elements.

{\em Maximum Clique} (MC). A clique is a {\em complete} subgraph of a graph (i.e., it is defined by a subset of vertices all adjacent to each other). Given a graph, the MC problem is to find a clique of maximum cardinality (or, if the graph is weighted, a clique of maximum weight). We refer the reader to \S \ref{sec:Graphs_and_networks} for a  detailed analysis.

{\em Cutting and Packing} (C\&P). Given a set of  ``small" items, and a set of ``large" containers, a problem in this area asks for an optimal arrangement of the items into the containers. Items and containers can be in one dimension ({\em Knapsack Problems} (KP), Bin Packing problems) or in more - usually two or three – dimensions (C\&P). See \S \ref{sec:Cutting_and_packing} for more details.

{\em Quadratic Variants of CO problems}. A currently hot research area concerns CO problems whose ``normal" linear objective function is replaced by a quadratic one. This greatly increases difficulty: in most cases problems which, in their linear formulation, can be solved in polynomial time (e.g., the AP) or in pseudo-polynomial time (e.g., the KP) become strongly ${\cal NP}$-hard.

\subsubsection*{Exact methods for ${\cal NP}$-hard problems}
For heuristic and approximation algorithms, we refer the reader to \S \ref{sec:Heuristics} and \S \ref{sec:Computational_complexity}. With the exception of DP methods (\S \ref{sec:Dynamic_programming}), most exact algorithms for $\mathcal{NP}$-hard CO problems, as well as most commonly used ILP solvers, are based on implicit enumeration. In the worst case, they can require the evaluation of all feasible solutions, and hence computing times growing exponentially with the problem size.
The most common methods can be classified as
\begin{itemize}[noitemsep]
\item {\em Branch-and-Bound} (B\&B);
\item {\em Branch-and-Cut} (B\&C);
\item {\em Branch-and-Price} (B\&P).
\end{itemize}
We will describe B\&B, the other methods (and their combinations, as B\&C-and-Price) being its extensions described in \S \ref{sec:Mixed_integer_programming}.

We consider a {\em maximisation} CO problem having an IP model with inequality constraints of `$\leq$' type. For a problem $P$, having feasible solution set $F(P)$, $z(P)$ denotes its optimal solution value, and $ub(P)$ an upper bound on $z(P)$. The main ingredients of B\&B are branching scheme and upper bound computation. \\

\noindent {\em Branching scheme}. The solution is obtained as follows:
\begin{itemize}[noitemsep]
\item[(i)] subdivide $P$ into $m$ subproblems, each having the same objective function as $P$ and a feasible solution set contained in $F(P)$, such that the union of their feasible solution sets is $F(P)$. The optimal solution of $P$ is thus given by the optimal solution of the subproblem having the maximum objective function value;
\item[(ii)] iteratively, if a subproblem cannot be immediately solved, subdivide it into additional subproblems.
\end{itemize}
The resulting method can be represented through a {\em branch-decision tree}, where the root node corresponds to $P$ and each node corresponds to a subproblem.

A node of the tree can be eliminated if the feasible solution set of the corresponding subproblem is empty, or its upper bound is not greater than the value of the best feasible solution to $P$ found so far.\\ 

\noindent {\em Upper bound computation}.
A valid upper bound $ub(P)$ can be computed as the optimal solution value of a {\em Relaxation} of the IP model of $P$, defined by:
\begin{itemize}[noitemsep]
\item[(i)] a feasible solutions set containing $F(P)$;
\item[(ii)] an objective function whose value is not smaller than that of $P$ for any solution in $F(P)$.
\end{itemize}
A relaxation is ``good" if the resulting upper bound $ub(P)$ is ``close" to $z(P)$ (i.e., if the {\em gap} between the two values, $z(P) - ub(P)$, is ``small''), and the relaxed problem is ``easy'' to solve, i.e., its optimal solution can be obtained with a computational effort much smaller than that required to solve $P$.

\subsubsection*{Relaxations}
The most common relaxation methods are:
\begin{itemize}[noitemsep]
\item {\em Constraint elimination}: a subset of constraints is removed from the IP model of $P$, so that the resulting problem is easy to to solve. The most widely used case is the linear relaxation;
\item {\em Linear relaxation}: when the model is an {\em Integer Linear Problem} (ILP), removing the constraints that impose integrality of the variables leads to a {\em Linear Program} (LP), which is polynomially solvable, commonly used in ILP solvers (see \S \ref{sec:Mixed_integer_programming});
\item {\em Surrogate relaxation}: a subset $\Sigma$ of inequality constraints is replaced by a single {\em surrogate} inequality, so that the corresponding relaxed problem is easy to solve. The surrogate inequality is obtained by multiplying both sides of each inequality of $\Sigma$ by a non-negative constant, and summing, respectively, the left-hand and right-hand sides of the resulting inequalities;
\item {\em Lagrangian relaxation}: a subset $\Lambda$ of inequality constraints is removed from the model and ``embedded'', in a Lagrangian fashion, into the objective function. For each inequality of $\Lambda$, the difference between left-hand and right-hand sides ({\em slack}) multiplied by a non-negative constant is added to the objective function.
\end{itemize}

The relaxations can be strengthened by adding one or more valid inequalities ({\em cuts}) to the IP model of $P$, such that they are redundant for the IP model, but can become active when the IP model is relaxed (see \S\ref{sec:Mixed_integer_programming}).

\subsubsection*{Further readings}
We refer the reader to the following selection of references for more details on the topics covered in this section. Well known, pre-1990 books are those by \citeauthor{GN72_SMPT} (\citeyear{GN72_SMPT}, IP), \citeauthor{C75_SMPT} (\citeyear{C75_SMPT}, algorithmic graph theory), \citeauthor{GJ79_SMPT} (\citeyear{GJ79_SMPT}, complexity), 
\citeauthor{BD80_SMPT} (\citeyear{BD80_SMPT}, AP), \citeauthor{LLRS85_SMPT} (\citeyear{LLRS85_SMPT}, TSP), and the CO specific volumes by  \cite{L76_SMPT}, \cite{CMTS79_SMPT}, \cite{PS82_SMPT},
 \cite{MLMR87_SMPT}, and \cite{NW88_SMPT}. We list more recent contributions in the order in which the topics were introduced:
\begin{itemize}[noitemsep]
\item CO: \cite{CCPS98_SMPT}, \cite{Schrijver:2003_IL};
\item AP: \cite{BDM12_SMPT} for linear and quadratic AP, \cite{C13_SMPT} for quadratic AP;
\item Network Flows: \cite{ahuja1993network_IL};
\item Scheduling: \cite{BEPSW01_SMPT,BEPSW07_SMPT}, \cite{P12_SMPT};
\item TSP: \cite{GP06_SMPT}, \cite{ABCC07_SMPT}, \cite{Cook:2011_IL};
\item VRP: \cite{GoldenRW2008_CA_MB}, \cite{TV14_SMPT};
\item Facility Location: \cite{MPL90_SMPT}, \cite{LNS15_SMPT};
\item Steiner trees: \cite{HRW92_SMPT}, \cite{PS12_SMPT}. Also see the recent survey by \cite{Ljubic:2021_IL};
\item Cutting and packing: \cite{MT90_SMPT}, \cite{KPP04_SMPT}. Also see the recent survey by \cite{CILM22A_SMPT,CILM22B_SMPT}.
\end{itemize}

\subsection[Computational complexity (Ulrich~Pferschy \& Clemens~Thielen)]{Computational complexity\protect\footnote{This subsection was written by Ulrich~Pferschy and Clemens~Thielen.}}
\label{sec:Computational_complexity}

Operational Research develops models and solution methods for problems arising from practical decision making scenarios. Often, these solution methods are \emph{algorithms}. The difficulty of a problem can be assessed empirically by evaluating the running times of corresponding algorithms, which requires careful implementations and meaningful test data. Moreover, this can be time-consuming and yields insights that depend on the skills of the programmer and are limited to the available test instances. 
\emph{Computational complexity} represents an alternative approach that allows for a more general assessment of a problem's difficulty that is independent of specific problem instances or solution algorithms. 

\subsubsection*{Problem encoding and running times of algorithms}

In complexity theory, the running time of an algorithm is expressed in terms of the size of the input, i.e., the amount of data necessary to encode an instance of the problem. Since computers store data in the form of binary digits (bits), the standard \emph{binary encoding} represents all data of a problem instance in the form of binary numbers. The number of required bits (the \emph{encoding length}) of an integer is roughly given by the binary logarithm of its absolute value. As an example, consider the binary encoding of instances of the well-known 0-1 knapsack problem (KP). An instance of KP consists of $n$~items -- each with a non-negative, integer weight and profit -- and a positive, integer knapsack capacity~$c$.
We can assume that all $n$~item weights are bounded by the capacity~$c$ and denote the value of the largest item profit by~$p_{\max}$. Then, the encoding length of a KP instance is bounded by $(n+1)\cdot\log_2(c) + n\cdot\log_2(p_{\max})\leq (2n+1)\cdot\log_2(\max\{c,p_{\max}\})$.

Rational numbers can be straightforwardly represented by their (integer) numerator and denominator, but their presence in the input might already influence a problem's computational complexity \citep{Woi18_UPCT}. Irrational numbers cannot be encoded in binary without rounding them appropriately, which means that a different kind of complexity theory is required when general real numbers are part of the input \citep[see][for details]{BCSS98_UPCT}. Hence, the following exposition is restricted to the case of integer inputs, where the encoding length of an instance can be bounded by the number of integers needed to represent it multiplied with the binary logarithm of the largest among their absolute values (see the bound for KP instances provided above as an example).

To allow universal running time analyses of algorithms that are independent of specific computer architectures, asymptotic running time bounds described using the so-called $\mathcal{O}$-notation \citep{Cor09_UPCT} are used. Informally, every polynomial in~$n$ with largest exponent~$k$ is in $\mathcal{O}(n^k)$. All terms with exponents smaller than~$k$ and the constant coefficient of~$n^k$ are ignored. One is then often interested in \emph{po\-ly\-no\-mi\-al-time algorithms} whose running time is in $\mathcal{O}(|I|^k)$ for some constant~$k$, where~$|I|$ denotes the encoding length of instance~$I$. A less preferred outcome would be a \emph{pseudopolynomial-time algorithm}, where the running time is only required to be polynomial in the number of integers in the input and the largest among their absolute values (or, equivalently, in the exponentially larger encoding length of the input when using \emph{unary encoding}, where the encoding length of an integer is roughly its absolute value). 

\subsubsection*{The complexity classes $\mathcal{P}$ and $\mathcal{NP}$}

Most application scenarios encountered in Operational Research finally lead to an optimisation problem (often a combinatorial problem -- see \S\ref{sec:Combinatorial_optimisation}), where a feasible solution is sought that minimises or maximises a given objective function. Every optimisation problem immediately yields an associated decision problem, asking a yes-no question. For example, a minimisation problem consisting of a set~$\mathcal{X}$ of feasible solutions and an objective function~$f$ can be written as $\min \{f(x): x \in \mathcal{X}\}$. For a given target value~$v$, the associated decision problem then asks: Does there exist a feasible solution~$x \in \mathcal{X}$ such that $f(x) \leq v$? Solving an optimisation problem to optimality trivially answers the associated decision problem for any given~$v$. On the other hand, every algorithm for the decision problem can be used to solve the underlying optimisation problem. Given upper and lower bounds, the optimal solution value can be identified in polynomial time by performing binary search between these bounds using the decision problem to answer the query in every iteration of the binary search (assuming that the range of objective function values and the encoding lengths of the bounds are polynomially bounded).

Motivated by the above, the computational complexity of an optimisation problem follows from the complexity of its associated decision problem. Here, the most relevant complexity classes in Operational Research are probably $\mathcal{P}$ and $\mathcal{NP}$, which are often used to draw the line between ``easy'' and ``hard'' problems in this context. Formally, the class $\mathcal{P}$ (``polynomial'') consists of all decision problems for which a polynomial-time solution algorithm exists on a deterministic Turing machine (or, equivalently, in any other ``reasonable'' deterministic model of computation), while the class $\mathcal{NP}$ (``nondeterministic polynomial'') consists of all decision problems for which the same holds on a \emph{nondeterministic} Turing machine. Equivalently, $\mathcal{NP}$ is the class of all decision problems such that, for any yes instance~$I$, there exists a certificate with encoding length polynomial in~$|I|$ and a deterministic algorithm that, given the certificate, can verify in polynomial time that the instance is indeed a yes instance. Since the most natural certificate is often a (sufficiently good) solution of the problem, $\mathcal{NP}$ can informally be defined as the class of decision problems for which solutions can be verified in polynomial time. For example, when considering the travelling salesman problem (TSP) on a given edge-weighted graph, the associated decision problem asks whether or not there exists a tour (Hamiltonian cycle) of at most a given length~$v$. While no polynomial-time algorithm for this decision problem is known to date, the problem can easily be seen to be in $\mathcal{NP}$ since the natural certificate to provide for a yes instance is simply a tour with length at most~$v$, whose feasibility and length can be easily verified in polynomial time.

Observe that these definitions directly imply that $\mathcal{P}\subseteq\mathcal{NP}$. Most researchers believe that $\mathcal{P}\subsetneq\mathcal{NP}$ or, equivalently, that there are problems in $\mathcal{NP}$ that do not admit polynomial-time solution algorithms. However, formally proving that $\mathcal{P}\neq\mathcal{NP}$ (or that $\mathcal{P}=\mathcal{NP}$) is still one of the most famous open problems in theoretical computer science to date.
            
This so-called \emph{$\mathcal{P}$ versus $\mathcal{NP}$ problem} can be equivalently expressed using the well-known notion of $\mathcal{NP}$-completeness \citep[see, e.g.,][]{GJ79_SMPT}. Intuitively, $\mathcal{NP}$-complete problems are the hardest problems in $\mathcal{NP}$ in the sense that, if one of these problems admits a polynomial-time solution algorithm, then so does every problem in $\mathcal{NP}$ (and, thus, we would obtain $\mathcal{P}=\mathcal{NP}$). A decision problem (not necessarily in $\mathcal{NP}$) with this property is also called $\mathcal{NP}$-hard. This means that a problem is $\mathcal{NP}$-complete if and only if it is both $\mathcal{NP}$-hard and contained in $\mathcal{NP}$. 
The first problem shown to be $\mathcal{NP}$-complete in Cook's famous theorem \citep{Cook:STOC71_UPCT} is the (Boolean) satisfiability problem (SAT). Shortly after, \cite{Karp:complexity_UPCT} gave a list of 21~fundamental problems that are $\mathcal{NP}$-complete. While Cook's proof that SAT is $\mathcal{NP}$-complete required considerable effort, proving that further problems are $\mathcal{NP}$-complete became significantly easier with this knowledge and hundreds -- if not thousands -- of problems were shown to be $\mathcal{NP}$-complete.

A decision problem is $\mathcal{NP}$-complete if and only if (1) it is contained in $\mathcal{NP}$ and (2) some $\mathcal{NP}$-complete problem (and, therefore, all problems in $\mathcal{NP}$) can be reduced to it via a \emph{polynomial-time reduction}. Such a polynomial-time reduction works as follows: For any instance of the known $\mathcal{NP}$-complete problem (e.g., SAT or TSP), one has to construct an instance of the investigated problem in polynomial time such that the two instances are equivalent, i.e., the constructed instance is a yes instance if and only if the given instance is a yes instance. Note that the requirement that the instance must be constructed in polynomial time (and, therefore, have encoding length polynomial in the encoding length of the original instance) is crucial. A common error in reductions is that the encoding length of the constructed instance depends polynomially on the size of numerical values in the given instance (instead of their encoding length).

The importance of the encoding can be illustrated by the 0-1 knapsack problem (KP), which is $\mathcal{NP}$-hard if binary encoding is used, but can be solved in polynomial time (via dynamic programming) if unary encoding is used (so $\mathcal{NP}$-hardness of the unary-encoded version would imply that $\mathcal{P}=\mathcal{NP}$). Problems like this, i.e., problems whose binary-encoded version is $\mathcal{NP}$-hard, but whose unary-encoded version can be solved in polynomial time, are called \emph{weakly $\mathcal{NP}$-hard}, while problems (such as SAT) that remain $\mathcal{NP}$-hard also under unary encoding are called \emph{strongly $\mathcal{NP}$-hard}. The existence of a pseudopolynomial-time algorithm is possible for weakly $\mathcal{NP}$-hard problems, but not for strongly $\mathcal{NP}$-hard problems (unless $\mathcal{P}=\mathcal{NP}$).
            
\subsubsection*{Approximation algorithms}

While some realistic-size instances of $\mathcal{NP}$-hard problems might still be solvable in reasonable time, this is not the case for all instances. In general, one can deal with 
$\mathcal{NP}$-hardness by relaxing the requirement of finding an optimal solution and instead settling for a ``good-enough'' solution. This leads to heuristics, whose aim is producing good-enough solutions in reasonable time (see \S\ref{sec:Heuristics} for details) and approximation algorithms \citep{Vazi10_UPCT, SW11_UPCT, Ausi13_UPCT}. Given~$\alpha\geq 1$, an \emph{$\alpha$-approximation algorithm} for an optimisation problem is a polynomial-time algorithm that, for each instance of the problem, produces a solution whose objective value is at most a factor~$\alpha$ worse than the optimal objective value. The factor~$\alpha$, which can be a constant or a function of the instance size, is then called the \emph{approximation ratio} or \emph{performance guarantee} of the approximation algorithm. While it is standard to use $\alpha\geq1$ for minimisation problems, there is no clear consensus in the literature as to whether $\alpha\geq1$ or $\alpha\leq1$ should be used for maximisation problems. For example, the simple extended greedy algorithm for the knapsack problem produces a solution with at least half of the optimal objective value on each instance, i.e., it is a $1/2$- or a $2$-approximation algorithm.

While \emph{inapproximability results} can be shown for some $\mathcal{NP}$-hard problems \citep[see][ch.~10]{Hoch96_UPCT}, others allow for approximation algorithms with approximation ratios arbitrarily close to one, i.e., they admit a \emph{polynomial-time approximation scheme} (PTAS). A PTAS is a family of algorithms that contains a $(1+\varepsilon)$-approximation algorithm for every $\varepsilon>0$. If the running time is additionally polynomial in $1/\varepsilon$, the PTAS is called a \emph{fully polynomial-time approximation scheme} (FPTAS). If all objective function values are integers, every FPTAS can be turned into a pseudopolynomial-time exact algorithm, so strongly $\mathcal{NP}$-hard problems do not admit an FPTAS (unless $\mathcal{P}=\mathcal{NP}$). Conversely, pseudopolynomial-time algorithms, in particular dynamic programming algorithms, often serve as a starting point for designing an FPTAS \citep{Woe00_UPCT,PW07_UPCT}.

\subsubsection*{Further complexity classes}
Theoretical computer science developed a wide range of complexity classes far beyond the $\mathcal{P}$ vs.\ $\mathcal{NP}$ dichotomy. Considering algorithms requiring polynomial space, i.e., for which the encoding length of the data stored at any time during the algorithm's execution is polynomial in the encoding length of the input (but no bound on the running time is required), gives rise to the class $\textsf{PSPACE}$. It is widely believed that $\mathcal{NP}\subsetneq\textsf{PSPACE}$, but even whether $\mathcal{P} \neq \textsf{PSPACE}$ holds is not known.

In the theoretical analysis of bilevel optimisation problems \citep[see, e.g.,][]{Lab16_UPCT} the complexity class~$\Sigma_2^p$ plays an important role \citep[see][]{Woe21_UPCT}. 
Here, a yes instance~$I$ is characterised by the existence of a certificate of encoding length polynomial in~$|I|$ such that a certain polynomial-time-verifiable property holds true \emph{for all elements of a given set~$\mathcal{Y}$}. As an example, consider the $2$-quantified (Boolean) satisfiability problem. Here, an instance consists of two sets~$X$ and~$Y$ of Boolean variables and a Boolean formula over $X\cup Y$. The question then is whether there exists a truth assignment of the variables in~$X$ such that the formula evaluates to true \emph{for all possible truth assignments of the variables in~$Y$}. This definition immediately sets the stage for a bilevel problem, where the decision~$x$ of the upper level (the leader) should guarantee a certain outcome \emph{for every possible decision~$y$ at the lower level} (the follower). It is widely believed that $\mathcal{NP}\subsetneq \Sigma_2^p$, although $\Sigma_2^p$-hardness does not rule out the existence of a PTAS \citep{Cap14_UPCT}. Under this assumption, $\Sigma_2^p$-hardness does, however, rule out the existence of a compact ILP-formulation, which can be a valuable finding for bilevel optimisation problems.

For some $\mathcal{NP}$-hard problems, one can construct algorithms with running time $\mathcal{O}(f(k)\cdot \mbox{poly}(|I|))$ for an arbitrary computable function~$f$, where the parameter~$k$ describes a property of the instance~$I$. Such problems are called \emph{fixed-parameter tractable}.
For example, the satisfiability problem SAT is fixed-parameter tractable with respect to the parameter~$k$ that represents the tree-width of the primal graph of the SAT instance.
In this graph, the vertices are the variables and two vertices are joined by an edge if the associated variables occur together in at least one clause, see \cite{Szei03}.
This parametric point of view is captured in the $W$-hierarchy of complexity classes -- see \cite{Nied06_UPCT} and the seminal book by \cite{DF99_UPCT}.

\subsection[Control theory (Xun~Wang)]{Control theory\protect\footnote{This subsection was written by Xun~Wang.}}
\label{sec:Control_theory}

Control theory deals with designing a control signal so that the state or output variables of the system meet certain criteria. It is a broad umbrella term that covers a variety of theories and techniques. Control theory has been widely applied in the studies of economics \citep{tustin1953_XW, grubbstrom1967_XW}, operations management \citep[][also see \citealt{Sarimveis2008_XW} for a recent review]{Simon1952_XW, Vassian1955_XW}, and finance \citep{Sethi1970_XW}. Here, we do not intend to provide an exhaustive or comprehensive review. Instead, we try to structurally organise the concepts and techniques commonly applied in operations research, which means that technical details will be omitted. We direct interested readers to a number of textbooks in the reference list, and an excellent review by \cite{aastrom2014_XW} for those interested in the development of control theory.

The major distinction between control theory and other optimisation theories is that the control variable to be designed is normally a time-varying, dynamic function. The control signal can either be dependent on the state variables (which is referred to as feedback control or closed-loop control) or independent (feedforward control or open-loop control). The design of control signals and control policies (defined as the function between the state of the system and the control, also known as ``control laws'' or ``decision rules'') is based on the structure of the system to be controlled (sometimes called the ``plant'' in the control engineering literature). Thus, the type of the dynamical system often define the type of control problem. In continuous systems, the time variable is defined on the real axis, suitable to describe continuous processes such as fluid processing and finance. In discrete systems, time is defined as integers, suitable in cases such as production and inventory control, where the production quantity is released every day. Linear systems are comprised of linear (or affine) state equations, while nonlinear systems contain nonlinear elements. Nonlinear systems are more difficult to analyse and control, and may lead to complex system behaviours such as bifurcation, chaos and fractals \citep{strogatz2018_XW}. But there are linearisation strategies which approximate the nonlinear system locally as linear systems \citep{Slotine1991_XW}. Based on whether random input is present, the dynamical system can be categorised into deterministic and stochastic.

There are two fundamental methods in the analysis of the system and control. The first relies on time-frequency transformations (Laplace transform for continuous systems and $z$-transform for discrete systems). A transfer function in the frequency domain can be used to represent and analyse the system \citep{towill1970_XW}. This method saves computational effort; however, it can only deal with linear system models and each transfer function only describes the relation between a single input and a single output (SISO). The second method directly tackles the state equations in the time domain and describes the movement of system state in the state space. It is suitable for nonlinear systems and multi-input-multi-output (MIMO) systems. With the advancement of computing technology, the computational burden faced by the time-based method becomes less significant. The literature refers to the frequency-based method as \emph{classic} control theory \citep{Ogata2010_XW} and the time-based method as \emph{modern} control theory.

The system under the effect of the control policy must be examined with respect to its properties and dynamic performance. \emph{Stability} is a property of the dynamical system, that the system can return to its steady state after receiving a finite external disturbance. Stability is a fundamental precondition that almost all control designs must meet, with few exceptions such as clocks and metronomes, where a periodic or cyclic response is desired. The stability criterion is straightforward to derive for linear systems, where both frequency-based \citep[e.g., Routh-Hurwitz stabibility criteria and Jury's inners approach,][]{jury1975_XW} and time-based (e.g., the eigenvalue approach) methods exist. However, stability analysis for nonlinear systems is more challenging \citep{Bacciotti2005_XW}. Other important properties of the control system include \emph{controllability}, defined as the ability to move the system to preferred state using only the control signal; and \emph{observability}, defined as the ability to infer the system state using the observable output signals \citep{Gilbert1963_XW}.

In addition to these intrinsic properties, the system can also be evaluated by the system's response to some characteristic input functions. The \emph{step function} (sometimes referred to as the Heaviside function) takes the value of zero before the reference time point, and one thereafter. The \emph{impulse function} (the Dirac $\delta$ function) takes the value of infinity at the reference time point and zero otherwise. These two input functions usually represent an abrupt change in the external environment. The \emph{sinusoidal function} can be used to describe the periodic and seasonal externalities. The Bode plot describes the amplitude and phase shift between the sinusoidal input and output. For stochastic environments, the \emph{white noise} is used to mimic random disturbances. It is a random signal that follows an identical and independent Gaussian distribution and has a constant power spectrum. The noise bandwidth of the system determines the ratio between output and input variances when the input is iid. The value of the noice bandwidth can be derived from either the transfer function or the state space representation. This concept is used in analysing the amplification phenomena in supply chains (see \S\ref{sec:Supply_chain_management}).

In practise, the system state and even the system structure may be unknown. Therefore, statistical techniques, known as \emph{state estimation} and \emph{system identification}, have been developed. State estimation uses observable output data to estimate the unobservable system states. A popular technique for this purpose is the Kalman filter \citep{kalman1960_XW}, essentially an adaptive estimator that can be applied not only in linear, time-invariant cases (LTI, where the system is linear and does not change over time), but also non-linear and time-variant cases. For example, it has been applied to estimate the demand process from the observed sales data \citep{Aviv2003_XW}. System identification attempts to ``guess'' the structure of the system from the input and output. 

Along with the development of control theory, various control strategies have been proposed. They are designed to fit the structure of the system, the objective of the control, and most importantly, to offer a paradigm to design the control policy. In what follows, we provide a brief summary of control strategies. \emph{Linear} control strategies can be represented linearly (in the form of transfer function). They offer great analytical tractability and satisfactory performance, especially when the open-loop system is also linear. Two widely adopted policies in this family are proportional-integral-derivative (PID) control and full-state feedback (FSF) control. In PID control, an error signal between the output and the reference input (e.g., a Heaviside function) is computed. The control signal is a linear combination of the error, the integral, and the derivative of the error. These three components can appear separately. The proportional control has been applied in mechanical and managerial mechanisms such as the centrifugal governor and production planning \citep{Chen2007_XW}. The full-state feedback control defines the control signal as a linear combination of the full system state vector, where the coefficient vector (the ``gain'') shares the same dimensionality as the state vector. By tuning the gain, the poles of the closed-loop system (the eigenvalue of the transition matrix or the roots of the characteristic equation) will change their position in the complex plane, adjusting the system performance. The full-state feedback policy can also be applied in production and inventory control \citep{Gaalman2006_XW}.

In contrast to the linear strategies, the \emph{nonlinear} control strategies are defined as policies where the control signal cannot be represented by a linear function of the system state \citep{Slotine1991_XW}. These policies are primarily used when the open loop system is also nonlinear. One such policy is sliding mode control, where the control signal is a switching function of the state, dependent on some switching rules. The system is then maintained near a hypersurface of the state space (sliding), where the dynamic behaviour of the system is desired. It should be ensured that the hypersurface is reachable from any initial state and that the system state can be maintained on the hypersurface by the policy. In practise, bang-bang control is adopted frequently as a special case of sliding mode control, where the control signal can take only two possible values. The rocket engines and domestic thermostats are examples of such (with on and off states). 

\emph{Optimal control} aims at finding the control signal or control policy that allows an objective function to reach its extreme point \citep{Sethi2009_XW, bertsekas2012dynamic_DLLD}. The objective function could be dependent on the state, output and/or control. Many control policies mentioned above, e.g. full-state feedback control and sliding mode control, have been proved to be the optimal control of some control problems. Optimal control in the special sense is based on Pontryagin's Maximum (or equivalently Minimum) Principle and mainly deals with the design of the open-loop control signal. When equipped with the Hamilton-Jacobi-Bellman (HJB) theory, it can be used to design optimal feedback control policies. Optimal control is closely connected with dynamic programming, which will be reviewed in \S [dynamic programming]. The optimal control technique has been widely applied in operations management \citep[e.g.,][]{Kumar2003_XW}.

When random external disturbances are present, the \emph{stochastic control} techniques are necessary \citep{aastrom2012_XW}. In these situations, objective functions are usually statistical functions of the state or the output, such as the absolute mean or variance. The most well-studied stochastic control problem is the Linear Quadratic Gaussian (LQG) problem, where the system is linear, the objective function is of quadratic form, and the noise signal follows a Gaussian distribution. The optimal control policy in this case is a linear one. Many supply chain management problems can be modelled in LQG form \citep{Lee1997_XW}. For more complex problems involving nonlinearity or an unspecified system structure, the \emph{model predictive control} (MPC) approach can be used \citep{Camacho2013_XW}. This approach transforms the infinite-horizon problem into a finite-horizon problem by focusing only on $T$ periods in the future, deriving the control signal for these $T$ periods, and adopting the most recent control. In the next period, the prediction is updated, and this process is repeated. MPC is not an optimal control method due to the finite-horizon approximation, yet it works very well in practise \citep{Doganis2008_XW}. To deal with parametrical uncertainties in the disturbance, \emph{robust control} provides guaranteed performance \citep{Zhou1998_XW}. The well-known $H_\infty$ control (H infinity) is one of such examples. It minimises the largest singular value of the transition matrix function, which in SISO systems equates to the peak value of the frequency response curve. This minimax strategy ensures that any frequency component in the input will not be amplified too much. Finally, if the system parameters vary over time, \emph{adaptive control} allows the control policy to update according to the estimated parameters \citep{aastrom2013_XW}. The difference between adaptive and robust control is that the policy is dynamic in the former and static in the latter. 

Recent development of control theory can be seen in the controlling of complex, large scale and network system; the use of artificial intelligence in control engineering; and the application of control theory in areas of physics, biology and economics.

\subsection[Data envelopment analysis (Sebastián~Lozano)]{Data envelopment analysis\protect\footnote{This subsection was written by Sebastián~Lozano.}}
\label{sec:Data_envelopment_analysis}

Data Envelopment Analysis (DEA) is a non-parametric frontier analysis methodology mainly used to assess the relative efficiency of a set of homogeneous operating units (termed Decision Making Units, DMUs). DMUs are assumed to consume inputs (i.e., resources) to produce outputs (e.g., goods and services). The production function that indicates the amount of outputs that can be produced from a given input vector is unknown. DEA does not make any assumption about the functional form of that dependency. Instead, DEA uses the observed data to infer the Production Possibility Set (PPS), also called the DEA technology, which contains all the operating points that are deemed feasible. This is achieved on the basis of a few assumptions (like envelopment of the observations, free disposability of inputs and outputs, convexity and returns to scale) and invoking the Minimum Extrapolation Principle. The resulting PPS contains all linear combinations of the observations along with all the operating points that they dominate. This leads to Linear Programming models whose main decision variables are the intensity variables used to compute the target operating point (projection). This target operating point must dominate the DMU being projected and represents maximal improvements (i.e., input reduction and output increase) with respect to the latter. Hence, the computed target belongs to the efficient frontier (which is the non-dominated subset of the PPS) and the efficiency score is a decreasing function of the distance from the DMU to the computed efficient target. There are different ways of measuring this distance, which, ultimately, depends on the potential input and output improvements (i.e., slacks) computed by DEA. Before diving into the DEA methodology note that, as \cite{Cook2014-xo_SL} point out, although DEA has a strong link with production theory in economics, it is often used to benchmark the performance of manufacturing and service operations. In such benchmarking exercises, the efficient DMUs, as defined by DEA, may not necessarily form a ``production frontier'', but rather a ``best-practice frontier''. Thus, the purpose of the performance measurement exercise affects the classification of the different variables considered into inputs or outputs. 

\subsubsection*{Efficiency assessment and target setting DEA models}
The seminal DEA models by \cite{Charnes1978-nm_SL} and \cite{Banker1984-no_SL}were oriented (i.e., gave priority to reducing the inputs or to increasing the outputs) and looked for a uniform (i.e., radial) improvement in all the input or output dimensions. The projection can also be estimated using a given direction, giving rise to Directional Distance Function (DDF) DEA models \citep{Wang2019-sv_SL}. However, most DEA approaches are non-radial and non-oriented \citep[e.g.,][]{Fukuyama2009-mo_SL}. Actually, because DEA aims at simultaneously improving inputs and outputs, it is inherently a multiobjective optimisation approach. Hence, taking into account the preferences of a decision maker, any Pareto optimal point can be selected as efficient target \citep{Soltani2020-xv_SL}.

Most DEA models compute targets that can be sometimes far away from the observed DMU. This increases the difficult and effort required to achieve the target. Hence, DEA models that compute closest efficient targets have been developed \citep{Aparicio2007-fb_SL}. An alternative is to compute stepwise efficiency improvement approaches that may eventually achieve ambitious efficient targets but after several gradual improvement steps \citep{Lozano2005-hz_SL}.

DEA models for handling non-discretionary variables \citep{Banker1986-ai_SL}, undesirable outputs \citep{Kuosmanen2005-og_SL}, integer variables \citep{Kazemi_Matin2009-ha_SL}, ratio variables \citep{Olesen2022-is_SL}, negative data \citep{Sharp2007-xj_SL}, and fuzzy data \citep{Arana-Jimenez2022-nv_SL} have also been proposed. Each of the above ``complications'' requires specific adaptations of the methodology and being capable of taking them into account is a proof of the power and flexibility of DEA.

The DEA models based on the PPS concept are labelled as envelopment formulations. There are also dual multiplier formulations in which the decision variables are not the intensity variables used to compute the target inputs and outputs but the corresponding input and output shadow prices. Multiplier formulations let each DMU choose these input and output weights so that its efficiency is maximised. This freedom often leads to DMUs choosing idiosyncratic or unreasonable weights. Imposing Assurance Regions (AR) and other types of weight restrictions has been proposed \citep{Allen1997-fa_SL} as well as measuring the efficiency of the DMUs as the average of the cross-efficiency scores computed with the input and output weights chosen by the different DMUs \citep{Doyle1994-of_DEA,Chen2022-ms_SL}. Another alternative that has been proposed is using a Common Set of Weights (CSW) instead of letting each DMU choose its own \citep{Salahi2021-eh_SL}.

In addition to computing efficiency scores, DEA can be used to rank the DMU. The problem here is that in conventional DEA all the DMUs labelled as efficient are tied and cannot be ranked. In addition to the CSW or cross-efficiency approaches mentioned above, there are other DEA-based full ranking methods, like the super-efficiency approach \citep{Tone2002-ha_SL}.  Alternatively, instead of fully ranking the DMUs, ranking intervals and dominance relations can be established \citep{Salo2011-io_DEA}.

\subsubsection*{Dynamic and Network DEA models}
DEA views DMUs as input-output black boxes. However, it is often the case that DMUs have an internal structure with different stages or processes (sometimes labelled subDMUs). Many different Network DEA (NDEA) models have been developed to address these scenarios \citep{Tone2009-ph_SL}. The key features of NDEA models are that each process has its own technology and that, except in the case of parallel processes, there exist intermediate product flows between the processes. Some NDEA models can compute an efficiency score for each process and relate the overall efficiency score to the scores of the individual processes \citep{Kao2016-sl_SL}. It must be noted that the NDEA configuration most studied and most commonly used in practice involves two stages in series \citep[see, e.g.,][]{Cook2010-bk_SL,Halkos2014-de_SL}.

Multi-period and dynamic scenarios can be modelled in a manner similar to NDEA simply by considering each time period as a subDMU. The difference between multi-period approaches \citep{Kao2014-dj_SL} and Dynamic DEA \citep{Tone2010-gh_SL} is that in the latter there are flows between consecutive periods (i.e., carryovers). Dynamic NDEA (DNDEA) models, in which there are carryovers between periods as well as intermediate product flows between the processes, have also been developed \citep{Tone2014-hp_SL}.

\subsubsection*{Centralised DEA models}
DEA generally projects each DMU separately onto the efficient frontier. There are situations in which the DMUs belong to the same organisation and there is a Central Decision Maker (CDM) that is interested in the overall system performance and therefore in projecting all the DMUs simultaneously. This type of Centralised DEA (CDEA) models are commonly used for resource allocation \citep{Lozano2004-mz_SL} and for centralised production planning \citep{Lozano2014-bb_SL}. Also, an approach to measure the centralised efficiency of the individual DMUs in CDEA scenarios has been proposed \citep{Davtalab-Olyaie2022-yb_SL}.

DEA models for allocating a fixed input or common revenue \citep{Li2021-ia_SL} or for fixed-sum-outputs \citep[FSO;][]{Zhu2017-ol_SL} also share with CDEA the need to project all the DMUs simultaneously to take into account their interrelationships. These models, same as CDEA, can use an envelopment or a multiplier formulation. While the key feature of the former is that all DMUs are projected simultaneously, that of the latter is that, same as in CSW, a single set of input and output weights is considered. 

\subsubsection*{DEA and Total Factor Productivity (TFP) growth}
DEA can be used to compute the Malmquist Productivity Index (MPI) by projecting the DMU in two consecutive periods onto the efficient frontier of each period and computing the geometric mean of the change in the corresponding radial efficiency scores between the two periods \citep{Fare1992-gf_SL}. The Malmquist-Luenberger Productivity Indicator (MLPI) is analogous but it employs the arithmetic average and an additive decomposition of DDF efficiency scores \citep{Chambers1996-zd_SL}. In both cases, the TFP growth of each DMU can be decomposed into an efficiency change and a technological change component. Other alternative decompositions of the MPI and MLPI have been developed \citep{Epure2011-yw_SL}.

Other approaches compute a global MPI \citep{Pastor2005-rg_SL,Kao2014-dj_SL}. These have the circularity property, missing in the adjacent-periods MPI. Changes in prices can be also incorporated to compute and decompose a global cost MPI \citep{Tohidi2012-nr_SL}. MPI variants that take into account the projections of all the observations or of different groups of observations as well as approaches to compute and decompose the aggregate productivity growth index of a whole industry and input-specific productivity growth indexes have also been proposed \citep{Aparicio2017-dh_SL,Kapelko2015-wf_SL}.

\subsubsection*{Metafrontier analysis}
In scenarios where the DMUs are heterogeneous and belong to different groups, not necessarily disjoint, the DMUs can be projected onto its group frontier as well as onto the metafrontier that results from enveloping all the group frontiers. Measuring the difference between the corresponding efficiency scores can be used to estimate the distance between both frontiers and hence the corresponding technology gap of each group. Although the group technologies are generally convex, the metatechnology is generally non-convex \citep{Afsharian2018-ah_SL}.

The metafrontier approach can be used in DNDEA \citep{See2021-yy_SL} and CDEA \citep{Gan2021-fz_SL} contexts. Also, using metafrontier concepts with each group of observations corresponding to a different time period, meta-MPI and meta-MLPI can be computed and appropriately decomposed \citep{Portela2010-zl_SL}.

\subsubsection*{Other DEA approaches}
There are other interesting DEA approaches that have not been covered above, like congestion \citep{Ren2021-lh_SL}, window analysis \citep{Peykani2021-qb_SL}, etc. Moreover, the field, although mature, is still expanding, with promising new developments, like Efficiency Analysis Trees (EAT) \citep{Esteve2020-ud_SL}, Support Vector Frontiers (SVF) \citep{Valero-Carreras2022-gd_SL}, or big data DEA \citep{Dellnitz2022-ei_SL}. This is not to mention the large and increasing number of DEA applications (see \S\ref{sec:Education}, \S\ref{sec:Environment}, and \S\ref{sec:Power_markets_and_systems}). For further learning on DEA the interested reader is referred to existing textbooks \citep{Cooper2007-pc_SL}, handbooks \citep{Cooper2011-sx_SL,Cook2014-fa_SL,Zhu2015-co_SL} and review papers \citep{Kao2014-jt_SL,Contreras2020-rn_SL,Peykani2020-zz_SL}.

\subsection[Decision analysis (Matthias~Ehrgott \& Salvatore~Greco)]{Decision analysis\protect\footnote{This subsection was written by Matthias~Ehrgott and Salvatore~Greco.}}
\label{sec:Decision_analysis}

The term decision analysis was introduced by \cite{howarddecision_MESG} as “a logical procedure for the balancing of the factors that influence a decision”, pointing out that “the procedure incorporates uncertainties, values and preferences in a basic structure that models the decision”.  According to \cite{keeney1982decision_MESG} decision analysis is a “formalisation of common sense for decision problems which are too complex for informal use of common sense” and, in more technical form “a philosophy, articulated by a set of logical axioms, and a methodology and collection of procedures, based upon those axioms, for responsibly analysing the complexities inherent in decision problems”. In a slighty different  perspective, \cite{roy1993decision_MESG} proposed the concept of decision aiding as “the activity of one who, in ways we call scientific, helps to obtain elements of answers to questions asked by actors involved in a decision-making process, elements helping to clarify this decision in order to provide actors with the most favourable conditions possible for that type of behaviour which will increase coherence between the evolution of the process, on the one hand, and the goals and/or systems of values within which these actors operate on the other”.

For \cite{howarddecision_MESG} “the essence of the procedure is the construction of a structural model of the decision in a form suitable for computation and manipulation”. For \cite{keeney1982decision_MESG}  “the foundations of decision analysis are provided by a set of axioms ... which provide principles for analysing decision problems”. Moreover, “the philosophical implications of the axioms are that all decisions require subjective judgements and that the likelihoods of various consequences and their desirability should be separately estimated using probabilities and utilities, respectively”. In this perspective, the key components of a decision problem are the set of alternatives to be taken into consideration; the set of consequences describing outcomes of alternatives, possibly in terms of a plurality of attributes or criteria; if the consequences are uncertain, the beliefs about their possible realisations expressed in terms of a probability distribution; the preferences of the decision maker.  The objective of the decision analysis is to construct a value function representing the preferences of the decision maker by assigning each alternative an evaluation of its desirability. In case of uncertainty of the consequences, the value function is expressed in terms of expected value with respect to the probability of the consequences. The basic methodology to induce the value function is based on the pioneering work of \cite{Vomo1944_GZ} that showed that a small set of axioms imply that the “utility” of an outcome $x$ is defined as the probability of getting the most-preferred outcome and otherwise the least-preferred outcome that would be indifferent to receiving outcome $x$ with certainty. For \cite{roy1993decision_MESG}, the decision aiding procedure should be developed in a constructive approach in which “concepts, models, procedures and results are here seen as suitable tools for developing convictions and allowing them to evolve, as well as for communicating with reference to the bases of these convictions”. In this perspective the “object is not to know or to approximate the best possible decision but to develop a corpus of conditions and means on which we can base our decisions in favour of what we believe to be most suitable”.

Decision Analysis is mainly based on concepts and tools related to the subjective probability of \cite{RePEc:hay:hetcha:ramsey1926_MESG} and \cite{de1937prevision_MESG}, the theory of expected utility of  \cite{Vomo1944_GZ} and subjective expected utility of \cite{savage1954foundations_MESG}, the Multiple Attribute Utility Theory (MAUT) of \cite{KeeneyRaiffa1976_MESG} and the psychology of judgement and decision-making of  \cite{tversky1974judgment_MESG}. The general idea is to try to evaluate each alternative by assigning a value based on the utilities of the outcomes obtained in each state of the nature multiplied by their probabilities. Delayed consequences may be discounted according to the time at which they are obtained. Each outcome may be evaluated by considering value trade-offs among multiple attributes. Decision analysis techniques include Utility Function Elicitation techniques, Probability Elicitation protocols, Net Present Value, Decision Trees, Influence Diagrams, and Monte Carlo simulation-based decision analysis \citep{clemen1996making_MESG}; Value-Focused Thinking \citep{keeney1996value_MESG};  Portfolio Decision Analysis \citep{salo2011portfolio_MESG}, Bayesian Networks \citep{pearl1988probabilistic_MESG}, and multi-stage decision optimization techniques such as dynamic rogramming and reinforcement learning.

Considering the distinction between normative, descriptive and prescriptive approaches \citep{bell1988descriptive_MESG}, the general perspective of the decision analysis is prescriptive rather than normative or descriptive \citep{edwards2007advances_MESG}. Descriptive analysis concerns the representation and prediction of observed decisions and normative analysis concerns the decisions that ideally coherent and rational individuals should take. Instead, prescriptive analysis tries to propose methods and techniques that will help real people make better decisions with lower regret and greater coherence of values and behaviors.  In this context, decision analysis takes a prescriptive approach that, focusing on the few basic axioms underlying subjective expected utility, adopts ``pragmatically''  the aspiration to the rationality of the normative approach, trying to correct all the heuristics and biases discovered and investigated by descriptive analysis \citep{tversky1974judgment_MESG}. The decision aiding approach \citep{roy1993decision_MESG} takes a different perspective that, criticising the idea that there is an objectively optimal decision to be discovered or at least approximated, aims to provide a recommendation consisting in a set of convictions constructed in the course of a decision process based on multiple interactions between the analyst and the decision maker. The decision aiding approach leads directly to a multi-criteria perspective \citep{belton2002multiple_MESG,greco2016multiple_MESG} taking explicitly into consideration the multiple attributes or criteria (e.g., related to finance, resources, time, and environmental impacts) to be considered in the decision problem at hand. This avoids the risk of a fictitious, not reasoned and arbitrary conversion of evaluations on different criteria  to a common unit, facilitating the discussion on the respective role of each criterion  \citep{roy2005paradigms_MESG,roy1996multicriteria_MESG}. To compare alternatives in a multicriteria decision procedure four main approaches can be adopted:
\begin{itemize}[noitemsep]
	\item aggregating criteria assigning a single value to each alternative: this is the case of above mentioned MAUT, as well as of some of the most well known multicriteria methods such as SMART \citep{edwards1994smarts_MESG}, and UTA \citep{jacquet1982assessing_MESG}; a specific mention deserves in this context the AHP approach \citep{saaty1977scaling_MESG}, that is probably the most adopted \citep[although controversial; see, e.g.,][]{dyer1990remarks} multicriteria method. It is based on the comparison of "importance" of criteria and of evaluation of alternatives with respect to considered criteria by means of a nine point qualitative scale; another specific class in this family are the distance-based methods which, following the main principle of TOPSIS \citep{Hwang1981multiple_MESG}, the first and most famous of these methods, evaluate each alternative on the basis of their distance from the positive ideal solution and the negative ideal solution (the fictitious alternatives that have the best and the worst evaluation on each criterion, respectively); two other well-known methods in this class are VIKOR \citep{opricovic2004compromise} and TODIM \citep{gomes1991todim}.
	\item aggregating criteria by means of one or more synthesising preference relations: the most well known methods based on this approach are the ELECTRE methods \citep{figueira2013overview_MESG,figueira2016electre_MESG}, that build a crisp or valued preference relation called outranking for which an alternative $a$ is at least as good as another alternative $b$ if a $a$ is not worse than $b$ for a majority of important criteria (concordance) and there is no criterion for which the advantage of $b$ over $a$ is so large that it prevents the possibility to declare $a$ at least as good as $b$ (non-discordance); 
	\item aggregating criteria through ``if $\ldots$, then $\ldots$'' decision rules \citep{greco2001rough_MESG}: the alternatives obtain an overall evaluation by matching decision rules with a syntax ``if the alternative is at least at level $l_{j_1}$ on criterion $g_{j_1}$ and $\ldots$ at least at level $l_{j_r}$ on criterion $g_{j_r}$, then the alternative is globally at least at level $l_{tot}$'', such as ``if the student has an evaluation at least good on mathematics and at least medium on literature, then the student is globally at least medium''; these rules are induced from a set of examples of decisions supplied by the decision maker. The advantage of this approach is its explainability due to the fact that the decision rules are expressed in natural language;  
	\item aggregating criteria through an interactive multiobjective optimisation \citep{branke2008multiobjective_MESG}: with this approach one can handle decision problems in which a set of objectives have to be optimised under given constraints \citep[see][]{sawnata85_MESG,steuer_MESG,miett99_MESG,ehrgott05_MESG}. In this context, the concept of  Pareto efficient solution is fundamental: it is a solution for which one cannot improve one objective without deteriorating some others. Several algorithms have been proposed for Pareto set generation and among them let us remember the weighted sum method, the lexicographic method, the achievement secularising function, the epsilon constraint method (for surveys see \citeauthor{marler2004survey}, \citeyear{marler2004survey}, or Chapters 18 and 19 in \citeauthor{greco2016multiple_MESG}, \citeyear{greco2016multiple_MESG}). Dealing with a multiobjective optimisation problem, it is important to discover the set of Pareto efficient solutions most preferred by the decision maker. Recently, beyond many exact methods, some heuristic methods have been proposed for these problems, such as some hybridisation between evolutionary multiobjective optimisation algorithms aiming to approximate the whole set of Pareto efficient solutions \citep{kdeb01_MESG} and some multicriteria preference elicitation methods to guide the optimisation algorithm toward the most interesting set of Pareto efficient solutions \citep[see, e.g.,][]{phelps2003interactive_MESG,branke2016using_MESG}.  
\end{itemize}

\subsection[Dynamic programming (Dong~Li \& Li~Ding)]{Dynamic programming\protect\footnote{This subsection was written by Dong~Li and Li~Ding.}}
\label{sec:Dynamic_programming}

Dynamic programming (DP) was the brainchild of Richard Bellman \citep{bellman1953introduction_DLLD}, who wrote ``DP is a mathematical theory devoted to the study of multistage processes''. Indeed, in the seven decades since his seminal work, the uses of DP have grown substantially thanks to its algorithmic nature in solving sequential decision-making problems, where the preceding actions and their realisation (in terms of consequences) will impact on the course of futures. Examples of such problems include multiperiod inventory management, or asset allocation (portfolio management) over a given time horizon. 
The central idea of DP is to break down the original multistage problem into a number of tail sub-problems by stages. For each stage, the tail sub-problem is a truncated version of the original problem starting from this stage. These tail sub-problems are then recursively solved one by one from the last stage backwards to the first one, at which point the original problem is solved. The solution of such a procedure is guaranteed to be optimal when the problem concerned satisfies a sufficient condition, i.e., the \textit{Principle of Optimality} \citep{bellman1953introduction_DLLD,Puterman2005-lr_HL}, which states “an optimal policy has the property that whatever the initial state and initial decision are, the remaining decisions must constitute an optimal policy with regard to the state resulting from the first decision” \citep{bellman1953introduction_DLLD}. Throughout this section, we focus our attention on discrete time systems. For continuous time dynamic systems, the readers are referred to the Hamilton-Jacobi-Bellman equations in optimal control \citep[see, for example,][]{bertsekas2012dynamic_DLLD}.

In particular, for a finite time horizon problem, the decisions are made over a number of \textit{stages} or \textit{decision epochs}, denoted by $t=0,\dots,T-1$. At each decision epoch, after observing the current system \textit{state} $x_t$ (comprised of one or more information variables that characterise how the system progresses), an \textit{action} $a_t$ is taken that leads to an immediate \textit{reward (cost)} of $r_t(x_t,a_t,w_t)$, where $w_t$ is the random disturbance at time $t$ with a known probability distribution. The system then evolves to state $x_{t+1}$ at the next decision epoch, following the \textit{transition function} $x_{t+1}=f_t(x_t, a_t, w_t)$ with the transition probability $p_t(x_{t+1}|x_t, a_t, w_t)$. After the last decision is made at epoch $T-1$, the system evolves to $x_T$ in the terminal stage with the salvage value $r_T(x_T)$. The objective of the problem is to find a \textit{policy} $\pi$, or a sequence of actions $(a_0, a_1,\cdots, a_{T-1})$ prescribed by $a_t = \pi(x_t)$, that maximises (minimises) the total expected reward (cost) across the entire time horizon. Note that for the expected total reward optimisation criterion (or additive reward functions) the Principle of Optimality is always satisfied \citep{Puterman2005-lr_HL}. To avoid the technical subtleties, in what follows we focus on discrete state space $S$ and action space $A$, and assume the random disturbance at an epoch is independent of those in the previous epochs. Define the dimension of the state space $S$ as the number of the information variables in the state. The mathematically inclined readers are referred to \cite{Puterman2005-lr_HL} for discussions on more general situations. Before proceeding, it is worth mentioning that when the random disturbance $w_t$ takes only a single value, the problem reduces to a deterministic problem. Perhaps the two most well known deterministic sequential decision-making problems solvable by DP are the Shortest Path problem \citep{dreyfus1969appraisal_DLLD} and the Knapsack problem \citep{KPP04_SMPT}.

Under the Principle of Optimality, the above-mentioned problem can be solved by backward induction. Denote by $V_t(x_t)$ the value function, or the optimal expected value-to-go from state $x_t$ at epoch $t$ until the end of the time horizon. The value function (for maximisation problems) satisfies the following optimality equations \citep[or Bellman Equations, see e.g.,][]{Puterman2005-lr_HL},
\begin{equation}
    V_t(x_t) = \max_{a_t \in A}\mathbb{E}\left[r_t(x_t,a_t,w_t) + V_{t+1}(f_t(x_t, a_t, w_t)) \right],
    \forall x_t \in S, t=0,\cdots,T-1, \label{Bellman}
\end{equation}
with the boundary condition $V_T(x_T) = r_T(x_T)$. 
By recursively solving the optimality equations from the last stage backwards to time zero, we obtain the optimal value functions and, at the same time, an optimal policy. For this method to work, however, at each stage one has to solve the value function for all states before proceeding to the previous stage. For problems with high dimensional state variables, the solution via this method is simply not practical due to the prohibitive amount of computational time and memory required. The recent development on DP research has been essentially trying to overcome this so called \textit{curse of dimensionality} \citep{Powell2011-km_HL}, which is discussed in the last paragraph of this section.

Many sequential decision-making problems in practice do not have a natural termination stage, leading to a rich body of literature studying infinite horizon problems, for which the total expected reward becomes unbounded as the time horizon tends to infinity. To this end, two alternative criteria have been widely used in the literature  \citep{Puterman2005-lr_HL, bertsekas2012dynamic_DLLD}. The first one applies a discount factor between 0 and 1, say $\beta$, to the future reward, which can be understood as the depreciation of monetary values over time. The total discounted reward is well defined as it is bounded by the sum of a decreasing infinite geometric sequence. In situations where discounting is not appropriate, a meaningful criterion is to consider the long run average reward, or the reward rate per stage. Assuming a stationary system (in which the transition function/probability, the reward function, and random disturbance do not change over time), the Bellman Equations for the total discounted reward criterion take the following form:
\begin{equation}
    V(x) = \max_{a \in A}\mathbb{E}\left[r(x,a,w) + \beta V(f(x, a, w)) \right], \forall x \in S,
    \label{BellmanInfinite}
\end{equation}
where the value function $V(x)$ is the optimal discounted value-to-go from state $x$ over an infinite time horizon. Note that there are no more boundary conditions. There is no more dependency on time either under the assumption of stationary systems, which is often satisfied in practice \citep{bertsekas2012dynamic_DLLD}. When such an assumption is not satisfied, a periodic or cyclic DP can be developed \citep{li2022bid_DLLD}. For brevity we do not include the Bellman equations for the long run average reward criterion but direct the readers to \cite{bertsekas2012dynamic_DLLD} and \cite{Puterman2005-lr_HL}.

There are mainly three solution algorithms \citep{tijms1994stochastic_DLLD, Puterman2005-lr_HL} for infinite horizon problems. The most widely used and understood algorithm is value iteration, or successive approximations as it was called in the early days. Starting from an arbitrary bounded value function vector (e.g., $V_0(x)=0, \forall x \in S$), this method iteratively updates value functions via the recursive equation below until the successive gaps between iterations $k+1$ and $k$ are within a predefined threshold. 
\begin{equation}
    V_{k+1}(x) = \max_{a \in A}\mathbb{E}\left[r(x,a,w) + \beta V_k(f(x, a, w)) \right], \forall x \in S.
    \label{VI}
\end{equation}
An alternative algorithm is policy iteration, which starts with an arbitrary policy and then iteratively improves it until no further improvements are possible. Each iteration includes two steps: firstly the expected value-to-go under the current policy is evaluated via a system of equations similar to (\ref{BellmanInfinite}) but for the actions prescribed by the policy; after that a policy improvement step is undertaken to find an improved action for each state that leads to a better value-to-go \citep{Puterman2005-lr_HL}. In the last algorithm, the system of Bellman Equations (\ref{BellmanInfinite}) are reformulated into a vary large scale linear program, which has one decision variable for each state and one constraint for each state-action pair. Regardless of the solution algorithms, just as in finite horizon problems, the curse of dimensionality remains the biggest hurdle for the implementation of DP. 

Various approximation methods have been proposed to improve the scalability of DP, leading to an important and thriving research field called Approximate Dynamic Programming (ADP). According to \cite{bertsekas2012dynamic_DLLD}, most of the ADP approaches fall into either the value space or policy space. We concentrate on the approaches in the value space (see also \S\ref{sec:Stochastic_models}) while we direct readers to \cite{bertsekas2012dynamic_DLLD} for the policy space counterparts. The basic idea of the value space approaches is to develop efficient methods to approximate the value functions or the expected value-to-go for a given policy. The most studied methods approximate the value functions via a linear or nonlinear combination of a set of handcrafted feature vectors (functions of the state) weighted by a set of parameters, which are calibrated by a suitable method \citep{bertsekas2012dynamic2_DLLD,ding2008allocation_DLLD}. Feature vectors are not always available, in which case Neural Networks have been used to construct feature vectors automatically \citep{Powell2011-km_HL,bertsekas2012dynamic_DLLD,he2018vehicle_DLLD}. Decomposition is also a popular method, which decomposes the original problem into a number of sub-problems each of which has a much smaller state space and can be solved efficiently by the exact algorithms mentioned above. The assembly of the value functions of these sub-problems provides an approximation to the original value functions \citep{kunnumkal2010new_DLLD, li2017dynamic_DLLD}. A distinct decomposition approach is \textit{Whittle's Restless Bandit} framework \citep{whittle1988restless_DLLD, glazebrook2014stochastic_DLLD, li2020switch_DLLD}, which decomposes the original problem via Lagrangian relaxation, calculates a state dependent index value for each sub-problem and uses these index values directly to derive policies for the original problem. Another method in the value space approximates the value functions of a specific policy via Monte Carlo simulation \citep{chang2007simulation_DLLD,bertsekas2012dynamic2_DLLD}, which are then used to find an improved policy. An alternative method is called \textit{Q-Learning} \citep{sutton2018reinforcement_LCAL}, which approximates the \textit{Q-factor} for each state-action pair. The Q-factor for $(x,a)$ is the expected value-to-go by taking action $a$ at state $x$ and then following either a given policy or the optimal policy thereafter. Due to the large number of combinations of state-action pairs, Q-Learning is more suitable for problems with a small state space \citep{bertsekas2012dynamic2_DLLD}. For an in-depth account on ADP we refer to two seminal books of \cite{Powell2011-km_HL} and \cite{bertsekas2012dynamic2_DLLD}. 

\subsection[Forecasting (Fotios~Petropoulos)]{Forecasting\protect\footnote{This subsection was written by Fotios~Petropoulos.}}
\label{sec:forecasting}

Forecasting is concerned with the prediction of unknown/future values of one or multiple variables of interest. If the values of these variables are collected over time, especially/in particular at regular intervals, the corresponding problem is referred to as time-series forecasting. The outputs of forecasting models include point estimates as well as expressions of uncertainty of such estimates in terms of probabilistic forecasts, prediction intervals, or path forecasts. Forecasting is applied in a wide range of applications. In this subsection, we offer an overview of established forecasting approaches that are useful in social settings \citep{Makridakis2020-eu_FP}, such as forecasts produced to support decision making in operations and supply chain management (\S \ref{sec:Inventory}; \S \ref{sec:Supply_chain_management}, finance (\S \ref{sec:Finance}), energy (\S \ref{sec:Power_markets_and_systems}), and other domains.

\textit{Exponential smoothing} is one of the most popular families of models for univariate time-series forecasting. The underlying principle of exponential smoothing models is that, at every step, the forecast is updated such that the most recent information is taken into account by exponentially discounting information from previous periods. The estimates for the exponential smoothing parameters are based on in-sample fits. The first and simplest exponential smoothing method, simple (or single) exponential smoothing, was developed by \cite{Brown56_FP}. This method was able to handle level-only data (no trend nor seasonal patterns). Soon after, it was extended to handle trended and seasonal data \citep{Holt20045_FP,Winters1960_FP}. Forty years later, \cite{Hyndman2002-jp_FP} introduced a fully fledged family of exponential smoothing models that are represented in a state-space framework. Usually, three states are considered: level, trend, and seasonality. The way that these three states interact to produce the final forecast determines the types of trend and seasonality (such as additive or multiplicative). Exponential smoothing models are fast to compute and perform well in a wide range of data \citep{Makridakis2019-oy_FP}, rendering them ideal benchmarks for forecasting applications. Detailed reviews of exponential smoothing models are offered by \cite{Gardner2006-bv_FP} and \cite{Hyndman2008-iu_FP}.

\textit{Autoregressive integrated moving average} (ARIMA) is another very popular family of univariate forecasting models \citep[for a seminal work on ARIMA, see][]{Box1976-af_FP}. In ARIMA, the data are first rendered stationary through transformations and differencing. The stationary data are then fitted in linear regression models (see also the next paragraph on regression models) in which the predictors are either past values of the data (autoregressive terms) or past errors (moving average terms). ARIMA models are theoretically appealing as they can depict a wide range of data generation processes. While manually identifying an optimal ARIMA model can be sometimes challenging, nowadays automated approaches exist \citep[see, for example,][]{Hyndman2008_FP,Franses2014-kl_FP}

When the variable of interest is known to be affected by other factors (also called ``exogenous variables''), then causal modelling can be applied. In its simplest form, causal models can be linear or nonlinear \textit{regression models} that regress the values of the dependent variable on the values of the independent variable(s). Apart from the ordinary least squares regression models, other types of regression models exist, such as the ordinal, logistic, Poisson, negative binomial regression models as well as the Generalised Linear Models (GLMs). The dependent variable (variable of interest) is usually continuous, however specific regression models exist for ordinal or binary dependent variables, such as the ordinal logistic regression model.” But of course there are also regression approaches for count data, like Poisson regression or negative binomial regression (Hilbe’s textbook of the same name is my go-to reference on this), or more generally Generalized Linear Models (GLMs). I would assume these to be more relevant to OR than binary or ordinal logistic regressions.

A common rule for using regression models for forecasting purposes is that the values of the independent variables are either known or can be predicted, as is very common in energy forecasting; see \cite{wer:14_DSRW} and \S\ref{sec:Power_markets_and_systems}. Transformations of the dependent or independent variables are sometimes necessary so that assumptions regarding normality of errors and constancy of the error variance are satisfied \citep{lag:mar:sch:wer:21_DSRW}. Another common issue in regression models is that of multicollinearity between independent variables. Linear regression models can also be used to produce time-series forecasts when no exogenous variables are available. In these cases, we can construct predictors for trend and seasonality and use these predictors as independent variables to model the time-series patterns. Finally, it is also worth mentioning that ARIMA models can be extended to ARIMAX models that can include the effects of exogenous variables, just like autoregressive (AR) models can be extended to ARX.

Instead of forecasting each time series separately, several approaches exist in order to forecast time series data as a collection. \textit{Multivariate models} (also known as structural models) are designed to model cross-sectional data, producing forecasts for many variables of interest at the same time. Such forecasts take into account interactions between all series. A common example is the vector autoregressive (VAR) models \citep{sims1980_FP,hasbrouck95_FP}. Another very popular cross-sectional approach is \textit{hierarchical forecasting} \citep{Athanasopoulos2020-cx_FP}. Hierarchical forecasting deals with time-series data that are naturally arranged in hierarchical structures (for example, product or geographical hierarchies). Forecasts for each node of the hierarchy are first produced independently using standard univariate forecasting approaches (such as exponential smoothing or ARIMA); then, forecasts across the hierarchy are reconciled to achieve coherency \citep{Wickramasuriya2019_FP,Hollyman2021-wj_FP}. Hierarchical forecasts offer better accuracy and are directly relevant for decision makers at multiple levels of an organisation. A different form of forecasting using multiple series, which is widely applied in machine-learning methods, is called \textit{cross-learning}. This approach implies learning \citep[usually through features;][]{Montero-Manso2020_FP,Wang2022comb_FP} from other series to be able to predict the variable of interest. Compared to other cross-sectional approaches, cross-learning requires access to a set of ``reference'' data which, though, do not have to be concurrent to the target data.

Given the plethora of available modelling options, we need ways to help us decide on the best approach for the target data. Two popular approaches for \textit{model selection} are information criteria and cross-validation. Information criteria select the best model amongst a pool of candidate models based on how well the in-sample forecasts fit the actual data (model fit), penalising at the same time for model complexity (Occam's razor). Information criteria are fast to compute and widely applied, mostly due to their implementations in open-source forecasting packages \citep{Hyndman2008_FP}. Cross-validation is based on the comparison of the out-of-sample performance between different models. To achieve this, the available data are split into ``training'' and ``validation'' data. The validation follows a rolling-origin process, where the forecasts of the candidate models are compared for multiple forecast origins \citep{tashman2000out_FP,Bergmeir2012use_FP}. A more recent approach to forecast selection is based on the concept of representativeness \citep{Petropoulos2022REP_FP}. Out-of-sample forecasts with higher representativeness to the past data patterns are preferred to ones with lower representativeness. Regardless of how one selects between forecasts and models, the values of the selection criteria can also be used to \textit{combine forecasts} \citep{kolassa2011combining_FP}. In fact, multiple studies have shown that combining forecasts, using equal or unequal weights, can significantly boost the forecasting performance of individual models \citep{Bates1969combination_FP,now:liu:wer:hon:16_DSRW,Wang2022_FP}. \cite{Claeskens2016-hr} offer a possible explanation on why the performance of forecast combinations is better than that of the individual forecasts.

Apart from statistical, algorithmic and computational approaches, the forecasting process can also be infused by \textit{judgement} (see, also, \S \ref{sec:Behavioural_OR} and \S \ref{sec:Soft_OR_and_problem_structuring_methods}). It is not unusual for forecasts to be directly produced in a judgemental way, without the support of any systematic approaches. Research suggests that such forecasts suffer from several biases \citep{Lawrence2006-ey_FP}. However, managers may sometimes have a unique appreciation of the situation, one that the hard data cannot communicate through models. In such cases, systematic approaches to elicit expert knowledge include the Delphi method \citep{Rowe1999-vs_FP}, structured analogies \citep{Green2007-ts_FP}, prediction markets \citep{Wolfers2004-mq_FP}, and interaction groups \citep{Van_de_Ven1971-vd_FP}; see also \cite{Graefe2011-ha_FP} and \cite{Nikolopoulos2015-mz_FP} for a comparison between these approaches. Apart from producing forecasts directly, judgement may also be used to adjust the formal/statistical forecasts. Judgemental interventions and their efficacy have been well-studied in the literature \citep[see, for example:][]{Fildes2009-tc_FP,Petropoulos2016-uk_FP,Fildes2019-uj_FP}. The main takeaways are: (\textit{i}) negative adjustments are generally more beneficial than positive ones; (\textit{ii}) larger adjustments should be preferred to smaller ones; and (\textit{iii}) the use of feedback and support will limit and improve the role of such judgemental adjustments. Finally, managerial judgement may be applied in other stages of the forecasting process, such as judgementally selecting between statistical models \citep{Petropoulos2018-mt_FP,De_Baets2020-du_FP} or setting their (hyper)parameters.

Forecasts produced in previous periods need to be evaluated once the corresponding actual values become available. Through feedback, \textit{forecast evaluation} allows analysts  to improve the forecasting process and, thus, forecasting performance. The main rule of forecast evaluation is that performance should be measured on data that were not used to fit the models or produce the forecasts. Solely measuring the in-sample performance will inevitably lead to over-fitting and the use of complex forecasting models. There exist a wide array of evaluation metrics. Some of them are suitable for measuring the accuracy or bias of the point forecasts, while others focus on how well the uncertainty around the forecasts is estimated. In the former category, popular metrics are the mean absolute error (MAE), the root mean squared error (RMSE), the mean absolute percentage error (MAPE, which is very popular in practice) and the mean absolute scaled error (MASE, which is theoretically more elegant and popular in academia). It should be noted that \cite{Kolassa2020-jh} showed that different error metrics are minimised by different (point) forecasts and that it makes little sense to evaluate one point forecast using multiple KPIs. For detailed overviews of forecasting metrics for point forecasts and their proper use, the reader is referred to \cite{HYNDMAN2006_FP}, \cite{Davydenko2013_FP}, and \cite{Koutsandreas2022_FP}. In the latter category, a popular metric is the interval score (IS), which is a proper scoring rule and considers both the calibration and the sharpness of the prediction intervals, as well as the pinball score, the continuous ranked probability score (CRPS), and the energy score. \cite{Gneiting2007-oj_FP} offer a review of (strictly) proper scoring rules. Finally, we should mention that nowadays it is common to go beyond strict forecasting performance and measure the performance of the forecasts on the utility \citep{hon:pin:etal:20_DSRW,Yardley2021_FP}.

For a detailed encyclopedic overview of the forecasting field, both in terms of theory and practice, we refer the reader to the work of \cite{PETROPOULOS2022705_FP}; a live version of this encyclopedia is available at \href{https://forecasting-encyclopedia.com/}{forecasting-encyclopedia.com}. \cite{HyndmanAthanasopoulos2021_FP} and \cite{FildesOrd2017_FP} have written comprehensive textbooks on forecasting and its applications. Notable open-source packages with implementations of the most popular forecasting models include the \texttt{forecast} \citep{Hyndman2022Rforecast_FP} and \texttt{smooth} \citep{Svetunkov2022smooth_FP} packages for R statistical software.

\subsection[Game theory (Georges~Zaccour)]{Game theory\protect\footnote{This subsection was written by Georges~Zaccour.}}
\label{sec:Game_theory}

Game theory is a branch of mathematics that studies strategic interactions between decision makers, called players. Strategic interactions means that a player's payoff depends not only on her own decision (action or choice), but also on the decisions made by the other players. The book by \cite{Vomo1944_GZ} is often considered the starting date of game theory, though some of its roots can be traced back to much earlier. Games can be classified along a series of features. In a static game, each player acts only once, whereas in a dynamic game, interactions are repeated over time. In a one-person game, the decision maker plays against a nonstrategic (or dummy) player, often referred to as ``nature'', whose action is the outcome of a probabilistic event with a fixed (known) distribution. Two-player games focus on one-on-one interactions. Duopolistic competition and management-union negotiations are situations that can be modelled as two-person games. Extending the model to $n>2$ players is conceptually easy but may be computationally challenging because each player needs to determine all the possible sequences of actions and reactions for all players. When the number of interacting players is very large, e.g., an economy with many small agents, the analysis shifts from individual-level decisions to understanding the group's behavioural dynamics. An illustration of this is traffic congestion: when an agent attempts to minimise her travel time on a route from A to B, her travel speed depends on the traffic density on that route. What matters is the number of drivers, not their identity. Population and evolutionary games \citep{Hosi1998_GZ,Cr2003_GZ,Sa2010_GZ} and mean-field games \citep{Huetal2003_GZ,Huetal2006_GZ,Huetal2007_GZ,Lali2006a_GZ,Lali2006b_GZ,Lali2007_GZ,Gosa2014_GZ} are branches of game theory that study situations with large numbers of players.

A game can be defined in three forms, namely, in strategic, extensive, or in coalitional form. To formulate a one-shot game in strategic form, we have to specify (\textit{i}) the set of \textit{players} and, for each player, (\textit{ii}) the set of \textit{actions}, and (\textit{iii}) a \textit{payoff function} measuring the desirability of the game's possible outcomes, which depends on the actions chosen by all players. The set of actions can be finite, e.g., to bid on a contract or not, or continuous, e.g., the amount bid. If the players intervene more than once in the game, then we should additionally define (\textit{iv}) the order of play, (\textit{v}) the information acquired by the players over time (stages), and (\textit{vi}) whether or not nature is involved in the game.

In a one-shot game, an action (move) and a strategy mean the same thing. In games where players intervene more than once, the two concepts no longer coincide. A strategy is then a \textit{decision rule} that associates a player's action with the information available to her at the time she selects her move. So an action, e.g., spending advertising dollars, is a result of the strategy. The word strategy comes from Greek (strategia) and has a military sense. An army general's main task is to design a plan that takes into account (adapts to) all possible contingencies. This is precisely the meaning of strategy in game theory. Whether in war, business or politics, it is never wise to allow yourself to be surprised by the enemy. This does not imply that a winning strategy always exists. Sometimes we must be content with a draw or even a reasonable loss.

One-shot games are a useful representation of strategic interactions when the past and the future are irrelevant to the analysis. However, if today's decisions also affect future outcomes and are dependent on past moves, then a dynamic game is needed. In a \textit{repeated game}, the agents play the same game in each round, that is, the set of actions and the payoff structures are the same in all stages \citep{Metal2015_GZ}. The number of stages can be finite or infinite, and this distinction matters in terms of achievable outcomes. In a \textit{stochastic game}, the transition between states depends on the players' actions \citep{Sh1953_GZ,Mene1981_GZ,Jano2018b_GZ,Jano2018a_GZ}. In a \textit{multistage game}, the players share the control of a discrete-time dynamic system (state equations) observed over stages \citep{Baol1999_GZ,En2005_GZ,Krpe2018_GZ}. Their choice of control levels, e.g., investments in production capacity or advertising, affects the evolution of the state variables (e.g., production capacity, reputation of the firm), as well as current payoffs. \textit{Differential games} are continuous-time counterparts of multistage games \citep{Is1975_GZ,Baetal2018_GZ}.

Information plays an important role in any decision process. In a game, the \textit{information structure} refers to what the players know about the game and its history when they choose an action. A player has \textit{complete information} if she knows who the players are, which set of actions is available to each one, what each player's information structure is, and what the players' possible outcomes can be. Otherwise, the player has \textit{incomplete information}. If, for instance, competing firms do not know their rivals' production costs, then the game is an incomplete-information game. The game can also have perfect or imperfect information. Roughly speaking, in a game of \textit{perfect information}, each player knows the other players' moves when she chooses her own action, as in, e.g., chess or a manufacturer-retailer game where the upstream player first announces a product's wholesale price, and then the downstream player reacts by selecting the retail price. The archetype of an imperfect-information game is the prisoner's dilemma, where (in the original story) the players have to simultaneously choose between confessing or denying a crime. A Cournot oligopoly, where each firm chooses its own production level without knowing its competitors' choices, is another instance of an imperfect-information game.

The outcome of a game depends on the players' behaviour. In a noncooperative game, e.g., R\&D competition to develop a vaccine, each player optimises her own payoff, whereas in a cooperative game, the players seek a collectively optimal solution. For instance, the members of a supply chain could agree to coordinate their strategies to maximise the chain's total profit. The
fundamental solution concept in a noncooperative game is the \textit{Nash equilibrium} \citep{Na1950_GZ,Na1951_GZ}. Let $I=\left\{ 1,\ldots ,n\right\} $ be the set of players, $S_{i}$ the set of strategies of player $i\in I$, and let her payoff function be given by $g_{i}\left( s_{1},\ldots,s_{n}\right)$ $:\prod\limits_{i\in I}S_{i}\rightarrow \mathbb{R}$, where $s=\left( s_{1},\ldots ,s_{n}\right) $. Assuming the players are maximisers, the strategy profile $s^{N}=\left(s_{1}^{N},\ldots ,s_{n}^{N}\right) $ is a Nash equilibrium if $g_{i}\left( s_{1}^{N},\ldots ,s_{n}^{N}\right) \geq $ $g_{i}\left( s_{1}^{N},\ldots s_{i-1}^{N},s_{i},s_{i+1}^{N}\ldots ,s_{n}^{N}\right) $ for all $s_{i}\in S_{i}$ and all $i\in I$. At an equilibrium, no player has an interest in deviating unilaterally to any other admissible strategy. Put differently, if all other players stick to their equilibrium values, then player $i$ does not regret implementing her equilibrium value too, which is obtained by best-replying to the choice of the others. That is, $s_{i}^{N}=\arg \max_{s_{i}\in S_{i}}g_{i}\left( s_{1}^{N},\ldots s_{i-1}^{N},s_{i},s_{i+1}^{N}\ldots ,s_{n}^{N}\right) .$ A Nash equilibrium does not always exist, and there may be multiple equilibria, raising the question of which one to select \citep{Selten1975-rl_GZ}. Existence and uniqueness conditions for Nash equilibrium typically rely on fixed-point theorems. If the game is one of incomplete information, then the solution concept is a \textit{Bayesian Nash equilibrium} \citep{Ha1967_GZ_1,Ha1967_GZ_2,Ha1967_GZ_3}. Another noncooperative equilibrium solution concept, which predates the Nash equilibrium, is the \textit{Stackelberg equilibrium}, introduced in a two-player framework by \citet{vo1934_GZ}. There is a hierarchy in decision-making between the two players: the \textit{leader} first announces her action, and next the follower makes a decision that takes the leader's action as given. Before announcing her action, the leader would of course anticipate the follower's response and selects the action that gives her the most favourable outcome. The framework has been extended to several followers and leaders \citep{Sh1984_GZ}.

In a \textit{cooperative game}, the players coordinate their strategies in view of optimising a collective outcome, e.g., a weighted sum of their payoffs, and must agree on how to share the dividend of their cooperation \citep{Mo1988_GZ,Ow1995_GZ}. Different solution concepts have been proposed, each based on some desirable properties, typically stated as axioms, such as fairness, uniqueness of allocation, and stability of cooperation. The most-used solutions in applications are the \textit{core} \citep{Gi1953_GZ}, and the \textit{Shapley value} \citep{Sh1953_GZ}. In any solution, the set of acceptable allocations only includes those that are individually rational. Individual rationality means that a player will agree to cooperate only if she can get a better outcome in the cooperative agreement than she would by acting alone. In a dynamic cooperative game, the agreement must specify, at the outset, the decisions that must be implemented by each player throughout the planning horizon. One concern in such games is the durability of the agreement over time. Clearly, it is rational for a player to leave the agreement at an intermediate time if she can achieve a better outcome. The literature on dynamic games has followed two streams in its quest to sustain cooperation over time, namely, building cooperative equilibria or defining time-consistent solutions. Through the implementation of some (punishing) strategies, the first stream seeks to make the cooperative solution an equilibrium of an associated noncooperative game. If this is achieved, then the result will be at once collectively optimal and stable, as no player will find it optimal to deviate unilaterally from the equilibrium. See \citet{Osru1994_GZ} for repeated games, \citet{Du1995_GZ} and \citet{Paza2015_GZ} for different types of stochastic games; and \cite{Hato1985_GZ}, \cite{Toetal1986_GZ}, and \citet{Hapo1987_GZ} for multistage and differential games. The second stream looks for time-consistent solutions, which are achieved by allocating the cooperative payoffs over time in such a way that, along the cooperative state trajectory, no player will find it optimal to switch to her noncooperative strategies. The idea was initiated in \cite{Pe1977_GZ} and has since been further developed \citep[see][]{Yepe2018_GZ,Peza2018_GZ}.

Game theory has found applications in biology, economics, engineering, management, Operational Research, and political and social sciences.

\subsection[Graphs and networks (Ivana~Ljubi\'c)]{Graphs and networks\protect\footnote{This subsection was written by Ivana~Ljubi\'c.}}
\label{sec:Graphs_and_networks}

Graphs and networks are used to represent interactions, connections or relationships between objects. In network optimisation problems, numerical attributes representing features such as costs, weights or capacities are assigned to objects (also called \emph{vertices}) or to connections between them. If connections are directed, we refer to them as \emph{arcs}, otherwise we call them \emph{edges}. Given an input graph with $n$ vertices and $m$ arcs (or edges), the goal is to find a subgraph that exhibits desired properties (described by a given set of constraints) and that optimises the given objective function (usually measured as the sum of edge or vertex ``weights'' of the solution's subgraph). In the following, we focus on some of the most fundamental and most studied problems in network optimisation. 

The \emph{shortest path problem} in arc-weighted graphs, for example, seeks to find a least costly path from the given source vertex $s$ to the given target $t$. When the arc costs are non-negative, one can use the algorithm of \cite{Dijkstra:1959_IL}, the efficient implementation of which uses Fibonacci heaps and runs in $\mathcal{O}(m + n\log n)$ time. For graphs with possible negative arc costs, in $\mathcal{O}(mn)$ time the Bellman-Ford algorithm either finds the shortest path from $s$ to all other vertices, or it proves that such a path does not exist due to the presence of a negative cost cycle reachable from $s$. The shortest path algorithms are explained in many textbooks, see e.g., \cite{cormen01introduction_IL,Kleinberg-Tardos:2006_IL,Schrijver:2003_IL,Williamson:2019_IL}. 
 
In the \emph{maximum flow problem} (MF), in a given network with arc capacities, we want to send as much flow as possible from the given source $s$ to the given sink $t$ without violating the arc capacities. The problem was motivated by the conflict between East and West during the Cold War \citep{schrijver2002history_JLYHK}. \cite{Ford-Fulkerson:1957_IL} develop the first exact algorithm that searches for augmenting paths in the \emph{residual network}. Their fundamental result, known as the \emph{max-flow/min-cut theorem} states that the maximum flow passing from the source to the sink is equal to the total capacity of the arcs in a \emph{minimum cut}, i.e., the network that indicates how much more flow is allowed in each arc., which is a subset of arcs of the smallest total capacity, the removal of which disconnects the source from the sink. The same result using the duality theory of LPs is given in \cite{Dantzig-Fulkerson:1955_IL}. The famous results from graph theory such as Menger's theorem, K\"{o}nig-Egeváry theorem, or Hall's theorem, follow from the max-flow/min-cut theorem \citep{Ford-Fulkerson:1962_IL}. The method of \cite{Ford-Fulkerson:1957_IL} is pseudo-polynomial when arc capacities are integral, however it may fail to find the optimal solution and need not terminate if some of the arc capacities are irrational \citep{Ford-Fulkerson:1962_IL}. An algorithm that overcomes this issue was independently discovered in the 1970s by \cite{Edmonds-Karp:1972_IL} and  \cite{Dinic:1970_IL}, see also \cite{Dinitz:2006_IL}.  Augmenting the flow along shortest paths (that is, along the paths with fewest edges) guarantees a polynomial-time complexity. Instead of augmenting the flow along a single augmenting path as in \cite{Edmonds-Karp:1972_IL}, the algorithm of \cite{Dinic:1970_IL} finds all shortest augmenting paths in a single phase. Another stream of MF algorithms exploits the \emph{preflow} idea of \cite{Karzanov:1974_IL} in which the vertices are overloaded with the excess flow (i.e., more incoming flow than the outgoing flow is allowed). Subsequent improvements are obtained in the following years. An important breakthrough is achieved by Goldberg \& Tarjan with the introduction of \emph{push-relabel} algorithms \citep{Goldberg-Tarjan:1988_IL}. A \emph{pseudoflow} algorithm for the maximum flow is introduced by  \cite{Hochbaum:2008_IL} and it is later improved in \cite{Hochbaum-Orlin:2013_IL}. The recent implementation by  \cite{Goldberg-et-al:2015_IL} is competitive with the  \cite{Boykov-Kolmogorov:2004_IL} method and the pseudoflow approach. Further historical details and a more complete list of references can be found in \cite{ahuja1993network_IL,Dinitz:2006_IL,Goldberg-Tarjan:2014_IL,Williamson:2019_IL}. Currently, the best strongly polynomial bounds are obtained by  \cite{Orlin:2013_IL} and  \cite{King-et-al:1994_IL}. However, new and improved MF algorithms continue to be discovered. The most recent trends use the idea of electrical flows for obtaining faster (exact or approximate) algorithms, see, e.g., Chapter 8 of \cite{Williamson:2019_IL}. 

In the \emph{minimum cost flow problem} (MCF), for each arc of the graph, a cost is incurred per unit of flow that traverses it. The goal is to send units of a good that reside at one or more \emph{supply} vertices to some other \emph{demand} vertices, without violating the given arc capacities at minimum possible cost. \cite{Edmonds-Karp:1972_IL} introduce the \emph{scaling technique} for the MCF. The technique is later improved by \cite{Orlin:1993_IL}. The algorithms of \cite{Vygen:2002_IL} and \cite{Orlin:1993_IL} have the best-known strongly polynomial complexity bounds for the MCF. \cite{Kovacs:2015_IL} provides a comprehensive literature overview and gives an experimental evaluation of MCF algorithms based on \emph{network-simplex, scaling} or \emph{cycle-cancelling} techniques. The MCF is treated in detail in many textbooks,  \cite{ahuja1993network_IL,Korte-Vygen:2008_IL,Williamson:2019_IL}. One of the important results is the  \emph{integrality of flow} property: if all demands/supplies and arc capacities are integers, then there exists an optimal MCF solution with integer flow on each arc. The result follows from the \emph{totally unimodular} property of the constraint matrix when the MCF is modelled as a linear program. 

In the \emph{minimum cut problem} (MC), one searches for a proper subset of vertices $S$ of a given arc-capacitated graph, such that the total capacity of arcs leaving $S$ is minimised.  For directed graphs, the algorithm of \cite{Hao-Orlin:1994_IL} is based on MF calculations between chosen pairs/subsets of vertices and exploits the push-relabel ideas. For undirected graphs, the \emph{Gomory–Hu tree}, which is a weighted tree that represents the minimum $s$-$t$ cuts for every $s$-$t$ pair in the graph, is introduced in \cite{Gomory-Hu:1961_IL}. This tree is constructed after $n-1$ MF computations, and a simpler procedure has been later given by  \cite{Gusfield:1990_IL}. The algorithm of   \cite{Padberg-Rinaldi:1990_IL} improves the ideas of \cite{Gomory-Hu:1961_IL}  and is widely used within branch-and-cut schemes for solving the travelling salesperson problem (TSP) and related problems. The \emph{maximum adjacency ordering} together with Fibonacci heaps is used in \cite{Nagamochi-et-al:1994_IL}. Randomised approaches can be found in \cite{Karger-Stein:1996_IL,Karger:2000_IL}. The method of   \cite{Karger:2000_IL} is de-randomised by   \cite{Li:2021_IL}. Practical performance of some of these algorithms is evaluated in \cite{Chekuri-et-al:1997_IL,Junger-et-al:2000_IL}.  For additional and more recent references, see the book by  \cite{Williamson:2019_IL}.

The problems mentioned so far all belong to the class $\mathcal{P}$, however most of the network optimisation problems that are relevant for practical applications are $\mathcal{NP}$-hard. We highlight two of them that serve as drivers for discovering new algorithms and methodologies that can be easily adapted to other difficult optimisation problems.

Given an undirected graph with non-negative edge costs, the \emph{Steiner tree problem in graphs} (STP) asks for finding a subtree that interconnects a given set of vertices (referred to as \emph{terminals}) at minimum cost. Two special cases can be solved in polynomial time: when all vertices are terminals (the minimum spanning tree problem), or when there are only two terminals (the shortest path problem). In general, however, the decision version of the STP is $\mathcal{NP}$-complete \citep{Karp:1972_IL}. Older surveys covering developments of first MIP formulations, Lagrangian relaxations, branch-and-bound methods and heuristics can be found in \cite{Maculan:1987_IL,Winter:1987_IL}. The research on the STP was marked by polyhedral studies in the 1990s \citep{Goemans:1994related_IL,Chopra:1994part1_IL}. Exact solution methods for the STP are based on a sophisticated combination of: \emph{reduction techniques} \citep{gamrathscip_IL,Rehfeldt-Koch:2021_IL}, \emph{dual and primal heuristics} \citep{Pajor:2018_IL} embedded within branch-and-cut or branch-and-bound frameworks, see \citep{Polzin:2003_IL,Daneshmand:2003_IL,Polzin:2009approaches_IL,gamrathscip_IL,Fischetti:2016_IL}. Currently best approximation ratio for the STP is 1.39 \citep{Goemans:2012matroids_IL}. A comprehensive survey of the results obtained in the last three decades is given by  \cite{Ljubic:2021_IL}. State-of-the-art computational techniques for the STP are due to   \cite{Rehfeldt:2021_IL}.  

The \emph{Travelling salesperson problem} (TSP) aims at finding the answer to the following question: If a travelling salesperson wishes to visit all $n$ cities from a given list exactly once, and then return to the home city, what is the cheapest route they need to take? For the history of the problem, see \cite{Applegate-et-al:2011_IL} and the book by \cite{Cook:2011_IL}. Since 1954, when \cite{DFJ:1954_IL} found a provably optimal solution for a 49-city problem instance, many important improvements in the development of exact methods have been achieved\footnote{\url{www.math.uwaterloo.ca/tsp/history/milestone.html}}. Facet-defining inequalities are investigated in  \cite{Padberg-Rinaldi:1990a_IL,JungerReineltRinaldi:1995_IL}. MIP formulations, including the famous \emph{subtour-elimination constraints} model by  \cite{DFJ:1954_IL},  are compared in \cite{Padberg-Sung:1991_IL}. Branch-and-cut methods are developed in   \cite{Applegate-et-al:2011_IL,JungerReineltRinaldi:1995_IL,Padberg-Rinaldi:1991_IL}. For the most recent overview on approximation algorithms for the TSP see \cite{Traub:2020_IL}. \cite{Helsgaun:2000_IL}\footnote{\url{akira.ruc.dk/~keld/research/LKH/}} provides an efficient implementation of the  \emph{$k$-opt heuristic} of \cite{lin1973effective_COIT}. \cite{Cook-et-al:2022_IL} extend the algorithm of \cite{Helsgaun:2000_IL} to deal with additional constraints in routing applications and win the Amazon Last Mile Routing Challenge in 2021. The TSP solver Concorde\footnote{\url{www.math.uwaterloo.ca/tsp/concorde}} \citep{Applegate-et-al:2011_IL} incorporates best algorithmic ideas from the past 60 years of  research on the topic. By combining techniques of \cite{Helsgaun:2000_IL} and \cite{Applegate-et-al:2011_IL}, instances with millions of vertices can be solved to within 1\% of optimality, see e.g., TSP solutions on graphs with up to 1.33 billion of vertices\footnote{\url{www.math.uwaterloo.ca/tsp/star}}.

\subsection[Heuristics (Ceyda~Oğuz \& İstenç~Tarhan)]{Heuristics\protect\footnote{This subsection was written by Ceyda~Oğuz and İstenç~Tarhan.}}
\label{sec:Heuristics}

Etymologically meaning to find/discover, heuristics make use of previous experience and intuition to solve a problem. A heuristic algorithm is designed to solve a problem in a shorter time than exact methods, by using different techniques ranging from simple greedy rules to complex structures, which could be dependent on the problem characteristics; however it does not guarantee to find the optimal solution (\S\ref{sec:Combinatorial_optimisation}; \S\ref{sec:Dynamic_programming}). Heuristics have been used in the operational research area extensively with respect to the applications (see, for example, \S\ref{sec:Inventory}, \S\ref{sec:Logistics}, \S\ref{sec:Manufacturing}, and \S\ref{sec:Transportation_Vehicle_routing}). In this subsection, we review the methods employed in the development of heuristics. 

Classifications and strategies provided in the literature guide us for the methods employed in heuristics. Below we provide a thorough classification and explain briefly the basic methods used under each class.

\textit{Induction}, being the simplest method to be applied with an analogy to the mathematical induction, is to solve the original complex problem by extending the results and insights obtained from small and simpler versions of the problem \citep{silver1980tutorial_COIT, silver2004overview_COIT,laguna2013heuristics_COIT}.

\textit{Restriction} methods primarily focus on explicitly eliminating some parts of the solution space so that the problem will be solved given a restricted set of solutions \citep{silver1980tutorial_COIT,zanakis1989heuristic_COIT,silver2004overview_COIT,laguna2013heuristics_COIT}. One way of doing this is to identify common attributes of the optimal solution and search among the solutions having these attributes only \citep{glover1977heuristics_COIT}. Another restriction can be applied by eliminating infeasible solutions considering a combination of decision variables which dictates incompatible values. Beam Search \citep{morton1993heuristic_COIT} is a good example of this class of heuristics which works with a truncated tree structure using strategies similar to a branch-and-bound algorithm (\S\ref{sec:Combinatorial_optimisation}). The trimming of the tree is utilised by a parameter called beam width to indicate how many nodes to have at every level of the tree. 

Heuristics using \textit{decomposition/partitioning} method employ different approaches to divide the problem into smaller and tractable parts, solve these parts separately and combine their solutions to give the solution to the original problem \citep{foulds1983heuristic_COIT,zanakis1989heuristic_COIT,silver2004overview_COIT,laguna2013heuristics_COIT}. The methods used to divide and then combine the solutions are usually dictated by the nature of the problem. For example, Hierarchical Planning proposed by \cite{hax1973hierarchical_COIT} considers the organisational level breakdown and the output of one decomposed problem becomes the input for the other. Rolling horizon also falls under this category \citep{stadtler2003multilevel_COIT}. A problem with a sequence of decisions that span a long planning horizon is solved by dividing the planning horizon into smaller planning intervals. The problem with these small planning intervals is solved continually by fixing the decisions for the first time period and moving into the next time period to solve the next problem. Another approach takes the characteristics of the input data into account and divides the problem such that each part includes only tractable amount of data. For example, data showing clusters of geographically close customers is suitable for this type of partitioning. The decomposed problems are solved independently, and their solutions are combined with a certain rule. Divide and Conquer algorithm heuristically clusters vertices on a given graph, generates a smaller graph for each cluster and solves the original problem for each cluster independently \citep{akhmedov2016divide_COIT}. Decomposition can be made based on an element of the problem, for example solving a logistics problem after dividing it into parts per vehicle. Other decomposition approaches benefit from the structure of the mathematical model developed for the problem. Examples of this sort are Lagrangian Relaxation \citep{fisher1981lagrangian_COIT}, in which complicating constraints are lifted to the objective function with a penalty, and Benders Decomposition \citep{Benders1962-sf_HL,rahmaniani2017benders_COIT}, in which once complicating variables are fixed, the remaining problem can be divided into problems to be solved independently.  

\textit{Approximation} methods focus on the mathematical models and utilise different strategies to make the problem tractable which results in a reduced size of the problem \citep{silver1980tutorial_COIT,silver2004overview_COIT}. One strategy widely used is the aggregation over variables or stages. Another common strategy is to modify the variables, the objective function or the constraints of the mathematical model in different ways, such as converting discrete variables into continuous variables, using a linear objective function instead of the non-linear objective function,, linearising nonlinear constraints, and either eliminating or weakening some of the constraints \citep{glover1977heuristics_COIT}. Kernel Search \citep{angelelli2010kernel_COIT}, which combines relaxation with decomposition over the decision variables, demonstrates that a heuristic may use more than one class of methods in its design. 

\textit{Constructive} heuristics start from an empty solution and build a complete solution by adding an element of the problem following a rule at every step, such as the nearest neighbour algorithm \citep{bellmore1968traveling_COIT} for the travelling salesman problem. Usually constructive heuristics are of greedy nature by making the decision for local optimum in every step. These algorithms can be enhanced by adding a look-ahead mechanism that is by estimating the future effects of a decision rather than just the current effect to avoid pitfalls of being greedy.

\textit{Improvement} heuristics start with a complete solution and improve it by modifying one or more elements of the solution in every iteration until a predetermined stopping condition is achieved. Improvement heuristics in their simplest form utilise a local search which is defined over a neighbourhood structure to express how the moves are performed from one iteration to the next. $k$-opt is an example of this sort which replaces $k$ elements of a solution with another set of $k$ elements in every step if it is beneficial \citep{lin1973effective_COIT}. The parameter $k$ determines the size of the local search and implicitly applies the restriction method discussed above. A neighbourhood is defined by a set of solutions which are reachable form the current solution. A local search is performed by moving from the current solution to another solution in the neighbourhood of it (next solution). Selection of the next solution is done by accepting either the one among random choices that improves the objective function value first (random descent if it is a minimisation problem) or the one resulting in the best objective function, i.e., the local optimum, with respect to that neighbourhood (steepest descent for a minimisation problem). This simple structure focuses on the local information (exploitation of the accumulated search experience) and is known as intensification \citep{glover1990tabu_COIT}. While it will be useful if the structure of the problem is appropriate, it may result in not good enough solutions otherwise. Hence, the improvement heuristic will benefit if it can explore other parts of the solution space, which is known as diversification \citep{glover1990tabu_COIT}. Two immediate strategies to be employed are either to start the search from different initial solutions and choose among the final solutions obtained (multi-start algorithms) or to allow moving to worse solutions if this direction will provide a better path for the future selections (hill-climbing strategy for a minimisation problem). 

Even though \textit{metaheuristics} \citep{glover1986future_COIT} are improvement methods, since they advance notably, considering them as a separate class is worthwhile. Metaheuristics utilise a local search together with intensification and diversification mechanisms and aim at eliminating the problem-dependent and domain-specific nature of other heuristics. Simulated Annealing \citep{kirkpatrick1983optimization_COIT} is one of the most popular metaheuristics which uses a single solution in its local search with a random descent and utilises hill-climbing strategy for diversification. Tabu Search \citep{glover1986future_COIT} is an example of deterministic metaheuristic working with a single solution throughout the search. It explicitly uses history of search in both intensification and diversification mechanisms. Genetic Algorithm \citep{Holland1975-qj_HL} is another popular metaheuristic comprising of random components for intensification and diversification but working with a set of solutions during the search. Variable Neighbourhood Search \citep{hansen1999introduction_COIT} is an excellent example of a design in which diversification is provided by systematically changing neighbourhood structures.     

\textit{Matheuristics} are heuristic approaches that exploit exact approaches (and their complementary strengths) without guaranteeing to find the optimal solutions. While matheuristics are designed with different strategies, we summarise the most widely used three strategies.

Those matheuristics which are originally exact approaches yet are implemented heuristically are overlapping with what is described under \textit{restriction} and \textit{decomposition/partitioning} methods in this subsection. Apart from those overlapping works, in the context of dynamic programming, the corridor method constructs neighbourhoods as corridors around the state trajectory of the incumbent solution \citep{sniedovich2006corridor_COIT}. Defined (preferably large) neighbourhoods can be searched with exact approaches. Dynasearch algorithm uses dynamic programming to search an exponential size neighbourhood stemming from compound moves in polynomial time \citep{congram2002iterated_COIT}.

Another group of matheuristics benefits from multiple exact models collectively within a heuristic mechanism. \cite{tarhan2022matheuristic_COIT} decompose the scheduling planning horizon into a set of buckets, solve a time-indexed model to generate a restricted model for each bucket and solve the restricted models sequentially to construct a complete feasible solution. \cite{della2014matheuristic_COIT} solve a restricted time-indexed model and a model with positional variables iteratively to search the neighbourhood of the incumbent solution. \cite{solyali2022effective_COIT} propose a matheuristic algorithm by sequentially solving different mixed integer linear programs.

Third strategy is to incorporate exact models into different components of the heuristics. This approach may have several variations. First version includes those matheuristics having a constant interaction between heuristics and mathematical programming models. \cite{manerba2014effective_COIT} use the Variable Neighbourhood Search to decide which variables to fix in their fix-and-optimise algorithm. \cite{Archetti2015_COIT} use different integer programming models in both the intensification and the diversification phases of their Tabu Search algorithm to improve the objective function value and/or restore feasibility. \cite{Adouani2022_COIT} apply exact and heuristic approaches respectively to change the value of so-called upper and lower level variables in the neighbourhood search. Other variations include matheuristics that sequentially call heuristics and the models; e.g., exact approaches following heuristics for post-optimisation \citep{Pillac2013_COIT}, exact approaches generating the initial solutions from heuristics \citep{macrina2019energy_COIT}, exact approaches supporting heuristics at their both beginning and end to provide an initial solution and to improve the final solution, respectively \citep{archetti2017matheuristic_COIT}.

We refer the reader for a detailed overview of heuristics to the works of \cite{muller1981heuristics_COIT} and \cite{silver2004overview_COIT}, of metaheuristics to the work of \cite{blum2003metaheuristics_COIT}, and of matheuristics to the work of \cite{boschetti2022matheuristics_COIT}. The most recent book by \cite{marte2018handbook_COIT} on heuristics is another invaluable resource. Finally, the progress of metaheuristics is discussed by \cite{swan2022metaheuristics_COIT}. This work provides a critical analysis of the current state of metaheuristics by focusing on cultural and technical barriers. 

For the future studies in the area of heuristics, new techniques and powerful mechanisms could be derived from practical problems to address complex systems of today's world. Another contribution can be to explore and integrate applications of artificial intelligence to deal with large scale data. Matheuristics are especially often applied for single-objective problems and accordingly, their implementation for multi-objective optimisation is a promising future research direction. For practical purposes, such as to be used within commercial solvers, it is also worthwhile to develop generic matheuristic frameworks that can address specific classes of optimisation problems. Parallel computing (i.e., parallel solution of mathematical models) and integration with machine learning (to, for example, manage the interaction with mathematical models and heuristics) are some other invaluable research directions for matheuristics.

\subsection[Linear programming (Jean-Marie~Bourjolly)]{Linear programming\protect\footnote{This subsection was written by Jean-Marie~Bourjolly.}}
\label{sec:Linear_programming}

Linear programming (LP) offers a framework for modelling the problem of extremising a linear economic function under a set of linear inequality constraints. Solving such models can be approached algebraically as well as geometrically: finding an extreme point of a polyhedron at which a given economic function is maximised or minimised. Since its inception in 1947 by Dantzig, the simplex method has been the standard algorithm for solving linear programs. A precursor, unbeknownst then to Dantzig, was a set of ideas exposed by Fourier in 1826 and 1827, and partly rediscovered by Motzkin in 1936, hence the now famous Fourier-Motzkin elimination method (\citeauthor{Dantzig1963-sy_JMB}, \citeyear{Dantzig1963-sy_JMB}, p. 84–85; \citeauthor{S98_SMPT}, \citeyear{S98_SMPT}, p. 155–157) that solves a set of linear inequalities by sequentially eliminating variables, at the cost though of exponentially increasing the number of constraints.

But since the 1930s, several researchers had been making a headway. Working independently from one another, they had grappled with specific problems: balancing the distribution of revenue (output) with the distribution of outlays (input) in the economic activity of a whole country \citep{Leontief1936-gg_JMB}; general economic equilibrium \citep{Neumann1945-wq_JMB}; production planning \citep{Kantorovich1960-ky_JMB}; transportation planning \citep{Hitchcock1941-mc_JMB,Kantorovitch1958-xh_JMB,Koopmans1949-qx_JMB}; deployment planning and logistics \citep{Dantzig1991-tc}. \cite{Dantzig1982-rx_JMB} said he had been “fascinated” by Leontief’s interindustry input-output model and wanted to generalise it by considering many alternative activities. He also credits von Neumann with the duality theory of linear programming, which parallels the work the latter did with Morgenstern on the theory of games. 

A linear program can always be expressed (in standard form) as \{minimise $cx$, subject to $Ax = b, x \geq 0$\}, where $x$ is an $n$-vector of decision variables, $A$ is an $m$ by $n$ constraint matrix that somehow weighs the variables, $b$ is an $m$-vector that puts limits on the possible values of $x$, and $cx$ is an economic function, called the objective function, that measures the quality of a given solution $x$. It is customary to assume, without loss of generality, that matrix $A$ is of rank $m$ and that $m$ is smaller than $n$ \citep[see, e.g.,][]{PS82_SMPT}. Since $\text{min}(cx) = -\text{max}(-cx)$, one can minimise a “cost” as well maximise a “profit”. Linear programs come in pairs: \{minimise $cx$, subject to $Ax = b, x \geq 0$\} and \{maximise $yb$, subject to $yA \leq c$\}. The former is called the primal, the latter is the dual. The duality theorem of linear programming has been proved to be equivalent to Farkas's lemma that was published in 1902 \citep[see, e.g.,][]{Dantzig1963-sy_JMB}. This implies that finding an optimal solution to a linear program is equivalent to finding a feasible solution to a system made of the primal and dual constraints with the additional inequality $cx \leq yb$. 

With the introduction of duality, we now have three algorithms for solving LPs: the (primal) simplex method that maintains (primal) feasibility throughout and tries to achieve optimality; the dual simplex method that maintains dual feasibility and moves toward primal feasibility; and the primal-dual algorithm that starts with a feasible solution to the dual and keeps improving it by solving an associated restricted primal. The primal-dual algorithm is the favoured simplex tool for solving most network flow problems, for instance, the famous algorithm of \cite{Ford-Fulkerson:1962_IL} for maximum flow. The dual simplex method together with column generation may come in handy when the number of constraints is huge in comparison with the number of variables \citep{DDS06_ALAL}.

A set of linear inequalities defines a convex polyhedron $P$. Therefore, since the objective function is linear, there are only three possibilities: no feasible solution (if and only if $P$ is empty); exactly one optimal solution, located at some extreme point of $P$; infinitely many optimal solutions, located at the points of a face of $P$ of dimension 1 or more, including its extreme points. The simplex method moves sequentially along the edges of $P$ from one extreme point to another. Algebraically, it moves from one set of $m$ linearly independent columns, called a basis, to another. Each basis induces a basic solution defined by setting to zero all the $n - m$ variables that do not correspond to its columns. A basic solution is feasible if all its components are non-negative. The move from one basis to another goes as follows: one column is dropped and replaced by a new one. This exchange, called a pivot, follows a set of rules for choosing the column that enters the basis and the one that exits. It is such that, barring degeneracy, the objective function decreases strictly in value at each pivot. 

Degeneracy is rooted in the fact that an extreme point of $P$ may correspond to several bases. The algebraic expression of this defectiveness is a basic feasible solution with more than $n - m$ zero components. This occurs when the number of hyperplanes intersecting at an extreme point is greater than the minimum necessary to define it. (Think of the tip of a pyramid that has a square base.) Pivoting in the presence of degeneracy may cause the simplex method to cycle. Several schemes have been devised to avoid cycling by carefully choosing the entering and leaving columns. Bland’s rule, considered as both simple and elegant, has been widely adopted \citep{Bland1977-dy_JMB}. As for finding an initial solution, if the problem is feasible, this can be done by introducing artificial non-negative variables that one then tries to drive down to zero.

Evidence shows that the simplex method is very fast in practice \citep[see][]{Shamir1987-pu_JMB}, but \cite{Klee1972-oq_JMB} designed an LP for which it must visit each one of the $2^n$ or so extreme points, which proved that it is not “good” in the sense of \cite{Edmonds1965-zj_JMB}. A “good algorithm” having been defined as one for which the worst-case complexity is polynomial with respect to the dimension of any instance, an important open question became “Is LP in $\mathcal{P}$?”. Khachiyan answered by the affirmative in 1979 when he adapted to the specific case of linear programming a known approach in convex optimisation that had been contributed to by several Soviet mathematicians \citep[see][]{Gacs1981-kz_JMB,Bland1981-fy_JMB,Chvatal1983-ri_JMB}. The argument goes as follows: given an LP, start with an ellipsoid that is big enough to contain the set $S$ of feasible solutions if it is not empty. At each iteration, check whether the centre of the ellipsoid is a solution. If it is not, there is a hyperplane $H$ separating it from $S$. Cut the ellipsoid in half by the hyperplane parallel to $H$ that goes through the centre. Then determine the smallest ellipsoid that contains the half-ellipsoid where one is trying to locate $S$, and repeat. Stop either with a solution (located at a centre) or with an ellipsoid that is too small to contain $S$. This is an important theoretical result \citep[see][]{GLS81_ALAL}, but with very little practical use as far as solving actual LPs goes. 

The same cannot be said, however, of the interior point algorithm introduced by \cite{Karmarkar84_EAY}, in which the moves happen strictly inside the set of feasible solutions instead of taking place on the envelope. Indeed, Karmarkar’s algorithm is polynomial and often competitive with the simplex method. It assumes a canonical form for linear programming in which the variables are constrained to $Ax = 0, x \geq 0$ and $x \in S = {x: x_1 + x_2 + \dots + x_n = 1}$; it further assumes, without loss of generality, that the point $e/n = (1/n, 1/n, \dots, 1/n)$ is feasible and that the minimum value of the objective function $cx$ is zero. As it seeks to stay away from the envelope of the solution polyhedron, the algorithm builds a sequence of strictly feasible solutions, i.e., that have strictly positive components, and makes a repetitive use of $e/n$. The gist of the algorithm is the following: given a strictly feasible solution $x^k$, one can define a simple bijective scaling function $f$ that maps $S$ onto itself so that $x^k$ is mapped onto $e/n$, away from the envelope, and so that $f$ has the following property: if, for any variable $x$, $f(x)$ is strictly feasible in the “new” space, then so too is $x$ in the initial space. In the “new” space, the gradient of the transformed objective function is projected on the null space of the transformed matrix $A$ augmented of a row of 1s, to account for $S$. If $p$ denotes that projection, one then moves in the direction of $-p$, i.e., in the direction of the steepest descent, while feasibility is maintained. The algorithm stops at a point $y^{k+1}$ before reaching the envelope of the feasible region. That point is transformed to $x^{k+1}$ by $f^{-1}$, and this is repeated with a new scaling bijection \citep[see][]{Strang1987-zh_JMB,Goldfarb1989-bz_JMB,Fang1993-bg_JMB,Winston2003-ls_JMB}. Important links between Karmarkar’s algorithm and the ellipsoid method have been pointed out \citep{Todd1988-uj_JMB,Ye1987-kd_JMB}. 

Integer Linear Programming (ILP), i.e., linear programs in which the variables are restricted to being integer-valued, is arguably the most challenging and beautiful expression of LP. Unfortunately, whereas LP is in $\cal{P}$, ILP is not, unless $\cal{P}=\cal{NP}$ \citep{Karp:complexity_UPCT}. However, there are classes of LPs for which, if there exists a solution at all, an integer solution is guaranteed without having to make it a requirement. This is the case, e.g., of most network flow models. And there are classes of ILPs for which any extreme point of the polyhedron of integer solutions can be obtained by “shaving off” non-integer extreme points of the outer polyhedron of real-valued solutions with hyperplanes the number of which is bounded by a polynomial in the dimension of the instance \citep{Edmonds1965-xd_JMB,Edmonds1965-zj_JMB,GLS81_ALAL,Grotschel_1988-ic_JMB,CCPS98_SMPT}. Furthermore, tackling $\cal{NP}$-complete problems has benefited greatly from this approach \citep[see, e.g.,][]{DFJ:1954_IL,ABCC07_SMPT}. 

\subsection[Mixed-integer programming (Adam~N.~Letchford \& Andrea~Lodi)]{Integer programming\protect\footnote{This subsection was written by Adam~N.~Letchford and Andrea~Lodi.}}
\label{sec:Mixed_integer_programming}

\emph{Mixed-integer programming} (MIP) is an $\cal{NP}$-hard generalisation of linear programming (LP; \S\ref{sec:Linear_programming}), in which some or all of the variables are required to take whole-number values. Way back in the late 1950s, it was already realised that a wide variety of important practical problems could be modelled as MIPs \citep{Da60_ALAL,MM57_ALAL}. Of course, at the time, there were no good algorithms, or indeed computers, to enable one to solve MIPs from real-life applications. Since then, however, dramatic progress has been made in theory, algorithms and software. Indeed, it is now possible to solve many real-life MIPs to proven optimality (or at least near-optimality) on a laptop. In this subsection, we review the main developments in this area. For more details, we refer the reader to the textbooks by \cite{CBD11_ALAL} and \cite{CCZ14_ALAL}.

In 1958, \cite{Go58_ALAL} developed the first finitely-convergent exact algorithm for pure IPs (i.e., MIPs in which all variables are restricted to whole-number values). His method was based on \emph{cutting planes}, i.e.
additional linear constraints which cut off fractional LP  solutions. Shortly after, \cite{LD60_ALAL} invented the \emph{branch-and-bound} method, in which a sequence of LP relaxations is embedded within a tree
structure. A few years later, \cite{Ba65_ALAL} devised a simpler branch-and-bound algorithm, for pure 0-1 LPs, which did not rely on LPs at all.

In the 1960s and 1970s, researchers invested considerable effort into deriving ``deep" cutting planes. This led to the discovery of Gomory mixed-integer cuts \citep{Go60_ALAL}, corner polyhedra \citep{Go69_ALAL}, intersection cuts \citep{Ba71_ALAL}, Chv\'atal-Gomory cuts \citep{Ch73_ALAL}, disjunctive cuts \citep{Ba79_ALAL,Ow73_ALAL}, and cuts derived from a study of the so-called knapsack polytope \citep{Ba75_ALAL,Wo75_ALAL}. These topics are still being studied to this day  \citep[see, e.g.,][]{CCZ14_ALAL,Co08_ALAL}.

In 1980, \cite{BM80_ALAL} developed a general-purpose heuristic for 0-1 LPs, called ``pivot-and-complement". This initiated a line of work on so-called ``primal heuristics", which also continues to this day. We will mention this again below.

A major step forward occurred in 1983, with the publication of an award-winning paper by \cite{CJP83_ALAL}. Basically, they did the following before running branch-and-bound: (\textit{i}) ``pre-process" the
formulation in order to make the LP relaxation stronger, (\textit{ii}) automatically generate  knapsack cuts to further improve the relaxation, (\textit{iii}) run a simple primal heuristic in order to obtain a feasible integer solution early on, and (\textit{iv}) permanently fix some variables to 0 or 1 based on reduced-cost arguments. In this way, they were able to solve ten real-life 0-1 LPs that had previously been regarded as unsolvable. The largest of these instances had 2756 variables and 756 constraints, a phenomenal achievement at the time. The approach of \cite{CJP83_ALAL} is now called ``cut-and-branch".

Around the same time, there were several major theoretical advances, such as the proof of the ``polynomial equivalence of separation and optimisation" \citep{GLS81_ALAL} and the development of a polynomial-time algorithm for pure IPs with a fixed number of variables \citep{Le83_ALAL}. For details, we recommend \cite{Sc86_ALAL}.

Coming back to a more practical perspective, several improvements were made to the basic cut-and-branch scheme in the 1980s and 1990s. For brevity, we just mention some highlights. Several authors proposed more powerful
pre-processing procedures \citep[e.g.,][]{DEC93_ALAL,HP91_ALAL,Sa94_ALAL}. \cite{GNS98_ALAL} developed more effective algorithms for generating knapsack cuts. Researchers also began to study cutting planes for \emph{mixed} 0-1 LPs  \citep[e.g.,][]{PVW84_ALAL,VW86_ALAL}, which eventually led to effective cut-and-branch algorithms for such problems \citep[e.g.,][]{VW87_ALAL}.

The next milestone was the invention of \emph{branch-and-cut} by \cite{PR87_ALAL}. In branch-and-cut, one has the option of generating cutting planes at \emph{any} node of the branch-and-bound tree, rather than only at the root node (as in cut-and-branch). Although this is a fairly simple idea, Padberg and Rinaldi added several ingredients to turn it into a highly effective tool. For example, (\textit{i}) care is taken to ensure that cutting planes generated at one node of the tree remain valid at all other nodes, (\textit{ii}) whenever a cutting plane is generated, it is stored in a so-called ``cut pool", (\textit{iii}) when visiting a new node of the tree, one can check the cut pool to see if it contains any useful cuts, (\textit{iv}) one uses a heuristic rule to decide when to stop cutting and start branching at any given node.

Several developments in the 1990s are also worth mentioning. First, there were some interesting works on methods to construct ``hierarchies" of relaxations for 0-1 and mixed 0-1 LPs \citep[e.g.,][]{BCC93_ALAL,LS91_ALAL,SA90_ALAL}. The method in \cite{BCC93_ALAL}, called \emph{lift-and-project}, turned out to be useful when embedded within a branch-and-cut algorithm for mixed 0-1 LPs \citep{BCC96_ALAL}. Shortly after that, \cite{Ba96_ALAL} obtained good results using Gomory mixed-integer cuts instead. This last result was a big surprise: up to then, researchers had thought that Gomory cuts were of theoretical interest only.

By the end of the 1990s, researchers were routinely solving real-life MIPs with thousands of variables and hundreds of constraints to proven optimality. Of course, MIP in general is ${\cal NP}$-hard, so one could not expect to solve \emph{all} instances so quickly. Indeed, \cite{CD99_ALAL} found a family of 0-1 LPs, called ``market split" problems, which proved to be especially challenging for branch-and-cut. This led to the development of a new class of specific algorithms called \emph{basis reduction} methods, see, e.g., \cite{Aa00_ALAL}.

In the period 2000-2010, there was a flurry of impressive works concerned with primal heuristics for MIP. For brevity, we mention just a few examples. \cite{FL03_ALAL} devised a method called \emph{local branching}, which is essentially a form of neighbourhood search in which the neighbourhoods -- being of exponential size -- are searched by solving auxiliary MIPs. Shortly after, \cite{DRL05_ALAL} presented \emph{relaxation-induced neighbourhood search} or RINS, which solves a series of small MIPs to search for integer solutions that are ``close" to the solution of the LP relaxation. Both local branching and RINS are improving heuristics, i.e., the neighbourhoods are defined with respect to a reference feasible solution to be improved. Remarkably, they solve auxiliary MIPs by simply calling a MIP solver in a black-box fashion (with work limits), thus witnessing the maturity of the field. In the same year of RINS, \cite{FGL05_ALAL} introduced the \emph{feasibility pump}, which is highly effective for MIPs where even finding a feasible solution is challenging.

The development of the branch-and-cut technology has been so impressive that many of the above-mentioned developments have been incorporated in software packages. This includes major commercial packages, such as {\tt CPLEX}, {\tt Gurobi} and {\tt FICO Xpress}, and non-commercial ones that are free to academics, such as {\tt SCIP}. We remark that this continual development in algorithms and software has been greatly enhanced by the creation and constant maintenance of {\tt MIPLIB}, a library of MIP instances on which all new methods are now routinely tested 
\citep[see][]{BBI92_ALAL,Gl17_ALAL}.

We end this section by briefly mentioning three other areas of constant development. First, there has been great progress on \emph{decomposition} approaches to MIPs that have special structure, with \emph{branch-and-price} being a particularly effective method \citep[e.g.,][]{DDS06_ALAL}. Second, there is also by now a substantial literature on \emph{stochastic} MIPs \citep[e.g.,][]{KS17_ALAL}. Third, considerable effort has been made to extend the MIP algorithmic technology to cope with nonlinearities, leading to the blossoming field of \emph{mixed-integer nonlinear programming} or MINLP \citep[e.g.,][]{LL12_ALAL}. Particularly effective algorithms and software packages are now available for \emph{convex} MINLP \citep[e.g.,][]{Kr19_ALAL}, and one of its important special cases, mixed-integer second order cone programming \citep[e.g.,][]{BS13_ALAL}.

\subsection[Nonlinear programming (E.~Alper~Y{\i}ld{\i}r{\i}m)]{Nonlinear programming\protect\footnote{This subsection was written by E.~Alper~Y{\i}ld{\i}r{\i}m.}}
\label{sec:Nonlinear_programming}

Nonlinear programming is a generalisation of linear programming (\S\ref{sec:Linear_programming}), in which the objective function or the constraints can be given by general nonlinear functions. Mathematically, a nonlinear programming problem is represented 
as 
\[
\textrm{(P)} \quad \min\{f(x): x \in \mathcal{S}\}.
\]
Here, $\mathcal{S} = \{x \in \mathbb{R}^n: g_i (x) \leq 0, \quad i = 1,\ldots,m\}$ denotes the feasible region, where, $g_i: \mathbb{R}^n \to \mathbb{R},~i = 1,\ldots,m$, and $f: \mathbb{R}^n \to \mathbb{R}$ denotes the objective function.

In comparison with linear programming, nonlinear programming problems have much more expressive power. As such, nonlinear programming problems naturally arise in almost every setting, ranging from investment planning to machine learning; from engineering to medicine; and from energy to sustainability (\S\ref{sec:Logistics}; \S\ref{sec:Ecommerce}; \S\ref{sec:Power_markets_and_systems}; \S\ref{sec:Finance}; \S\ref{sec:Healthcare}; \S \ref{sec:Location}).

In this subsection, we will give a brief overview of theory and algorithms. While we will not cover the modelling aspect, we will mention some classes of optimisation problems with desirable properties, which should imply that using an optimisation model from such classes would significantly increase the likelihood of solving it. 

The difficulty of the generic optimization problem (P) is largely determined by the properties of the objective function $f: \mathbb{R}^n \to \mathbb{R}$ and of the functions $g_i,~i = 1,\ldots,m$ that define the feasible region $\mathcal{S} \subseteq \mathbb{R}^n$. Generally speaking, increasingly more restrictive assumptions on $f$ and on $g_i,~i = 1,\ldots,m$ give rise to increasingly more structured optimisation problems with stronger and more desirable properties. For instance, the special case in which each of $f$ and $g_i,~i = 1,\ldots,m$ is a linear function, referred to as \emph{linear programming} (\S\ref{sec:Linear_programming}), is arguably the most structured class of optimization problems with very appealing theoretical properties, which lay the groundwork for several effective solution methods such as the simplex method \citep[see, e.g.,][]{Dantzig1990_EAY} and interior-point methods \citep{Karmarkar84_EAY,Wright97_EAY,Ye1997_EAY}. In contrast, general nonlinear programming problems usually enjoy fewer desirable properties.

The class of \emph{convex optimisation problems} is comprised of optimisation problems in which each of $f: \mathbb{R}^n \to \mathbb{R}$ and $g_i,~i = 1,\ldots,m$ is a convex function, which implies that $\mathcal{S} \subseteq \mathbb{R}^n$ is a convex set, and includes linear programming as a special case. Any optimisation problem that does not belong to this class is a \emph{nonconvex optimisation problem}. On the other hand, (P) is called an \emph{unconstrained optimisation problem} if $\mathcal{S} = \mathbb{R}^n$, and a \emph{constrained optimisation problem} otherwise. 

A useful notion in nonlinear programming is that of \emph{local optimality}. A point $\hat x \in \mathbb{R}^n$ is said to be a \emph{local minimiser} of (P) if there exists an open ball $\mathcal{B} \subset \mathbb{R}^n$ of positive radius centred at $\hat x$ such that $\hat x$ is a minimiser of $f$ over the potentially smaller feasible region $\mathcal{B} \cap \mathcal{S}$. In contrast, $\hat x$ is a global minimiser of (P) if $\hat x$ is a minimiser of $f$ over the entire feasible region $\mathcal{S}$. Note that a global minimiser is also a local minimiser.

We next briefly give an overview of optimality conditions for each aforementioned class of optimisation problems. We start with unconstrained optimisation problems in the one-dimensional setting (i.e., $n = 1$). If $\hat x \in \mathbb{R}$ is a local minimiser of (P), then $f$ should be neither decreasing nor increasing at $\hat x$. Assuming that $f$ is a continuously differentiable function, we therefore obtain $f^\prime(\hat x) = 0$. This geometric interpretation carries over to the higher-dimensional setting (i.e., $n \geq 2$) by simply viewing a mutivariate function as a collection of one-dimensional functions along \emph{feasible directions} at each $\hat x \in \mathbb{R}^n$, i.e., directions along which one can move starting from $\hat x \in \mathbb{R}^n$ and still remain in the feasible region. In the unconstrained case, every direction $d \in \mathbb{R}^n$ is a feasible direction at every $x \in \mathbb{R}^n$. Using the result from the one-dimensional case, if $\hat x \in \mathbb{R}^n$ is a local minimiser of (P), then the partial derivatives of $f$ with respect to each variable should be zero, or equivalently, that $\nabla f(\hat x) = 0 \in \mathbb{R}^n$, where $\nabla f: \mathbb{R}^n \to \mathbb{R}^n$ is the gradient of $f$. Such a point is called a \emph{stationary point}.

For the special case of convex unconstrained optimisation problems, the convexity of the objective function $f: \mathbb{R}^n \to \mathbb{R}$ implies that the aforementioned necessary conditions are also sufficient, i.e., a point is a local minimiser if and only if it is a stationary point. Furthermore, for convex functions, every local minimiser is, in fact, a global minimiser. Therefore, we obtain the equivalence between global minimisers and stationary points. On the other hand, for a nonconvex optimisation problem, there may be stationary points that may not correspond to a local minimiser of (P) (e.g., if $f(x) = x^3$, then $\hat x = 0$ is a stationary point but not a local minimiser). As illustrated by this example, the complete characterisation of global optimality does not carry over from convex optimisation to nonconvex optimisation, even in the unconstrained setting.

For constrained optimisation problems, we first consider the convex optimisation case. By the convexity of the feasible region $\mathcal{S} \subseteq \mathbb{R}^n$, for any $\hat x \in \mathcal{S}$, the set of all feasible directions is given by $\tilde x - \hat{x} \in \mathbb{R}^n$, where $\tilde x \in \mathcal{S}$. Arguing similarly to the unconstrained case and using the convexity of $f$, a point $\hat x \in \mathcal{S}$ is a global minimiser of (P) if and only if $f$ does not decrease along any feasible direction, i.e., if and only if $\nabla f(\hat x)^T (\tilde x - \hat x) \geq 0$ for all $\tilde x \in \mathcal{S}$. Therefore, as in the unconstrained case, we once again have the equivalence between local and global minimisers.

Next, consider a nonconvex constrained optimisation problem. If the feasible region $\mathcal{S} \subseteq \mathbb{R}^n$ is a convex set but $f$ is a nonconvex function, a similar argument as in the convex case gives rise to the following necessary condition: If $\hat x \in \mathcal{S}$ is a local minimiser of (P), then $\nabla f(\hat x)^T (\tilde x - \hat x) \geq 0$ for all $\tilde x \in \mathcal{S}$. As in the unconstrained case, simple examples show that this condition is no longer sufficient for local optimality. If, on the other hand, $\mathcal{S} \subseteq \mathbb{R}^n$ is a nonconvex set, then we instead rely on a more general notion of \emph{tangent directions} to the feasible region $\mathcal{S}$ at $\hat x$. Therefore, if $\hat x \in \mathcal{S}$ is a local minimiser, then $\nabla f(\hat x)^T d \geq 0$ for every tangent direction $d \in \mathbb{R}^n$ to the feasible region $\mathcal{S}$ at $\hat x$. In general, the set of such tangent directions may not be easy to characterise. Under certain additional assumptions about the geometry of the feasible region $\mathcal{S} \subseteq \mathbb{R}^n$, referred to as \emph{constraint qualifications} \cite[Chapter 5]{BSS2005_EAY}, explicit necessary optimality conditions can be derived. 

Having reviewed optimality conditions, we finally give a brief overview of methods for solving optimisation problems. Nonlinear optimisation algorithms are generally \emph{iterative} in nature, i.e., they generate a sequence of points $x_k \in \mathbb{R}^n,~k = 1,2,\ldots$ that satisfies certain properties. For instance, the sequence may either converge to a local or global minimiser of an optimisation problem, or may simply have a limit point that satisfies the necessary conditions for local optimality. As illustrated by the discussion on optimality conditions, one can establish considerably weaker properties for nonconvex optimisation problems in comparison with convex optimisation problems. In fact, most classes of nonconvex optimisation problems are provably difficult in a formal complexity sense, even when restricted to minimising a quadratic function over a polyhedron \citep{MK87_EAY,PV91_EAY}. As such, it would not be reasonable to expect an algorithm to solve every optimisation problem to global optimality in a reasonable amount of time. 

Therefore, different performance metrics are employed for assessing algorithms for different classes of optimisation problems. While, for convex optimisation problems, one usually expects a ``good'' algorithm to compute a global optimal solution, an algorithm for a nonconvex optimisation problem could be deemed ``effective'' if it always converges to a local (rather than a global) optimal solution. 

In the unconstrained case, given an iterate $x_k \in \mathbb{R}^n$, the main idea is to identify a feasible direction $d \in \mathbb{R}^n$ along which the objective function will decrease. Such a direction $d \in \mathbb{R}^n$, called a \emph{descent direction}, would necessarily satisfy $\nabla f(x_k)^T d < 0$. Then, a step size in this direction is determined according to certain criteria that would guarantee a decrease in the objective function. Therefore, this family if algorithms is referred to as \emph{gradient descent methods} and includes \emph{steepest descent} as a special case (i.e., the case where $d = - \nabla f(x_k)$). Under mild assumptions, this class of algorithms converges to a stationary point of $f$. Recall that such a point is a global minimiser if $f$ is a convex function. Other methods in this class are Newton methods, conjugate gradient methods, and quasi-Newton methods, each of which generates iterates that converge to a stationary point under appropriate assumptions. 

Considering the constrained case, while general convex optimisation problems do not retain all desirable properties of the simpler class of linear programming problems, they still have a sufficiently rich structure that pave the way for provably efficient solution algorithms. In fact, every convex optimisation problem, in theory, can be solved to global optimality by the ellipsoid method \citep{YN76_EAY,Shor77_EAY} or by interior-point methods in polynomial time \citep{NN94_EAY}. Furthermore, a variety of highly effective commercial and non-commercial solvers are available for solving several classes of convex optimisation problems such as linear programming, second-order cone programming, and semidefinite programming that frequently arises
in applications (see, e.g., https://neos-server.org/neos/solvers/index.html).

For the nonconvex constrained case, one approach is based on approximating a constrained optimisation problem by a sequence of unconstrained optimisation problems by either using a penalty function, based on penalising violation of constraints (\emph{penalty methods}), or using a barrier function, based on preventing the violation of constraints by keeping the iterates strictly in the relative interior of the feasible region $\mathcal{S} \subseteq \mathbb{R}^n$ (\emph{barrier methods}). Other methods include Augmented Lagrangian methods, based on combining Lagrangian relaxation with penalty methods, and Sequential Quadratic Programming methods, based on approximating the optimisation problem by a quadratic programming problem.

Finally, various real-life applications in machine learning and data science give rise to very large-scale problems that are beyond the capability of current solvers and computing platforms. For such problems, there exist a variety of heuristic optimisation methods that can be employed to find a good solution in a reasonable amount of time (\S\ref{sec:Heuristics}). However, in contrast with exact methods, such methods usually do not provide any guarantees on the quality of the solution.

Nonlinear optimisation is a very active area of research. The reader is referred to excellent textbooks for further information \cite[e.g.,][]{FM68_EAY,Man94_EAY,BSS2005_EAY,JW06_EAY,Bert16_EAY,LY2016_EAY}.

\subsection[Queueing (Hayriye~Ayhan \& Tuğçe~Işık)]{Queueing\protect\footnote{This subsection was written by Hayriye~Ayhan and Tuğçe~Işık.}}
\label{sec:Queueing}

Queueing systems arise in many real life applications including production, service systems, finance, logistics and transportation. As mentioned in \cite{stidham_HATI}, many queueing models have been studied even before the introduction of Operational Research in the 1950s. We provide a brief overview of methodologies used in queueing systems analysis. We start with exact methods and then continue with approximations and asymptotic analysis. 

The classical analysis of queueing systems involved modelling the single stage Markovian queues as birth-death processes and computing their steady state performances using Markov chain theory. These earlier theoretical contributions were initially summarised in Feller's two volume books \cite{feller1_HATI,feller2_HATI}, then in classical textbooks such as \cite{cooper_HATI}, \cite{grossharris_HATI} and  Kleinrock's two volumes \cite{kleinrock_HATI,kleinrock_HATI2}, and more recently in \cite{gautam_HATI}, \cite{mor_HATI}, and many other books. \cite{takas_HATI} focused on using transforms and generating functions for steady state and transient behaviour of queueing systems. In the early days, transforms and generating functions were considered to be exact expressions but one had to invert these generating functions in order to obtain the actual performance measures which is in general difficult. Marcel Neuts was the first one who approached this inversion problem algorithmically. In his 1981 book, \cite{neuts_HATI} focused on queues that generalise the G/M/1 structure, whereas in his second book, \cite{neuts2_HATI} generalised the structure of the M/G/1 queue. The main idea in these books is to approximate the non-exponential distributions with a phase type distribution (convolution and mixture of exponentials) which yields a continuous time Markov chain model for the original system that could be analysed, at least numerically. This line of research resulted in many contributions on the so-called {\em matrix-geometric} methods \citep[see also][]{latouche_HATI}. Arguably the most well known  result in queueing theory is Little's law ($L=\lambda W$ or its generalisation $H=\lambda G$) which provides a relationship between the mean steady state number of customers and the mean sojourn time in a system. For a thorough survey of the Little's result and its extensions, the reader is referred to \cite{whitt_HATI}. There are numerous proofs of Little's law but \cite{eltaha_HATI} provide an elegant sample path proof. On the other hand, \cite{bertsimas_HATI} relate the steady-state distribution of the number in the system (or in the queue) to the steady state distribution of the time spent in the system (or in the queue) in a queueing system under FIFO (First In First Out). 

While there has been a lot of interest in stationary queues, Massey's 1981 dissertation \cite{massey_HATI} drew attention to the analysis of non-stationary queues (i.e., queues with time dependent arrival and service processes). Massey's dissertation started with the analysis of M(t)/M(t)/1 queue and then extended to other non-stationary Markovian systems. Many subsequent papers, such as \cite{massey4_HATI}, focused on queueing models with time-dependent arrival rates, especially infinite-server ``offered-load'' models which describe the load that would be on the system if there were no limit to the available resources. The main idea of these papers is to provide algorithms (approximations) to solve the Poisson equation. On the other hand, \cite{bertsimas2_HATI} derived a set of transient distributional laws that relate the number of customers in the system (queue) at time t to the system (waiting) time of a customer that arrived to the system (queue) at time $t$.

Networks of queues have been of interest to researchers since 1950s. \cite{jackson1_HATI} was the first one to observe that joint steady state distribution of the number of customers at the nodes of a network of Markovian queues with single server (at each node) is the product of individual distribution of M/M/1 queues. \cite{jackson2_HATI} generalised this result to networks of queues with multiple servers at the nodes. \cite{gordon_HATI} discovered that the stationary distribution again has a product form in closed Markovian networks but in this case a normalisation constant is required. \cite{bcmp_HATI} proved that the product form is insensitive to the service time distribution if the service discipline satisfies certain assumptions. This and other insensitivity results in networks were also considered by \cite{kelly_HATI} and \cite{dick_HATI} which also has results on other networks such as those with blocking and rerouting. \cite{daduna_HATI} focused on obtaining explicit expressions for the steady behaviour of discrete time queueing networks and gave a moderately positive answer to the question of whether there can be a product form calculus in discrete time. In recent years, a number of models involving different compatibilities between jobs and servers in queueing systems, or between agents and resources in matching systems, have been studied, and, under Markovian assumptions and appropriate stability conditions, the stationary distributions were again shown to have product forms \citep[see][and the references therein]{rhonda_HATI}. \cite{francois_HATI} modelled a class of networks using the so-called (max,+) linear systems. In their pioneering work, using (max,+) algebra techniques, \cite{baccellischmidt_HATI} derived Taylor series expansions for the mean waiting times in Poisson driven queueing networks  that belong to the class of (max,+) linear systems. Even though these expansions are sometimes referred to as light traffic approximations, in some cases all coefficients of the series expansion can be computed yielding an exact expression. These results were generalised to transient performance measures by \cite{sven_HATI} and  joint characteristics by \cite{ayhan1_HATI}.

Exact analysis of general queueing systems is often challenging, making the characterisation of performance measures difficult. Thus, asymptotic analyses are commonly carried out via various approximation methods. We next provide an overview of such methods.

Many of the earlier works on the asymptotic analysis of queueing systems focused on heavy traffic and many server approximations for single stage queues. In his pioneering work, \cite{kingman1961single_HATI,kingman1962queues_HATI,kingman1965heavy_HATI} have asymptotically characterised the waiting time distribution for the single server queue with general interarrival and service time distributions under heavy traffic conditions (i.e., when the traffic intensity $\rho \to 1$). Several others have developed heavy traffic approximations for the G/G/s queue, where a sequence of systems with fixed number of servers and traffic intensities $\{\rho_n\}$ approaching one are considered \citep[see, for example,][]{kollerstrom1974heavy_HATI}. In these approximations, the sequence of normalised (i.e., scaled) queue length processes converge to a reflected Brownian motion with negative drift  \citep[see][]{whitt2002stochastic_HATI}, and the associated sequence of scaled stationary queue-length distributions (i.e., the stationary distribution of the limiting diffusion process) converges to an exponential distribution. We refer the reader to \cite{harrison1985brownian_HATI} for a detailed technical treatment of heavy traffic limits and diffusion approximations. Asymptotic analysis was also considered for  multi-class and multi-stage queueing networks. Defining the stability region of these networks using fluid limit analysis was considered in \cite{chen1995fluid_HATI, dai1995positive_HATI,dai1995stability_HATI}. Many of the works that considered heavy traffic analysis of multi-class queueing networks focus on achieving the so-called state space collapse. \cite{bramson1998state_HATI} demonstrated the state space collapse for first-in first-out queueing networks of Kelly type and head-of-the-line proportional processor sharing queueing networks. His framework has been used to prove state space collapse results in several other works including \cite{stolyar2004maxweight_HATI} and \cite{mandelbaum2004scheduling_HATI}. For a more comprehensive review of heavy traffic analysis of multi-class queueing networks, we refer the reader to \cite{chen2001fundamentals_HATI}.

Many-server approximations were also considered for asymptotic analyses of  queueing systems. In these approximations, the traffic intensity can be kept constant while letting the arrival rate and the number of servers go to infinity. \cite{iglehart1965limiting_HATI} showed that the resulting sequence of normalised queue length processes converges to an Ornstein-Uhlenbeck process in the many server setting when the service time distributions are exponential. Later on, \cite{whitt1982heavy_HATI} generalised this result for systems with non-exponential service times.  For a more comprehensive overview of results in this area, see \cite{whitt2002stochastic_HATI}. In their seminal work, \cite{halfin1981heavy_HATI} defined the so called Halfin-Whitt regime for the GI/M/s queue where the traffic intensities converge to one from below, the number of servers and arrival rates tend to infinity, but steady-state probability that all servers are busy remains fixed. They showed that under the appropriate scaling, the queue length processes converge to a diffusion process. In the past decades, many other asymptotic results have been obtained for many server queues in the Halfin-Whitt regime. \cite{reed2009g_HATI} studied the G/GI/s queue and obtained fluid and diffusion limit results for the queue length process.  We refer the reader to \cite{van2019economies_HATI} for a further review of the various asymptotic results obtained in the Halfin-Whitt regime. 

Although heavy traffic approximations for queues have been popular in recent decades, light traffic (as the traffic intensity $\rho\rightarrow 0$) and interpolation approximations have also been developed. \cite{bloomfield1972low_HATI} developed light traffic approximations for a single server queue. \cite{burman1983light_HATI} developed approximations for the expected delay in M/G/s queue both for heavy and light traffic, and showed that as traffic intensity goes to zero, probability of delay depends only on mean service time distributions. \cite{daley1992light_HATI} used light traffic approximations to study the limiting properties of the waiting time in many-server queues. Light traffic approximations have also been used to study the limiting processes in queueing networks  \citep[see, for example,][]{simon1992simple_HATI}. 

As mentioned earlier, approximation methods were commonly used in the asymptotic analysis of time-varying (i.e., non-stationary) queues. \cite{mandelbaum1999time_HATI} developed a fluid approximation for the queue length process in time-varying multiserver queue with abandonments and retrials. \cite{pang2010two_HATI} have developed heavy traffic approximations for infinite server queues with time-varying arrivals. The reader is referred to \cite{whitt2018time_HATI} for a recent review of the literature on non-stationary queues. 

Due to the interest in communication/telecommunication systems, in late 1990s and early 2000s, there was a lot of research on queues with heavy tailed interarrival and/or service times. Intuitively, heavy tailed distributions decay slower than an exponential distribution \citep[see][for a thorough discussion]{resnick2007heavy_HATI}. \cite{boxma2000single_HATI} provided an overview of results for single service queues with heavy tailed interarrival and/or service time distributions. \cite{sabine_HATI} and \cite{ayhan2_HATI} showed that the asymptotics of response time was dominated by the station with the heavy tailed service time in a class of open and closed networks, respectively. \cite{foss2012large_HATI} developed upper and lower bounds on the tail distribution of the stationary waiting time in the GI/GI/s queue with heavy tailed service times.

\subsection[Risk analysis (Louis~Anthony~Cox,~Jr)]{Risk analysis\protect\footnote{This subsection was written by Louis~Anthony~Cox,~Jr.}}
\label{sec:Risk_analysis}

Risk analysis is a discipline that seeks to inform people about what might happen and how to reduce the probability and severity of undesired outcomes. It draws on decision analysis, game theory, and other areas of Operational Research but is distinguished from them by the questions it asks, the frameworks it provides for answering them, and the uses to which its answers are put \citep{Aven2020-tq_TC,Greenberg2020-xl_TC}. Where decision analysis focuses on principles for identifying logically coherent choices that make preferred outcomes more likely based on a decision maker’s beliefs and value trade-offs, risk analysis seeks to inform \textit{analytic-deliberative decision-making} by multiple stakeholders –- possibly with conflicting worldviews, values, and beliefs –- for managing critically important matters ranging from the safe operation of nuclear power plants to priority-setting for public and occupational health and safety measures. 

Risk analysis is often subdivided into \textit{risk perception}, \textit{risk assessment}, \textit{risk communication}, \textit{risk management}, and \textit{risk governance and policy-making} \citep{Greenberg2021-vu_TC}. The following sections describe these components.

\subsubsection*{Risk Perception}
Public concerns and political appetite to address them are shaped by \textit{perceived} risks, whether or not they are accurate.  Several frameworks have been developed to help understand the technical, psychological, and social drivers of risk perceptions \citep{Siegrist2020-ul_TC}.  The \textit{psychometric paradigm} \citep{Slovic2000-zr_TC} explains many aspects of risk perceptions in terms of a few underlying factors such as as \textit{dread risk} (associated with a lack of control, dreaded consequences, catastrophic potential, inequity in the distribution of risks, risks increasing over time, and fatal consequences) and \textit{unknown risk} (associated with unobservability, novelty, unknown exposure, being unknown to science, and delayed consequences). The cognitive \textit{heuristics and biases} literature positions risk perceptions within a “dual process” framework in which rapid emotional evaluations (“System 1”) can be modified by slower, more effortful cognition (“System 2”) \citep{Kahneman2011-vb_TC,Skagerlund2020-cs_TC}. The \textit{cultural theory of risk} \citep{Douglas1983-td_TC,McEvoy2017-zp_TC,Bi2021-gl_TC} posits that individual perceptions of risk are shaped by social and ideological processes that emphasise or suppress perceptions of risks depending on the respondent’s values and preferred form of social order. The \textit{social amplification of risk framework} (SARF) \citep{Kasperson2022-dv_TC} describes the social amplification or attenuation of perceived risks as risk information is communicated among people with different worldviews. 

Major lessons from the study of risk perception are that experts and members of the public often view risks quite differently; that experts often focus on the probability and consequence severity dimensions of risk while members of the public consider many other aspects; that most people tend to overestimate the frequencies of rare but vivid events (e.g., terrorist attacks, murders) and underestimate the frequencies of common but familiar ones (e.g., car accidents, heart attack fatalities); and that risk perceptions of both experts and lay people are predictably shaped and distorted by cognitive heuristics and biases and are amplified or attenuated by media reports and other communications in ways that reflect the recipients’ worldviews. System 1 tends to be innumerate, responding emotionally to possibilities and categories of harm while underweighting or ignoring relevant frequencies and magnitudes. System 2 often fails to sufficiently adjust or correct the promptings of System 1 leading to decisions with predictable regrets. These findings help to explain why expert and actuarial assessments of risk often differ from lay perceptions of risk. In a democratic society, perceptions affect decisions. A major challenge for risk analysis is to assess and communicate risks to help inform and improve collective decisions in ways that understand and respect the realities of risk perception. 

\subsubsection*{Risk Assessment}
\textit{Risk assessment} addresses how large and uncertain risks are. It begins with qualitative questions about what might go wrong and proceeds to quantitative assessments of how likely adverse events are to occur and what their possible consequences and their probabilities would be \citep{Kaplan1981-vg_TC}.  Probabilistic risk assessment (PRA) and quantitative risk assessment (QRA) methods apply probability models and statistical methods to data and modelling assumptions to quantify or bound the predicted frequencies and severities of losses and to estimate how their joint probability distribution would be changed by different risk management policies or interventions. Quantitative measures of risk can be derived from the full probability distribution or stochastic process descriptions of uncertain outcomes \citep{Smidts1997-ma_TC}, including dynamic coherent risk measures used in financial risk analysis \citep{Bielecki2017-kr_TC}. Stochastic process models of the occurrence frequencies of adverse events (such as accidents at power plants or tornadoes in cities) and the probability distribution of losses for each event can also be used to estimate entire cumulative probability distributions for losses over a stated time interval for different scenarios or sets of assumptions \citep[][see also \S\ref{sec:forecasting} and \S\ref{sec:Power_markets_and_systems}]{Kaplan1981-vg_TC}.

In the past decade, PRA techniques such as causal Bayesian networks (BNs), dynamic Bayesian networks (DBNs), and related probabilistic graphical models have increasingly been used to predict the probabilistic effects caused by interventions in engineering systems \citep{Ruiz-Tagle2022-ts_TC} and public health applications \citep{Butcher2021-zo_TC}. They have largely supplanted older and less general PRA techniques such as fault tree analysis and Markov decision processes \citep{Hanea2022-ls_TC,Cox2018-ui_TC}. Together with discrete-event stochastic simulation models and continuous (systems dynamics) simulation, they provide constructive methods to predict how risk management interventions would change the probabilities of outcomes over time. This information can be used for simulation-optimisation of risk management decisions \citep{Better2011-mq_TC}.  

PRA techniques have been extended to address adversarial risks in which intelligent adaptive adversaries rather than chance events threaten the safety and values that a risk manager seeks to protect \citep{Banks2022-fi_TC}; and \textit{unknown risks} (or risks under radical uncertainty, sometimes called \textit{Knightian uncertainty} by economists), in which relevant probabilities are unknown, e.g., by using uncertainty sets that replace precise probability distributions by (usually convex) sets of possible probability distributions or by scenarios of possibilities that are not necessarily exhaustive \citep{Gilboa2017-hp_TC}.  Recent artificial intelligence and machine learning (AI/ML) methods are now being applied to natural hazards and disasters \citep{Guikema2020-pw_TC}, cybersecurity \citep{Nifakos2021-jk_TC}, power markets \citep{mar:nar:wer:zie:22_DSRW}, and financial portfolio risk management problems where new, changing, and unknown conditions make it necessary to learn effective risk prediction and management decision rules from data and experience without the guidance of well-validated PRA models \citep{Cox2020-sj_TC}.  

\subsubsection*{Risk Management, Governance, Communication, and Risk-Cost-Benefit Analysis}
Given public perceptions and technical estimates of risks, what should be done about them? Who should decide, and how? Managing risks to human health, safety, or the environment often involves “wicked” decision problems and “deep” uncertainties, meaning that there are no clear, widely agreed-to definitions of the decision problem and solutions to it \citep{Lempert2021-os_TC}. Although multiobjective and risk-sensitive or risk-constrained optimisation problems can be formulated for some risk management problems, such as routing hazardous cargo, for many wicked risk management problems, relevant decision variables, constraints, possible outcomes, and objective functions may be unknown or not widely agreed to. Risk management in such challenging cases usually involves issues of causation (what can be done and how much difference in outcome probabilities would different feasible choices make?), \textit{collective choice} (how should the disparate perceptions and preferences of individuals be resolved or aggregated for purposes of collective decision-making?) and \textit{risk governance} (who should be responsible for making, implementing, obeying, enforcing, and revising risk management decisions; how should stakeholders participate in risk management decisions; what institutions and processes should guide, restrain, and integrate collective risk management; and how should conflicts be resolved and collective decisions be made at individual, organisational, community, local, national, and international levels?) \citep{Klinke2021-zn_TC}. 

The most immediate decisions for risk managers responding to potential or actual crises are often about \textit{risk communication}. For example, if a pandemic or natural disaster such as a tsunami, volcanic eruption or hurricane, seems possible but not necessarily imminent or certain, then what should scientists and government officials tell policy-makers and the public about the uncertain risks? Who should say what to whom, and how soon?  Different risk communication goals such as informing and empowering individual and community decisions, persuading individuals to change their behaviours, instructing citizens what to do, and informing or shaping policy deliberations and decisions, require different communication approaches. Risk communication frameworks for addressing these challenges overlap with risk perception frameworks but also emphasise the roles of trust in information sources and of \textit{outrage} in mobilising public engagement and changing behaviours \citep{Malecki2020-cm_TC}.   

\textit{Risk-cost-benefit analysis} provides a simple-sounding approach to collective risk management: take risk management actions to maximise expected social utility or net societal benefits \citep[expressed as expected net present value, possibly with a risk-adjusted discount rate;][]{Eliashberg1981-dt_TC,Hammond1992-gv_TC}. However, mathematical impossibility theorems have shown that when different people have sufficiently different beliefs and preferences, there may be no coherent way to aggregate them to make collective decisions that respect normative principles such as Pareto efficiency \citep{Nehring2007-zh_TC}. Trade-offs between measures of accuracy and fairness have recently been identified for machine-learning algorithms used in risk assessment in areas such as mortgage lending and criminal justice \citep{Corbett-Davies2017-ls_TC}. Societal risk management is now often viewed less as a top-down or centralised decision and control process in which experts provide estimated probabilities and social utilities or net benefits to use in risk-cost-benefit or social utility maximisation calculations than as a participatory democratic deliberative process \citep{Rad2021-yr_TC}. Experts in risk analysis can provide useful technical information in this process but should not dominate it \citep{Pellizzoni2000-fl_TC,Greenberg2021-vu_TC,Klinke2021-zn_TC}.

Risk analysis poses intellectual, technical, and practical implementation challenges that are likely to engage and challenge Operational Research and risk analysis professionals for the foreseeable future. A more detailed review of the accomplishments, current state, and remaining challenges for much of the field of risk analysis can be found in the 40\textsuperscript{th} Anniversary Special Issue of \textit{Risk Analysis} \citep{Greenberg2020-xl_TC}, in specialised books such as \cite{Aven2015-wq_TC}, and in online resources\footnote{For example, \url{https://www.sra.org/risk-analysis-introduction/}}.

\subsection[Simulation (Christine~S.M.~Currie)]{Simulation\protect\footnote{This subsection was written by Christine~S.M.~Currie.}}
\label{sec:Simulation}
A simulation aims to reproduce the important behaviour of a real system. Our focus here is on the use of computer simulation models within operational research (OR), whilst acknowledging that the field is much wider and ranges from computer models of sub-atomic particles to simulations involving real human actors, particularly prevalent in medicine and health sciences. Three main \textit{flavours} of simulation are used in this context: discrete event simulation, agent-based modelling and system dynamics. After discussing the uses of simulation, we continue this subsection by introducing these three main flavours before going on to discuss four important areas in simulation research: conceptual modelling, input modelling and parameterisation, simulation optimisation, and finally the newer area of data-driven simulation, linking to Industry 4.0 and digital twins. A selection of open source tools for simulation are given in \S\ref{sec:Open_source_software_for_OR}.

A simulation model, built on a computer, has a number of potential functions. Principally it is used for experimentation because testing out new settings or ways of working on a simulation results in fewer negative implications than experimenting with the real system. This can allow simulation to be used for optimisation of complex stochastic systems and there has been considerable research in this area in recent years, as we discuss below. Simulation can also be used for predicting future behaviour, and the COVID-19 pandemic showcased the predictive power of simulation modelling in a very high-profile situation \citep[e.g., the agent-based model used to advise the UK government and described in][]{Ferguson2020_CC}. The process of building a simulation model results in a better understanding of the real system because of the need to identify and model the important relationships between different entities. Running the model can also help with estimating the sensitivity of model outputs to system parameters. Another use of simulation models is for training. Within the OR context, this most often takes the form of strategic game-playing \citep[e.g., the beer game developed by MIT and described in][]{Sterman1989_XW} to practice decision-making under different scenarios in a safe environment.

\textit{Discrete event simulation} (DES) is typically used to model systems in which entities move through a set of activities. Where these activities require resources, entities will queue until the resource becomes available. Such simulations are described as discrete event because the system state only varies at discrete time points, known as events. For example, a DES model might be used to describe a manufacturing line and in this case the events could include an item starting or finishing processing by a machine on the line. Usually in DES, the simulation clock will jump from one event to the next rather than moving in equal time steps. 

\textit{System dynamics} (SD; \S\ref{sec:Systems_dynamics}) was first developed in the 1950s by Jay Forrester to help with the understanding of industrial problems. SD models deal with stocks and flows, where the dynamics are dictated by a set of connected differential equations. Describing a system using an SD model can help with detecting feedback loops and delay effects and SD modelling is useful for strategic decision-making. 
 
\textit{Agent based modelling} (ABM) describes the behaviour of individual entities or \textit{agents} within a population. As \cite{Macal2016_CC} states, one of the key differences between ABM and both DES and SD is that it takes an agent perspective of a system. Each agent in the simulation will follow its own set of rules dictating its behaviour and how it interacts with the environment and other agents. Agent behaviour is typically stochastic, allowing natural variability to be included in the model. An ABM can be used as a \textit{bottom up} approach to determine \textit{emergent behaviour} where individual actions lead to a system level response.  

Regardless of the simulation modelling technique, the first part of any simulation project is to gain an understanding of the system being modelled, the objectives of the work, and the key components that should be included, referred to as \textit{conceptual modelling}. There is some discussion of the exact definition of conceptual modelling in \cite{Robinson2008_CC} but some key points are made to support the process. First, the conceptual model can be thought of as separate from the final computer model that is built and serves as an abstraction of the real system that describes what is going to be modelled and gives an indication of how that might be done. Second, the development of the conceptual model requires input from the modeller and the system owners. Third, a conceptual model does not remain constant through a simulation project but is revisited and adapted as the project continues. Recent research in conceptual modelling is reviewed by \cite{Robinson2020_CC} and has focused on designing modelling frameworks. The work described tends to be related to DES models but the core principles can also apply to the building of both ABS and SD models.

Like any other model, the utility of a simulation is very much dependent on its inputs: the garbage-in-garbage-out principle holds true here. Setting up the probability distributions used for inputs of a stochastic simulation model or parameterisation of a deterministic simulation model is referred to as \textit{input modelling} and is typically achieved through fitting statistical models to available data and eliciting expert opinion. When estimating the inputs for a simulation model from data there is some uncertainty in their true values. With a different set of data, the estimates of the inputs would likely be different. Any uncertainties setting the model inputs will propagate through to the model outputs, resulting in \textit{input uncertainty}. This is influenced both by the accuracy of the estimates of the inputs and the sensitivity of the model output to that particular distribution or parameter. \cite{Corlu2020_CC} provide a review of the current state of the art in input uncertainty research for simulation, while \cite{SongWSC2014_CC} provide practical suggestions on how to estimate the impact of unput uncertainty on the output results.

Often simulation models are used to experiment with different system set-ups. \textit{Simulation optimisation}, sometimes referred to as \textit{optimisation via simulation} describes the use of a simulation model to find the optimal value for one or more decision variables. Typically it is used in the design of stochastic systems that are too complex to be effectively described by an analytical model. Practical examples of problems that can be solved using simulation optimisation include finding the optimal number and configuration of beds in a hospital ward; determining the appropriate number of repair staff on a production line; choosing between a selection of different configurations for a system.

The problem can be represented mathematically as 
\[
\min g(\mathbf{x}), \mbox{  } \mathbf{x} \in \mathbf{\Theta},
\]
where the function we are optimising $g(\mathbf{x})$ is generally the expected value of the output of a stochastic simulation model, $g(\mathbf{x}) = E[Y(\mathbf{x}, \xi)]$; $\mathbf{x}$ is a vector of decision variables; $\mathbf{\Theta}$ denotes the feasible region for $\mathbf{x}$; and $\xi$ indicates the randomness inherent in the model. The majority of research in simulation optimisation aims to improve the efficiency of the optimisation algorithms by reducing the number of simulation replications needed to estimate the optimal values of $\mathbf{x}$. Where a complex and slow-running simulation model is used to generate the $Y(\mathbf{x}, \xi)$ this efficiency is particularly important. \cite{HongNelsonWSC2009_CC}  classify simulation optimisation problems into three categories:
\begin{enumerate}[noitemsep]
\item \textbf{Small number of solutions:} $\mathbf{\Theta}$ contains only a small number of solutions and the decision variable $\mathbf{x}$ might define a particular system configuration. In this case the problem is one of \textit{ranking and selection}.
\item \textbf{Decision variables are continuous:} $\mathbf{\Theta}$ is a convex subset of the set of $d$-dimensional real numbers and the problem is \textit{continuous optimisation via simulation}.
\item \textbf{Decision variables $\mathbf{x}$ are discrete and ordered:} $\mathbf{\Theta}$ is a subset of the set of $d$-dimensional integers and the problem is \textit{discrete optimisation via simulation}.
\end{enumerate}
A set of algorithms exists for solving each category of problem. There has also been significant work on multi-objective optimisation via simulation; for example, see \cite{Hunter2019_CC} for a detailed description of the problem and different solution approaches.

In recent years, sensing has become more widespread and the transfer of data from physical systems to control systems is now happening in close to real time. This has allowed modellers to design simulation models that are automatically fed data from the real system allowing them to either predict the future \citep[e.g., prediction of emergency departments crowding][]{Hoot2008_CC} or to use the simulation models to optimise system parameters as part of a dynamic control process. Such models are sometimes referred to as a \textit{digital twin} or \textit{symbiotic simulation} \cite[see][for a definition]{Onggo_chapter_2019_CC}. \cite{Xu2016_CC} describe how simulation can be incorporated into the Industry 4.0 framework using the example of a semiconductor fab operation. The use of simulation in dynamic control is still in its infancy and requires fast and reliable simulation optimisation algorithms as well as mechanisms for enabling the simulation model to evolve based on new input data. Industry 5.0 is intended to complement Industry 4.0 by putting societal goals at the heart of industrial decision-making\footnote{\url{https://msu.euramet.org/current_calls/documents/EC_Industry5.0.pdf}}. This is a potential growth area for simulation optimisation, particularly multi-objective simulation optimisation \citep[e.g., see ][]{Hunter2019_CC} which enables solutions to be found that support several competing objectives.

Being a huge topic with many different facets, there is no single article that provides an overview of simulation but there are several excellent textbooks covering simulation techniques including \cite{Law2015_CC} and \cite{Banks2004_CC}. The archive of the Winter Simulation Conference\footnote{\url{www.wintersim.org}} is also an extensive resource in understanding the state-of-the-art in the field and its tutorial papers provide more basic tuition in the effective implementation of simulation. Recently, the history track at the conference has also provided an overview of the evolution of the simulation field. 

\subsection[Soft OR and problem structuring methods (Mike~Yearworth \& Leroy~White)]{Soft OR\protect\footnote{Henceforth in this subsection we just refer to PSMs to avoid the obvious dual of ‘Hard OR’ and therefore to manufacture, or at least to continue to propagate, an unhelpful distinction. We retain its use here for signposting.} and problem structuring methods\protect\footnote{This subsection was written by Mike~Yearworth and Leroy~White.}}
\label{sec:Soft_OR_and_problem_structuring_methods}

Problem Structuring Methods (PSMs) are concerned with addressing problem formulation in OR. Following definitions of \cite{Mingers2004-qf_MYLW} and \cite{Rosenhead1996-xa_MYLW} they consist of a set of rigorous but not mathematical methods based on qualitative, diagrammatic modelling. They allow for a range of distinctive stakeholder views of a problem to be expressed, explored and accommodated. They encourage active participation of stakeholders in the modelling process, through facilitated workshops and the cognitive accessibility of the modelling approach. PSMs afford negotiating a joint agenda and ownership of actions. The aim is for exploration, learning and commitment from stakeholders, rather than optimisation or prediction. PSMs thus are vital and constitute a significant developmental direction for OR. See \cite{Smith2019-pw_JL} and \cite{Franco2022-mn_AFRH} for recent reviews.

Understanding the contribution of PSMs to OR requires some knowledge of their evolution. We characterise the development of the field into three phases: (\textit{i}) origins, (\textit{ii}) growth (only noted here through the increased publication rate of PSM related articles), and (\textit{iii}) maturity, covering the diffusion of PSMs to fields outside of OR, and re-integration of problem structuring into mainstream OR. Looking at the last first, we see PSMs at an important turning point, as recent work by \cite{Dyson2021-sb_MYLW} specifically identify the centrality of problem structuring in the origins of OR and lead us towards the important question of why PSMs are not seen as an essential element of every OR engagement. 

The origins of PSMs as a set of formal methods in OR arose as a consequence of the broad critique of the process of OR in the 1970-80s; the label itself was pioneered by \cite{Woolley1981-go_MYLW}. Ackoff became a trenchant critic of the sole pursuit of objectivity and optimisation in OR describing it as an “\textit{opt-out}” \citep{Ackoff1977-kh_MYLW} and set out an agenda for reconceptualising OR practice \citep{Ackoff1979-zx_GL}. \cite{Dando1981-hv_MYLW} described the situation as a “\textit{Kuhnian crisis in management science}”. In Rosenhead’s “\textit{Rational Analysis for a Problematic World: Problem Structuring Methods for Complexity, Uncertainty and Conflict}” their prescription in OR engagements was associated with dealing with problem contexts identified variously as wicked, messy, or swampy \citep[][pp. 3-11]{Rosenhead1989-nu_MYLW}. These can be summarised as problem situations that are not well-defined, involving many interested parties with different perspectives (worldviews), where there is difficulty agreeing objectives and the meaning of success, and that require creating agreement amongst the parties involved for actions to be taken. The implication of the dichotomous framing of problem contexts – i.e., wicked/tame, swamp/high ground, hard/soft, tactical/strategic – was to set out a clear critique for the whole field of OR and to suggest that to retain its relevance in dealing with the messiness of real-world problems PSMs were required to bring some rigour – and indeed a reminder of the importance – to the process of problem formulation. Importantly, the pioneers of PSMs were concerned that traditional (i.e., ‘Hard OR’) processes for problem formulation were practitioner-free \citep{Checkland1983-st_MYLW,Rosenhead1986-ti_MYLW}. 

The main PSMs set out by \cite{Rosenhead1989-nu_MYLW} were Strategic Options Development and Analysis \citep[SODA;][]{Eden1989-rz_MYLW}, arising from cognitive mapping; Soft Systems Methodology \citep[SSM;][]{Checkland1989-fc_MYLW}, emerging from the failure of Hard Systems Thinking approaches (e.g., Systems Engineering, RAND-style Systems Analysis) when applied to messy problems; and the Strategic Choice Approach \citep[SCA;][]{Friend1989-yj_MYLW}, arising from planning. In addition, Robustness Analysis, Metagame Analysis, and Hypergame Analysis were also included. However, setting the boundary of PSMs has always been an open question \citep{mingers_soft_2011_JHLS}. In the main, the core methods (SODA, SSM, SCA) are seen as exemplary and provide sufficient coherence for OR scholars and practitioners to be provided with a clear view of a common theme. 

The methodology of PSMs has long been associated with contextual matters through Systems Thinking \citep[\S\ref{sec:Systems_thinking}][]{Checkland1983-st_MYLW,Mingers2010-ei_MYLW}, Community OR \citep{Johnson2018-sj,Jones1981-xt_MYLW,Parry1991-px_MYLW}, and large group processes \citep{Shaw2004-gw_MYLW,White2002-ck_MYLW}. Methodological individualism has been addressed through Behavioural OR \citep[\S\ref{sec:Behavioural_OR};][]{Franco2016-qk_MYLW,White2016-jg_MYLW}. There has also been long-standing relevance to Multi-Criteria Decision Analysis \citep[MCDA;][]{Marttunen2017-yk_JL}, value-focused thinking \citep{Keeney1996-yo}, policy analysis \citep{Eden2004-ws_MYLW}, and strategy making \citep{Ackermann2011-as_MYLW,Dyson2000-gz_MYLW}. Bridging between PSMs and other techniques in OR has been developed as multimethodology \citep{Mingers1997-nd_MYLW}; for example, integration with Simulation \citep[\S\ref{sec:Simulation};][]{Kotiadis2006-gv_MYLW,Tako2015-wi_MYLW}. Some approaches to using the Viable Systems Model (VSM), System Dynamics (\S\ref{sec:Systems_dynamics}), and Decision Analysis (\S\ref{sec:Decision_analysis}) would also be considered as PSMs too \citep[pp. 267-288]{Rosenhead2001-ai_JL}, e.g., VSM \citep{Lowe2020-ai_MYLW} and System Dynamics \citep{Lane1998-ff_MYLW}.  We also see developments in Group Model Building (GMB) from the System Dynamics community making a significant contribution to PSMs \citep{Andersen2007-pg}. In their growth and mature phase, applications of SSM, SCA, and SODA have extended the reach of PSMs into, e.g., project management \citep{Franco2004-qi_MYLW} and environment, sustainability, and energy policy, e.g., SCA \citep{Fregonara2013-yd_MYLW}, SODA \citep{Hjortso2004-aa_MYLW}, SSM \citep{Pahl-Wostl2007-br_MYLW}, and the Drivers, Pressures, State, Impact and Response framework \citep[DPSIR][]{Bell2012-gn_MYLW}.

The state of the art and research agenda for PSMs has been the subject of periodic reflection e.g., reviews by \cite{Rosenhead1996-xa_MYLW} and \cite{Mingers2004-qf_MYLW}. A Special Issue of JORS in 2006 questioned where PSMs were heading \citep{Rosenhead2006-ih_MYLW} – variously argued as a “grassroots revolution” \citep{Westcombe2006-vm_MYLW}, an appeal to common principles \citep{Eden2006-as_MYLW}, and observations that “form and content have evolved through interaction between the ideas and their practical use” \citep{Checkland2006-qc_MYLW}. A more recent viewpoint debate “whither PSMs” again questioned their direction of travel \cite{Harwood2019-kp_MYLW,Lowe2019-gm_MYLW}. 

The qualitative nature of PSM methods raises questions about evaluating both effectiveness and value. \cite{Midgley2013-ec_MYLW}, \cite{White2006-wo_MYLW}, and \cite{Franco2022-mn_AFRH} have addressed the question of effectiveness, and whilst White goes some way towards defining the value of PSMs it is important to note the reservations expressed by \citeauthor{Checkland1990-ip_MYLW} (\citeyear{Checkland1990-ip_MYLW}, p. 299) – that measuring value in a unique problem context, the ‘messy’ realm of PSMs, is unlikely to be meaningful. \cite{Tully2019-dd_MYLW} examine this conundrum in depth from the perspective of a consulting business and make some practical suggestions for its resolution.

Theory provides an important basis for PSM development. The range of PSM practice reported has been explained by the constitutive rules that underpin specific methods. First articulated by \citeauthor{Checkland1981-op_MYLW} (\citeyear{Checkland1981-op_MYLW}, pp. 252-254), constitutive rules are generative of method rather than prescriptive and account for the range of practices that emerge, even when adopting a specific methodology such as SSM i.e., adaptation is always necessary to address the specifics of the application context. The idea was developed further by \citeauthor{Jackson2003-ln_MYLW} (\citeyear{Jackson2003-ln_MYLW}, pp. 307-311) and then by \cite{Yearworth2014-xf_MYLW} into a generic constitutive definition for PSMs. Another significant development has been a focus on PSMs as practice and drawing on practice theories. These theories provide a means of understanding OR practices by “zooming-in” to the detailed, fine-grained, scale and by “zooming-out”, looking at how specific practices affect the broader context \citep{Ormerod2022-qp_MYLW}. Together these theoretical strands provide sufficient basis on the one hand, to liberate PSMs from the pigeon-hole of the dichotomous framing of their origins, and on the other, to address OR practice as a whole and to see problem structuring as a normal, indeed necessary, part of every OR intervention. For instance, Actor Network Theory (ANT) provides a lens to look at the processes of problematisation (i.e., problem formulation) in OR practice \citep{White2009-di_MYLW}. \cite{Callon1981-tn_MYLW} draws specific attention to the “abundance of problematisations” facing expert practitioners – that there is no single specific way of problematising. Strands of ANT focus on the performative idiom; \citep{Ormerod2014-ax_MYLW} draws attention to the “mangle of practice” and the need for more informative case studies of OR practice. Other theoretical underpinnings are relevant to PSM developments e.g., PSMs as technology \citep{Keys1998-bo_MYLW}, Critical Realism \citep{Mingers2000-oa_MYLW}, Activity Theory \citep{White2016-sv_MYLW}, and the specific role of models as boundary objects \citep{Alberto_Franco2013-yw_MYLW} in facilitated workshops \citep{Franco2010-pf_MYLW}.

From a practitioner point of view, the recent report “Reinvigorating Soft OR for Practitioners” by \cite{Ranyard2021-kw_MYLW} to the Heads of OR and Analytics Forum (HORAF) and the inclusion of the knowledge requirement “How to select and apply, a range of problem structuring methods to understand complex problems” in the Operational Research Specialist Degree Apprenticeship specification by the \cite{Institute_for_Apprenticeships_Technical_Education_undated-av_MYLW} are a welcome development. 

In conclusion, for PSMs we see a return to the roots of OR as a discipline – encompassing both practice and academic scholarship – through the centrality of problem formulation to the process of OR \citep[][p. 13]{Churchman1957-ax_MYLW} and a reminder that the seeds of problem structuring can be seen in the work of the ‘founders’ of OR as uncovered by \cite{Dyson2021-sb_MYLW}. We have identified a number of research gaps that indicate future research directions for the development of PSMs. In the area of the impact of new digital technologies, \cite{Yearworth2019-qi_MYLW} propose greater use of online “same time/different places” problem structuring workshops in order to meet the requirements for fast meeting setup times, reducing carbon emissions, enabling the scale-up to large group participation, and supporting new post-pandemic working patterns. The need to address complex policy issues in the context of wicked problems is highlighted by \cite{Howick2017-of_MYLW} and \cite{Ferretti2019-qt_MYLW}, who argue for a re-invigorated engagement for PSM practice in policy analysis. Finally, \cite{Ormerod2014-kb_MYLW}, \cite{Ranyard2015-xx_MYLW} and \cite{Ormerod2023-qp} remind us that we need to see a renewed practitioner-led orientation for OR scholarship that grounds future development in solid empirical work.

\subsection[Stochastic models (Haitao~Li)]{Stochastic models\protect\footnote{This subsection was written by Haitao~Li.}}
\label{sec:Stochastic_models}

Many decision problems involve uncertainty, e.g., network design with disruption risk, portfolio selection with uncertain return, resource planning with unknown resource availability, crop planting with uncertain yield, inventory control with varying demand, and project scheduling with random task duration, etc. While the effect of uncertain parameters on the optimal solution and objective value can be studied through the well-known sensitivity analysis, or what-if analysis, in a deterministic optimisation approach, such post-optimality analysis does not prescribe solutions under uncertainty \textit{a priori}. This subsection provides an overview of a suite of optimisation models and methods that seek to obtain optimal or near-optimal solutions for the class of decision problems where some parameters are \textit{stochastic} with known probability distribution\footnote{Decision problems with \textit{incomplete information} about random parameter’s probability distribution go beyond the scope of this review.}.

Originated in the seminal work of \cite{Dantzig1955-um_HL}, stochastic programming is one of the earliest and most prominent approaches to deal with optimisation problems with stochastic parameters. The basic stochastic programming model has a two-stage framework, called two-stage stochastic programming with recourse \citep[2S-SPR;][]{Birge2011-ut_HL}. In the first stage, the \textit{here-and-now} decision is made. Then in the second stage, the \textit{recourse} decision is prescribed for each scenario of stochastic parameters after their realisation. The objective function minimises the total cost as the summation of the first-stage cost and the expected second-stage cost given the probability distribution of the stochastic parameters. It is often insightful to compute the value of stochastic solution \citep[VSS;][]{Birge1982-ix_HL} as the difference between the optimal objective function value of the deterministic counterpart (by substituting the stochastic parameters with their point estimates) and that of the stochastic programming model. We refer to \cite{Birge2011-ut_HL} and \cite{Shapiro2021-wq_HL} for a systematic and updated treatment on the modelling and theory of stochastic programming, and to \cite{Wallace2005-lf_HL} for a collection of applications of stochastic programming. Recent applications include disaster relief management \citep{Grass2016-rg_HL}, transit network design \citep{Zhao2017-jt_HL}, portfolio selection \citep{Masmoudi2018-ya_HL}, treatment plant placement in drinking water systems \citep{Schwetschenau_S_E2019-ws_HL}, process systems \citep{Li2021-qy_HL}, multi-product aggregate planning \citep{Gomez-Rocha2022-kp_HL}, and resource allocation for infrastructure planning \citep{Zhang2022-su_HL}, among others. 

A stochastic programming model may also include a constraint that is satisfied with a probability. This model is is known as the chance constrained programming model introduced by \cite{Charnes1958-nv_HL}. The probabilistic constraint can often be transformed into a deterministic constraint given the known probability distribution of a stochastic parameter. Detailed coverage on the chance constrained programming models and methods is available in \cite{Prekopa2013-du_HL}. Notable applications include farm management \citep{Moghaddam2011-yh_HL}, broadband wireless network design \citep{Clasen2014-ki_HL}, supply chain network design \citep{Shaw2016-sh_HL}, equity trading server allocation \citep{Sun2019-lq_HL}, and power system planning \citep{Geng2019-ph_HL}, among others.  

A stochastic programming model can be formulated as a deterministic mathematical programming model by associating its decision variables with the scenarios of stochastic parameters, an approach often referred to as deterministic equivalent formulation (DEF). Solving a stochastic programming model via its DEF can be computationally challenging as the size of DEF grows rapidly with the number of scenarios of stochastic parameters. Thus custom designed algorithms are often needed to obtain quality solutions for medium- and large-size stochastic programming models. Assuming the set of random parameters has finite support, the DEF of a 2S-SPR has a block structure with L-shape, which motivates the well-known L-shape method \citep{Van_Slyke1969-qz_HL} based on Benders decomposition \citep{Benders1962-sf_HL}. For problems with a large number of random scenarios, it can be computationally intractable for the exact decomposition method to obtain an optimal solution. One may resort to various sampling-based methods to obtain approximate solutions. The stochastic decomposition method proposed by \cite{Dantzig1991-yh_HL} and \cite{Higle1991-lz_HL} employs Monte Carlo simulation and importance sampling to compute sampling cuts instead of generating the exact cuts in the L-shape method. The other successful approach is sample average approximation \citep[SAA;][]{Kleywegt2002-mi_HL,Shapiro2003-sm_HL}, which approximates the second-stage objective function by an expected value function corresponding to a set of scenarios of the random parameters. Numerical experiments and results of the SAA method on various benchmark instances are available in \cite{Linderoth2006-ff_HL}. 

Another well-known stochastic modelling and solution approach is the integrated simulation-optimisation \citep{Fu2005-qb_HL}, especially used for solving problems involving discrete decision variables, widely encountered in applications in management science, operations and supply chain management. A typical integrated simulation-optimisation framework consists of two inter-related components: \textit{search} and \textit{sampling}. The search component deals with the solution space, often combinatorial in nature with large size, for which various metaheuristics \citep{Glover2003-cz_HL} can be applied. These include local search based metaheuristics such as simulated annealing \citep{kirkpatrick1983optimization_COIT}, tabu search \citep{Glover1997-bi_HL} and scatter search \citep{Glover1998-hd_HL}, as well as population-based metaheuristics, e.g., genetic algorithm \citep{Holland1975-qj_HL}. The sampling component evaluates a candidate solution via simulation, e.g., Monte Carlo or discrete event simulation. Thus an integrated simulation-optimisation approach can be viewed as an augmented deterministic metaheuristic that employs simulation to evaluate/estimate solutions in an uncertain environment. Recent applications include maritime logistics \citep{Zhou2021-vs_HL}, pooled ride-hailing operators \citep{Bischoff2018-fr_HL}, and staffing for service operations \citep{Solomon2022-ou_HL}. 

Many real-world applications need decisions to be made sequentially under uncertainty, e.g., production planning, inventory control, resource allocation, and project scheduling, etc. One approach to this type of applications is the multi-stage stochastic programming \citep{Birge2011-ut_HL}, which is a generalisation of the 2S-SPR. In a typical multi-stage stochastic programming framework, a decision is made in a stage, based on the observed realisation of random parameters and the decisions made in the previous stage, to minimise the total expected future cost. The random parameters are assumed to evolve according to some known stochastic process. We refer to \cite{Zhang2023-xu_HL} for an updated and comprehensive treatment on various stochastic processes.  A nested decomposition as a generalisation of the L-shape method can be applied to obtain exact solutions to a multi-stage stochastic programming model \citep{Birge1985-hx_HL}. Conceptually, it applies Benders decomposition or the L-shape method recursively to a series of nested two-stage subproblems. Although theoretically sound, it can be computationally challenging to handle reasonably large instances as the number of scenarios grows exponentially with the number of stages and random parameters. Thus one often resorts to various approximation algorithms for obtaining quality solutions efficiently, including value function approximation, constraint relaxation, scenario reduction, and Monte Carlo methods, among others \citep{Birge2011-ut_HL}. 

An alternative approach to sequential decision making under uncertainty is stochastic dynamic programming \citep{Ross1983-yc_HL} or as a Markov decision processes \citep[MDP;][]{Puterman2005-lr_HL}. See \S\ref{sec:Dynamic_programming} for more details.

There are two general approaches for solving an MDP model: open-loop and closed-loop \citep{bertsekas2012dynamic_DLLD}. An open-loop approach obtains a solution to all the decision variables upfront, which is \textit{static} in nature without updating during execution of the sequential decision-making process. The integrated simulation-optimisation approach introduced above is a successful way to obtain an open-loop policy, e.g., using genetic algorithms \citep{Ballestin2007-ma_HL}, tabu search \citep{Tsai1998-tz_HL}, or the greedy randomised adaptive search procedure \citep[GRASP;][]{Ballestin2009-rp_HL} with simulation for the Stochastic Resource-Constrained Project Scheduling Problem (SRCPSP). 

Instead of optimising the entire problem upfront, a closed-loop approach seeks to obtain an optimal decision rule (\textit{policy}) to map the state at a stage to a decision, given the information available to the decision-maker at the current stage. A closed-loop policy is \textit{dynamic} and \textit{adaptive} in nature, thus is more flexible than an open-loop policy \citep{Dreyfus1977-sl_HL,bertsekas2012dynamic_DLLD}. Although theoretically attractive, to obtain an optimal closed-loop policy through the well-known Bellman equation in recursive way \citep{Bellman1957-ue_HL} is computationally intractable due to the curse-of-dimensionalities of MDP in state space, solution space of decision variables, and sample space of random parameters. 

Recent advances advocate the design and implementation of approximate dynamic programming (ADP) for solving large-scale MDPs. ADP has its roots in neural dynamic programming \citep[NDP;][]{Bertsekas1996-qb_HL} and reinforcement learning \citep[RL;][]{sutton2018reinforcement_LCAL}. Its key idea is to replace the exact cost-to-go function with some sort of approximation. We refer to \cite{Si2004-mp_HL} and \cite{Powell2011-km_HL} for comprehensive coverage on ADP and its applications. There are two approximation paradigms for the design of an ADP algorithm. The value function approximation approach works directly on the cost-to-go function to replace it with an alternative functional form that is computationally tractable. Using sample path simulation, a forward iteration procedure can be implemented to solve a deterministic sub-problem with the approximated objective function subject to the set of constraints corresponding to the state of the current stage \citep{Powell2011-km_HL}. This approach has been successfully applied to the multicommodity network flow problem \citep{Topaloglu2006-bn_HL}, dynamic fleet management \citep{Simao2009-ag_HL}, and dynamic resource planning \citep{Solomon2019-cu_HL}. 

While the value function approximation approach works well for problems with structures amenable to efficient mathematical programming methods such as linear programming or network optimisation, many combinatorial optimisation problems are $\mathcal{NP}$-hard themselves, and can be computationally demanding for mathematical programming to handle. We refer to \S\ref{sec:Computational_complexity} for a review on the topic of computational complexity and $\mathcal{NP}$-hardness. This calls for an alternative approximation paradigm known as the rollout policy \citep{Bertsekas1997-qm_HL}. A rollout policy estimates the cost-to-go function using some heuristic via simulation, which can be either an efficient problem-specific heuristic or a custom-designed metaheuristic for the problem at hand. It can be viewed as a look-ahead policy that estimates the cost of a decision-state pair under uncertainty about the future, which can be in contrast to the lookup table approach in RL \citep{sutton2018reinforcement_LCAL} where the cost of a decision-state pair is learned through simulation in a look-back fashion. A hybrid look-ahead and look-back ADP algorithm has been developed by \cite{Li2015-ey_HL} to take advantage of the complementary strengths of the pure rollout approach and the lookup table approach alone. Successful applications of rollout policy have been reported for stochastic vehicle routing \citep{Secomandi2001-lw_HL,Goodson2013-fd_HL}, SRCPSP with stochastic activity durations \citep{Li2015-ey_HL}, RCPSP with multiple-overlapping modes \citep{Chu2019-bw_HL}, ride-hailing system planning \citep{Al-Kanj2020-rn_HL}, and attended home delivery \citep{Koch2020-bf_HL}. 

All the aforementioned models and methods assume that the probability distribution of random parameters is known or can be properly estimated. This assumption may not hold in some situations where there is lack of knowledge about the uncertain parameters, or error in measurement or implementation. Optimisation with uncertain parameters without probability distribution calls for the robust optimisation (RO) approach. Although the origin of RO can be dated back to the 1970s \citep{Soyster1973-hb_HL}, RO has been growing as an active research field since the last two decades. In an RO model, one assumes that uncertain parameters are within a user-specified uncertainty set. A robust feasible solution satisfies the set of uncertain constraints for all realisations of the uncertain parameters in the uncertainty set. One main technique to solve an RO model is the robust reformulation approach to obtain a computationally tractable robust counterpart (RC) with a finite number of deterministic constraints \citep{Bertsimas2011-fv_HL}. When choosing the type of uncertainty set for the model, one often needs to trade-off between robustness against realisations of the uncertain parameters and computational tractability, i.e., size of the uncertainty set \citep{Gorissen2015-zv_HL}. We refer to \cite{Ben-Tal2002-as_HL} and \cite{Ben-Tal2009-lm_HL} for systematic treatment on robust optimisation. RO has been applied in various fields including finance \citep{Georgantas2021-nj_HL}, energy and utility \citep{Sun2021-pm_HL}, supply chain \citep{Ben-Tal2005-bz_HL,Pishvaee2011-hu_HL}, healthcare \citep{Meng2015-dg_HL}, and marketing \citep{Wang2012-qi_HL}. 

\subsection[System dynamics (Martin~Kunc \& John~D.W.~Morecroft)]{System dynamics\protect\footnote{This subsection was written by Martin~Kunc and John~D.W.~Morecroft.}}
\label{sec:Systems_dynamics}

System Dynamics (SD), founded by \cite{Forrester1961-px_MCJM}, is a ``rigorous method for qualitative description, exploration and analysis of complex systems in terms of their processes, information, organisational boundaries and strategies; which facilitates quantitative simulation modelling and analysis for the design of system structure and control'' \citep{Wolstenholme1990-df_MCJM}. SD modelling involves (as extracted from the System Dynamics Society website - \href{www.systemdynamics.org}{www.systemdynamics.org}):
\begin{itemize}[noitemsep]
   \item ``Defining problems dynamically, in terms of graphs over time.
   \item Striving for an endogenous, behavioural view of the significant dynamics of a system, a focus inward on the characteristics of a system that themselves generate or exacerbate the perceived problem.
   \item Thinking of all concepts in the real system as continuous quantities interconnected in loops of information feedback and circular causality.
   \item Identifying independent stocks or accumulations (levels) in the system and their inflows and outflows (rates).
   \item Formulating a behavioural model capable of reproducing, by itself, the dynamic problem of concern. The model is usually a computer simulation model expressed in nonlinear equations or can be left without quantities as a diagram capturing the stock-and-flow/causal feedback structure of the system.
   \item Deriving understandings and applicable policy insights from the resulting model.
   \item Implementing changes resulting from model-based understandings and insights.''
\end{itemize}

SD can be employed for both qualitative and quantitative modelling. On the one hand, tools and methods employed for qualitative SD modelling are also considered Soft Operational Research or Problem Structuring methods. On the other hand, quantitative SD modelling shares many aspects of traditional simulation methods or Hard Operational Research. Using SD quantitatively implies the development of a 5-steps process \citep{Sterman2000_SMD} that starts with a dynamic hypothesis about a structure responsible for the performance over time observed in the system followed by the model formulation, testing and experimentation. The next section discusses both approaches in detail.

One interesting characteristic of SD models is the spectrum of model fidelity they cover \citep{Morecroft2012-zf_MCJM}. Figure 1 illustrates a spectrum of model fidelity and realism. Models range in size from large and detailed to small and metaphorical.  On the left-hand side are analogue, high-fidelity models epitomised by aircraft flight simulators used to train pilots and to rehearse crisis scenarios.  They are constructed with realistic detail and accurate scaling to provide a vivid and lifelike experience of flying the aircraft they represent. People typically expect business and social models to be similarly realistic; the more realistic the better. Realistic high-fidelity models are discussed later in this subsection. But very often smaller models are extremely useful, particularly when their purpose is to aid communication and to build shared understanding of contentious problem situations in business and society. As Figure 1 suggests, the spectrum of useful models can include illustrative models (of limited detail yet plausible scaling) or even tiny metaphorical models (of minimal detail yet transferable insight).

\begin{figure}[ht!]
\begin{center}
\includegraphics[width=10cm]{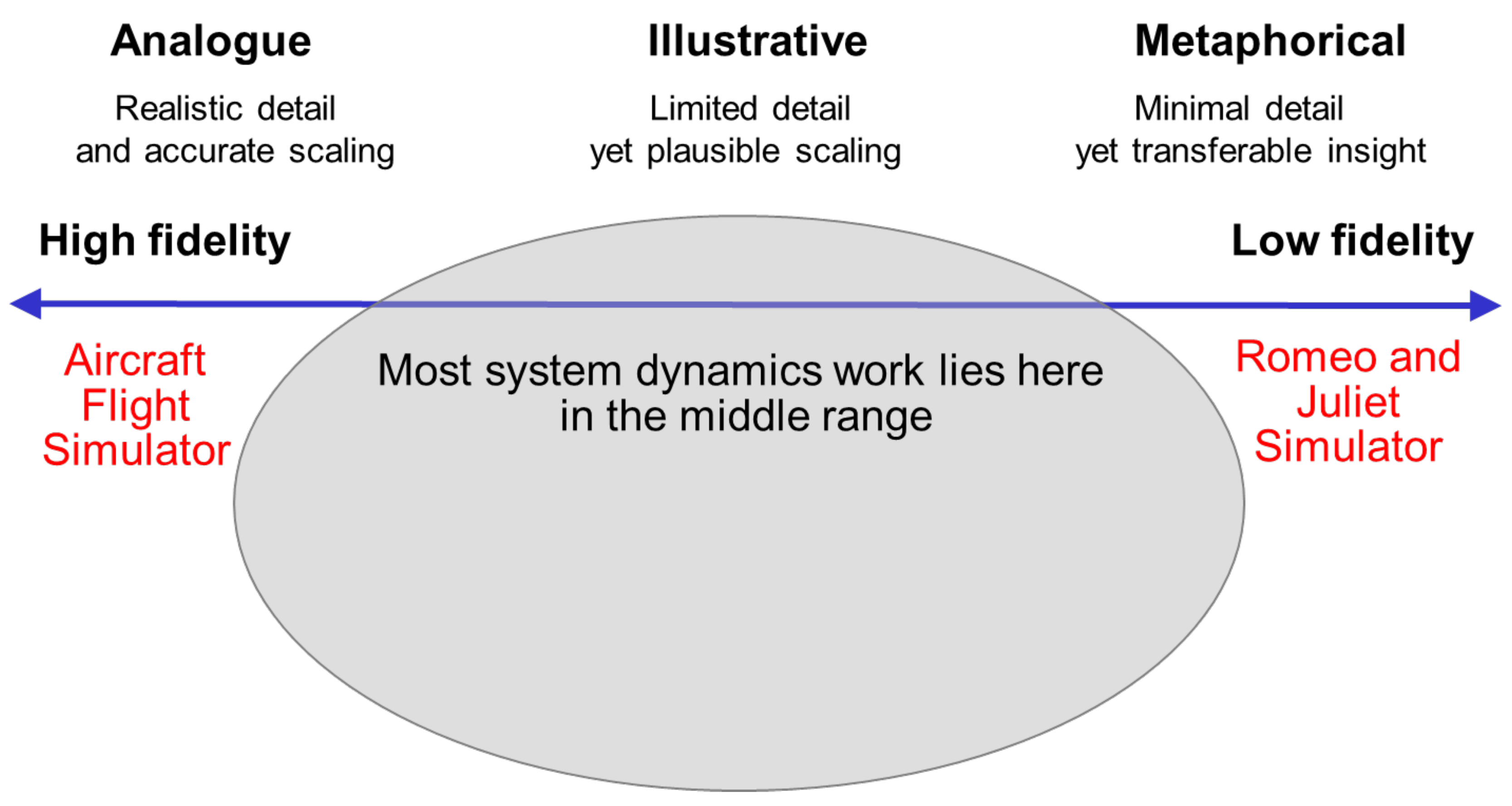} 
\caption{Modelling and Realism: A Spectrum of Model Fidelity. Adapted from \cite{Morecroft2015-vr_MCJM}, Chapter 10.}
\label{fig:MCJM_fig1}
\end{center} 
\end{figure}

At the other end of the spectrum, on the far right, is a low fidelity Romeo and Juliet simulator \citep{Morecroft2010-ui_MCJM}. This particular simulation model contains just four main concepts: Romeo's love for Juliet, Juliet’s love for Romeo and the corresponding rates of change of their love.  It is used as a metaphorical model or transitional object to help undergraduates and high school students to better understand something complex and abstract -- differential equations or even Shakespeare's play. Clearly, a simulator cannot possibly replicate Shakespeare’s play, but it can encourage students to study the play more closely than they otherwise would. It is this metaphorical property of small models -- to attract people's attention, to encourage them to reflect and debate -- that often underpins their value to model users. Sometimes metaphorical models enable client engagement. One could say that `small is beautiful' in the world of policy and strategy modelling. Over the years SD studies have included models and simulators that cover the entire range. See \cite{Kunc2017-pw_MCJM} for a sample of papers published in the Journal of the Operational Research Society and \cite{Kunc2018-uj_MCJM} for a study on published SD models. 

\subsubsection*{Qualitative System Dynamics}
The main objective of qualitative SD involves discovering the structure, in terms of feedback loops, driving the dynamic behaviour of key variables, usually with clients through facilitated workshops. The main tool employed in qualitative SD modelling is causal loop diagram (CLD). The steps for developing a CLD are \citep[based on][]{Kunc2017-bq_MCJM}:
\begin{enumerate}[noitemsep]
   \item Understanding the direction of causality between two variables. Interestingly, it is a source of important discussion among participants in facilitated qualitative SD modelling. 
   \item Defining polarities involves identifying the relationships between two variables as either positive (same sense of direction of change) or negative (opposite sense of direction of change).
   \item Identifying feedback processes responsible for the dynamic behaviour of variables. They originate from connecting variables in a circular chain of cause-and-effect. There are two types of feedback process: reinforcing and balancing.
\end{enumerate}

Finally, table \ref{tab:MVJM_t1} shows description of the modelling process

\begin{table}[h]
\centering
\caption{Qualitative  SD modelling process based on \cite{Kunc2017-bq_MCJM}.}
\label{tab:MVJM_t1}
\medskip
\begin{tabular}{ p{2cm} p{13cm} }
 \hline 
\textbf{Modelling process}	& \textbf{Qualitative SD} \\
 \hline
\textit{Objective} & Understand the feedback structure of the system. \\
\textit{Inputs} &	Text data obtained through facilitated face-to-face meetings, interview or the interpretation of causal mechanism in reports and from theories. \\
\textit{Process}	& The modelling process implies the construction of CLD to represent individual and/or group-level interpretations of causal links. Facilitation processes are critical to uncover the causal links. \\
\textit{Outputs}	& There are three main outputs: learning about the structure of the system, changes in participants perspectives, and agreement on future policies. \\
 \hline
\end{tabular}
\end{table}

\subsubsection*{Quantitative System Dynamics}
Quantitative System Dynamics characterises the system behaviour using a set of accumulation processes linked through feedback processes. The structure of the model is represented through stocks and flows diagrams. The numerical results, which are deterministic and continuous, aim to replicate past system behaviour through calibration and testing processes before the model is used to test interventions in the system. Table \ref{tab:MVJM_t2} presents a summary of the modelling process.

\begin{table}[h]
\centering
\caption{Quantitative SD modelling process based on \cite{Kunc2017-bq_MCJM}.}
\label{tab:MVJM_t2}
\medskip
\begin{tabular}{ p{2cm} p{13cm} }
 \hline 
\textbf{Modelling process}	& \textbf{Quantitative SD} \\
 \hline
\textit{Objective}	& Test a hypothesis about the structure driving the reference mode of a variable.  \\
\textit{Inputs}	& Text data obtained through facilitated face-to-face meetings, interview or the interpretation of causal mechanism in reports and from theories to determine the structure.
Numerical data for the model can come from three sources: judgement from experts or managers, numerical data sets and facilitation processes for nonlinear functions. \\
\textit{Process}	& After defining the boundary of the model, a stock and flow diagram is developed to represent the structure of the system. Equations are formulated for each variable and parameters entered. Testing of the structure and outputs are performed to confirm the model structure replicates the behaviour observed in the key variables. \\
\textit{Outputs}	& There are three outputs: time series showing performance over time; performance over time of policies or interventions in the system; and learning about the dynamic behaviour of the system.  \\
 \hline
\end{tabular}
\end{table}

\subsubsection*{Application areas}
\begin{enumerate}[noitemsep]
\item[a)] \textit{Behavioural modelling}: There are three main areas of application. Firstly, research in decision making under dynamic complexity focused on identifying and documenting systematic misperceptions of feedback in decision making processes across multiple industries and environmental conditions using SD models \citep{Gary2008-uu_MCJM,Atkinson2016-wl_MCJM}. Secondly, experimental studies explore decision making and performance using management flight simulators or microworlds based on SD models \citep{Gary2008-uu_MCJM,Sterman2014-du_MCJM}. Thirdly, individual experimental work using SD models examines how differences in mental model accuracy and decision rules lead to differences in the performance \citep{Torres2017-vm_MCJM}. Recently, scholars have advocated for a practice of behavioural system dynamics \citep{Lane2022-md_MCJM}.
\item[b)] \textit{Group model building}: There is a wide body of research on  model conceptualisation in groups; see \cite{Rouwette2010-sc_MCJM}, where the outcome is either qualitative or quantitative SD models. Researchers have assessed the effects on communication, learning, consensus and commitment in the behaviour of groups, as well as measuring the changes in mental models and understanding the impact of group model building in terms of persuasion and attitudes \citep{Rouwette2016-su_MCJM}. 
\item[c)] \textit{Multi-scale high-fidelity systems modelling}: SD modellers are embracing new approaches to improve the scale and fidelity of their models to move from aggregate conceptual models into realistic detailed models supported with specific data. There are multiple considerations to develop high-fidelity models \citep{Sterman2018-cg_MCJM}. Firstly, these models represent heterogeneous actors in the system, which involves disaggregating single stocks into multiple stocks reflecting their differences across dimensions (e.g., age). While this solution increases granularity in the model, it also implies increasing computational burden and long simulation times, which can limit the ability to perform sensitivity analysis, change structure and test interventions by stakeholders. Secondly, high-fidelity models reflect business and social processes in detail fitting the data. Therefore, models may move from the traditional representation of time, continuous, to discrete; from state continuous variables or discrete variables; and include uncertainty using stochastic variables. Consequently, SD models can employ ordinary differential equations, stochastic differential equations, discrete event simulations, agent-based models and dynamic network models. Thirdly, multiscale modelling involves integrating models working at different temporal scales (e.g., fast and slow dynamics). Fourthly, since SD models tend to employ qualitative data (e.g., decision making rules), modellers should identify and mitigate biases in sample selection and data elicitation to collect robust qualitative data. Fifthly, quantitative data should have a clear purpose in terms of the model, so specific data related to the problem the model is solving has to be collected rather than accepting only available numerical data. Sixthly, high-fidelity models should consider parameter estimation and model analysis extremely necessary to replicate historical data.
\item[d)] \textit{SD and Artificial Intelligence (AI)}: Since abundant information is available in different forms (images, text, and numbers), there is a need for technologies that not only predict data but also learn from the environment such as AI \citep{Baryannis2019-vp_MCJM}. AI can be used for cognitive thinking, learning from behaviour, recalling, and drawing inferences \citep{Min2010-tg_MCJM}. SD models use inferences of the casual structure in system to predict future trends or test interventions (e.g., new policies). SD can combine with AI to generate AI-driven simulations based on machine-learned and mathematical rules to make more accurate models \citep{Li2022_MC}. Another use is the employment of AI methods to interpret the results of simulations, especially feedback loop dominance in complex SD models. 
\end{enumerate}

\subsubsection*{Future of System Dynamics research} 
The future of SD may be driven by developments in several on different areas. Firstly, SD can be used as a problem structuring or systems thinking method (in terms of qualitative SD) so improvements in terms of facilitation will be critical. Secondly, when SD is an aggregated simple model that helps modeller and client to learn about dynamic complexity, improvements in terms of impact on behaviour from using the model \citep{Kunc2020-fn_MCJM} will be expected. Thirdly, SD can be high fidelity systems models using all the toolkit available in terms of simulation methods and AI. The next section on systems thinking (\S\ref{sec:Systems_thinking}) looks at other systems methodologies for different purposes.

\subsection[Systems thinking (Gerald~Midgley)]{Systems thinking\protect\footnote{This subsection was written by Gerald~Midgley.}}
\label{sec:Systems_thinking}

`Systems thinking' involves us viewing complex problem situations and possible human responses to them using systems theories, methodologies, methods and concepts. We will start this section by presenting a contemporary understanding of what a ‘system’ is. An explanation of how systems thinkers use this understanding to support action to address or prevent complex problems will then follow. Subsequently, we will review 70$+$ years of systems thinking to show how we got to this contemporary understanding via three ‘waves’ of methodological development. 

\subsubsection*{What is a System?} 
A system is made up of a set of \textit{interrelated parts}, with \textit{emergent properties} \citep{Emmeche1997-cj_GM}. An emergent property is a feature that cannot be traced back to any single part of the system, so can only be understood as arising from the whole (all the parts and interrelations together). Systems have \textit{boundaries}: we can say what is inside and outside the system \citep{Ulrich1994-fw_GM}, although some interactions may cross these boundaries \citep{Von_Bertalanffy1968-yb_GM}. However, systems are always seen from the \textit{perspective} of an observer/participant \citep{Churchman1979_GM,Cabrera2015-rd_GM}. Indeed, there can be multiple perspectives on the boundaries of the system, what interrelationships (within the system and with its environment) need to be considered, what emergent properties matter, and what other perspectives should be heard. 

\subsubsection*{What is Systems Thinking?} 
Based on the above understanding of systems, we can now explain \textit{systems thinking}. It is taking a systems approach to \textit{rethinking the taken-for-granted assumptions} of decision makers, OR practitioners and stakeholders on what perspectives, boundaries, interrelationships and/or emergent properties matter in a given situation, and what the implications are for action. Many systems thinkers use the adjective ‘critical-systemic’: thinking is \textit{systemic} because of the use of the above systems concepts, but it is also \textit{critical} because it involves rethinking options for understanding and action in relation to the deployment of these concepts \citep{Ulrich1994-fw_GM,Gregory2020-zf_JL}.

\subsubsection*{Three Waves of Systems Methodology} 
Since the 1950s, there have been three ‘waves’ (or successive paradigms) of systems methodology, although the second and third waves didn't fully replace their predecessors: some groups of practitioners stuck with earlier ideas. The first wave was typified by early work in systems engineering \citep[e.g.,][]{Hall1962-by_GM,Jenkins1969-vq_GM}, systems analysis \citep[e.g.,][]{Miser1985-yx_GM,Miser1988-ev_GM}, system dynamics \citep[e.g.,][see also \S\ref{sec:Systems_dynamics}]{Forrester1961-px_MCJM}, and organisational cybernetics \citep[e.g.,][]{Beer1966-aw_GM,Beer1981-fe_GM}. The first wave emphasised quantitative computer modelling by experts serving clients. These experts explained emergent properties of systems by understanding interrelatedness, and then deployed these explanations to make recommendations to clients on the possible consequences of strategic and tactical decisions.

In terms of the definition of a system presented earlier, systems were seen as real-world entities; the emphases were on interrelationships and emergent properties; boundaries were relevant because modelling had to account for all the parts and interrelationships in a system that are needed to understand given emergent properties; but multiple perspectives were often bypassed, rather than listened to, in the interests of objectivity or impartiality.

Almost all the first-wave methodologies regarded models as representations of reality, with people often being viewed as deterministic parts of systems being modelled rather than self-conscious actors who can change their purposes \citep{Ackoff1979-zx_GL}. Indeed, stakeholder purposes can differ significantly from those of the systems modeller and his/her client, and ignoring this can create conflict that undermines an OR project \citep{Checkland1985-mi_GM}. Some critics \citep[e.g.,][]{Hoos1972-cq_GM,Lee1973-us_GM} argued that massive investments in large-scale modelling were wasted because systems practitioners tried to be comprehensive (e.g., modelling all interacting problems at the city scale), yet they didn’t sufficiently account for the actual questions that decision makers wanted to address –- more modest modelling for specific purposes would have been better. Worse, the typical response to project failures was to say that the models were \textit{not comprehensive enough}, so the ideal of comprehensiveness remained unquestioned \citep{Lilienfeld1978-kd_GM}. 

These criticisms led to a second wave of systems methodologies focused on stakeholder participation, qualitative modelling and dialogue for collaborative learning. The idea of producing expert recommendations was replaced by a facilitation role for the practitioner, so multiple stakeholders could develop and integrate their ideas into proposals for change. Modelling shifted from a focus on real-world systems to understanding stakeholder perspectives, which could help people develop better mutual understanding and agree broadly-acceptable ways forward. Second-wave methodologies included soft systems methodology \citep{Checkland1981-op_MYLW}, strategic assumption surfacing and testing \citep{Mason1981-wa_GM}, interactive planning \citep{Ackoff1981-rt_GM} and interactive management \citep{Warfield1994-fa_GM}. Several earlier, first-wave methodologies were transformed in the second wave to become more participative, most notably system dynamics \citep[e.g.,][]{Vennix1996-ei_GM} and organisational cybernetics \citep[e.g.,][]{Espejo1989-fd_GM}. 

It was during the second wave that the definition of a system was expanded to recognise that all systems are understood \textit{from a perspective}. Boundaries were no longer considered the real-world edges of systems, but instead marked what people include in or exclude from their deliberations \citep{Churchman1970-hv_GM}. There was a shift away from seeing systems as real-world entities to viewing them as useful \textit{ways of thinking} to structure interpretations, either of the world or of prospective actions to change that world \citep{Checkland1981-op_MYLW}.

However, this second wave came to be critiqued by a third wave of systems thinkers. Two issues came to the fore. First, a bitterly-entrenched paradigm war between first- and second-wave systems thinkers was sparked by the emergence of the second wave \citep{Jackson1984-zt_GM}. In response, there were many third-wave proposals for \textit{methodological pluralism}: drawing creatively from both first- and second-wave methodologies, and reinterpreting methods through new frameworks or guidelines for choice. The idea was to refuse the forced choice between first- and second-wave thinking, and embrace the best of both. This gave us a more flexible and responsive practice than either of the previous two waves could deliver \citep[e.g.,][]{Jackson1991-pi_GM,Mingers1997-vb_GM}. Much of the work on methodological pluralism was developed under the banner of ‘critical systems thinking’ \citep{Flood1991-ej_GM,Flood1996-pg_GM,Jackson2019-eb_GM}.

The second issue identified in the third wave was that earlier approaches were relatively naïve with respect to power relations. The first-wave assumption that the practitioner and/or client knows best could result in the coercive imposition of ‘solutions’ and/or a lack of stakeholder buy-in, which would frustrate implementation of recommendations for change \citep[e.g.,][]{Jackson1991-pi_GM,Rosenhead2001-ai_JL}. Also, there was a second-wave, practice-limiting belief that stakeholder participation in dialogue, in and of itself, allows the better argument to prevail. This overly minimises problems of bias, coercion, groupthink, deceit, ideological framing and disempowerment \citep{Mingers1980-ya_GM,Mingers1984-yt_GM,Jackson1982-rr_GM}. 

A seminal, third-wave response to the power issue was Ulrich’s (\citeyear{Ulrich1987-hh_GM,Ulrich1994-fw_GM}) critical systems heuristics. Ulrich’s central idea is being critical of the boundary judgements made by decision makers, including the OR practitioner him/herself. Nobody can have a comprehensive viewpoint, so boundaries are inevitably set with reference to the purposes and values of decision makers. However, too often, boundary judgements are taken for granted, so decision makers (often unknowingly) foist their normative assumptions on those affected by their decisions, and the latter’s voices are not heard. Ulrich encourages those involved in and affected by an OR project to reach agreement in dialogue on the key assumptions upon which that project should be based. However, when dialogue is avoided by decision makers, those affected by their ideas have the right to make a ‘polemical’ case to embarrass the decision makers into accepting discussion. The key principle is preventing powerful stakeholders (decision makers and ‘experts’, including the OR practitioner) from simply taking their boundaries and values for granted and imposing them on others. 

Following \cite{Ulrich1994-fw_GM}, \cite{Midgley1998-rw_GM} then reviewed all the second- and third-wave work on boundaries, and proposed a broader theory and practice of \textit{boundary critique}. This encourages the practitioner to explore different possible boundaries, purposes and values in an OR project, and also to uncover conflicts \citep{Midgley2011-gm_GM} and processes of marginalisation \citep{Midgley1992-cr_GM}. \cite{Midgley1998-rw_GM} argue that boundary critique is necessary in all projects dealing with complex issues, as there are likely to be initially-hidden elements of the situation that need to be accounted for. Indeed, even deciding whether a problem situation should be viewed as complex or not requires some up-front boundary critique.

In terms of the definition of a system given earlier and its implications for systems thinking, this work deepened our understanding of boundaries: taken-for-granted boundaries can reflect the structural entrenchment of power relations in our organisations, institutions and wider society \citep{Jackson1985-ws_GM}, which can cause major socio-political and environmental issues \citep{Midgley1994-cj_GM}. Therefore, third-wave systems thinkers started talking about evolving stakeholder perspectives \textit{and} structural relationships: doing either without the other can result in systemic resistance to change \citep{Gregory2000-mo_GM}. However, the starting point for intervention (following an initial boundary critique) is usually stakeholder perspectives because it is the stakeholders themselves who can then turn their attention to structural reform \citep{Boyd2007-am_GM}. Here we see the co-existence of both the first-wave understanding of real-world systems and the second-wave emphasis on stakeholder perspectives. Methodological pluralism makes perfect sense in this context, as some approaches are particularly useful for evolving stakeholder perspectives \citep[e.g.,][]{Checkland1981-op_MYLW}, and others support intervention in organisational and institutional structures \citep[e.g.,][]{Beer1981-fe_GM}. Both can be integrated into an OR project design \citep[e.g.,][]{Sydelko2021-qp_GM,Sydelko2023-pa_GM}.

Eventually, research on methodological pluralism and boundary critique was integrated into a new ‘systemic intervention’ approach by \cite{Midgley2000-gr_GM}. He recognised that boundary critique could support deep diagnoses of problem contexts, and these diagnoses could then inform the design of OR projects, drawing creatively upon methods from both previous waves of systems thinking \textit{and} from other traditions. This work unified the different strands of third-wave methodology.

Recently, however, there have been discussions about whether a fourth wave is forming. Current research foci include whether a universal theory of systems thinking is possible and necessary \citep{Cabrera2023-nw_GM}; how to construct a simple narrative of systems thinking to effectively communicate our work \citep{Midgley2021-hy_GM}; how arts-based methods can enhance practice \citep{Rajagopalan2021-cv_GM}; and what we can learn from neuroscience to inform methodological development \citep{Lilley2022-ud_GM}. It remains to be seen whether addressing these issues will extend the third wave or launch a fourth wave of systems thinking. 

\subsection[Visualisation (Martin~J.~Eppler)]{Visualisation\protect\footnote{This subsection was written by Martin~J.~Eppler.}}
\label{sec:Visualisation}

Visualisation, the graphic (and often interactive) display of quantitative or qualitative information, has established itself not only as a powerful working modality for many management and engineering contexts \citep{Basole2021-ma_MJE,Lindner2022-lk_MJE}, but also as a research field (and research method) in its own right \citep{Eppler2008-gh_MJE}. In this subsection we briefly review the visualisation field, its relevance for Operational Research (including its benefits and caveats), its various types and application contexts, its theoretical perspectives and approaches, as well as its likely future evolution.

Why care about the graphic representation of information, especially for Operational Research? The short answer that research has provided over the last decades to this question is that it provides numerous cognitive, emotional, and social benefits and thus improves our individual and collective ability to make use of information. These benefits include a quicker comprehension of information \citep{Kress2006-ug_MJE}, the detection of important patterns \citep{Bendoly2016-cx_MJE}, the ability to better discuss information \citep{Bendoly2016-cx_MJE,Meyer2018-ui_MJE}, or the greater recall of information \citep{Paivio1990-mw_MJE,Childers1984-iq_MJE}.

The visual representation of information is not, however, without risks or potential disadvantages \citep[see][]{Bresciani2015-ac_MJE,Basole2021-ma_MJE}. Visualisations can be misleading, manipulating, oversimplified, biased, or simply confusing or overwhelming. If, for example, the y-axis of a line chart has been cropped, a small improvement may be mistakenly interpreted as a substantial one. The correct interpretation of information may also require what is often referred to as visual literacy \citep[including data and information literacy, see][]{Locoro2021-dc_MJE} on behalf of the viewers. A diagram is sometimes worth ten thousand words \citep{Larkin1987-se_MJE}, but at times it requires that many words of explanation to properly understand it.

To avoid such risks, professionals need to choose the right visualisation format for the task at hand and use it diligently and in line with existing guidelines \citep[such as those made popular by][]{Tufte2001-bh_MJE} and our perceptual preferences \citep{Ware2020-sj_MJE}. There is research on both of these questions, i.e., on the available types of visualisation \citep{Shneiderman1996-an_MJE,Chi2000-og_MJE} and on their proper use \citep[see for example][]{Ware2020-sj_MJE}.

In terms of segmenting the different kinds of graphic representations for operations management contexts, one can, at the highest level, differentiate between quantitative and qualitative information visualisation. This distinction is based on the type of information that is represented: in the case of numbers or data this is referred to as quantitative visualisation. Typical examples of this genre of visualisation are business intelligence dashboards or simple overhead slides with bar and pie charts. Pie charts, however, are perceptually problematic, as we cannot visually distinguish pie section sizes accurately, let alone compare them in different pie charts. In the case of concepts, arguments, ideas, or issues this is often labelled as qualitative visualisation. Argument mapping \citep{Bresciani2018-xj_MJE} is one approach within this group that is already used in different management contexts. Whereas quantitative visualisation is mostly software-based, qualitative visualisation can be done on paper, walls, flipcharts, and other physical media.

There are, of course, also instances of mixed visualisations that combine quantitative information and qualitative insights in a single image \citep[see][for such combined representations]{Eppler2016-tb_MJE}. An example of such a hybrid visualisation would be a business intelligence dashboard (consisting of charts) that reveals conceptual diagrams through mouse-over comments (or vice versa). 

The aforementioned distinctions are part of one tradition of visualisation research, namely the classificatory or taxonomic approach \citep[see, for example,][]{Shneiderman1996-an_MJE}. This research stream or visualisation perspective aims at providing a systematic and comprehensive overview on all forms of information visualisation that are useful for the engineering or management sector. 

Another theoretical framing of the visualisation field comes from the literature on graphic representations as boundary objects that span professional frontiers and connect expertise across disciplines -- through the help of joint visual displays \citep{Black2012-ym_MJE}. This stream of literature emphasises the dual nature of visualisations to be simultaneously fixed and fluid, clear and open to multiple interpretations or functions (for example the blueprint chart of a building or the Gantt chart for a project). 

A third influential approach to make sense of the use and impact of visualisation in workplace settings is the cognitive or collaborative dimensions approach \citep{Green2006-cm_MJE,Bresciani2018-xj_MJE}. This theoretical lens sheds light on the different qualities of graphic representation that make them more or less suited to be collaboration catalysts –- for example based on their (procedural or representational) clarity, unevenness or facilitated insight. 

A similar theoretical perspective is the affordance approach \citep{Meyer2018-ui_MJE}, that highlights the different cognitive `invitations' or incentives that visualisations can provide, such as their attention grabbing effect, their interpretive flexibility or their story telling potential. 

Another influential research stream focuses not on how images are best designed for their application contexts, but on how they are appropriated and interpreted. Many researchers with this research stream employ a semiotic approach to the study of visualisation based on the seminal work by \cite{Kress2006-ug_MJE}. This approach is also informed by (research-based) insights into our perception of visual information, but additionally enriched by the conventions and (cultural) traditions that govern our interpretation processes of graphic symbols.

There are of course many other research streams discussing the design and use of visualisation in management or engineering contexts. While some of them focus on particular visualisation formats, such as diagrams, maps, 3D models, or sketches, others focus on certain application contexts, such as (big data) analytics, creativity and innovation, production simulation, or planning. This brings us to the actual application contexts of visualisation.

When are the visualisation formats and perspectives described above mostly used? Typical application contexts for information visualisation are strategising \citep{Eppler2009-qw_MJE} and planning sessions of managers and experts (for example with the help of Gantt charts or technology roadmaps \citep{Blackwell2008-gs_MJE}, risk analysis \citep{Eppler2009-qy_MJE}, ideation and problem solving workshops (using mind mapping, argument mapping, or simple sketching) in product development and business model innovation contexts (using canvases and other visual framerworks), training and development (including knowledge transfer and retention), as well as for performance management, simulations and forecasting or scenario workshops. Last but not least, visual methods are also used as research tools in their own right \citep{Comi2014-rt_MJE} to enable better access to practitioners' expectations, experiences, or priorities \citep{Bell2013-eb_MJE}. 

Many new application contexts are currently emerging within the realm of Operational Research and management, including new forms of visualisation. These novel forms include trainings and simulations in three dimensional immersive settings such as the Metaverse or augmented reality visualisations for simulations or assisted on-site decision making or operations. Another fascinating recent phenomenon consists of images created by artificial intelligence based on user instructions (such as DALL-E or similar systems). Such artificially created, at times photo-realistic images, can help (for example) in the ideation, service innovation, or product development context. The rise of artificial intelligence also impacts the interpretation of information visualisation: A case in point are data visualisation packages (such as Tableau or PowerBI) that (through AI) already assist the user in the exploration and interpretation of the provided data charts and suggest areas for deeper analysis. The visualisation field is thus a highly dynamic area with great promise, both in terms of its methodological repertoire as well as its application scope.

\clearpage

\section{Applications}
\label{sec:applications}

\subsection[Auctions and bidding (Bo~Chen)]{Auctions and bidding\protect\footnote{This subsection was written by Bo~Chen.}}
\label{sec:Auctions_and_bidding}

The 2020 Nobel Prize in Economics was awarded to Paul Milgrom and Robert Wilson for their improvements to auction theory and inventions of new auction formats. Their theoretical discoveries have improved auctions in practice and benefited sellers, buyers and taxpayers around the world \citep{RSAS20_BC}.

An \emph{auction} is usually a process of selling and/or buying goods or services that are up for bids. A \emph{bid} is a competitive offer of a price and/or quantity tag for a good or service. Auction is a particular way to determine prices and allocation of goods or services.

Auctions have been used since antiquity for the sale of a variety of objects. Today, both the range and the value of objects sold by auction have grown to staggering proportions \citep{Krishna10_BC}. The contexts within which auctions are applied include art objects, antiques, rare collectibles, expensive wines, numerous kinds of commodities, livestock, radio spectrum, used cars, real estate, online advertising, vacation packages, wholesale electricity and emission trading, and many more.

In the basic economic model, the price of a good or service is obtained when the supply and demand meet and it is normally an equilibrium value after adjustments over time. However, in some situations such adjustments cannot be made to reach an equilibrium. As \citet{Haeringer18_BC} points out, auctions are commonly used when (\textit{a}) sellers and/or buyers have little knowledge of what would be the ``right'' price (e.g., a tract of land with an unknown amount of oil underground); (\textit{b}) the supply is scarce (e.g., an art painting); (\textit{c}) the quantity or quality of the good changes very frequently (e.g., electricity or fish); and (\textit{d}) transaction frequency is low (e.g., radio spectrum).

Bidders behave strategically. Based on the available information, what they know themselves and what they believe other bidders to know, it is difficult to analyse the outcomes of different bidding rules. This is where auction theory comes in, which is closely linked to many other domains of operational research, such as game theory (\S\ref{sec:Game_theory}), behavioural OR (\S\ref{sec:Behavioural_OR}), combinatorial optimisation (\S\ref{sec:Combinatorial_optimisation}), computational complexity (\S\ref{sec:Computational_complexity}), linear programming (\S\ref{sec:Linear_programming}) and integer programming (\S\ref{sec:Mixed_integer_programming}).

\subsubsection*{Key concepts and results}
As detailed by \citet{Haeringer18_BC}, an auction consists of the following component rules: (\textit{a}) bidding format (e.g., a price, a price and a quantity, a quantity only, or a list of items if more than one item are for sale); (\textit{b}) bidding process (e.g., auction stopping criteria and information for bidders); and (\textit{c}) price and allocation (i.e., auction winner(s) and the final price(s)).

In studying auctions, it is important to bear in mind the underlying model of \emph{valuation}, the values attached to the objects by individual buyers and/or sellers. If the value, though unknown at the time of bidding, of an object is the same for all bidders, then the evaluation is of a \emph{common value}. More generally, in situations of \emph{private values}, the value of an object varies from one bidder to another. These values can be \emph{independent} or \emph{interdependent}.

If there is only one item to be sold, we have the most basic auction. Some common forms of such simple auctions are well known. In an open-outcry auction, an auctioneer takes bids from the participants and at some point of time a winner is declared, who will then pay for the item at some price related to the bids. If all bids follow the dynamics of ascending prices and the winner is the highest bidder, who pays his bidding price, then we have an \emph{English auction}. In the case of private values, the English auction is strategically equivalent to the \emph{second-price sealed-bid auction} \citep{Krishna10_BC}, in which bidders submit written bids without knowledge of other bids. The highest bidder wins but pays the price that is the second highest in the auction. On the other hand, if the auctioneer in an open-outcry auction begins with a high asking price (in the case of selling) and lowers it until some participant accepts the price (or until it reaches a predetermined reserve price), then we have a \emph{Dutch auction}. This type of open-outcry descending-price auction is most commonly used for goods that are required to be sold quickly such as flowers, or fresh produce \citep{Mishra09_BC}, as it has the advantage of speed since a sale never requires more than one bid. The Dutch auction is strategically equivalent to the \emph{first-price sealed-bid auction} \citep{Krishna10_BC}, which is the same as the second-price sealed-bid auction except that the winner pays his bidding price.

If there are multiple homogeneous (resp., heterogeneous) items to be sold, we have a \emph{multiunit} (resp.\ \emph{combinatorial}) auction.

One of the most important results in auction theory is the \emph{revenue equivalence theorem} \citep{Heydenreich09_BC,Nisan07_BC}, which in its simple form states that when bidders' valuations are private and uniformly distributed, the expected revenue of the seller is the same in the English (or second-price) and Dutch (or first-price) auctions.

In a \emph{forward auction}, a number of buyers compete to obtain goods or services from one seller  (e.g., spectrum auction). In contrast, in a \emph{reverse auction}, a number of sellers compete to obtain business from one buyer (e.g., electricity capacity market). In a \emph{double auction}, there are multiple sellers and multiple buyers (e.g., wholesale electricity market). Potential buyers submit their bid prices and potential sellers submit their ask prices to the market institution, which then chooses the price that clears the market. At this price $p$, all the sellers who asked no more than $p$ sell and all buyers who bid at least $p$ buy.

The main issues that guide auction theory involve a comparison of the performance of different auction formats. Naturally \emph{revenue} is by far the most common yardstick from the seller's perspective. However, if the auction concerns the sale of a publicly held asset to the private sector, such as the case of spectrum auction, \emph{efficiency} may be more important -- the object ends up in the hands of whoever values it most \emph{ex post}, or in the more general case of multiple items, the \emph{sum} of realised values for all participants is maximised. Besides, \emph{simplicity} and \emph{susceptibility} to collusion among bidders are among other criteria for the choice of an auction format \citep{Krishna10_BC}.

\subsubsection*{Some best practices}
One of the most important applications of auction theory is the implementation of \textit{spectrum auctions} to allocate licenses to mobile phone carriers, who act as buyers in the forward auction. One of the auction formats introduced by \citet{Milgrom87_BC} and \citet{Wilson98_BC} was first used in 1994 by the US authorities to sell radio frequencies. This practice has since spread globally and led to great benefit to society.

There can be many ways for allocating licences in general. In addition to an auction, it can proceed either with a lottery in which any interested party would just have to sign up possibly with an entry fee, or with a beauty contest in which all those wishing to obtain a licence are required to present a case and the final winners would be selected by a committee. \citet{Haeringer18_BC} provides a detailed argument why a lottery or a beauty contest is inappropriate
in the case of radio frequencies and why an auction offers a more attractive solution. There are a number of issues in selling licences of spectrum, such as those concerning collusion, demand reduction and lack of entry. Of particular relevance for common-value auctions is the so-called \emph{winner's curse} -- the winner pays too much and loses out. \citet{Haeringer18_BC} discusses how a suitable format of an auction can be used to address these issues.

A \textit{wholesale electricity market} exists when competing generators offer their electricity output to retailers. Double auctions are normally used for such a market \citep{may:tru:18_DSRW}. By its nature electricity is difficult to store and has to be available on demand. Consequently, unlike other products, it is not possible, under normal operating conditions, to keep it in stock, ration it or have customers queue for it, so the supply should match the demand very closely at any time despite the continuous variations of both \citep[][\S\ref{sec:Power_markets_and_systems}]{wer:14_DSRW}. The supply uncertainty becomes particularly relevant with an increased use of green energy (such as solar, tidal, wind energy).

The \textit{electricity capacity market} becomes necessary to build and maintain electricity capacity that may be called upon in time of need to maintain the grid balance. In the UK’s system for purchasing Short Term Operating Reserve (STOR) for electricity supply \citep{NG22_BC}, the National Grid maintains a reserve generation ability in case of sudden demand or supply variations. Part of the operating reserve is made up by contracts through auctions. In this market, the bids come as electricity capacity, so the National Grid determines the right amount of capacity to reserve from a competitive tender. Tenders are assessed on the basis of availability prices and utilisation prices together with a consideration of response times and geographical locations. In this reverse auction, the National Grid acts as the buyer, while individual electricity operators act as sellers. Extensive studies on such auctions can be found in \citet{Chao02_BC} and \citet{Schummer03_BC}. More recently, \citet{Anderson17_BC, Anderson22_BC} investigate the problem under more general settings. They show that a natural equilibrium is not only efficient but also optimal for individual bidders.

The Internet is a new exciting venue for auctions and eBay is certainly the most well-known auction place on the Internet. Auctions on eBay face new challenges due to the nature of the Internet, where an auction can take days or even weeks and potential buyers can bid whenever they want. In response, eBay uses \emph{proxy bidding} wherein a computer programme is used to bid on behalf of the bidder, who enters an auction effectively with a \emph{maximum} bid. The computer programme raises rival bids by the minimum increment set beforehand as long as it is below the maximum bid. It is easy to see that such an auction is effectively a second-price auction in which the amount entered by the bidder serves as the bidding amount. \cite{Ariely2003-eb_BC} propose an analytical framework for studying bidding behaviour in online auctions. \cite{Chothani2015-eb_BC} provide an overview of online auctions. \citet{Hickman08_BC} analyses significant differences between electronic auctions and non-electronic auctions.

\subsubsection*{Closing remarks}
There are many excellent surveys of auction theory and applications. \citet{Milgrom85_BC} and \citet{McAfee87_BC} provide a cogent account of the theory of single-object auctions and explain many extensions and applications of the theory. \citet{Milgrom04_BC} provides a comprehensive introduction to modern auction theory and its important new applications. \citet{Samuelson14_BC} examines the use of auctions, paying equal attention to theory and practice. \citet{Haeringer18_BC} and \citet{Kagel20_BC} give respectively an overview of empirical and experimental studies on auctions. \citet{Cassady67_BC} provides a colourful and insightful overview of real-world auction institutions.

\subsection[Community Operational Research (Amanda~J.~Gregory)]{Community Operational Research\protect\footnote{This subsection was written by Amanda~J.~Gregory.}}
\label{sec:Community_Operational_Research}

Community Operational Research (COR) reflects the aspirations of OR's early theorists and practitioners of “science helping society” (\citeauthor{Cook1973-fz_AJG}, \citeyear{Cook1973-fz_AJG}/\citeyear{Cook1984-cu_AJG}, p.36). There is a long tradition of COR practice that includes Ackoff's \citeyear{Ackoff1970-if_AJG} engagement with members of the Mantua ghetto in Philadelphia, Cook's projects with inner-city community organisations \citep{Cook1973-fz_AJG,Luck1984-sj_AJG}, Beer's work with the Allende Government in Chile \citep{Beer1981-fe_GM}, and numerous projects undertaken from the University of Bath \citep{Jones1981-xt_MYLW,Sims1982-sn_AJG}. See \cite{Jackson2004-qf_AJG} and \cite{Rosenhead1993-dj_AJG} for a discussion of such work. Although these early examples of COR are significant, they were far from the norm as a focus on “science helping the establishment” predominated \citep[p.36]{Cook1973-fz_AJG,Cook1984-cu_AJG}. In recognition of this, \cite{Rosenhead1986-ti_MYLW} posed the question of “who O.R. worked for (“custom”)” (p.335) in his inaugural address as President of the UK's Operational Research Society. Rosenhead answered his own question in stating that the customers were, in the main, “big business, public utilities, the military and central government departments, with a thin scatter of local governments and health and other public authorities” (p.336) to the neglect of other groups “located outside the power structure” (p.337). \cite{Rosenhead1986-ti_MYLW} not only discussed the custom of OR, but also tackled the related issue of practice in asserting that “The evolved forms of tools reflect the circumstances of their use” (p.338). Hence, mainstream OR's focus on quantification and modelling reflected its customers’ privileging of technical matters over dialogue between stakeholders and issues of emancipation (\citeauthor{Rosenhead1993-dj_AJG}, \citeyear{Rosenhead1993-dj_AJG}, drawing on \citeauthor{Habermas1972-yr_AJG}, \citeyear{Habermas1972-yr_AJG}), and involved the use of OR methods “beyond the comprehension of most people” \citep[p.339]{Rosenhead1986-ti_MYLW}, effectively masking the social and value-laden nature of much decision making. In contrast, concerns for better mutual understanding in society and freedom from oppressive power relationships inspired the call for a more transparent OR to support “a more lively, complex and elaborate social process of decision-making” \citep[p.339]{Rosenhead1986-ti_MYLW}.

Such was the impact of Rosenhead's inaugural speech and his efforts within the OR Society that engagement with non-traditional clients quickly became legitimised and formalised through the founding of a research centre, the Community Operational Research Unit, located at Northern College, which later moved to its present location at the University of Lincoln, UK. In the 1980s, the OR Society also provided support for the Centre for Community OR at the University of Hull (later to be merged into the Centre for Systems Studies), and the Community OR Network of around 300 OR practitioners. The universities of Lincoln and Hull continue to actively practice and promote COR. More recently, in 2011, the OR Society created a Pro Bono OR scheme that connects volunteer analysts with good causes\footnote{\url{https://www.theorsociety.com/get-involved/pro-bono-or/}, accessed 2023-01-19}\textsuperscript{,}\footnote{\url{https://www.theorsociety.com/get-involved/society-groups/special-interest-groups-and-networks/or-in-the-third-sector/}, accessed 2023-01-19}.

Given the multi-faceted and often complex nature of COR projects, there appears to be no one particular OR approach that has emerged as dominant. There are, though, three streams of complementary, sometimes overlapping, approaches that have proven useful in multiple reported cases of COR:

\begin{enumerate}[noitemsep]
\item Problem Structuring Methods (PSMs) are a collection of approaches that offer decision support by “way of representing the situation (that is, a model or models) that will enable participants to clarify their predicament, converge on a potentially actionable mutual problem or issue within it, and agree commitments that will at least partially resolve it” \citep[p.531]{Mingers2004-qf_MYLW}. The modelling effort may involve clarification of normative agendas through dialogue, as PSMs are largely founded on interpretivist or social constructivist epistemologies \citep{Jackson2006-cb_AJG}. For more on PSMs see \S\ref{sec:Soft_OR_and_problem_structuring_methods}.
\item Critical Systems Thinking (CST) and Critical Systems Practice (CSP) focus on the distinction of a broad range of problem contexts and the development of systems-based methods appropriate to those contexts \citep{Flood1991-fb_AJG,Flood1991-ej_GM,Mingers1997-nd_MYLW}. Having a broad range of systems methodologies to draw on is necessary but not sufficient for good practice. Consequently, \cite{Jackson2000-ei_AJG} encapsulated the notion of good practice in his statement of three commitments of CSP: critical awareness, relating to critique of the different systems methodologies, and social awareness of the societal and organisational context; improvement, referring to the achievement of something beneficial, reflecting a cautious approach to the aspiration of universal liberation; and pluralism, the need to work with multiple paradigms without recourse to some unifying metatheory. For more on systems thinking, see \S\ref{sec:Systems_thinking}.
\item Systemic Intervention (SI) developed out of CST and took as its two primary concerns critical reflection on boundaries of inclusion and exclusion \citep{Churchman1970-hv_GM,Ulrich1983-ki_AJG,Ulrich1987-hh_GM,Midgley2000-gr_GM} and methodological pluralism. \cite{Midgley2000-gr_GM} defines SI thus: “If intervention is purposeful action by an agent to create change, then systemic intervention is purposeful action by an agent to create change in relation to reflection on boundaries” (p.129). He shows how exploring boundaries informs the methodological design of a project, with the meaningful engagement of communities built in. For more on SI see section \S\ref{sec:Systems_thinking}.
\end{enumerate}

These three streams of approaches have much in common with action research \citep[AR;][]{Levin1994-yg_AJG,Midgley2000-gr_GM,Mingers2004-qf_MYLW} and, perhaps not surprisingly, AR has been a focus of a lot of COR work. Indeed, the Community Operational Research Unit explicitly articulated a working philosophy of AR following the traditions established in Latin America and Scandinavia \citep{Thunhurst1992-vg_AJG}. Over the years, a considerable and diverse body of COR work has amassed, with some contemporary and notable examples including conference papers \citep[e.g.,][]{Wong2020-so_AJG}, case-based research papers \citep[e.g.,][]{Rosenhead1996-zi_AJG,Deutsch2022-ns_AJG,Paucar-Caceres2022-dy_AJG,Pinzon-Salcedo2022-im_AJG,Chowdhury2023-th_AJG}, journal special issues \citep[e.g.,][]{Johnson2018-sj}, project reports \citep[e.g.,][]{Stephens2018-ov_AJG} and edited books \citep[e.g.,][]{Bowen1995-to_AJG,Ritchie1994-rw_AJG,Midgley2004-xx_AJG,Johnson2012-jm_AJG}.

What counts as COR is not a simple matter though, and there have been several papers over the years that have critically discussed this (see for example the different understandings reflected in \citeauthor{Midgley2018-jv_AJG}, \citeyear{Midgley2018-jv_AJG}, and \citeauthor{White2018-rh_AJG}, \citeyear{White2018-rh_AJG}). Importantly, \cite{Johnson2012-vr_AJG} suggest that some examples of COR might be more appropriately classed as capacity-building instead of “applications based on analytic models intended to provide specific policy and operational guidance to decision-makers in a way that extends existing theory and methods” (p.39). While some COR might indeed be classed as capacity building (for example, \citeauthor{Boyd2007-am_GM}, \citeyear{Boyd2007-am_GM}, are explicit that capacity building was part of their project), it is important not to confuse such interventions with those that are based on the use of models of a qualitative rather than quantitative nature. Indeed, the commitment to knowledge being embedded within the client organisation \citep{Klein2007-dj_AJG}, handing over tools and techniques \citep{Gregory1992-ab_AJG,Gregory1992-yb_AJG,Boyd2007-am_GM,Gregory2015-co_AJG} and self-organised learning \citep{Herron2018-ff_AJG} serve to bring about capability-building alongside model building and the use of analytical approaches at the local level, which does not rule out modelling and data analysis. The needs and skills of citizens and associated groups have moved on since the 1980s, such that the tools of OR (data and models) are not so incomprehensible as they might previously have been regarded, and are familiar to most if not all citizens \citep{Caulkins2008-xg_AJG}. Indeed, \citeauthor{Hindle2018-zo_JEB}'s (\citeyear{Hindle2018-zo_JEB}) work with the Trussell Trust on mapping food bank data demonstrates that charities can make good use of big data and data visualisation.

Although there are examples of COR projects being undertaken world-wide, sustained organised support has been most evident in the UK and US. \cite{Johnson2012-jm_AJG,Johnson2012-iq_AJG} brought renewed interest to the field in the US with his promotion of a stream of activity that goes by the title of Community-Based Operations Research (CBOR). \cite{Johnson2012-vr_AJG} define CBOR as “a subfield of public-sector OR... that emphasizes most strongly the needs and concerns of disadvantaged human stakeholders in well-defined neighborhoods. Within these neighborhoods, localized characteristics vary over space and exert a strong influence over relevant analytic models and policy or operational prescriptions” (p.38). Complementary to the remit of CBOR is the Institute for Operations Research and the Management Sciences (INFORMS) Pro Bono Analytics® initiative\footnote{\url{https://connect.informs.org/probonoanalytics/home}, accessed 2023-02-08}.

Whilst COR, CBOR and pro-bono OR may be said to have a related remit, it is worth mentioning a key difference, “COR takes as its remit to work with (i.e., to take as its clients) disadvantaged community groups themselves” \citep[p.610]{Rosenhead2013-zp_AJG}, whereas CBOR and pro-bono OR are more focused on making OR and analytics available to third sector and public organisations. This distinction is not undisputed \citep{Midgley2018-jv_AJG}, but the important thing is that such efforts, geared to meaningful community engagement, have not only enabled community access to OR, but have also provided a strong impetus for its theoretical and methodological development in a way that honours the legacy of OR's early founders.

As we have an increasing number of ways to connect with others and form communities, it would be easy to assume that, going forwards, COR merely needs to develop new forms of practice to support communities in these different realms. But, in a VUCA (volatile, uncertain, complex and ambiguous) world \citep{Bennis1986-ia_AJG}, we must be alert to the need to challenge simple assumptions. Rather, there is a good argument for a critical turn in COR involving the explicit examination of underpinning values and ethics \citep{Cordoba2006-tv_AJG,Jackson2006-cb_AJG}. \cite{Midgley2004-xx_AJG} have already claimed that “if practitioners do not reflect on the different visions that it is possible to promote, then there is a danger that they will default to the understanding of community that is implicit in the liberal/capitalist tradition” (p.259). This brings with it missed opportunities to pursue more challenging and empowering practices that enable political activism and give some of the most marginalised people in our society a meaningful voice in OR projects. Much of COR has arguably been quite tame, doing good in local communities without challenging the political status quo \citep{Wong1994-lb_AJG}, but in an era of climate change, biodiversity loss, rising nationalism, insecure employment, mass migration, and increasing wealth inequality, a new, more critical agenda for COR is urgently needed.

\subsection[Cutting and packing (Julia~A.~Bennell)]{Cutting and packing\protect\footnote{This subsection was written by Julia~A.~Bennell.}}
\label{sec:Cutting_and_packing}

Cutting and packing (C\&P) problems are geometric assignment problems, in which \textit{small items} are assigned to \textit{large objects} such that a given objective function is optimised and two basic geometric feasibility conditions hold, specifically containment and non-overlap. They appear in a wide range of settings, but are most commonly investigated for applications in manufacturing and transportation. For example, cutting pattern pieces from material or packing boxes into containers. These are combinatorial optimisation problems and $\mathcal{NP}$-hard. Depending on the size or geometry of the problem, there exists strong formulations that can be solved using exact methods. However, there remains many open problems that have instances that cannot be solved to optimality, or computational times are impractical for applications in practice. Moreover, there are problems where bounds are weak and only toy instances can be solved to optimality. As a result, heuristics remain an important tool in C\&P.

Given the wide variety of C\&P problems, \cite{Dyckhoff1990-lr_JB} and later \cite{Wascher2007-bg_JB} defined a typology of problems using the following dimensions:\\

\noindent Objective function
\begin{itemize}[noitemsep]
\item Output maximisation: packing the greatest value of items in a given fixed dimension finite number of large object(s).
\item Input minimisation: pack all items using the minimal number of fixed dimension objects or the minimum size large object with at least one unconstrained dimension.
\end{itemize}

\noindent Assortment of small items (items to be packed)
\begin{itemize}[noitemsep]
\item All items are identical. 
\item Weakly heterogeneous: few item types given the total number of items. 
\item Strongly heterogeneous: many item types that are unique or have few copies.
\end{itemize}

\noindent Assortment of large items
\begin{itemize}[noitemsep]
\item Single large object: fixed dimension for output maximisation, open dimension(s) for input minimisation.
\item Multiple large objects: fixed dimension, either identical or heterogeneous.
\end{itemize}

These distinctions lead to named problem types, e.g., bin packing problems (BPP) are input minimisation problems, with strongly heterogeneous small items and multiple large objects, while a knapsack problem (KP) shares the same characteristics but is an output maximisation problem. Note that problem names and their definitions are not universally accepted or consistently used, so researchers should check the articulated problem definition in the paper when selecting literature. 

The following focuses on two-dimensional (2D) and three-dimensional packing (3D) as these include the unique challenges of the geometric constraints associated with C\&P problems. One-dimensional (1D) problems remain interesting and challenging \citep[see][]{Martinovic2018-by_JB,Munien2021-sh_JB}. For an introduction to C\&P, see \cite{Scheithauer2018-xm_JB}.

\subsubsection*{Geometry}
Handling the geometric characteristic of C\&P problems adds significantly to the computational burden and the number of variables needed in a model. These increase with the spatial dimensions and with the irregularity of the shape of the small items. For 1D problems, the geometric constraints of overlap and containment are trivial. Regular shapes (rectangles, boxes, circles, spheres) add complexity through additional item location variables $x$, $y$ (and $z$) co-ordinates, and dimensions: length, width (and depth). However, the common characteristics of the shape mean these are straightforward to model. Pairwise constraints between items and between each item and the boundary of the large object ensure feasibility.

In the case of irregular shaped items, accurate non-overlap constraints cannot be reduced to comparing a set of common dimensions. While the item location is still determined by a defined origin, the arbitrary nature of the shape significantly increase the complexity of assessing geometric feasibility. At a basic level, it requires testing for edge intersections between items and containment of one item inside another. Methods to reduce the complexity include the nofit polygon, raster method and phi-functions in 2D and voxels and phi-functions in 3D. \cite{Bennell2008-nv_JB} provide a tutorial in geometric methods for 2D nesting problems. \cite{Lamas-Fernandez2022-fe_JB} describe approaches for modelling 3D geometry. Developing solution methods for irregular packing problems requires a comprehensive, fast and robust geometry library.

\subsubsection*{Constraints}
There exists a wide range of practical constraints arising from the applications. These may relate to the material being cut having defects or quality variability, the cutting tool requiring space between items or constraints on the types of cut. There may be sequencing constraints or assignment constraints that include precedence or prevent/require the packing of items together. In 2D rectangle C\&P, a common requirement is guillotine cuts where all cuts must be orthogonal and span the entire width or breadth of the rectangular material sheet. Furthermore, the number of alternate cuts (e.g., a switch from vertical to horizontal cuts) may be restricted.  Applications in 3D container loading provide a challenging set of constraints on the arrangement of boxes to ensure weight distribution, horizontal and vertical stability of load and consider the weight baring strength of the stacked boxes. \cite{Bortfeldt2013-wd_JB} describe the different types of constraints and their adoption by researchers.

\subsubsection*{Two-dimensional problems}
In two dimensions, research has focused on rectangle packing problems, and irregular shape packing problems, often called nesting problems. There is also a smaller body of literature on circle packing \citep{Hifi2009-ce_JB}.

Exact solution approaches to the 2D rectangle packing problem are reviewed by \cite{Iori2021-hs_JB} and cover the main problem types. The paper evidences the recent advances in exact methods for these problems while identifying a number of open problems that remain very challenging. Specifically, they identify the open-dimension problem, some specific instances of BPP, and problems with multiple heterogeneous large objects. Moreover, nesting problems remain a rich area of research for developing high fidelity and scalable exact methods.

Heuristics and metaheuristics are a natural choice, particularly for problems with a large number of small items. Early research focused on placement heuristics such as First Fit and Next Fit for BPP \citep{Coffman1980-vj_JB}, Bottom Left and Bottom Left Fill for open dimension problems \citep{Albano1980-jw_JB,Chazelle1983-fp_JB,Burke2004-ia_JB}. These place items into the large object in a given sequence according to the placement rule and observing any additional placement constraints. A natural evolution of this approach is to apply a metaheuristic to re-sequence the packing order to obtain better solutions, of which there are many examples. 

The 2D rectangle \textit{identical item packing} problem is known as the manufacturer's pallet loading problem. Research is mature with exact methods and heuristics that perform well across benchmark instances. \cite{Silva2016-jj_JB} comprehensively review these problems. Fewer papers have looked at the case where the small item is irregular, for example cutting metal blanks \citep[e.g.,][]{Costa2009-il_JB}.

Output maximisation problems focus on the 2D rectangle \textit{knapsack problem}, with few articles considering problems with weakly heterogeneous data. Problems cover guillotine and non-guillotine cutting and item values may be equivalent to area or have an assigned value. Furthermore, the constrained variant places upper and lower bounds on the number of each item type placed. The guillotine variant where area and value align has a fast exact method \citep{Oliveira1990-cd_JB}. Non-aligned item values and non-guillotine cutting is still challenging. \cite{CILM22B_SMPT} includes a summary of 2DKP. Packing a single large object is often a component of a larger practical problem, where the decision problem of whether a set of rectangles will fit into a fixed dimension rectangle, referred to as a 2D orthogonal packing problem, is of interest \citep{Clautiaux2007-ti_JB}. 

\textit{Cutting stock problems} have been studied for over 60 years with the seminal paper by \cite{Gilmore1965-mm_JB} that described the column generation approach for 2D rectangle cutting stock problem with guillotine cuts. Most papers focus on ILP/MILP approaches. A notable feature of these problems and how they differ from the BPP, is the way solutions are constructed arising from the data instances. Since there are few item types, but many items of each type, the solutions are composed of pattern types that are repeated across multiple stock sheets leaving a residual problem of unmet demand. 

\textit{Bin packing problems} have been extensively studied and include the guillotine and non-guillotine variant. \cite{Lodi2014-jq_JB} provide a review of BPPs. Early heuristic approaches \citep{Berkey1987-lb_JB} include two-phase algorithms that pack multiple bin width strips and then solve a 1D BPP where the item size is the height of the strip, while single-phase algorithms pack directly into the bins either using a level packing approach or a placement heuristic such as bottom-left. Increasingly, researchers are focusing on exact methods; see for example \cite{Pisinger2007-rx_JB} who use branch and price for variable size and fixed size bins. There are very few examples of bin packing with irregular shapes, where one example is glass cutting \citep{Bennell2018-fm_JB}.

The \textit{open dimension problem} variant is often called the strip packing problem. This can be formulated as a linear mixed integer programme, although practical size problems are still very challenging. \cite{Martello2003-mo_JB} develop bounds by relaxing the problem so it can be solved as a one dimensional BPP. Placement heuristics (bottom left and bottom left fill) based on sequencing of items within a (meta)heuristic framework are widely used. \cite{Hopper2001-dz_JB} undertook an extensive analysis of this type of approach. 

\textit{Nesting problems}, where the small items are irregular, are commonly formulated as open dimension problems. Solution approaches are dominated by the use of heuristics and metaheuristics. \cite{Bennell2009-kx_JB} provide a review of these methodologies including using exact models to improve local optima. This approach is also used by \cite{Stoyan2016-fk_JB} who use phi-functions, which allow orientation as a decision variable. In the last decade, researchers have developed formulations that can be solved to a global optimum for small problems. \cite{Toledo2013-xy_JB} approximates the items and the packing area to a discrete set of points allowing the problem to be solved as a MIP model. \cite{Alvarez-Valdes2013-fs_JB} used the nofit polygon to define a finite set of convex spaces and used binary variables to activate constraints. 

\subsubsection*{Three dimensional problems}
These problems are solved across the range of problem types, largely considering single container output maximisation, multi-container input minimisation and the open dimension problem. For packing boxes, the mix of constraints addressed across the literature is inconsistent and frequently not congruent with industry standards. Solutions to the problem focus on building walls, layers or blocks of identical boxes. See \cite{Zhao2016-nt_JB} for a comparative review of algorithms including exact methods. Recent papers are now looking at the additional constraints arising from a vehicle, such as axle load and stability under breaking and acceleration \cite[see][]{Ali2022-ow_JB}. 3D packing of irregular shapes is an open problem that had increasing relevance in areas such as additive manufacturing. Efficient handling of the geometry is a significant factor in the solution approach along with the level of fidelity required for the application. 

\subsubsection*{Data}
Across all problem types there are standard data sets and data generators that provide a useful means to test the effectiveness of solution approaches. Many of these are listed on the EURO Special Interest Group on Cutting and Packing (ESICUP) website\footnote{\url{https://www.euro-online.org/websites/esicup/data-sets/}}.

\subsection[Disaster relief and humanitarian logistics (Bahar~Y.~Kara \& Özlem~Karsu)]{Disaster relief and humanitarian logistics\protect\footnote{This subsection was written by Bahar~Y.~Kara and Özlem~Karsu.}}
\label{sec:Disaster_relief_and_humanitarian_logistics}

Humanitarian logistics (HL) is one of the key application areas that Operational Research (OR) has been offering solutions to improve the welfare of the society under difficult circumstances. Humanitarian logistics problems are highly relevant in today's world due to various challenges including but not limited to, climate change and its consequences (increases in extreme weather events such as heatwaves, floods), natural disasters (e.g., earthquakes, tsunamis), man-made conflicts (e.g., Syria and Ukraine crises) and health-related catastrophic events (e.g., pandemics). Humanitarian logistics operations involve  complex systems with multiple stakeholders such as victims, planners, public/private service providers, volunteers, general public and media, each with their own preferences and priorities; and the inherently challenging decisions of scarce resource allocation have to be made over a long time span, under high uncertainty. We use humanitarian logistics as an umbrella term, which covers relief logistics, disaster logistics, and development logistics. The humanitarian literature uses all these terms interchangeably. Actually, disaster logistics and relief logistics should refer to the cases where a disaster is/was/is expected to be in action whereas development logistics refers to cases which aim to improve daily life. 

In relief logistics, the operations require advanced planning, hence the authorities are constantly facing challenges in the four main stages of: \textit{mitigation, preparedness, response, and recovery} \citep{altay2006or_BYKOK,ccelik2012humanitarian_BYKOK,kara2017humanitarian_BYKOK}. In Table \ref{tab:prob_BYKOK}, we list some of the most frequently considered problems, categorised based on the phase that they arise. As seen in the table, mitigation and preparedness phases mostly consist of activities related to planning, which involves network design, location, allocation and routing operations as well as provisioning processes that include inventory and other supply chain-related decisions \citep[see, e.g.,][for applications in location and prepositioning, respectively]{balcik2008facility_BYKOK,rawls2010pre_BYKOK}. 

Response activities occur after the crisis or the disaster hits. In this phase, the aim is providing a rapid response, prioritising the survival needs. This, however, does not preclude considering efficiency in these operations as the system requires scarce resource allocation such as personnel, equipment and supplies across demand points, invoking a need for good decision support mechanisms. In line with this need, a large body of work is devoted to the application problems arising in this phase. Finally, recovery phase focuses mainly on debris management and infrastructure repair and restoration. Most of the mentioned operations involve additional decisions regarding workforce planning and scheduling and require structured methods for data management, information sharing and coordination, which are key for effective response \citep{altay2014challenges_BYKOK}. Some of these models require quantification of human suffering due to lack of services or goods and the deprivation cost can be used for this purpose, as discussed in detail in \cite{holguin2013appropriate_BYKOK}. 

\begin{table}[htbp]
	\centering
	\caption{Problems in Relief Logistics}
\resizebox{\textwidth}{!}{	\begin{tabular}{clll}
		\textbf{Phase} & \multicolumn{1}{l}{\textbf{Main Problems considered}} & \multicolumn{1}{l}{\textbf{Main decisions involved}} & \multicolumn{1}{l}{\textbf{Main concerns}} \\
				\hline
		Preparedness and mitigation & \multicolumn{1}{l}{Infrastructure  and network structuring } & Network design & \multicolumn{1}{l}{Connectivity} \\
				&  &  (mainly for strengthening purposes)  &  \\
		\cline{2-4}
		& \multicolumn{1}{l}{Risk assessment } & Prioritization  & \multicolumn{1}{l}{Data preparation } \\
			&  & (of roads, buildings, arcs)  & for further analysis \\
		\cline{2-4}
		& \multicolumn{1}{l}{Evacuation planning } & Location of gathering points & \multicolumn{1}{l}{Traffic} \\
		&       & Allocation of evacuees & \multicolumn{1}{l}{Accessibility} \\
		&       & Routing  & \multicolumn{1}{l}{Evacuation Time} \\
		&       &       & \multicolumn{1}{l}{Behavioral factors} \\
		\cline{2-4}
		& \multicolumn{1}{l}{Shelter location (and allocation)} & Selecting among potential locations  & \multicolumn{1}{l}{Ensuring accesibility} \\
		&       &       & \multicolumn{1}{l}{Shelter utilization} \\
		&       &       & \multicolumn{1}{l}{Infrastructure (Risk)} \\
		&       &       & \multicolumn{1}{l}{Behavioral factors} \\
		\cline{2-4}
		& \multicolumn{1}{l}{Prepositioning} & Locating point of distributions  & \multicolumn{1}{l}{Budget} \\
		&       & Provisioning (of supplies) & \multicolumn{1}{l}{Accessibility} \\
		&       & Related supply chain decisions  & \multicolumn{1}{l}{Fairness} \\
		&       &    (supplier selection, allocation and routing)   & \multicolumn{1}{l}{Speed} \\
		&       &       & \multicolumn{1}{l}{Efficiency} \\
		\cline{2-4}
		& \multicolumn{1}{l}{Supply chain (procurement) planning } & Contract design  & \multicolumn{1}{l}{Quality} \\
			&  &   &  Efficiency\\
		\hline
		Response & \multicolumn{1}{l}{Damage assessment } & Demand assessment & \multicolumn{1}{l}{Speed} \\
		&       & Infrastructure assessment & \multicolumn{1}{l}{Accuracy} \\
		\cline{2-4}
		& \multicolumn{1}{l}{Search and rescue operations} & Team formation and allocation  & \multicolumn{1}{l}{Speed} \\
		&       &       & \multicolumn{1}{l}{Fairness} \\
		&       &       & \multicolumn{1}{l}{Capacity} \\
		\cline{2-4}
		& \multicolumn{1}{l}{Evacuation } & Location & \multicolumn{1}{l}{Speed} \\
		&       & Allocation & \multicolumn{1}{l}{Risk} \\
		&       & Routing &  \\
		\cline{2-4}
		& \multicolumn{1}{l}{Shelter management} & Location & \multicolumn{1}{l}{Physical environment} \\
		&       & Allocation & \multicolumn{1}{l}{Risk} \\
		&       & Related supply chain decisions  & \multicolumn{1}{l}{Accessibility} \\
		&       &    (distribution, routing)   & \multicolumn{1}{l}{Fairness} \\
		&       &       & \multicolumn{1}{l}{Utilization} \\
		\cline{2-4}
		& \multicolumn{1}{l}{Donation management/ } & Allocation & \multicolumn{1}{l}{Fairness } \\
		&   Resource allocation    & Routing  & \multicolumn{1}{l}{Efficiency} \\
		&       & Inventory Management  &  \\
		\hline
		Recovery & \multicolumn{1}{l}{Debris management } & Network design & \\
		&       & Prioritization and scheduling  &  \\
			&  & (deciding which nodes/arcs to clean first)  &  \\
		\hline
	\end{tabular}}%
	\label{tab:prob_BYKOK}%
\end{table}%

Not all humanitarian operations are triggered by challenges stemming from a single well-defined event. There are crises that can not be attributed to a single cause e.g., famine in under-developed countries. There are well-established efforts in the development logistics literature to alleviate the effects of such crises. Some examples are global health projects for increasing access to health coverage  and fighting diseases that occur in low and middle-income countries on a wide-scale, such as malaria and AIDS, through distribution of effective tools and/or medication.  Vaccine development and distribution to poorer regions as well as distribution of other basic needs to the deprived populations: food aid distribution \citep{foodaid2015_BYKOK, mahmoudi2022decision_BYKOK}, clean water network design and distribution \citep{LAPORTE2022_BYKOK}, energy, education and hygiene provision are also widely considered. A recent trend is utilising cash and voucher distribution whenever possible, since it is a method that respects human dignity, avoids complications of relief item logistics and supports the local market \citep{karsu2019refugee_BYKOK}. 

The recent COVID-19 pandemic has also motivated a wide range of HL applications such as personnel protective equipment allocation and distribution, frontline workforce planning and system design for testing, tracing and vaccination \citep{FARAHANI20231_BKOK}.

The HL literature also integrates newer technologies to the delivery systems: There has been recent attempts to  use drones in the last mile distribution as they constitute a convenient tool to reach remote areas in short time \citep[examples include delivery of blood samples, vaccines and food aid; see also][]{GENTILI2022108057_BKOK,GHELICHI2021105443_BKOK,alfandari2022tailored_BKOK}.

OR offers decision support for humanitarian settings based on a wide range of quantitative and qualitative tools. Mathematical modelling and optimisation is used in almost all problems arising in HL to make the related location, allocation, routing and network design decisions. The models are shaped by the priorities in the phase and constraints imposed by the physical infrastructure, resource availability as well as the social, economic and cultural environment. 
As the underlying technical problems are difficult to solve, various mathheuristic and metaheuristic approaches are employed. There is also an increasing trend in using system dynamics \citep{besiou2021system_BYKOK} and empirical analysis \citep{pedraza2016empirically_BYKOK}.

There is no one-size-fits-all methodology but some key properties require specific methods to be used. In most relief logistics problems the environment is highly stochastic, calling for applications of forecasting and stochastic programming. The uncertain factors include but are not limited to the number of affected individuals, the extent of the effect, types of needs, and usability status of the infrastructure and other resources. Multiple stakeholders and conflicting criteria are involved, requiring multicriteria decision making \citep{ferrer2018multi_BYKOK} approaches. Unlike commercial logistics systems, fairness is a key concern in humanitarian settings. Fairness or equity, however, is hard-to-quantify and is context dependent: a rule that is considered fair under some circumstances may not be deemed so in others. The policy makers may want to prioritise beneficiaries based on attributes such as socio-economic status and hence ensure vertical equity or may consider them indistinguishable and seek horizontal equity \citep{KARSU2015_BYKOK}. 
  
HL has been receiving attention with an increasing rate, which led to many review studies that the interested reader  can refer to: see, e.g., \cite{luis2012disaster_BYKOK, CELIK201647_BYKOK, aringhieri2017emergency_BYKOK, besiou2018or_BYKOK,8935364_BYKOK,DONMEZ20211_BYKOK}. See also \cite{kunz2017relevance_BYKOK} for a discussion on how to make humanitarian research more impactful for humanitarian organisations and beneficiaries. 
	
Next we categorise the future of the humanitarian applications to motivate and direct new researchers:
\begin{itemize}[noitemsep]
\item 	OR is responsive to the difficulties the world is facing and humanitarian challenges are no exception to this. The recent COVID-19 pandemic has shown that relief logistics can be applied in health-related crisis management to provide quick, effective, efficient and fair responses to health care problems.  WHO ``urges countries to build a fairer, healthier world post-COVID-19" and this is doable with good humanitarian practices relying on OR.  The recent efforts in designing fair and efficient systems using OR would contribute in addressing inequities in health and welfare, which have been exacerbated by the pandemics. We believe that there is still room for improvement in adopting a holistic approach and conducting multidisciplinary work when designing such systems. One example is the vaccine implementation and roll-out problem: conceptualising  this problem as a sole logistics problem may not be the best practice as the success of any design highly depends on human behaviour. People have different views, risk attitudes and preferences over available options, which affects how any proposed policy will perform.  Incorporating such behavioural factors is an important yet scarcely studied issue. 
\item The underlying technical problems in the HL domain are hard to solve due to uncertainty in various parts of the system, lack of (reliable) data and multiple criteria that are involved. Moreover, a significant portion of these problems are combinatorial optimisation problems, i.e. they require choosing from a prohibitively large set of solutions that are implicitly defined by constraints of the system. Therefore advances in OR methodology to obtain better, quicker solutions to optimisation problems, and in data analytics on handling big data such as the one obtained through geographic information systems (GIS), would pave the way for quicker and better response. Effective data analysis would especially help when learning from past practice. Indeed, lessons learned from humanitarian supply chain practice can also be used in managing supply chain disruptions in other sectors, as discussed in \citep{kovacs2021lessons_BYKOK}. 
\item UN's Sustainable Development Goals emphasise the global challenges faced including  poverty, inequality, climate change, environmental degradation, peace and justice \citep{SDG_BYKOK}. As stated in \citep{report_BYKOK}, ``Increases in extreme weather events and  climate change can compound risks of international food shocks, water insecurity, conflict and other humanitarian emergencies and crises. Difficulty of access to critical resources such water and food may trigger migrations or	exacerbate conflict risks.''  All these areas are, by definition, related to humanitarian operations, hence humanitarian logistics has a lot to offer in these domains. 
\item The Turkey/Syria earthquakes in February 2023 have clearly demonstrated the importance of effective coordination and strategic planning. Thus, we would like to emphasise the need for collaborative research that brings field expertise (of e.g,. municipalities, NGOs and volunteers) and academic know-how together. 
\end{itemize}

\subsection[E-commerce (Charlotte~Köhler \& Tom~Van~Woensel)]{E-commerce\protect\footnote{This subsection was written by Charlotte~Köhler and Tom~Van~Woensel.}}
\label{sec:Ecommerce}

\subsubsection*{What is E-Commerce about?}
E-commerce deals with the transactions of goods and services through online communications (computers, tablets, smartphones, etc.). Both business-to-business (B2B) and business-to-consumer (B2C) realisations are observed in practice. In B2B, companies operate their supply chains through online networks. In B2C, products and services are sold directly to consumers. E-commerce sales steadily increased for years and amounted to \$5,211 billion worldwide in 2021, with the pandemic being a major contributor.\footnote{https://bit.ly/3dfiwDW, 2022-08-08}$^{,}$\footnote{https://bit.ly/3SzAWzf, 2022-08-08}

\textit{E-fulfilment} describes all fulfilment activities for e-commerce. All necessary steps for a customer to receive an order after placing are thus referred to as the e-fulfilment process. Due to the nature of the e-commerce domain, these e-fulfilment activities often occur in a city context \citep{savelsbergh201650th_CKTVW}. E-fulfilment processes are planning intensive, and creating a profitable business in this environment is challenging. Customer service expectations are high, however, and the customer is more and more in the lead on how and where their orders need to be delivered \citep[the ``logsumer'' takes an active role in time, price, quality, and sustainability decisions of logistic services][]{DHL_CKTVW}.

The e-fulfilment process can be divided into three steps, namely (\textit{i}) order acceptance, (\textit{ii}) order assembly, and (\textit{iii}) order delivery \citep{campbell2005decision_CKTVW}. For most online companies, these steps take place separately, one after the other. However, new on-demand companies have considerably shortened lead times and perform these steps simultaneously \citep{wassmuth2022demand_CKTVW}.

During \textit{order acceptance}, customer requests arrive on a retailer's website and ask for service. As fulfilment capacities are limited (for example, delivery capacities), the retailer wants to accept the most profitable subset of all customer requests. However, customer requests arrive one at a time. Thus, the retailer does not know the total delivery costs until all customers are accepted and the final delivery route is planned. In addition, when a request is accepted, the retailer does not know whether requests with higher revenues will arrive afterward for which capacity should have been reserved. To estimate costs, vehicle routing methods are adapted for usage as customer acceptance mechanisms \citep[e.g.,][]{ehmke2014customer_CKTVW,kohler2019evaluation_CKTVW}). Revenue management methods are used to allocate capacities to high revenue requests \citep[e.g.,][]{cleophas2014deliveries_CKTVW,klein2019differentiated_CKTVW}). However, since decisions in the online environment must be made instantly, the use of complex and, thus, computationally intensive solution methods is limited. 

The warehouse picking and consolidating ordered goods are summarised under \textit{order assembly}. Before this, the retailer must decide on the location and design of the warehouses. Choosing the location is closely linked to the fulfilment capacity of the retailer and must be well-planned. The design of the warehouse determines the efficiency of picking the ordered items. Finding efficient picking strategies to reduce retailer costs is studied in, for example, \cite{Schiffer2022_CA_MB}. Lastly, the retailer must determine the optimal stock level of items. Given the short lead times in e-commerce, this task must be completed before customer requests arrive. The task is closely linked to the research field of inventory management, where techniques such as forecasting \citep{ulrich2021distributional_TVWCK} or artificial intelligence \citep{albayrak2023applications_TVWCK} are commonly used to address this challenge effectively.

For \textit{order delivery}, routes are planned for all accepted orders. For e-commerce, the last-mile delivery is usually towards the customer's location, i.e., the consumer's home or company site \citep{agatz2008fulfillment_CKTVW}, leading to a magnitude of fragmented delivery locations with small drop sizes. Significant challenges arise from how these last-mile deliveries (routes) are designed. Delivery route planning is closely related to the established field of vehicle routing, and approaches are being adapted for use in e-fulfilment \citep[e.g.,][]{emec2016adaptive_CKTVW}. Two-echelon routing systems are often considered to maintain economies of scale and satisfy the emission zone requirements in the cities \citep{sluijk2022chance_CKTVW,sluijk2022two_CKTVW}). In most cases, delivery is made by conventional delivery vehicles. However, individual retailers are also starting to bring orders to customers in the city centre using bikes. We also see drones \citep[e.g.,][]{ulmer2018same_CKTVW,dayarian2020crowdshipping_CKTVW}) and robots \citep[e.g.,][]{simoni2020optimization_CKTVW}.

\subsubsection*{E-fulfilment Challenges}
E-fulfilment processes present several challenges. 
For \textit{unattended deliveries}, delivery is possible without the customer being present. Pick-up point delivery enhances the efficiency of the delivery operations via consolidation opportunities. Consumers can also find it a more convenient delivery option than waiting for the delivery at home. There is a need for incentive mechanisms to increase the attractiveness of pick-up points (e.g., reduced delivery price). \cite{Galiullina2022_CKTVW} study this problem as a trade-off between routing cost savings gained from steering the customer demand and the investments required to influence customer behaviour. Another challenge is to find the optimal locations for pick-up points, such that delivery costs are minimised and customers still have convenient access, which is, for example, considered in  \cite{lin2020last_CKTVW} and \cite{wang2020locating_CKTVW}. The customer must accept the delivery herself for \textit{attended deliveries}, e.g., to prevent grocery spoiling. To avoid delivery failures, the customer and the retailer usually agree on a delivery time window.

Customers expect short time windows, which increase the retailer's delivery costs  \citep{kohler2020flexible_CKTVW}. As the time windows assignment to orders is crucial for the retailer's profitability, several approaches consider balancing demand along the offered time windows. One possibility is to withhold specific time windows from customers and only offer a subset of beneficial time windows. \cite{campbell2005decision_CKTVW} and \cite{cleophas2014deliveries_CKTVW} consider routing costs and customer value and only offer time windows to customers that are expected to maximise the profit. Another possibility is to assign prices to time windows to nudge customers to specific time window options \citep{campbell2006incentive_CKTVW,yang2016choice_CKTVW,klein2019differentiated_CKTVW}. Some approaches consider adapting the time window design to increase routing flexibility. \cite{kohler2020flexible_CKTVW} only offer short time windows to customers when it does not impact the routing costs too much, and \cite{strauss2021dynamic_CKTVW} hand out time window bundles to customers that are only narrowed down to one option once more customers requests are known. 

Recently, many online retailers began offering \textit{on-demand deliveries}, so customers can receive their orders the same day (some grocery stores promise delivery times within a few minutes\footnote{https://bit.ly/3B1H1ww, 2022-09-09}). Shortening lead times poses another challenge as there is almost no time for planning or consolidation of orders available. The approach presented by \cite{klapp2020request_CKTVW} hence supports retailers in deciding which customers can be promised an immediate delivery and which can only be served from the next delivery day. \cite{ulmer2018same_CKTVW} investigate how the number of same-day deliveries can be increased if delivery is not only done by vehicles but additionally by drones. In \cite{banerjee2022fleet_CKTVW}, the authors examine how retailers must allocate their delivery capacity to cover same-day delivery needs per service area. 

For delivery in an urban context, high demand in densely populated areas often goes hand in hand with high traffic and unreliable travel times, and vice versa. \cite{ehmke2014customer_CKTVW} therefore create acceptance mechanisms that present the customer with a time window offer that is as reliable as possible so that the customer does not notice an unforeseen change in travel times. \cite{kohler2019evaluation_CKTVW} test the suitability of customer acceptance mechanisms for more and less densely populated areas to derive how well different routing mechanisms approximate delivery times. 

Another ongoing challenge is the increasing prevalence of customers being granted the option to \textit{return ordered items} free of charge by many companies. As a result, the e-fulfillment process expands beyond the three steps outlined earlier to include the management of returns. Despite the typically high return rates that result in substantial additional costs for retailers, offering a return option is still profitable due to the subsequent improvement in customer satisfaction and retention \citep{rintamaki2021customers_TVWCK}. The management of returns can be perceived as reverse order delivery, leading to routing challenges related to those presented earlier. To mitigate costs, several studies, such as \citet{mahar2017store_TVWCK} and \citet{yan2022whether_TVWCK}, explore the implementation of in-store returns.

\subsubsection*{Operational Research challenges: Time, Timing, and Data}

The \textit{time} dimension involves all dimensions to how key elements are (conceptually) modelled with regards to the time (e.g., travel times or handling times). Identifying the time features in modelling and solution methodologies are essential qualifiers for realistic model representations. 

The \textit{timing} dimension involves all actions at a particular point or in a period when something happens (e.g., a new order arrives). Timing considers synchronisation issues where, for example, vehicles need to meet at a certain point in time and geographical location. \cite{drexl2012synchronization_CKTVW} presents a survey of vehicle routing problems with multiple synchronisation constraints. Synchronisation requirements between the vehicles relate to spatial, temporal, and load aspects. Synchronisation is a challenge, for example, in heterogeneous fleets \citep{ulmer2018same_CKTVW} or, in the case of battery-powered vehicles that must be charged in time. 
\begin{itemize}[noitemsep]
    \item \textit{Offline} means that we do the planning and scheduling before the execution, often assigned to tactical planning. Data is estimated (forecast) based on past observations, and the operations are planned based on that. For example, \cite{agatz2011time_CKTVW} use expected demand to decide which time windows should be offered within different parts of the delivery area. \cite{lang2021anticipative_CKTVW} propose a preparation offline phase that serves as input to speed up decisions during later online customer acceptance. The data considered could be either time-independent (i.e., independent of time) or time-dependent (i.e., the data has a time-stamp). For example, travel times can be modelled time-independent (i.e., constant speed) or time-dependent \citep[e.g.,][]{spliet2018time_CKTVW}.
    \item \textit{Online} refers to the optimisation in real-time, where revealing new data and planning and scheduling operations happen simultaneously. The terms ``dynamic'' or ``operational planning'' is also often used. As time is critical in online planning, methods are always limited by their solution time. Instead of finding a routing solution, delivery costs are approximated \citep[e.g.,][]{yang2017approximate_CKTVW,lebedev2021dynamic_CKTVW} or a simple routing heuristic is applied \citep[e.g.,][]{mackert2019choice_CKTVW,klein2018model_CKTVW}). Alternatively, customer choice is estimated simply \citep[e.g.,][]{campbell2006incentive_CKTVW} instead of complex and time-consuming customer choice modelling. \cite{vanderhagen2022machine_CKTVW} uses a machine learning approach to fasten up feasibility checks of time windows offered during order acceptance. 
\end{itemize}

The \textit{data dimension} refers to how the data and observations are modelled. The data can be handled deterministic or stochastic, or we observe the realised data. Most models assume deterministic data and build their solution approach around this notion. More and more researchers, however, recognise the challenge of adequately representing reality in their models. \cite{yang2016choice_CKTVW} use booking data of an online grocer to estimate realistic customer behaviour. \cite{kohler2022data_CKTVW} investigate how to accept high revenue requests by applying a sampling procedure with booking data from an e-grocer in Germany. 

\subsubsection*{Relevant literature}
\cite{agatz2013revenue_CKTVW} provide the first overview of how retailers can manage e-fulfilment processes. A recent review on e-fulfilment for attended home deliveries can be found in \cite{wassmuth2022demand_CKTVW}. We refer the reader to \cite{fleckensteinRecentAdvancesIntegrating2022_AKSJF} and \cite{snoeck2020revenue_CKTVW} for a focus on routing and revenue management methods in e-fulfilment, respectively.

\subsection[Education (Jill~Johnes)]{Education\protect\footnote{This subsection was written by Jill~Johnes.}}
\label{sec:Education}

Education spans activity from kindergarten, through primary and secondary schooling, to higher education. The earlier years of education are often compulsory reflecting the premise that an educated workforce is crucial to economic performance. The extent to which education is publicly funded varies from one level of education to another, as well as from one country to another depending on the local view concerning the social return on investment. Public funding for education alongside the role education and training play in the performance of an economy therefore make education a prime context for application of operational research (OR) tools. This section provides a brief overview of some of the main areas.

Many OR methods can be useful to the policy maker for macro-planning and financial allocation purposes. Forecasting student numbers can be done using Markov chain models \citep{Nicholls2009-vc_JJ,Brezavscek2017-cj_JJ} or machine learning (ML) and artificial intelligence (AI) \citep{Yan2021_JJ} – the importance of AI applications to education will be further expanded later. Allocation of finances is typically supported by multi-objective decision analysis \citep{Cobacho2010-mn_JJ}. 

One important aspect of resource allocation relates to the efficient use of resources. Availability of published education data in many countries provides an opportunity to examine the “black box” of education production. Consequently, there is a long-standing literature surrounding efficiency in education, typically a not-for-profit context where conventional measures of performance are inappropriate. 

Early studies of efficiency in higher education applied deterministic ordinary least squares methods to university-level data to examine efficiency in the production of specific outputs \citep{Jauch1975-hd_JJ,Johnes1990-zt_JJ} while schools adopted multilevel modelling methods to derive performance insights from pupil-level as opposed to school-level data \citep{Woodhouse1988-fr_JJ}. But the multi-product nature of production in education establishments means that looking at inputs separately provides only a partial picture. The tools of multiple-criteria decision analysis such as principle components, the analytic hierarchy process and co-plot have therefore been adopted to examine and visualise the many dimensions more easily \citep{Johnes1996-pk_JJ,Paucar-Caceres2005-kz_JJ,Mar-Molinero2007-pi_JJ}.

Two frontier estimation approaches to analysing efficiency, both of which derive from Farrell (1957), have evolved to address various shortcomings of early approaches. The non-parametric data envelopment analysis (DEA) easily handles the multi-input multi-output nature of production observed in education and provides easily-interpreted measures of efficiency \citep{Charnes1978-nm_SL}. DEA shows each observation in its best possible light (in efficiency terms) by computing a distinct set of input and output weights. This permits the derivation of benchmark observations for each inefficient institution, i.e., the establishment(s) the observation should be looking to emulate to become more efficient. Non-parametric frontier estimation techniques have been applied in the context of education at all levels, providing management information at the institution level, and policy insights at the macro-level \citep{Thanassoulis2011-yg_JJ,Portela2012-bw_JJ,Burney2013-iw_JJ}.

Network DEA provides a more forensic examination of the “black box” \citep{Fare2000-wo_JJ} by breaking down the production process into its component parts, and overall efficiency can be decomposed into efficiency in each of the stages \citep{Wang2019-qo_JJ,Lee2022-an_JJ}. 

When longitudinal data are available, DEA can be used to analyse changes in efficiency using the \cite{Malmquist1953-ve_JJ} productivity index which decomposes productivity change into efficiency and technological change components \cite{Wolszczak-Derlacz2018-xo_JJ}. The method can be used to make comparisons between groups rather than (or as well as) between time periods \citep{Aparicio2017-dh_SL}.

The deterministic non-parametric nature of DEA has been addressed in numerous extensions including by introducing bootstrapping and significance tests \citep{Johnes2006-id_JJ,Essid2010-at_JJ,Papadimitriou2019-ap_JJ}. Second stage analyses which examine the determinants of efficiency also abound \citep{Haug2017-bv_JJ}. This approach is only valid if the hypothesis of separability holds i.e. the variables used in the second stage should only influence the efficiency scores and not the determination of the efficiency frontier \citep{Simar2011-oz_JJ}. The development of separability tests \citep{Daraio2018-ut_JJ} and the robust conditional estimation approach address these issues \citep{Daraio2007-qr_JJ}; their application in education provide more robust and insightful results \citep{Lopez-Torres2021-yl_JJ}.

Stochastic frontier analysis (SFA) provides both parameter estimates (with significance tests) and efficiency estimates which allow for stochastic errors \citep{Aigner1977-df_JJ,Meeusen1977-gq_JJ,Jondrow1982-fq_JJ}. Compared to DEA it is more difficult to model multi-input, multi-output production; most SFA applications in education therefore focus on cost efficiency \citep{Agasisti2016-bg_JJ}, or a single output model \citep{Kirjavainen2012-cr_JJ}, although there are some exceptions \citep{Abbott2009-xv_JJ,Johnes2014-jm_JJ}. The parameter estimates have made SFA popular in the cost function context as scope and scale economies can be estimated and these have useful policy implications \citep{Johnes2005-ue_JJ,Johnes2013-eg_JJ}.

In its basic form, SFA parameter estimates apply to every observation in the dataset. Extensions of the technique include latent class SFA and random parameters SFA which allow the parameters to vary by specific groups (latent class) or by each observation (random parameter). These approaches benefit from the advantages of DEA and SFA although are computationally demanding but have been applied in education to interesting effect \citep{Johnes2011-hd_JJ,Johnes2016-di_JJ}.

The interested reader is referred to comprehensive reviews of the relevant literature \citep{Kao2014-jt_SL,Thanassoulis2016-du_JJ,Witte2017-xg_JJ,Johnes2022-jv_JJ}.

All these OR methods can be applied in the contexts of macro- and micro-level planning and budget allocation. One area at the micro-level for which OR techniques are useful is timetabling. 
Timetabling of examinations and/or teaching is most complex at secondary and tertiary levels and can be viewed as a scheduling problem whereby resources, limited in supply, are allocated to a constrained number of times and locations, with the allocation satisfying stated objectives. Timetabling differs from scheduling in that the resources (staff members) are typically specified in advance rather than being a part of the allocation problem; and while scheduling aims to minimise costs, the objective of timetabling is to realise desirable objectives (e.g., no clashes) as closely as possible
\citep{Petrovic2004-ly_JJ}. Timetablers face both hard and soft constraints in constructing the timetable \citep{Asmuni2009-yb_JJ} and this is therefore a problem which lends itself to solution by various possible OR techniques in the field of combinatorial optimisation. The main approaches are briefly summarised below.

Mathematical programming (particularly integer linear programming) is commonly used in timetabling \citep{Cataldo2017-la_JJ} but often leads to  computationally demanding problems. Heuristics (see below) are introduced for increased efficiency \citep{Dimopoulou2001-rz_JJ}. Case-based reasoning approaches use a past solution (stored in the case base) as the starting point for a new timetable and use similarity measures to identify the optimal solutions \citep{Burke2006-wo_JJ}. These approaches are often problem-specific making them non-transferable. Their computational demands can be addressed by using heuristics \citep{Petrovic2007-ls_JJ}. The multi-criteria approach assumes that there are solutions to the timetabling problem satisfying the hard constraints and then the quality of these solutions is assessed on the basis of how well each one satisfies the soft constraints \citep{Burke2002-dp_JJ}. As with other methods it is often combined with heuristics. 

Heuristics are an increasingly common method for application to timetabling either on their own or in combination with other methods. Low level construction heuristics include largest degree, largest weighted degree, largest colour degree, largest enrolment, saturation degree and random. Extensions include meta-heuristics which work in the search space guiding neighbourhood moves to a solution \citep{Qu2015-oi_JJ}; fuzzy heuristics which can find a best approach in the initial timetable construction phase \citep{Asmuni2009-yb_JJ}; and hyper-heuristics which find or generate appropriate heuristics to solve complex search problems as encountered in timetabling \citep{Qu2015-oi_JJ}. Given their focus, hyper-heuristics have the potential to provide more generalised solutions to timetabling problems than other approaches \citep[see][for a review]{Pillay2016-pp_JJ}.

The interested reader is referred to reviews of educational timetabling approaches \citep{Oude_Vrielink2019-qh_JJ,Tan2021-ba_JJ}.

Finally, an emerging area of interest is the application of AI and ML to education. AI and ML are, as already highlighted,  useful for forecasting as they can analyse rich data on, for example, student numbers, retention, achievement, teaching and quality to derive better predictions and/or understanding of the challenges \citep{Alyahyan2020-eh_JJ,Bates2020-bc_JJ,Teng2022-ra_JJ}. They can also be used in the teaching and learning process itself by personalising each student’s experience for example through use of chatbots, by creating exercises for students which address their weaknesses, and by reviewing assessments highlighting strengths and weaknesses \citep{Teng2022-ra_JJ}. In the growing distance learning education arena where it is more difficult to manage participants who have more freedom to learn when they want and may encounter more distractions, AI can be used to support teachers in gauging student engagement. Thus AI algorithms can be used to develop an online education classroom management system \citep{Wang2021-rm_JJ}. AI and ML have much to offer in education but their potential across all disciplines has yet to be properly explored \citep{Bates2020-bc_JJ}. See \cite{Zawacki-Richter2019-ki_JJ} for further literature.

\subsection[Environment (Judit~Lienert)]{Environment\protect\footnote{This subsection was written by Judit~Lienert.}}
\label{sec:Environment}

Environmental problems are at the centre of societal concerns, and of many research activities, also in Operational Research (OR). It is impossible to comprehensively present this literature. Instead, we first introduce characteristics of environmental problems, then present some insights from specific OR fields, mainly citing review articles. Thereafter, we discuss \textit{Decision Analysis} (\S\ref{sec:Decision_analysis}) methods applied to environmental problems.

Environmental problems are usually multi-faceted and complex \citep{French2005-pc_JL,Gregory2012-zh_JL,Reichert2015-jd_JL}. Since 50 years, such public policy issues are known as ``wicked problems'' \citep{Rittel1973-wn_JL}. In many environmental cases, uncertainties are high. It may be difficult to establish scientific knowledge and adequately model environmental systems. They usually span all sustainability dimensions, which requires making trade-offs between achieving environmental, economic, and societal objectives. Various decision-makers and stakeholders with different world-views are affected, sparking conflicts of interest. Any action may have irreversible or far-reaching consequences over long time horizons. Additionally to the temporal dimension, spatial considerations over varying geographic regions may be important. As wicked problems are typically unique, we might need to find new solutions in each case. OR methods can be highly suitable to disentangle and structure complex environmental problems, and can certainly contribute to problem solving. Below, we present some viewpoints.

\textit{Soft OR} methodologies, and \textit{problem structuring methods} (PSMs; see also \S\ref{sec:Soft_OR_and_problem_structuring_methods}) have been developed to tackle complex real-world problems in interaction with stakeholders \citep{Rosenhead2001-ai_JL,Smith2019-pw_JL}. However, most (review) articles are not specific to environmental problems. Using an applied example, \cite{White2009-dd_JL} explored the potential of soft OR for a city development case. \cite{Marttunen2017-yk_JL} reviewed the combination of PSMs with Multi-Criteria Decision Analysis (MCDA) methods. More complex PSMs seem to be under-utilised, suggesting that their benefits cannot sufficiently inform real-world issues, including environmental decision-making. Similarly, \cite{French2022-kr_JL} argued that literature of quantitative and qualitative OR approaches has developed in silos, and that an intertwined, cyclic understanding of soft and hard OR methods is needed to address complex problems. This author was also concerned that behavioural issues are less well understood in qualitative compared to quantitative model building. Related to problem structuring, \textit{stakeholder analysis and participation} is central to environmental problems. Such research is recently gaining increased interest by OR \citep{De_Gooyert2017-yt_JL,Gregory2020-zf_JL,Hermans2009-rp_JL}. \textit{Behavioural OR} (BOR; \S\ref{sec:Behavioural_OR}) is also gaining momentum \citep{Franco2021-sm_JL}. BOR strongly focuses on interventions, and could increase the understanding of societal and psychological issues in environmental problems. However, to date an environmental perspective is rarely taken. One exception is a conceptual paper about behavioural issues in environmental modelling \citep{Hamalainen2015-fv_JL}. A meta-analysis of 61 environmental and energy cases analysed patterns and \textit{biases} that may occur in the problem structuring phase of decision-making \citep{Marttunen2018-qs_JL}.

\textit{Sustainable supply chains} (see also \S\ref{sec:Supply_chain_management}) have been recently reviewed by \cite{Barbosa-Povoa2018-dw_JL}. These authors took a multi-stakeholder perspective along the supply chain to achieve sustainability goals. They found a predominance of \textit{optimisation} methods applied to strategic decision levels. Most of the 220 reviewed articles focused on economic and environmental aspects, leaving behind the social aspects. Similarly, another review focused on \textit{combinatorial optimisation} (\S\ref{sec:Combinatorial_optimisation}), integrating \textit{reverse logistics} (see also \S\ref{sec:Logistics}) and waste management \citep{Van_Engeland2020-qq_JL}. Among other aspects, the authors emphasised the importance of environmental, social and performance indicators, and stakeholder integration, when dealing with flows of waste products. Taking a life-cycle perspective, usually addressed with \textit{life cycle sustainability assessment} (LCSA), \cite{Thies2019-by_JL} reviewed advanced OR methods for \textit{sustainability assessment of products}. While most articles used ecological indicators, the integration of economic and social indicators is emerging. They concluded that improved systematic procedures for uncertainty treatment are needed, and better integration of qualitative social indicators as well as spatially explicit data.

Other authors reviewed specific OR methods. For instance, \cite{Zhou2018-et_JL} reviewed \textit{Data Envelopment Analysis} (DEA; \S\ref{sec:Data_envelopment_analysis}) for sustainability assessments. Again, economic and environmental measures were well included, but the literature lacked social measures such as customer satisfaction. New DEA methods should be developed that include \textit{social network relationships}. \textit{Mathematical programming} and \textit{optimisation} methods to support biodiversity protection were reviewed by \cite{Billionnet2013-hw_JL}. Some of these difficult combinatorial optimisation problems were well solved, but further research is needed to satisfactorily address real-world biodiversity issues. For conservation management, spatial aspects are central, for example creating biological corridors in the landscape to increase biodiversity. Future research should include the temporal dimension and needs of practitioners. \textit{Robust optimisation} (\S\ref{sec:Stochastic_models}) could be a research avenue to handle uncertainty. A review of invasive species also took a mathematical perspective \citep{Buyuktahtakin2018-bb_JL}. Among other conclusions, research should develop more realistic models to capture spatial and temporal dynamics of invasive species, improve uncertainty treatment and coordination among stakeholders, and include holistic approaches for addressing trade-offs between conservation management and costs of such programs.

\textit{Multi-criteria decision analysis} (MCDA; \S\ref{sec:Decision_analysis}) provides a rich literature addressing environmental decision problems. \cite{French2005-pc_JL} discussed properties of wicked environmental problems from a conceptual point of view (introduced above), and implications for decision support. \cite{Cinelli2014-gf_JL} analysed MCDA methods for sustainability assessments. They voiced some concern that choosing the MCDA methods is rather based on preferences, not analytic considerations. Indeed, text-mining of 3,000 articles provided little evidence that particular environmental application fields used certain methods more frequently, possibly because researchers are unaware of specific method merits \citep{Cegan2017-kf_JL}. To overcome this, \cite{Cinelli2014-gf_JL} classified five MCDA methods using ten criteria important for sustainability assessments, e.g., uncertainty management and testing robustness of results, software, and user-friendliness. We know of two general articles for systematically choosing a suitable MCDA method \citep{Cinelli2020-lq_JL,Roy2013-ii_JL}. Several articles reviewed \textit{decision support systems} (DSS) to identify features and best practices for supporting environmental problems \citep{Mustajoki2017-eq_JL,Walling2020-kn_JL}. Moreover, there are many reviews of MCDA applied to a specific environmental field, but only few were published in OR journals \citep[e.g.,][]{Colapinto2020-mc_JL,Kandakoglu2019-pm_JL}. There is a pronounced increase of articles applying MCDA in all environmental areas \citep[e.g., water, air, energy, natural resources, and waste management;][]{Cegan2017-kf_JL,Huang2011-hy_JL}. Below, we introduce some important findings from \textit{decision analysis}.

Some authors defined frameworks for environmental assessments taking a method perspective. \cite{Gregory2012-zh_JL} proposed \textit{structured decision making} (SDM) to tackle real-world environmental decision problems. Based on \textit{multi-attribute value theory} (MAVT), SDM can be applied without much (mathematical) formalisation. This textbook discusses many practical environmental issues, highlighting solutions from international decision cases. \cite{Reichert2015-jd_JL} proposed a framework for environmental decisions that emphasises uncertainty of scientific knowledge and societal preferences. They argued that theoretical requirements are best met by combining \textit{multi-attribute utility theory} (MAUT) with \textit{scenario planning} and \textit{probability theory}, illustrated with a river management case. Scenario planning has been advocated by various authors for tackling wicked problems \citep{Wright2019-av_JL}. The combination of scenario planning with MCDA has been reviewed by \cite{Stewart2013-ez_JL}, and applied to e.g., nuclear remediation management \citep{Geldermann2009-cz_JL}, coastal engineering under climate change \citep{Karvetski2011-sk_JL}, or water infrastructure planning \citep{Scholten2015-ja_JL}. Scenario analysis has also been combined with probabilistic statements and \textit{mathematical optimisation} for \textit{risk assessment} (see also \S\ref{sec:Risk_analysis}) of nuclear waste repositories \citep{Salo2022-yz_JL}. A climate policy review illustrates the importance of integrating various OR methods to effectively support decision-making \citep{Doukas2020-vl_JL}. The currently predominant evaluation of policy strategies with \textit{climate-economy} or \textit{integrated assessment models} (IAMs) fails to incorporate all relevant uncertainties and stakeholders, and sufficiently address system complexity. These authors proposed integrated approaches, including participatory stakeholder processes with \textit{fuzzy cognitive maps}, combined with MCDA and \textit{portfolio analysis} (PA). PA is especially useful as meta-analysis, and has been reviewed by \cite{Liesio2021-zv_JL}. A PA-framework for environmental decision-making has been proposed by \cite{Lahtinen2017-sh_JL}.

To address \textit{spatial aspects} of environmental problems, \textit{geographic information systems} (GIS) are often combined with MCDA, sometimes also developing DSSs \citep{Keenan2019-mp_JL}. \textit{Risk analysis} (\S\ref{sec:Risk_analysis}) and OR research increasingly focuses on spatial planning \citep{Ferretti2019-jm_JL,Malczewski2020-gs_JL}. One example is the axiomatic foundation of spatial multi-attribute value functions \citep{Harju2019-mb_KVRH,Keller2019-uq_JL}.

Many reviews found that \textit{stakeholder integration} throughout the decision-making process was insufficiently considered, e.g., in flood risk management \citep{Madruga_de_Brito2016-ul_JL} or nature conservation \citep{Esmail2018-wp_JL}. This reflects generally found deficits in \textit{problem structuring} (\S\ref{sec:Soft_OR_and_problem_structuring_methods}), for instance insufficient consideration of social objectives \citep{Kandakoglu2019-pm_JL}, or systematic underestimation of the importance of economic objectives \citep{Marttunen2018-qs_JL,Walling2020-kn_JL}. Moreover, there is a tendency to choose too many objectives in environmental cases \citep{Diaz-Balteiro2017-tl_JL}, potentially inducing \textit{biases} in later stages of MCDA \citep{Marttunen2019-nl_JL}. 

Many reviews emphasised the importance of \textit{uncertainty analyses} in environmental decisions, but this is strongly ignored in practice. One review found that only 19\% of 271 articles included uncertainty analysis, 17\% using \textit{fuzzy techniques} to capture imprecise numbers \citep{Diaz-Balteiro2017-tl_JL}. In another review, 34\% of 343 articles dealt with the imprecision of predictions, 70\% using fuzzy sets, and 20\% \textit{stochastic modelling} \citep{Kandakoglu2019-pm_JL}. In both reviews, only 20\%–30\% of the articles performed \textit{sensitivity analysis}. Additionally, only 5\% of 343 reviewed papers included \textit{temporal aspects} of the environmental decision \citep{Kandakoglu2019-pm_JL}.

As conclusion, OR researchers are widely engaging in environmental problems. Environmental problems are intriguingly complex, thus offering opportunities for inspiring research. Although our evaluation is neither comprehensive nor systematic, some general research needs appear across all OR fields. Many articles emphasised the importance of better integrating practitioners and stakeholders in environmental problems, and of better considering societal objectives. Various fields require improved methods to address the complexities of environmental problems, including appropriately dealing with many types of uncertainties, time, and space. Combining soft with hard OR, improving problem structuring, and integrating questions from behavioural OR will increase the chances of finding sustainable solutions for our worlds' environmental problems. This can also spark cross-disciplinary research over different fields of OR.

\subsection[Ethics and fairness (John~N.~Hooker)]{Ethics and fairness\protect\footnote{This subsection was written by John~N.~Hooker.}}
\label{sec:Ethics_and_fairness}

There is substantial literature on the ethical practice of operational research, surveyed in \cite{Brans2007-xi_JH}, \cite{Ormerod2013-km_JH}, \cite{Tsoukias2021-pz_JH}, and \cite{Bellenguez2023-hc_JH}.  While this is a vitally important discussion, it is useful to consider how the \mbox{science} of operational research can contribute to ethics, as well as how ethics can contribute to the practice of operational research.  It has accomplished this primarily through the development of modelling techniques and algorithms that embody ethical concepts, notably distributive justice.

An operational research model that aims simply to minimise total cost or maximise total benefit may unfairly distribute costs or benefits across stakeholders.  This concern arises in a number of application areas, including healthcare (\S\ref{sec:Healthcare}), disaster relief (\S\ref{sec:Disaster_relief_and_humanitarian_logistics}), facility location (\S\ref{sec:Location}), task assignment, telecommunications (\S\ref{sec:Telecommunications}), and machine learning (\S\ref{sec:Artificial_intelligence_machine_learning_data_science}). It
poses the problem of finding a suitable formulation of equity or fairness that can be incorporated into a mathematical model.

For example, if donated organs are allocated in the most economically efficient fashion, patients with certain medical conditions may wait far longer for a transplant than other patients \citep{McElfresh2018_JH}.  If earthquake shelters are located so as to minimise average distance from residents, persons living in less densely populated areas may have much further to travel \citep{Sibel2019inequity_JH}.  If a machine learning algorithm awards mortgage loans so as to maximise expected earnings, members of a \mbox{minority} group may find themselves unable to obtain loans even when they are financially responsible 
\citep{SaxHuaDefRadParLiu20_JH}. If traffic signals at intersections are timed to maximise traffic throughput, motorists on side streets may have to wait forever for a green light \citep{CheShaTse13_JH}.  

We provide here a brief overview of mathematical formulations of fairness that have been proposed for OR and AI models. Comprehensive treatments can be found in \cite{KARSU2015_BYKOK} and \cite{CheHoo23_JH}. In addition, \cite{OgrLusPioNacTom14_JH} review formulations developed for telecommunications and facility location, two major users of fairness models.  Recent years have seen an enormous surge of interest in fairness criteria for machine learning, many of which are surveyed in \cite{MehMorSaxLerGal22_JH}.

We suppose that the model into which one wishes to incorporate fairness allocates {\em utilities} to a collection of {\em stakeholders}, and we are concerned about the fairness of this allocation.  Utility could take the form of wealth, resources, negative cost, health outcomes, or some other type of benefit.  Stakeholders can be individuals, organisations, demographic groups, geographic regions, or other entities for which distributive justice is a concern.

Fairness models can be divided into three broad categories.  
{\em Inequality measures} are normally used to constrain the degree of inequality in solutions obtained by maximising total benefit or minimising total cost.   Some of these focus on inequalities across individuals, and others on inequalities across groups.  Various statistics for measuring the former are discussed in \cite{Cow00_JH} and \cite{JenVanKer11_JH}.  Perhaps the best known is the {\em Gini coefficient}, widely used to measure income or wealth inequality \citep{Gin12_JH,yitzhaki2013more_JH}. The {\em Hoover index} \citep{Hoov36_JH} is proportional to the relative mean deviation of utilities and represents the fraction of total utility that must be redistributed to achieve perfect equality.  Both the Gini coefficient and the Hoover index can be given linear formulations (\S\ref{sec:Linear_programming}) in an optimisation model by means of linear-fractional programming \citep{ChaCoo62_JH}. {\em Jain's index} \citep{JaiChiHaw84_JH}, well known in telecommunications, is a strictly monotone function of the coefficient of variation.

Inequality between groups, generally referred to as {\em group disparity}, is by far the most discussed type of inequality metric in the machine learning field \citep[\S\ref{sec:Artificial_intelligence_machine_learning_data_science};][]{VerRub18_JH,MehMorSaxLerGal22_JH}.  It assesses whether AI-based decisions (e.g., mortgage loan awards, job interviews, parole, college admission) are biased against a designated group, perhaps defined by race, ethnic background, or gender. Fairness implementations in machine learning typically strive to minimise loss (due to defaults on loans, etc.) while placing a bound on some measure of resulting group disparities. The best known measures are {\em demographic parity} \citep{dwork2012fairness_JH}, {\em equalised odds} \citep{HarPriSre16_JH}, and {\em predictive rate parity} \citep{dieterich2016compas_JH,chouldechova2017fair_JH}, and {\em counterfactural fairness} \citep{KusLofRusSil17_JH,RusKusLofSil17_JH}.  The first two have mixed \mbox{integer}/\mbox{linear} programming (MILP) formulations (\S\ref{sec:Mixed_integer_programming}), and the third a mixed \mbox{integer}/nonlinear formulation.  Weaknesses of group parity measures include a lack of consensus on which one is suitable for a given application \citep{CasCruGreRegPenCos22_JH}, as well as on which groups should be monitored for bias.

A second category of models is concerned with {\em fairness for the disadvantaged}.  They strive for equality, but with greater emphasis on the lower end of the distribution.  The {\em maximin} criterion, based on the famous difference principle of John Rawls, maximises the welfare of the worst-off individual or social class \citep{Raw99_JH}.  It is defended with a social contract argument that has been intensely discussed in the philosophical literature \citep[as surveyed in][]{Fre03_JH,RicWei99_JH}. A more sophisticated form of the principle is lexicographic maximisation ({\em leximax}), which maximises the worst-off, then the second worst-off, and so forth. The {\em McLoone index} compares the total utility of stakeholders at or below the median utility to the utility they would enjoy of all were brought up to the median.  It is based on a concern that no one be disadvantaged but tolerates inequality in the top half of the distribution.  It has been used to assess the allocation of public services, particularly education \citep{Ver96_JH} and can be given an MILP formulation \citep{CheHoo23_JH}.

Criteria that {\em balance efficiency and fairness} can be placed in three categories: convex combinations of efficiency and fairness, criteria from classical social choice theory, and threshold criteria.  Convex combinations provide the simplest approach, as for example a combination of total utility and a fairness measure \citep[e.g.,][]{Sibel2019inequity_JH}.  Other formulations are given by \cite{Yager1997_JH}, \cite{ogryczak2003solving_JH}, and \cite{ReaFroMasSteFerPan21_JH}.  Convex combinations and other weighted averages pose the general problem of justifying a choice of weights, particularly when utility and equity are measured in different units, although \cite{ArgKarYav22_JH} propose a means of avoiding this issue.

The task of balancing fairness and efficiency gave rise to one of the oldest research streams in social choice theory, beginning with the {\em Nash bargaining solution}, also known as {\em proportional fairness} \citep{Nas50_JH}.  Proportional fairness has seen application in such engineering contexts as telecommunication and traffic signal timing \citep{Mazumdar1991_JH,Kelly1998_JH} and elsewhere.  Nash gave an axiomatic argument for the criterion, while \cite{Har77_JH}, \cite{Rub82_JH}, and \cite{BinRubWol86_JH} have shown that it is the outcome of certain bargaining procedures.  {\em Alpha fairness} generalises proportional fairness by introducing a parameter $\alpha$ that governs the importance of fairness, where $\alpha=0$ corresponds to a purely utilitarian criterion, $\alpha=1$ to proportional fairness, and $\alpha=\infty$ to the maximin criterion \citep{MoWal00_JH,VerAyeBor10_JH}.  Alpha fairness has been derived from a set of axioms \citep{Lan2010_JH,LanChi11_JH}, including an ``axiom of partition'' that is largely responsible for the result. It provides an objective function to be maximised that is nonlinear but concave (\S\ref{sec:Nonlinear_programming}). Another criterion, {\em Kalai-Smorodinsky bargaining}, likewise has an axiomatic defence \citep{KalSmo75_JH} and addresses what one might see as a weakness in Nash bargaining, namely that it can result in reduced utility for some stakeholders when the feasible set is enlarged.  The Kalai-Smorodinsky criterion can be viewed as a kind of normalised maximin, as it calls for allocating to each stakeholder the largest possible fraction of his or her potential utility (ignoring other stakeholders) on the condition that this fraction be the same for everyone. This criterion has received support from \cite{Tho94_JH} as well as the ``contractarian'' ethical philosophy of \cite{Gau87_JH} and has been recommended for wage negotiations and similar applications \citep{Ale92_JH}.

{\em Threshold criteria} are of two types.  One, based on an {\em efficiency threshold}, imposes a maximin objective until the efficiency cost becomes unacceptably great, at which point some stakeholders are switched to a utilitarian criterion.  The other, based on an {\em \mbox{equity} threshold}, imposes a utilitarian criterion until inequity becomes unacceptably great, at which point a maximin criterion is introduced.  Originally proposed for two stakeholders \citep{WilliamsCookson2000_JH}, the threshold criteria have been extended to $n$ persons, using an MILP formulation for the former \citep{HooWil12_JH} and a linear programming model for the latter \citep{ElcHooZha22_JH}.  A parameter $\Delta$ regulates the equity/efficiency trade-off in both models, in that stakeholders with utility within $\Delta$ of the worst-off are given special priority.  Thus, the parameter $\Delta$ may be interpretable in a practical situation in a way that $\alpha$ in the alpha fairness criterion is not.  Both threshold criteria inherit a weakness of the maximin criterion, namely that they may be insensitive to the equity position of disadvantaged stakeholders other than the very worst-off.  This has been addressed for the efficiency threshold by combining a utilitarian criterion with a {\em leximax} rather than a maximin criterion.  \cite{McElfresh2018_JH} accomplish this by assuming there is a pre-existing priority ordering of stakeholders.  \cite{CheHoo21_JH} avoid this assumption by giving greater priority to stakeholders with utilities closer to the lowest, and by solving a sequence of MILP models to balance the leximax element with total utility.

Fairness modelling is a relatively recent research program in operational research that may forge new connections with other fields.  Much as interactions between OR and economics, management, and engineering have been mutually beneficial on both a theoretical and practical level, collaboration with ethicists on the precise formulation of fairness concepts may bring similar benefits to both ethical philosophy and operational research.

\subsection[Finance (Michèle~Breton)]{Finance\protect\footnote{This subsection was written by Michèle~Breton.}}
\label{sec:Finance}

The use of mathematical models and numerical algorithms to solve an extensive range of problems in finance is widespread, by both researchers and practitioners. In this subsection, we offer an overview of some established models and discuss a selection of the corresponding OR approaches and techniques.

\subsubsection*{Resource allocation models}
As in any other industry, the optimal allocation of resources to activities is a central problem in finance. Prototype models include short-term cash flow management (a linear program), portfolio dedication and immunisation (linear programs), capital budgeting (knapsack problem), asset/liability management (stochastic program with recourse), and portfolio selection (quadratic program).

The \emph{portfolio selection} model introduced in \cite{Markowitz52_MB} and discussed in \cite{markowitz00_MB} is one of the best known optimisation models in finance. This mean-variance model consists of determining the composition of a portfolio of risky assets -- a vector of weights -- where the performance (to be maximised) is measured by the expected portfolio return, a linear function of the assets' weights, while the risk (to be minimised) is measured by the variance of the portfolio return, a quadratic function of the weight vector. The resulting optimisation problem gives rise to a convex quadratic program. This model and its analytical properties led to a formalisation of diversification as a strategy to mitigate risk and to important developments in financial theory.

While the Markowitz model represents a considerable simplification of the portfolio management problem, mean-variance optimisation models are still very much applied in practice. Straightforward variations of the Markowitz model can account for various constraints on the asset weights (e.g., bounds, minimum participation, regulatory or operational restrictions, logical constraints, etc.), yielding mixed integer quadratic programs.

Mean-variance models rely on sets of parameters describing the expected returns and their correlation matrix in the universe of the set of considered assets. Various forecasting approaches (\S\ref{sec:forecasting}) have been proposed to obtain estimates of these parameters, often relying on some assumptions about the correlation structure. One important issue related to the use of mean-variance optimisation models is the sensitivity of their solutions to the estimated parameter values \citep{michaud89_MB}, specifically when the feasible region is relatively unconstrained. Robust optimisation (see also \S\ref{sec:Stochastic_models}) is increasingly used to limit the estimation risk of mean-variance portfolio solutions \citep{ismail19_MB, yin21_MBl, blanchet22_MB}.

Another related limitation of mean-variance models is the fact that they are static models, that is, expectations and correlations of asset returns are assumed to be known and constant over the planning horizon. In practice, estimations are updated periodically to reflect changes in data, and portfolios are rebalanced to the optimal composition corresponding to the new set of estimates. Small perturbations in the values of the input parameters may lead to significant changes in the composition of the portfolio from one period to the next (for instance, when groups of assets have similar characteristics). When the costs associated with changing the composition of the portfolio are significant, a static model may be far from optimal. The portfolio selection problem can be readily extended to a multi-period context, allowing to account for transaction costs and/or to use a dynamic model of the evolution of asset prices over time \citep{li00_MB}. Dynamic models can also account for additional frictions, such as taxes on capital gains or losses \citep{dammon01_MB}. The resulting dynamic portfolio selection problem may be a large-scale stochastic dynamic program (\S\ref{sec:Dynamic_programming}; \S\ref{sec:Stochastic_models}). Moreover, risk measures based on portfolio variance are not additively separable, precluding the efficient use of dynamic programming. \cite{Steinbach01_MB} proposes a solution approach based on scenario decomposition.

\subsubsection*{Risk management}
While, in OR, the classical way to deal with decisions under risk is utility theory, finance models usually take a different approach by directly measuring and/or pricing risk. Various measures, such as variance, semi-variance, Value at Risk (VaR) or Conditional value at risk (CVaR) have been proposed to characterise risk\footnote{VaR is a quantile of the distribution of investment losses over a specified period, while CVaR is the amount of the expected losses, conditional to being above the VaR threshold}. VaR is effectively concerned with computing quantiles of the predictive distrubution (see also \S\ref{sec:forecasting} and \S\ref{sec:Power_markets_and_systems}). In the following paragraphs, we present two contrasting families of approaches to financial risk management.

\emph{Diversification and hedging }approaches are closely related to the resource allocation models presented above. They consist in setting up and managing portfolios of securities with desirable properties. Diversification is effective in reducing risk that is uncorrelated across securities, while hedging is used to reduce systematic risk, for instance by holding securities exposed to the same risk factors to eliminate uncertainty, or by buying insurance in the form of derivative contracts. In general, hedging positions must be continuously adjusted to account for the time evolution of risk factors and security prices. In addition, investment portfolios are
often required to satisfy institutional or regulatory constraints. Risk mitigation portfolio planning problems give rise to dynamic stochastic mathematical programs. In recent years, CVaR has become prominent for measuring portfolio risk; CVaR is well-suited to measure down-side risk in skewed distribution and, as shown in \cite{artzner99_MB}, it has the desirable properties of a \emph{coherent} risk measure. Moreover, the use of CVaR in optimisation models gives rise to convex or linear programs, allowing to efficiently solve the large-scale problems encountered in practice \citep{rockafellar00_MB, andersson01_MB, rockafellar02_MB}.

\emph{Risk pricing} approaches rather seek to evaluate the consequence of unpredictable events and are notably used for the management of credit and counterparty risk, that is, the risk that the issuer of a security (for instance, a corporate bond) will not be able to meet its future obligations. A variety of models have been proposed to evaluate the VaR of debt instruments, mainly for the purpose of assessing regulatory requirements ensuring that financial institutions put aside sufficient capital to sustain eventual losses. \cite{crouhy00_MB} presents a review of methodologies currently proposed by the industry to evaluate the probability and consequences of default events. Most approaches used in the industry to price credit and counterparty risk are based on probabilistic models or
Monte-Carlo simulation (\S\ref{sec:Simulation}) and, as such, cannot account for strategic behaviour by the debtor or the lender \citep{breton18_1_MB}.

\subsubsection*{Asset pricing}
Most asset pricing models are founded on an absence of arbitrage assumption, which is usually motivated by the efficiency of markets. Under this assumption, the value of a financial asset is equal to the expected value of its future payoffs, under a suitable probability measure. One specific application is the valuation and replication of contingent claims, such as financial options. The contribution of OR to this area lies in the development and implementation of efficient numerical pricing methods for complex financial securities.

Starting from the binomial tree model of \cite{cox79_MB}, \emph{numerical methods} for option pricing include Monte-Carlo (\S\ref{sec:Simulation}) and quasi-Monte-Carlo approaches \citep{acworth98_MB, lEcuyer09_MB}; dynamic programming (\S\ref{sec:Dynamic_programming}) and approximate dynamic programming models accounting for optimal exercise strategies \citep{ben02_MB,longstaff01_MB}; and robust control models (\S\ref{sec:Stochastic_models}) accounting for transaction costs and model uncertainty \citep{davis93_MB,bernhard05_MB, bandi14_MB}. Numerical algorithms developed for option pricing have also been applied to the valuation of numerous instruments, including corporate bonds, credit derivatives, contracts, and, under the designation of \emph{real options}, managerial flexibility \citep{trigeorgis96_MB, schwartz04_MB, dixit09_MB}.

In the context of \emph{algorithmic trading}, asset pricing algorithms have been revisited using artificial intelligence approaches, for instance by using machine learning to identify factor models or reinforcement learning to compute optimal exercise strategies \citep{dixon20_MB, gu20_MB}, or by augmenting the set of covariates with textual data \citep{algaba20_MB}.

\subsubsection*{Strategic interactions}
Decisions made by investors, firms, financial institutions, and regulators have a direct impact on asset values, returns and risk. Players in the financial sector have competing interests and interact strategically over time, and these interactions are recognised in many game-theoretic models of investment and corporate finance (\S\ref{sec:Game_theory}). Important issues include market impact and market manipulation, option games, strategic exercise of real options, agency conflicts, corporate investment, dividend and capital structure policies, financial distress, and mergers and acquisitions.

\emph{Optimal execution} refers to the determination of a trading strategy minimising the expected cost of trading a given volume over a fixed period, accounting for the impact of the trades on the price of the security. This problem is addressed in \cite{bertsimas98_MB}, using a stochastic dynamic program minimising the execution costs, and in \cite{almgren01_MB}, where a combination of volatility risk and transaction costs is minimised. Optimal execution and market impact are particularly significant issues in the context of algorithmic trading and have been addressed by the recently developed mean-field game theory, acknowledging the fact that price is impacted by the trades of many atomic players \citep{Firoozi17_MB, cardaliaguet18_MB, Mojtaba19_MB}.

\emph{Option games} appear in asset pricing models when a security gives interacting optional rights to more than one holder, that is, when the exercise of an optional right by one holder modifies those of the others. Examples include callable, putable and convertible bonds, warrants, and, especially, instruments subject to credit or counterparty risk. In general, the pricing of such financial instruments corresponds to the solution of a non-zero-sum stochastic game where players use feedback strategies \citep{ben2007_MB}.

\emph{Financial distress} models are used to price corporate debt, according to various assumptions about strategic default, debt service and bankruptcy procedures \citep{fan00_MB, broadie07_MB,annabi12_MB}.

Finally, a large literature in corporate finance uses game-theoretic models to deal with financial decisions made by firms, such as the choice between debt and equity when financing operations, the amount of dividends paid out to shareholders, and decisions about whether to invest in risky projects.

\subsubsection*{Further readings}
The recognition of finance as an thriving application area for OR methods developed about thirty years ago \citep[see, for instance, ][for an introduction to optimisation problems underlying risk management strategies and instruments]{dahl93I_MB,dahl93II_MB}. A review of practical applications of OR methods in finance appeared in \cite{board03_MB}. For a comprehensive textbook covering optimisation models in finance, we refer the reader to \cite{cornuejols06_MB}. A unified framework for asset pricing can be found in \cite{cochrane09_MB} and a review of applications of dynamic games in finance in \cite{breton18_2_MB}. Recent discussions about the interface of operations, risk management and finance as a promising research area are presented in \cite{wang21_MB} and \cite{babich21_MB}.

\subsection[Government and public sector (Katherine~Kent \& Sam~Rose)]{Government and public sector\protect\footnote{This subsection was written by Katherine~Kent and Sam~Rose.}}
\label{sec:Government_and_public_sector}

This subsection will present some OR applications within the UK’s government operational research service (GORS). GORS represents over 26 departments and agencies across Great Britain and Northern Ireland with analysts working in multi-disciplinary teams to find workable solutions to real life problems. The outbreak of the Coronavirus pandemic in 2020 introduced a new global backdrop and we were faced with the challenge of producing appropriate analysis to answers questions during an ever-changing landscape where time was of the essence. This led to collaborations across a wide range of departments across the nations.

A few examples of where this collaborative approach was adopted successfully are highlighted by the work carried out by the Department for Transport (DfT) and the Office for National Statistics (ONS). The ONS worked with other government departments such as Department of Health and Social Care (DHSC) and schools across the UK to monitor infection rates. They also applied their expertise in artificial intelligence (AI) in the form of semantic maps to gather insight into the pandemic. Additionally, the DfT along with other government departments used agent based modelling and discrete event simulation to unpick the issues around border disruptions and international travel. 

\subsubsection*{Coronavirus (COVID-19) Infection survey and Schools Infection Survey}
The Office for National Statistics (ONS) played a vital role during the pandemic in monitoring infection rates. The Coronavirus (COVID-19) infection survey estimates how many people across England, Wales, Northern Ireland, and Scotland would have tested positive for a COVID-19 infection, regardless of whether they report experiencing symptoms. This study was a collaboration with academic partners and funded by Department of Health and Social Care.
This major study involved asking people up and down the country to provide nose and throat swabs on a regular basis. These are analysed to see if they have contracted COVID-19. In addition, some adults are also asked to provide blood samples to determine what proportion of the population has antibodies to COVID-19. Further details of the methodology can be found in \cite{Office_for_National_Statistics2022-ey}.

Estimates of the total national proportion of the population testing positive for COVID-19 are weighted to be representative of the population that live in private-residential households in terms of age (grouped), sex and region. The analysis for the infection study is complex, the model generates estimated daily rates of people testing positive for COVID-19 controlling for age, sex, and region. This technique is known as dynamic Bayesian multi-level regression post-stratification (MRP). Details about the methodology are also provided by \cite{Pouwels2021-xu}.

Estimates from the ONS survey are published weekly, a critical element was how best to communicate the uncertainty, for dissemination estimates were translated into, for example, 1 in 50 people, with appropriate visuals including the ONS insights tool \citep{Office_for_National_Statistics2022-rx}. A complementary piece of work was monitoring transmission and antibody levels within schools, enabling the government to accurately assess the risk of different policy options.

The Schools Infection Survey (SIS) was a longitudinal study which collected data through polymerase chain reaction (PCR) tests, antibody tests and questionnaires. As well as monitoring transmission within a school environment, this data was used to assess the wider impacts of the pandemic and repeat lockdowns on our children and young people, including long covid, mental health and physical activity levels. Further detail can be seen in \cite{Office_for_National_Statistics2022-gy} and \cite{Hargreaves2022-yt}.

The Daily Contact Testing (DCT) trial was a blind medical trial which compared infection rates across two groups subject to different policies: the control group where children were in “bubbles”, and after one child testing positive the entire bubble would be sent home from school, and the intervention group where after a child tested positive, close contacts would then test daily and were allowed to remain in school as long as their results were negative. The study was a success and led to a policy change that resulted in schools being kept open for longer. Further detail can be seen in \cite{Young2021-ws}.

\subsubsection*{Semantic Maps and their use for understanding regional disparities}
Semantic Maps are a type of knowledge graph, championed in the world of robotics and Artificial Intelligence as a way to provide infrastructure to exploit all kinds of potentially even crowdsourced data and information in such a way as to provide dynamic, online, interactive visualisations that support the controlled and secure use of live data. The can be geospatial in nature, but they can also reflect connections through semantic relationships. These maps populated with data would provide users with many different ways to consume the underlying data and help inspire citizens about the potential power of data to drive understanding and generate insights. 

During the development of the Levelling Up White Paper the evidence base needed to be developed across government in order to define the key metrics and measures to focus policies in areas that would drive change. The white paper itself was delivered by the Levelling Up Taskforce in the Cabinet Office, however, ONS worked with the geospatial commission in convening a group of chief analysts from all departments on a regular basis. ONS and the Levelling Up Taskforce used the group to commission and collate existing evidence and then worked with officials in His Majesty’s Treasury (HMT) to develop a systems thinking model from that evidence. This systems mapping was the basis of the theory of change that underpins the white paper. Subsequently, the metrics and missions were developed and refined with this group in recognition of the fact that there needs to be a focus across the system to reduce disparities that are often larger within areas such as local authorities or regions than they are between them. 

ONS took this and developed a semantic map that identifies potential data sources for various aspects of the knowledge graph, using this to both prioritise filling evidence gaps where data and evidence do not currently exist and developing an integrated data asset for Levelling Up. This data is in the process of being acquired and engineered to be able to be easily linked through a set of linking ‘spines’ referred to as the reference management database. This engineering and architecture is key for supporting the sharing of information in a way that ensures privacy. It is intended that in future this asset will be made available across government and in a secure integrated data service. 

\subsubsection*{Agent Based Modelling}
The COVID-19 pandemic presented challenges for the international travel community. Government officials in transport and health needed to model the preventative effect of the various policy options including testing and isolation on importation of infections from international travel.  

The approach chosen in the Department for Transport to support this fast-moving policy area was agent based modelling which built upon the more scientific epidemiological modelling undertaken by colleagues in Health and academia. This allowed for the incorporation of the various differing parameters of international travellers including, where they were coming from and their risk of being infected and infectious, the uncertainty over incubation and infectious periods, and their likely behavioural response to various isolation and testing regulations. 

Whilst not designed to be a scientific forecast, the modelling allowed the cross-government community to estimate the relative effectiveness of policy options. This work supported policy making during a highly uncertain and changing environment when Government had to balance risk with the wider impacts on the aviation sector and the second order impacts on the economy.

\subsubsection*{Discrete Event Simulation}
As part of the EU Exit preparations, it was important for both the Department for Transport, the Home office and regional resilience teams to understand the impact of the expected border disruption at UK ports for roll-on roll-off freight traffic travelling to EU Member States. As the issue was around changes to the time taken and resource available to carry out additional border processes, the natural choice was discrete event simulation. Analysts in government developed a detailed model of the Short Strait crossings 
\citep[Port of Dover and the Channel Tunnel to France accounting for 84\% of accompanied heavy goods vehicles travelling to continental Europe in 2019;][]{Department_for_Transport2018-vw}. Regional models were developed to cover other ports. These allowed government officials to understand the likely queues and flow of vehicles and to understand the impact of changes to the system, which was vital to supporting contingency planning.

\subsubsection*{Statistical analysis and forecasting}
The COVID-19 lockdowns of 2020 accelerated the uptake of new and novel data sources for understanding mobility. Government analysts rapidly ingested new data sources such as that provided by Google mobility as open source data as well as procuring in additional anonymised and aggregated mobile network operator data. By analysing these new data sets alongside traditional demographic and geographic data sets it was possible to generate insights into the changes in mobility being seen across the country as a result of the various national and regional restrictions. Regression analysis was undertaken to produce a predictive model. It was then possible to forecast the impact of later changes to restrictions on population mobility.

\subsubsection*{Net zero - system thinking}
In 2020 the UK Prime Minister’s Council for Science and Technology advised on the following: ‘a whole systems approach can provide the framework that government requires to lead change across public and private sectors and \dots enables decision makers to understand the complex challenges posed by the net-zero target and devise solutions and innovations that are more likely to succeed’ \citep{GOVUK_2020-yk}. The Prime minister agreed.

As transport represents a huge portion of the challenge, ORs have run participatory systems mapping workshops in the Department for Transport with modal subject matter experts to identify the key causes and effects in the Transport Net Zero system. This aims to enable those working on Transport policy to explore the evidence, gain new insights and visibility of interdependencies within the system, and to understand the likely wider impact of their policy choices.

\subsubsection*{Conclusions}
All these examples help to illustrate the breadth of analysis undertaken across central government during the global pandemic to tackle real life issues; and the ranges of techniques we have as OR analysts to find workable solutions in an ever changing world.

\subsection[Healthcare (Christos~Vasilakis)]{Healthcare\protect\footnote{This subsection was written by Christos~Vasilakis.}}
\label{sec:Healthcare}

Why is the organisation and delivery of health and care services so difficult to manage, plan for and improve? Difficulties and delays in accessing care services, cancellations and increasing costs have a negative impact on all of us: patients, carers, and care professionals. Despite the attention and resources invested in addressing these problems, many health systems face increasing pressure to improve the effectiveness and efficiency of their operations. Part of the problem is the complexity inherent in the organisation of care services and our limited understanding of how changes will affect their delivery. Another problem is the intrinsic uncertainty and variability in many aspects of care service delivery. Add the multifaceted dynamics arising in this very complex socio-technical system involving professionals, patients and existing and new technologies, against a background of increased demand and budgetary constraints, and it is no surprise the effort to improve healthcare has been termed ‘rocket science’ \citep{Berwick2005-hf_CV}.

Operational research has a long established history in this area with the first application (scheduling outpatient hospital appointments) reported in the early 1950s \citep{Bailey1952-xt_CV}. Since then, there has been a proliferation of OR applications reported in the literature \citep{Katsaliaki2010-ae_CV}, and evidence of use to support policy making and care delivery \citep{Royston2009-ae_CV}. This is not a surprise given the importance of healthcare in our lives and that many of the problems faced by those managing and delivering care services are amenable to the methods and ethos of OR \citep{utley_crowe_pagel_2022_CV}. In the short review that follows, which is by no means exhaustive and draws primarily (but not exclusively) from the UK academic community and the National Health Service (NHS), I have attempted to give examples of review and individual studies grouped in a few broad areas of healthcare. 

\subsubsection*{Applications to hospital settings}
Hospital care has been the setting of a large number of OR studies \citep{Jun1999-oz_CV}. Hospitals typically survive reorganisations and funding cuts (unlike management, policy and other statutory bodies), and are large enough to be able to engage meaningfully in research projects (unlike, for example,  many primary care practices). Some (mainly teaching) hospitals host large biomedical research centres, with many of the professionals working in them active researchers.  

A specific area that has attracted the attention of operational researchers is the Emergency Department (ED). A recent review article identified 21 studies that used a computer simulation method to capture patient progression through the ED of an established UK NHS hospital, mainly focusing on service redesign \citep{Mohiuddin2017-zh_CV}. Individual studies have addressed the micro (single hospital) level \citep{Baboolal2012-ci_CV}, as well as the meso-level of emergency and on-demand healthcare within a region \citep{Brailsford2004-ib_CV}. The study by \cite{Lane2000-qz_CV} used System Dynamics to model the interaction of demand patterns and resources deployed in ED and other parts of the hospital to examine the link between emergency and elective operations in hospitals.

Another hospital area that has been the focal point of OR is peri-operative care. \cite{Sobolev2011-xv_CV} in their systematic review identified 34 studies modelling the flow of surgical patients. Various forms of optimisation have also been applied to surgical scheduling problems including operating room \citep{Fairley2019-eu_CV}, staffing \citep{Bandi2020-xy_CV} and nurse rostering \citep{Xiang2015-qe_CV} among others. \cite{Cardoen2010-rw_CV} identified almost 250 papers, with the rate of published studies accelerating at around the start of the new millennium (similar trends have been observed across many disciplines). The review revealed that most of the research was directed towards the planning and scheduling of elective patients in highly stylised scenarios – although many operational challenges are triggered by factors such as the arrival of non-elective (emergency) patients. More recently, the problems tackled have become more realistic to include considerations of downstream resource availability such as critical care and general ward beds \citep{Fugener2014-jk_CV}, and scheduling elective operations in such a way that randomly arriving emergency patients can be accommodated without excessive delays \citep{Jung2019-ls_CV}.

An area OR has demonstrably made a beneficial impact is the organisation of acute stroke services. Several studies have attempted to address the rate and speed with which patients with suspected acute ischemic stroke go through the initial diagnostic steps and receive treatment \citep{Meretoja2014-he_CV}. \cite{Monks2012-ec_CV} made a number of recommendations for improving treatment rates in a rural hospital. In a follow-up study that evaluated the results of their recommendations, mean door-to-needle times (a key performance metric with direct impact on patient survival and recovery) fell from 100 min to 55 min while thrombolysis rates increased to 14.5\% \citep{Monks2015-yy_CV}. More recently, the focus has shifted to supporting decision around the centralisation of regional acute stroke services \citep{Wood2020-jf_CV,Wood2022-ot_CV}, as well as the supporting the introduction of endovascular thrombectomy, a new and very effective treatment for ischemic stroke \citep{Maas2022-ty_CV}.

\subsubsection*{Applications to non-acute hospital care settings}
Much of healthcare is delivered outside of large hospital facilities. Primary care, home and community care, social care are significant components of the healthcare ecosystem. Primary care, whether provided by physicians, nurse practitioners or pharmacists, is typically concerned with providing a first contact and principal point of continuing care and/or coordinating other specialist care. There are early examples of theoretical work to assist primary care planning by estimating the coverage achieved by staff and facilities, using antenatal care as an example \citep{Kemball-Cook1983-yd_CV}. More recently in the area of maternity care provided in community as well as hospital facilities, \cite{Erdogan2019-ms_CV} developed and empirically tested an open source facility location solver to assist with a decision on the number and location of regional maternity facilities.

Home-based care has attracted considerable attention from operational researchers \citep{Grieco2021-le_CV}. This review identified studies proposing models and solution methods for operational decisions on staff rostering, the allocation of staff to patient visits, the scheduling of visits and the routing of staff. An example of impactful OR project is the Swedish study by \cite{eveborn_operations_2009_JHLS}, where a set of algorithms and accompanying software tool were developed to provide solutions to staff-to-patient allocations, staff scheduling and staff routing problems. Having deployed the tool to more than 200 units/organisations, operational efficiency was increased by up to 15\%, resulting in annual savings of 20-30 million euros. More recently, modelling work has supported the effort to address the timely discharging of hospital patients by using a combination of home-based and bedded ‘step-down’ community care \citep{Harper2021-qk_CV}.

Mental health, one of the leading causes of disease burden internationally, has also received the attention of operational researchers \citep{Long2018-ji_CV}. Specific areas of application varied from psychiatric ICUs \citep{Moss2022-ao_CV} to system design \citep{Smits2010-dq_CV} and planning \citep{Vasilakis2013-sn_CV}, medical decision making \citep{Haji_Ali_Afzali2012-jh_CV}, and epidemiology \citep{Ciampi2011-fi_CV}.

\subsubsection*{Public health, health system preparedness and resilience, and pandemic response}
Public health, the science and practice of helping people stay healthy and protecting them from threats to their health, is another area of OR applications. The review article by \cite{Fone2003-uq_CV}, identified OR studies of infection and communicable disease, screening, and several epidemiological and health policy studies. Microsimulation, a type of simulation which models individual life trajectories through a number of healthy and disease states, has found wide applicability in the area of public health \citep{Krijkamp2018-iz_CV}, such as forecasting the long-term care needs of the older population in England \citep{Kingston2018-fi_CV}. Multicriteria decision analysis (MCDA) methods have also been used extensively to address questions of health policy or health technology assessment \citep{Glaize2019-jz_CV}. 

An area that has seen increased attention over the last two decades is that of emergency preparedness and health system resilience \citep{Tippong2022-ec_CV}. Emergency preparedness studies include a study of red blood cell provision following mass casualty events \citep{Glasgow2018-th_CV}. Examples of health system resilience studies include the paper by \cite{Crowe2014-gr_CV}, which examined the feasibility of using modelling to assess the capacity of a care system to continue operating in the face of major disruption. 

The COVID-19 pandemic not only gave rise to a large number of modelling studies, it also raised the profile of mathematical modelling with the general public and policy makers. \cite{Pagel2022-vg_CV}, in their excellent article on the role of modelling in the pandemic response, discussed the early lessons learnt from this experience including the poor understanding of policy makers and the public of key concepts such as exponential growth. They argue that infection disease modelling, which generated much of the evidence used to support decisions of pandemic response \citep{Brooks-Pollock2021-lw_CV}, is intrinsically difficult given the complex relationships between the model parameters, and the difficulties associated with quantifying these parameters. 

The possible benefits of modelling in addressing the challenges presented by the pandemic were outlined by \cite{Currie2020-tz_CV}. Indeed, several studies emerged early in the pandemic including, for example, an attempt to forecast the number of infected and recovered cases used univariate time series models \citep{Petropoulos2020-qm_CV}. \cite{Wood2020-wx_CV} published one of the first OR studies that examined the likely impact of increases in critical care capacity as a means to reduce the COVID-19 death toll. In a follow-up study, the sophistication of the model was increased to capture notions of triaging access of patients to critical care beds during periods of intense demand \citep{Wood2021-jv_CV}. The operation of large vaccination centres was also the topic of several modelling studies, both theoretical \citep{Franco2022-dx_CV} and empirical \citep{Wood2021-jt_CV,Valladares2022-yk_CV}.

\subsubsection*{Concluding remarks}
Despite the large body of literature, the role and impact of OR on improving care systems is less clear. Hospitals have “largely failed to use one of the most potent methods currently available for improving the performance of complex organisations” \citep{Buhaug2002-ck_CV} and “staff may be largely unaware of the potential applications and benefits of OR” \citep{utley_crowe_pagel_2022_CV}. A systematic review found that only half of the included studies reported models that were constructed to address the needs of policy-makers, and only a quarter reported some involvement of stakeholders \citep{Sobolev2011-xv_CV}. Recent positive developments include the introduction of guidelines to improve the reporting of OR studies \citep[e.g.,][]{Monks2019-wt_CV}, studies that recognise the importance of behavioural factors in attempts to influence practice and decision making with OR \citep{Crowe2022-bc_CV} and attempts to systematically generate evidence on the value and impact of OR on patient and system outcomes \citep{Monks2015-yy_CV,Soorapanth2022-jt_CV}. The research agenda should continue to evolve with the aim of addressing the challenges around engagement, implementation and evidencing the impact of OR applied to healthcare problems.

\subsection[Inventory (Jing-Sheng~Song)]{Inventory\protect\footnote{This subsection was written by Jing-Sheng~Song.}}
\label{sec:Inventory}

Inventories are the materials, parts, and finished goods held by an organisation for future use or sale. Not having enough inventory is costly. Shortages of materials and parts cause interruptions in production processes, delays in product delivery, and stockouts of finished goods. On the other hand, carrying inventory is costly, too, involving the cost of capital due to tied-up capital, storage cost, insurance, taxes, and spoilage and obsolescence costs.  

Inventory theory studies analytical models and solution techniques to help organisations meet the service requirement most cost-effectively or minimise the total expected costs of ordering, inventory-holding, and shortage. It does so by quantifying the tradeoffs driven by economies of scale, lead time (the time it takes to receive the ordered quantity after placing an order), and supply and demand uncertainties. It prescribes effective inventory-control policies that govern when to order an item (called \textit{reorder point}) and how much to order (called \textit{order quantity}). 

Inventory models distinguish from each other along several features: single or multiple planning periods, discrete- or continuous-time inventory monitoring, single- or multi-product, single- or multi-stage (or location), demand nature (deterministic or stochastic, stationary or nonstationary, distribution known or unknown), product perishability, lost sales or backlogging when shortages occur, deterministic or stochastic lead time, supply system (single- or dual-source, exogenous or endogenous, a finite or infinite capacity), and cost structure (with or without a fixed ordering cost, etc.). The following research-based textbooks and handbooks offer more detailed coverage and references:  \cite{arrow1958studies_JSS}, \cite{axsater2006inventory_JSS}, \cite{de2003supply_JSS}, \cite{graves1993logistics_JSS},
\cite{hadley1963within_JSS}, \cite{nahmias2011perishable_JSS}, \cite{porteus2002foundations_JSS}, \cite{silver1988inventory_JSS}, \cite{simchi2014logic_JSS}, \cite{snyder2019fundamentals_JSS}, \cite{song2023handbook_JSS}, and \cite{Zipkin2000_JSS}.

One class of models focuses on characterising the optimal inventory-control policy under a given supply and demand environment and cost structure. A common approach is formulating a multi-period inventory decision problem as a dynamic program and transforming the original formulation into a simpler one through state reduction. Next, identify the structural properties of the single-period cost function to determine the optimal policy form for a single-period problem. Then, show that these properties are preserved by the (Bellman) optimality equation, so the policy form is optimal for each period. The optimal policy parameters may not be easy to compute; hence some works develop efficient algorithms to calculate the optimal policy parameters. 

Another class of models focuses on developing efficient performance evaluation tools for a given type of inventory policy that is either commonly used in practice or of simple structure and easy to implement. This is particularly important for systems where state reduction is not viable and the dimension of the system state grows exponentially in the number of periods (the so-called curse of dimensionality), so the optimal policy has no simple form. Typically, this type of work analyses a continuous-review system in which demand follows a stochastic process and derives steady-state performance measures of any given policy, such as average inventory, average backorders, and stockout rate, as well as the long-run average cost. Then, optimisation tools can be developed to find the optimal policy parameters that minimise the long-run average cost. 

The third class of models conducts asymptotic analysis to establish asymptotic optimality of some simple-structured policies for less tractable inventory systems with unknown and complex optimal policies. 

The following are several classic models where the optimal policies are shown to have simple forms. Unless otherwise stated, the models assume a single stage, a single source, and a single nonperishable product. 

\textit{The EOQ (Economic Order Quantity) Model} was first developed by \cite{harris1913many_JSS} \citep[see the reprint][]{harris1990many_JSS} and popularised by \cite{wilson1934scientific_JSS}. It concerns the balancing of holding and ordering costs due to economies of scale in procurement or production. It is a continuous-review model over an infinite planning horizon, assuming the annual demand for the stocked item is a constant $\lambda$. There is a fixed procurement cost $k$ independent of the order size, accounting for administrative, material handling, and transportation-related costs. The annual per-unit inventory-holding cost is $h$. The optimal order quantity (EOQ) that minimises the annual order and holding costs equals $\sqrt{2k\lambda/h}$, which is insensitive to small perturbations of the model parameters. Variations of this model can accommodate finite production rate, planned backlogs, random yield, a quantity discount, and time-varying demand \citep[also known as the dynamic lot sizing problem; see][]{wagner1958dynamic_JSS,silver1973heuristic_JSS}. It also forms the basis for the development of efficient multi-item joint replenishment policies and multiechelon coordinated replenishment policies such as the power-of-two policies; see \cite{roundy198598_JSS}, \cite{roundy198698_JSS}. For major developments and references, see \cite{axsater2006inventory_JSS}, \cite{muckstadt1993analysis_JSS}, \cite{silver1988inventory_JSS}, \cite{simchi2014logic_JSS}, and \cite{Zipkin2000_JSS}.

\textit{The Newsvendor Model}, which originated from \cite{edgeworth1888mathematical_JSS} in a banking application, was formalised by \cite{arrow1951optimal_JSS} in the general inventory context. It optimises the tradeoff between too much and too little inventory caused by demand uncertainty for a seasonal product. It is a single-period model with only one ordering opportunity before the selling season, assuming an estimated demand distribution. The fixed order cost is negligible. After the ordered quantity arrives, the selling season begins, and demand realises. At the end of the season, there will be either unsold units (overage) or unmet demand (underage). The unit overage cost ($o$) = purchasing cost less - salvage value, while the unit underage cost ($u$) is the lost profit. The optimal \textit{newsvendor order quantity} equals to the fractile of the demand distribution at the critical ratio $u/(u+o)$. The model can be generalised in many ways, including random yield, different cost structures, pricing, and distribution-free bounds \citep{gallego1993distribution_JSS, petruzzi1999pricing_JSS, porteus1990stochastic_JSS, qin2011newsvendor_JSS} and multi-location with risk-pooling effect \citep{bimpikis2016inventory_JSS,eppen1979note_JSS}.

\textit{Dynamic Backlogging Models}. The most tractable and developed setting for multi-period models with stochastic demand and a constant lead time is full backlogging. When stockouts are rare, this model is a reasonable approximation for the lost-sales system. An important concept (due to state reduction) is \textit{inventory position}, which is the sum of the on-hand inventory plus total pipeline inventory minus backorders. This is the total system inventory available to satisfy future demand if we do not order again. 

Assume demand is independent over time. A \textit{base-stock policy} is optimal if the order cost is linear (no fixed order cost). Each period has a target inventory position called the \textit{base-stock level}. If the inventory position before ordering is below this level, order up to this level; otherwise, do not order. If the demand is stationary, the myopic base-stock level that minimises a single-period expected cost is optimal. The base-stock level has the same form as the newsvendor quantity, with the holding cost as the overage cost, the backorder cost as the underage cost, and the demand during a lead time replacing the single-period demand. For nonstationary demand, as long as the myopic base-stock levels are nondecreasing in time, the myopic base-stock level is still optimal. See \cite{veinott1965optimal_JSS} and \cite{porteus1990stochastic_JSS}. 

When the order cost is linear plus a fixed cost $k$, the optimal policy is an \textit{$(s, S)$ policy}. In each period, if the inventory position before ordering is below a threshold $s$, order up to $S$; otherwise, do not order. The key enabler of this result is that the single-period cost is $k$-convex, a property discovered by \cite{scarf1960optimality_JSS}. When the demand is stationary, the policy is also stationary. In a continuous-review system with Poisson demand, the optimal policy is an \textit{$(r,q)$ policy}: When the inventory position reaches $r$, order $q$ units. It is equivalent to the $(s,S)$ policy with $r=s$ and $q=S-s$. A simple yet effective heuristic policy is to use the optimal base-stock level to approximate $r$, and use the EOQ formula to approximate $q$; see \cite{zheng1992properties_JSS} and \cite{axsater1996using_JSS}.

These policy structures have been extended to more complex models, such as Markov modulated demand \citep{iglehart1960optimal_JSS, song1993inventory_JSS, sethi1997optimality_JSS}, exogenous and sequential stochastic lead times \citep{kaplan1970dynamic_JSS, nahmias1979simple_JSS, ehrhardt1984s_JSS, song1994effect_JSS, song1996inventory_JSS}, capacity constraints \citep{federgruen1986ainventory_JSS, federgruen1986binventory_JSS}, unknown demand distribution \citep{scarf1959bayes_JSS, scarf1960some_JSS,  azoury1985bayes_JSS}, and a dual-source problem where the lead times of the two sources differ by one period \citep{fukuda1964optimal_JSS} or the lead times are stochastic and endogenous \citep{song2017optimal_JSS}. See \cite{veinott1966status_JSS}, \cite{perera2022asurvey_JSS}, \cite{perera2022bsurvey_JSS}, \cite{porteus1990stochastic_JSS}, and \cite{Zipkin2000_JSS} for more detail. 

Multiechelon (or multi-stage) inventory systems are common in supply chains where the stages are interrelated, such as production facilities, warehouses, and retail locations. The literature focuses on understanding three basic system structures: series, assembly, and distribution systems. 

In a series system with $N$ stages and backlogging, random customer demand arises at stage 1, stage 1 orders from stage 2, and so on, and stage N orders from an outside supplier with ample supply. There is a constant transportation time between two consecutive stages. Define the echelon inventory of each stage to be the inventory at the stage plus all downstream inventories (including those in transit). Assuming no fixed order costs, \cite{clark1960optimal_JSS} establish that an echelon base-stock policy is optimal for all stages. That is, we can treat each echelon as a single location and order the echelon inventory position up to a target base-stock level. \cite{axsater1993installation_JSS} show that for any echelon base-stock policy, there is an equivalent local base-stock policy; therefore, the implementation of the optimal policy is simple. \cite{federgruen1984computational_JSS} find that the optimal echelon base-stock policy for the infinite horizon problem can be efficiently obtained. \cite{rosling1989optimal_JSS} proves that under certain mild conditions, an assembly system can be transformed into an equivalent series system, so the Clark-Scarf result applies. \cite{chen1994lower_JSS} further streamline the proofs of these results. \cite{shang2003newsvendor_JSS} construct effective single-stage newsvendor solutions to approximate the optimal echelon base-stock levels. \cite{chen2001optimal_JSS} show that a state-dependent echelon base-stock policy is optimal for Markov-modulated demand. See \cite{axsater1993continuous_JSS}, \cite{axsater2003supply_JSS}, \cite{axsater2006inventory_JSS}, \cite{angelus2023generalization_JSS}, \cite{federgruen1993centralized_JSS}, \cite{kapuscinski2023capacitated_JSS}, and \cite{shang2023single_JSS} for more developments, including batch ordering, capacity limits, distribution systems, transshipment, and expediting. 

Many other features are much less tractable, such as lost-sales systems \citep{bijvank2023lost_JSS}, censored demand data, perishable products \citep{li2023perishable_JSS}, general dual-sourcing systems \citep{xin2023dual_JSS}, distribution systems, and assemble-to-order systems \citep{atan2017assemble_JSS,song2003supply_JSS, de2023assemble_JSS}. Nonetheless, significant progress has been made in recent years on structural properties of the optimal policy, asymptotic optimal policies, and effective heuristics, thanks to more analytical tools such as discrete convexity, asymptotic analysis, and machine learning algorithms. See \cite{chao2023online_JSS}, \cite{cheung2023statistical_JSS}, \cite{shi2023approximation_JSS}, and other chapters in \cite{song2023handbook_JSS}.

\subsection[Location (Sibel~A.~Alumur)]{Location\protect\footnote{This subsection was written by Sibel~A.~Alumur.}}
\label{sec:Location}

In the domain of operations research, location problems are concerned with determining the location of a facility or multiple facilities to optimise one or more objective functions under constraints. 
Location problems seek answers to questions such as how many facilities should be located, where should each facility location be, how large should each facility be, and how should the demand for the facilities’ services be allocated to these facilities \citep{daskin1995network_SAA}. 
An example of a facility to be located is a factory, distribution centre, warehouse, cross-dock, or hub, where demand can be for raw materials, components, products, passengers, data, etc. 

Location decisions arise in a variety of public and private sector decision-making problems. Some examples from different sectors include locating landfills where demand is for disposal of household waste \citep{ERKUT1989275_SAA},  ambulances where demand is for transporting emergency patients to hospitals \citep{brotcorne2003ambulance_SAA}, warehouses where demand is for storing products arriving from factories \citep{aghezzaf2005capacity_SAA}, schools where demand is for students \citep{haase2013management_SAA}, regenerators in optical networks where demand is for data \citep{yildiz2017regenerator_SAA}, shelter sites where demand is for refugees \citep{bayram2018shelter_SAA}, and charging stations where demand is for electric vehicles that need to charge \citep{kinay2021full_SAA}. More applications of location problems from practice can be found in \cite{eiselt2015applications_SAA}.

Location decisions refer to the placement of a facility considering its interactions with demand points (e.g., customers, suppliers, retailers, households) and possibly with other facilities to be located. It includes selecting the location and determining how this location supports meeting a decision-maker or organisation's objective. It is important to note that facility location decisions are different from facility \textit{design} decisions. Facility design decisions usually consist of facility layout and material handling systems design. The layout entails all equipment, machinery, and furnishing within the building, whereas material handling systems comprise the mechanism needed to satisfy the required facility interactions. Facilities planning and design are extensively discussed in \cite{tompkins2010facilities_SAA}.

Several factors influence facility location decisions, the most prominent ones being transportation costs and the availability of the transportation infrastructure. Among other important factors are the availabilities and costs of land, market, labour, materials, equipment, energy, government incentives, and competitors as well as geographical factors and weather conditions. 
Distance is usually considered to be one of the most important criteria in facility location models. Several distance metrics can be used in location models such as Euclidean (straight-line), rectilinear (Manhattan), Cheybyshev (Tchebychev), and network distance. Network distance is the distance that is calculated on an existing transportation network, for example, through using Google or Bing maps. 

An important criterion to be considered in location problems is how demands are to be satisfied by the facilities to be located. In some applications, the whole demand of each customer must be satisfied from a single facility (``non-divisible" demand) which is referred to as a \textit{single source} or \textit{single allocation}. Single allocation location problems are also referred to as \textit{location-allocation} problems as each demand point is allocated to a single facility. On the other hand, in \textit{multiple source}/\textit{allocation} problems, the demand of a single customer can be served from several facilities.  

Location decisions are usually classified according to their decision space. In \textit{continuous} or \textit{planar} location problems, the facilities can be located anywhere in the decision space. The search is for the optimal coordinates; i.e, latitude and longitude. In \textit{discrete} location problems, a finite set of potential facility locations is provided, possibly determined through a pre-selection process. In \textit{network} location problems, on the other hand, there is a given network and the facilities are to be located on this network. In network location problems, facilities can further be restricted to be placed only on the vertices or nodes of this network and not on the edges or arcs, referred to as \textit{vertex-} or \textit{node-restricted} location problems.

Continuous location problems focus on minimising some function related to the distance between the facilities to be located and the existing facilities or demand points, such as suppliers and customers, where minisum (minimising the total weighted distance) and minimax (minimising the maximum or worst weighted distance) are among the most commonly employed objectives. Special cases of continuous single-facility location problems with commonly used distance metrics (e.g., rectilinear and Euclidean) are well-studied and polynomial time solution algorithms exist \citep{francis2004facility_SAA}. In the case of multi-facility continuous location problems, the facilities to be located can be homogeneous or non-homogeneous; in the latter, there are different types (e.g., a factory and a warehouse) or sizes of facilities to locate. 

One of the most studied discrete location problems is the \textit{$p$-median} problem. The goal is to pick a subset $p$ of (homogeneous) facilities to open from among a given set of potential locations that minimise the total transportation cost of satisfying each demand point from the (nearest) facility it takes service from. There is a well-known node optimality theorem by \cite{hakimi1965optimum_SAA} for the $p$-median problem on networks that proves that at least one optimal solution to the $p$-median problem consists of locating the facilities only on the nodes of the network (even though a facility is allowed to be located anywhere on the network including any point on an edge between the nodes). Possible applications of the $p$-median problem are clustering, transit network timetabling and scheduling, placement of cache proxies in a computer network, diversity management, cell formation and much more \citep{marin2019p_SAA}. 

An important related problem is the \textit{uncapacitated facility location problem} (UFLP) which is also referred to as the \textit{simple plant location problem}. Unlike the $p$-median problem, in  UFLP, the number of facilities to be located is no longer known and determined by optimising an objective function that considers the trade-off between the fixed costs of locating facilities and the transportation costs. Numerous extensions of UFLP with uncertainty, multiple commodities (e.g., products or services), multi-period planning horizon, multiple objectives, and network design decisions have been studied with applications in several domains such as supply chain and distribution systems design. 

A nice structure of $p$-median and UFLP is that since the facilities to be located are assumed to have enough capacity (e.g., space or labour hour), all demands of each customer can be served from a single facility with minimum allocation costs. This is no longer the case for \textit{capacitated} versions of the facility location problems, where single- and multiple-allocation versions are both extensively studied. For multiple allocation capacitated (fixed-charge) facility location problems, when the set of open facilities is given, the resulting subproblem of finding the best allocations is a transportation problem. In the single allocation case, on the other hand, when the set of open facilities is pre-determined the resulting allocation subproblem is a generalised assignment problem \citep{fernandez2015fixed_SAA}.  

When the worst-case is more important than the average, it might be better to consider the furthest or most disadvantaged demand point to ensure equity in servicing the demand. Accordingly, the $p$-centre problem aims to locate $p$ facilities such that the maximum distance (or travel time/cost) from a demand point to its nearest facility is minimised (minmax). The $p$-centre problem can be used to locate public schools and various emergency service facilities such as police stations, hospitals, and fire stations. Different variations of this problem have been studied such as the capacitated, conditional, continuous, fault-tolerant, and probabilistic $p$-centre problems \citep{ccalik2019p_SAA}. 

In covering location problems, the aim is to locate facilities so as to cover demand. Typically, a demand point is considered to be covered if it is within a certain distance or travel time of a facility. Unlike in the previous models, the demand points are not assigned to facilities in covering location problems. The two most common covering location problems are \textit{set covering} and \textit{maximal covering} location problems. In the set covering location problem, the aim is to minimise the total cost of locating facilities to cover all demand points, whereas, in the maximal covering location problem, the aim is to maximise the total demand covered subject to a budget constraint or a constraint on the total number of facilities to locate. Continuous variants of these covering location problems are also studied \citep{plastria2002continuous_SAA}. Several different versions of covering location problems have been studied in the literature, including, but not limited to weighted, redundant, hierarchical, backup covering problems with applications in emergency services, crew scheduling, mail advertising, archaeology, metallurgy, and nature reserve selections \citep{garcia2015covering_SAA}. 

In general, facility location problems consider and model only a single echelon; i.e., either the flows of commodities  (e.g., products, customers) coming into or out from the facilities to be located are negligible, for instance, when one of those transportation costs is borne by another decision-maker and somehow not related to the current decision-making problem. An example would be a manufacturing company determining the location of its new factory for delivering products to its customers with minimum total cost, where the company is not directly involved in the delivery of raw materials from their suppliers to the factory. When the flow of commodities coming into the facilities to be located as well as the flow going out of those facilities are simultaneously considered in the models, these location problems are referred to as \textit{two-echelon} location problems. For example, while locating a distribution centre, the transportation cost of products from the factory to this distribution centre as well as the transportation costs from the distribution centre to the retailers may need to be considered in the model. Sometimes there are facilities to be located at several echelons where flows of commodities in and out of all those facilities need to be considered. These \textit{multi-echelon} types of location problems are encountered for several applications of supply chain network design \citep{melo2009facility_SAA}. Another related category is when there is a hierarchical network structure among the facilities to be located, referred to as \textit{hierarchical} facility location problems \citep{csahin2007review_SAA}. An example of a hierarchical location problem is designing a postal delivery network where the locations of the sorting centres as well the locations of the post offices that are to be allocated to those sorting centres need to be determined. 

There might also be interactions among the facilities to be located. This is the case, for example, for \textit{hub location} problems where the demand is defined between pairs of demand points (origin-destination pairs) as opposed to having the demand of an individual point. 
In that case, to satisfy the demand from an origin to a destination point, flow can be transported between the facilities to be located en route to the destination, where those facilities can act as switching, transshipment, sorting, connection, consolidation, or break-bulk points. Hub location models have several applications in passenger and freight airlines, express shipment, postal delivery, trucking, public transit, and telecommunication network design \citep{alumur2021perspectives_SAA}.

Location problems have been a testbed for many algorithmic and methodological advances in operations research. Most discrete location problems commonly belong to a class of $\mathcal{NP}$-hard decision problems (\S\ref{sec:Computational_complexity}) and they can usually be formulated with mixed-integer programming (MIP) models (see \S\ref{sec:Mixed_integer_programming}). In addition to using commercial MIP solvers, several exact and (meta)heuristic algorithms (\S\ref{sec:Heuristics}) have been developed and tested on benchmark instances from the literature. Some of those benchmark instances can be obtained from \cite{beasley1990or_SAA}, \cite{posta2014exact_SAA}, and \cite{fischetti2017redesigning_SAA}.

Location science is a very broad field of research that encompasses geography, continuous and discrete optimisation (\S\ref{sec:Combinatorial_optimisation}), graph theory (\S\ref{sec:Graphs_and_networks}), logistics (\S\ref{sec:Logistics}), and supply chain management (\S\ref{sec:Supply_chain_management}). This section only highlights the basic and most well-known location models. For a more detailed overview of the field of location science, we refer the reader to several books written in this field, such as \cite{drezner2004facility_SAA}, \cite{eiselt2011foundations_SAA}, and \cite{LNS15_SMPT}. 

\subsection[Logistics (Janny~Leung \& Yong-Hong~Kuo)]{Logistics\protect\footnote{This subsection was written by Janny~Leung and Yong-Hong~Kuo.}}
\label{sec:Logistics}

Logistics refers to the organisation and implementation of the processes related to the procurement, transport and maintenance of materials, personnel and facilities. The application of operational research to logistics dates back to 1930 \citep{schrijver2002history_JLYHK}, where \cite{tolstoi1930metody_JLYHK} solved to optimality the problem of transporting salt, cement, and other cargo on the railway network of the Soviet Union. In general, the objective of logistics management can be summed up as ``getting the right thing/people to the right place at the right time in the right quantity at the right cost''. For materials, logistics operations require the co-ordination of forecasting, purchasing, inventory control, warehousing, distribution, transportation, delivery and installation. Logistics management of personnel involves, in addition, skills matching,   capabilities training, labour rules and worker preferences. At a strategic level, logistics involves the design of the transport network and facilities. In this subsection, we discuss several major domains of logistics applications, namely, military, inventory, time-sensitive, reverse and humanitarian logistics. We also mention some new technologies for emerging logistics applications.

\subsubsection*{Military Logistics}

Logistics play an important role in military operations. Indeed, the word ``logistics" itself is derived from the position {\it{Mar\'echal des logis}} created in the French army in the 17\textsuperscript{th} century, whose responsibilities of establishing camps and arranging transport/supplies were referred as ``la logistique" \citep{Jomini1862_JLYHK}. Many historians credited logistics as the success factor in wars from ancient to modern times. During World War II, the need for large-scale logistics planning accelerated the development of operational research. The ability to sustain the convoy of supply ships was a major factor in the Battle of the Atlantic \citep{Kirby2003_JLYHK}. The 1948-1949 Berlin Airlift, where over 2.3 million tons of goods were flown to besieged West Berlin, is well-known as the first use of logistics as a military and political strategy \citep{Tine2005_JLYHK}.

The North Atlantic Treaty Organisation defines logistics as the science of planning and carrying out the movement and maintenance of forces, and covers acquisition, transport, maintenance and evacuation of materiel, personnel and facilities, and provision of services and medical support. Operational research methodologies are extensively used \citep{Scala2020_JLYHK}. Reliability and operability of the supply lines are a major concern in military logistics \citep{McConnell2021_JLYHK}, and simulation is much utilised. \cite{Cioppa2004_JLYHK} review agent-based simulation for military applications. Emerging technologies -- such as additive manufacturing \citep{Boer2020_JLYHK} and unmanned transport \citep{Jotrao2021_JLYHK} -- have also sparked research in smart military logistics \citep{Schutz2017_JLYHK}. The reader is also referred to \S\ref{sec:Military_and_Homeland_security}.

\subsubsection*{Inventory Logistics}

In modern logistics, most activities are related to products and goods, where their availability to customers or users is a key concern. Inventory, thus, plays an important role in this respect. A classic problem related to inventory logistics is the inventory-routing problem (IRP), introduced by \cite{bell1983improving_JLYHK} for the distribution of industrial gases. Since then, various IRP applications have been studied, including those related to automobile components \citep{blumenfeld1985analyzing_JLYHK}, groceries \citep{gaur2004periodic_JLYHK}, cement \citep{christiansen2011maritime_JLYHK}. Typically, IRP arises in vendor-managed inventory systems as the supplier monitors the inventory and makes replenishment decisions for its retailers \citep{archetti2007branch_JLYHK}. Because inventory can be carried from one period to the next, IRP considers {\it joint} decisions of inventory and routing across multiple periods and aims to minimise the total transportation and inventory holding costs over the planning horizon, subject to all demands being satisfied. \cite{speranza1994minimizing_JLYHK} extended the IRP to settings with multiple products. When demands are uncertain, IRP becomes stochastic IRP \citep{federgruen1984combined_JLYHK,trudeau1992stochastic_JLYHK}, where the objective function includes additional shortage cost. \cite{coelho2014heuristics_JLYHK} investigated the stochastic dynamic IRP where decisions are made as customers' demand become realised. Inventory logistics is even more timely today due to e-commerce \citep{archetti2021recent_JLYHK}. The main challenge of these inventory logistics problems is due to the computational complexity of solving multiple $\mathcal{NP}$-hard problems simultaneously. The reader is also referred to \S\ref{sec:Inventory}.

\subsubsection*{Time-sensitive Logistics}

The quality and functionality of items, in storage or transit, deteriorate over time. For items such as fresh vegetables, the value degrades continuously. Other perishables, such as blood, have fixed lifetimes and cannot be used beyond expiry. Logistics management of time-sensitive goods must consider production, distribution and transport jointly. \cite{Federgruen1986_JLYHK} was one of the first papers to consider jointly inventory allocation and transportation for fixed-lifetime perishables with probabilistic demand. Since then, there has been much research exploring additional issues, such as freight consolidation \citep{Hu2018_JLYHK}, storage/transport capacities \citep{Crama2022_JLYHK} and environmental concerns \citep{GOVINDAN2014_JLYHK}. \cite{Shaabani2022_JLYHK} gives a comprehensive literature review.

For continually decaying food items, delivery costs must be traded off with freshness-upon-arrival which may lead to lost sales or revenue \citep{Mirzaei2015_JLYHK}. The overall network design -- especially decisions on where along the supply chain processing occurs -- is important, since deterioration rates differ for unprocessed vs.~finished/packaged goods, and for items in transport vs.~in storage \citep{DeKeizer2017_JLYHK}.

An important category of perishable goods is blood. Integer-programming models were developed by \cite{Hemmel2009_JLYHK} for collection and distribution of blood products to Austrian hospitals, and by \cite{Araujo2020_JLYHK} for blood delivery in south Portugal. \cite{Piraban2019_JLYHK} survey research on blood supply chain management. 

\subsubsection*{Reverse Logistics}

Due to increased sustainability awareness and legislation, reducing the environmental impact of production and distribution has become important. Twenty years ago, \cite{Beamon1999_JLYHK} advocated that supply chains must be extended from one-way to a closed loop, where used products and materials are recovered for {\em re-use, recycle or re-manufacture}. Reverse logistics, thus, refer to the material flow from the point of consumption back upstream for regenerating value \citep{Rogers2001_JLYHK}. Compared to a forward supply chain, reverse logistics processes are more complicated. Firstly, the source, quality and quantity of recoverable used products/materials from end-users are highly unpredictable. There is an added decision-stage for inspection, evaluation and sorting of the collected materials, and streaming them into various processes (re-use, re-manufacture, disassembly, disposal, etc.). These re-purposing processes may be expensive, so trade-offs must be made between recovery cost and salvage value.

Re-manufacturing, where items are repaired to serviceable (like-new) condition, is an important aspect of reverse logistics. \cite{Simpson1978_JLYHK} was the first to address a multi-period repairable inventory problem with random demand and returns supply; using dynamic programming, he found the optimal policy structure which specified the repair, purchase and scrap levels for each period. Later, the model was extended to consider side-sales \citep{Calmon2017-gn_AM} and warranty demands \citep{Lin2020_JLYHK}.  Nowadays, the concept of reverse logistics is broadened holistically to closed-loop supply chains and the circular economy \citep{SG2019_JLYHK}. See \cite{Van_Engeland2020-qq_JL} for a recent review.

\subsubsection*{Humanitarian Logistics}

When disasters strike, speedy evacuation and prompt delivery of resources to affected areas are critical.
From some sparse early studies \citep{Sherali1991_JLYHK}, humanitarian logistics research grew rapidly since 2000. The research stream yielded insights that have changed how humanitarian agencies plan and manage disaster relief. A key concept is {\em inventory pre-positioning} where depots are set up already stocked with supplies {\em in anticipation} of disaster occurrences, instead of scrambling for procurement in the aftermath. \cite{Duran2011_JLYHK} developed a facility-location and supply pre-positioning plan for CARE. See also \cite{Rawls2011_JLYHK}. Many of the models used are large-scale mixed-integer-programming models.

Humanitarian logistics involve multiple objectives: costs, response urgency and fairness are all important.
\cite{Huang2012_JLYHK} considered {\em equity} in last-mile distribution;  \cite{Sheu2014_JLYHK} incorporated perceptions of people awaiting rescue. Other researchers considered decision under uncertainty:  \cite{Mete2010_JLYHK} developed a stochastic model for location and delivery of medical supplies. Yet other research took an interdiction approach and anticipated post-disaster deployment \citep{OHanley2011_JLYHK, Irohara2013_JLYHK}. Recent technological advances have stimulated new research and practices. \cite{Maharjan2020_JLYHK} investigated pre-positioning of mobile logistics/ telecommunications hubs for Nepal. See \cite{Behl2019_JLYHK} for a survey. The reader is also referred to \S\ref{sec:Disaster_relief_and_humanitarian_logistics}.

\subsubsection*{Emerging Technologies}

As technologies advance, the role of logistics has become more important in the Industry 4.0 era \citep{tang2019strategic_JLYHK}. Tracking and locating technologies (RFID, GPS, IoT, etc.) enable organisations and companies to acquire information in real time. Powerful computing facilities can perform analytics of massive volumes of historical data to support near real-time solutions for large-scale problems -- essential for city logistics involving thousands of orders to fulfil within a day, or even an hour. Exciting emerging applications include TSP/VRP for routing of drones \citep{masmoudi2022vehicle_JLYHK} and/or autonomous vehicles \citep{reed2022impact_JLYHK}, risk analysis powered by block-chain technology \citep{choi2019mean_JLYHK}, flow-based optimisation for crowd-sourcing logistics \citep{sampaio2020delivery_JLYHK} and cargo hitching \citep{fatnassi2015planning_JLYHK}, demand-driven optimisation for car/bike sharing systems \citep{wang2022dynamic_JLYHK}, and queuing and simulation for robotic warehouses \citep{fragapane2021planning_JLYHK}.

These more complex and larger-scale problems with tighter response times require new solution methodologies. Most of these new approaches are combinations of operational research and data science techniques, for example, robust optimisation \citep{zhang2021robust_JLYHK}, reinforcement learning \citep{yan2022reinforcement_JLYHK} and other machine learning-based optimisation approaches \citep{bengio2021machine_JLYHK}.

\subsection[Manufacturing (Kathryn~E.~Stecke \& Xuying~Zhao)]{Manufacturing\protect\footnote{This subsection was written by Kathryn~E.~Stecke and Xuying~Zhao.}}
\label{sec:Manufacturing}

Manufacturing is the production process from materials to goods. Such goods can be finished goods sold to end consumers or components sold to other manufacturers for the production of other more complex products. Manufacturing has gone through several different phases (Industry 2.0 to 4.0) in the twentieth and twenty-first centuries. Here we offer an overview of important manufacturing topics in different time periods. 

In Industry 2.0 (from the end of nineteenth century to the 1980s), demand was relatively stable. Important manufacturing systems include the Toyota production system (TPS) and cellular manufacturing. The aim of these systems is to increase productivity with lower production cost, which fits the needs of a stable market during this time period. 

Taichi Ohno published Toyota Seisan Hoshiki describing the TPS in 1978. TPS is an integrated production system that can supply products to meet both requirements of product volumes and product varieties. Research and practical papers, reports, and books were published in various media to describe TPS. The underlying management principles and theoretical mechanisms of TPS are well-known. A TPS is an integrated production system that generates products to satisfy requirements of volumes and varieties simultaneously with minimum resource waste. A large number of TPS enablers have been reported and include just-in-time material system (JIT-MS), seven wastes, heijunka, multi-skilled workers, quick set-up and changeover, and keiretsu. Excellent analysis and review papers on the TPS are \cite{Treville2006-ak_KESXZ}, \cite{Hines2004-fo_KESXZ}, and \cite{Narasimhan2006-qt_KESXZ}. 

Cellular manufacturing (CM) uses group technology to efficiently produce a high variety of parts. Cells are converted from job shops with functional layouts to improve efficiency \citep{Yin2006-vm_KESXZ}. A cell consists of a machine group and a part family. The first step in CM system design is cell formation. Part families and machine groups are identified to form manufacturing cells such that the intercell movements of parts are minimised. Parts in the same family have similar machining sequences. Machines in a cell are arranged to follow this sequence. In this way, parts flow from machine to machine in their processing sequence, resulting in an efficient machining flow that is similar to an assembly line. For each part family, the volume of any particular part type may not be high enough to utilise a dedicated cell. The total volume of all part types in a part family should be high enough to utilise a machine cell well. CM attempts to flexibly accommodate high variety and simultaneously efficiently take advantage of flow lines \citep{Celikbilek2015-dz_KESXZ}.

In Industry 3.0 (from the 1980s to today), demand is relatively volatile because of technological innovations, higher product variety, and shorter product life cycles. The important manufacturing topics include flexible manufacturing systems (FMSs) and \textit{seru} production system. The theme of these topics is to meet the increased demand for high variety and short delivery time. Product life cycles decreased during this time period, which drives manufacturers to focus on responsiveness and delivery time. Short changeover time between different product types is useful. 

An FMS is an integrated, computer-controlled manufacturing system of automated material handling and computer numerically controlled machine tools that can simultaneously process medium-sized volumes of a variety of part types. A fully automated FMS can attain the efficiency of well-balanced, machine-paced transfer lines, while utilising the flexibility that job shops have to simultaneously machine multiple part types \citep{Stecke1981-eb_KESXZ,Stecke1983-lc_KESXZ,Browne1984-ec_KESXZ}.

A seru production system is an assembly system that has been adopted by many Japanese electronics companies. 
\cite{Yin2008-jw_KESXZ} is the first English language paper on seru production. They describe and analyse the success of seru production systems in Canon and other Japanese companies. It is more flexible than TPS, which cannot achieve the required responsiveness in this innovative time period. A seru production system consists of one or more serus. Serus within a seru system are quickly reconfigurable, i.e. they can be constructed, modified, dismantled, and reconstructed frequently in a short time. There are three types of serus, called divisional serus, rotating serus, and yatais. They represent the evolutionary development of serus. A divisional seru is a short, often U-shaped, assembly line staffed with several partially cross-trained workers. Tasks within a divisional seru are divided into different sections. One worker is in charge of each section. A rotating seru is often arranged in a U-shaped short line with several workers. Each worker performs all required tasks from start to finish without interruption. Tasks are performed on fixed stations, so workers walk from station to station. A worker follows the worker ahead of her or him, and is also followed by the worker behind her or him. A seru with only one worker is called a yatai. An important performance of the seru production system is that it can quickly respond to product varieties with fluctuated volumes. By applying seru, delivery time is reduced. Variety and volume are easily handled. 

The TPS-based assembly line became inefficient because of an inability to change very frequently to produce small-volume demands. The typical seru creation process in Sony and Canon can be summarised as follows \citep{Yin2017-gh_KESXZ}. Assembly lines were dismantled and replaced with divisional seru systems through resource co-location and removal/replacement, cross training, and autonomy. The technique of karakuri (involves procedures to discover and appropriate the useful functions of expensive equipment into inexpensive self-made equipment) is applied to replace expensive dedicated equipment by inexpensive self-made and/or general-purpose equipment that can be duplicated and redeployed as needed by serus. As cross-training progresses, divisional serus evolve into rotating serus and yatais. 

More details about the underlying management and control principles of seru can be found in \cite{Stecke2012-rg_KESXZ}, \cite{Yin2008-jw_KESXZ}, \cite{Yin2017-gh_KESXZ}, and \cite{Liu2014-mk_KESXZ}. \cite{Roth2016-pg_KESXZ} reviews the last 25 years of OM research and provides eight promising research directions, one of which is seru production systems. 

Today manufacturing has entered a new age (Industry 4.0) because of the emergence of disruptive technological innovations. Examples of important manufacturing topics include smart manufacturing, mass-customisation, sustainable manufacturing, and additive layer manufacturing. \cite{Strozzi2017-ff_KESXZ} examines the evolution, trends, and emerging topics of a smart factory and provides topics for future research. \cite{Hughes2022-jg_KESXZ} provides a review for manufacturing in the Industry 4.0 era. 

Smart manufacturing refers to flexible and adaptable manufacturing processes through integrated systems and using advanced technologies such as sensors, IoT, cloud computing, big data, artificial intelligence, automation, robots, cyber-physical systems, and additive layer manufacturing. Some detailed discussions can be found in \cite{Ivanov2016-qo_KESXZ}, \cite{Kersten2017-zq_KESXZ}, \cite{Liao2017-qi_KESXZ}, \cite{Theorin2017-ae_KESXZ}, \cite{Thoben2017-il_KESXZ}, and \cite{Hughes2022-jg_KESXZ}.   

One important benefit of smart manufacturing is that it aids the capability of mass customisation and short lead time to quickly meet changing demands. \cite{Zawadzki2016-ek_KESXZ} suggested smart product design and production control for efficient operations in a smart factory to enable mass customisation. \cite{Brettel2014-ho_KESXZ} show that self-improving smart manufacturing systems can utilise data and quickly react (e.g., reconfigure) to personalised customer orders, which helps realise mass customisation. Some efficient mathematical models that use big data to manage and control manufacturing processes can be applied in smart factories \citep{Ivanov2016-qo_KESXZ,Ivanov2017-rv_KESXZ}. 

Sustainable manufacturing aims to minimise negative environmental impacts while conserving energy and natural resources. Sustainable manufacturing also enhances employee, community, and product safety. The emergence of blockchain technology and its potential disruption within the manufacturing and supply chain industries present opportunities for greater levels of sustainability in Industry 4.0. The immutability and smart contract capability of blockchain technology allow the provenance and integrity of products to be monitored more effectively. These factors contribute to reducing verification costs and the provision of real-time status information on the quality of materials throughout the supply chain \citep{Ko2018-pi_KESXZ}. The disintermediation attributes of blockchain can directly contribute to manufacturing sustainability by effectively reducing complexity, and improving efficiency with less waste via the streamlining of the supply chain \citep{Hughes2019-iv_KESXZ}. 

Additive layer manufacturing may generate a disruptive and revolutionary impact on manufacturing \citep{Garrett2014-on_KESXZ}. It enables a manufacturer to further increase responsiveness by reducing lead time and increasing customisation levels. \cite{Long2017-ie_KESXZ} provide a definition, characteristics, and mainstream technologies of 3D printing. \cite{Dong2016-aq_KESXZ} compared the optimal assortment strategies under traditional flexible technology and 3D printing to find that 3D printing may allow a larger set of product assortment. \cite{Song2020-ti_KESXZ} and \cite{Ivan2017-ux_KESXZ} examined the use of 3D printing on a logistics system for spare parts inventory design. 3D printing tends to be slower than other manufacturing methods, which currently limits its use in practice.  

For a detailed encyclopedic overview of the manufacturing field, both in terms of theory and practice, see \cite{Yin2017-gh_KESXZ}. They discuss and compare production systems from Industry 2.0 to Industry 4.0. The demand dimensions of each industry era are analysed and provided as the driving force for the changes in the production systems over time.

\subsection[Military and homeland security (Kai~Virtanen \& Raimo~P.~Hämäläinen)]{Military and Homeland security\protect\footnote{This subsection was written by Kai~Virtanen and Raimo~P.~Hämäläinen.}}
\label{sec:Military_and_Homeland_security}

The birth of OR is related to the use of optimisation modelling for military operations and resource planning during the Second World War. The early linear programming (\S\ref{sec:Linear_programming}) problems ranged from the efficient use of weapon systems to logistics and strategy planning. Today, the arena of defence has expanded extensively with new areas including information and cyber warfare. The need to counter terrorism has created the field of homeland security. OR has a role in all these emerging topics. One can say that all OR methods are applied in military and homeland security problems.

Optimisation methods are used in a wide range of defence and security applications. For instance, assigning weapons to targets \citep{Kline2019-fx_KVRH} using integer programming (\S\ref{sec:Mixed_integer_programming}; \S\ref{sec:Combinatorial_optimisation}) has been addressed with a variety of optimisation algorithms. Other integer programming studies include, for example scheduling of training for military personnel \citep{Fauske2016-fi_KVRH} as well as military workforce and capital planning \citep{Brown2004-mg_KVRH}. Mixed integer linear programming is utilised in diverse applications such as path planning of unmanned ground and areal vehicles including drones, mission planning, acquisition decisions of military systems as well as load planning in transportation. Optimisation of vehicles’ routes is also carried out by solving network optimisation problems (\S\ref{sec:Graphs_and_networks}) with shortest path algorithms \citep{Royset2009-ps_KVRH}. Network optimisation is also used, e.g., in developing military countermeasures. Examples of bilevel and robust optimisation (\S\ref{sec:Stochastic_models}) formulations cover positioning of defensive missile interceptors \citep{Brown2005-ao_KVRH} and design of a supply chain for medical countermeasures against bioattacks \citep{Simchi-Levi2019-pv_KVRH}. Multiobjective optimisation has been applied, for example in optimising boat resources of coast guard \citep{Wagner2012-rt_KVRH} and planning of airstrikes against terrorist organisations \citep{Dillenburger2019-ir_KVRH}. Inherent structures of specific military optimisation problems have motivated the development of new solution techniques \citep{Boginski2015-iz_KVRH} including, for example, metaheuristics (\S\ref{sec:Heuristics}). Such techniques are used, e.g., for solving nonlinear military optimisation tasks (\S\ref{sec:Nonlinear_programming}) such as design of projectiles.

Game theoretic modelling (\S\ref{sec:Game_theory}) is used in many defence studies. Information related topics include misinformation in warfare \citep{Chang2022-xf_KVRH} and public warnings against terrorist attacks \citep{Bakshi2018-hw_KVRH}. Examples of game theoretic strategy design problems cover the optimal use of missiles and the validation of combat simulations \citep{Poropudas2010-jx_KVRH}. Designing security and counter strategies against enemies, terrorists and adversarial countries naturally lead to the use of game models. Interdiction network game models arise in security applications \citep{Holzmann2021-ib_KVRH}, and they are used, e.g., in route planning through a minefield.

Military simulation models (\S\ref{sec:Simulation}) are classified into constructive, virtual and live simulations \citep{Tolk2012-fx_KVRH}. Constructive simulations do not involve real-time human participation. They are based on well-known modelling methodologies such as Monte Carlo, discrete event and agent-based simulations. Applications of constructive models cover, e.g., the development and use of weapons, sensor and communications systems, planning of operations and campaigns, improving maintenance processes of military systems, and evaluating effects of fire. In addition, cyber-defence analyses have been conducted \citep{Damodaran2020-us_KVRH}. Constructive simulations have also been used in simulation-optimisation studies such as scheduling maintenance activities of aircraft, military workforce planning, and aircraft fleet management \citep{Mattila2014-mm_KVRH,Jnitova2017-vp_KVRH}.

The complexity of modelling human behaviour generates a major challenge for constructive simulation. This issue is avoided in virtual simulations, i.e. simulators in which real people operate simulated systems and in live simulations where real people operate real systems with simulated weapon effects. These practices are typically used, e.g., in military exercises and training of personnel. An emerging trend is to combine live, virtual and constructive simulations into a single simulation activity \citep{Mansikka2021-is_KVRH}. Applications of this simulation type vary from training to testing large-scale systems and mission rehearsals \citep{Hodson2014-ft_KVRH}. In a combined simulation, new ways to measure performance are introduced \citep{Mansikka2021-rw_KVRH} by complementing traditional measures such as loss exchange or kill ratio by human measures such as participants' situation awareness, mental workload and normative performance \citep{Mansikka2019-ol_KVRH}.

Features of virtual simulation can be recognised in wargaming \citep{Turnitsa2021-dz_KVRH} that has been used for military training and educating since the early 19\textsuperscript{th} century. Other wargaming areas are, for example, examination of warfighting tactics as well as evaluation of military operations and scenarios. Nowadays, wargames are also applied in studies of international relations and security as well as in analyses of government policy, international trade, and supply-chain mechanics \citep{Reddie2018-zs_KVRH}. The implementations of wargames range from manual tabletop map exercises to computer-supported setups in which different OR and artificial intelligence techniques are utilised \citep{Davis2022-kx_KVRH}.

Dynamic phenomena regarding military and defence are often represented with differential or difference equations. Examples of these models are Lanchester attrition equations that describe the evolution of strengths of opposing forces in gunfire combat \citep[e.g.,][]{Jaiswal2012-md_KVRH}. There are also several modifications of these equations aiming to model, e.g., asymmetrical combat, tactical restrictions and even morale issues. Another example of simple combat models is the salvo model that represents naval combat of warships involving missiles \citep{Hughes1995-za_KVRH}. Optimal control (see also \S\ref{sec:Control_theory}) has been utilised, for example in planning optimal paths of military vehicles as well as in guidance systems of unmanned aerial vehicles, drones and missiles \citep{Karelahti2007-jb_KVRH}. For a recent overview, see for example \cite{Israr2022-woKAI}. Another type of optimal control application is the assignment of resources to counter-terror policies and measures \citep{Seidl2016-ga_KVRH}. Markov decision processes  and approximate dynamic programming (\S\ref{sec:Dynamic_programming}) have recently emerged as important techniques for analysing dynamic military decision-making problems related to, e.g., missile defence interceptors and military medical evacuation \citep{Jenkins2021-nx_KVRH}.

The need for multicriteria evaluation is common in military decision-making. Example applications of multi-criteria decision analysis (MCDA; see also \S\ref{sec:Decision_analysis}) are acquisition of military systems and equipment procurement, military unit realignments and base closures, locating military bases, and assessment of future military concepts and technologies \citep{Ewing2006-as_KVRH,Geis2011-zd_KVRH,Harju2019-mb_KVRH}. Public procurement even for the military is regulated in many countries, and directives require it to consider multiple criteria \citep{Lehtonen2022-yc_KVRH}. It is interesting to notice that the recent acquisition decision of a 5\textsuperscript{th} generation multirole fighter aircraft in Finland was, indeed, supported by MCDA \citep{Keranen2018-ok_KVRH}. MCDA weighting methods have also been used to create measures of mental workload in military tasks \citep{Virtanen2021-ft_KVRH}. In portfolio decision analysis problems, the goal is to find the best set of components, e.g., in weapons systems or in force mix for reconnaissance, with respect to multiple criteria \citep{Burk2011-ns_KVRH}. The evaluation of the effectiveness of military systems calls also for the use of cost-benefit analysis  \citep[\S\ref{sec:Risk_analysis};][]{Melese2015-zf_KVRH}. Data envelopment analysis (\S\ref{sec:Data_envelopment_analysis}) is a multicriteria approach helping to seek efficiency also in military problems such as personnel planning.

MCDA studies in homeland security is a broad area ranging from the design of countermeasure portfolios to threat analysis and cyber-security \citep{Wright2006-ug_KVRH}. The questions of interest include, e.g., identification of terrorists' goals and preferences, estimation of attacks' consequences, and comparison of countermeasure actions \citep{Abbas2017-sv_KVRH}. Cost-benefit models are also relevant in terrorism research \citep{Hausken2018-gv_KVRH}.

Today, we are witnessing the vast growth of the use of machine learning and artificial intelligence (\S\ref{sec:Artificial_intelligence_machine_learning_data_science}) in military and security problems \citep{Dasgupta2022-mj_KVRH,Galan2022-xx_KVRH}. Such problem areas are, e.g., wargaming and simulation, command and control of autonomous unmanned vehicles including drones, air surveillance, and cyber-security only to mention a few. Data analytics (see also \S\ref{sec:Business_analytics}) is naturally also used in military OR \citep{Hill2020-sv_KVRH}, e.g., for supporting logistics planning. Considering uncertainty is essential, e.g., in intelligence analysis and risk analysis related to terrorism (see also \S\ref{sec:Risk_analysis}) . Adversarial risk analysis \citep{Rios_Insua2021-yx_KVRH} uses Bayesian approaches (see also \S\ref{sec:Risk_analysis}) for taking into account information, beliefs and goals of adversaries. A similar approach is also applied in the modelling of pilot decision-making in air combat with influence diagrams \citep{Virtanen2004-su_KVRH}. Markov models and Bayesian networks are used to evaluate risks and conduct time dependent probabilistic reasoning related to military missions \citep{Poropudas2011-ns_KVRH}. \cite{Kaplan2010-dt_KVRH} studies the infiltration and interdiction of terror plots using queueing theory (\S\ref{sec:Queueing}).

In the future, combat models need to include socio-cultural and behavioural factors \citep{Numrich2012-ov_KVRH}. We are also likely to see an increase in the modelling of individual and group behaviour as well as the consideration of behavioural issues in military and homeland security contexts. Behavioural game theory can give insights into military strategy and conflict situations. Behavioural OR (\S\ref{sec:Behavioural_OR}), which studies the impacts of the human modeller and model users including cognitive biases in decision support, is likely to receive increasing attention in military applications as well.

For further readings, we refer to the military OR textbooks by \cite{Fox2019-ky_KVRH} and \cite{Jaiswal2012-md_KVRH}. The recently edited volume by \cite{Scala2020_JLYHK} describes various OR methods and how to apply them in military problems. \cite{Abbas2017-sv_KVRH} and \cite{Herrmann2012-cg_KVRH} focus on homeland security modelling.

\subsection[Natural Resources (Emel~Aktas)]{Natural Resources\protect\footnote{This subsection was written by Emel~Aktas.}}
\label{sec:Natural_Resources}

Climate change and natural resource management require different quantitative and qualitative models that support public policy \citep{Ackermann2021_EA}. One of the early papers on the use of modelling for natural resource utilisation describes a resource analysis simulation procedure to assess the environmental impact of human activities \citep{Bryant1978_EA}. The procedure comprises a structural model to express the complex network of interacting human activity systems and a parametric model to determine the scale of the activity being modelled. 

An integrated decision support system for water distribution and management was built to generate alternative water allocation and agricultural production scenarios for a semi-arid region \citep{Datta1995_EA}. The model considers ground and surface water sources as the supply. The water demand is a combination of the need for drinking, irrigation, household and public utility, natural vegetation, industrial use, and ecological balance. The decision support tool visualises water allocation to competing crops under a range of simulation scenarios, providing a wider set of options to the department taking decisions how water is distributed. 

A web-based decision support system developed for the US Fish and Wildlife Service and the US Geological Survey initiative facilitates cross-organisational data sharing and performs analyses to improve conservation delivery \citep{Hunt2016_EA}. Situation-specific management actions such as controlled burn or prescribed graze required by this decision support tool improves ecological outcomes of other conservation efforts. Buffelgrass is an invasive species that causes significant damage to the native desert ecosystem. A multi-period multi-objective integer programming model was proposed to find optimal treatment strategies to control the buffelgrass population in the Arizona desert \citep{Buyuktahtakin2014_EA}. The multiple objectives minimise damage to threatened resources: a native cactus species (saguaros), buildings, and vegetation subject to budget and labour constraints. The results show the necessity of cooperation between different interest groups to establish reasonable treatment strategies and the need for a policy change because current resources cannot stop an ecological disaster in the future.

A mixed-integer programming model is developed to evaluate fishery management policies over an infinite horizon by incorporating steady state levels of variables into a multi-period analysis framework \citep{Glen1997_EA}. This model is intended to be used annually with updated stock estimates to set a dynamic total allowable catch per year depending on the stock estimates over several years. Statistical and fuzzy multiple criteria analysis establishes which materials contribute the most to the presence and the abundance of species in artificial reef structures \citep{Shipley2018_EA}. Managers of fisheries can use this model to screen different species without loss of rigour and validity of results. Multiple-criteria analysis is used in conjunction with integer programming to assist complex management plans in ecology and natural resources \citep{miranda2020_EA}. A case study on the Mitchell River catchment (Australia) shows the trade-offs between ecological, spatial, and cost criteria, enabling decision-makers to explore and analyse a broad range of conservation plans. The use of catastrophe theory in management of natural resource systems are described with cases on forestry and fishery management \citep{Wright1983_EA}. Catastrophe theory applies the mathematical theory of structural stability to practical systems. It allows modelling of ecosystems with low and high levels of predator and prey populations. It helps model a catastrophic jump from one level to the other, supporting decision making for management of natural resources. 

As a natural resource, wind provides clean, renewable, and sustainable energy. A multi-objective model minimises cost and idle time under reliability thresholds, maintenance priority, and opportunism \citep{Ma2022_EA}. Reliability thresholds trigger maintenance activity. Maintenance priority indicates which maintenance tasks need to be performed under limited maintenance resources. Opportunistic approach indicates whether additional maintenance should be performed when a maintenance team is already out to service several turbines. The proposed multi-objective optimisation model is tested using a stochastic simulation model of a wind farm and confirmed to keep the wind system at a higher performance level with lower cost and higher availability.

Natural resource exploration is frequently subject to real options analysis \citep{Nishihara2012_EA, Martzoukos2009_EA}. A stochastic mixed integer nonlinear programme is proposed to incorporate geological and market uncertainty into mineral value chain optimisation \citep{Zhang2018_EA}. Simulation of mine deposits and commodity market informs the profitability of strategic and tactical plans. A range of real options applied to natural resource management include investments in infrastructure, use of land, and management of natural resources \citep{Trigeorgis2018_EA}. Firms require high output price levels to invest in environmental technologies, because they would not want to commit to an investment that could turn out to be unprofitable in the event of a price fall \citep{Cortazar1998_EA}.

Several papers are published on the use of operational research for natural resource management. Typical operational research problems and actors in agricultural supply chains informs strategic investment and operations management under increased pressure on natural resources \citep{pla2014_EA}. The contribution of operational research applications to agricultural value chain sustainability and resilience call for applications of complex systems methods such as agent-based modelling, systems modelling, and network analysis \citep{Higgins2010_EA}. A review of environmental management and sustainability papers in major management science/operational research and systems journals revealed dominance of hard optimisation methods \citep{paucar2011_EA}. 

The environment-development problem concerns reconciling industrial development and environmental protection. A methodological framework is proposed to model the environmental impact of development under uncertainty arising from the degree of unpredictability arising from decision makers and environmental processes \citep{Dzidonu1993_EA}. Natural resource development contracts depend on the bargaining power of transnational corporations and host country governments \citep{Anandalingam1987_EA}. Contracts that stipulate sharing of the net income from resource development between the developing corporation and the government show that government would receiver higher income if many corporations are involved and if the government agrees to contribute to production costs. A review of Operational Research in mine planning reports optimisation and simulation applied to surface and underground mine planning problems, including mine design, long- and short-term production scheduling, and equipment selection \citep{Newman2010_EA}. The operational research on mining is evolving to solve larger and more detailed and realistic models.

A series of cases studies from Asia, Africa, and Latin America presents principles and applications of an integrated approach to natural resources management, including the complexity of systems and redirecting research towards participatory approaches, multi-scale analysis, and tools for systems analysis, information management, and impact assessment \citep{Campbell2003_EA}. A specialised book on the Baltic region presents scientific research on activities depleting natural resources, emissions from energy use, pollution, and strategies for environmental management \citep{Fenger1991_EA}. Stochastic Models and Option Values: Applications to Resources, Environment and Investment Problems presents a collection of research papers on the use of control theoretic methods to address problems that arise in natural resource development \citep{Lund1991_EA}. Strategic Planning in Energy and Natural Resources contain innovative and methodologically rigorous operational research applications \citep{Lev1987_EA}. The Handbook of Operations Research in Natural Resources collate research papers that address optimal allocation of scarce resources in agriculture, fisheries, forestry, mining, and water resources \citep{Weintraub2007_EA}. Operations Research and Environmental Management organises its content by regional and global policies \citep{Haurie2012_EA}. Models help local and regional authorities optimise their energy distribution and minimise natural resource waste.  

\subsection[Open-source software for OR (Cihan~Tugrul~Cicek \& Güneş~Erdoğan)]{Open-source software for OR\protect\footnote{This subsection was written by Cihan~Tugrul~Cicek and Güneş~Erdoğan.}}
\label{sec:Open_source_software_for_OR}

Commercial solvers for solving Operational Research (OR) problems have been used for several decades and have provided both practitioners and academics with access to the state-of-the-art OR techniques. Mathematical programming solvers IBM ILOG CPLEX \citep{Cplex_CTCGE}, Gurobi \citep{Gurobi_CTCGE}, BARON \citep{Baron_CTCGE}, and discrete event simulation software Arena \citep{Arena_CTCGE} and Simul8 \citep{Simul8_CTCGE} are among the best-known examples. There also exists ad-hoc software for particular problems raising in manufacturing \citep[e.g., AIMMS;][]{Aimms_CTCGE}, healthcare \citep[e.g., SimCAD Pro Health;][]{Simcad_CTCGE}, and logistics \citep[e.g., AnyLogistix][]{Anylogistix_CTCGE}. The strength of commercial software is primarily based on the fact that they provide users with a simple interface to declare a problem, utilise state-of-the-art solution algorithms, and visualise the result with minimal effort. 

These software do not only solve problems but also provide modelling, debugging, and scenario analysis to improve the solution process \citep{Dagkakis2016_CTCGE}. However, the lack of access to the source code and knowledge of how these tools work internally inhibit users from customisation. It is difficult to contribute to the development of commercial software as it is a black-box to the end-user. The high licence costs of those software has been one of the most prominent factors blocking many companies, especially small and medium enterprises, from integrating them into their tactical and operational planning \citep{Linderoth2005_CTCGE}. \cite{Dagkakis2016_CTCGE} argued that the lack of reusability and modularity have been the additional factors impeding the use of commercial OR software.

Open-source software, on the other hand, enable users to solve OR problems without a significant initial investment. Although using open-source software does not require licensing fees, the effort to deploy it may require a significant amount of effort and time. Nevertheless, the opportunity to access the core components of a solver (or simulator) and ease of development has driven the OR community to shift from a strict, slow-pace black-box software development to modular, flexible, and quick open-source software development. 

In this section, we discuss open-source OR software, categorising them into two main groups: (\textit{i}) open-source solvers and (\textit{ii}) open-source simulators. The former category covers the software focusing on solving mathematical programming problems. The latter includes software for simulating a real-world environment and helping decision makers to understand and analyse the system without consuming physical resources. Note that we neither provide the specific features of such software nor the characteristics in terms of programming languages, etc. Interested readers are referred to the comprehensive reviews in \cite{Linderoth2005_CTCGE} and \cite{Dagkakis2016_CTCGE} for some of the software we mention below.

\subsubsection*{Open-source Solvers}
A solver can be defined as a set of computationally efficient analytical tools that can find optimal (or near-optimal) solutions to a mathematical programming model. In 2000, a public initiative was built by the IBM Research Division \citep{Coinor_CTCGE} to promote and support community-driven development of open-source solvers that utilise the state-of-the-art research in OR. Subsequently, a public project called COIN-OR \citep{Coinor2_CTCGE} has been initiated to host a range of open-source solvers with an open-source interface that enables contributors, users, and developers to implement their own algorithms. The repository has been expanded to provide different open-source solvers for different programming problems such as CLP \citep{CLP_CTCGE} and HiGHS \citep{Huangfu2018_CTCGE} for linear programming (LP); ABACUS \citep{Abacus_CTCGE}, BCP \citep{Bcp_CTCGE}, CBC \citep{Cbc_CTCGE}, Pyomo \citep{Pyomo_CTCGE}, and SYMPHONY \citep{Symphony_CTCGE} for mixed integer linear programming (MILP); Bonmin \citep{Bonmin_CTCGE}, Couenne \citep{Couenne_CTCGE}, DisCO \citep{Disco_CTCGE}, Ipopt \citep{Ipopt_CTCGE}, and SHOT \citep{Shot_CTCGE} for linear and mixed integer nonlinear programming (MINLP); SMI \citep{Smi_CTCGE} and Pyomo \citep{Pyomo_CTCGE} for Stochastic Programming (SP). COIN-OR also includes several other projects that would help users to improve their experience with modelling like PuLP \citep{Pulp_CTCGE} and visualising such as GiMPy \citep{Gimpy_CTCGE}.

SCIP Optimization Suite \citep{Scipopt_CTCGE} can be used as a framework for mixed-integer linear or nonlinear programming as well as a standalone solver for such problems. A recent initiative commenced by the introduction of Julia language \citep{Julia_CTCGE}, which is a high-level, high-performance dynamic language for technical computing, is JuMP \citep{Jump_CTCGE} which helps users to solve a variety of problem classes including linear, mixed-integer linear, and nonlinear programming. It allows developers to use its framework and introduce new open-source solvers for particular problem classes.

We would like to also mention GLPK \citep{Glpk_CTCGE}, which is the default linear programming solver behind some of the aforementioned mixed integer linear programming solvers. GLPK can also be used as a standalone linear programming solver. Finally, a suite of open-source solvers has been developed by Google \citep{ORTools_CTCGE} to tackle integer programming and constraint programming problems. The OR-Tools provide a modelling interface and allow users to select different commercial or open-source solvers to generate solutions.

\subsubsection*{Open-source Simulation Software}
Simulation software can be categorised into three based on the methods that they use to define the system and its resources. We should note that we cover the software that has been applied particularly in OR domains. We opt to omit open-source software that focus on specific domains, e.g. OMNet++ \citep{Omnet_CTCGE} for communication networks, for the sake of brevity. We refer interested readers to the works of \citep{Dagkakis2016_CTCGE} and \citep{Lang2021_CTCGE} and references therein for a broader review. An experimental comparison of some of the software presented here can be found in \cite{Kristiansen2022_CTCGE}.

The first method, Discrete Event Simulation (DES), is based on the processes of a system. In DES, the processes are defined as hosts of resources that run different operations on entities. For instance, one can define a part to be manufactured as an entity and create a manufacturing process with three machines to shape the part. DES software can be used to model and visualise complex queuing systems in order to help decision makers better understand the interactions between entities and processes.

JaamSim \citep{Jaamsim_CTCGE} is one of the most popular open-source DES with its user-friendly interface, easy-to-use drag and drop facilities, and continuous maintenance support. JaamSim provides a standalone executable which allows users to start using the software without technical knowledge on installations. Another DES framework is SimPy \citep{Simpy_CTCGE} which is based on standard Python functions. Its simple structure enables users to quickly obtain results for their simulation problems. SimPy has also initiated two other DESs, SimSharp \citep{Simsharp_CTCGE} and SimJulia \citep{SimJulia_CTCGE}, which are the implementations of SimPy on C\# and Julia languages, respectively. The last DES we would like to mention is Facsimile \citep{Facsimile_CTCGE} which uses Scala as its basis scripting language. The main purpose of Facsimile, is to provide a high-quality discrete-event simulation library that can be used for industrial projects.

The second method, known as System Dynamics (SD), is based on representing a system as a \emph{causal loop diagram} to define interactions between different components of a system. Some of the open-source SD software are PySD \citep{Pysd_CTCGE}, InsightMaker \citep{Insightmaker_CTCGE}, SysDyn \citep{Sysdyn_CTCGE}, and OpenModelica \citep{Openmodelica_CTCGE}. PySD can convert the well-known commercial SD software Vensim \citep{Vensim_CTCGE} input and allow user to configure. SysDyn uses the OpenModelica environment for simulation but provides an alternative built-in environment to speed up the simulation process. All these software have their own visualisation and reporting tools. 

The third method is called Agent Based Simulation (ABS) and focuses on the autonomous individuals, i.e., agents, in a system. Each agent in ABS has its own characteristics and its way to interact with the other agents and the surrounding environment can differ. One of the open-source ABS software is Gama \citep{Gama_CTCGE}, which provides users a modelling language, a cross-platform to reproduce simulations, as well as a visualisation tool. InsightMaker \citep{Insightmaker_CTCGE} is another open-source ABS software that allows users to create their own model on a web-based interface. Lastly, NetLogo \citep{NetLogo_CTCGE} provides a modelling environment together with different applications to interact with other scripting languages.

We would like to complete this section with a brief summary of the application areas of both solvers and simulation software. Table~\ref{tab:application_areas} provides examples of areas on which the OR software can be used.

\begin{table}[!htb]
    \centering
    \caption{Application areas of solvers and simulation software.}
    \begin{tabular}{p{0.12\textwidth}p{0.12\textwidth}p{0.65\textwidth}}
        \toprule
        Subject & Methodology & Application Areas \\
        \midrule
        Solver & LP & Transportation, agriculture, manufacturing \\
        & MILP & Logistics, healthcare, network design, pooling, disaster response  \\
        & MINLP & Scheduling, telecommunication, energy systems, layout design, network design, portfolio optimization, water systems \\
        & SP & Supply chain planning, production planning, process control and optimisation \\
        \midrule
        Simulation & DES & Manufacturing, network design, healthcare operations, financial systems, inventory management \\
        & SD & Telecommunication, macro- and micro-economics, social systems, ecological systems \\
        & ABS & Stock markets, robotics, epidemiology, game theory, evacuation planning \\
        \bottomrule
    \end{tabular}
    \label{tab:application_areas}
\end{table}

\subsubsection*{Discussion}
For the sake of completeness, we should also mention that there are also several ad-hoc software that address specific OR problems. For instance, OptaPlanner \citep{Optaplanner_CTCGE} can solve staff rostering, scheduling, timetabling, and Vehicle Routing Problems (VRPs). Another example is VRP Spreadsheet Solver \citep{Erdogan2017_CTCGE}, which is an Excel-based solver. Although these software provide easy and fast access to solutions, the lack of generalisation to more complex OR problems and limited development opportunities can be seen as barriers to widespread impact.

As a final discussion point, we would like to list some of the essential features of an open-source OR software. First and foremost is the performance of the software. A user would expect a comparable level of performance from an open-source software with respect to commercial software. Secondly, the scalability of a solver, i.e. its performance when the problem size increases, is one of the factors desired by practitioners. Finding an optimal solution to a VRP instance with 20 customers does not guarantee that a VRP solver will achieve the same performance when the number of customers increases to 2000. Thirdly, technical support for a software has a crucial role in attracting users. Continuous development, documentation, and clear descriptions to change requests are some of the aspects that an open-source software should address to improve its maintainability. Finally, integration with existing libraries would help an open-source software widen its community and attract more developers to contribute.

\subsection[Power markets and systems (Dimitrios~Sotiros \& Rafał~Weron)]{Power markets and systems\protect\footnote{This subsection was written by Dimitrios~Sotiros and Rafał~Weron.}}
\label{sec:Power_markets_and_systems}

The energy industry relies on forecasts (\S\ref{sec:forecasting}) and decision support tools (\S\ref{sec:Decision_analysis}) for operations and planning. While long-term demand forecasts -- with lead times measured in months, quarters or years -- have been used for planning purposes for over a century, contemporary energy forecasting literature focuses more on the short- (minutes, hours) and mid-term (days, weeks) horizons \citep{hon:pin:etal:20_DSRW}. 
Since the late 1990s, the workhorse of power trading and a typically used reference point for long-term contracts is the day-ahead market \citep{may:tru:18_DSRW}, where prices for all load periods of the next day are determined at the same time during a uniform-price auction \citep[][see also \S\ref{sec:Auctions_and_bidding}]{wer:14_DSRW}. No wonder, the majority of studies focus on predicting intermittent generation from \textit{renewable energy sources} (RES), electric load (or demand) and prices for the 24 hours of the next day \citep{mac:nit:wer:21_DSRW}. Two classes of approaches dominate: regression-based models and \textit{artificial neural networks} \cite[ANN;][see also \S\ref{sec:Artificial_intelligence_machine_learning_data_science}]{lag:mar:sch:wer:21_DSRW}.

Regression and ANN models of the 1990s and 2000s were built on expert knowledge, often independently for each hour of the day. Their inputs included past values of (depending on the context) RES generation, loads or prices from the last few days, day-ahead forecasts of exogenous variables (e.g., temperatures for load, load for prices) and a seasonal component captured by sinusoidal functions or weekday dummies \citep{hon:14_DSRW,wer:14_DSRW}. Their sub-optimal  performance could be readily improved by combining forecasts across different models \citep{bor:bun:lis:nan:13_DSRW}, calibration sets \citep{now:liu:wer:hon:16_DSRW} or calibration windows \citep{hub:mar:wer:19_DSRW}. Interestingly, combining is not only a remedy for time-varying point forecasting performance. Together with quantile regression it provides a simple, yet powerful tool for probabilistic predictions -- \textit{Quantile Regression Averaging} \cite[QRA;][]{now:wer:18_DSRW}. During the Global Energy Forecasting Competition 2014, teams using variants of QRA were ranked 1\textsuperscript{st} and 2\textsuperscript{nd} in the price track \citep{gai:gou:ned:16_DSRW,mac:now:16_DSRW}. QRA can be also used to construct dynamic strategies aiming at finding the optimal trade-off between risk and return when trading the intraday and day-ahead markets \citep{jan:woj:22_DSRW}.

With the advent of easily accessible computational power, the models became more complex and expert knowledge was no longer enough to handle them. A major breakthrough came with the introduction of regularisation methods to energy forecasting in the 2010s. Although regularisation is a much older concept, its use in load \citep{cha:hor:hwa:lee:16_DSRW,zie:liu:16_DSRW}, price \citep{uni:now:wer:16_DSRW,zie:16_DSRW,zie:wer:18_DSRW}, wind \citep{mes:pin:19_DSRW} and solar forecasting \citep{yan:18_DSRW} began only recently. Ridge regression has not seen many applications in energy forecasting, however, the \textit{least absolute shrinkage and selection operator} (LASSO) and elastic nets \citep{has:tib:wai:15_DSRW} have been shown to yield extremely competitive predictive models. \textit{LASSO-estimated autoregressive} (LEAR) models often have hundreds of inputs, e.g., spanning all hours of the past week, but LASSO can shrink redundant coefficients to zero and, thus, perform variable selection. Despite their ability to handle only linear relationships between variables, LEAR models tend to be only slightly inferior to the much more complex and much harder to estimate deep ANNs \citep{lag:mar:sch:wer:21_DSRW}. 

The availability of high-performance GPUs and advances in optimisation algorithms  made it possible to efficiently train ANNs with hundreds of inputs and outputs, multiple hidden layers and recurrent connections (\S\ref{sec:Artificial_intelligence_machine_learning_data_science}). This led to a wave of deep learning and hybrid energy forecasting models in the late 2010s \citep{gao:etal:19_DSRW,wan:etal:17_DSRW}. A prominent, yet relatively parsimonious example is the deep neural network (DNN) model proposed for price forecasting \citep{lag:rid:sch:18_DSRW}. It uses a feedforward architecture with two hidden layers, 24 outputs (one for each hour of the next day) and ca.\ 250 inputs: past prices from the previous week, day-ahead forecasts of fundamental variables (demand, RES generation) and weekday dummies. To decrease the computational burden, its hyperparameters (number of neurons per layer, activation functions, optimisation algorithm, etc.) and inputs (treated as binary hyperparameters -- either selected or discarded) are jointly optimised once every few weeks, while the weights are recalibrated every day to account for the most recent market data. Despite this simplification, daily recalibration of the DNN model is two orders of magnitude slower than of the LEAR model with the same inputs \cite[minutes vs.\ seconds on a quadcore i7 CPU; see][]{lag:mar:sch:wer:21_DSRW}.

The increased complexity of deep ANNs is a major obstacle in understanding the underlying processes. Partial remedy provide recently proposed architectures like the neural basis expansion analysis for interpretable time series forecasting \citep[NBEATS;][]{ore:dud:pel:tur:21_DSRW,oli:cha:mar:wer:dub:22_DSRW}, which project the time series onto basis functions in the fundamental blocks of the network structure. The final forecasts can be decomposed into interpretable components returned by groups of blocks (called stacks). Separate stacks can account for the trend, seasonality and exogenous variables. Another recent innovation in energy forecasting is a distributional ANN \citep{mas:etal:21_DSRW}. It only requires a small change in the architecture -- instead of the 24 hourly forecasts, the network can return the parameter sets of 24 probability distributions (e.g., the mean and standard deviation for the Gaussian). The benefits are clear. The downside, however, is that the distribution itself has to be estimated  (it is a hyperparameter). Somewhat surprisingly, distributional ANNs not only can yield more accurate probabilistic predictions, but also better point forecasts \citep{mar:nar:wer:zie:22_DSRW}.

For horizons beyond the next 48 hours other  approaches have been proposed \citep{wer:14_DSRW}, not necessarily forecasting per se. Structural models define the functional relationships between physical (weather, generation, consumption, etc.) and economic (bidding, trading) variables that set the price, then utilise -- typically -- parsimonious statistical or machine learning techniques to provide the stochastic inputs. Due to the nature of fundamental data, often of weekly or monthly granularity, such models are more suitable for medium-term risk management, portfolio optimisation and derivatives pricing \citep{kie:kus:16_DSRW}, than for short-term forecasting \citep{mah:gir:kar:22_DSRW}.
In the class of multi-agent approaches, \cite{ven:bai:ram:riv:05_DSRW} identify three trends:  equilibrium, simulation and optimisation models. The former (Nash-Cournot framework, supply function equilibrium, strategic production-cost models) have seen limited application in oligopoly markets \citep{rui:maz:08_DSRW}. Agent-based simulations are used when the problem is too complex to be addressed within an equilibrium framework \citep{fra:kra:kel:fic:21_DSRW}. 

Optimisation models address profit maximisation from the point of view of a firm competing in the market. One of the simplest settings is that of the price clearing process being exogenous to electricity generation optimisation -- as the price is fixed, the market revenue is a linear function of the production and linear programming (\S\ref{sec:Linear_programming}) or \textit{mixed integer linear programming} (MILP) can be employed \citep{ven:bai:ram:riv:05_DSRW}. On the other hand, \textit{Virtual Power Plant} (VPP) operations constitute a more complex problem of decision-making under uncertainty. A VPP is a cluster of dispersed generating units (e.g., intermittent rooftop solar panels on residential houses), flexible loads and battery storage that operates as a single entity. Robust optimisation and stochastic programming can be used to derive the optimal VPP trading strategy
\citep{mor:con:mad:pin:zug:14_DSRW}.

To support broader regulatory decisions at the firm or country level, frontier analysis methods are employed. Such methods aim to estimate the efficient frontier of the evaluated production units, to measure their relative efficiency (against the frontier) and to provide targets that can support policymakers. The benchmarking nature of these methods has established them as a flexible and multifaceted decision-making tool. In particular, \textit{Data Envelopment Analysis} (DEA, \S\ref{sec:Data_envelopment_analysis}) has been employed in a wide spectrum of energy applications. Early DEA studies relied only on a few factors (labour, fuel, capital, electricity production) to assess the technical efficiency of electric utilities \citep{Fare_1983_DSRW}. Later studies took into account sustainable practices by including environmental variables. Such factors are commonly treated as undesirable outputs that arise as by-products of the production process \citep{Fare_1996_DSRW} or as non-controllable variables, which reflect external factors that the unit under evaluation cannot control \citep{Hattori_2005_DSRW}. When price information is available, DEA allocation models can be used to evaluate revenue, cost and profit efficiency. \cite{Ederer_2015_DSRW} argued that sophisticated cost efficiency assessment methods should be employed to evaluate RES, and relied on DEA models to assess the capital and the operating cost efficiency of offshore wind farms. Notably, DEA is often combined with multi-criteria decision-making techniques to incorporate decision maker's preferences into the assessment \citep{Lee_2011_DSRW,Wang_a_2022_DSRW} and econometric techniques to study causal effects \citep{Shah_2022_DSRW}.

For a review and outlook into the future of energy (load, price, wind, solar) forecasting see
\cite{hon:pin:etal:20_DSRW}. \cite{hon:fan:16_DSRW} offer a tutorial review on probabilistic load forecasting. The standard reference for electricity price forecasting is \cite{wer:14_DSRW}. \cite{lag:mar:sch:wer:21_DSRW} offer a more recent viewpoint focusing on deep learning and hybrid models. They also provide a set of guidelines/best practices and make freely available the \textit{epftoolbox} with Python codes for two highly competitive benchmark models (LEAR, DNN). Two thorough treatments of probabilistic price forecasting are \cite{now:wer:18_DSRW} and \cite{zie:ste:18_DSRW}. \cite{swe:bes:bro:pin:20_DSRW} present a brief overview of the state-of-the-art in RES forecasting, whereas \cite{yan:etal:22_DSRW} jointly discuss atmospheric science and power system engineering in the context of solar forecasting. Finally, for detailed literature reviews on energy related applications of DEA see \cite{Mardani_2017_DSRW}, \cite{Sueyoshi_2017_DSRW} and \cite{Yu_2020_DSRW}.

\subsection[Project management (Willy~Herroelen \& Erik~Demeulemeester)]{Project management\protect\footnote{This subsection was written by Willy~Herroelen and Erik~Demeulemeester.}}
\label{sec:Project_management}

Operational Research methods play a fundamental role in managing a portfolio of projects, in project selection and in the management of each individual project. \textit{Project portfolio management} is concerned with the optimal mix and prioritisation of proposed projects in order to maximise the organisation’s overall goals \citep{Levine2007-bi_WH_ED}. At the strategic level, \textit{project selection} deals with the selection of and resource allocation among a group of projects \citep{Kavadias2012-mr_WH_ED}. Static models rely on mathematical programming, scoring and sorting, financial modelling, graphical and charting techniques. Dynamic models for selecting projects from a stream of arrivals may rely on queueing theory (\S\ref{sec:Queueing}), simulation (\S\ref{sec:Simulation}), decision theory (\S\ref{sec:Decision_analysis}) and stochastic dynamic programming (\S\ref{sec:Dynamic_programming}; \S\ref{sec:Stochastic_models}).

At the tactical and operational levels,  \textit{project management} \citep{Meredith2011-uu_WH_ED} basically involves the planning, scheduling and control of project activities to achieve performance, cost and time objectives for a given scope of work, while using resources efficiently and effectively. The scope of a project is the magnitude of the work that must be performed to make sure that the product or items to be provided (the project result or the project deliverables) meet the requirements or acceptance criteria agreed upon at the onset of a project. Once the project is properly defined in terms of its scope and objectives, the planning phase may start through the identification of the project activities, the estimation of time and resources, the identification of relationships and dependencies between the activities and the identification of the schedule constraints. The activities can be graphically portrayed in the form of a project network showing the necessary interdependencies of the activities. Based on the type and quantities of resources required, cost estimates can be made. \textit{Project scheduling} \citep{Demeulemeester2002-an_WH_ED} then involves the construction of a project base plan which specifies for each activity the precedence and resource feasible start and completion dates, the amounts of the various resource types that will be needed during each time period, and as a result the budget. Once a baseline schedule has been established, it must be implemented. This involves performing the work according to the plan and controlling the work by monitoring the progress and taking necessary corrective action when the project is on its way to run behind schedule, to overrun the budget, or to violate the original technical specifications.

\subsubsection*{Construction of the project network}
A project network is a graph consisting of a set of nodes and a set of arcs. In the activity-on-arc representation (AoA), the nodes represent the events and the arcs represent the activities. AoA networks form the basis of the Project Evaluation and Review Technique \citep[PERT;][]{Malcolm1959-yn_WH_ED} and the Critical Path Method \citep[CPM;][]{Kelley1961-ko_WH_ED}. The precedence relationship used is the finish-start relationship with a zero time lag: an activity can start as soon as all its predecessor activities have finished. In the mostly used activity-on-node representation (AoN) the nodes represent the activities and the arcs denote the precedence relations. The AoN representation allows for the specification of generalised precedence relations of four types: start-start, start-finish, finish-start and finish-finish with minimal and/or maximal time lags. A minimal time lag specifies that an activity can only start (finish) when its predecessor activity has already started (finished) for a certain time period, whereas a maximal time lag specifies that an activity should be started (finished) at the latest a number of time periods beyond the start (finish) of another activity.

\subsubsection*{Temporal analysis for deterministic unconstrained project scheduling}
In this case a single deterministic duration estimate is used for the activities. Basically, the temporal analysis then involves the computation of the activity start times under the objective of minimising the project duration. In the presence of strict finish-start precedence relations, this can be achieved by simple forward and backward pass calculations. Generalised precedence relations with maximal time lags call for the use of graph algorithms for computing the longest path (critical path) in networks. 

The temporal analysis may also be performed under the objective of maximising the net present value of the project. The deterministic max-npv problem can be formulated as a nonlinear problem. An efficient recursive solution procedure has been developed for AoN networks and has been extended to deal with the case of time-dependent cash flows \citep{Vanhoucke2001-kp_WH_ED}. 

Another non-regular performance measure is the minimisation of the weighted earliness-tardiness penalty costs of the project activities, where activities have an individual due date with associated unit earliness and unit tardiness penalty costs. The problem can be solved by an exact recursive search procedure \citep{Vanhoucke2001-nb_WH_ED}.

\subsubsection*{The deterministic resource-constrained project scheduling problem}
Project activities require resources for their execution. Renewable resources (e.g., manpower, machines) are available on a per-period basis. Their introduction into the analysis complicates matters considerably. Computing a precedence- and resource-feasible deterministic schedule that minimises the project duration, the resource-constrained project scheduling problem (RCPSP) is $\mathcal{NP}$-hard in the strong sense (\S\ref{sec:Computational_complexity}). Both exact and suboptimal procedures have been presented in the literature.

Many mathematical programming formulations (\S\ref{sec:Mixed_integer_programming}), either binary or mixed integer linear programs, have been developed \citep{Demassey2010-pt_WH_ED}. The RCPSP may also be solved through constraint-based scheduling \citep{Laborie2010-li_WH_ED}. Also a number of branch-and-bound algorithms have been presented for optimally solving the RCPSP \citep{Brucker1998-qg_WH_ED,Demeulemeester1992-xs_WH_ED}. 

Heuristic procedures broadly fall into two categories: constructive heuristics and improvement heuristics. Constructive heuristics start from an empty schedule and add activities one by one until a feasible schedule is obtained. Activities are ranked by priority rules which determine the order in which the activities are added to the schedule. Improvement heuristics start from a feasible schedule that was obtained using a constructive heuristic. Operations are then performed on a schedule which transforms a solution into an improved one. These operations are repeated until a locally optimal solution is obtained.

Project scheduling metaheuristics come in a wide variety and broadly include tabu search \citep{Baar1999-ni_WH_ED}, simulated annealing \citep{Bouleimen2003-ud_WH_ED}, genetic algorithms \citep{Hartmann2002-vg_WH_ED}, ant colony optimisation \citep{Merkle2002-xi_WH_ED}, scatter search and electromagnetic approaches \citep{Debels2006-di_WH_ED}. 

\subsubsection*{Resource problem variants and generalisations}
Branch-and-bound may be used for solving the RCPSP with generalised precedence relations \citep{Demeulemeester1997-bt_WH_ED,De_Reyck1998-sh_WH_ED}, when activities may be preempted \citep{Demeulemeester1996-yl_WH_ED}, when the problem has to be solved under the objective of maximising the net present value \citep{Vanhoucke2001-kp_WH_ED} or with the earliness-tardiness objective \citep{Vanhoucke2001-fb_WH_ED}.

The resource levelling problem aims at completing the project within its deadline with a resource usage that is leveled as much as possible over the entire project duration. Exact solution procedures based on integer or dynamic programming and branch-and-bound as well as heuristic procedures have been developed \citep{Neumann2000-yl_WH_ED}. The resource availability cost problem, that consists of scheduling the project activities such that the total cost of acquiring the necessary resources is minimised, assuming that a resource is assigned to the project for the total project duration, can be solved optimally \citep{Demeulemeester1995-gy_WH_ED}. The resource renting problem \citep{Nubel2001-ou_WH_ED} which assumes that resources can be added or removed from the resource pool over the project life, can be solved optimally using branch-and-bound or heuristically using genetic algorithms and scatter search \citep{Ballestin2007-xx_WH_ED,Kerkhove2017-et_WH_ED}. 

The multi-mode RCPSP assumes limited availability of renewable and nonrenewable (e.g., money) resource types and assumes that a project activity may be executed in multiple modes, where an activity mode corresponds to the assignment of a mode-specific number of units of a (non)renewable resource type to the activity with correspondingly resulting activity duration. The project decisions then involve the decisions to start and perform the activities in a specific mode in order to minimise the project duration. Branch-and-bound \citep{Hartmann1998-ny_WH_ED}, branch-and-cut, local search \citep{De_Reyck1999-oz_WH_ED} and genetic algorithms \citep{Hartmann2001-bu_WH_ED} are available. 

For projects with a flexible project structure, where activities to be performed are not known in advance, decisions for the implementation of optional activities can be made using genetic algorithms and tabu search \citep{Kellenbrink2015-av_WH_ED,Servranckx2019-zs_WH_ED,Servranckx2019-vw_WH_ED}.

\subsubsection*{Dealing with uncertainty}
Risk analysis involves the identification of the qualitative and quantitative assessment of the risk factors for the project through the estimation of the probability of the risk factors (activity duration, cost and resource requirement increases, start time delays) as well as their potential impact. The impact of each risk is best assessed individually and mapped to the duration of a project activity \citep{Creemers2014-lb_WH_ED}. Risk responses may then involve risk avoidance by performing an alternative approach without the risk, taking actions to reduce the risk, and risk impact reduction by switching to a different execution mode, adding additional resources, etc.

Stochastic scheduling does not generate a baseline schedule before the start of the project, but deals with time uncertainty by viewing the scheduling problem as a multi-stage decision process where scheduling policies are used to decide at each of the stages which occur serially through time at random decision points, which activities selected from the set of precedence and resource feasible activities have to be started under the objective of minimising the expected project duration \citep{Demeulemeester2011-eb_WH_ED}. 

Proactive project scheduling generates a robust baseline schedule through solving the RCPSP and subsequently tries to protect it as well as possible against time and resource disruptions that may occur during project execution. This protection can be achieved by deciding on a clever way to transfer the renewable resources between the activities \citep{Leus2004-kp_WH_ED}.  Both branch-and-bound and heuristics are available for the minimisation of the weighted sum of the difference between the planned and the realised activity start times \citep{Van_de_Vonder2008-hv_WH_ED,Lambrechts2008-fp_WH_ED}. Another way involves the insertion of time buffers that should prevent the propagation of distortions throughout the schedule. The critical chain methodology \citep{Goldratt1997-tq_WH_ED,Herroelen2001-vl_WH_ED,Newbold1998-xg_WH_ED} defines the critical chain as that set of tasks which determines the overall project duration. Protection is then realised through the insertion of feeding buffers and resource buffers in combination with a project buffer at the end of the critical chain. 

When during the actual execution of the project disruptions occur that cause deviations from the protected baseline schedule or even render this schedule infeasible, reactive scheduling procedures may be deployed.

For reviews and comprehensive textbooks on project management and scheduling we refer the reader to \cite{Demeulemeester2002-an_WH_ED},
\cite{Demeulemeester2011-eb_WH_ED},
\cite{Hartmann2010-dc_WH_ED},
\cite{Herroelen2007-hh_WH_ED},
\cite{Herroelen2005-me_WH_ED},
\cite{Meredith2011-uu_WH_ED},
\cite{Neumann2002-ve_WH_ED},
\cite{Shtub2004-ob_WH_ED}, and
\cite{Vanhoucke2018-vv_WH_ED}.

\subsection[Revenue management (Arne~K.~Strauss \& Jens~Frische)]{Revenue management\protect\footnote{This subsection was written by Arne~K.~Strauss and Jens~Frische.}}
\label{sec:Revenue_management}

The discipline of revenue management (RM) deals, in the widest sense, with demand management decisions to improve overall revenue or profit. Demand management decisions aim at influencing demand, such as pricing and availability control. Occasionally, demand management decisions can also take different forms like ranking lists (e.g., when showing customers of a meal delivery platform a rank order of restaurants that offer home delivery) or green van icons to denote which time slots for grocery home delivery are more environment-friendly (because there is a delivery planned to take place already anyway). RM is about IT-supported decision-making, mostly on the operational level, in contrast to strategic pricing theory encountered in the marketing domain. 
            
Such decision support systems, referred to as RM systems, have been first developed in the airline industry in the 1970s when deregulation introduced competition in the US airline market. They were so successful that the practice of RM soon spread to other industry domains, particularly to those that sell services or perishable goods (perishability creates pressure to sell within a limited selling horizon).  Examples include restaurant tables, hotel rooms, rental cars, or airline seats, among many other applications. In these industries, the supply is usually fairly inflexible, fixed costs are high and variable costs are relatively low (which also makes revenue maximisation mostly equivalent to profit maximisation, hence the name \emph{revenue} management).

In this section, we first outline recent research trends on demand models (and their estimation) that are required to provide an input to the RM optimisation system. Then, we present recent research on efficient optimisation of demand management decisions. Finally, we outline further reading suggestions including some current popular application areas. We mostly use the passenger airline industry throughout this subsection to illustrate developments in the field of RM.
        
\subsubsection*{Modelling demand}
In order to make good demand management decisions, we first need to  have a model of demand to describe the response to specific RM actions (such as pricing or changing the availability of products or services). The first demand models used in RM assumed that demand for a given product is independent of what else is offered. These so-called `independent demand models' are relatively easy to estimate. However, this independence assumption usually only holds in applications where rate fences (such as advance purchase requirements) strongly limit customers' abilities to substitute a product with one another. 

One way to relax the -- often unrealistic -- independent demand assumption is by considering that the customer looks at all alternatives available and chooses one. The requires modelling of customers’ choice behaviour; the seminal paper by \cite{talluriRevenueManagementGeneral2004_AKSJF} introduced discrete choice modelling to the domain of revenue management. In choice-based RM, demand for a product is assumed to depend on the available purchase alternatives and their attributes. These models tend to be more accurate in predicting demand if the independent demand assumption is not met, at the cost of being more difficult to estimate and implement \citep{kleinReviewRevenueManagement2020_AKSJF}. Much research has been carried out on choice-based RM since 2005; for a recent review, see \cite{straussReviewChoicebasedRevenue2018_AKSJF}.
            
Among the most recent trends -- building on the aforementioned choice-based RM literature with fixed choice model parameters -- is a stream of work on personalisation and choice model parameter learning. For example, \cite{cheungThompsonSamplingOnline2017_AKSJF} solve an online personalised assortment optimisation problem formulated as a multi-armed bandit problem. Demand learning models balance the trade-off between gathering new samples (and thereby learning more about the true customer behaviour) and applying the RM decision that, based on the current belief of customer behaviour, looks to be the best. In the short term, this means that we occasionally make decisions that seem not very promising, yet that will gain us insights into customer behaviour (for instance, by offering price points that were never offered before). For our airline example, a potential application is the ongoing learning of model parameters governing the choice of ancillary products (like seat upgrades, extra luggage, etc.). Models like the one by \cite{agrawalMNLbanditDynamicLearning2019_AKSJF} can be used for this purpose.

Demand models in RM may be biased if they are estimated on constrained data, meaning that the sales data does not necessarily reflect the actual demand. For example, if a flight is fully booked, we observe no further sales transactions for that flight. Yet demand may well exceed flight capacity and, as such, should be estimated somehow. Methods to statistically unconstrained demand data are reviewed in \cite{guoUnconstrainingMethodsRevenue2012_AKSJF}. 

Another dimension of demand modelling is represented by strategic versus myopic customer behaviour. One of the earliest papers on this topic is the work by \cite{suIntertemporalPricingStrategic2007_AKSJF}. He considers customers who may delay their purchase when they expect lower-priced offers in the future. With RM mostly focusing on myopic customers (meaning customers who do not anticipate future developments in their purchase decision), the behaviour of strategic customers leads to inefficiency. \cite{suIntertemporalPricingStrategic2007_AKSJF} proposes an intertemporal pricing component to adjust the offering based on the market composition between these customer types, and more work has built on this since.
            
\subsubsection*{Optimisation advances}
A central element of an RM system is the decision of what to offer whenever a customer arrives. Decision policies (essentially a mapping from the state space of available information to the action space) are used to determine which products are made available (i.e., availability control), or at which price (i.e., dynamic pricing) -- and sometimes, a combination of both. Using dynamic prices to manage demand can be very similar to availability control: when there are products defined with identical features but different price tags such that there is a discrete set of prices to choose from for a product, it can be considered a special instance of the aforementioned availability control \citep{straussReviewChoicebasedRevenue2018_AKSJF}. This can be observed in airline's implementations of differently priced booking classes for the same seat such that a customer can only purchase that seat for the fare of the booking class made available to them. The groundbreaking papers by \cite{gallegoOptimalDynamicPricing1994_AKSJF,gallegoMultiproductDynamicPricing1997_AKSJF} also featured a dynamic pricing concept for a single-commodity and a network-level problem, respectively. Both papers also considered the effect of significant cancellations and no-shows (meaning bookings that are not actually being used in the end). In that context, it can be valuable to accept more reservations than physical capacity would allow. This practice is called overbooking and is widespread in many industries where the risk of having to reject a customer with a valid reservation is not overly costly \citep[with examples such as][showing it being applied already before RM but now usually integrated into systems to manage demand for available capacity]{simonAlmostPracticalSolution1968_AKSJF}. For an overview of recent contributions on this matter, see \cite{kleinReviewRevenueManagement2020_AKSJF}.
            
Decision policies in RM are trading off the immediate reward of having a customer buy a product versus the so-called opportunity cost associated with this purchase, stemming from having to commit some resources to a given sale. For example, a resource might be a flight with a specific seat capacity. Selling a ticket for a seat on this flight requires us to commit a seat to this customer, which otherwise might have still gotten sold in the future at a potentially larger fare. Therefore, by having a customer buy the product, we incur the cost of losing the opportunity to sell the associated resource units in the future. This value (or at least an approximation thereof) is sometimes used as a revenue threshold defining which products shall be shown as available; such special versions of availability control are known as bid price policies, with the bid price being this threshold, and only products with revenues that exceed the bid price being shown.
            
There are two major challenges in deriving optimal decision policies: first, we need to somehow obtain the opportunity costs involved with having a customer book a particular product at a given point in time; second, we need to solve the resulting optimisation problem to give us the actual decision to be implemented. 
            
The latter decision problem, given opportunity cost, may be as simple as a comparison of two numbers (traditionally used in independent-demand settings), but can be non-trivial in the presence of sophisticated models of customer behaviour (dependent demand settings). Much research over the past few years has been devoted to studying properties of choice models that can be exploited to efficiently solve -- or at least closely approximate -- the online RM decision problem. This work is further motivated by the need to solve these RM decision problems quickly to ensure an acceptable user experience. Within availability control, assortment optimisation under various choice models has received particularly much attention because this problem becomes $\mathcal{NP}$-hard for many customer choice models unless a certain structure can be exploited; for a review, see  \cite{straussReviewChoicebasedRevenue2018_AKSJF}. 
            
The other challenge in deriving optimal decision policies is the computation of opportunity costs. This is usually the more difficult task for real applications because the opportunity costs depend on time, the current state (of inventories), and future demand and actions. Dynamic programming (DP; \S\ref{sec:Dynamic_programming}) is usually being applied to solve or at least characterise the optimal decision policy over a given booking horizon. 

However, obtaining opportunity costs using DP is often only possible when dealing with problems that have a single resource \citep[like optimising for a single flight only, for example in][]{wollmerAirlineSeatManagement1992_AKSJF}. When there are products that use more than one resource (like a itinerary of multiple flights connecting in a hub), then we speak of network RM problems. These require much more effort to solve (so as to get opportunity cost estimates for our decision policy) due to the fact that decisions on one product may affect many others that are using the same resource. Therefore, to reach at least an approximate solution for a network RM problem, one usually needs to resort to deterministic linear programming \citep[][\S\ref{sec:Linear_programming}]{liuChoicebasedLinearProgramming2008_AKSJF} or approximate dynamic programming (\citeauthor{gallegoRevenueManagementPricing2019_AKSJF}, \citeyear{gallegoRevenueManagementPricing2019_AKSJF}, describe how approximate dynamic programming can be used in RM). In practical applications, the network-level optimisation problem is often decomposed into a collection of single-resource problems such as in \cite{kemmerDynamicSimultaneousFare2012_AKSJF} who were motivated by methods deployed by Lufthansa Systems in their RM optimisation module. In older RM systems, booking control was typically implemented using versions of the so-called Expected Marginal Seat Revenue (EMSR) heuristics \citep{belobabaAirTravelDemand1987_AKSJF}, which are in turn rooted in the work of \cite{littlewoodSpecialIssuePapers2005_AKSJF}, originally written in 1972.
            
Once a decision policy has been obtained (by first obtaining opportunity costs and then solving the corresponding decision problem), we then need to evaluate the decision policy using simulation or even in real-world trials. \cite{bertsimas2005simulation_VLVV} give an overview of different decision rules for airline RM that are evaluated with a simulation study. Further details on simulation techniques can be found in \S\ref{sec:Simulation}. An example of testing a dynamic pricing policy in live trials is the work by \cite{fisherCompetitionbasedDynamicPricing2018_AKSJF}.
        
\subsubsection*{Further reading}
RM is also applied in retailing, both for e-commerce and offline shopping. \cite{agatz2013revenue_CKTVW} provide a practical overview of ways online retailers implement RM in their business. But there are also retail RM applications outside of online shopping. For example,  \cite{caroClearancePricingOptimization2012_AKSJF} describe how brick-and-mortar fashion stores optimise their price markdowns during season clearing sales.
                     
In particular, linking RM to general transportation problems has received significant attention over recent years. An overview of advances in that field is given by \cite{fleckensteinRecentAdvancesIntegrating2022_AKSJF}. Applications thereof can be seen here especially in business models with delivery constraints, such as same-day deliveries \citep{ulmerDynamicPricingRouting2020_AKSJF} or for attended home delivery (AHD) which is common for online grocery shopping. Customers' desire for short and guaranteed time windows in AHD leads to less than optimal routings. \cite{yang2016choice_CKTVW} show how RM can be used to steer demand for delivery time slots towards a routing solution closer to the optimum, thereby increasing overall profit for businesses shipping goods that require AHD. Another example of applying RM ideas to new transportation problems is the work by \cite{kunnenValueFlexibleFlighttoroute2022_AKSJF}. They analyse how an air traffic network manager could reduce overall delays for all airspace users (i.e., airlines) by offering dynamically priced flight trajectories.
            
For more detailed readings about the development of the RM domain, the techniques being used, and more applications, we refer the reader to the books by \cite{talluriTheoryPracticeRevenue2004a_AKSJF} and \cite{gallegoRevenueManagementPricing2019_AKSJF}.

\subsection[Service industries (Jan~Holmström \& Lauri~Saarinen)]{Service industries\protect\footnote{This subsection was written by Jan~Holmström and Lauri~Saarinen.}}
\label{sec:Service_industries}

\textit{Service industry from the perspective of operations}: Service industry is a concept from economics originally defined by what it is not. It is not a manufacturing industry that produces tangible goods (cars, clothes, equipment), but industry that provides intangible outputs, such as hospitality, healthcare, and education.  In service research, services are also defined by additional characteristics. In addition to intangibility, the so-called IHIP characteristics, recognises heterogeneity, inseparability, and perishability as the defining characteristics of service industries.

\begin{figure}[ht!]
\begin{center}
\includegraphics[width=13.5cm]{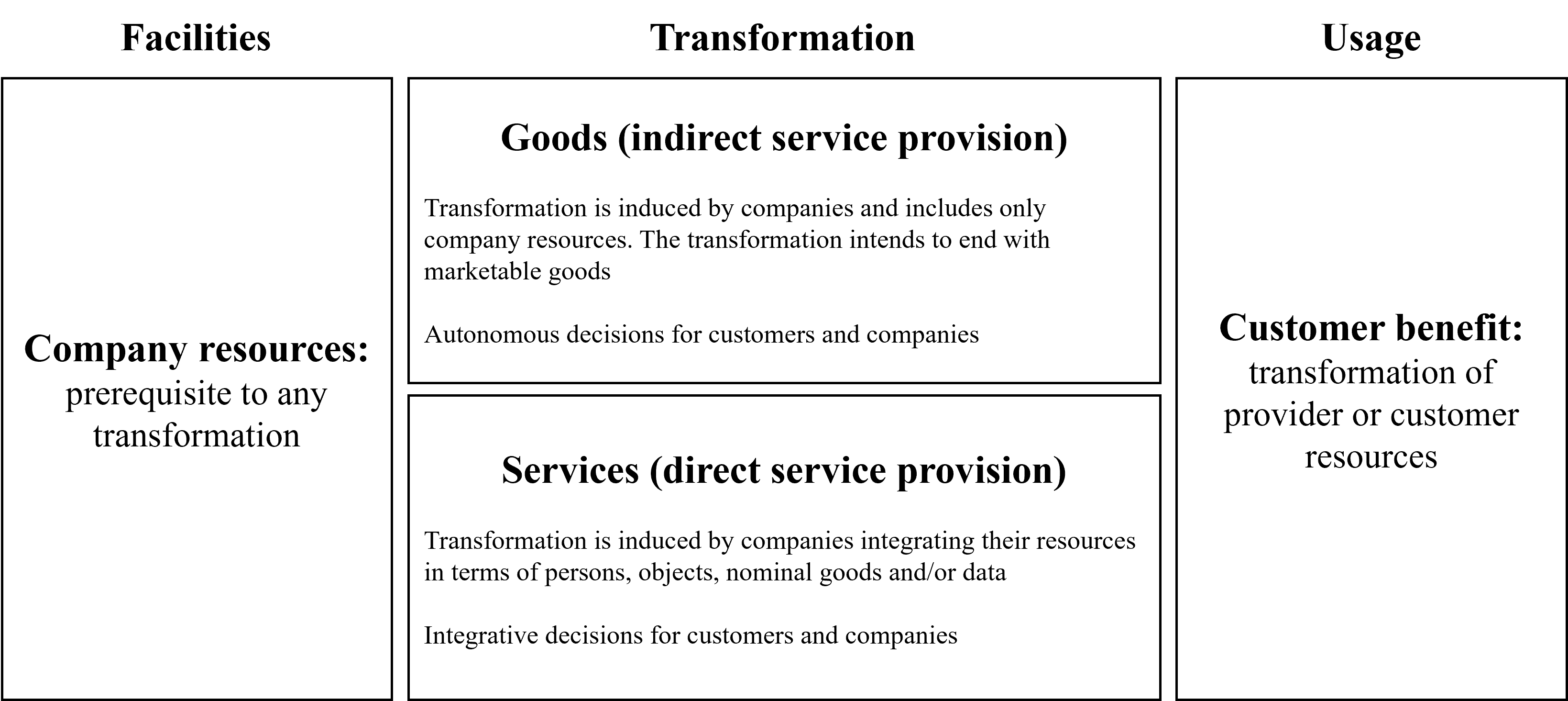} 
    \caption{An actionable framework for service industries: Facilities-Transformation-Usage \citep[adapted from ][]{moeller_characteristics_2010_JHLS}}
\label{fig:FTU_JHLS}
\end{center} 
\end{figure}

In operations and Operational Research, service industry is approached not from its characteristics but operationally to support actionable insights \citep{burger_developing_2019_JHLS}. For Operational Research applications, service industry can be approached through the FTU-framework, defining services as a particular type of transformation (Figure \ref{fig:FTU_JHLS}). Service industries are distinguished from goods industries through the direct provision and integrative decision making. In services the decision making of customers and providing companies is intertwined, while in goods industries customer and providers make autonomous decisions. In service industries the value is directly provided to the customer, while in goods industries indirectly through the product.

However, operational reality is not this clear-cut. In manufacturing industries, through servitisation \citep{kohtamaki_practices_2018_JHLS}, some companies seek to make their products more like services to differentiate their offering, and directly create value to their customers. In service industries, managers seek to make services more like goods, to be able to run service facilities more like factories and improve productivity \citep{levitt_production-line_1972_JHLS, schmenner_service_2004_JHLS}. OR methods that were initially developed in manufacturing industries (e.g., forecasting, queuing, scheduling, simulation), are increasingly applied in service industries to make service facilities operate more like factories \citep[cf. ][]{eveborn_operations_2009_JHLS}. For servitisation, OR presents a more limited range of methods. Methods supporting the servitisation of products are for example, value constellation modelling 
\citep{holmstrom_comparing_2010_JHLS, brax_meta-model_2017_JHLS}, and ecosystem modelling to support innovation and new business model design in an open environment \citep{talmar_mapping_2020_JHLS}.

The challenge in service industries is that service systems tend to be open, problems wicked, and optimising solutions difficult to develop and apply. Value provision often requires interaction with customers (customisation) limiting the situations where facilities can be organised for flow and efficiency as service factories. Also, servitisation of products occurs in an open systems environment, requiring responsiveness to influences and disturbances from the outside, as will be seen for our application examples from industrial services and home healthcare.

\textit{Service industry applications}: In the following we will present two examples on the use of OR methods for creative insight and novel solutions in service provision. The examples are homecare of elderly patients \citep{groop_improving_2017_JHLS}, and line maintenance of commercial aircraft \citep{ohman_frontlog_2021_JHLS}. In the first example systems thinking, in the form of soft OM methods from Theory of Constraints \citep{davies_theory_2005_JHLS}, is used in combination with design science research (implementation and evaluation). In the second example design and simulation are used in combination, uncovering an unexpected new way of simultaneously improving resilience and reducing costs in a commercial airline.

\textit{Homecare of elderly patients} \citep{groop_improving_2017_JHLS}: Nurses, team leaders, and healthcare management had distinct and diverging views on what is the problem with the homecare operations. Strongly held and divergent views are an indication of a possibly wicked problem \citep{Sydelko2021-qp_GM}. The divergent views in the case were uncovered through engagement (following actors in their work) and interviewing, with the purpose of articulating what different stakeholders identified as undesirable effects (UDE) of the current solution.  UDE is thinking tools terminology from Theory of Constraints \citep{dettmer_goldratts_1997_JHLS}. These UDEs of the current operation were pruned for overlaps and narrowed down to a list of seven (including the seeming contradiction between low utilisation of full-time employed nurses, stressed-out nurses, and chronic under-staffing requiring frequent use of temporary staff). Using effect-cause logic, the interconnections between effects, and mechanisms behind the effects were specified and then evaluated by all stakeholder groups in joint workshops. In this case, the first effect-cause analysis pointed towards a core problem, a contradiction, which when addressed, would improve efficiency. The needed change was in the way the nurse visits are scheduled to improve effectiveness. Instead of scheduling nearby patients after each other to save travelling time, the home care organisation should focus on only scheduling nurses for time-critical visits (visits that must be performed at a specific time) at the peak-demand in the morning. This way the time of full-time nurses will be more effectively used.

However, when implemented, the scheduling change had next to no effect. With the initial solution a failure, the evaluation of the implemented change pointed to issues with the problem framing. Going back, considering the stated problems (UDEs) and initial solution, the field researchers found that they had missed an important undesirable effect originating from the way the organisation operated. In the mapping nobody had raised as a problem that full-time nurses, when not busy, do not help-out in other teams. When nurses stay within their teams, work is evenly divided between everybody in the team, which is not a problem for nurses, nor for team leaders. Instead, when there is need for more nurses, outside temporary nurses are called in, and they move between teams if needed, but not the full-time employed nurses.

Management, for whom the low utilisation of full-time employed nurses is a cost issue (with payment of salaries both for idling full-time nurses and busy temporary nurses), were not aware of how nurses staying with their team was a mechanism behind the low utilisation. Nor had the researchers working in the field realised this before failing with the first solution design. Re-framing the problem once more, another solution changing the scheduling for employed nurses was proposed. Instead of dividing work equally between all nurses in their teams, the team leaders should seek to schedule work so that one, or even two nurses in their teams have no work, and can be made available to help-out in other teams. 

\textit{Line-maintenance of commercial aircraft} \citep{ohman_frontlog_2021_JHLS}:  The second example illustrates the use of simulation in problem reframing and finding a new type of solution. The service operations are line-maintenance of aircraft in an airline. Initially the problem was framed by management as improving departure reliability without increasing the number of maintenance technicians. The intended solution was introducing lean in the turn-around of aircraft. 

However, in line maintenance there are no material and time buffers for which lean approaches have been so impactful in manufacturing. The minimum frequency and content of maintenance tasks are regulated. Departures are delayed by technical problems that add unplanned tasks, which need to be carried out. Here, lean principles can increase productivity but not reduce the unplanned tasks. To reduce the delays caused by unplanned tasks a resource buffer of maintenance technicians appears necessary.

In this example, the same method of engagement was applied as in the homecare example. Observing and interviewing different actors involved, field researchers sought to articulate a set of undesirable effects of the current way of operating. However, no agreement on a core problem to address could here be reached. Instead, problem framing ended with a question and a puzzling response that pointed in a new direction. Line-maintenance scheduling assumed that maximising the interval for planned tasks is optimal, also when there are unplanned tasks and constrained resources. Engaging and interviewing maintenance planners for both long-haul and short-haul fleets and operations, production, and resource planning, field researchers began to gain an in-depth understanding of the airline maintenance planning function. Heuristics and principles related to dealing with over-maintenance not visible in the operational documentation were encountered.

To explore the implications, the researchers first modelled the relationship between over-maintenance and planned workload variance in a deterministic setting, focusing solely on scheduled maintenance. The model indicated a promising relationship: an increase of one percent in the total planned workload (over-maintenance) could result in up to a six percent reduction in workload variance. Next, simulation of the airline operation and maintenance included the unplanned events according to their historical distribution. The simulation surprisingly indicated that increasing over-maintenance could reduce over-all costs and improve departure reliability, if combined with a re-scheduling solution for maintenance task. Re-scheduling introduces a new type of time buffer, a frontlog of planned maintenance tasks that can be postponed to allow technicians to address unplanned tasks without disruptions to departure schedule.

\textit{Summary and conclusion}: In service industry applications problem framing methods are particularly important. The openness of service operations and wicked problems often require the Operational Researcher in service industry applications to go outside their comfort zone regarding methods \citep{mingers_soft_2011_JHLS, mingers_helping_2015_JHLS} to search for actionable insight \citep{burger_developing_2019_JHLS}. In the examples provided, a combination of approaches, tools, and methods were contingently employed in the search for a good problem framing as the basis for an effective solution design. For the application of OR methods in service industries, the homecare example illustrates the use of a soft OR method in framing the problem \citep[from the practice of Theory of Constraints, ][]{davies_theory_2005_JHLS}, the use of scheduling from hard OR as a solution component, and implementation as a method of design science for evaluation and redesign \citep{holmstrom_bridging_2009_JHLS, sein_action_2011_JHLS}. The second example illustrates the use of simulation as a method of explorative design. In the empirical grounding of the simulation model we encountered the good problem, which is the key to success in simulation projects \citep{law_how_2003_JHLS}. Through simulation we developed and explored the effect of the dynamic re-scheduling and buffer management approach, with surprising outcomes. Before the simulation study nobody knew about the opportunity to both improve departure reliability and reduce costs.

The example multi-method approach combined soft OR, simulation, and systems thinking for framing the problem. As in cross-agency problem solving in government and public administration, service industry problem solving benefits from mapping different actor perspectives, as the purposes, perspectives and values of the service supply chain actors can easily be in conflict \citep{Sydelko2021-qp_GM}. However, in addition to methods for actionable insight, methods for turning insights into solution proposals are also needed. For proposing and developing a solution design, the two examples relied on explorative design science \citep{holmstrom_bridging_2009_JHLS}, relying on OR methods in evaluation when implementation is possible, and simulation for substituting implementation. In the search of effective solution designs, OR methods such as scheduling, and forecasting were applied as potential solution components in both examples. 

\subsection[Sports (Ian~G.~McHale)]{Sports\protect\footnote{This subsection was written by Ian~G.~McHale.}}
\label{sec:Sports}

Moneyball \citep{Lewis2003-hl_IM} told the story of how the Oakland Athletics Major League Baseball team was able to leverage an inefficiency in the labour market for baseball players, and perform above expectations (given the team’s salary spend). Its impact on how quantitative analysis is viewed within sport and wider society is unprecedented. We have moved from an age when society tended to undervalue quantitative skills to a post-Moneyball era where analytics is generally accepted as being ``cool''. Told in both a best-selling book \citep{Lewis2003-hl_IM}, and a Hollywood movie of the same name, Moneyball has driven a rapid expansion of interest in the field of sports analytics. For an analysis of Moneyball, see, for example, \cite{Hakes2006-tt_IM}.

The history of quantitative analysis in sports dates back to centuries before the Moneyball story, and to the conception of probability itself. The concepts of chance are as old as the first dice games, but they did not evolve into the mathematical principles of probability until the 17\textsuperscript{th} century when Pascal and Fermat exchanged ideas in a series of letters during 1654. The letters were written in response to the following problem: \textit{two players, A and B, each stake 32 pistoles on a first-to-three-point game. When A has 2 points and B has 1 point, the game is interrupted and cannot continue. How should the stakes of 64 pistoles be fairly distributed?} Fast forward three centuries and the similarities of this problem with the problem encountered in limited-overs cricket, when a match is cut short because of rain, are uncanny. Indeed, the solution offered by \cite{Duckworth1998-cu_IM} is one of the great success stories of OR in sport, or arguably OR in any field. That sports fans routinely use the names of a statistician, and an operational researcher should be the source of great pride to the OR community.

The field of sports analytics now boasts specialist journals, regular special issues in top-rated mainstream journals, large departments in sports teams, and many stories of success and over-achievement in professional sport.

\subsubsection*{What is `sports analytics'?}
Analytics is largely an umbrella term for data science, statistics, operational research, and nowadays, machine learning. A simple definition of sports analytics is \textit{the use of analytics to gain a competitive edge in sport}. A wider definition would be \textit{the use of analytics to improve decision making in sport}.

Research has been published on the use of analytics in almost all popular sports including: football, tennis, cricket, golf, American football, baseball, motor sport, martial arts, and many more. Rather than review the field by sport, it is more logical to consider the field by task. The following is not a comprehensive list of such tasks, but provides an overview of the more common objectives of sports analytics.

\subsubsection*{Ranking and Rating}
Ranking of competitors is, to a large extent, the entire purpose of organised sport, and rating is a popular area of research. There are several families of models used to rate competitors. Paired comparisons models are used when two competitors play in each contest. For example, Elo ratings were first developed for use in chess, but have since been used by, for example, \cite{Hvattum2010-nn_IM} for forecasting football results. \cite{Glickman2001-vc_IM} presented a more general Elo model based on a Bayesian updating system and applied it to the problem of dynamic ratings of chess players. Another paired comparisons model is the Bradley-Terry model and this was used in \cite{McHale2011-qs_IM} to forecast the results of tennis matches. Multiple comparisons models are used when several competitors play in each contest and \cite{Baker2015-cv_IM} use a time-varying multiple comparisons model to rate golfers from different eras. \cite{Langville2013-cf_IM} provide an excellent overview of rankings models.

Rating individuals in team sports is a somewhat more complex task than the examples given above, especially when the individuals have different objectives, as is the case in football for example, where some players are mainly responsible for defending, whilst others are mainly responsible for attacking. Basketball, ice-hockey and football all fall into this category. In such circumstances plus-minus ratings are useful. At its most basic level, a player’s plus-minus rating is a comparison of a team’s performances with and without the player. \cite{Rosenbaum2004-ut_IM} presented a method for calculating plus-minus player ratings in basketball, before extensions were added by \cite{Macdonald2012-xc_IM} and \cite{Kharrat2020-gg_IM} to account for the intricacies of ice-hockey and football, respectively.

The availability of more granular data, such as event data (detailing each and every event in a game, e.g., the timing, coordinates and players involved in a pass) and player tracking data (the coordinates of all players on the field of play recorded at several times per second), has enabled more advanced measures of player performance to be calculated. One such measure is that of \textit{expected value of possession} (EVP) for valuing individual actions in team sports. The concept of EVP was first presented in \cite{Cervone2016-gk_IM} and asks the question “what is the probability of the objective happening before an action, compared to the probability after an action?”. The objective may be to score a goal in football. If an action is a good one – the probability of a goal should increase, whilst if it is a bad one, the probability of a goal will likely decrease. The change in the probability of the objective occurring is then the value of the action. Recent applications of deep reinforcement learning have seen EVP calculated for football \citep[e.g.,][]{Liu2020-ps_IM,Fernandez2021-do_IM}. Indeed, it is likely that the EVP concept will be used in many sports in the future.

\cite{Akhtar2015-lj_IM} uses the change in probability of a team winning a Test match to rate cricketers. The idea is similar to the EVP idea proposed by \cite{Cervone2016-gk_IM}, but uses multinomial regression to calculate probabilities of a team winning/drawing/losing the match, before and after each ball of the match.

The idea of monitoring the change in an expected value is also used in golf’s ‘strokes gained’ metric \citep{Broadie2012-eo_IM}. Strokes gained measures how good an individual shot is, and by aggregating over many shots, one can identify how good a player is overall, or how good certain areas (e.g., putting and driving) of a player’s game are.

\subsubsection*{Decision making}
A core tenet of sports analytics is that it should drive improvement, indeed improving decision making is central to the OR paradigm. There are many papers looking to use analytics to improve decision making in a sports context.

Perhaps the costliest consequences of decision making in sport concern the recruitment of athletes. Indeed, the Moneyball premise is built on the idea of avoiding overpaying for talent.

Football clubs exchange huge sums of money to acquire the services of players. These transfer fees were studied in \cite{Coates2022-sm_IM} who consider the issue of the wisdom of the crowd in estimating the fees. \cite{McHale2022-vx_IM} use machine learning techniques to model transfer fees as a function of performance metrics and contract status, amongst other things.

In lucrative team sports such as American football, football and basketball, recruitment of young talent with high potential is of potentially great value, but it appears a relatively little researched area. In one of only a handful of papers on this issue, \cite{Craig2021-fc_IM} present a model to predict the potential of college quarterbacks to one day play in the NFL. 

In addition to making good decisions around player recruitment, sports teams must make good decisions about their coaches. \cite{Peeters2020-yb_IM} consider the impact of coaches on the performance of Major League Baseball teams, whilst \cite{Muehlheusser2018-no_IM} rate coaches in German football. Identifying good coaches is just one dimension of decision-making surrounding running a sports team, and it is often the case that team owners are faced with the decision of whether or not to fire a coach. The impact of managerial dismissals has been the focus of attention in the economics literature. In football, \cite{De_Dios_Tena2007-dp_IM} measure the consequences of mid-season managerial dismissals on a team’s performance and find that there is a short-term improvement in results, but only in home matches.

The final area of decision making we note is that of team selection. In cricket the ordering of the batting line-up was considered in \cite{Perera2016-iw_IM}, whilst \cite{Watson2021-od_IM} use machine learning to optimise team selection in rugby union. \cite{Cao2022-dj_IM} look at optimising team selection in soccer.

\subsubsection*{Other areas of sports analytics}
Sport has attracted the attention of quantitative analysis in numerous other areas, though some do not have the objective of improving performance and/or decision making. For example, OR has been used to inform scheduling of tournaments (see also \S\ref{sec:Timetabling}). 

The popularity of sports betting means forecasting results has received a great deal of attention in the literature. As the sport with the largest global betting market, football has attracted the most attention in the forecasting literature. A notable contribution was that of \cite{Dixon1997-cf_IM}, whose Poisson model has been used as the basis of subsequent work for over two decades. More recently, machine learning techniques have begun to outperform Poisson-type models. See \cite{Dubitzky2019-kc_IM} for details of the results of the ‘Soccer Prediction Challenge’.

Tournament design has been the subject of research in, for example, \cite{Scarf2009-lc_IM}. The idea is that tournaments should maintain excitement. On a similar theme, \cite{Friesl2017-ji_IM} and \cite{Scarf2019-mk_IM} looked at the rules of ice-hockey and rugby and considered how they might be adjusted to increase excitement. By lowering scoring rates, the outcome of a game is more uncertain, and according to the uncertainty of outcome hypothesis this is what drives interest. However, there is conflicting evidence on the uncertainty of outcome hypothesis \citep[see, for example,][]{Forrest2002-rg_IM}. Understanding what drives the interest of fans was the subject of \cite{Buraimo2020-ur_IM} who looked at how suspense, surprise and shock during a match drives in-match television viewing figures.

To find more articles on sports analytics, the interested reader has several options including specialist journals (the \textit{Journal of Quantitative Analysis in Sports}, the \textit{Journal of Sports Economics}, and the \textit{Journal of Sports Analytics}), and discipline journals such as the \textit{European Journal of Operational Research}, the \textit{Journals of the Royal Statistical Society}, the \textit{Journal of the Operational Research Society}, and the \textit{International Journal of Forecasting}, together with a plethora of blogs and websites all focused on sports analytics.

\subsection[Supply chain management (Stephen~M.~Disney)]{Supply chain management\protect\footnote{This subsection was written by Stephen~M.~Disney.}}
\label{sec:Supply_chain_management}

The field of supply chain management (SCM) is concerned with the information, material, and cash flows within and between supply chain members. Materials generally flow down a supply chain (like water in a river); information and money flow up the supply chain. The way we design, source, produce, move, store, schedule, communicate, collaborate, and compete are important factors in SCM.

\subsubsection*{Lean production} 
SCM is built on the foundations of good \emph{industrial engineering}. The pioneering industrial engineers Frank and Lilian Gilbreth provided us with \emph{time and motion studies} \citep{Gilbreth1911_SMD}, \emph{human factors}, and \emph{scientific management}. During the 1920s scientific management techniques were imported into Japan's Imperial Navy's shipyards and factories to improve efficiency and quality. Initially, some table top management games learnt from the Gilbreths in the United States were taken to the Kure Navel Arsenal (a navel shipyard) in 1923 \citep{Robinson1994_SMD}. The table top management games were used demonstrate the efficient flow, and organisation, of work. \cite{Robinson1994_SMD} claims these table top games facilitated Japan in general, and Toyota in particular, to become highly efficient at producing high quality, low cost, reliable products. The \emph{Toyota Production System} (TPS) became the world standard in the highly efficient \emph{lean production} technique \citep{Ohno1988_SMD}. Western companies soon sought to emulate the success of TPS \citep{Womack1990_SMD}, hunting high and low for the \emph{seven lean wastes} \citep{Hines1997_SMD}. \cite{Holweg2007_SMD} provides an excellent summary of the genealogy of lean production.

\subsubsection*{Value stream mapping} 
One of the best ways to document and understand a supply chain is to draw a \emph{value stream map} \citep[VSM;][]{Rother1999_SMD}. VSMs detail how the material flow is controlled by the information flow and decision-making activities.  Key is to determine the point in the material flow where the customer order directly regulates the production cadence. This point is known as the \emph{pacemaker process} or the \emph{customer order decoupling point} \citep[CODP;][]{Olhager2010_SMD}. The pacemaker is often the process that separates the work that is pulled through the system by a \emph{Kanban} system, and the work that flows out to the customer in a \emph{first-in-first-out} (FIFO) queue. 

\subsubsection*{Agile and leagile supply chains}
Lean supply chains are characterised by \emph{just-in-time} inventories and \emph{high capacity utilisation}. But not all supply chains should be lean. Some supply chains need to be responsive, with extra inventory and spare capacity held in reserve so the system can quickly respond to unexpected demand \citep{Fisher1997_SMD}. This has become known as \emph{agile production}. The lean and agile paradigms can integrated in together in a concept known as \emph{leagility} \citep{Naylor1999_SMD}. In leagile supply chains, the material flow is set up to follow lean principles upstream from the CODP; downstream from the CODP, agile principles are followed.

\subsubsection*{Bullwhip and supply chain dynamics}
The \emph{bullwhip effect} is one of the biggest areas of SCM research. The moniker, coined by \cite{Lee1997_XW}, refers to the tendency of the slowly changing consumer demand (the bullwhip handle) to create wildly fluctuating fast moving demand at the raw material processors (the bullwhip popper). This \emph{variance amplification} effect is caused by the decision-making activities  \citep{Forrester1958_SMD}. The seminal paper by \cite{Lee1997_XW} highlights four causes of the bullwhip effect: \emph{demand signal processing}, \emph{order batching}, \emph{shortage gaming}, and \emph{price fluctuations}.

Demand signal processing has been the most studied cause of the bullwhip effect. Demand signal processing refers to the activity of \emph{forecasting} the demand over the lead time (and review period), so that one may determine production and/or replenishment order quantities to maintain \emph{finished goods inventory} and \emph{raw material} levels close to a target. Setting target inventory levels is a problem similar to the \emph{newsvendor problem} \citep{Churchman1957-ax_MYLW}. As orders eventually turn into the inventory here is a \emph{feedback loop} in the decision; there is also a \emph{work-in-progress} feedback loop in the system \citep{Sterman2000_SMD}.  Both these feedback loops contain delays.  This creates a complex system whose dynamics are in part driven by the external demand, but are mostly an internally generated effect caused by the fundamental structure of the supply chain \citep{Sterman2000_SMD}.

Control engineers have developed a large toolkit to understand and manipulate the dynamics of feedback systems. \cite{Towill1982_SMD} and \cite{John1994_SMD} studied the dynamics of continuous time replenishment rules with the \emph{Laplace transform}. \cite{Dejonckheere2003_SMD} studied discrete time replenishment rules via the \emph{z-transform} and the \emph{Fourier transform}. They showed the \emph{order-up-to replenishment policy} with \emph{moving average} and \emph{exponential smoothing} forecasts, for all lead-times and all possible demand patterns, always created bullwhip. 

\cite{Michna2020_SMD} studied \emph{stochastic lead times}, revealing the forecasting of lead times is an important cause of the bullwhip effect. \cite{Gaalman2022_SMD} explores the interaction between the lead time and bullwhip under \emph{general order auto-regressive moving average} demand. They reveal the interaction between demand, lead times, and bullwhip is complex and subtle; bullwhip does not always increase in the lead time. \cite{Wang2016_SMD} provides a recent review of the bullwhip effect, its causes, solution approaches, and thoughts on future research directions.

\subsubsection*{Location and localisation}
The number, and location, of distribution centres (DC) is an important problem in \emph{distribution network design}.  Too few DCs result in longer travel distances (and times) to customers; too many result in high amounts of distribution inventory. The \emph{square root law for inventory}  \citep{Maister1976_SMD} shows the amount of inventory in a distribution network falls by $1\slash\sqrt{n}$ when $n$ DC's are consolidated into a single DC.  The transportation costs involved in delivering customer demand from $n$ DCs can be accurately modelled using \emph{transportation planning software}  \citep{Hammant1999_SMD}. This software typically includes: road maps, speed limits, tolls, congestion, as well as various methods for modelling transport costs.

\emph{Postponement} can also reduce inventory in supply chains. Postponement involves delaying final assembly until demand reveals itself; products are then quickly customised to meet the consumer's desires. For example, HP build generic printers in Mexico to ship to Europe. Upon arrival, they are assigned to a country and the correct power pack is ``assembled'' into the product  \citep{Feitzinger1997_SMD}. With postponement HP holds less generic inventory to buffer the shipping lead time compared to the amount of country specific inventory it would need if the power packs were assembled in Mexico. 

Another important SCM decision is where to produce?  Should you produce locally where perhaps labour cost is high, or should you \emph{outsource}, or \emph{off-shore}, to a low labour cost country? Sometimes, offshore production is supplemented with a local factory or a \emph{near-shored supplier} in a \emph{dual sourcing} arrangement \citep{Allon2010_SMD}. A \emph{tailored base surge} policy sends constant orders to the offshore supplier with the long lead time, while the near-shore supplier flexes production quantities with a short lead time. A small local \emph{SpeedFactory} may be able to correct for the forecast errors and gain enough of an inventory advantage to offset the increased local labour costs \citep{Boute2022_SMD}. 

\subsubsection*{Information flows in supply chains}
Changing the information used in replenishment decisions can improve the dynamics of supply chains. The sharing of demand information with upstream suppliers is often referred to as the \emph{information sharing} \citep{Lee2000_SMD}, or \emph{information enrichment} strategy \citep{Dejonckheere2004_SMD}. Knowing the end consumer demand allows upstream members to base their demand forecasts on the real demand information, removing one of the potential causes of the bullwhip effect. Indeed, information sharing allows for a linear, rather than a geometrically, increasing bullwhip effect as orders go echelon-to-echelon up the supply chain \citep{Chen2000_SMD}. \cite{Kaipia2017_SMD} considers the practicalities of implementing the information sharing strategy.

Sharing both demand and inventory information with your supplier can enable the \emph{vendor managed inventory} (VMI) strategy  \citep{Dong2002_SMD}. In the VMI strategy, the consumer demand and downstream inventory information is used by the supplier (the vendor) to make replenishment  decisions on behalf of his customer. This allows two supply chain echelons to behave dynamically as one echelon, removing a bullwhip generating decision from the supply chain \citep{Holweg2005_SMD}.

\subsubsection*{Coordinating supply chain contracts}
Supply chains often consist of many different organisations, each operating to maximise their own profit. Due to the \emph{double marginalisation} problem, if each player acts solely in their own interests, the supply chain will not be able to reach the \emph{first best} solution; money will be left on the table. Sometimes, the first best solution can be reached by a \emph{centralised decision-maker} coordinating the supply chain; at other times the altruistic behaviour of one supply chain member, in return for a \emph{transfer payment}, can coordinate. 

There are many different types of contracts \citep{Cachon2005_SMD}: revenue sharing, buy-back, price-discount, quantity-flexibility, sales-rebate, franchise, and quantity discount contracts to name just a few. All have their strengths and weakness and are applicable in different settings. Many contracts are based on \emph{newsvendor} principles \citep{Lariviere2016_SMD}. Another important concept in contract design is the idea of \emph{Pareto improving} contracts, where no player is worse off than the (locally optimised) \emph{base case}, but at least one other player is better off. Other contracts allow for the \emph{arbitrary allocation of profits} between players, and for the delegation of decision-making activities to others \citep{Chintapalli2017_SMD}.

\subsubsection*{Emerging topics in the field of supply chain management}
Emerging topics in the SCM field include: 
\begin{itemize}[noitemsep]
    \item The distributed ledger technology behind \emph{block chains} \citep{Babich2020_SMD} and \emph{cryptocurrencies} \citep{Choi2020_SMD} can be used to create a permanent record of provenance and ownership. Ensuring your cotton has not been produced by slaves, your diamonds are not conflict, and children did not mine your Lithium is vital now as UK Directors can face prison time under the Proceeds of Crime Act for crimes committed in their supply chains.
    \item \emph{Opaque pricing} is a technique used to sell last minute travel industry inventory (e.g., hotel rooms) at discounted prices.  The traveller books a room without knowing the exact hotel brand. Cost sensitive travellers are happy because they get a bargain. The hotel is happy because they get extra income without damaging their brand. Opaque pricing can be used for products as well; for example, a red pen sells for \$10, and a blue pen sells for \$10, but if you don't care which colour you have, a red-or-blue pen is offered at the opaque price of \$8. The vendor is able to use the customer's lack of preference to reduce inventory requirements \citep{Ren2022_SMD}.
    \item \emph{Quantum computing} allows one to solve $\mathcal{NP}$-hard problems (such as the \emph{travelling salesman} problem) to optimality instantaneously, rather than waiting for months with regular computers \citep{Srinivasan2018_SMD}. This technology has the potential to make supply chains more efficient. 
\end{itemize}

\subsection[Sustainability (Akshay~Mutha)]{Sustainability\protect\footnote{This subsection was written by Akshay~Mutha.}}
\label{sec:Sustainability}

In this subsection, we focus on the area of sustainable operations from the perspective of \textit{closed-loop supply chains (CLSC)}. We consider literature that focuses on product-, module/part-, and material-level recovery and reuse activities. These activities provide economic and environmental benefits. CLSC entail transportation and acquisition of used products; sorting, grading and disposition for different recovery methods; disassembly and reassembly (i.e., remanufacturing operations); and marketing of remanufactured products. \cite{Guide2003-qd_AM} and \cite{Ferguson2010-pb_AM} provide comprehensive overviews of the strategic, tactical, and operational aspects of CLSC. 

The supply side in CLSC differs from traditional supply chains in the following ways. The quantity of used products being returned is uncertain; the timing of when they are returned is uncertain; and the condition (quality) in which they are returned is also uncertain. These differences lead to uncertain recovery rates and processing times, uncertain cost of recovery, and imperfect matching between supply of used products and demand for remanufactured products, and hence the subsequent demand for new parts needed to make the remanufactured (finished) product. Below, we provide a brief overview of some of the methods used to optimise the different activities in CLSC, while managing these uncertainties. 

The \textit{reverse logistics (RL) network} (see also \S\ref{sec:Logistics}) handles the collection of used products from end-users, and their transportation between collection points, consolidation centres, testing, sorting, and grading facilities, and recovery (e.g., remanufacturing, reuse, recycling) facilities and landfill locations. Stylised and game-theoretic models are developed to determine the optimal collection strategy for producers (if they choose to, or are required to collect used products). The collection strategy includes decisions on whether producers should collect directly from end-users, or use the retail network as collection points, or use third-party collectors \citep[e.g.,][]{Savaskan2006-ou_AM}. In further analysis of the collection strategy, the continuous approximation method is used to determine whether the producer (or business) should offer to pick-up, or have end-users drop-off the used product \citep[e.g.,][]{Fleischmann2003-fx_AM}. 

Several quantitative models are developed to determine the optimal RL network design. An extensive discussion of these models and solution approaches can be found in \cite{Akcali2009-zg_AM}. Linear programming, mixed-integer linear programming (MILP), and stochastic programming are widely used to determine optimal network structures. \cite{Fleischmann2004-wv_AM} provide and excellent overview of MILP and stochastic programming models for facility location and network design for dedicated reverse, and integrated (forward and reverse) logistics networks. Mixed-integer nonlinear programming models are also sometimes used to determine the optimal RL network structure \citep[e.g.,][]{De_Figueiredo2008-rt_AM}. In addition to optimal network design, vehicle routing models (\S\ref{sec:Transportation_Vehicle_routing}) are used to determine optimal collection and pick-up routes. These vehicle routing problems are often $\mathcal{NP}$-hard, and are based on location of demand: node, arc, and general. The models are extended to include vehicle routing with backhaul, routing with simultaneous delivery and pick-up, and routing with partially mixed deliveries and pick-ups \citep[see][and references therein]{Beullens2004-nf_AM}.  

One way of managing the supply uncertainty in CLSC is to forecast the return of used products. Different methods are used to compute the product return probability. These include modelling returns as a function of past sales (via a known delay distribution), regression models \citep{Samorani2019-rw_AM}, simulation models, and queueing models \citep[e.g.,][]{Toktay2000-iy_AM}

Buyers of used products (i.e., producers or their contract-remanufacturers, and third-party remanufacturers) actively manage the supply uncertainty (timing, quantity, and quality) by using incentive mechanisms such as quality-based pricing, trade-ins, and buybacks. Buyers acquire used products either in sorted (i.e., known quality-levels) or unsorted (i.e., unknown quality-levels) form. Lot-sizing models are developed to determine the optimal acquisition quantity when the used products are available in unsorted form \citep[e.g.,][]{Galbreth2009-ar_AM}; and when they are available in sorted (continuous or discrete quality-levels) form \citep[e.g.,][]{Mutha2016-dc_AM}. The acquisition process has also been analysed in the context of buyer-supplier contracts. The objective of these analyses is to determine the optimal contract structure with known or unknown quality levels, e.g., quality-dependant acquisition costs and quantities \citep{Mutha2019-nv_AM}, and coordination mechanisms \citep{Debo2004-xh_AM,Vedantam2021-zj_AM,Li2022-ad_AM}. Several models are also developed to determine the optimal acquisition cost for used products and selling price for remanufactured products for an exogenous set of discrete quality-levels of used products \citep[e.g.,][]{Guide2003-nv_AM}. 

Testing, sorting, and grading of the acquired used products are important activities in product recovery operations. \cite{Hahler2017-qt_AM} provide a detailed description of these operations for used consumer electronics. Sorting defective, and economically and technically infeasible-to-remanufacture units from the acquired quantity streamlines the subsequent operations (e.g., transportation, disassembly, and reassembly). Knowing the quality of the incoming units before scheduling recovery operations significantly improves the performance of the system. The benefit of yield information (i.e., information on the quality distribution of the incoming units) has been analysed using lot-sizing models and simulation models \citep[e.g.,][]{Ketzenberg2003-ek_AM}. Several models are developed to optimise the different decisions in the grading process, e.g., the optimal number of grades \citep[e.g.,][]{Ferguson2009-vt_AM}, the resulting optimal grade-wise remanufacturing cost and selling price \citep[e.g.,][]{Mutha2022-ry_AM}, and the optimal location and timing of the sorting and grading process, for example at the point of collection/return or at the disassembly stage \citep[e.g.,][]{Guide2006-eo_AM,Zikopoulos2008-sq_AM}. 

The disposition of sorted and graded used products typically involve a problem of optimal assignment of the economically and technically recoverable units to different recovery options, e.g., product-level recovery (i.e., remanufacturing); module/part-level recovery (i.e., reuse for making remanufactured and new products, or for spares); and the non-recoverable products for material-level recovery (i.e., recycling). The assignment decisions are usually based on considerations of supply (yield information, processing times, and costs) and demand (revenue, opportunity cost, and inventory cost). Optimal control models \citep[e.g.,][]{Inderfurth2001-ye_AM} and revenue management-based models \citep[e.g.,][]{Pince2016-ol_AM,Calmon2017-gn_AM,Calmon2021-vf_AM} are widely used to determine optimal disposition decisions. Depending on the type of the product, single-period models (for products with short lifecycles, e.g., cellphones) and multiperiod models (for products with long lifecycles, e.g., engines) are used in the disposition analyses. For example, \cite{Ozdemir-Akyildirim2014-by_AM} formulate the optimisation problem as a multiperiod Markov decision process (MDP) and provide a linear-programming model for solving the deterministic approximation of the MDP model. 

Within the production planning and control literature in CLSC, a relatively small part has focused on disassembly planning and sequencing, and material requirement planning (MRP). \cite{Inderfurth2004-qp_AM} provide an extensive overview of the various optimisation models developed to optimise these elements, including shop floor control rules, in remanufacturing-only and hybrid (joint manufacturing and remanufacturing) systems. Disassembly sequencing is mainly analysed using direct graphs \citep[see][and \S\ref{sec:Graphs_and_networks}]{Lambert2003-au_AM}, and MRP decisions are analysed from an inventory control perspective \citep[e.g.,][]{Inderfurth1999-kv_AM,Ferrer2001-yn_AM}. A significant part of the literature on CLSC has focused on inventory management. Optimal inventory control policies are derived using periodic-review models \citep[e.g.,][]{R_Teunter2004-cp_AM,Zhou2011-zd_AM} and continuous-review models \citep[e.g.,][]{Van_der_Laan1999-kp_AM,Toktay2000-iy_AM,Jia2016-vy_AM}. The single-period newsvendor-like models are largely analysed as acquisition lot-sizing models (discussed in the preceding paragraphs). 

The research on market (selling)-related aspects of CLSC is focused around understanding the profit and pricing implications due to the co-existence of new and remanufactured products in the (same) market, and on understanding the customers of remanufactured products. Optimisation models are developed to determine the pricing and profitability of remanufactured products \citep[e.g.,][]{Ovchinnikov2011-qd_AM,Abbey2015-dl_AM}. Game-theoretic models are developed to determine optimal market-segments (based on pricing) for new and remanufactured products \citep[e.g.,][]{Debo2005-vu_AM,Atasu2008-za_AM}. The market for- and customers of- remanufactured products are mostly analysed using empirical methods, e.g., using sales data from websites selling used and remanufactured products, usually accompanied by customer surveys \citep[e.g.,][]{Guide2010-ey_AM,Subramanian2012-om_AM}. Behavioural experiments are used to understand consumer perceptions (e.g., quality, functionality), their acceptance (and rejection), and willingness to pay for remanufactured products \citep[e.g.,][and references therein]{Abbey2017-cn_AM}.

\subsection[Telecommunications (Bernard~Fortz)]{Telecommunications\protect\footnote{This subsection was written by Bernard~Fortz.}}
\label{sec:Telecommunications}

Operational Research plays a key role in the design and management of telecommunication networks. A large variety of applications of both exact methods and heuristics can be found in the literature. We focus here on the applications for wired networks. 

\subsubsection*{Topological network design}
The earliest works on telecommunication networks focused on wired fixed line telephony. For the long-term planning of these networks, clients' demands are not known in advance, or with a lot of uncertainty.  This often gives rise to two-stage approaches where only the fixed cost of opening links are considered first, and the decisions on routing and capacity allocation taken in a second (later) stage. This approach is relevant when the fixed costs are very high compared to routing and capacity costs, and/or when topological decisions do not affect capacity decisions too much . For example, digging a trench to install fiber optic cables is very costly, while increasing capacity can be done by adding or upgrading equipment into nodes, which is relatively simple and cheap. The objective is to build a network at minimum cost, considering only the fixed cost associated with opening a link, ignoring capacity and routing costs.

Two main issues appear in the planning process of such networks: economy and survivability.  Economy refers to the construction cost, while survivability refers to the restoration of services in the event of equipment failure. A network is called a tree if it is connected (i.e., there exists a path between all pairs of nodes), and removing any link disconnects at least one pair of nodes. Trees satisfy the primary goal of minimising the total cost while connecting all nodes. The minimum cost spanning tree problem therefore received a lot of attention, see e.g.,\citet{magnanti.wolsey:95_BF}.

However, only one node or edge breakdown causes a tree network to become disconnected and therefore to fail in its main objective of enabling communication between all pairs of nodes. This means that some survivability constraints have to be considered while building the network. Usually, these constraints come in the form with $k$-connectivity requirements, i.e. the ability to restore network service in the event of a failure of at most $k-1$ components of the network. In their earliest work on the subject, \cite{grotschel.monma:90_BF} introduced a general model for survivability requirements, and studied the polytope associated with an integer programming formulation of the problem. 

The minimum-cost two-connected spanning network problem, that consists in finding a network with minimal total cost for which two node-disjoint paths are available between every pair of nodes, was studied extensively, starting with the work of \citet{monma.shallcross:89_BF}. Such networks have been found to provide a sufficient level of survivability in most cases, but it turns out that the optimal solution of this problem is often very sparse. In such a topology, primary routing paths and re-routing paths in case of failure might become very long, introducing large delays in the network.

Two kinds of solutions have been proposed to remedy this problem: The first one imposes a constraint on the length of the paths (in terms of number of links crossed), the so-called \emph{hop-constrained} models. The second approach consists of imposing that each edge belongs to at least one cycle (or {\em ring}) whose length is bounded by a given constant.

Hop-constraints were first considered by \cite{balakrishnan.altinkemer:92_BF} in order  to generate alternative solutions for a network design problem. Later on, \citet{gouveia:98_BF} presented a layered network flow reformulation that has since been used in many network design applications involving hops-constraints.

The second approach to avoid long re-routing paths in case of failure is based on the technology of {\em self-healing rings}. These are cycles in the network equipped in such a way that any link failure in the ring is automatically detected and the traffic rerouted by the alternative path in the cycle. Many problems involve setting a bound on the length of the ring including each edge. Network design problems with bounded rings were first studied in \citet{fortz.labbe.ea:00_BF}.

\subsubsection*{Location problems}

Location problems play a central role in telecommunications network design. We focus here on problems arising in wired (optical) telecommunications networks. These problems are mostly concerned with decisions related to the placement of specific equipments into nodes of the network, and are closely related to hub location problems \citep{alumur.kara:08_BF}.

The \emph{Concentrator Location Problem} is probably the most basic application of equipment placement. The problem consists of determining the number and location of concentrators that are used to aggregate end-user demands before sending them on the backbone network. The allocation of end-users demands to the concentrators has also to be determined such that the capacities of the concentrators are not exceeded. This problem has received much attention in the literature, starting with the work of 
\citet{pirkul:87_BF}.

Another classical problem arises with the replacement of an old technology by a new one, e.g., when telecommunications companies replace outdated copper twisted cable connections by fiber optic connections. The \emph{Connected Facility Location Problem} (ConFL) aims at optimising the building cost for networks involving the two technologies, which is modeled as \emph{tree-star} networks: the core network, made of fiber optic connections, has a tree topology and interconnects multiplexers that switch traffic between fiber optic and copper connections. Each multiplexer is the centre of a star-network of copper connections to the customers. Early work on ConFL concentrated on approximation algorithms, such as the primal-dual procedures proposed by \cite{swamy.kumar:04_BF}. The currently best-known constant approximation ratio is given by the 4-approximation algorithm of \cite{eisenbrand.grandoni.ea:10_BF}. Heuristic approaches have been proposed by \cite{ljubic:07_BF} and \cite{bardossy.raghavan:10_BF}. Different Mixed Integer Programming models for ConFL were proposed by \cite{gollowitzer.ljubic:11_BF}. 

In addition to these long-term design problems, operational short-term decisions are related to the routing of demands in the network, with a focus on avoiding congestion. Most networks nowadays operate the \emph{Internet Protocol}. The internet is a collection of inter-connected networks called autonomous systems, that operates under a hierarchy of layers. An \emph{Autonomous System} (AS) is defined as a set of routers under a single
technical administration, such as an internet service provider or a country. As of July 2022, over $100,000$ ASes were registered\footnote{\url{https://www-public.imtbs-tsp.eu/~maigron/RIR_Stats/RIR_Delegations/World/ASN-ByNb.html}}, connecting over 5 billion internet users worldwide\footnote{\url{https://www.statista.com/statistics/617136/digital-population-worldwide/}}. 

\subsubsection*{Traffic engineering}
Traffic engineering (TE) addresses the problem of efficiently allocating resources in the network so that user constraints are met. Several criteria can be used to measure the effectiveness of a routing configuration. The selection of the objective function may drastically change the quality of the resulting routing. This distinction has been illustrated in \citet{pioro.medhi:04_BF}. \citet{balon.skivee.ea:06_BF} discuss various TE objective functions and evaluate how well these objective functions meet TE requirements.

The internet routing protocols can be clustered into two main groups: \emph{inter-domain} and \emph{intra-domain}. While inter-domain are used to route traffic between ASes, \emph{Interior Gateway Protocols} (IGPs) handle the routing within ASes. As inter-domain protocols are mostly governed by administrative and political considerations, there is not much room for Operational Research techniques to be applied for performing TE. On the other hand, the optimisation of IGPs have received a lot of attention. The most popular IGPs
are based on shortest path routing: shortest paths are calculated using a \emph{link metric system}, which corresponds to the set of link weights or link metrics that belong to the same AS. The network operator controls the routing of the traffic indirectly by setting the link metrics. This gives rise to very challenging optimisation that have mostly been tackled heuristically by many authors, starting with the seminal work of \citet{fortz.thorup:00_BF}. Some exact models have also been proposed, e.g., by \citet{pioro.szentesi.ea:00_BF}.

Recently, \cite{filsfils.kumar-nainar.ea:15_BF} proposed Segment Routing (SR), a new routing protocol developed to address known limitations of traditional routing protocols in IP networks. SR offers the possibility to deviate from the shortest path by using detours in the form of nodes or links respectively called node segments and adjacency segments. Optimisation of SR is a very active field of research and has been already addressed in
\cite{bhatia.hao.ea:15_BF,hartert.vissicchio.ea:15_BF,jadin.aubry.ea:19_BF}.

\subsubsection*{Further readings}
For surveys on survivable network design, we refer the reader to \cite{christofides.whitlock:81_BF,kerivin.mahjoub:05_BF,fortz.labbe:06_BF,fortz:21_BF}. Location problems in telecommunications are surveyed in \cite{skorin-kapov.skorin-kapov.ea:06_BF,fortz:15_BF} and a unified view on location and network design problems was proposed by \cite{contreras.fernandez:12_BF}. For a detailed survey on the Concentrator Location Problem, see Chapter~2 in \cite{yaman:05*1_BF}. Traffic engineering with shortest paths routing protocols is covered in the surveys of \cite{bley.fortz.ea:09_BF,fortz:11_BF,altn.fortz.ea:13_BF}.

\subsection[Timetabling (Greet~Vanden~Berghe \& Sanja~Petrovic)]{Timetabling\protect\footnote{This subsection was written by Greet~Vanden~Berghe and Sanja~Petrovic.}}
\label{sec:Timetabling}

Timetabling  represents a particular subgroup of scheduling problems, namely the set of problems for which activities must be assigned to resources within a set of fixed timeslots. Nevertheless, the two disciplines, scheduling and timetabling, are tightly related and benefit from mutual advancements in both modelling and method development. 

Practical timetabling problems appear in many sectors, for example, in education, healthcare, sports and public transportation. They have been drawing academic attention for a few decades, partly because they are easy to grasp but challenging to solve. The timetabling community gathered at its first international conference in Edinburgh in 1995, one year before the Association of European Operational Research Societies (EURO) established a EURO Working Group on the Practice and Theory of Automated Timetabling \citep{EWG_PATAT_GVBSP}.  Ever since the third conference, which took place in 2000, the timetabling community has gathered every two years\footnote{The 25th anniversary of PATAT conferences had to be postponed till 2022, due to the COVID-19 pandemic.} to share ideas on both theoretical and practical aspects of timetabling.

This subsection provides a brief overview of timetabling history, while highlighting what makes timetabling problems computationally challenging, which initiatives have boosted timetabling research and how state-of-the-art knowledge, models and algorithms can be applied in practice. We restrict the discussion to timetabling problems involving human resources, such as students, teachers, healthcare workers and sports teams. 

\subsubsection*{Problem definition}
Let us consider a set of timeslots $T = \{1,...,|T|\}$, a set of activities $A = \{1, ..., |A|\}$ and a set of resources $R = \{1,...,|R|\}$. A timetabling problem then consists in assigning (all) the activities in $A$ to resources in $R$ and timeslots in $T$ in such a way that a set of constraints is met.  Constraints may apply to resources, timeslots and activities.  They usually restrict the number of assignments to certain resources within subsets of $T$. 

Constraints are usually divided into two categories: hard constraints, which must be strictly satisfied, and soft constraints, for which violations may be tolerated but should be avoided if possible. Weights may be set on the soft constraints, denoting their relative importance. A common timetabling objective is to minimise the weighted sum of soft constraint violations.  This objective sometimes has to be combined with other timetabling objectives, for example, to minimise the cost associated with the employed resources. 

\subsubsection*{Educational timetabling} 
Educational timetabling problems can be split into three major groups: university examination timetabling, university course timetabling and high-school timetabling. In {\it examination timetabling}, the task is to assign examinations in $A$ to a limited number of timeslots in $T$ and rooms in $R$ such that no student has more than one exam at a time.  Each student's exams should be spread out in time as much as possible. Additional constraints may include precedence constraints between exams, special room requirements, and limited room capacities. {\it Course timetabling} involves the assignment of course sections  (lectures, tutorials, lab sessions, seminars) to specific days of the week and times of the day. Real-world problems may require sectioning, when students have to be split into separate subgroups for different sections. Typically, the objective is to minimise the number of students' conflicts. {\it High-school timetabling} assumes that students are split into classes and each class has to take a set of resources. Given a set of timeslots, each activity (involving both a student group and a teacher) must be assigned to a timeslot so that no teacher and no student group are participating in more than one activity at a time. Most practical problems have additional constraints;  for example, teachers may have limited  availability and some activities may require more than one timeslot. In general, educational timetabling problems are $\mathcal{NP}$-hard \citep{DeWerra_GVBSP}. Additionally, the constraints often pose a feasibility challenge. 

The educational timetabling community made a considerable effort to create rich sets of benchmark instances to be used for comparing methods. The first set of examination timetabling instances was defined by \citet{carter1996_GVBSP}.  Four competitions on educational timetabling, entitled  ITC-2002 \citep{ITC-2002_GVBSP}, ITC-2007 \citep{McCollumEtal2010_GVBSP}, ITC-2011 \citep{post2016_GVBSP} and ITC-2019 \citep{muller_GVBSP}, further advanced the development of timetabling algorithms. \cite{post2012_GVBSP} developed a general format and benchmark instances for high-school timetabling, which were extended later by \cite{post2014_GVBSP}. \citet{CESCHIA2022_GVBSP} published a review of educational timetabling, presenting detailed characteristics of all benchmark instances and state-of-the-art results. \textsc{OptHub}\footnote{\hyperlink{https://OPTHUB.uniud.it}{https://OPTHUB.uniud.it}} provides a common platform for storing problem instances and solutions to selected optimisation problems, including educational timetabling. 

\subsubsection*{Personnel timetabling}
Personnel timetabling, also referred to as employee timetabling or rostering, concerns the construction of a timetable for personnel in $R$ in such a way as to satisfy coverage constraints throughout a time horizon \citep{Ernst_GVBSP}. The timeslots in $T$ often represent shifts, which correspond to tasks or duties in $A$. Some activities may require certain skills, and hence can only be conducted by a subset of $R$. Many {\it work-rest}-related objectives are formulated in terms of {\it time-related constraints}, restricting, for example, the number of hours worked, the number of weekends worked, the number of consecutive night shifts \citep{JOS_GVBSP}. Additionally, personnel rostering problems typically consider {\it personal preferences} as regards working time or days off. Whereas the problem is generally considered $\mathcal{NP}$-hard, \cite{smet2016polynomially_GVBSP} showed that some personnel rostering problems are polynomially solvable, provided they do not contain a particular class of  constraints. \cite{categorisation_GVBSP} developed a categorisation of personnel rostering problems, based on the characterisation of resources, objectives and constraints.  \cite{kingston2018unifed_GVBSP} complemented this work by providing a unified notation for nurse rostering problems. 

The Practice and Theory of Automated Timetabling (PATAT) community organised two International Nurse Rostering Competitions, entitled INRC I and INRC II. The problem definition of INRC I \citep{INRCI_GVBSP} was based on the instances collected by \cite{curtois_instances_GVBSP}. INRC II \citep{ceschia2019second_GVBSP} incorporated real-world constraints concerning subsequent rostering horizons. The competition datasets have been collected and published\footnote{\url{https://patatconference.org/}}.

Apart from the constraints and objective functions considered in the two INRCs, some sectors expect their personnel rosters to be cyclic \citep{musliu_GVBSP,cyclic_GVBSP}. Recent trends also include  objectives related to  fairness \citep{fairness_GVBSP} and well-being \citep{well-being_GVBSP}.  Objective priorities set by the users may lead to unwanted solutions. To address this issue, \cite{weights_GVBSP} developed an approach to automatically set acceptable weights which avoid conflicting objectives from leading to poor solutions. 

\subsubsection*{Sports timetabling}
Sports timetabling problems often address tournament or competition scheduling. They require assigning sports activities in $A$, represented by {\it pairs} of teams in $R$, to timeslots in $T$ in such a way that each team meets all the other teams. Constraints depend on the competition's rules, which may differ in different  parts of the world \citep{ribeiro_GVBSP,duran2021sports_GVBSP}. Specific sports timetabling constraints prescribe that teams must not meet the same opponent within consecutive timeslots, or that the number of consecutive home or away games is restricted. The travelling tournament problem (TTP), introduced by \cite{Easton2001_GVBSP}, is an academic adaptation of the Major League Baseball competition in the United States.  The objective of the TTP is to minimise the sum of travelling distances for each  team.  Travelling umpire scheduling \citep{TUP_GVBSP} is subject to similar constraints, but it assumes that the tournament is fixed and that each game is assigned an umpire.

\cite{Rasmussen2008_GVBSP} provided a review on round robin sports timetabling, where each team plays against each other team twice, once at home and once away. \citeauthor{Drexl2007_GVBSP}'s (\citeyear{Drexl2007_GVBSP}) review focused on graph-theoretical approaches to sports timetabling.  \cite{Briskorn2010_GVBSP} investigated the complexity of several variants of the round-robin tournament problem, and similarly, \cite{DEOLIVEIRA2015101_GVBSP} studied the complexity of travelling umpire scheduling problems. The characteristic sports timetabling constraints, which forbid the assignment of activities to subsets of $T$, can be  challenging in terms of feasibility. 

\cite{TTP_GVBSP} and \cite{TUP_Website_GVBSP} boosted sports timetabling research by publishing challenging benchmark instances and monitoring best known and/or optimal results. \cite{sports_competition_GVBSP} organised the first international sports timetabling competition, for which the instances are available at the website of \cite{UGENT_STT_competition_GVBSP}.

\subsubsection*{Timetabling and related problems}
Academic timetabling problems are often considered in isolation from other problems. However, many real-world situations face timetabling entangled with other optimisation problems. Solutions for one of them have an impact on the solution for the other problems. For example, the {\it staffing} problem is concerned with optimising a group of human resources and their characteristics such as skills and contracts in an organisation, across a relatively large time horizon.  From a staffing perspective, the personnel structure should adequately cover the organisation’s anticipated workload while respecting the available budget. On the other hand, from a rostering perspective, the personnel structure should enable computing good quality rosters across many subsequent rostering periods \citep{Komarudin_GVBSP}. Similarly, {\it task scheduling} usually assumes personnel rosters are fixed, but both problems can also be addressed in an integrated manner \citep{paul2015classification_GVBSP}. The {\it workforce routing and scheduling} problem is related to vehicle routing. Apart from scheduling a fleet of vehicles to serve a set of customers, timetabling issues, such as temporal constraints, contracts and skills are also imposed on the problem \citep{Dario_GVBSP}. Some {\it production scheduling and inventory} problems are subject to additional timetabling restrictions which apply to their employees \citep{altachem_GVBSP}.

\subsubsection*{Where do we stand and what is the future}
Academic timetabling has made good progress and instances, models and algorithms have been shared and published.  For example, the heuristic search strategies Step  Counting Hill-climbing \citep{step_GVBSP} and Late Acceptance Hill-climbing \citep{LA_GVBSP} were initially developed for solving timetabling problems.  Due to their simplicity and effectiveness, they continue to be used in a much wider application domain by many computational experts. 

So long as some instances remain unsolved, or solutions for instances have not been proven optimal, algorithm development  remains open for improvement. Future challenges may also apply to new combinatorial optimisation problems encompassing a timetabling component. They may not necessarily map to any of the three timetabling categories detailed in this chapter. However, they may gain importance due to either increased practical need or academic initiatives, such as the publication of benchmarks or the organisation of competitions.

Apart from these future computational challenges, timetabling research should also focus on how to address {\it human resources'} considerations. Besides the traditional work-rest constraints and objectives, academia should also reconcile personnel well-being with their perception of fair workload within a team and with their level of autonomy in determining their personal timetables.  Research should also focus on how to address the increasing personnel resignation in human-centric working environments such as education and healthcare. Robust timetabling, for example, has a lot of potential and at the same time induces scientifically interesting modelling questions.

\subsection[Transportation: Rail (David~Canca)]{Transportation: Rail\protect\footnote{This subsection was written by David~Canca.}}
\label{sec:Transportation_Rail}

The transportation of goods and passengers by rail has played an important role in the evolution of industrialised societies, contributing to their development and prosperity. Rail freight transport still holds critical importance in supporting the economic growth of many countries around the world due to its contribution to guaranteeing an efficient flow of goods internally and across borders. Furthermore, rail transportation is also essential for the movement of people, being the preferred transportation mode for commuters in many large urban areas. This preponderant role also affects the internal mobility of cities. First, a differentiation must be made between freight and passenger transport. Freight trains are longer and heavier than passenger trains, and can often have multiple propulsion units. Compared to that, passenger trains are much lighter and have more horsepower per tonne. There are also important planning and operational differences, whereas passengers decide freely where they will travel, each load of freight must be managed and routed from a specific origin to its destination. These differences originate very different problems in both areas. Even in passenger transportation, different problems arise depending on the type of service; long- and medium-distance, commuter rail, urban rapid transit, and scenic and sightseeing train transportation; see, for instance, \citet{Caprara2007_DC}.

Despite all these differences, a set of common hierarchical stages can be highlighted in the process of planning and operating a rail transportation system  \citep{Bussieck1997a_DC}: network design and/or line planning, timetabling, platforming, rolling stock circulation, shunting, and crew planning. 

At a strategic level, the problems are characterised by long planning horizons and typically involve resource acquisition. This level includes network design and line planning problems. The first  refers to the construction or modification of existing railway infrastructure and mainly concerns urban rapid transit systems.  For a railway company or agency, the line planning problem consists of defining a set of lines and determining their frequencies, and it is usually the first stage in planning medium and long-distance passenger rail networks. 

\citet{Bussieck2004_DC} considered the design of line plans in public transport with the objective of minimising the total cost.  \citet{Goossens2006_DC} presented several models for solving line planning problems in which lines can have different halting patterns. \citet{Laporte2007_DC} proposed a first railway rapid transit network design model to maximise the expected trip coverage. \citet{GutierrezJarpa2013_DC} presented a model to minimise travel cost while maximising the captured demand. See also \cite{Laporte2015_DC} for an extension where  the idea consists of first building a set of segments within broad corridors connecting some vertex sets to later assemble the segments into lines.

A different set of works pays attention to the formulation of network design models from scratch. Starting from an underlying network, these models construct lines by joining edges, incorporating topological constraints to guarantee connectivity between consecutive edges of each line.  This approach gives rise to complex models which are quite difficult to solve using exact procedures; see, for instance, the work by \citet{Szeto2014_DC}, or the recent works by \citet{Canca2017_DC} and \citet{Canca2019b_DC} which concern the design of a railway rapid transit network.

For a comprehensive review of the different methodologies used in practice to solve this problem, the readers can consult the review of \citet{Guihaire2008_DC}. More recent reviews of \citet{Shoebel2012_DC} and \citet{IbarraRojas2015_DC} present a systematic classification of problem variants, considered objectives and solving methodologies. 

At the tactical level, the next stage in planning a railway system consists of several problems, starting with scheduling and timetabling, followed by rolling stock planning, crew rostering, and crew scheduling. The timetabling problem concerns the determination of the arrival and departure times of trains to stations. When overtaking and overlapping are allowed, the timetabling problem becomes a train scheduling problem. Timetables can be cyclic, regular, hybrid, and demand-driven. 
Concerning the design of cyclic timetables, \citet{Caprara2002_DC} proposed a graph-theoretic formulation for the train timetabling problem using a directed multigraph in which nodes correspond to departures and arrivals at a certain station at a given time instant.  \citet{Liebchen2002_DC} used a Periodic Event Scheduling problem (PESP) with several add-ons concerning problem reduction and strengthening. \citet{Chierici2004_DC} extended the classical timetabling model to take into account the reciprocal influence between the quality of a timetable and the transport demand captured by the railway with respect to alternative means of transport. \citet{Cacchiani2008_DC} proposed heuristic and exact algorithms for the (periodic and non-periodic) train timetabling problem on a corridor.

Regular timetables have been commonly used in the case of railway rapid transit systems, especially at relatively short time planning horizons where demand can be considered approximately constant.  \citet{Canca2016a_DC} proposed a sequential optimisation approach to determine the best regular timetable for a railway rapid transit network where lines share tracks. \citet{Canca2017a_DC} incorporated aspects of energy consumption in the design of a two-way rapid rail transit line. Later, \citet{Canca2018_DC} extended the previous work to a full network, taking into account transfers between lines. \citet{Robenek2017_DC} proposed a new type of timetable combining both the regularity of the cyclic timetables and the flexibility of the non-cyclic ones. 

During recent years, starting from the works of \citet{Canca2014_DC} and \citet{Niu2013_DC} many researchers have paid attention to the design of demand-driven timetables \citep[see, for instance,  ][]{Barrena2014a_DC, Barrena2014_DC}. The design of a specific train timetable can be combined by using different acceleration strategies such as stop-skipping and short-turning. For example, given predetermined train skip-stop patterns, \citet{Niu2015_DC} proposed a quadratic integer programming model with linear constraints to synchronise effective passenger loading and train arrival and departure times at stations. \citet{Zhou2022_DC} proposed a mixed integer linear programming model to jointly optimise the train timetable and the rolling stock circulation plan, allowing rolling stock to change its composition through coupling/decoupling operations at the terminal stations of a metro line. \cite{Yuan2022_DC} introduced a new integrated optimisation model for train timetabling that also considered rolling stock assignment and incorporated a short turn strategy on a bidirectional metro line. 

Several authors have also proposed methods to increase the transport capacity of a given timetable \citep[see ][]{Burdett2009_DC}. \citet{Cacchiani2010_DC} studied the problem of incorporating freight transport trains in railway networks, where both passenger and freight trains are running. To finish this description of the OR contributions for the train timetabling problem, a special mention of the work by \citet{Kroon2009_DC} is convenient. In this research, the authors generated several real timetables using Operational Research techniques for the Dutch railway network.

Rolling stock management is probably the most complex stage in the classical sequential railway planning process and plays a key role in the efficient operation of railway networks. At a tactical level, the rolling stock circulation plan consists of a set of interrelated subproblems such as train composition decisions (coupling and decoupling operations involving locomotives and carriages), selection of rest locations, the design of vehicle circulations (specific paths that vehicles must follow to guarantee an efficient and safe operation), and the definition of maintenance policies \citep{Caprara2007_DC}. In a general rolling stock circulation problem, every train circulation has a variable length (distance and number of days) and incorporates information about the allowed specific rolling stock types, composition, coupling/decoupling operations, maintenance and cleaning activities. \citep{Maroti2005_DC, Maroti2007_DC}. Other practical considerations such as rolling stock availability, depot capacity \citep{Lai2015_DC}, coupling and decoupling activities \citep{Fioole2006_DC}, turnaround times, maintenance \citep{Maroti2007_DC}, and track and platform capacities are simultaneously considered depending on the specific problem. Given the importance of this topic within the set of planning tasks, other contributions have been proposed for different problems concerning rolling stock management, as, for instance, determining a set of minimum cost equipment cycles such that the most convenient rolling material is assigned to each planned trip \cite{Cordeau2000_DC} or obtaining the optimal circulation of rolling stock considering order in train compositions \citep{Alfieri2006_DC, Peeters2008_DC}. Maintenance also plays an important role in several rolling stock management contributions, see, for instance, the works by \citet{Maroti2005_DC, Giacco2014_DC} and \citet{Dariano2019_DC}. Robustness is another topic of interest in the related literature. Interested readers can consult the works by \citet{Cacchiani2008a_DC, Cacchiani2012_DC}.

After rolling stock management, the crew scheduling process determines the set of duties that covers all programmed services \citep{Caprara1998_DC}. Finally, the crew is assigned to serve the crew schedule and the corresponding train services \citep{Huisman2005_DC}. The rostering process aims at determining an optimal sequencing of a given set of duties into rosters satisfying operational constraints deriving from union contract and company regulations \citep{Caprara2003_DC}. 

To finish this section, two important problems of rail freight transportation are briefly commented. The first concerns the strategic design of freight transport networks and the second concerns the tactical operation of marshalling yards. Concerning the design of service networks, \citet{Crainic1984_DC} analysed the problems of routing freight traffic, scheduling train services, and allocating classification activities between yards on a rail network. \cite{Crainic1990_DC} developed a model of rail freight transportation adapted for the strategic planning of freight traffic considering other transportation modes. \citet{Zhu2014_DC} addressed the problem of scheduled service network design for freight rail transportation integrating service selection and scheduling, car classification and blocking, train makeup, and shipments routing based on a three-layer cyclic space-time network representation. 

Shunting yards, also known as marshalling or classification yards, play a key role in rail freight transport networks, acting as hubs where inbound trains are first disassembled and the carriages are then to form new convoys, generating new trains which transport the load towards the correct destinations. This procedure allows carriages to be sent through the network according to their destinations without the need for many connections. Therefore, time savings in shunting operations \citep{Jaehn2015_DC} have a great impact on cost savings in the movement of freight through the rail network \citep{Boysen2012_DC}. In passenger transportation, shunting operations focus on train units that are not necessary to operate a schedule and must be parked at shunt yards. Since different types of trains use the rail infrastructure, the specific type of a unit restricts the set of shunt tracks where they can be parked. The aim of this problem is to assign train locations to the shunt tracks while minimising routing costs from platforms to the corresponding shunt tracks \citep{Huisman2005_DC, Kroon2008_DC}. For a more detailed description of the optimisation problems involved in shunting operations, we refer the reader to the works by \citet{Jaehn2016_DC} and \cite{Ruf2021_DC}.

\subsection[Transportation: Maritime (Harilaos~N.~Psaraftis)]{Transportation: Maritime\protect\footnote{This subsection was written by Harilaos~N.~Psaraftis.}}
\label{sec:Transportation_Maritime}

Maritime transportation carries more than 80\% of the world’s trade and some 70\% of the value of that trade \citep{Unctad2022-sd_HP}. The spectrum of Operational Research (OR) applications in maritime transportation is broad. Following the classification of \cite{Christiansen2013-cg_HP}, these problems can be broken down into three levels: \textit{strategic}, \textit{tactical} and \textit{operational}. Some typical problems in each of these levels will be described in this section.

It is important to note that, in much of the OR maritime transportation literature, traditional economic criteria such as cost minimisation or profit maximisation are the norm, and environmental criteria (for instance emissions minimisation) are less frequent. However, with the quest to decarbonise shipping \citep{Imo2018-zj_HP}, the body of knowledge that includes environmental criteria is growing very fast in recent years. Sometimes environmental criteria map directly into economic criteria: if for instance \textit{fuel cost} is the criterion, and since it is directly proportional to ship \textit{emissions}, if fuel cost is to be minimised as an objective, so will emissions, and the solution is \textit{win-win}. However, for other objectives this direct relationship may cease to exist and one would need to look at environmental criteria in their own right. 

In conceptual terms, if $x$ is a vector of the decision variables of the problem at hand, $f(x)$ is the fuel cost associated with $x$, $c(x)$ is the cost other than fuel and $m(x)$ are the associated maritime emissions (carbon, sulphur, or other), then a generic optimisation problem is the following:
\begin{gather*} 
\text{Minimise } \alpha(f(x)+c(x))+\beta m(x) \\
\text{s.t. }  x \in X
\end{gather*}
where $\alpha$ and $\beta$ are user-defined weights (both  $\geq0$) representing the relative importance the decision maker assigns to cost versus emissions, and $X$ represents the feasible solution space, usually defined by a set of constraints.

One can safely say and without loss of generality that if $d(x)$ is the amount of fuel consumed, $p$ is the fuel price, and $e$ is the emissions coefficient (kg of emissions per kg of fuel), then $f(x)=pd(x)$ and $m(x)= ed(x)$. Therefore $f(x)= km(x)$ with $k=p/e$, as both $f(x)$ and $m(x)$ are proportional to the amount of fuel consumed $d(x)$. The cases that different fuels are used onboard the ship, for instance in the main engine vs the auxiliary engines, or if fuel is switched from high to low sulphur along the ship's trip, represent straightforward generalisations of the above formulation. Then the above problem can also be written as
\begin{gather*} 
\text{Minimise } \alpha c(x)+(\alpha k+\beta)m(x) \\
\text{s.t. }  x \in X
\end{gather*}

\noindent The following special cases of the above problem are important:

\begin{enumerate}[noitemsep]
   \item The case $\alpha=0$, $\beta>0$, in which the problem is to minimise emissions.
   \item The case $\alpha>0$, $\beta=0$, in which the problem is to minimise total cost.
   \item The case $c(x)=0$, in which fuel cost is the only component of the cost.
\end{enumerate}

\noindent A solution $x$* is called \textit{win-win} if both case (1) and case (2) have $x$* as an optimal solution. It is important to realise that such a solution may not necessarily exist.

It is also straightforward to see that in case (3), cost and emissions are minimised at the same time and we have a \textit{win-win} solution. It is clear that $c(x)=0$ is a \textit{sufficient condition} for a win-win solution. But this is not a \textit{necessary} condition, as it is conceivable to have the same solution being the optimal solution under two different objective functions. An interesting question is to what extent policy makers can introduce either (\textit{a}) a Market Based Measure (MBM) such as a fuel tax and/or (\textit{b}) a set of constraints, that would make win-win solutions possible. 

Let us now examine some typical OR problems in the 3-level hierarchy.

\textit{Strategic level problems} involve planning horizons of several years (from 1 to 25). Among them, \textit{fleet size and mix} problems involve the basic questions, what is the best mix for a shipping company's fleet in the years ahead? How large should these ships be? How many should they be, and how fast they should go? See \cite{Alvarez2011-do_HP}, \cite{Zeng2007-xf_HP} and \cite{Pantuso2014-bk_HP} for some work in this area.

\textit{Network design} problems also belong to the strategic problem category and are special to liner shipping. They involve the design of a liner company's network, which comprises the ports it would serve, the routes it will use, which ports will be chosen as hub ports, how are the company's feeder networks configured, and whether the company will use the hub-and-spoke concept or direct calls. See \cite{Agarwal2008-vx_HP}, \cite{Reinhardt2012-do_HP}, and \cite{Brouer2014-ue_HP} for more on these problems.

\textit{Tactical level problems} involve intermediate planning horizons, from a few days to a year. Among them, \textit{ship routing and scheduling} is perhaps the most important problem class, mainly for tramp shipping, with works by \cite{Christiansen2013-cg_HP}, \cite{Andersson2011-qg_HP}, \cite{Fagerholt2010-lr_HP}, and \cite{Lin2011-yp_HP}. Routing and scheduling of offshore supply vessels belongs also to this area \citep{Halvorsen-Weare2011-fm_HP,Norlund2013-qg_HP}. All of these problems call for the determination of the best set of ship routes under some predefined criteria.

\textit{Fleet deployment} is also included in the class of tactical level problems, calling for the allocation of ships to routes \citep[see][among others]{Meng2011-of_HP,Andersson2015-rj_HP,Lai2022-nz_HP}. \textit{Speed optimisation} problems are also tactical level problems and have received increased attention in recent years, due to the pivotal role of ship speed with regard to both economic and environmental criteria. Due to the fact that fuel consumption is a nonlinear function of ship speed, these problems are typically nonlinear. Related formulations attempt to find best vessel speeds along the legs of the route, according to specific criteria \citep[see][]{Psaraftis2013-zv_HP,Fagerholt2013-on_HP,Magirou2015-bq_HP}. These problems may also involve \textit{flexible frequencies} \citep{Giovannini2019-fb_HP}.

Speed and route decisions may also be combined \citep{Psaraftis2014-lu_HP,Wen2017-ku_HP}. One of the perhaps counter-intuitive results of these combined scenarios is that sailing the minimum distance route at minimum speed does not necessarily minimise fuel consumption and hence emissions. This may be so whenever the minimum distance route involves a heavier load profile for the ship. Depending on ship type, the difference in fuel consumption between a fully loaded and a ballast (empty) condition can be up to 40\%. A result that is less surprising is that expensive cargoes sail faster and hence induce more emissions. This is to be expected if cargo in-transit inventory costs are taken into account.

\textit{Modal split/discrete choice models} examine scenarios in which shippers may choose a transportation mode that is alternative to the maritime mode as a result of unfavourable time, cost, or other considerations. As a result, cargoes from the Far East to Europe may prefer the rail vs the maritime mode, or cargoes in European short sea trades may choose the road mode as opposed to shipping. Such modal shifts may increase the overall level of CO2 and may warrant mitigation measures by the shipping lines and the policy makers. Papers that look into this problem include \cite{Psaraftis2010-rt_HP} and \cite{Zis2017-ka_HP,Zis2019-ne_HP}. A multi-commodity network flow formulation in the context of China’s Belt and Road initiative is given by \cite{Qi2022-le_HP}. 

\textit{Operational level problems} concern problems with planning horizons from a few hours to a few days. Among them, a very important class of problems concerns \textit{weather routing} scenarios. The important difference vis-à-vis the ship routing and scheduling problems described earlier is that weather routing problems are typically \textit{path} problems defined as trying to optimise a ship's track from a specified origin to a specified destination, under a prescribed objective and under time varying and maybe also stochastic weather conditions. Decision variables include the selection of the ship’s path and the speeds along the path, and typical objectives include minimum transit time and minimum fuel consumption. Several constraints such as time windows, or constraints to accommodate a feasible envelope on ship motions, vertical and transverse accelerations and ship loads such as shear forces, bending moments and torsional moments can be introduced. The influence of currents, tides, winds and waves, which may be varying in both time and space should be taken into account. See \cite{Perakis1989-zi_HP}, \cite{Lo1998-ix_HP}, and \cite{Zis2020-uk_HP} for some references on this topic.

\textit{Disruption management} is also another important operational level problem class and typically refers to liner shipping. It entails actions that can help the shipping company manage its recovery from possible disruptions of its schedule. Such disruptions may be the result of bad weather, port strikes, equipment malfunction, or more recently, the COVID-19 pandemic that caused massive congestion in many ports worldwide or the Ever Given incident that disrupted traffic in the Suez Canal and the Far East to Europe route in 2021. See \cite{Qi2015-rp_HP} and \cite{Asghari2022-pq_HP} for work in this area.

\textit{Terminal management, berth allocation, and stowage planning} problems also belong to the class of operational level problems, as they deal with an important part of the overall maritime supply chain, that of the coordination between a ship and a port. See \cite{Moccia2006-cl_HP}, \cite{Goodchild2007-rb_HP}, and \cite{Zhen2015-mh_HP} for some related work. 

To conclude, maritime transportation constitutes an important application area for OR, and the related problems are interesting and significant, both from a methodological perspective and from a business and policy perspective. This is so both for traditional economic performance criteria and for environmental criteria, the importance of the latter getting higher in recent years. 

\subsection[Transportation: Aviation (Virginie~Lurkin \& Vikrant~Vaze)]{Transportation: Aviation\protect\footnote{This subsection was written by Virginie~Lurkin and Vikrant~Vaze.}}
\label{sec:Transportation_Aviation}
According to the Air Transport Action Group, in 2019, the world's 1,478 airlines transported 4.5 billion passengers to 3,780 airports, generating 11.3 million direct jobs. Today's airlines are sophisticated businesses making aviation a worldwide economic engine. Yet, aviation is a competitive industry, vulnerable to exogenous shocks, e.g., oil prices, infectious diseases or terrorism. This leads to high costs, and low profit margins, even in the best of times. To tackle these challenges, the industry relies heavily on Operational Research (OR) for decision-making. Prominent OR application domains within aviation include revenue management, airline schedule planning, airline operations recovery, airport flight scheduling, and air traffic flow management. Additionally, some recent OR studies focus on modelling delay propagation through aviation networks

\subsubsection*{Revenue management (RM)}
RM is broadly defined as the strategies and tactics to increase revenues by optimally matching demand for products/services with the available capacity. Seat allocation and pricing are the two main decisions to \textit{control} ticket sales of different fare-classes. Models using capacity allocation as the control variable are called quantity-based RM models. They allocate seats to fare-classes with exogenously determined prices. In contrast, price-based RMs uses pricing policies to maximise revenues. Early RM models focused on overbooking -- the practice of selling more tickets than seats to hedge against cancellations or no-shows. Though various static and dynamic models have been presented since the pioneering work of \cite{rothstein1971airline_VLVV}, airlines mostly use simpler static policies in practice.

Static and dynamic models have been proposed for both single-leg and network-wide seat allocation. Static models optimise seat allocation at a certain time, typically the beginning of the booking period. Dynamic models monitor and adjust to the booking process over time. The earliest static leg-based approach \citep{littlewood1972forecasting_VLVV} considered two fare-classes. \cite{brumelle1990allocation_VLVV} relaxed the assumption of statistical independence between demands. For the multi-class problem, \cite{belobaba1987survey_VLVV} introduced the Expected Marginal Seat Revenue heuristic, a widely used approach in practice. Many studies \citep[e.g.,][]{brumelle1993airline_VLVV} provided optimality conditions for static models, while others developed methods to compute optimal protection levels in the absence of demand information, using optimality conditions \citep{van2000revenue_VLVV} or stochastic approximations \citep{kunnumkal2009stochastic_VLVV}. Dynamic formulations allow time-based controls, but require restrictive demand assumption for tractability, limiting practical impact. Solving network models exactly is computationally hard. Accordingly, most studies on network models use approximations, based on deterministic linear programming \citep{talluriRevenueManagementGeneral2004_AKSJF}, randomised linear programming \citep{talluri1999randomized_VLVV} or decomposition into single-resource problems, as well as solutions using simulation-based optimisation \citep{bertsimas2005simulation_VLVV}. Seat inventory control usually assumes capacity to be fixed, an assumption relaxed by \cite{busing2019capacity_VLVV} integrating capacity uncertainty in leg-based RM. Others integrated inventory control and pricing \citep{you1999dynamic_VLVV}.

Simplest deterministic pricing models are price-sensitive versions of the well-known newsvendor problem \citep{gallego1994optimal_VLVV}. This allows mathematical derivations of optimal prices. Several studies, such as, \cite{feng1995optimal_VLVV}, generalised this problem to include demand dynamics and/or multiple products. Stochastic dynamic programming is a natural way to tackle dynamic pricing. Dynamic models depict reality more accurately, but are harder to solve \citep{gallego1994optimal_VLVV}. Interestingly, solutions to deterministic models are usually good approximations for their stochastic counterparts, and are often used in practice. Traditional RM assumed independent demand, ignoring product substitutability. With the seminal paper of \cite{talluriRevenueManagementGeneral2004_AKSJF}, the RM field has shifted toward including customer choice behaviors within pricing and capacity decisions. \S\ref{sec:Revenue_management} provides a detailed overview of RM concepts and trends beyond aviation.

\subsubsection*{Airline schedule planning (ASP)}
ASP is the process of designing airline schedules maximising profits subject to resource constraints. Taking demand, airport and aircraft characteristics, and maintenance and personal requirements as inputs, ASP outputs selected flight timetables, aircraft schedules and crew duty plans. Most ASP steps typically occur before RM actions and thus constrain the set of decisions available to RM systems. Key ASP steps include fleet planning, route planning, frequency planning, timetable design, fleet assignment, aircraft routing and crew scheduling.  Fleet planning involves decisions regarding purchasing, selling, and leasing of aircraft fleet, while route planning selects airport pairs to operate nonstop flights. Early studies, e.g., \cite{hane1995fleet_VLVV}, matched a predetermined set of flights with aircraft types, developing a fleet assignment model (FAM). The FAM specifically focuses on fleet assignment, which is one particular step within the overall ASP process. The basic FAM, a mixed-integer linear program, minimised costs of operating aircraft and passengers unserved, given passenger demand for individual flight legs. This leg-based approach ignores that passengers often fly on multiple flights in connecting itineraries. 

\cite{barnhart2002itinerary_VLVV} overcame this limitation via an itinerary-based FAM to explicitly model network effects. Some studies developed tractable solution approaches. \cite{barnhart2009airline_VLVV} proposed a subnetwork-based decomposition for capturing FAM's revenue implications, an approach recently extended by \cite{yan2022choice_VLVV} to solve a FAM incorporating passenger choice. Others extended FAM by incorporating incremental timetable design decisions, e.g., changes to flight timings \citep{desaulniers1997daily_VLVV} or selection of optional flights \citep{lohatepanont2004airline_VLVV}. \cite{wei2020airline_VLVV} developed a \textit{clean slate} heuristic optimising entire timetables and fleet assignments under choice-based demand. Frequency planning, which optimises the number of flights operated during a day or part of a day, rather than deciding exact timetables, has also received attention, with an emphasis on capturing affects of competition from other airline and high-speed rail operators \citep[e.g.,][]{cadarso2017integrated_VLVV}.

The last two steps in schedule planning are conceptually similar. Aircraft routing assigns individual aircraft to flights while ensuring that each aircraft undergoes periodic maintenance, and crew scheduling assigns crew to operate flights while satisfying a myriad of crew regulations. Early studies individually optimised aircraft routing \citep{gopalan1998aircraft_VLVV} or crew scheduling \citep{graves1993flight_VLVV}. \cite{lavoie1988new_VLVV} used column generation, an effective solution approach for both problems, to crew scheduling, while \cite{cordeau2001benders_VLVV} used Benders decomposition to jointly solve both problems.

Good schedules not only minimise planned costs, but are also robust to disruptions, to keep the actual costs low. Researchers in early 2000s optimised robustness proxies, e.g., station purity, short cycles, crew swapping opportunities, and crew schedule slack \citep{schaefer2005airline_VLVV}. Later studies directly minimised total planned and unplanned costs of aircraft routing \citep{lan2006planning_VLVV} and crew scheduling \citep{yen2006stochastic_VLVV} separately, and also jointly \citep{dunbar2012robust_VLVV}. Recent studies have used robust optimisation to solve the aircraft routing \citep{yan2018robust_VLVV} and crew scheduling \citep{antunes2019robust_VLVV} problems.

\subsubsection*{Airline operations recovery (AOR)}
AOR encompasses the actions undertaken to repair schedules, when disruptive events such as inclement weather, equipment failures, etc., take place. \cite{rosenberger2003rerouting_VLVV} developed a model and a solution heuristic for repairing aircraft routing, whereas \cite{lettovsky2000airline_VLVV} tackled crew recovery. For the integrated recovery problem, \cite{petersen2012optimization_VLVV} developed a decomposition strategy, while \cite{maher2016solving_VLVV} used column-and-row-generation. Recent recovery studies incorporated other key elements, including flight planning \citep{marla2017integrated_VLVV} and passenger no-shows \citep{cadarso2022passenger_VLVV}.

\subsubsection*{Airport flight scheduling (AFS)}
Beyond airline decision-making, OR is also used to improve decision-making of central authorities and air traffic managers. Research over the past decade demonstrated the potential for enhancing social welfare by constraining schedules at busy airports via slot-control mechanisms \citep{swaroop2012more_VLVV}. Some studies balanced strategic cost of scheduling changes against tactical cost of delays, for a single airport \citep{jacquillat2015integrated_VLVV} or multiple airports \citep{wang2020stochastic_VLVV}. \cite{zografos2012dealing_VLVV} used an integer program for allocating slots to airlines under administrative controls. \cite{fairbrother2020slot_VLVV} attempted to balance the often-conflicting goals of efficiency, equity and the incorporation of airline preferences in optimising slot-scheduling mechanism.

\subsubsection*{Air traffic flow management (ATFM)}
The tactical side of airport and airspace capacity management has received considerable OR attention since the 1990s. ATFM is a broad term used to define key interventions, such as ground holding of airplanes, that ensure safe and efficient flight operations by restricting flow of aircraft into congested airspaces. \cite{terrab1993strategic_VLVV} and \cite{vranas1994multi_VLVV} proposed the single-airport and multi-airport ground holding problems, respectively. The latter was extended to include enroute capacities by \cite{bertsimas1998air_VLVV}. \cite{bertsimas2011integer_VLVV} additionally incorporated flight rerouting and solved larger-scale problems. Adoption of the collaborative decision-making (CDM) paradigm in practice ushered in a new era of research. Advocating increased agency to airlines, \cite{vossen2006slot_VLVV} provided an integer program for slot trading mechanism design under CDM. Recent studies \citep[e.g.,][]{starita2020air_VLVV} are increasingly focused on explicit handling of uncertainty on both demand and capacity side within the ATFM optimisation problems.

\subsubsection*{Modelling delay propagation}
Tightly coupled aviation networks make disruption management particularly challenging. Delays and disruptions in one part of the network propagate to other parts, through aircraft, crew and passenger connections. Recent studies quantified these propagation effects. First, \cite{pyrgiotis2013modelling_VLVV} proposed an analytical queuing and network decomposition model for aircraft-based delay propagation. \cite{barnhart2014modeling_VLVV} presented discrete choice models for passenger itinerary estimation and a reaccommodation heuristic for passenger delay calculations. \cite{wei2018modeling_VLVV} solved inverse optimisation for estimating crew itineraries and crew-based delay propagation. These studies attempted bridging the gap between sparse and aggregate public datasets, and the detailed and disaggregated data needs for aviation OR research.

\subsubsection*{Further reading}
Readers interested in aviation OR are referred to the second edition of the book by \cite{belobaba2015global_VLVV}. In particular, Chapters 4 and 5 focus on pricing and RM, Chapters 8 and 10 on schedule optimisation, robustness and recovery, and Chapter 14 on air traffic management and control. Looking ahead, it is apparent that OR will keep finding natural applications within aviation, especially given the exciting disruptive innovations within urban air mobility. Rapidly growing fields of passenger air taxi operations and drone operations for parcel deliveries are giving rise to new variants of well-known OR problems, e.g., network design \citep{wang2022vertiport_VLVV}, travelling salesperson \citep{roberti2021exact_VLVV}, vehicle routing \citep{dayarian2020same_VLVV}, and facility location \citep{chen2022scalable_VLVV}. 

\subsection[Transportation: Network design (Mike~Hewitt)]{Transportation: Network design\protect\footnote{This subsection was written by Mike~Hewitt.}}
\label{sec:Transportation_Network_design}

In a transportation context, the term \emph{Network Design} \citep{magnanti1984network_MH} generally refers to planning the \emph{supply} side of a transportation system so that it efficiently satisfies some estimate of \emph{demand} within the quality standards of the customers using the system. The planning decisions typically prescribe the movements of vehicles, or convoys (e.g., a railroad train or tug and barges), between stations/terminals in the network to transport people or goods. Network design is typically undertaken for situations wherein what is transported, be it people or goods, is small relative to vehicle capacity. Thus, one primary measure of efficiency is vehicle utilisation, with high utilisation achieved through consolidation. Quality is typically measured based on on-time delivery. 

Network design is relevant to passenger transportation systems such as urban public-transport \citep{Mauttone2021_MH} by bus \citep{ceder1986bus_MH} or light rail \citep{farahani2013review_MH}, as well as systems providing interurban transport by train \citep{hooghiemstra1999decision_MH} or airplane \citep{franke2017network_MH}. It is also relevant to a wide range of goods transportation markets, such as parcel and small-package \citep{barnhart1996air_MH} and less-than-truckload freight \citep{powell/sheffi:1989_MH}. A network design case study for a postal carrier can be found in \cite{winkenbach2016strategic_MH}. Transportation carriers serving these markets may rely on one or more modes, including motor carrier \citep{Bakir2021_MH}, rail \citep{Chouman2021_MH}, ocean \citep{Christiansen2021_MH}, and inland waterway \citep{konings2003network_MH}. The planning of vehicle and goods movements by each mode and synchronisation of goods moving from one mode to the next (e.g., intermodal) can be assisted by network design \citep{arnold2004modelling_MH}.

For different modes the scope of design decisions prescribed by network design models may be broadened in different ways. For example, modes such as rail and inland waterway involve multiple layers of consolidation. For rail \citep{Zhu2014_DC}, goods are consolidated into rail cars, which are then consolidated into blocks that are transported by the same locomotive. For motor carriers, vehicles can not yet move without a driver, whose movements and schedules are restricted by governmental safety regulations and potentially labour management practices that dictate the driver return periodically to a specific physical location in the network (e.g., his/her domicile). Network design models for motor carriers may build schedules for drivers that observe safety regulations \citep{Crainic2018_MH} as well as determine how many drivers should be associated with each physical location \citep{HEWITT2019324_MH}.

The network design problem is typically modelled as a Mixed Integer Program (MIP) formulated on a directed graph \citep{Crainic2021MIP_MH}. Nodes in such a graph model physical locations, potentially at different points in time. Directed edges between such nodes model transportation that begins in one physical location and ends at another. Edges may encode a scheduling dimension, such as when a vehicle departs from one location and arrives at another, that depends in part on the travel time required for the physical move \citep{erera2013improved_MH}. Associated with an edge is a function that maps the amount of vehicle capacity made available on that edge to cost. Typically, it is a step function with each step modelling an increase in capacity due to dispatching an extra vehicle. Commodities model people or goods that are to be transported; associated with each commodity is an origin node, a destination node, and a size. 

The classical network design problem seeks to find a path for each commodity that begins at its origin node, ends at its destination node, and potentially visits one or more intermediate nodes. The problem evaluates these paths with respect to the total cost of capacity made available to support them and seeks to minimise that total cost. Some network design models \citep{Frangioni2021_MH} instead minimise costs that are a function of the amount of goods transported on an edge, as opposed to the capacity made available to transport them. Network design is an optimisation problem that has received significant attention both for its practical relevance and the computational challenges \citep{johnson1978complexity_MH} associated with solving it. 

Most MIP formulations of the network design problem involve commodity flow variables that model the transportation of goods within the network and another set of edge-based variables that model the transportation of vehicles. Typically, commodity flow variables are continuous when a shipper's goods can be divided and routed on multiple paths or binary when they cannot.  Commodity flow variables are typically edge-based, but some models involve paths from shipment origin to shipment destination. The use of a path formulation typically necessitates column generation \citep{HEWITT2019324_MH}. However, unlike the vehicle routing problem, extended, path-based formulations of the network design problem do not provide stronger linear relaxations than compact, arc-based formulations. Depending on the context and mode the vehicle edge variables may either be binary or integer. Linking constraints are included in the formulation to ensure sufficient vehicle capacity travels on an edge to carry the commodities making that transportation move.  Typically, much larger cost coefficients are associated with vehicle edge variables than commodity flow variables.

The majority of literature on network design focuses on deterministic models wherein it is presumed all parameter values (costs, capacities, demands) are known with certainty. However, given that network design models are often solved as part of a tactical planning exercise, uncertainty has been studied \citep{Hewitt2021SND_MH}. Much of that work focuses on uncertainty in commodity sizes and models such problems as two stage stochastic programs wherein vehicle movements are planned in the first stage and commodities are routed in the second stage given the vehicle movements prescribed in the first. There has been limited work on robust optimisation models \citep{Koster2021_MH} or those that view network design in a dynamic context \citep{ALHAJJHASSAN2022102885_MH}. 

Both exact \citep{Crainic2021Exact_MH} and heuristic \citep{Crainic2021Heuristic_MH} solution methods for deterministic network design models have been proposed.  One challenge associated with solving MIP formulations of network design problems is that the linking constraints often lead to fractional vehicle edge variables. Thus, the linear programming relaxations of network design MIPs often yield weak bounds on the objective function value of an optimal solution to the MIP. As a result, much of the literature that focuses on speeding up the solution of MIP formulations of the network design problem focuses on strengthening formulations with \textit{valid inequalities} \citep{NW88_SMPT}. Such inequalities are typically either based on classical ideas such as flow covers from integer programming \citep{gu1999lifted_MH} or leveraging the network structure of the problem \citep{raack2011cut_MH}. Another approach taken to solve network design problems is Benders decomposition \citep{benders2005partitioning_MH,costa2005survey_MH}, particularly when second stage variables are continuous and the optimisation problem resulting from fixing the network design is a linear program.

Another challenge associated with solving MIP formulations of the network design problem is due to the size of the network on which the MIP is formulated  when that network encodes time. The classical approach to representing time in network design is to formulate a MIP on a network wherein multiple nodes represent the same physical location, albeit at different points in time \citep{HewittVuCrainic2016_MH}. Similarly, multiple edges represent the same physical transportation move, albeit at different departure and arrival times. Such networks are typically referred to as time-expanded networks and the overall solution procedure in contexts that require the modelling of time is to construct such a network, formulate a MIP on that network, and then solve that MIP. \cite{BOLAND2019195_MH} study the impact on solution quality of modelling time at different granularities and observe that the finer the representation the higher the quality of the resulting solution. However, such an approach can be computationally challenging when long planning horizons must be modelled or fine representations of time are used, as both cases lead to networks and resulting MIPs that are very large. An alternate approach, called \textit{Dynamic Discretisation Discovery} \citep{DDD2017_MH,EnhancedDDD_MH} proposed to instead generate time-expanded networks in a dynamic and iterative manner.
 
Heuristic methods for deterministic network design models can be classified into one of two categories. The first category focuses on metaheuristics \citep{hussain2019metaheuristic_MH} and  neighbourhood structures. Early heuristics \citep{POWELL1983471_MH} proposed for network design models searched neighbouring solutions by reducing the capacity on one edge in the network and, if necessary, increasing the capacity on another. However, more recent and  effective methods have proposed more complex neighbourhood structures such as cycles or paths \citep{ghamlouche2003cycle_MH}. The second category focuses on what is generally called matheuristics \citep{maniezzo2021matheuristics_MH}. In these heuristics, a neighbourhood of a solution is searched by formulating and solving the MIP of the network design problem, albeit with the values of subsets of variables fixed to their values in the solution at hand \citep{hewitt2010combining_MH}. This is repeatedly done and with different mechanisms used for selecting subsets of variables to fix. 

Similarly, both exact and heuristic solution methods have been proposed for stochastic network design models that take the form of scenario-based two stage stochastic programs. The vast majority of such stochastic programs studied to date involve continuous commodity flow variables in the second stage. As a result, the second stage subproblems are linear programs and the overall stochastic program is amenable to Benders decomposition \citep{Birge2011-ut_HL}. Thus, much of the methodological work on solving such stochastic programs has focused on techniques for speeding up or rendering more impactful different steps in the Benders scheme \citep{magnanti1981accelerating_MH,crainic2021partial_MH}. While Progressive Hedging \citep{RockWetsPH_MH} is an exact method for two stage stochastic programs with continuous variables in both stages, it has been used as the basis of heuristic methods for stochastic network design \citep{crainic2011progressive_MH,crainic2014scenario_MH}. 

\cite{crainic2021network_MH} contains deeper dives into the subjects touched on here as well as discussions of those not discussed. 

\subsection[Transportation: Vehicle routing (Claudia~Archetti \& Maria~Battarra)]{Transportation: Vehicle routing\protect\footnote{This subsection was written by Claudia~Archetti and Maria~Battarra.}}
\label{sec:Transportation_Vehicle_routing}

The Capacitated Vehicle Routing Problem (CVRP) was first proposed by \cite{DantzigR1959_CA_MB}, and named the Truck Dispatching Problem. The goal was that of routing a fleet of identical gasoline delivery trucks from a central depot to service stations (often referred as `customers'). Each truck had to return to the central depot, after visiting an ordered subset of the customers. All customers had to be visited once by a vehicle delivering all their gasoline requirements in the one delivery. The objective was the minimisation of the routing costs, as the sum of the travelling distances of every truck. The CVRP classical definition is the same as that proposed by \cite{DantzigR1959_CA_MB} more than 60 years ago. Introducing a capacitated fleet of vehicles makes the CVRP for a much harder generalisation of the Travelling Salesman Problem \citep{Flood1956_CA_MB}. 
    
The CVRP definition has been enriched over the decades to take into account all the delivery requirements of the customers and of the transportation providers, as well as the characteristics of the available fleet of vehicles, and the increasing availability of technology (i.e., GIS and real time mapping, autonomous vehicles, shared mobility systems and so on). The research literature has flourished with new variants, as well as more  sophisticated  and flexible solution approaches. This chapter aims at providing pointers to key milestones achieved in the last 60 years of the CVRP literature, identifying the latest and most successful exact and metaheuristic algorithms, as well as referencing the most famous online challenges and standard techniques for benchmarking CVRP solution algorithms. 
    
The CVRP `classical' variants and solution approaches are well summarised in  \cite{TothV2002_CA_MB}. 
This book provides key references and definitions for critical application features, as for the CVRP with Time Windows, the CVRP with Backhauls and the CVRP with Pickup and Delivery, the CVRP with vehicle/site dependencies, the CVRP with inventory and the stochastic CVRP. \cite{GoldenRW2008_CA_MB} extends the definition of the classical variants to routing problems with heterogenous fleets, periodic routing problems, split routing problems, dynamic and online routing problems. 
\cite{TV14_SMPT} further widen the remit of application of routing algorithms to maritime applications, disaster relief distribution problems, and considers up-to-date objective functions different than minimising the distance travelled. More recently, fleets of electric vehicles \citep{PelletierJL2016_CA_MB}, problems over time \citep{mor2022vehicle_CA_MB}, drones \citep{OttoACGP2018_CA_MB}, cargo boats \citep{Christiansen2013-cg_HP} and warehouse pickers  \citep{Schiffer2022_CA_MB} have been embedded in routing settings. The new dynamic environment inspired research on stochastic \citep{Gendreau2016_CA_MB}, dynamic \citep{soeffker2021stochastic_CA_MB} and time-dependent \citep{gendreau2015time_CA_MB} routing problems.
    
An up-to-date survey on recent trends can be found in \cite{VIDAL2020401_CA_MB}, in which the CVRP extensions due to richer objective functions, the integration with other optimisation problems, and application-oriented transportation requirements are surveyed. \cite{PartykaH2014_CA_MB} discuss routing algorithms from the practitioners' perspective, and surveys which are the requirements of a logistics company when they acquire a  routing software. 
    
Next the most successful CVRP solution algorithms are summarised, first discussing exact methods. Formulations with a polynomial number of variables and constraints were the first proposed mathematical models, as for the two-commodity formulation  by \cite{laporte1992vehicle_CA_MB} and \cite{Baldacci2Flow2004_CA_MB}. They have the advantage of being easy to use (as they just require encoding in the syntax of the solver). The disadvantage of them however is their poor performance due to high dimension of the formulations, and the weakness of the continuous relaxation. Better results were obtained from formulations with an exponential number of constraints, such as those in which subtour elimination constraints are added dynamically to the formulations in a branch\&cut fashion \citep{Padberg-Rinaldi:1991_IL}. The CVRPSEP library by \cite{LysgaardLE2004_CA_MB} provides separation procedures for subtour elimination constraints, as well as other strengthening additional inequalities. The most successful exact solution framework is up-to-date the branch\&cut\&price \citep{DDS06_ALAL,laporte2009fifty_CA_MB}. This method is based on the Dantzig-Wolfe decomposition \citep{Desrosiers2005_CA_MB}. Binary variables model if a route is used or not in the solution, thus their corresponding set is exponential in size. As a consequence, a restricted set of variables is used to initiate the formulation and only profitable routes are iteratively generated solving a subproblem, called the pricing problem. The CVRP pricing problem is a shortest path with resource constrains, and it is typically solved through dynamic programming \citep{Irnich2005_CA_MB}.  Some of the most relevant milestones in developing branch\&cut\&price algorithms for the CVRP are combining branch\&cut and column generation into the first branch\&cut\&price \citep{Fukasawa2006_CA_MB}, applying bi-directional search in the subproblem \citep{RighiniS2008_CA_MB}, introducing subset row cuts \citep{JepsenPSP2008_CA_MB},   using ng-routes to speed up the subproblem solution \citep{Roberti2011_CA_MB},  using stabilisation techniques for dual values \citep{Gschwind2016_CA_MB,Pessoa2018_CA_MB}, and  proposing primal heuristics based on the restricted master problem \citep{Sadykov2019_CA_MB}. The reader might refer to \cite{desaulniers2002accelerating_CA_MB} for the most widely used acceleration techniques for the solution of the pricing problem.
   
Lately, the work of \cite{PessoaSUV2020_CA_MB} provides an impressive open-source branch\&cut\&price algorithm, based on \cite{Pecin2017_CA_MB}. This algorithm provides state-of-the-art exact solutions for the CVRP and, using a flexible solution representation, for most of the well known routing variants and other sequencing problems. The tool incorporates the algorithmic components previously mentioned, as well as other recent developments \cite[see for example,][]{Sadykov2021_CA_MB}, and compares favourably to other  branch\&cut\&price implementations.  Some of the most powerful exact algorithms for the CVRP, available in different programming languages, are publicly available at \cite{Sadykov2022_CA_MB}.
   
Metaheuristics are  capable of solving very large CVRP instances in limited computing time, however there is no proof of optimality for the solutions found. They are typically initialised with solutions generated by constructive heuristics (the Clarke and Wright is a famous example, \citealp{ClarkeW1964_CA_MB}). Metaheuristics rely heavily on local search procedures to improve the solution quality and intensify the search, and on a metaheuristic framework to obtain a good balance of diversification and intensification \citep{gendreau2010handbook_CA_MB}. In chronological order, popular CVRP frameworks have been the Tabu Search \citep{Cordeau2005_CA_MB}, the Adaptive Large Neighbourhood Search \citep{PISINGER20072403_CA_MB}, the Iterated Local Search \citep{SUBRAMANIAN20132519_CA_MB}, and the Hybrid Genetic algorithm \citep{VIDAL2022105643_CA_MB}. 
The latter two examples of metaheuristic frameworks are particularly relevant to the CVRP literature  due to their high performance,  their flexibility in solving effectively many VRP variants, and because their code had been made publicly available to the research community (the code presented in \citealp{VIDAL2022105643_CA_MB} is, for example, available at \citealp{Vidal2022_CA_MB}). \cite{VIDAL20131_CA_MB} provide a very good summary of the features that make a CVRP metaheuristic successful.

More recently, examples of algorithms producing very high quality solutions for the CVRP have been:\begin{itemize}[noitemsep]
    \item \cite{arnold2019makes_CA_MB}: data mining is used to identify solution features, and these features are used to effectively guide the search algorithms;
    \item \cite{christiaens2020slack_CA_MB}: SISRs is a ruin and recreate algorithm based on an innovative string removal operator;
    \item \cite{QUEIROGA2021105475_CA_MB}: POPMUSIC is a matheuristic that iteratively solves smaller subproblems by means of the branch\&cut\&price by \cite{PessoaSUV2020_CA_MB};
    \item  \cite{AccorsiVigo2021_CA_MB}: FILO is an Iterated Local Search with acceleration techniques and annealing-based neighbour acceptance criteria;
    \item    \cite{MAXIMO20211108_CA_MB}: AILS-PR is an Iterated Local Search metaheuristic hybridised with Path Relinking; and, 
    \item \cite{CavaliereBF2022_CA_MB}: a refinement heuristic using a penalty-based extension of the Lin and Kerninghan heuristic is combined with a restricted column generation to iteratively select meaningful routes.
\end{itemize} 
  
Clear standards have been set by the CVRP community around which benchmark instances should be used for testing the performance of an algorithm, and which are ways of testing a computer code for a fair comparison with other previously proposed algorithms. \cite{UCHOA2017845_CA_MB} discuss the most widely used instances and provides a link to the repository, in which the input data, as well as the best known solutions, are provided and kept up-to-date by the authors. A more recent set of instances and best known solutions is available in \cite{queiroga2022_CA_MB}, where the authors provide data enabling the use of machine learning approaches to solve the CVRP. \cite{ACCORSI2022229_CA_MB} present the standard practices to test CVRP algorithms: how to determine computing time (typically on a single thread), common ways of tuning parameters, and providing best and average solutions on a specified number of executions, among others.
   
Finally, another popular and flourishing avenue for boosting research on the development of effective solutions approaches for the CVRP and variants is represented by competitions.  Some of the most famous CVRP and routing challenges are:
\begin{itemize}[noitemsep]
   \item the DIMACS challenge \citep{DIMACS_CA_MB}, where the goal was to promote research on challenging routing problem variants;
   \item the Amazon Last Mile Routing Research Challenge \citep{AMAZON_CA_MB}, where a specific problem was tackled, namely, the challenge of embedding driver knowledge into route optimisation;
   \item the recently launched EURO Meets NeurIPS 2022 Vehicle Routing Competition \citep{EURO_CA_MB}, with the goal of developing and comparing machine learning techniques for the CVRP.
\end{itemize}

The Vehicle Routing problem has inspired an incredible amount of research. This is due to the challenges it poses when it comes to solving it, to the many variants related to it and to the relevant practical applications. Despite the decades of research efforts and achievements, interest continues to grow mainly thanks to the emerging topics raised by the ever changing application environment. This chapter provides a brief, but hopefully sufficiently comprehensive overview of the techniques, problem variants and emerging trends which will inspire further research.

\clearpage

\section[Conclusions (Daniele~Vigo \& Said~Salhi)]{Conclusions\protect\footnote{This subsection was written by Daniele~Vigo and Said~Salhi.}}
\label{sec:Conclusions}

This encyclopedic article, dedicated to the 75\textsuperscript{th} anniversary of the \textit{Journal of the Operational Research Society}, is made up of an \textit{Introduction} and two distinct though related sections: \textit{Methods} and \textit{Applications}. The introduction section gives an interesting overview of OR with an emphasis on its origin in the UK and highlights the methods and applications that are covered in this paper. A brief summary of the two sections is given below.

In the first main section (\S\ref{sec:methods}), 24 OR-based methods are presented by experts in their respective areas. These methods, which are given in alphabetical order, are concisely described, each starting with the basics, then moving to advanced and contemporary aspects. The authors also pinpoint challenging limitations while highlighting promising research directions.

As OR is rooted in the need to solve decision problems either through optimisation, statistics, visualisation and information technology tools, or through soft system methodologies, we aim to retain this historical flavor in summarising these methods by adopting a simple three-group categorisation.  

The first category covers optimisation-related topics and includes 10 out of the 24 subsections. It ranges from the original optimisation model of linear programming (LP; \S\ref{sec:Linear_programming}) in the late 1940s to its various extensions. One is obtained  by restricting the decision variables to discrete elements including  binary ones (\S\ref{sec:Mixed_integer_programming}), allowing uncertainty in the input (\S\ref{sec:Stochastic_models}), or relaxing the objective function or the constraints not to be  necessarily linear (\S\ref{sec:Nonlinear_programming}). An interesting area that had been dormant for more than 30 years was revived in the late 1970s and early 1980s by studying  a special case of fractional LP which defines relative efficiency and is  known as data envelopment analysis (\S\ref{sec:Data_envelopment_analysis}). Combinatorial optimisation (\S\ref{sec:Combinatorial_optimisation}), a topic that has fascinated and intrigued many mathematicians of the 18\textsuperscript{th} century, seeks an optimal subset or values from a large finite set of elements. These problems can be defined and solved through graphs and networks (\S\ref{sec:Graphs_and_networks}),  some of which are relatively more difficult  than others. To measure the performance of algorithms in terms of time and space complexity, computational complexity (\S\ref{sec:Computational_complexity}) emerged as a solid foundation for distinguishing between classes defined as $\mathcal{P}$ and $\mathcal{NP}$ and studying $\alpha$-approximation algorithms. One methodology  can be traced back to the Ancient Greek times, and is based on the ‘find and discover principle’, now known as ‘heuristic search’ (\S\ref{sec:Heuristics}), which has experienced a phenomenal growth in the late 1980s and early 1990s. This is a major development since  these methodologies provide the best way  to reduce  not only the risk of getting stuck at poor local optima, but also have the power to  yield practical solutions for complex discrete and global optimisation problems that could not have been solved otherwise. A methodology that  is free from restrictions of linearity and convexity is the study of multi-stage process, as given in \S\ref{sec:Dynamic_programming}. 

The next category includes statistics and decision-based tools and also covers 10 of the 24 subsections. For example, business analytics (\S\ref{sec:Business_analytics}), decision analysis (\S\ref{sec:Decision_analysis}) and visualisation (\S\ref{sec:Visualisation}), though they previously existed under different names, have grown significantly while retaining their simplicity. Machine learning, including artificial intelligence (\S\ref{sec:Artificial_intelligence_machine_learning_data_science}), which borrowed its principles from heuristic search and statistics, has taken off very rapidly in teaching, research and applications. This is mainly due to computer power, sophisticated algorithms, freely available  computer languages such as R and Python, and their ability to handle massive amount of data that are now easily available to the public. Other older topics, though still relevant and widely applicable, have also seen a surge in new developments. These include queueing (\S\ref{sec:Queueing}), forecasting (\S\ref{sec:forecasting}), control theory (\S\ref{sec:Control_theory}),  and game theory (\S\ref{sec:Game_theory}). Given the uncertainty and risk involved in many decisions, risk analysis (\S\ref{sec:Risk_analysis}) is evolving fast so as to handle such environments alongside computer simulation (\S\ref{sec:Simulation}), especially discrete event simulation.  The latter, which has a wide spectrum of applications in both the private and public sectors, has recently been enriched by incorporating  multi-objective optimisation within its evaluation component. 

The last category covers the remaining four subsections.  Although some of these research areas existed in other fields such as system engineering in the 1950s, they have become contemporary OR topics especially in the UK in the late 1970s. Soft OR and problem structuring methods (\S\ref{sec:Soft_OR_and_problem_structuring_methods}) question the problem definition and aim to involve stakeholders for a better understanding, with system thinking (\S\ref{sec:Systems_thinking}) analysing the interactions between people, machines and systems  while also questioning the system boundaries.  A related area is system dynamics (\S\ref{sec:Systems_dynamics})  where the dynamism is incorporated throughout and found to suit better applications with limited but plausible scaling. An interesting, though relatively recent OR area, but with a long history rooted in social psychology, is behavioural OR (\S\ref{sec:Behavioural_OR}), where people's behaviour and culture are incorporated into the decision making process. Although the methodologies included in this category usually do not directly aim to solve problems, they can be complementary  to the harder OR techniques.

The second section covers applications that have been, since the very beginning, strongly interconnected with the development of OR methodologies. This section is very rich in examples coming from many fields. For the sake of brevity, we will not refer to each subsection individually but mention just a few. By reading the section it is evident that, on the one hand OR provides appropriate modelling and solution tools to practical problems that  arise in the real world and are nowadays crucial in the design and management of most systems, from healthcare and other public services, to transportation and manufacturing. On the other hand, the complexity and size of practical problems has always stimulated the progress of OR  towards more efficient and flexible techniques which are capable to cope with the challenges posed by the applications. This mutual and virtuous connection is well reflected by the richness of the Applications section of this work. It highlights not only the traditional areas which saw tremendous research efforts and successful implementations, as the traditional fields of transportation, manufacturing, cutting and packing, and inventory management, but also relatively new and interesting sectors such as sports and education.

It is worth noting that in the \textit{Applications} section (\S\ref{sec:applications}), several dimensions of OR impact in the real world clearly emerge. The first one is the broad range of fields to which OR techniques have already successfully been applied and offer an even larger potential yet to be exploited. These range from \textit{vertical} sectors,  such as supply chain management, disaster relief and recovery, or military applications, where a wide array of  problems are defined and solved through appropriate and varied methodologies, to more \textit{horizontal} domains which may impact  several vertical sectors, like vehicle routing or facility location, for which highly specialised methods have been developed. The second dimension is related to the great variety of methodologies applied to the different contexts. These span the whole tool set of OR, including exact and heuristic methods developed to solve specific optimisation problems, to techniques created to handle uncertainty and multi-criteria and, more recently, integrating artificial intelligence methods. Indeed, the great improvements achieved in the last decades in integer and nonlinear programming now allow  to effectively model and solve many problems arising at the operational and tactical levels, where data are more available and reliable. The uncertainty in the data and the modelling typical of strategic decisions are  successfully handled by a variety of methodologies that have proven to be effective in the solution of real applications which are well reviewed in this work. A third very interesting dimension is represented by the development of new broad research perspectives which may have a strong impact in all fields of OR and  are deeply motivated by  applications. An excellent example is the inclusion of fairness and ethics  in  optimisation which, on the one hand allow for considering important issues favouring the acceptability and usability of the results, and on the other hand pose new methodological challenges.

As a general conclusion, thanks to the advances in computer technology, the availability of massive amount of live data, and novel developments, in both optimisation and statistics,   effective optimisation software, powerful machine learning techniques and visualisation tools now exist to solve problems that were considered practically unsolvable just a decade ago.  Applications have always been a main driver for OR development, and the successes achieved increase the appetite for further improvements. 

In the more classical area of exact and heuristic techniques, there is clearly a need to improve the capability of handling efficiently large and very large-scale instances to cope with more complex and demanding scenarios. This  increase in scale is not only generated by the need to solve larger problems, but also to incorporate various steps of the planning processes into integrated and more comprehensive methods. A  field that still deserves further research efforts is the  consideration of uncertainty in OR methods. Important methodological obstacles have yet to be surmounted and there is clearly a need for the development of simple and pragmatic methods, possibly resulting from the integration of artificial intelligence techniques, which can be applied to the solution of large-scale problems arising in  several important application domains. However, it is also worth stressing that these advances, though they are welcome, may  suffer from shortcomings, such as the local optimality trap, biased data,  and impractical assumptions. These hidden aspects could yield poor outcomes on which  academics and practitioners ought to keep an open eye.

\clearpage

\section*{Disclaimer}
The views expressed in this paper are those of the authors and do not necessarily reflect the views of their affiliated institutions and organisations.

\section*{Acknowledgements}
Fotios Petropoulos would like to thank all the co-authors of this article for their very enthusiastic response and participation in this initiave. He is also indebted to his lead advisor for this project, Gilbert Laporte, as well as Christos Vasilakis, Güneş~Erdoğan, Stephen Disney and Maria Battarra for their help and suggestions. Finally, Fotios is grateful to John Boylan and the other Editors-in-Chief of the \textit{Journal of the Operational Research Society} for inviting this paper to be part of the 75\textsuperscript{th} issue of the journal. Fotios dedicates this article to Professor John Boylan: John, your kindness will always be remembered.

Maria Battarra's work reported in this paper was undertaken as part of the Made Smarter Innovation: Centre for People-Led Digitalisation, at the University of Bath, University of Nottingham, and Loughborough University.

David Canca's work was supported by the University of Sevilla, the Regional Government of Andalucia (Spain) and the European Regional Development Fund (ERDF) under grant US-1381656.

Laurent Charlin and Andrea Lodi would like to thank Didier Chételat and Mizu Nishikawa-Toomey for reading and commenting on drafts of their subsection (\S\ref{sec:Artificial_intelligence_machine_learning_data_science}) and the CIFAR AI Chair and the CERC programs for funding.

Salvatore Greco wishes to acknowledge the support of the Ministero dell'Istruzione, dell'Universita e dellaRicerca (MIUR) - PRIN 2017, project ``Multiple Criteria Decision Analysis and Multiple Criteria Decision Theory'', grant 2017CY2NCA.

Katherine Kent and Sam Rose thank Mithu Norris for the help and coordination with \S\ref{sec:Government_and_public_sector} but also Emma Hickman and Ffion Lelii for their contribution.

Silvano Martello, Paolo Toth and Daniele Vigo were supported by Air Force Office of Scientific Research under Grants no. FA8655-20-1-7012, FA8655-20-1-7019, FA9550-17-1-0234 and FA8655-21-1-7046.

Dimitrios Sotiros's work was partially supported by the National Science Center (NCN, Poland) grant no. 2020/37/B/HS4/03125.

Greet~Vanden~Berghe and Sanja~Petrovic acknowledge the advice provided by Andrea Schaerf (University of Udine).

Rafał Weron's work was partially supported by the National Science Center (NCN, Poland) grant no. 2018/30/A/HS4/00444.

\clearpage

\appendix

\section{List of acronyms}\label{sec:acronyms}

\noindent 
1D: One-Dimensional \\
2D: Two-Dimensional \\
2DKP: Two-Dimensional Knapsack Problem \\
2S-SPR: Two-Stage Stochastic Programming with Recourse \\
3D: Three-Dimensional \\
ABM: Agent Based Modelling \\
ABS: Agent Based Simulation \\
ADP: Approximate Dynamic Programming \\
AFS: Airport Flight Scheduling \\
AHD: Attended Home Delivery \\
AHP: Analytic Hierarchy Process \\
AI: Artificial Intelligence \\
ANN: Artificial Neural Network \\
ANT: Actor Network Theory \\
AoA: Activity-on-Arc \\
AoN: Activity-on-Node \\
AOR: Airline Operations Recovery \\
AP: Assignment Problem \\
ARIMA: AutoRegressive Integrated Moving Average (model) \\
AR: Assurance Region \\
AR: Action Research \\
ARIMAX: AutoRegressive Integrated Moving Average with eXogenous variables (model) \\
AS: Autonomous System \\
ASP: Airline Schedule Planning \\
ATFM: Air Traffic Flow Management \\
B2B: Business-To-Business \\
B2C: Business-To-Consumer \\
B\&B: Branch-and-Bound \\
B\&C: Branch-and-Cut \\
B\&P: Branch-and-Price \\
BN: Bayesian Network \\
BOR: Behavioural OR \\
BPP: Bin Packing Problem \\
C\&P: Cutting and Packing \\
CBOR: Community-Based Operations Research \\
CDEA: Centralised DEA \\
CDM: Central Decision Maker or Collaborative Decision-Making \\
CLD: Causal Loop Diagram \\ 
CLSC: Closed-Loop Supply Chains \\
CM: Cellular Manufacturing \\
CNN: Convolutional Neural Network \\
CO: Combinatorial Optimisation \\
CODP: Customer Order decoupling point \\
ConFL: Connected Facility Location Problem \\
COR: Community Operational Research \\
CPM: Critical Path Method \\
CRPS: Continuous Ranked Probability Score \\
CST: Critical Systems Thinking \\
CSW: Common Set of Weights \\
CVaR: Conditional Value at Risk \\
CVRP: Capacitated Vehicle Routing Problem \\
DBN: Dynamic Bayesian Network \\
DC: Distribution Centre \\
DCT: Daily Contact Testing \\
DDF: Directional Distance Function \\
DEA: Data Envelopment Analysis \\
DEF: Deterministic Equivalent Formulation \\
DES: Discrete Event Simulation \\
DfT: Department for Transport \\
DHSC: Department of Health and Social Care \\
DMU: Decision Making Unit \\
DNDEA: Dynamic Network DEA \\
DNN: Deep Neural Network \\
DP: Dynamic Programming \\
DPSIR: Drivers, Pressures, State, Impact and Response \\
DS: Data Science \\
DSS: Decision Support Systems \\
EAT: Efficiency Analysis Trees \\
ED: Emergency Department \\
EMSR: Expected Marginal Seat Revenue \\
EOQ: Economic Order Quantity \\
ERP: Enterprise Resource Planning \\
ESICUP: EURO Special Interest Group on Cutting and Packing \\
EURO: European Operational Research Societies \\
EVP: Expected Value of Possession \\
FAM: Fleet Assignment Model \\
FIFO: First-In-First-Out \\
FMS: Flexible Manufacturing Systems \\
FPTAS: Fully Polynomial-Time Approximation Scheme \\
FSF: Full-State Feedback \\
FSO: Fixed-sum output \\
FTU: Facilities-Transformation-Usage (framework) \\
GIS: Geographic Information Systems \\
GLM: Generalised Linear Model \\
GMB: Group Model Building \\
GNN: Graphical Neural Network \\
GORS: Government Operational Research Service \\
GP: Gaussian Process \\
GPS: Global Positioning System \\
GPU: Graphics Processing Unit \\
GRASP: Greedy Randomised Adaptive Search Procedure \\
HJB: Hamilton-Jacobi-Bellman \\
HMT: His Majesty’s Treasury \\
HL: Humanitarian Logistics \\
HORAF: Heads of OR and Analytics Forum \\
IAM: Integrated Assessment Model \\
ICU: Intensive Care Unit \\
IGP: Interior Gateway Protocol \\
IHIP: Intangibility, Heterogeneity, Inseparability, and Perishability \\
IID: Independently and Identically Distributed \\
INFORMS: Institute for Management Science and Operations Research \\
INRC: International Nurse Rostering Competition \\
ILP: Integer Linear Problem \\
ILP: Integer Linear Programming \\
IoT: Internet of Things \\
IP: Integer Programming \\
IRP: Inventory-Routing Problem \\
JIT-MS: Just-In-Time Material System \\
KP: Knapsack Problem \\
LASSO: Least Absolute Shrinkage and Selection Operator \\
LCSA: Life Cycle Sustainability Assessment \\
LEAR: LASSO-Estimated AutoRegressive (model)\\
LP: Linear Programming \\
LQG: Linear Quadratic Gaussian \\
MAE: Mean Absolute Error \\
MAPE: Mean Absolute Percentage Error \\
MASE: Mean Absolute Scaled Error \\
MAUT: Multi-Attribute Utility Theory \\
MAVT: Multi-Attribute Value Theory \\
MBM: Market Based Measure \\
MC: Maximum Clique or Minimum Cut (problem)\\
MCDA: Multi-Criteria Decision Analysis \\
MCF: Minimum Cost Flow (problem) \\
MDP: Markov Decision Process \\
MF: Maximum Flow (problem) \\
MILP: Mixed-Integer Linear Programming \\
MINLP: Mixed-Integer NonLinear Programming \\
MIMO: Multi-Input-Multi-Output \\
MIP: Mixed-Integer Programming \\
ML: Machine Learning \\
MLPI: Malmquist Luenberger Productivity Indicator \\
MPC: Model Predictive Control \\
MPI: Malmquist Productivity Index \\
MRP: Material Requirement Planning \\
MRP: Multi-level Regression Post-stratification \\
NBEATS: Neural Basis Expansion Analysis for interpretable Time Series forecasting \\
NDEA: Network DEA \\
NDP: Neural Dynamic Programming \\
NFL: National Football League \\
NHS: National Health Service \\
NGO: Non-Governmental Organisation \\
OM: Operations Management \\
ONS: Office for National Statistics \\
OR: Operational (or Operations) Research \\
PA: Portfolio Analysis \\
PATAT: Practice and Theory of Automated Timetabling \\
PCR: Polymerase Chain Reaction \\
PERT: Project Evaluation and Review Technique \\
PESP: Periodic Event Scheduling Problem \\
PID: Proportional Integral Derivative \\
POMDP: Partially Observable Markov Decision Process \\
PPS: Production Possibility Set \\
PRA: Probabilistic Risk Assessment \\
PSM: Problem Structuring Method \\
PTAS: Polynomial-Time Approximation Scheme \\
QRA: Quantile Regression Averaging or Quantitative Risk Assessment \\
R\&D: Research and Development \\
RCPSP: Resource-Constrained Project Scheduling Problem \\
RES: Renewable Energy Sources \\
RFID: Radio-Frequency IDentification \\
RINS: Relaxation-Induced Neighbourhood Search \\
RL: Reinforcement Learning or Reverse Logistics\\
RM: Revenue Management \\
RMSE: Root Mean Squared Error \\
RNN: Recurrent Neural Network \\
SAA: Sample Average Approximation \\
SARF: Social Amplification of Risk Framework \\
SAT: SATisfiability (problem) \\
SCA: Strategic Choice Approach \\
SCM: Supply Chain Management \\
SD: Systems Dynamics \\
SDM: Structured Decision Making \\
SFA: Stochastic Frontier Analysis \\
SI: Systemic Intervention \\
SIS: Schools Infection Survey \\
SISO: Single-Input-Single-Output \\
SODA: Strategic Options Development and Analysis \\
SR: Segment Routing \\
SRCPSP: Stochastic Resource-Constrained Project Scheduling Problem \\
SSM: Soft Systems Methodology \\
SST: Shortest Spanning Trees \\
STP: Steiner Tree Problem (in graphs) \\
SVF: Support Vector Frontiers \\
SVM: Support Vector Machine \\
TE: Traffic Engineering \\
TFP: Total Factor Productivity \\
TPS: Toyota Production System \\
TSP: Travelling Salesman Problem \\
TTP: Travelling Tournament Problem \\
UDE: UnDesirable Effects \\
UFLP: Uncapacitated Facility Location Problem \\
VaR: Value at Risk \\
VAR: Vector AutoRegressive (model) \\
VMI: Vendor Managed Inventory \\
VPP: Virtual Power Plant \\
VRP: Vehicle Routing Problems \\
VSM: Viable Systems Model or Value Stream Map \\
VSS: Value of Stochastic Solution \\
VUCA: Volatile, Uncertain, Complex and Ambiguous \\
WHO: World Health Organisation \\

\clearpage

\bibliographystyle{elsarticle-harv}
\bibliography{refs}

\begin{thebibliography}{2196}
\expandafter\ifx\csname natexlab\endcsname\relax\def\natexlab#1{#1}\fi
\expandafter\ifx\csname url\endcsname\relax
  \def\url#1{\texttt{#1}}\fi
\expandafter\ifx\csname urlprefix\endcsname\relax\def\urlprefix{URL }\fi

\bibitem[{Aardal et~al.(2000)Aardal, Bixby, Hurkens, Lenstra, and
  Smeltink}]{Aa00_ALAL}
Aardal, K., Bixby, R., Hurkens, C., Lenstra, A., Smeltink, J., 2000. Market
  split and basis reduction: Towards a solution of the {C}ornu\'ejols-{D}awande
  instances. INFORMS Journal on Computing 12, 192--202.

\bibitem[{Abadi et~al.(2015)Abadi, Agarwal, Barham, Brevdo, Chen, Citro,
  Corrado, Davis, Dean, Devin, Ghemawat, Goodfellow, Harp, Irving, Isard, Jia,
  Jozefowicz, Kaiser, Kudlur, Levenberg, Man\'{e}, Monga, Moore, Murray, Olah,
  Schuster, Shlens, Steiner, Sutskever, Talwar, Tucker, Vanhoucke, Vasudevan,
  Vi\'{e}gas, Vinyals, Warden, Wattenberg, Wicke, Yu, and
  Zheng}]{tensorflow2015-whitepaper_LCAL}
Abadi, M., Agarwal, A., Barham, P., Brevdo, E., Chen, Z., Citro, C., Corrado,
  G.~S., Davis, A., Dean, J., Devin, M., Ghemawat, S., Goodfellow, I., Harp,
  A., Irving, G., Isard, M., Jia, Y., Jozefowicz, R., Kaiser, L., Kudlur, M.,
  Levenberg, J., Man\'{e}, D., Monga, R., Moore, S., Murray, D., Olah, C.,
  Schuster, M., Shlens, J., Steiner, B., Sutskever, I., Talwar, K., Tucker, P.,
  Vanhoucke, V., Vasudevan, V., Vi\'{e}gas, F., Vinyals, O., Warden, P.,
  Wattenberg, M., Wicke, M., Yu, Y., Zheng, X., 2015. {TensorFlow}: Large-scale
  machine learning on heterogeneous systems. Software available from
  tensorflow.org.

\bibitem[{Abbas et~al.(2017)Abbas, Tambe, and von
  Winterfeldt}]{Abbas2017-sv_KVRH}
Abbas, A.~E., Tambe, M., von Winterfeldt, D., 2017. Improving Homeland Security
  Decisions. Cambridge University Press, Cambridge.

\bibitem[{Abbey et~al.(2015)Abbey, Blackburn, and Guide}]{Abbey2015-dl_AM}
Abbey, J.~D., Blackburn, J., Guide, V. D.~R., 2015. Optimal pricing for new and
  remanufactured products. Journal of Operations Management 36, 130--146.

\bibitem[{Abbey et~al.(2017)Abbey, Kleber, Souza, and Voigt}]{Abbey2017-cn_AM}
Abbey, J.~D., Kleber, R., Souza, G.~C., Voigt, G., 2017. The role of perceived
  quality risk in pricing remanufactured products. Production and Operations
  Management 26~(1), 100--115.

\bibitem[{Abbott and Doucouliagos(2009)}]{Abbott2009-xv_JJ}
Abbott, M., Doucouliagos, C., 2009. Competition and efficiency: overseas
  students and technical efficiency in {A}ustralian and {N}ew {Z}ealand
  universities. Education Economics 17~(1), 31--57.

\bibitem[{Accorsi et~al.(2022)Accorsi, Lodi, and Vigo}]{ACCORSI2022229_CA_MB}
Accorsi, L., Lodi, A., Vigo, D., 2022. Guidelines for the computational testing
  of machine learning approaches to vehicle routing problems. Operations
  Research Letters 50~(2), 229--234.

\bibitem[{Accorsi and Vigo(2021)}]{AccorsiVigo2021_CA_MB}
Accorsi, L., Vigo, D., 2021. A fast and scalable heuristic for the solution of
  large-scale capacitated vehicle routing problems. Transportation Science
  55~(4), 832--856.

\bibitem[{Acito and Khatri(2014)}]{Acito2014-my_JEB}
Acito, F., Khatri, V., 2014. Business analytics: Why now and what next?
  Business Horizons 57~(5), 565--570.

\bibitem[{Ackermann and Eden(2011)}]{Ackermann2011-as_MYLW}
Ackermann, F., Eden, C., 2011. Making Strategy: Mapping Out Strategic Success,
  2nd Edition. SAGE, London.

\bibitem[{Ackermann and Howick(2022)}]{Ackermann2021_EA}
Ackermann, F., Howick, S., 2022. Experiences of mixed method or practitioners:
  Moving beyond a technical focus to insights relating to modelling teams.
  Journal of the Operational Research Society 73~(9), 1905--1918.

\bibitem[{Ackoff(1956)}]{Ackoff1956-qh}
Ackoff, R.~L., 1956. The development of operations research as a science.
  Operations Research 4~(3), 265--295.

\bibitem[{Ackoff(1970)}]{Ackoff1970-if_AJG}
Ackoff, R.~L., 1970. A black ghetto's research on a university. Operations
  Research 18~(5), 761--771.

\bibitem[{Ackoff(1977)}]{Ackoff1977-kh_MYLW}
Ackoff, R.~L., 1977. Optimization + objectivity = optout. European Journal of
  Operational Research 1~(1), 1--7.

\bibitem[{Ackoff(1979{\natexlab{a}})}]{Ackoff1979-zx_GL}
Ackoff, R.~L., 1979{\natexlab{a}}. The future of operational research is past.
  Journal of the Operational Research Society 30~(2), 93--104.

\bibitem[{Ackoff(1979{\natexlab{b}})}]{Ackoff1979-vt_GL}
Ackoff, R.~L., 1979{\natexlab{b}}. Resurrecting the future of operational
  research. Journal of the Operational Research Society 30~(3), 189--199.

\bibitem[{Ackoff(1981)}]{Ackoff1981-rt_GM}
Ackoff, R.~L., 1981. Creating the Corporate Future: Plan Or be Planned For.
  Wiley.

\bibitem[{Acworth et~al.(1998)Acworth, Broadie, and Glasserman}]{acworth98_MB}
Acworth, P.~A., Broadie, M., Glasserman, P., 1998. A comparison of some {Monte
  Carlo} and quasi {Monte Carlo} techniques for option pricing. In:
  Niederreiter, H., Hellekalek, P., Larcher, G., Zinterhof, P. (Eds.), Monte
  Carlo and Quasi-Monte Carlo Methods 1996. Springer-Verlag, New York, NY, pp.
  1--18.

\bibitem[{Adouani et~al.(2022)Adouani, Jarboui, and
  Masmoudi}]{Adouani2022_COIT}
Adouani, Y., Jarboui, B., Masmoudi, M., 2022. A matheuristic for the 0{–}1
  generalized quadratic multiple knapsack problem. Optimization Letters 16,
  37--58.

\bibitem[{Afsharian and Podinovski(2018)}]{Afsharian2018-ah_SL}
Afsharian, M., Podinovski, V.~V., 2018. A linear programming approach to
  efficiency evaluation in nonconvex metatechnologies. European Journal of
  Operational Research 268~(1), 268--280.

\bibitem[{Afzali et~al.(2012)Afzali, Karnon, and
  Gray}]{Haji_Ali_Afzali2012-jh_CV}
Afzali, H. H.~A., Karnon, J., Gray, J., 2012. A critical review of model-based
  economic studies of depression: modelling techniques, model structure and
  data sources. PharmacoEconomics 30~(6), 461--482.

\bibitem[{Agarwal and Ergun(2008)}]{Agarwal2008-vx_HP}
Agarwal, R., Ergun, {\"O}., 2008. Ship scheduling and network design for cargo
  routing in liner shipping. Transportation Science 42~(2), 175--196.

\bibitem[{Agasisti(2016)}]{Agasisti2016-bg_JJ}
Agasisti, T., 2016. Cost structure, productivity and efficiency of the
  {I}talian public higher education industry 2001--2011. International Review
  of Applied Economics 30~(1), 48--68.

\bibitem[{Agatz et~al.(2011)Agatz, Campbell, Fleischmann, and
  Savelsbergh}]{agatz2011time_CKTVW}
Agatz, N., Campbell, A., Fleischmann, M., Savelsbergh, M., 2011. Time slot
  management in attended home delivery. Transportation Science 45~(3),
  435--449.

\bibitem[{Agatz et~al.(2013)Agatz, Campbell, Fleischmann, Van~Nunen, and
  Savelsbergh}]{agatz2013revenue_CKTVW}
Agatz, N., Campbell, A.~M., Fleischmann, M., Van~Nunen, J., Savelsbergh, M.,
  2013. Revenue management opportunities for internet retailers. Journal of
  Revenue and Pricing Management 12~(2), 128--138.

\bibitem[{Agatz et~al.(2008)Agatz, Fleischmann, and
  Van~Nunen}]{agatz2008fulfillment_CKTVW}
Agatz, N., Fleischmann, M., Van~Nunen, J.~A., 2008. E-fulfillment and
  multi-channel distribution--a review. European Journal of Operational
  Research 187~(2), 339--356.

\bibitem[{Aghezzaf(2005)}]{aghezzaf2005capacity_SAA}
Aghezzaf, E., 2005. Capacity planning and warehouse location in supply chains
  with uncertain demands. Journal of the Operational Research Society 56~(4),
  453--462.

\bibitem[{Agrawal et~al.(2019)Agrawal, Avadhanula, Goyal, and
  Zeevi}]{agrawalMNLbanditDynamicLearning2019_AKSJF}
Agrawal, S., Avadhanula, V., Goyal, V., Zeevi, A., 2019. {{MNL-bandit}}: A
  dynamic learning approach to assortment selection. Operations Research
  67~(5), 1453--1485.

\bibitem[{Ahuja et~al.(1993)Ahuja, Magnanti, and Orlin}]{ahuja1993network_IL}
Ahuja, R.~K., Magnanti, T.~L., Orlin, J.~B., 1993. Network Flows: Theory,
  Algorithms, and Applications. Prentice Hall, Englewood Cliffs, NJ.

\bibitem[{Aigner et~al.(1977)Aigner, Lovell, and Schmidt}]{Aigner1977-df_JJ}
Aigner, D., Lovell, C. A.~K., Schmidt, P., 1977. Formulation and estimation of
  stochastic frontier production function models. Journal of Econometrics
  6~(1), 21--37.

\bibitem[{{AIMMS}(2022)}]{Aimms_CTCGE}
{AIMMS}, 2022. {AIMMS}.
\newline\urlprefix\url{https://www.aimms.com/}

\bibitem[{Ak{\c c}al{\i} et~al.(2009)Ak{\c c}al{\i}, {\c C}etinkaya, and
  {\"U}ster}]{Akcali2009-zg_AM}
Ak{\c c}al{\i}, E., {\c C}etinkaya, S., {\"U}ster, H., 2009. Network design for
  reverse and closed-loop supply chains: An annotated bibliography of models
  and solution approaches. Networks 53~(3), 231--248.

\bibitem[{Akhmedov et~al.(2016)Akhmedov, Kwee, and
  Montemanni}]{akhmedov2016divide_COIT}
Akhmedov, M., Kwee, I., Montemanni, R., 2016. A divide and conquer matheuristic
  algorithm for the prize-collecting steiner tree problem. Computers \&
  Operations Research 70, 18--25.

\bibitem[{Akhtar et~al.(2015)Akhtar, Scarf, and Rasool}]{Akhtar2015-lj_IM}
Akhtar, S., Scarf, P., Rasool, Z., 2015. Rating players in test match cricket.
  Journal of the Operational Research Society 66~(4), 684--695.

\bibitem[{{Al Hajj Hassan} et~al.(2022){Al Hajj Hassan}, Hewitt, and
  Mahmassani}]{ALHAJJHASSAN2022102885_MH}
{Al Hajj Hassan}, L., Hewitt, M., Mahmassani, H.~S., 2022. Daily load planning
  under different autonomous truck deployment scenarios. Transportation
  Research Part E: Logistics and Transportation Review 166, 102885.

\bibitem[{Al-Kanj et~al.(2020)Al-Kanj, Nascimento, and
  Powell}]{Al-Kanj2020-rn_HL}
Al-Kanj, L., Nascimento, J., Powell, W.~B., 2020. Approximate dynamic
  programming for planning a ride-hailing system using autonomous fleets of
  electric vehicles. European Journal of Operational Research 284~(3),
  1088--1106.

\bibitem[{Albano and Sapuppo(1980)}]{Albano1980-jw_JB}
Albano, A., Sapuppo, G., 1980. Optimal allocation of {Two-Dimensional}
  irregular shapes using heuristic search methods. IEEE Transactions on
  Systems, Man, and Cybernetics 10~(5), 242--248.

\bibitem[{Albayrak~{\"U}nal et~al.(2023)Albayrak~{\"U}nal, Erkayman, and
  Usanmaz}]{albayrak2023applications_TVWCK}
Albayrak~{\"U}nal, {\"O}., Erkayman, B., Usanmaz, B., 2023. Applications of
  artificial intelligence in inventory management: A systematic review of the
  literature. Archives of Computational Methods in Engineering 30~(4),
  2605--2625.

\bibitem[{Alexander(1992)}]{Ale92_JH}
Alexander, C., 1992. The {Kalai--Smorodinsky} bargaining solution in wage
  negotiations. Journal of the Operational Research Society 43, 779--786.

\bibitem[{Alfandari et~al.(2022)Alfandari, Ljubi{\'c}, and
  da~Silva}]{alfandari2022tailored_BKOK}
Alfandari, L., Ljubi{\'c}, I., da~Silva, M. D.~M., 2022. A tailored benders
  decomposition approach for last-mile delivery with autonomous robots.
  European Journal of Operational Research 299~(2), 510--525.

\bibitem[{Alfieri et~al.(2006)Alfieri, Groot, Kroon, and
  Schrijver}]{Alfieri2006_DC}
Alfieri, A., Groot, R., Kroon, L., Schrijver, A., 2006. Efficient circulation
  of railway rolling stock. Transportation Science 40~(3), 378--391.

\bibitem[{Algaba et~al.(2020)Algaba, Ardia, Bluteau, Borms, and
  Boudt}]{algaba20_MB}
Algaba, A., Ardia, D., Bluteau, K., Borms, S., Boudt, K., 2020. Econometrics
  meets sentiment: An overview of methodology and applications. Journal of
  Economic Surveys 34~(3), 512--547.

\bibitem[{Ali et~al.(2022)Ali, Ramos, Carravilla, and Oliveira}]{Ali2022-ow_JB}
Ali, S., Ramos, A.~G., Carravilla, M.~A., Oliveira, J.~F., 2022. On-line
  three-dimensional packing problems: A review of off-line and on-line solution
  approaches. Computers \& Industrial Engineering 168, 108122.

\bibitem[{Allen et~al.(1997)Allen, Athanassopoulos, Dyson, and
  Thanassoulis}]{Allen1997-fa_SL}
Allen, R., Athanassopoulos, A., Dyson, R.~G., Thanassoulis, E., 1997. Weights
  restrictions and value judgements in data envelopment analysis: Evolution,
  development and future directions. Annals of Operations Research 73, 13--34.

\bibitem[{Allon and Van~Mieghem(2010)}]{Allon2010_SMD}
Allon, G., Van~Mieghem, J.~A., 2010. Global dual sourcing: Tailored base-surge
  allocation to near- and offshore production. Management Science 56~(1),
  110--124.

\bibitem[{Almgren and Chriss(2001)}]{almgren01_MB}
Almgren, R., Chriss, N., 2001. Optimal execution of portfolio transactions.
  Journal of Risk 3, 5--40.

\bibitem[{Altay and Green~III(2006)}]{altay2006or_BYKOK}
Altay, N., Green~III, W.~G., 2006. {OR/MS} research in disaster operations
  management. European Journal of Operational Research 175~(1), 475--493.

\bibitem[{Altay and Labonte(2014)}]{altay2014challenges_BYKOK}
Altay, N., Labonte, M., 2014. Challenges in humanitarian information management
  and exchange: evidence from {H}aiti. Disasters 38~(s1), S50--S72.

\bibitem[{Alt{\i}n et~al.(2013)Alt{\i}n, Fortz, Thorup, and
  {\"U}mit}]{altn.fortz.ea:13_BF}
Alt{\i}n, A., Fortz, B., Thorup, M., {\"U}mit, H., 2013. Intra-domain traffic
  engineering with shortest path routing protocols. Annals of Operations
  Research 204~(1), 65--95.

\bibitem[{Alumur and Kara(2008)}]{alumur.kara:08_BF}
Alumur, S., Kara, B.~Y., 2008. Network hub location problems: The state of the
  art. European Journal of Operational Research 190, 1--21.

\bibitem[{Alumur et~al.(2021)Alumur, Campbell, Contreras, Kara, Marianov, and
  O’Kelly}]{alumur2021perspectives_SAA}
Alumur, S.~A., Campbell, J.~F., Contreras, I., Kara, B.~Y., Marianov, V.,
  O’Kelly, M.~E., 2021. Perspectives on modeling hub location problems.
  European Journal of Operational Research 291~(1), 1--17.

\bibitem[{Alvarez et~al.(2011)Alvarez, Tsilingiris, Engebrethsen, and
  Kakalis}]{Alvarez2011-do_HP}
Alvarez, J.~F., Tsilingiris, P., Engebrethsen, E.~S., Kakalis, N. M.~P., 2011.
  Robust fleet sizing and deployment for industrial and independent bulk ocean
  shipping companies. INFOR: Information Systems and Operational Research
  49~(2), 93--107.

\bibitem[{{\'A}lvarez-Miranda et~al.(2020){\'A}lvarez-Miranda, Salgado-Rojas,
  Hermoso, Garcia-Gonzalo, and Weintraub}]{miranda2020_EA}
{\'A}lvarez-Miranda, E., Salgado-Rojas, J., Hermoso, V., Garcia-Gonzalo, J.,
  Weintraub, A., 2020. An integer programming method for the design of
  multi-criteria multi-action conservation plans. Omega 92, 102147.

\bibitem[{Alvarez-Valdes et~al.(2013)Alvarez-Valdes, Martinez, and
  Tamarit}]{Alvarez-Valdes2013-fs_JB}
Alvarez-Valdes, R., Martinez, A., Tamarit, J.~M., 2013. A branch \& bound
  algorithm for cutting and packing irregularly shaped pieces. International
  Journal of Production Economics 145~(2), 463--477.

\bibitem[{Alyahyan and D{\"u}steg{\"o}r(2020)}]{Alyahyan2020-eh_JJ}
Alyahyan, E., D{\"u}steg{\"o}r, D., 2020. Predicting academic success in higher
  education: literature review and best practices. International Journal of
  Educational Technology in Higher Education 17, 3.

\bibitem[{{Amazon last mile routing}(2021)}]{AMAZON_CA_MB}
{Amazon last mile routing}, 2021. Research challenge.
  \url{https://routingchallenge.mit.edu/}, accessed on 2021-09-14.

\bibitem[{Anandalingam(1987)}]{Anandalingam1987_EA}
Anandalingam, G., 1987. Asymmetric players and bargaining for profit shares in
  natural resource development. Management Science 33~(8), 1048--1057.

\bibitem[{Andersen et~al.(2007)Andersen, Vennix, Richardson, and
  Rouwette}]{Andersen2007-pg}
Andersen, D.~F., Vennix, J. A.~M., Richardson, G.~P., Rouwette, E. A. J.~A.,
  2007. Group model building: problem structuring, policy simulation and
  decision support. Journal of the Operational Research Society 58~(5),
  691--694.

\bibitem[{Anderson et~al.(2017)Anderson, Chen, and Shao}]{Anderson17_BC}
Anderson, E., Chen, B., Shao, L., 2017. Supplier competition with option
  contracts for discrete blocks of capacity. Operations Research 65~(4),
  952--967.

\bibitem[{Anderson et~al.(2022)Anderson, Chen, and Shao}]{Anderson22_BC}
Anderson, E., Chen, B., Shao, L., 2022. Capacity games with supply function
  competition. Operations Research 70~(4), 1969--1983.

\bibitem[{Andersson et~al.(2001)Andersson, Mausser, Rosen, and
  Uryasev}]{andersson01_MB}
Andersson, F., Mausser, H., Rosen, D., Uryasev, S., 2001. Credit risk
  optimization with conditional value-at-risk criterion. Mathematical
  Programming 89~(2), 273--291.

\bibitem[{Andersson et~al.(2011)Andersson, Duesund, and
  Fagerholt}]{Andersson2011-qg_HP}
Andersson, H., Duesund, J.~M., Fagerholt, K., 2011. Ship routing and scheduling
  with cargo coupling and synchronization constraints. Computers \& Industrial
  Engineering 61~(4), 1107--1116.

\bibitem[{Andersson et~al.(2015)Andersson, Fagerholt, and
  Hobbesland}]{Andersson2015-rj_HP}
Andersson, H., Fagerholt, K., Hobbesland, K., 2015. Integrated maritime fleet
  deployment and speed optimization: Case study from {RoRo} shipping. Computers
  \& Operations Research 55, 233--240.

\bibitem[{Andras(2010)}]{Omnet_CTCGE}
Andras, V., 2010. Omnet++. In: Klaus, W., Mesut, G., James, G. (Eds.), Modeling
  and Tools for Network Simulation. Springer, Berlin, Heidelberg, pp. 35--59.

\bibitem[{Angelelli et~al.(2010)Angelelli, Mansini, and
  Speranza}]{angelelli2010kernel_COIT}
Angelelli, E., Mansini, R., Speranza, M.~G., 2010. Kernel search: A general
  heuristic for the multi-dimensional knapsack problem. Computers \& Operations
  Research 37~(11), 2017--2026.

\bibitem[{Angelus(2023)}]{angelus2023generalization_JSS}
Angelus, A., 2023. Generalizations of the clark-scarf model and analysis. In:
  Song, J.-S. (Ed.), Research Handbook on Inventory Management. Edward Elgar
  Publishing.

\bibitem[{Annabi et~al.(2012)Annabi, Breton, and Fran{\c{c}}ois}]{annabi12_MB}
Annabi, A., Breton, M., Fran{\c{c}}ois, P., 2012. Resolution of financial
  distress under chapter 11. Journal of Economic Dynamics and Control 36~(12),
  1867--1887.

\bibitem[{Antunes et~al.(2019)Antunes, Vaze, and
  Antunes}]{antunes2019robust_VLVV}
Antunes, D., Vaze, V., Antunes, A.~P., 2019. A robust pairing model for airline
  crew scheduling. Transportation Science 53~(6), 1751--1771.

\bibitem[{Aparicio et~al.(2017)Aparicio, Crespo-Cebada, Pedraja-Chaparro, and
  Sant{\'\i}n}]{Aparicio2017-dh_SL}
Aparicio, J., Crespo-Cebada, E., Pedraja-Chaparro, F., Sant{\'\i}n, D., 2017.
  Comparing school ownership performance using a pseudo-panel database: A
  {M}almquist-type index approach. European Journal of Operational Research
  256~(2), 533--542.

\bibitem[{Aparicio et~al.(2007)Aparicio, Ruiz, and
  Sirvent}]{Aparicio2007-fb_SL}
Aparicio, J., Ruiz, J.~L., Sirvent, I., 2007. Closest targets and minimum
  distance to the {P}areto-efficient frontier in {DEA}. Journal of Productivity
  Analysis 28~(3), 209--218.

\bibitem[{Applegate et~al.(2007)Applegate, Bixby, Chv\'atal, and
  Cook}]{ABCC07_SMPT}
Applegate, D., Bixby, R., Chv\'atal, V., Cook, W., 2007. The Traveling Salesman
  Problem: A Computational Study. Princeton University Press, Princeton.

\bibitem[{Applegate et~al.(2011)Applegate, Bixby, Chvátal, and
  Cook}]{Applegate-et-al:2011_IL}
Applegate, D.~L., Bixby, R.~E., Chvátal, V., Cook, W.~J., 2011. The Traveling
  Salesman Problem. Princeton University Press, Princeton.

\bibitem[{Arana-Jim{\'e}nez et~al.(2022)Arana-Jim{\'e}nez, S{\'a}nchez-Gil,
  Lozano, and Younesi}]{Arana-Jimenez2022-nv_SL}
Arana-Jim{\'e}nez, M., S{\'a}nchez-Gil, M.~C., Lozano, S., Younesi, A., 2022.
  Efficiency assessment using fuzzy production possibility set and enhanced
  {R}ussell {G}raph measure. Computational and Applied Mathematics 41~(2), 79.

\bibitem[{Ara{\'u}jo et~al.(2020)Ara{\'u}jo, Santos, Marques, and
  Barbosa-Povoa}]{Araujo2020_JLYHK}
Ara{\'u}jo, A.~M., Santos, D., Marques, I., Barbosa-Povoa, A., 2020. Blood
  supply chain: a two-stage approach for tactical and operational planning. OR
  Spectrum 43, 1023--1053.

\bibitem[{Archetti and Bertazzi(2021)}]{archetti2021recent_JLYHK}
Archetti, C., Bertazzi, L., 2021. Recent challenges in routing and inventory
  routing: E-commerce and last-mile delivery. Networks 77~(2), 255--268.

\bibitem[{Archetti et~al.(2007)Archetti, Bertazzi, Laporte, and
  Speranza}]{archetti2007branch_JLYHK}
Archetti, C., Bertazzi, L., Laporte, G., Speranza, M.~G., 2007. A
  branch-and-cut algorithm for a vendor-managed inventory-routing problem.
  Transportation Science 41~(3), 382--391.

\bibitem[{Archetti et~al.(2017)Archetti, Boland, and
  Speranza}]{archetti2017matheuristic_COIT}
Archetti, C., Boland, N., Speranza, M.~G., 2017. A matheuristic for the
  multivehicle inventory routing problem. INFORMS Journal on Computing 29~(3),
  377--387.

\bibitem[{Archetti et~al.(2015)Archetti, Corber{\' a}n, Plana, Sanchis, and
  Speranza}]{Archetti2015_COIT}
Archetti, C., Corber{\' a}n, A., Plana, I., Sanchis, J.~M., Speranza, M.~G.,
  2015. A matheuristic for the team orienteering arc routing problem. European
  Journal of Operational Research 245, 392--401.

\bibitem[{Argyris et~al.(2022)Argyris, Karsu, and Yavuz}]{ArgKarYav22_JH}
Argyris, N., Karsu, {\"{O}}., Yavuz, M., 2022. Fair resource allocation:
  {Using} welfare-based dominance constraints. European Journal of Operational
  Research 297~(2), 560--578.

\bibitem[{Ariely and Simonson(2003)}]{Ariely2003-eb_BC}
Ariely, D., Simonson, I., 2003. Buying, bidding, playing, or competing? value
  assessment and decision dynamics in online auctions. Journal of consumer
  psychology: the official journal of the Society for Consumer Psychology
  13~(1), 113--123.

\bibitem[{Aringhieri et~al.(2017)Aringhieri, Bruni, Khodaparasti, and van
  Essen}]{aringhieri2017emergency_BYKOK}
Aringhieri, R., Bruni, M.~E., Khodaparasti, S., van Essen, J.~T., 2017.
  Emergency medical services and beyond: Addressing new challenges through a
  wide literature review. Computers \& Operations Research 78, 349--368.

\bibitem[{Arnold and S{\"o}rensen(2019)}]{arnold2019makes_CA_MB}
Arnold, F., S{\"o}rensen, K., 2019. What makes a {VRP} solution good? {T}he
  generation of problem-specific knowledge for heuristics. Computers \&
  Operations Research 106, 280--288.

\bibitem[{Arnold et~al.(2004)Arnold, Peeters, and
  Thomas}]{arnold2004modelling_MH}
Arnold, P., Peeters, D., Thomas, I., 2004. Modelling a rail/road intermodal
  transportation system. Transportation Research Part E: Logistics and
  Transportation Review 40~(3), 255--270.

\bibitem[{Arrow et~al.(1951)Arrow, Harris, and Marschak}]{arrow1951optimal_JSS}
Arrow, K.~J., Harris, T., Marschak, J., 1951. Optimal inventory policy.
  Econometrica 19~(3), 250--272.

\bibitem[{Arrow et~al.(1958)Arrow, Karlin, and Scarf}]{arrow1958studies_JSS}
Arrow, K.~J., Karlin, S., Scarf, H.~E., 1958. Studies in the mathematical
  theory of inventory and production. Stanford University Press.

\bibitem[{Artzner et~al.(1999)Artzner, Delbaen, Eber, and Heath}]{artzner99_MB}
Artzner, P., Delbaen, F., Eber, J.-M., Heath, D., 1999. Coherent measures of
  risk. Mathematical Finance 9~(3), 203--228.

\bibitem[{Arvan et~al.(2019)Arvan, Fahimnia, Reisi, and
  Siemsen}]{Arvan2019-xy_AFRH}
Arvan, M., Fahimnia, B., Reisi, M., Siemsen, E., 2019. Integrating human
  judgement into quantitative forecasting methods: A review. Omega 86,
  237--252.

\bibitem[{Asghari et~al.(2023)Asghari, Jaber, and
  Mirzapour~{Al-e-hashem}}]{Asghari2022-pq_HP}
Asghari, M., Jaber, M.~Y., Mirzapour~{Al-e-hashem}, S. M.~J., 2023.
  Coordinating vessel recovery actions: Analysis of disruption management in a
  liner shipping service. European Journal of Operational Research 307~(2),
  627--644, , DOI: 10.1016/j.ejor.2022.08.039.

\bibitem[{Asmuni et~al.(2009)Asmuni, Burke, Garibaldi, McCollum, and
  Parkes}]{Asmuni2009-yb_JJ}
Asmuni, H., Burke, E.~K., Garibaldi, J.~M., McCollum, B., Parkes, A.~J., 2009.
  An investigation of fuzzy multiple heuristic orderings in the construction of
  university examination timetables. Computers \& Operations Research 36~(4),
  981--1001.

\bibitem[{{\AA}str{\"o}m(2012)}]{aastrom2012_XW}
{\AA}str{\"o}m, K.~J., 2012. Introduction to Stochastic Control Theory. Courier
  Corporation.

\bibitem[{{\AA}str{\"o}m and Kumar(2014)}]{aastrom2014_XW}
{\AA}str{\"o}m, K.~J., Kumar, P.~R., 2014. Control: A perspective. Automatica
  50~(1), 3--43.

\bibitem[{{\AA}str{\"o}m and Wittenmark(2013)}]{aastrom2013_XW}
{\AA}str{\"o}m, K.~J., Wittenmark, B., 2013. Adaptive control. Courier
  Corporation.

\bibitem[{Atan et~al.(2017)Atan, Ahmadi, Stegehuis, de~Kok, and
  Adan}]{atan2017assemble_JSS}
Atan, Z., Ahmadi, T., Stegehuis, C., de~Kok, T., Adan, I., 2017.
  Assemble-to-order systems: A review. European Journal of Operational Research
  261~(3), 866--879.

\bibitem[{Atasu et~al.(2008)Atasu, Sarvary, and
  Van~Wassenhove}]{Atasu2008-za_AM}
Atasu, A., Sarvary, M., Van~Wassenhove, L.~N., 2008. Remanufacturing as a
  marketing strategy. Management Science 54~(10), 1731--1746.

\bibitem[{Athanasopoulos et~al.(2020)Athanasopoulos, Gamakumara, Panagiotelis,
  Hyndman, and Affan}]{Athanasopoulos2020-cx_FP}
Athanasopoulos, G., Gamakumara, P., Panagiotelis, A., Hyndman, R.~J., Affan,
  M., 2020. Hierarchical forecasting. In: Fuleky, P. (Ed.), Macroeconomic
  Forecasting in the Era of Big Data: Theory and Practice. Springer, pp.
  689--719.

\bibitem[{Atkinson and Gary(2016)}]{Atkinson2016-wl_MCJM}
Atkinson, S., Gary, M.~S., 2016. Mergers and acquisitions: Modeling decision
  making in integration projects. In: Kunc, M., Malpass, J., White, L. (Eds.),
  Behavioral Operational Research: Theory, Methodology and Practice. Palgrave
  Macmillan, London, pp. 319--336.

\bibitem[{Ausiello et~al.(1999)Ausiello, Crescenzi, Gambosi, Kann,
  Marchetti-Spaccamela, and Protasi}]{Ausi13_UPCT}
Ausiello, G., Crescenzi, P., Gambosi, G., Kann, V., Marchetti-Spaccamela, A.,
  Protasi, M., 1999. Complexity and Approximation: Combinatorial Optimization
  Problems and Their Approximability Properties. Springer, Berlin.

\bibitem[{Aven(2015)}]{Aven2015-wq_TC}
Aven, T., 2015. Risk Analysis. John Wiley \& Sons.

\bibitem[{Aven(2020)}]{Aven2020-tq_TC}
Aven, T., 2020. Three influential risk foundation papers from the 80s and 90s:
  Are they still state-of-the-art? Reliability Engineering \& System Safety
  193, 106680.

\bibitem[{Aviv(2003)}]{Aviv2003_XW}
Aviv, Y., 2003. A time-series framework for supply-chain inventory management.
  Operations Research 51~(2), 210--227.

\bibitem[{Axs{\"a}ter(1993)}]{axsater1993continuous_JSS}
Axs{\"a}ter, S., 1993. Continuous review policies for multi-level inventory
  systems with stochastic demand. Handbooks in Operations Research and
  Management Science 4, 175--197.

\bibitem[{Axs{\"a}ter(1996)}]{axsater1996using_JSS}
Axs{\"a}ter, S., 1996. Using the deterministic eoq formula in stochastic
  inventory control. Management Science 42~(6), 830--834.

\bibitem[{Axs{\"a}ter(2003)}]{axsater2003supply_JSS}
Axs{\"a}ter, S., 2003. Supply chain operations: Serial and distribution
  inventory systems. In: Graves, S.~C., de~Kok, A.~G. (Eds.), Handbooks in
  Operations Research and Management Science. Vol.~11. Elsevier, pp. 525--559.

\bibitem[{Axs{\"a}ter(2006)}]{axsater2006inventory_JSS}
Axs{\"a}ter, S., 2006. Inventory Control. Springer.

\bibitem[{Axs{\"a}ter and Rosling(1993)}]{axsater1993installation_JSS}
Axs{\"a}ter, S., Rosling, K., 1993. Installation vs. echelon stock policies for
  multilevel inventory control. Management Science 39~(10), 1274--1280.

\bibitem[{Ayhan and Baccelli(2001)}]{ayhan1_HATI}
Ayhan, H., Baccelli, F., 2001. Expansions for joint laplace transform of
  stationary waiting times in (max,+)-linear systems with poisson input.
  Queueing Systems 37~(1), 291--328.

\bibitem[{Ayhan et~al.(2004)Ayhan, Palmowski, and Schegel}]{ayhan2_HATI}
Ayhan, H., Palmowski, Z., Schegel, S., 2004. Cyclic queueing networks with
  subexponential service trimes. Journal of Applied Probability 41~(3),
  291--301.

\bibitem[{Azoury(1985)}]{azoury1985bayes_JSS}
Azoury, K.~S., 1985. Bayes solution to dynamic inventory models under unknown
  demand distribution. Management Science 31~(9), 1150--1160.

\bibitem[{Baar et~al.(1999)Baar, Brucker, and Knust}]{Baar1999-ni_WH_ED}
Baar, T., Brucker, P., Knust, S., 1999. Tabu search algorithms and lower bounds
  for the {Resource-Constrained} project scheduling problem. In: Vo{\ss}, S.,
  Martello, S., Osman, I.~H., Roucairol, C. (Eds.), {Meta-Heuristics}: Advances
  and Trends in Local Search Paradigms for Optimization. Springer US, Boston,
  MA, pp. 1--18.

\bibitem[{Babich et~al.(2021)Babich, Birge, et~al.}]{babich21_MB}
Babich, V., Birge, J.~R., et~al., 2021. The interface of finance, operations,
  and risk management. Foundations and Trends{\textregistered} in Technology,
  Information and Operations Management 15~(1--2), 1--203.

\bibitem[{Babich and Hilary(2020)}]{Babich2020_SMD}
Babich, V., Hilary, G., 2020. {OM Forum—D}istributed ledgers and operations:
  {W}hat operations management researchers should know about blockchain
  technology. Manufacturing \& Service Operations Management 22~(2), 223--240.

\bibitem[{Baboolal et~al.(2012)Baboolal, Griffiths, Knight, Nelson, Voake, and
  Williams}]{Baboolal2012-ci_CV}
Baboolal, K., Griffiths, J.~D., Knight, V.~A., Nelson, A.~V., Voake, C.,
  Williams, J.~E., 2012. How efficient can an emergency unit be? {A} perfect
  world model. Emergency Medicine Journal 29~(12), 972--977.

\bibitem[{Baccelli et~al.(1992)Baccelli, Cohen, Olsder, and
  Quadrat}]{francois_HATI}
Baccelli, F., Cohen, G., Olsder, G.~J., Quadrat, J.-P., 1992. Synchronization
  and linearity: an algebra for discrete event systems. John Wiley \& Sons Ltd.

\bibitem[{Baccelli et~al.(1997)Baccelli, Hasenfuss, and Schmidt}]{sven_HATI}
Baccelli, F., Hasenfuss, S., Schmidt, V., 1997. Transient and stationary
  waiting times in (max,+)-linear systems with poisson input. Queueing Systems
  26~(3), 301--342.

\bibitem[{Baccelli et~al.(1999)Baccelli, Schlegel, and Schmidt}]{sabine_HATI}
Baccelli, F., Schlegel, S., Schmidt, V., 1999. Asymptotics of stochastic
  networks with subexponential service times. Queueing Systems 33, 205--232.

\bibitem[{Baccelli and Schmidt(1996)}]{baccellischmidt_HATI}
Baccelli, F., Schmidt, V., 1996. Taylor series expansions for poisson-driven
  $(max,+ )$-linear systems. The Annals of Applied Probability 6~(1), 138--185.

\bibitem[{Bacciotti and Rosier(2005)}]{Bacciotti2005_XW}
Bacciotti, A., Rosier, L., 2005. Liapunov functions and stability in control
  theory. Springer Science \& Business Media.

\bibitem[{Ba\c{s}ar et~al.(2018)Ba\c{s}ar, Haurie, and Zaccour}]{Baetal2018_GZ}
Ba\c{s}ar, T., Haurie, A., Zaccour, G., 2018. Nonzero-sum differential games.
  In: Ba\c{s}ar, T., Zaccour, G. (Eds.), Handbook of Dynamic Game Theory.
  Springer, Cham, pp. 61--110.

\bibitem[{Ba\c{s}ar and Olsder(1999)}]{Baol1999_GZ}
Ba\c{s}ar, T., Olsder, G., 1999. Dynamic noncooperative game theory, 2nd
  Edition. SIAM, Philadelphia.

\bibitem[{Bailey(1952)}]{Bailey1952-xt_CV}
Bailey, N. T.~J., 1952. A study of queues and appointment systems in hospital
  out-patient departments, with special reference to waiting-times. Journal of
  the Royal Statistical Society 14~(2), 185--199.

\bibitem[{Baker and McHale(2015)}]{Baker2015-cv_IM}
Baker, R.~D., McHale, I.~G., 2015. Deterministic evolution of strength in
  multiple comparisons models: Who is the greatest golfer? Scandinavian Journal
  of Statistics, Theory and Applications 42~(1), 180--196.

\bibitem[{Bakir et~al.(2021)Bakir, Erera, and Savelsbergh}]{Bakir2021_MH}
Bakir, I., Erera, A., Savelsbergh, M., 2021. Motor carrier service network
  design. In: Crainic, T.~G., Gendreau, M., Gendron, B. (Eds.), Network Design
  with Applications to Transportation and Logistics. Springer, pp. 427--467.

\bibitem[{Bakshi and Pinker(2018)}]{Bakshi2018-hw_KVRH}
Bakshi, N., Pinker, E., 2018. Public warnings in counterterrorism operations:
  Managing the {``Cry-Wolf''} effect when facing a strategic adversary.
  Operations Research 66~(4), 977--993.

\bibitem[{Balakrishnan and Altinkemer(1992)}]{balakrishnan.altinkemer:92_BF}
Balakrishnan, A., Altinkemer, K., 1992. Using a hop-constrained model to
  generate alternative communication network design. {ORSA} Journal on
  Computing 4~(2), 192--205.

\bibitem[{Balas(1965)}]{Ba65_ALAL}
Balas, E., 1965. An additive algorithm for solving linear programs with
  zero-one variables. Operations Research 13, 517--546.

\bibitem[{Balas(1971)}]{Ba71_ALAL}
Balas, E., 1971. Intersection cuts---a new type of cutting planes for integer
  programming. Operations Research 19, 19--39.

\bibitem[{Balas(1975)}]{Ba75_ALAL}
Balas, E., 1975. Facets of the knapsack polytope. Mathematical Programming 8,
  146--164.

\bibitem[{Balas(1979)}]{Ba79_ALAL}
Balas, E., 1979. Disjunctive programming. Annals of Discrete Mathematics 5,
  3--51.

\bibitem[{Balas et~al.(1993)Balas, Ceria, and Cornu{\'e}jols}]{BCC93_ALAL}
Balas, E., Ceria, S., Cornu{\'e}jols, G., 1993. A lift-and-project cutting
  plane algorithm for mixed 0--1 programs. Mathematical Programming 58,
  295--324.

\bibitem[{Balas et~al.(1996{\natexlab{a}})Balas, Ceria, and
  Cornu{\'e}jols}]{BCC96_ALAL}
Balas, E., Ceria, S., Cornu{\'e}jols, G., 1996{\natexlab{a}}. Mixed 0-1
  programming by lift-and-project in a branch-and-cut framework. Management
  Science 42, 1229--1246.

\bibitem[{Balas et~al.(1996{\natexlab{b}})Balas, Ceria, Cornu{\'e}jols, and
  Natraj}]{Ba96_ALAL}
Balas, E., Ceria, S., Cornu{\'e}jols, G., Natraj, N., 1996{\natexlab{b}}.
  Gomory cuts revisited. Operations Research Letters 19, 1--9.

\bibitem[{Balas and Martin(1980)}]{BM80_ALAL}
Balas, E., Martin, C., 1980. Pivot and complement---a heuristic for 0-1
  programming. Management Science 26, 86--96.

\bibitem[{Balcik and Beamon(2008)}]{balcik2008facility_BYKOK}
Balcik, B., Beamon, B.~M., 2008. Facility location in humanitarian relief.
  International Journal of Logistics 11~(2), 101--121.

\bibitem[{Baldacci et~al.(2004)Baldacci, Hadjiconstantinou, and
  Mingozzi}]{Baldacci2Flow2004_CA_MB}
Baldacci, R., Hadjiconstantinou, E., Mingozzi, A., 2004. An exact algorithm for
  the capacitated vehicle routing problem based on a two-commodity network flow
  formulation. Operations Research 52~(5), 723--738.

\bibitem[{Baldacci et~al.(2011)Baldacci, Mingozzi, and
  Roberti}]{Roberti2011_CA_MB}
Baldacci, R., Mingozzi, A., Roberti, R., 2011. New route relaxation and pricing
  strategies for the vehicle routing problem. Operations Research 59~(5),
  1269--1283.

\bibitem[{Ballest{\'\i}n(2007{\natexlab{a}})}]{Ballestin2007-xx_WH_ED}
Ballest{\'\i}n, F., 2007{\natexlab{a}}. A genetic algorithm for the resource
  renting problem with minimum and maximum time lags. In: Cotta, C., Hemert, J.
  (Eds.), Evolutionary Computation in Combinatorial Optimization. Springer,
  Berlin, Heidelberg, pp. 25--35.

\bibitem[{Ballest{\'\i}n(2007{\natexlab{b}})}]{Ballestin2007-ma_HL}
Ballest{\'\i}n, F., 2007{\natexlab{b}}. When it is worthwhile to work with the
  stochastic {RCPSP}? Journal of Scheduling 10~(3), 153--166.

\bibitem[{Ballest{\'\i}n and Leus(2009)}]{Ballestin2009-rp_HL}
Ballest{\'\i}n, F., Leus, R., 2009. Resource-constrained project scheduling for
  timely project completion with stochastic activity durations. Production and
  Operations Management 18~(4), 459--474.

\bibitem[{Balon et~al.(2006)Balon, Skiv\'ee, and Leduc}]{balon.skivee.ea:06_BF}
Balon, S., Skiv\'ee, F., Leduc, G., 2006. How well do traffic engineering
  objective functions meet {TE} requirements? In: Proceedings of IFIP
  Networking 2006, Coimbra. Vol. 3976. Springer LNCS, pp. 75--86.

\bibitem[{Bandi and Bertsimas(2014)}]{bandi14_MB}
Bandi, C., Bertsimas, D., 2014. Robust option pricing. European Journal of
  Operational Research 239~(3), 842--853.

\bibitem[{Bandi and Gupta(2020)}]{Bandi2020-xy_CV}
Bandi, C., Gupta, D., 2020. Operating room staffing and scheduling.
  Manufacturing \& Service Operations Management 22~(5), 958--974.

\bibitem[{Banerjee et~al.(2022)Banerjee, Erera, and
  Toriello}]{banerjee2022fleet_CKTVW}
Banerjee, D., Erera, A.~L., Toriello, A., 2022. Fleet sizing and service region
  partitioning for same-day delivery systems. Transportation Science 56~(5),
  1327--1347.

\bibitem[{Banker et~al.(1984)Banker, Charnes, and Cooper}]{Banker1984-no_SL}
Banker, R.~D., Charnes, A., Cooper, W.~W., 1984. Some models for estimating
  technical and scale inefficiencies in data envelopment analysis. Management
  Science 30~(9), 1078--1092.

\bibitem[{Banker and Morey(1986)}]{Banker1986-ai_SL}
Banker, R.~D., Morey, R.~C., 1986. Efficiency analysis for exogenously fixed
  inputs and outputs. Operations Research 34~(4), 513--521.

\bibitem[{Banks et~al.(2022)Banks, Gallego, Naveiro, and
  R{\'\i}os~Insua}]{Banks2022-fi_TC}
Banks, D., Gallego, V., Naveiro, R., R{\'\i}os~Insua, D., 2022. Adversarial
  risk analysis: An overview. Wiley Interdisciplinary Reviews: Computational
  Statistics 14, 1530.

\bibitem[{Banks et~al.(2004)Banks, Carson, Nelson, and Nicol}]{Banks2004_CC}
Banks, J., Carson, J., Nelson, B.~L., Nicol, D., 2004. {Discrete-Event} System
  Simulation, 4th Edition. Prentice Hall.

\bibitem[{Barbosa-P{\'o}voa et~al.(2018)Barbosa-P{\'o}voa, da~Silva, and
  Carvalho}]{Barbosa-Povoa2018-dw_JL}
Barbosa-P{\'o}voa, A.~P., da~Silva, C., Carvalho, A., 2018. Opportunities and
  challenges in sustainable supply chain: An operations research perspective.
  European Journal of Operational Research 268~(2), 399--431.

\bibitem[{Bardossy and Raghavan(2010)}]{bardossy.raghavan:10_BF}
Bardossy, M., Raghavan, S., 2010. Dual-based local search for the connected
  facility location and related problems. INFORMS Journal on Computing 22,
  584--602.

\bibitem[{Barnhart et~al.(2009)Barnhart, Farahat, and
  Lohatepanont}]{barnhart2009airline_VLVV}
Barnhart, C., Farahat, A., Lohatepanont, M., 2009. Airline fleet assignment
  with enhanced revenue modeling. Operations Research 57~(1), 231--244.

\bibitem[{Barnhart et~al.(2014)Barnhart, Fearing, and
  Vaze}]{barnhart2014modeling_VLVV}
Barnhart, C., Fearing, D., Vaze, V., 2014. Modeling passenger travel and delays
  in the national air transportation system. Operations Research 62~(3),
  580--601.

\bibitem[{Barnhart et~al.(2002)Barnhart, Kniker, and
  Lohatepanont}]{barnhart2002itinerary_VLVV}
Barnhart, C., Kniker, T.~S., Lohatepanont, M., 2002. Itinerary-based airline
  fleet assignment. Transportation Science 36~(2), 199--217.

\bibitem[{Barnhart and Schneur(1996)}]{barnhart1996air_MH}
Barnhart, C., Schneur, R.~R., 1996. Air network design for express shipment
  service. Operations Research 44~(6), 852--863.

\bibitem[{Barocas et~al.(2019)Barocas, Hardt, and Narayanan}]{barocas19_LCAL}
Barocas, S., Hardt, M., Narayanan, A., 2019. Fairness and Machine Learning:
  Limitations and Opportunities. fairmlbook.org,
  \url{http://www.fairmlbook.org}.

\bibitem[{Barrena et~al.(2014{\natexlab{a}})Barrena, Canca, Coelho, and
  Laporte}]{Barrena2014a_DC}
Barrena, E., Canca, D., Coelho, L.~C., Laporte, G., 2014{\natexlab{a}}. Exact
  formulations and algorithm for the train timetabling problem with dynamic
  demand. Computers \& Operations Research 44, 66--74.

\bibitem[{Barrena et~al.(2014{\natexlab{b}})Barrena, Canca, Coelho, and
  Laporte}]{Barrena2014_DC}
Barrena, E., Canca, D., Coelho, L.~C., Laporte, G., 2014{\natexlab{b}}.
  Single-line rail rapid transit timetabling under dynamic passenger demand.
  Transportation Research Part B: Methodological 70, 134--150.

\bibitem[{Baryannis et~al.(2019)Baryannis, Validi, Dani, and
  Antoniou}]{Baryannis2019-vp_MCJM}
Baryannis, G., Validi, S., Dani, S., Antoniou, G., 2019. Supply chain risk
  management and artificial intelligence: state of the art and future research
  directions. International Journal of Production Research 57~(7), 2179--2202.

\bibitem[{Baskett et~al.(1975)Baskett, Chandy, Muntz, and Palacios}]{bcmp_HATI}
Baskett, F., Chandy, K.~M., Muntz, R.~R., Palacios, F.~G., 1975. Open, closed,
  and mixed networks of queues with different classes of customers. Journal of
  the ACM 22~(2), 248--260.

\bibitem[{Basole et~al.(2022)Basole, Bendoly, Chandrasekaran, and
  Linderman}]{Basole2021-ma_MJE}
Basole, R., Bendoly, E., Chandrasekaran, A., Linderman, K., 2022. Visualization
  in operations management research. INFORMS Journal on Data Science 1~(2),
  172--187.

\bibitem[{Bates and Granger(1969)}]{Bates1969combination_FP}
Bates, J.~M., Granger, C.~W., 1969. The combination of forecasts. Journal of
  the Operational Research Society 20~(4), 451--468.

\bibitem[{Bates et~al.(2020)Bates, Cobo, Mari{\~n}o, and
  Wheeler}]{Bates2020-bc_JJ}
Bates, T., Cobo, C., Mari{\~n}o, O., Wheeler, S., 2020. Can artificial
  intelligence transform higher education? International Journal of Educational
  Technology in Higher Education 17~(1), 1--12.

\bibitem[{Baxter et~al.(2020)Baxter, Wilborn~Lagerman, and
  Keskinocak}]{8935364_BYKOK}
Baxter, A.~E., Wilborn~Lagerman, H.~E., Keskinocak, P., 2020. Quantitative
  modeling in disaster management: A literature review. IBM Journal of Research
  and Development 64~(1/2), 3:1--3:13.

\bibitem[{Bayram and Yaman(2018)}]{bayram2018shelter_SAA}
Bayram, V., Yaman, H., 2018. Shelter location and evacuation route assignment
  under uncertainty: A {B}enders decomposition approach. Transportation Science
  52~(2), 416--436.

\bibitem[{Bazaraa et~al.(2005)Bazaraa, Sherali, and Shetty}]{BSS2005_EAY}
Bazaraa, M.~S., Sherali, H.~D., Shetty, C.~M., 2005. Nonlinear Programming -
  Theory and Algorithms, 3rd Edition. Wiley, Hoboken, NJ.

\bibitem[{Beamon(1999)}]{Beamon1999_JLYHK}
Beamon, B.~M., 1999. Designing the green supply chain. Logistics Information
  Management 12~(6), 332--342.

\bibitem[{Beasley(1990)}]{beasley1990or_SAA}
Beasley, J.~E., 1990. {OR-Library}: distributing test problems by electronic
  mail. Journal of the Operational Research Society 41~(11), 1069--1072.

\bibitem[{Beer(1966)}]{Beer1966-aw_GM}
Beer, S., 1966. Decision and Control. Wiley, Chichester.

\bibitem[{Beer(1981)}]{Beer1981-fe_GM}
Beer, S., 1981. Brain of the firm, 2nd Edition. Wiley, Chichester.

\bibitem[{Beham(2020)}]{Simsharp_CTCGE}
Beham, A., 2020. {SimSharp}.
\newline\urlprefix\url{https://github.com/heal-research/SimSharp}

\bibitem[{Behl and Dutta(2019)}]{Behl2019_JLYHK}
Behl, A., Dutta, P., 2019. {Humanitarian supply chain management: a thematic
  literature review and future directions of research}. Annals of Operations
  Research 283~(1), 1001--1044.

\bibitem[{Bell et~al.(1988)Bell, Raiffa, and
  Tversky}]{bell1988descriptive_MESG}
Bell, D., Raiffa, H., Tversky, A., 1988. Descriptive, Normative, and
  Prescriptive Interactions in Decision Making. Cambridge University Press,
  Cambridge.

\bibitem[{Bell and Davison(2013)}]{Bell2013-eb_MJE}
Bell, E., Davison, J., 2013. Visual management studies: Empirical and
  theoretical approaches. International Journal of Management Reviews 15~(2),
  167--184.

\bibitem[{Bell and O'Keefe(1995)}]{Bell1995-lw_AFRH}
Bell, P.~C., O'Keefe, R.~M., 1995. An experimental investigation into the
  efficacy of visual interactive simulation. Management Science 41~(6),
  1018--1038.

\bibitem[{Bell(2012)}]{Bell2012-gn_MYLW}
Bell, S., 2012. {DPSIR=A} problem structuring method? {A}n exploration from the
  ``{I}magine'' approach. European Journal of Operational Research 222~(2),
  350--360.

\bibitem[{Bell et~al.(1983)Bell, Dalberto, Fisher, Greenfield, Jaikumar, Kedia,
  Mack, and Prutzman}]{bell1983improving_JLYHK}
Bell, W.~J., Dalberto, L.~M., Fisher, M.~L., Greenfield, A.~J., Jaikumar, R.,
  Kedia, P., Mack, R.~G., Prutzman, P.~J., 1983. Improving the distribution of
  industrial gases with an on-line computerized routing and scheduling
  optimizer. Interfaces 13~(6), 4--23.

\bibitem[{Bellenguez et~al.(2023)Bellenguez, Brauner, and
  Tsouki{\`a}s}]{Bellenguez2023-hc_JH}
Bellenguez, O., Brauner, N., Tsouki{\`a}s, A., 2023. Is there an ethical
  operational research practice? and what this implies for our research? EURO
  Journal on Decision Processes 11, 100029.

\bibitem[{Bellman(1953)}]{bellman1953introduction_DLLD}
Bellman, R., 1953. An introduction to the theory of dynamic programming. Tech.
  rep., RAND Corp Santa Monica CA.

\bibitem[{Bellman(1957)}]{Bellman1957-ue_HL}
Bellman, R., 1957. Dynamic Programming. Princeton University Press.

\bibitem[{Bellmore and Nemhauser(1968)}]{bellmore1968traveling_COIT}
Bellmore, M., Nemhauser, G.~L., 1968. The traveling salesman problem: {A}
  survey. Operations Research 16~(3), 538--558.

\bibitem[{Belobaba et~al.(2015)Belobaba, Odoni, and
  Barnhart}]{belobaba2015global_VLVV}
Belobaba, P., Odoni, A., Barnhart, C., 2015. The global airline industry, 2nd
  Edition. John Wiley \& Sons.

\bibitem[{Belobaba(1987{\natexlab{a}})}]{belobabaAirTravelDemand1987_AKSJF}
Belobaba, P.~P., 1987{\natexlab{a}}. Air travel demand and airline seat
  inventory management. Thesis, Massachusetts Institute of Technology.

\bibitem[{Belobaba(1987{\natexlab{b}})}]{belobaba1987survey_VLVV}
Belobaba, P.~P., 1987{\natexlab{b}}. Survey paper—airline yield management an
  overview of seat inventory control. Transportation Science 21~(2), 63--73.

\bibitem[{Belotti et~al.(2009)Belotti, Berthold, Bonami, Cafieri, Margot,
  Megaw, Vigerske, and Wächter}]{Couenne_CTCGE}
Belotti, P., Berthold, T., Bonami, P., Cafieri, S., Margot, F., Megaw, C.,
  Vigerske, S., Wächter, A., 2009. {Couenne}.
\newline\urlprefix\url{https://projects.coin-or.org/Couenne}

\bibitem[{Belton and Stewart(2002)}]{belton2002multiple_MESG}
Belton, V., Stewart, T., 2002. Multiple Criteria Decision Analysis: An
  Integrated Approach. Springer-Verlag, Berlin.

\bibitem[{Ben-Ameur et~al.(2007)Ben-Ameur, Breton, Karoui, and
  L’Ecuyer}]{ben2007_MB}
Ben-Ameur, H., Breton, M., Karoui, L., L’Ecuyer, P., 2007. A dynamic
  programming approach for pricing options embedded in bonds. Journal of
  Economic Dynamics and Control 31~(7), 2212--2233.

\bibitem[{Ben-Ameur et~al.(2002)Ben-Ameur, Breton, and L'Ecuyer}]{ben02_MB}
Ben-Ameur, H., Breton, M., L'Ecuyer, P., 2002. A dynamic programming procedure
  for pricing {A}merican-style {A}sian options. Management Science 48~(5),
  625--643.

\bibitem[{Ben-Tal et~al.(2009)Ben-Tal, El~Ghaoui, and
  Nemirovski}]{Ben-Tal2009-lm_HL}
Ben-Tal, A., El~Ghaoui, L., Nemirovski, A., 2009. Robust Optimization.
  Princeton University Press.

\bibitem[{Ben-Tal et~al.(2005)Ben-Tal, Golany, Nemirovski, and
  Vial}]{Ben-Tal2005-bz_HL}
Ben-Tal, A., Golany, B., Nemirovski, A., Vial, J.-P., 2005. {Retailer-Supplier}
  flexible commitments contracts: A robust optimization approach. Manufacturing
  \& Service Operations Management 7~(3), 248--271.

\bibitem[{Ben-Tal and Nemirovski(2002)}]{Ben-Tal2002-as_HL}
Ben-Tal, A., Nemirovski, A., 2002. Robust optimization -- methodology and
  applications. Mathematical Programming 92~(3), 453--480.

\bibitem[{Benders(1962)}]{Benders1962-sf_HL}
Benders, J.~F., 1962. Partitioning procedures for solving mixed-variables
  programming problems. Numerische Mathematik 4~(1), 238--252.

\bibitem[{Benders(2005)}]{benders2005partitioning_MH}
Benders, J.~F., 2005. Partitioning procedures for solving mixed-variables
  programming problems. Computational Management Science 2~(1), 3--19.

\bibitem[{Bendoly and Clark(2016)}]{Bendoly2016-cx_MJE}
Bendoly, E., Clark, S., 2016. Visual Analytics for Management: Translational
  Science and Applications in Practice. Routledge.

\bibitem[{Bengio et~al.(2021)Bengio, Lodi, and
  Prouvost}]{bengio2021machine_JLYHK}
Bengio, Y., Lodi, A., Prouvost, A., 2021. Machine learning for combinatorial
  optimization: a methodological tour d’horizon. European Journal of
  Operational Research 290~(2), 405--421.

\bibitem[{Bennell et~al.(2018)Bennell, Cabo, and
  Mart{\'\i}nez-Sykora}]{Bennell2018-fm_JB}
Bennell, J.~A., Cabo, M., Mart{\'\i}nez-Sykora, A., 2018. A beam search
  approach to solve the convex irregular bin packing problem with guillotine
  cuts. European Journal of Operational Research 270~(1), 89--102.

\bibitem[{Bennell and Oliveira(2008)}]{Bennell2008-nv_JB}
Bennell, J.~A., Oliveira, J.~F., 2008. The geometry of nesting problems: A
  tutorial. European Journal of Operational Research 184~(2), 397--415.

\bibitem[{Bennell and Oliveira(2009)}]{Bennell2009-kx_JB}
Bennell, J.~A., Oliveira, J.~F., 2009. A tutorial in irregular shape packing
  problems. Journal of the Operational Research Society 60~(1), S93--S105.

\bibitem[{Bennis and Nanus(1985)}]{Bennis1986-ia_AJG}
Bennis, W.~G., Nanus, B., 1985. Leaders: Strategies for Taking Charge. Harper
  \& Row, New York, NY.

\bibitem[{Benson and Sa{\u{g}}lam(2013)}]{BS13_ALAL}
Benson, H., Sa{\u{g}}lam, {\"U}., 2013. Mixed-integer second-order cone
  programming: a survey. In: Topaloglu, H. (Ed.), Theory Driven by Influential
  Applications. INFORMS, Catonsville, MD, pp. 13--36.

\bibitem[{Bergmeir and Ben{\'\i}tez(2012)}]{Bergmeir2012use_FP}
Bergmeir, C., Ben{\'\i}tez, J.~M., 2012. On the use of cross-validation for
  time series predictor evaluation. Information Sciences 191, 192--213.

\bibitem[{Berkey and Wang(1987)}]{Berkey1987-lb_JB}
Berkey, J.~O., Wang, P.~Y., 1987. {Two-Dimensional} finite {Bin-Packing}
  algorithms. Journal of the Operational Research Society 38~(5), 423--429.

\bibitem[{Bernhard(2005)}]{bernhard05_MB}
Bernhard, P., 2005. The robust control approach to option pricing and interval
  models: An overview. In: Breton, M., Ben-Ameur, H. (Eds.), Numerical Methods
  in Finance. Springer, New York, NY, pp. 91--108.

\bibitem[{Bertsekas(2012{\natexlab{a}})}]{bertsekas2012dynamic_DLLD}
Bertsekas, D., 2012{\natexlab{a}}. Dynamic programming and optimal control
  (Volume I). Athena Scientific, Belmont, MA.

\bibitem[{Bertsekas(2012{\natexlab{b}})}]{bertsekas2012dynamic2_DLLD}
Bertsekas, D., 2012{\natexlab{b}}. Dynamic programming and optimal control
  (Volume II: Approximate dynamic programming). Athena Scientific, Belmont, MA.

\bibitem[{Bertsekas and Tsitsiklis(1996)}]{Bertsekas1996-qb_HL}
Bertsekas, D., Tsitsiklis, J.~N., 1996. {Neuro-Dynamic} Programming. Athena
  Scientific.

\bibitem[{Bertsekas(2016)}]{Bert16_EAY}
Bertsekas, D.~P., 2016. Nonlinear programming, 3rd Edition. Athena Scientific,
  Belmont, MA.

\bibitem[{Bertsekas et~al.(1997)Bertsekas, Tsitsiklis, and
  Wu}]{Bertsekas1997-qm_HL}
Bertsekas, D.~P., Tsitsiklis, J.~N., Wu, C., 1997. Rollout algorithms for
  combinatorial optimization. Journal of Heuristics 3~(3), 245--262.

\bibitem[{Bertsimas et~al.(2011{\natexlab{a}})Bertsimas, Brown, and
  Caramanis}]{Bertsimas2011-fv_HL}
Bertsimas, D., Brown, D.~B., Caramanis, C., 2011{\natexlab{a}}. Theory and
  applications of robust optimization. SIAM Review 53~(3), 464--501.

\bibitem[{Bertsimas and De~Boer(2005)}]{bertsimas2005simulation_VLVV}
Bertsimas, D., De~Boer, S., 2005. Simulation-based booking limits for airline
  revenue management. Operations Research 53~(1), 90--106.

\bibitem[{Bertsimas and Lo(1998)}]{bertsimas98_MB}
Bertsimas, D., Lo, A.~W., 1998. Optimal control of execution costs. Journal of
  Financial Markets 1~(1), 1--50.

\bibitem[{Bertsimas et~al.(2011{\natexlab{b}})Bertsimas, Lulli, and
  Odoni}]{bertsimas2011integer_VLVV}
Bertsimas, D., Lulli, G., Odoni, A., 2011{\natexlab{b}}. An integer
  optimization approach to large-scale air traffic flow management. Operations
  Research 59~(1), 211--227.

\bibitem[{Bertsimas and Mourtzinou(1997)}]{bertsimas2_HATI}
Bertsimas, D., Mourtzinou, G., 1997. Transient laws of non-stationary queueing
  systems and their applications. Queueing Systems 25~(1), 115--155.

\bibitem[{Bertsimas and Nakazato(1995)}]{bertsimas_HATI}
Bertsimas, D., Nakazato, D., 1995. The distributional {L}ittle's law and its
  applications. Operations Research 43~(2), 298--310.

\bibitem[{Bertsimas and Patterson(1998)}]{bertsimas1998air_VLVV}
Bertsimas, D., Patterson, S.~S., 1998. The air traffic flow management problem
  with enroute capacities. Operations Research 46~(3), 406--422.

\bibitem[{Berwick(2005)}]{Berwick2005-hf_CV}
Berwick, D.~M., 2005. The {John Eisenberg} lecture: health services research as
  a citizen in improvement. Health Services Research 40~(2), 317--336.

\bibitem[{Besiou et~al.(2018)Besiou, Pedraza-Martinez, and
  Van~Wassenhove}]{besiou2018or_BYKOK}
Besiou, M., Pedraza-Martinez, A.~J., Van~Wassenhove, L.~N., 2018. {OR} applied
  to humanitarian operations. European Journal of Operational Research 269~(2),
  397--405.

\bibitem[{Besiou and Van~Wassenhove(2021)}]{besiou2021system_BYKOK}
Besiou, M., Van~Wassenhove, L.~N., 2021. System dynamics for humanitarian
  operations revisited. Journal of Humanitarian Logistics and Supply Chain
  Management 11~(4), 599--608.

\bibitem[{Bestuzheva et~al.(2021)Bestuzheva, Besan{\c{c}}on, Chen, Chmiela,
  Donkiewicz, van Doornmalen, Eifler, Gaul, Gamrath, Gleixner, Gottwald,
  Graczyk, Halbig, Hoen, Hojny, van~der Hulst, Koch, L{\"u}bbecke, Maher,
  Matter, M{\"u}hmer, M{\"u}ller, Pfetsch, Rehfeldt, Schlein, Schl{\"o}sser,
  Serrano, Shinano, Sofranac, Turner, Vigerske, Wegscheider, Wellner, Weninger,
  and Witzig}]{Scipopt_CTCGE}
Bestuzheva, K., Besan{\c{c}}on, M., Chen, W.-K., Chmiela, A., Donkiewicz, T.,
  van Doornmalen, J., Eifler, L., Gaul, O., Gamrath, G., Gleixner, A.,
  Gottwald, L., Graczyk, C., Halbig, K., Hoen, A., Hojny, C., van~der Hulst,
  R., Koch, T., L{\"u}bbecke, M., Maher, S.~J., Matter, F., M{\"u}hmer, E.,
  M{\"u}ller, B., Pfetsch, M.~E., Rehfeldt, D., Schlein, S., Schl{\"o}sser, F.,
  Serrano, F., Shinano, Y., Sofranac, B., Turner, M., Vigerske, S.,
  Wegscheider, F., Wellner, P., Weninger, D., Witzig, J., 2021. {The SCIP
  Optimization Suite 8.0}.
\newline\urlprefix\url{https://www.scipopt.org/}

\bibitem[{Better and Glover(2011)}]{Better2011-mq_TC}
Better, M., Glover, F., 2011. Simulation optimization in risk management. In:
  Cochran, J.~J., Cox~Jr., L.~A., Keskinocak, P., Kharoufeh, J.~P., Smith,
  J.~C. (Eds.), Wiley Encyclopedia of Operations Research and Management
  Science. John Wiley \& Sons, Hoboken, NJ, pp. 1--9.

\bibitem[{Beullens et~al.(2004)Beullens, Van~Oudheusden, and
  Van~Wassenhove}]{Beullens2004-nf_AM}
Beullens, P., Van~Oudheusden, D., Van~Wassenhove, L.~N., 2004. Collection and
  vehicle routing issues in reverse logistics. In: Dekker, R., Fleischmann, M.,
  Inderfurth, K., Van~Wassenhove, L.~N. (Eds.), Reverse Logistics: Quantitative
  Models for {Closed-Loop} Supply Chains. Springer, Berlin, pp. 95--134.

\bibitem[{Bezanson et~al.(2017)Bezanson, Edelman, Karpinski, and
  Shah.}]{Julia_CTCGE}
Bezanson, J., Edelman, A., Karpinski, S., Shah., V.~B., 2017. Julia: {A} fresh
  approach to numerical computing. SIAM Review 59, 65--98.

\bibitem[{Bhatia et~al.(2015)Bhatia, Hao, Kodialam, and
  Lakshman}]{bhatia.hao.ea:15_BF}
Bhatia, R., Hao, F., Kodialam, M., Lakshman, T., 2015. Optimized network
  traffic engineering using segment routing. 2015 IEEE Conference on Computer
  Communications (INFOCOM), 657--665.

\bibitem[{Bi et~al.(2021)Bi, Zhang, Fan, Tang, Guan, Gao, Cui, Ma, Wu, Hao,
  Ning, and Liu}]{Bi2021-gl_TC}
Bi, X., Zhang, Q., Fan, K., Tang, S., Guan, H., Gao, X., Cui, Y., Ma, Y., Wu,
  Q., Hao, Y., Ning, N., Liu, C., 2021. Risk culture and {COVID-19} protective
  behaviors: A {Cross-Sectional} survey of residents in china. Frontiers in
  Public Health 9, 686705.

\bibitem[{Bielecki et~al.(2017)Bielecki, Cialenco, and
  Pitera}]{Bielecki2017-kr_TC}
Bielecki, T.~R., Cialenco, I., Pitera, M., 2017. A survey of time consistency
  of dynamic risk measures and dynamic performance measures in discrete time:
  {LM-measure} perspective. Probability Uncertainty and Quantitative Risk
  2~(3).

\bibitem[{Bijvank et~al.(2023)Bijvank, Huh, and
  Janakiraman}]{bijvank2023lost_JSS}
Bijvank, M., Huh, W.~T., Janakiraman, G., 2023. Lost-sales inventory systems.
  In: Song, J.-S. (Ed.), Research Handbook on Inventory Management. Edward
  Elgar Publishing.

\bibitem[{Billionnet(2013)}]{Billionnet2013-hw_JL}
Billionnet, A., 2013. Mathematical optimization ideas for biodiversity
  conservation. European Journal of Operational Research 231~(3), 514--534.

\bibitem[{Bimpikis and Markakis(2016)}]{bimpikis2016inventory_JSS}
Bimpikis, K., Markakis, M.~G., 2016. Inventory pooling under heavy-tailed
  demand. Management Science 62~(6), 1800--1813.

\bibitem[{Binmore et~al.(1986)Binmore, Rubinstein, and
  Wolinsky}]{BinRubWol86_JH}
Binmore, K., Rubinstein, A., Wolinsky, A., 1986. The {Nash} bargaining solution
  in economic modeling. {RAND} Journal of Economics 17, 176--188.

\bibitem[{Birge(1982)}]{Birge1982-ix_HL}
Birge, J.~R., 1982. The value of the stochastic solution in stochastic linear
  programs with fixed recourse. Mathematical Programming 24~(1), 314--325.

\bibitem[{Birge(1985)}]{Birge1985-hx_HL}
Birge, J.~R., 1985. Decomposition and partitioning methods for multistage
  stochastic linear programs. Operations Research 33~(5), 989--1007.

\bibitem[{Birge and Louveaux(2011)}]{Birge2011-ut_HL}
Birge, J.~R., Louveaux, F., 2011. Introduction to Stochastic Programming.
  Springer.

\bibitem[{Bischoff et~al.(2018)Bischoff, Kaddoura, Maciejewski, and
  Nagel}]{Bischoff2018-fr_HL}
Bischoff, J., Kaddoura, I., Maciejewski, M., Nagel, K., 2018. Simulation-based
  optimization of service areas for pooled ride-hailing operators. Procedia
  Computer Science 130, 816--823.

\bibitem[{Bixby et~al.(1992)Bixby, Boyd, and Indovina}]{BBI92_ALAL}
Bixby, R., Boyd, E., Indovina, R., 1992. {MIPLIB}: a test set of mixed integer
  programming problems. SIAM News 25, 16.

\bibitem[{Black and Andersen(2012)}]{Black2012-ym_MJE}
Black, L.~J., Andersen, D.~F., 2012. Using visual representations as boundary
  objects to resolve conflict in collaborative model-building approaches.
  Systems Research and Behavioral Science 29~(2), 194--208.

\bibitem[{Blackwell et~al.(2008)Blackwell, Phaal, Eppler, and
  Crilly}]{Blackwell2008-gs_MJE}
Blackwell, A.~F., Phaal, R., Eppler, M., Crilly, N., 2008. Strategy roadmaps:
  New forms, new practices. In: Stapleton, G., Howse, J., Lee, J. (Eds.),
  Diagrammatic Representation and Inference. Springer, Berlin, Heidelberg, pp.
  127--140.

\bibitem[{Blanchet et~al.(2022)Blanchet, Chen, and Zhou}]{blanchet22_MB}
Blanchet, J., Chen, L., Zhou, X.~Y., 2022. Distributionally robust
  mean-variance portfolio selection with wasserstein distances. Management
  Science 68~(9), 6382--6410.

\bibitem[{Bland(1977)}]{Bland1977-dy_JMB}
Bland, R.~G., 1977. New finite pivoting rules for the simplex method.
  Mathematics of Operations Research 2~(2), 103--107.

\bibitem[{Bland et~al.(1981)Bland, Goldfarb, and Todd}]{Bland1981-fy_JMB}
Bland, R.~G., Goldfarb, D., Todd, M.~J., 1981. Feature {Article---The}
  ellipsoid method: A survey. Operations Research 29~(6), 1039--1091.

\bibitem[{B{\l}a{\.z}ewicz et~al.(2001)B{\l}a{\.z}ewicz, Ecker, Pesch, Schmidt,
  and Weglarz}]{BEPSW01_SMPT}
B{\l}a{\.z}ewicz, J., Ecker, K., Pesch, E., Schmidt, G., Weglarz, J., 2001.
  Scheduling Computer and Manufactoring Processes. Springer, Berlin.

\bibitem[{B{\l}a{\.z}ewicz et~al.(2007)B{\l}a{\.z}ewicz, Ecker, Pesch, Schmidt,
  and Weglarz}]{BEPSW07_SMPT}
B{\l}a{\.z}ewicz, J., Ecker, K., Pesch, E., Schmidt, G., Weglarz, J., 2007.
  Handbook on Scheduling: {F}rom Theory to Applications. Springer, Berlin.

\bibitem[{Bley et~al.(2009)Bley, Fortz, Gourdin, Holmberg, Klopfenstein,
  Pi{\'o}ro, Tomaszewski, and {\"U}mit}]{bley.fortz.ea:09_BF}
Bley, A., Fortz, B., Gourdin, E., Holmberg, K., Klopfenstein, O., Pi{\'o}ro,
  M., Tomaszewski, A., {\"U}mit, H., 2009. Optimization of {OSPF} routing in
  {IP} networks. In: Koster, A., Muñoz, X. (Eds.), Graphs and Algorithms in
  Communication Networks. Springer, pp. 199--240.

\bibitem[{Bloomfield and Cox(1972)}]{bloomfield1972low_HATI}
Bloomfield, P., Cox, D., 1972. A low traffic approximation for queues. Journal
  of Applied Probability 9~(4), 832--840.

\bibitem[{Blum and Roli(2003)}]{blum2003metaheuristics_COIT}
Blum, C., Roli, A., 2003. Metaheuristics in combinatorial optimization:
  {O}verview and conceptual comparison. ACM Computing Surveys (CSUR) 35~(3),
  268--308.

\bibitem[{Blum et~al.(1998)Blum, Cucker, Shub, and Smale}]{BCSS98_UPCT}
Blum, L., Cucker, F., Shub, M., Smale, S., 1998. Complexity and Real
  Computation. Springer, New York, NY.

\bibitem[{Blumenfeld et~al.(1985)Blumenfeld, Burns, Diltz, and
  Daganzo}]{blumenfeld1985analyzing_JLYHK}
Blumenfeld, D.~E., Burns, L.~D., Diltz, J.~D., Daganzo, C.~F., 1985. Analyzing
  trade-offs between transportation, inventory and production costs on freight
  networks. Transportation Research Part B: Methodological 19~(5), 361--380.

\bibitem[{Board et~al.(2003)Board, Sutcliffe, and Ziemba}]{board03_MB}
Board, J., Sutcliffe, C., Ziemba, W.~T., 2003. Applying operations research
  techniques to financial markets. Interfaces 33~(2), 12--24.

\bibitem[{Boginski et~al.(2015)Boginski, Pasiliao, and
  Shen}]{Boginski2015-iz_KVRH}
Boginski, V., Pasiliao, E.~L., Shen, S., 2015. Special issue on optimization in
  military applications. Optimization Letters 9~(8), 1475--1476.

\bibitem[{Boland et~al.(2017)Boland, Hewitt, Marshall, and
  Savelsbergh}]{DDD2017_MH}
Boland, N., Hewitt, M., Marshall, L., Savelsbergh, M., 2017. The
  continuous-time service network design problem. Operations Research 65~(5),
  1303--1321.

\bibitem[{Boland et~al.(2019)Boland, Hewitt, Marshall, and
  Savelsbergh}]{BOLAND2019195_MH}
Boland, N., Hewitt, M., Marshall, L., Savelsbergh, M., 2019. The price of
  discretizing time: a study in service network design. EURO Journal on
  Transportation and Logistics 8~(2), 195--216.

\bibitem[{Bommasani et~al.(2021)Bommasani, Hudson, Adeli, Altman, Arora, von
  Arx, Bernstein, Bohg, Bosselut, Brunskill, Brynjolfsson, Buch, Card,
  Castellon, Chatterji, Chen, Creel, Davis, Demszky, Donahue, Doumbouya,
  Durmus, Ermon, Etchemendy, Ethayarajh, Fei-Fei, Finn, Gale, Gillespie, Goel,
  Goodman, Grossman, Guha, Hashimoto, Henderson, Hewitt, Ho, Hong, Hsu, Huang,
  Icard, Jain, Jurafsky, Kalluri, Karamcheti, Keeling, Khani, Khattab, Koh,
  Krass, Krishna, Kuditipudi, Kumar, Ladhak, Lee, Lee, Leskovec, Levent, Li,
  Li, Ma, Malik, Manning, Mirchandani, Mitchell, Munyikwa, Nair, Narayan,
  Narayanan, Newman, Nie, Niebles, Nilforoshan, Nyarko, Ogut, Orr,
  Papadimitriou, Park, Piech, Portelance, Potts, Raghunathan, Reich, Ren, Rong,
  Roohani, Ruiz, Ryan, R'e, Sadigh, Sagawa, Santhanam, Shih, Srinivasan,
  Tamkin, Taori, Thomas, Tram{\`e}r, Wang, Wang, Wu, Wu, Wu, Xie, Yasunaga,
  You, Zaharia, Zhang, Zhang, Zhang, Zhang, Zheng, Zhou, and
  Liang}]{Bommasani2021FoundationModels}
Bommasani, R., Hudson, D.~A., Adeli, E., Altman, R., Arora, S., von Arx, S.,
  Bernstein, M.~S., Bohg, J., Bosselut, A., Brunskill, E., Brynjolfsson, E.,
  Buch, S., Card, D., Castellon, R., Chatterji, N.~S., Chen, A.~S., Creel,
  K.~A., Davis, J., Demszky, D., Donahue, C., Doumbouya, M., Durmus, E., Ermon,
  S., Etchemendy, J., Ethayarajh, K., Fei-Fei, L., Finn, C., Gale, T.,
  Gillespie, L.~E., Goel, K., Goodman, N.~D., Grossman, S., Guha, N.,
  Hashimoto, T., Henderson, P., Hewitt, J., Ho, D.~E., Hong, J., Hsu, K.,
  Huang, J., Icard, T.~F., Jain, S., Jurafsky, D., Kalluri, P., Karamcheti, S.,
  Keeling, G., Khani, F., Khattab, O., Koh, P.~W., Krass, M.~S., Krishna, R.,
  Kuditipudi, R., Kumar, A., Ladhak, F., Lee, M., Lee, T., Leskovec, J.,
  Levent, I., Li, X.~L., Li, X., Ma, T., Malik, A., Manning, C.~D.,
  Mirchandani, S.~P., Mitchell, E., Munyikwa, Z., Nair, S., Narayan, A.,
  Narayanan, D., Newman, B., Nie, A., Niebles, J.~C., Nilforoshan, H., Nyarko,
  J.~F., Ogut, G., Orr, L., Papadimitriou, I., Park, J.~S., Piech, C.,
  Portelance, E., Potts, C., Raghunathan, A., Reich, R., Ren, H., Rong, F.,
  Roohani, Y.~H., Ruiz, C., Ryan, J., R'e, C., Sadigh, D., Sagawa, S.,
  Santhanam, K., Shih, A., Srinivasan, K.~P., Tamkin, A., Taori, R., Thomas,
  A.~W., Tram{\`e}r, F., Wang, R.~E., Wang, W., Wu, B., Wu, J., Wu, Y., Xie,
  S.~M., Yasunaga, M., You, J., Zaharia, M.~A., Zhang, M., Zhang, T., Zhang,
  X., Zhang, Y., Zheng, L., Zhou, K., Liang, P., 2021. On the opportunities and
  risks of foundation models. arXiv:2108.07258.

\bibitem[{Bonami et~al.(2005)Bonami, Biegler, Conn, Cornu{\'e}jols, Grossmann,
  Laird, Lee, Lodi, Margot, Sawaya, and Waechter}]{Bonmin_CTCGE}
Bonami, P., Biegler, L., Conn, A., Cornu{\'e}jols, G., Grossmann, I., Laird,
  C., Lee, J., Lodi, A., Margot, F., Sawaya, N., Waechter, A., 2005. {An
  Algorithmic Framework for Convex Mixed Integer Nonlinear Programs}.
\newline\urlprefix\url{https://projects.coin-or.org/Bonmin}

\bibitem[{Bordignon et~al.(2013)Bordignon, Bunn, Lisi, and
  Nan}]{bor:bun:lis:nan:13_DSRW}
Bordignon, S., Bunn, D.~W., Lisi, F., Nan, F., 2013. Combining day-ahead
  forecasts for {British} electricity prices. Energy Economics 35, 88--103.

\bibitem[{Bor\r{u}vka(1926)}]{B26_SMPT}
Bor\r{u}vka, O., 1926. O jist\'em probl\'emu minim\'aln\'im. Pr\'ace Moravsk'e
  P\v{r}\'{\i}rodov\v{e}\-deck\'{e} Spole\v{c}nosti, 37--58.

\bibitem[{Bortfeldt and W{\"a}scher(2013)}]{Bortfeldt2013-wd_JB}
Bortfeldt, A., W{\"a}scher, G., 2013. Constraints in container loading -- a
  state-of-the-art review. European Journal of Operational Research 229~(1),
  1--20.

\bibitem[{Boschetti and Maniezzo(2022)}]{boschetti2022matheuristics_COIT}
Boschetti, M.~A., Maniezzo, V., 2022. Matheuristics: using mathematics for
  heuristic design. 4OR - A Quarterly Journal of Operations Research 20~(2),
  173--208.

\bibitem[{Bouleimen and Lecocq(2003)}]{Bouleimen2003-ud_WH_ED}
Bouleimen, K., Lecocq, H., 2003. A new efficient simulated annealing algorithm
  for the resource-constrained project scheduling problem and its multiple mode
  version. European Journal of Operational Research 149~(2), 268--281.

\bibitem[{Boute et~al.(2022)Boute, Disney, Gijsbrechts, and
  Mieghem}]{Boute2022_SMD}
Boute, R., Disney, S.~M., Gijsbrechts, J., Mieghem, J. A.~V., 2022. Dual
  sourcing and smoothing under non-stationary demand time series: {R}e-shoring
  with {S}peed{F}actories. Management Science 68, 1039--1057.

\bibitem[{Bowen(1995)}]{Bowen1995-to_AJG}
Bowen, K. (Ed.), 1995. In at the Deep End: {MSc} Student Projects in Community
  Operational Research. Community Operational Research Unit Publications.

\bibitem[{{Box} et~al.(1976){Box}, George, {Jenkins}, and
  {Gwilym}}]{Box1976-af_FP}
{Box}, George, E.~P., {Jenkins}, {Gwilym}, 1976. Time Series Analysis
  Forecasting and Control. {Holden-Day}, San Francisco, {CA}.

\bibitem[{Boxma and Cohen(2000)}]{boxma2000single_HATI}
Boxma, O., Cohen, J., 2000. The single server queue: Heavy tails and heavy
  traffic. In: Park, K., Willinger, W. (Eds.), Self-Similar Network Traffic and
  Performance Evaluation. Wiley Online Library, pp. 143--169.

\bibitem[{Boyd et~al.(2007)Boyd, Geerling, Gregory, Kagan, Midgley, Murray, and
  Walsh}]{Boyd2007-am_GM}
Boyd, A., Geerling, T., Gregory, W.~J., Kagan, C., Midgley, G., Murray, P.,
  Walsh, M.~P., 2007. Systemic evaluation: a participative, multi-method
  approach. Journal of the Operational Research Society 58~(10), 1306--1320.

\bibitem[{Boykov and Kolmogorov(2004)}]{Boykov-Kolmogorov:2004_IL}
Boykov, Y., Kolmogorov, V., 2004. An experimental comparison of
  min-cut/max-flow algorithms for energy minimization in vision. {IEEE}
  Transactions on Pattern Analysis and Machine Intelligence 26~(9), 1124--1137.

\bibitem[{Boysen et~al.(2012)Boysen, Fliedner, Jaehn, and
  Pesch}]{Boysen2012_DC}
Boysen, N., Fliedner, M., Jaehn, F., Pesch, E., 2012. Shunting yard
  operations{: T}heoretical aspects and applications. European Journal of
  Operational Research 220, 1--14.

\bibitem[{Brailsford and Schmidt(2003)}]{Brailsford2003-rs_AFRH}
Brailsford, S., Schmidt, B., 2003. Towards incorporating human behaviour in
  models of health care systems: An approach using discrete event simulation.
  European Journal of Operational Research 150~(1), 19--31.

\bibitem[{Brailsford et~al.(2004)Brailsford, Lattimer, Tarnaras, and
  Turnbull}]{Brailsford2004-ib_CV}
Brailsford, S.~C., Lattimer, V.~A., Tarnaras, P., Turnbull, J.~C., 2004.
  Emergency and on-demand health care: modelling a large complex system.
  Journal of the Operational Research Society 55~(1), 34--42.

\bibitem[{Bramson(1998)}]{bramson1998state_HATI}
Bramson, M., 1998. State space collapse with application to heavy traffic
  limits for multiclass queueing networks. Queueing Systems 30~(1), 89--140.

\bibitem[{Branke et~al.(2016)Branke, Corrente, Greco, S{\l}owi{\'n}ski, and
  Zielniewicz}]{branke2016using_MESG}
Branke, J., Corrente, S., Greco, S., S{\l}owi{\'n}ski, R., Zielniewicz, P.,
  2016. Using {C}hoquet integral as preference model in interactive
  evolutionary multiobjective optimization. European Journal of Operational
  Research 250~(3), 884--901.

\bibitem[{Branke et~al.(2008)Branke, Deb, Miettinen, and
  Slowi{\'n}ski}]{branke2008multiobjective_MESG}
Branke, J., Deb, K., Miettinen, K., Slowi{\'n}ski, R. (Eds.), 2008.
  Multiobjective Optimization: Interactive and Evolutionary Approaches. Vol.
  5252 of Lecture Notes in Computer Science. Springer-Verlag, Berlin.

\bibitem[{Brans and Gallo(2007)}]{Brans2007-xi_JH}
Brans, J.-P., Gallo, G., 2007. Ethics in {OR/MS}: past, present and future.
  Annals of Operations Research 153~(1), 165--178.

\bibitem[{Brax and Visintin(2017)}]{brax_meta-model_2017_JHLS}
Brax, S.~A., Visintin, F., 2017. Meta-model of servitization: {The} integrative
  profiling approach. Industrial Marketing Management 60, 17--32.

\bibitem[{Breiman(2001)}]{Breiman01_LCAL}
Breiman, L., 2001. Random forests. Machine Learning 45~(1), 5--32.

\bibitem[{Bresciani and Eppler(2015)}]{Bresciani2015-ac_MJE}
Bresciani, S., Eppler, M.~J., 2015. The pitfalls of visual representations: A
  review and classification of common errors made while designing and
  interpreting visualizations. SAGE Open 5~(4), 2158244015611451.

\bibitem[{Bresciani and Eppler(2018)}]{Bresciani2018-xj_MJE}
Bresciani, S., Eppler, M.~J., 2018. The collaborative dimensions of argument
  maps: A socio-visual approach. Semiotica 2018~(220), 199--216.

\bibitem[{Breton(2018)}]{breton18_2_MB}
Breton, M., 2018. Dynamic games in finance. In: Duan, J.-C., H\"{a}rdle, W.~K.,
  Gentle, J.~E. (Eds.), Handbook of Dynamic Game Theory. Springer, pp.
  827--863.

\bibitem[{Breton and Marzouk(2018)}]{breton18_1_MB}
Breton, M., Marzouk, O., 2018. Evaluation of counterparty risk for derivatives
  with early-exercise features. Journal of Economic Dynamics and Control 88,
  1--20.

\bibitem[{Brettel et~al.(2014)Brettel, Friederichsen, Keller, and
  Rosenberg}]{Brettel2014-ho_KESXZ}
Brettel, M., Friederichsen, N., Keller, M., Rosenberg, M., 2014. How
  virtualization, decentralization and network building change the
  manufacturing landscape: An industry 4.0 perspective. International Journal
  of Mechanical, Industrial Science and Engineering 8~(1), 37--44.

\bibitem[{Brezav{\v s}{\v c}ek et~al.(2017)Brezav{\v s}{\v c}ek, Bach, and
  Baggia}]{Brezavscek2017-cj_JJ}
Brezav{\v s}{\v c}ek, A., Bach, M.~P., Baggia, A., 2017. Markov analysis of
  students' performance and academic progress in higher education. Organizacija
  50~(2), 83--95.

\bibitem[{Briskorn et~al.(2010)Briskorn, Drexl, and
  Spieksma}]{Briskorn2010_GVBSP}
Briskorn, D., Drexl, A., Spieksma, F. C.~R., 2010. Round robin tournaments and
  three index assignments. 4OR - A Quarterly Journal of Operations Research
  8~(4), 365--374.

\bibitem[{Broadie(2012)}]{Broadie2012-eo_IM}
Broadie, M., 2012. Assessing golfer performance on the {PGA} {TOUR}. Interfaces
  42~(2), 146--165.

\bibitem[{Broadie et~al.(2007)Broadie, Chernov, and Sundaresan}]{broadie07_MB}
Broadie, M., Chernov, M., Sundaresan, S., 2007. Optimal debt and equity values
  in the presence of {C}hapter 7 and {C}hapter 11. The Journal of Finance
  62~(3), 1341--1377.

\bibitem[{Brooks-Pollock et~al.(2021)Brooks-Pollock, Danon, Jombart, and
  Pellis}]{Brooks-Pollock2021-lw_CV}
Brooks-Pollock, E., Danon, L., Jombart, T., Pellis, L., 2021. Modelling that
  shaped the early {COVID-19} pandemic response in the {UK}. Philosophical
  Transactions of the Royal Society of London. Series B, Biological Sciences
  376~(1829), 20210001.

\bibitem[{Brotcorne et~al.(2003)Brotcorne, Laporte, and
  Semet}]{brotcorne2003ambulance_SAA}
Brotcorne, L., Laporte, G., Semet, F., 2003. Ambulance location and relocation
  models. European Journal of Operational Research 147~(3), 451--463.

\bibitem[{Brouer et~al.(2014)Brouer, Alvarez, Plum, Pisinger, and
  Sigurd}]{Brouer2014-ue_HP}
Brouer, B.~D., Alvarez, J.~F., Plum, C. E.~M., Pisinger, D., Sigurd, M.~M.,
  2014. A base integer programming model and benchmark suite for
  {Liner-Shipping} network design. Transportation Science 48~(2), 281--312.

\bibitem[{Brown et~al.(2005)Brown, Carlyle, Diehl, Kline, and
  Wood}]{Brown2005-ao_KVRH}
Brown, G., Carlyle, M., Diehl, D., Kline, J., Wood, K., 2005. A two-sided
  optimization for theater ballistic missile defense. Operations Research
  53~(5), 745--763.

\bibitem[{Brown et~al.(2004)Brown, Dell, and Newman}]{Brown2004-mg_KVRH}
Brown, G.~G., Dell, R.~F., Newman, A.~M., 2004. Optimizing military capital
  planning. Interfaces 34~(6), 415--425.

\bibitem[{Brown(1956)}]{Brown56_FP}
Brown, R.~G., 1956. Exponential smoothing for predicting demand. Little,
  Cambridge, MA.

\bibitem[{Browne et~al.(1984)Browne, Dubois, Sethi, Rathmill, and
  Stecke}]{Browne1984-ec_KESXZ}
Browne, J., Dubois, D., Sethi, S., Rathmill, K., Stecke, K., 1984.
  Classification of flexible manufacturing systems. The FMS Magazine 2~(2),
  114--117.

\bibitem[{Brucker et~al.(1998)Brucker, Knust, Schoo, and
  Thiele}]{Brucker1998-qg_WH_ED}
Brucker, P., Knust, S., Schoo, A., Thiele, O., 1998. A branch and bound
  algorithm for the resource-constrained project scheduling problem. European
  Journal of Operational Research 107~(2), 272--288.

\bibitem[{Brumelle and McGill(1993)}]{brumelle1993airline_VLVV}
Brumelle, S.~L., McGill, J.~I., 1993. Airline seat allocation with multiple
  nested fare classes. Operations Research 41~(1), 127--137.

\bibitem[{Brumelle et~al.(1990)Brumelle, McGill, Oum, Sawaki, and
  Tretheway}]{brumelle1990allocation_VLVV}
Brumelle, S.~L., McGill, J.~I., Oum, T.~H., Sawaki, K., Tretheway, M., 1990.
  Allocation of airline seats between stochastically dependent demands.
  Transportation Science 24~(3), 183--192.

\bibitem[{Bryant(1978)}]{Bryant1978_EA}
Bryant, J.~W., 1978. Modelling for natural resource utilization analysis.
  Journal of the Operational Research Society 29~(7), 667--676.

\bibitem[{Buhaug(2002)}]{Buhaug2002-ck_CV}
Buhaug, H., 2002. Long waiting lists in hospitals. BMJ 324~(7332), 252--253.

\bibitem[{Bulut et~al.(2019)Bulut, Xu, Ralphs, and Vigerske}]{Disco_CTCGE}
Bulut, A., Xu, Y., Ralphs, T., Vigerske, S., 2019. {Discrete Conic Optimization
  Solver}.
\newline\urlprefix\url{https://projects.coin-or.org/Disco}

\bibitem[{Buraimo et~al.(2020)Buraimo, Forrest, McHale, and
  Tena}]{Buraimo2020-ur_IM}
Buraimo, B., Forrest, D., McHale, I.~G., Tena, J.~D., 2020. Unscripted drama:
  Soccer audience response to suspense, surprise, and shock. Economic Inquiry
  58~(2), 881--896.

\bibitem[{Burdett and Kozan(2009)}]{Burdett2009_DC}
Burdett, E.~L., Kozan, E., 2009. Techniques for inserting additional trains
  into existing timetables. Transportation Research Part B: Methodological
  43~(8-9), 821--836.

\bibitem[{Burger et~al.(2019)Burger, White, and
  Yearworth}]{burger_developing_2019_JHLS}
Burger, K., White, L., Yearworth, M., 2019. Developing a smart operational
  research with hybrid practice theories. European Journal of Operational
  Research 277~(3), 1137--1150.

\bibitem[{Burk and Parnell(2011)}]{Burk2011-ns_KVRH}
Burk, R.~C., Parnell, G.~S., 2011. Portfolio decision analysis: Lessons from
  military applications. In: Salo, A., Keisler, J., Morton, A. (Eds.),
  Portfolio Decision Analysis: Improved Methods for Resource Allocation.
  Springer, New York, NY, pp. 333--357.

\bibitem[{Burkard et~al.(2012)Burkard, Dell'Amico, and Martello}]{BDM12_SMPT}
Burkard, R., Dell'Amico, M., Martello, S., 2012. Assignment problems, revised
  reprint. SIAM, Philadelphia.

\bibitem[{Burkard and Derigs(1980)}]{BD80_SMPT}
Burkard, R., Derigs, U., 1980. Assignment and Matching Problems: Solution
  Methods with {FORTRAN} Programs. Springer-Verlag, Berlin, Heidelberg.

\bibitem[{Burke and Bykov(2008)}]{LA_GVBSP}
Burke, E.~K., Bykov, Y., 2008. A late acceptance strategy in hill-climbing for
  exam timetabling problems. In: Proceedings of the International Conference of
  the Practice and Theory of Automated Timetabling (PATAT 2008). pp. 1--7.

\bibitem[{Burke and Curtois(2014)}]{curtois_instances_GVBSP}
Burke, E.~K., Curtois, T., 2014. New approaches to nurse rostering benchmark
  instances. European Journal of Operational Research 237, 71–81.

\bibitem[{Burke et~al.(2004{\natexlab{a}})Burke, {De Causmaecker}, {Vanden
  Berghe}, and {Van Landeghem}}]{JOS_GVBSP}
Burke, E.~K., {De Causmaecker}, P., {Vanden Berghe}, G., {Van Landeghem}, H.,
  2004{\natexlab{a}}. The state of the art of nurse rostering. Journal of
  Scheduling 7~(6), 441--499.

\bibitem[{Burke et~al.(2004{\natexlab{b}})Burke, Kendall, and
  Whitwell}]{Burke2004-ia_JB}
Burke, E.~K., Kendall, G., Whitwell, G., 2004{\natexlab{b}}. A new placement
  heuristic for the orthogonal {Stock-Cutting} problem. Operations Research
  52~(4), 655--671.

\bibitem[{Burke et~al.(2006)Burke, MacCarthy, Petrovic, and
  Qu}]{Burke2006-wo_JJ}
Burke, E.~K., MacCarthy, B.~L., Petrovic, S., Qu, R., 2006. Multiple-retrieval
  case-based reasoning for course timetabling problems. Journal of the
  Operational Research Society 57~(2), 148--162.

\bibitem[{Burke and Petrovic(2002)}]{Burke2002-dp_JJ}
Burke, E.~K., Petrovic, S., 2002. Recent research directions in automated
  timetabling. European Journal of Operational Research 140~(2), 266--280.

\bibitem[{Burman and Smith(1983)}]{burman1983light_HATI}
Burman, D.~Y., Smith, D.~R., 1983. A light-traffic theorem for multi-server
  queues. Mathematics of Operations Research 8~(1), 15--25.

\bibitem[{Burney et~al.(2013)Burney, Johnes, Al-Enezi, and
  Al-Musallam}]{Burney2013-iw_JJ}
Burney, N.~A., Johnes, J., Al-Enezi, M., Al-Musallam, M., 2013. The efficiency
  of public schools: the case of {K}uwait. Education Economics 21~(4),
  360--379.

\bibitem[{B{\"u}sing et~al.(2019)B{\"u}sing, Kadatz, and
  Cleophas}]{busing2019capacity_VLVV}
B{\"u}sing, C., Kadatz, D., Cleophas, C., 2019. Capacity uncertainty in airline
  revenue management: Models, algorithms, and computations. Transportation
  Science 53~(2), 383--400.

\bibitem[{Bussieck et~al.(2004)Bussieck, Lindner, and
  Lubbecke}]{Bussieck2004_DC}
Bussieck, M.~R., Lindner, T., Lubbecke, M.~E., 2004. A fast algorithm for near
  cost optimal line plans. Mathematical Methods of Operations Research 59,
  205--220.

\bibitem[{Bussieck et~al.(1997)Bussieck, Winter, and
  Zimmerman}]{Bussieck1997a_DC}
Bussieck, M.~R., Winter, T., Zimmerman, U.~T., 1997. Discrete optimization in
  public rail trainsport. Mathematical Programming 79, 415-- 444.

\bibitem[{Butcher et~al.(2021)Butcher, Huang, Robinson, Reffin, Sgaier,
  Charles, and Quadrianto}]{Butcher2021-zo_TC}
Butcher, B., Huang, V.~S., Robinson, C., Reffin, J., Sgaier, S.~K., Charles,
  G., Quadrianto, N., 2021. Causal datasheet for datasets: An evaluation guide
  for {Real-World} data analysis and data collection design using bayesian
  networks. Frontiers in Artificial Intelligence 4, 612551.

\bibitem[{B{\"u}y{\"u}ktahtak{\i}n and Haight(2018)}]{Buyuktahtakin2018-bb_JL}
B{\"u}y{\"u}ktahtak{\i}n, {\.I}.~E., Haight, R.~G., 2018. A review of
  operations research models in invasive species management: state of the art,
  challenges, and future directions. Annals of Operations Research 271~(2),
  357--403.

\bibitem[{Bykov and Petrovic(2016)}]{step_GVBSP}
Bykov, Y., Petrovic, S., 2016. A step counting hill climbing algorithm applied
  to university examination timetabling. Journal of Scheduling 19, 479--492.

\bibitem[{Bynum et~al.(2021)Bynum, Hackebeil, Hart, Laird, Nicholson, Siirola,
  Watson, and Woodruff}]{Pyomo_CTCGE}
Bynum, M.~L., Hackebeil, G.~A., Hart, W.~E., Laird, C.~D., Nicholson, B.~L.,
  Siirola, J.~D., Watson, J.-P., Woodruff, D.~L., 2021. Pyomo--optimization
  modeling in {P}ython, 3rd Edition. Vol.~67. Springer Science \& Business
  Media.
\newline\urlprefix\url{http://www.pyomo.org/}

\bibitem[{Bö$\eth$varsdóttir et~al.(2021)Bö$\eth$varsdóttir, Smet, {Vanden
  Berghe}, and Stidsen}]{weights_GVBSP}
Bö$\eth$varsdóttir, E.~B., Smet, P., {Vanden Berghe}, G., Stidsen, T.~J.,
  2021. Achieving compromise solutions in nurse rostering by using
  automatically estimated acceptance thresholds. European Journal of
  Operational Research 292~(3), 980--995.

\bibitem[{Büyüktahtakin et~al.(2014)Büyüktahtakin, Feng, and
  Szidarovszky}]{Buyuktahtakin2014_EA}
Büyüktahtakin, I.~E., Feng, Z., Szidarovszky, F., 2014. A multi-objective
  optimization approach for invasive species control. Journal of the
  Operational Research Society 65~(11), 1625--1635.

\bibitem[{Cabrera et~al.(2023)Cabrera, Cabrera, and
  Midgley}]{Cabrera2023-nw_GM}
Cabrera, D., Cabrera, L., Midgley, G., 2023. The four waves of systems
  thinking. In: Cabrera, D., Cabrera, L., Midgley, G. (Eds.), Routledge
  Handbook of Systems Thinking. Routledge, London.

\bibitem[{Cabrera et~al.(2015)Cabrera, Cabrera, and Powers}]{Cabrera2015-rd_GM}
Cabrera, D., Cabrera, L., Powers, E., 2015. A unifying theory of systems
  thinking with psychosocial applications. Systems Research and Behavioral
  Science 32~(5), 534--545.

\bibitem[{Cacchiani et~al.(2008{\natexlab{a}})Cacchiani, Caprara, Galli, Kroon,
  and Mar{\'o}ti}]{Cacchiani2008a_DC}
Cacchiani, V., Caprara, A., Galli, L., Kroon, L., Mar{\'o}ti, G.,
  2008{\natexlab{a}}. {Recoverable Robustness for Railway Rolling Stock
  Planning}. In: Fischetti, M., Widmayer, P. (Eds.), 8th Workshop on
  Algorithmic Approaches for Transportation Modeling, Optimization, and Systems
  (ATMOS'08). Vol.~9 of OpenAccess Series in Informatics (OASIcs). Schloss
  Dagstuhl--Leibniz-Zentrum fuer Informatik, Dagstuhl, Germany, pp. 9--10.

\bibitem[{Cacchiani et~al.(2012)Cacchiani, Caprara, Galli, Kroon, Mar{\'o}ti,
  and Toth}]{Cacchiani2012_DC}
Cacchiani, V., Caprara, A., Galli, L., Kroon, L., Mar{\'o}ti, G., Toth, P.,
  2012. Railway rolling stock planning: Robustness against large disruptions.
  Transportation Science 46~(2), 217--232.

\bibitem[{Cacchiani et~al.(2008{\natexlab{b}})Cacchiani, Caprara, and
  Toth}]{Cacchiani2008_DC}
Cacchiani, V., Caprara, A., Toth, P., 2008{\natexlab{b}}. A column generation
  approach to train timetabling on a corridor. 4OR - A Quarterly Journal of
  Operations Research 6~(2), 125--142.

\bibitem[{Cacchiani et~al.(2010)Cacchiani, Caprara, and
  Toth}]{Cacchiani2010_DC}
Cacchiani, V., Caprara, A., Toth, P., 2010. Scheduling extra freight trains on
  railway networks. Transportation Research Part B: Methodological 44~(2),
  215--231.

\bibitem[{Cacchiani et~al.(2022{\natexlab{a}})Cacchiani, Iori, Locatelli, and
  Martello}]{CILM22A_SMPT}
Cacchiani, V., Iori, M., Locatelli, A., Martello, S., 2022{\natexlab{a}}.
  Knapsack problems -- an overview of recent advances. {Part I}: Single
  knapsack problems. Computers \& Operations Research 143, 105692.

\bibitem[{Cacchiani et~al.(2022{\natexlab{b}})Cacchiani, Iori, Locatelli, and
  Martello}]{CILM22B_SMPT}
Cacchiani, V., Iori, M., Locatelli, A., Martello, S., 2022{\natexlab{b}}.
  Knapsack problems -- an overview of recent advances. {Part II}: Multiple,
  multidimensional, and quadratic knapsack problems. Computers \& Operations
  Research 143, 105693.

\bibitem[{Cachon and Lariviere(2005)}]{Cachon2005_SMD}
Cachon, G.~P., Lariviere, M.~A., 2005. Supply chain coordination with
  revenue-sharing contracts: Strengths and limitations. Management Science
  51~(1), 30--44.

\bibitem[{Cadarso and Vaze(2022)}]{cadarso2022passenger_VLVV}
Cadarso, L., Vaze, V., 2022. Passenger-centric integrated airline schedule and
  aircraft recovery. Transportation Science 56~(6), 1410--1431.

\bibitem[{Cadarso et~al.(2017)Cadarso, Vaze, Barnhart, and
  Mar{\'\i}n}]{cadarso2017integrated_VLVV}
Cadarso, L., Vaze, V., Barnhart, C., Mar{\'\i}n, {\'A}., 2017. Integrated
  airline scheduling: Considering competition effects and the entry of the high
  speed rail. Transportation Science 51~(1), 132--154.

\bibitem[{{\c{C}}al{\i}k et~al.(2019){\c{C}}al{\i}k, Labb{\'e}, and
  Yaman}]{ccalik2019p_SAA}
{\c{C}}al{\i}k, H., Labb{\'e}, M., Yaman, H., 2019. p-center problems. In:
  Laporte, G., Nickel, S., Saldanha~da Gama, F. (Eds.), Location Science.
  Springer, pp. 51--65.

\bibitem[{Callon(1981)}]{Callon1981-tn_MYLW}
Callon, M., 1981. Struggles and negotiations to define what is problematic and
  what is not. In: Knorr, K.~D., Krohn, R., Whitley, R. (Eds.), The Social
  Process of Scientific Investigation. Springer Netherlands, Dordrecht, pp.
  197--219.

\bibitem[{Calmon and Graves(2017)}]{Calmon2017-gn_AM}
Calmon, A.~P., Graves, S.~C., 2017. Inventory management in a consumer
  electronics {Closed-Loop} supply chain. Manufacturing \& Service Operations
  Management 19~(4), 568--585.

\bibitem[{Calmon et~al.(2021)Calmon, Graves, and Lemmens}]{Calmon2021-vf_AM}
Calmon, A.~P., Graves, S.~C., Lemmens, S., 2021. Warranty matching in a
  consumer electronics {Closed-Loop} supply chain. Manufacturing \& Service
  Operations Management 23~(5), 1314--1331.

\bibitem[{Camacho and Bordons(2013)}]{Camacho2013_XW}
Camacho, E.~F., Bordons, C., 2013. Model Predictive Control. Springer, London.

\bibitem[{Campbell and Savelsbergh(2006)}]{campbell2006incentive_CKTVW}
Campbell, A.~M., Savelsbergh, M., 2006. Incentive schemes for attended home
  delivery services. Transportation Science 40~(3), 327--341.

\bibitem[{Campbell and Savelsbergh(2005)}]{campbell2005decision_CKTVW}
Campbell, A.~M., Savelsbergh, M.~W., 2005. Decision support for consumer direct
  grocery initiatives. Transportation Science 39~(3), 313--327.

\bibitem[{Campbell and Sayer(2003)}]{Campbell2003_EA}
Campbell, B.~M., Sayer, J.~A., 2003. Integrated natural resources management:
  linking productivity, the environment and development. CABI Publishing,
  Cambridge, MA.

\bibitem[{Canca et~al.(2018)Canca, Andrade-Pineda, De-los Santos, and
  Calle}]{Canca2018_DC}
Canca, D., Andrade-Pineda, J.~L., De-los Santos, A., Calle, M., 2018. The
  railway rapid transit frequency setting problem with speed-dependent
  operation costs. Transportation Research Part B: Methodological 117~(A),
  494--519.

\bibitem[{Canca et~al.(2014)Canca, Barrena, Algaba, and Zarzo}]{Canca2014_DC}
Canca, D., Barrena, E., Algaba, E., Zarzo, A., 2014. Design and analysis of
  demand‐adapted railway timetables. Journal of Advanced Transportation
  48~(2), 119--137.

\bibitem[{Canca et~al.(2016)Canca, Barrena, De-Los-Santos, and
  Andrade-Pineda}]{Canca2016a_DC}
Canca, D., Barrena, E., De-Los-Santos, A., Andrade-Pineda, J.~L., 2016. Setting
  lines frequency and capacity in dense railway rapid transit networks with
  simultaneous passenger assignment. Transportation Research Part B:
  Methodological 93(A), 251--267.

\bibitem[{Canca et~al.(2017)Canca, De-Los-Santos, Laporte, and
  Mesa}]{Canca2017_DC}
Canca, D., De-Los-Santos, A., Laporte, G., Mesa, J.~A., 2017. An adaptive
  neighborhood search metaheuristic for the integrated railway rapid transit
  network design and line planning problem. Computers \& {O}perations
  {R}esearch 78, 1--14.

\bibitem[{Canca et~al.(2019)Canca, De-Los-Santos, Laporte, and
  Mesa}]{Canca2019b_DC}
Canca, D., De-Los-Santos, A., Laporte, G., Mesa, J.~A., 2019. Integrated
  railway rapid transit network design and line-planning problem with maximum
  profit. Transportation Research Part E: Logistics and Transportation Review
  127, 1--30.

\bibitem[{Canca and Zarzo(2017)}]{Canca2017a_DC}
Canca, D., Zarzo, A., 2017. Design of energy{-E}fficient timetables in
  two{-w}ay railway rapid transit lines. Transportation Research Part B:
  Methodological 102, 142--161.

\bibitem[{Cao et~al.(2022)Cao, Lan, Xie, Chen, Zhang, Zhang, and
  Wu}]{Cao2022-dj_IM}
Cao, A., Lan, J., Xie, X., Chen, H., Zhang, X., Zhang, H., Wu, Y., 2022.
  {Team-Builder}: Toward more effective lineup selection in soccer. IEEE
  Transactions on Visualization and Computer Graphics, DOI:
  10.1109/TVCG.2022.3207147.

\bibitem[{Cappart et~al.(2021)Cappart, Ch{\'e}telat, Khalil, Lodi, Morris, and
  Veli{\v c}kovi{\'c}}]{cappart2021combinatorial_LCAL}
Cappart, Q., Ch{\'e}telat, D., Khalil, E.~B., Lodi, A., Morris, C., Veli{\v
  c}kovi{\'c}, P., 8 2021. Combinatorial optimization and reasoning with graph
  neural networks. In: Zhou, Z.-H. (Ed.), Proceedings of the Thirtieth
  International Joint Conference on Artificial Intelligence, {IJCAI-21}.
  International Joint Conferences on Artificial Intelligence Organization, pp.
  4348--4355, survey Track.

\bibitem[{Caprara et~al.(2014)Caprara, Carvalho, Lodi, and
  Woeginger}]{Cap14_UPCT}
Caprara, A., Carvalho, M., Lodi, A., Woeginger, G.~J., 2014. A study on the
  computational complexity of the bilevel knapsack problem. {SIAM} Journal on
  Optimization 24~(2), 823--838.

\bibitem[{Caprara et~al.(2002)Caprara, Fischetti, and Toth}]{Caprara2002_DC}
Caprara, A., Fischetti, M., Toth, P., 2002. Modeling and solving the train
  timetabling problem. Operations Research 50~(5), 851--861.

\bibitem[{Caprara et~al.(2007)Caprara, Kroon, Monaci, Peeters, and
  Toth}]{Caprara2007_DC}
Caprara, A., Kroon, L., Monaci, M., Peeters, M., Toth, P., 2007. Passenger
  railway optimization. In: Barnhart, C., Laporte, G. (Eds.), Handbooks in
  Operations Research and Management Science. Vol. 14. Transportation.
  Elsevier, pp. 129--187.

\bibitem[{Caprara et~al.(2003)Caprara, Monaci, and Toth}]{Caprara2003_DC}
Caprara, A., Monaci, M., Toth, P., 2003. Models and algorithms for a staff
  scheduling problem. Mathematical Programming 98~(1), 445--476.

\bibitem[{Caprara et~al.(1998)Caprara, Toth, Vigo, and
  Fischetti}]{Caprara1998_DC}
Caprara, A., Toth, P., Vigo, D., Fischetti, M., 1998. Modeling and solving the
  crew rostering problem. Operations Research 46~(6), 820--830.

\bibitem[{Cardaliaguet and Lehalle(2018)}]{cardaliaguet18_MB}
Cardaliaguet, P., Lehalle, C.-A., 2018. Mean field game of controls and an
  application to trade crowding. Mathematics and Financial Economics 12~(3),
  335--363.

\bibitem[{Cardoen et~al.(2010)Cardoen, Demeulemeester, and
  Beli{\"e}n}]{Cardoen2010-rw_CV}
Cardoen, B., Demeulemeester, E., Beli{\"e}n, J., 2010. Operating room planning
  and scheduling: A literature review. European Journal of Operational Research
  201~(3), 921--932.

\bibitem[{Caro and Gallien(2012)}]{caroClearancePricingOptimization2012_AKSJF}
Caro, F., Gallien, J., 2012. Clearance pricing optimization for a fast-fashion
  retailer. Operations Research 60~(6), 1404--1422.

\bibitem[{Carraro and Haurie(1996)}]{Haurie2012_EA}
Carraro, C., Haurie, A., 1996. Operations Research and Environmental
  Management. Kluwer Academic, Dordrecht/Boston/London.

\bibitem[{Carter et~al.(1996)Carter, Laporte, and Lee}]{carter1996_GVBSP}
Carter, M.~W., Laporte, G., Lee, S.~Y., 1996. Examination timetabling:
  Algorithmic strategies and applications. Journal of the Operational Research
  Society 74, 373--383.

\bibitem[{Cassady(1967)}]{Cassady67_BC}
Cassady, R., 1967. Auctions and Auctioneering. University of California Press,
  Berkeley, CA.

\bibitem[{Castelnovo et~al.(2022)Castelnovo, Crupi, Greco, Regoli, Penco, and
  Cosentini}]{CasCruGreRegPenCos22_JH}
Castelnovo, A., Crupi, R., Greco, G., Regoli, D., Penco, I.~G., Cosentini,
  A.~C., 2022. A clarification of the nuances in the fairness metrics
  landscape. Scientific Reports 12, article 4209.

\bibitem[{Castillo-Salazar et~al.(2016)Castillo-Salazar, Landa-Silva, and
  Qu}]{Dario_GVBSP}
Castillo-Salazar, A.~J., Landa-Silva, D., Qu, R., 2016. Workforce scheduling
  and routing problems: literature survey and computational study. Annals of
  Operations Research 239, 39--67.

\bibitem[{Cataldo et~al.(2017)Cataldo, Ferrer, Miranda, Rey, and
  Saur{\'e}}]{Cataldo2017-la_JJ}
Cataldo, A., Ferrer, J.-C., Miranda, J., Rey, P.~A., Saur{\'e}, A., 2017. An
  integer programming approach to curriculum-based examination timetabling.
  Annals of Operations Research 258~(2), 369--393.

\bibitem[{Caulkins et~al.(2008)Caulkins, Eelman, Ratnatunga, and
  Schaarsmith}]{Caulkins2008-xg_AJG}
Caulkins, J.~P., Eelman, E., Ratnatunga, M., Schaarsmith, D., 2008. Operations
  {R}esearch and public policy for {A}frica: harnessing the revolution in
  management science instruction. International Transactions in Operational
  Research 15, 151--171.

\bibitem[{Cavaliere et~al.(2022)Cavaliere, Bendotti, and
  Fischetti}]{CavaliereBF2022_CA_MB}
Cavaliere, F., Bendotti, E., Fischetti, M., 2022. An integrated
  local-search/set-partitioning refinement heuristic for the capacitated
  vehicle routing problem. Mathematical Programming Computation 14, 749--779.

\bibitem[{Cavallo et~al.(2019)Cavallo, Ishizaka, Olivieri, and
  Squillante}]{Cavallo2019-sj_AFRH}
Cavallo, B., Ishizaka, A., Olivieri, M.~G., Squillante, M., 2019. Comparing
  inconsistency of pairwise comparison matrices depending on entries. Journal
  of the Operational Research Society 70~(5), 842--850.

\bibitem[{Ceder and Wilson(1986)}]{ceder1986bus_MH}
Ceder, A., Wilson, N.~H., 1986. Bus network design. Transportation Research
  Part B: Methodological 20~(4), 331--344.

\bibitem[{Cegan et~al.(2017)Cegan, Filion, Keisler, and
  Linkov}]{Cegan2017-kf_JL}
Cegan, J.~C., Filion, A.~M., Keisler, J.~M., Linkov, I., 2017. Trends and
  applications of multi-criteria decision analysis in environmental sciences:
  literature review. Environment Systems and Decisions 37~(2), 123--133.

\bibitem[{Cela(2013)}]{C13_SMPT}
Cela, E., 2013. The quadratic assignment problem: theory and algorithms.
  Springer Science \& Business Media, Berlin, Heidelberg.

\bibitem[{{\c{C}}elik(2016)}]{CELIK201647_BYKOK}
{\c{C}}elik, M., 2016. Network restoration and recovery in humanitarian
  operations: Framework, literature review, and research directions. Surveys in
  Operations Research and Management Science 21~(2), 47--61.

\bibitem[{{\c{C}}elik et~al.(2012){\c{C}}elik, Ergun, Johnson, Keskinocak,
  Lorca, Pekg{\"u}n, and Swann}]{ccelik2012humanitarian_BYKOK}
{\c{C}}elik, M., Ergun, {\"O}., Johnson, B., Keskinocak, P., Lorca, {\'A}.,
  Pekg{\"u}n, P., Swann, J., 2012. Humanitarian logistics. In: Mirchandani,
  P.~B. (Ed.), New Directions in Informatics, Optimization, Logistics, and
  Production. INFORMS TutORials in Operations Research. INFORMS, pp. 18--49.

\bibitem[{Celikbilek and S{\"u}er(2015)}]{Celikbilek2015-dz_KESXZ}
Celikbilek, C., S{\"u}er, G.~A., 2015. Cellular design-based optimisation for
  manufacturing scheduling and transportation mode decisions. Asian Journal of
  Management Science and Applications 2~(2), 107--129.

\bibitem[{Cervone et~al.(2016)Cervone, D'Amour, Bornn, and
  Goldsberry}]{Cervone2016-gk_IM}
Cervone, D., D'Amour, A., Bornn, L., Goldsberry, K., 2016. A multiresolution
  stochastic process model for predicting basketball possession outcomes.
  Journal of the American Statistical Association 111~(514), 585--599.

\bibitem[{Ceschia et~al.(2019)Ceschia, Dang, {De Causmaecker}, Haspeslagh, and
  Schaerf}]{ceschia2019second_GVBSP}
Ceschia, S., Dang, N., {De Causmaecker}, P., Haspeslagh, S., Schaerf, A., 2019.
  The second international nurse rostering competition. Annals of Operations
  Research 274~(1), 171--186.

\bibitem[{Ceschia et~al.(2022)Ceschia, {Di Gaspero}, and
  Schaerf}]{CESCHIA2022_GVBSP}
Ceschia, S., {Di Gaspero}, L., Schaerf, A., 2022. Educational timetabling:
  Problems, benchmarks, and state-of-the-art results. European Journal of
  Operational Research, DOI: 10.1016/j.ejor.2022.07.011.

\bibitem[{Chae et~al.(2016)Chae, Horesh, Hwang, and
  Lee}]{cha:hor:hwa:lee:16_DSRW}
Chae, Y., Horesh, R., Hwang, Y., Lee, Y., 2016. Artificial neural network model
  for forecasting sub-hourly electricity usage in commercial buildings. Energy
  and Buildings 111, 184--194.

\bibitem[{Chambers et~al.(1996)Chambers, F{\"a}re, and
  Grosskopf}]{Chambers1996-zd_SL}
Chambers, R.~G., F{\"a}re, R., Grosskopf, S., 1996. Productivity growth in
  {APEC} countries. Pacific Economic Review 1~(3), 181--190.

\bibitem[{Chang et~al.(2007)Chang, Hu, Fu, and
  Marcus}]{chang2007simulation_DLLD}
Chang, H.~S., Hu, J., Fu, M.~C., Marcus, S.~I., 2007. Simulation-based
  algorithms for Markov decision processes. Springer.

\bibitem[{Chang et~al.(2022)Chang, Keblis, Li, Iakovou, and
  White}]{Chang2022-xf_KVRH}
Chang, Y., Keblis, M.~F., Li, R., Iakovou, E., White, C.~C., 2022.
  Misinformation and disinformation in modern warfare. Operations Research
  70~(3), 1577--1597.

\bibitem[{Chao and Wilson(2002)}]{Chao02_BC}
Chao, H.-P., Wilson, R., 2002. Multi-dimensional procurement auctions for power
  reserves: Robust incentive-compatible scoring and settlement rules. Journal
  of Regulatory Economics 22~(2), 161--183.

\bibitem[{Chao et~al.(2023)Chao, Chen, and Zhang}]{chao2023online_JSS}
Chao, X., Chen, B., Zhang, H., 2023. Online learning in inventory and pricing
  optimization. In: Song, J.-S. (Ed.), Research Handbook on Inventory
  Management. Edward Elgar Publishing.

\bibitem[{Charnes and Cooper(1962)}]{ChaCoo62_JH}
Charnes, A., Cooper, W.~W., 1962. Programming with linear fractional
  functionals. Naval Research Logistics Quarterly 9, 181--186.

\bibitem[{Charnes et~al.(1978)Charnes, Cooper, and Rhodes}]{Charnes1978-nm_SL}
Charnes, A., Cooper, W.~W., Rhodes, E., 1978. Measuring the efficiency of
  decision making units. European Journal of Operational Research 2~(6),
  429--444.

\bibitem[{Charnes et~al.(1958)Charnes, Cooper, and Symonds}]{Charnes1958-nv_HL}
Charnes, A., Cooper, W.~W., Symonds, G.~H., 1958. Cost horizons and certainty
  equivalents: An approach to stochastic programming of heating oil. Management
  Science 4~(3), 235--263.

\bibitem[{{Chazelle}(1983)}]{Chazelle1983-fp_JB}
{Chazelle}, 1983. The {Bottomn-Left} {Bin-Packing} heuristic: An efficient
  implementation. IEEE Transactions on Computers C-32~(8), 697--707.

\bibitem[{Checkland(1981)}]{Checkland1981-op_MYLW}
Checkland, P., 1981. Systems Thinking, Systems Practice. Wiley, Chichester.

\bibitem[{Checkland(1983)}]{Checkland1983-st_MYLW}
Checkland, P., 1983. {O.R}. and the systems movement: Mappings and conflicts.
  Journal of the Operational Research Society 34~(8), 661--675.

\bibitem[{Checkland(1985)}]{Checkland1985-mi_GM}
Checkland, P., 1985. From optimizing to learning: A development of systems
  thinking for the 1990s. Journal of the Operational Research Society 36~(9),
  757--767.

\bibitem[{Checkland(1989)}]{Checkland1989-fc_MYLW}
Checkland, P., 1989. Soft systems methodology. In: Rosenhead, J. (Ed.),
  Rational analysis for a problematic world: Problem structuring methods for
  complexity, uncertainty, and conflict. Wiley, Chichester, pp. 71--100.

\bibitem[{Checkland(2006)}]{Checkland2006-qc_MYLW}
Checkland, P., 2006. Reply to {E}den and {A}ckermann: Any future for problem
  structuring methods? Journal of the Operational Research Society 57~(7),
  769--771.

\bibitem[{Checkland and Poulter(2006)}]{Checkland2006-tv_JEB}
Checkland, P., Poulter, J., 2006. Learning For Action: A Short Definitive
  Account of Soft Systems Methodology, and Its Use for Practitioners, Teachers
  and Students. Wiley.

\bibitem[{Checkland and Scholes(1990)}]{Checkland1990-ip_MYLW}
Checkland, P., Scholes, J., 1990. Soft Systems Methodology in Action. Wiley,
  Chichester.

\bibitem[{Chekuri et~al.(1997)Chekuri, Goldberg, Karger, Levine, and
  Stein}]{Chekuri-et-al:1997_IL}
Chekuri, C., Goldberg, A.~V., Karger, D.~R., Levine, M.~S., Stein, C., 1997.
  Experimental study of minimum cut algorithms. In: Saks, M.~E. (Ed.),
  Proceedings of the Eighth Annual {ACM-SIAM} Symposium on Discrete Algorithms,
  5-7 January 1997, New Orleans, Louisiana, {USA}. {ACM/SIAM}, pp. 324--333.

\bibitem[{Chen et~al.(2011)Chen, Batson, and Dang}]{CBD11_ALAL}
Chen, D.-S., Batson, R., Dang, Y., 2011. Applied Integer Programming. Wiley,
  Hoboken, NJ.

\bibitem[{Chen et~al.(2000)Chen, Drezner, Ryan, and Simchi-Levi}]{Chen2000_SMD}
Chen, F., Drezner, Z., Ryan, J.~K., Simchi-Levi, D., 2000. Quantifying the
  bullwhip effect in a simple supply chain: The impact of forecasting, lead
  times, and information. Management Science 46~(3), 436--443.

\bibitem[{Chen and Song(2001)}]{chen2001optimal_JSS}
Chen, F., Song, J.-S., 2001. Optimal policies for multiechelon inventory
  problems with markov-modulated demand. Operations Research 49~(2), 226--234.

\bibitem[{Chen and Zheng(1994)}]{chen1994lower_JSS}
Chen, F., Zheng, Y.-S., 1994. Lower bounds for multi-echelon stochastic
  inventory systems. Management Science 40~(11), 1426--1443.

\bibitem[{Chen(1995)}]{chen1995fluid_HATI}
Chen, H., 1995. Fluid approximations and stability of multiclass queueing
  networks: work-conserving disciplines. The Annals of Applied Probability
  5~(3), 637--665.

\bibitem[{Chen and Yao(2001)}]{chen2001fundamentals_HATI}
Chen, H., Yao, D.~D., 2001. Fundamentals of Queueing Networks: Performance,
  Asymptotics, and Optimization. Springer.

\bibitem[{Chen et~al.(2022)Chen, Wandelt, Dai, and Sun}]{chen2022scalable_VLVV}
Chen, L., Wandelt, S., Dai, W., Sun, X., 2022. Scalable vertiport hub location
  selection for air taxi operations in a metropolitan region. INFORMS Journal
  on Computing 34~(2), 834--856.

\bibitem[{Chen and Wang(2022)}]{Chen2022-ms_SL}
Chen, L., Wang, Y.-M., 2022. Data envelopment analysis cross-efficiency method
  of non-homogeneous decision-making units. Journal of the Operational Research
  Society, DOI: 10.1080/01605682.2022.2056535.

\bibitem[{Chen et~al.(2013)Chen, Sharma, and Tseng}]{CheShaTse13_JH}
Chen, L.-W., Sharma, P., Tseng, Y.-C., 2013. Dynamic traffic control with
  fairness and throughput optimization using vehicular communications. IEEE
  Journal on Selected Areas in Communications 31, 504--512.

\bibitem[{Chen and Guestrin(2016)}]{Chen16_LCAL}
Chen, T., Guestrin, C., 2016. {XGBoost}: A scalable tree boosting system. In:
  Proceedings of the 22nd ACM SIGKDD International Conference on Knowledge
  Discovery and Data Mining. KDD '16. ACM, New York, NY, USA, pp. 785--794.

\bibitem[{Chen and Hooker(2022{\natexlab{a}})}]{CheHoo21_JH}
Chen, V., Hooker, J.~N., 2022{\natexlab{a}}. Combining leximax fairness and
  efficiency in an optimization model. European Journal of Operational Research
  299, 235--248.

\bibitem[{Chen and Hooker(2022{\natexlab{b}})}]{CheHoo23_JH}
Chen, V., Hooker, J.~N., 2022{\natexlab{b}}. A guide to formulating fairness in
  optimization models. Annals of Operations Research.

\bibitem[{Chen and Disney(2007)}]{Chen2007_XW}
Chen, Y.~F., Disney, S.~M., 2007. The myopic order-up-to policy with a
  proportional feedback controller. International Journal of Production
  Research 45~(2), 351--368.

\bibitem[{Cheung and
  {Simchi-Levi}(2017)}]{cheungThompsonSamplingOnline2017_AKSJF}
Cheung, W.~C., {Simchi-Levi}, D., 2017. Thompson sampling for online
  personalized assortment optimization problems with multinomial logit choice
  models. SSRN Electronic Journal 3075658.

\bibitem[{Cheung and Simchi-Levi(2023)}]{cheung2023statistical_JSS}
Cheung, W.~C., Simchi-Levi, D., 2023. Statistical learning in inventory
  management. In: Song, J.-S. (Ed.), Research Handbook on Inventory Management.
  Edward Elgar Publishing.

\bibitem[{Chi(2000)}]{Chi2000-og_MJE}
Chi, E.~H., 2000. A taxonomy of visualization techniques using the data state
  reference model. In: Proceedings of the {IEEE} Symposium on Information
  Vizualization 2000 (INFOVIS 2000). IEEE Computer Society, USA, pp. 69--75.

\bibitem[{Chierici et~al.(2004)Chierici, Cordone, and Maja}]{Chierici2004_DC}
Chierici, A., Cordone, R., Maja, R., 2004. The demand-dependent optimization of
  regular train timetable. Electronic Notes in Discrete Mathematics 17,
  99--104.

\bibitem[{Childers and Houston(1984)}]{Childers1984-iq_MJE}
Childers, T.~L., Houston, M.~J., 1984. Conditions for a {Picture-Superiority}
  effect on consumer memory. The Journal of Consumer Research 11~(2), 643--654.

\bibitem[{Chintapalli et~al.(2017)Chintapalli, Disney, and
  Tang}]{Chintapalli2017_SMD}
Chintapalli, P., Disney, S.~M., Tang, C.~S., 2017. Coordinating supply chains
  via advance-order discounts, minimum order quantities, and delegations.
  Production and Operations Management 26~(12), 2175--2186.

\bibitem[{Choi(2020)}]{Choi2020_SMD}
Choi, T.~M., 2020. Supply chain financing using blockchain: {I}mpacts on supply
  chains selling fashionable products. Annals of Operations Research, DOI:
  10.1007/s10479-020-03615-7.

\bibitem[{Choi et~al.(2019)Choi, Wen, Sun, and Chung}]{choi2019mean_JLYHK}
Choi, T.-M., Wen, X., Sun, X., Chung, S.-H., 2019. The mean-variance approach
  for global supply chain risk analysis with air logistics in the blockchain
  technology era. Transportation Research Part E: Logistics and Transportation
  Review 127, 178--191.

\bibitem[{Chopra and Rao(1994)}]{Chopra:1994part1_IL}
Chopra, S., Rao, M.~R., 1994. The {S}teiner tree problem {I}: Formulations,
  compositions and extension of facets. Mathematical Programming 64~(1),
  209--229.

\bibitem[{Chothani et~al.(2015)Chothani, Patel, Dekavadiya, and
  Patel}]{Chothani2015-eb_BC}
Chothani, R.~G., Patel, N.~A., Dekavadiya, A.~H., Patel, P.~R., 2015. A review
  of online auction and its pros and cons. International Journal of Advance
  Engineering and Research Development 2~(1), 8--11.

\bibitem[{Chouldechova(2017)}]{chouldechova2017fair_JH}
Chouldechova, A., 2017. Fair prediction with disparate impact: A study of bias
  in recidivism prediction instruments. Big data 5~(2), 153--163.

\bibitem[{Chouman and Crainic(2021)}]{Chouman2021_MH}
Chouman, M., Crainic, T.~G., 2021. Freight railroad service network design. In:
  Crainic, T.~G., Gendreau, M., Gendron, B. (Eds.), Network Design with
  Applications to Transportation and Logistics. Springer, pp. 383--426.

\bibitem[{Chowdhury et~al.(2023)Chowdhury, Gregory, and
  Queah}]{Chowdhury2023-th_AJG}
Chowdhury, R., Gregory, A.~J., Queah, M., 2023. Creative and flexible
  deployment of systems methodologies for child rights and child protection
  through (the conceptual lens of) holistic flexibility. Systems Research and
  Behavioral Science 40~(4).

\bibitem[{Christiaens and Vanden~Berghe(2020)}]{christiaens2020slack_CA_MB}
Christiaens, J., Vanden~Berghe, G., 2020. Slack induction by string removals
  for vehicle routing problems. Transportation Science 54~(2), 417--433.

\bibitem[{Christiansen et~al.(2011)Christiansen, Fagerholt, Flatberg, Haugen,
  Kloster, and Lund}]{christiansen2011maritime_JLYHK}
Christiansen, M., Fagerholt, K., Flatberg, T., Haugen, {\O}., Kloster, O.,
  Lund, E.~H., 2011. Maritime inventory routing with multiple products: A case
  study from the cement industry. European Journal of Operational Research
  208~(1), 86--94.

\bibitem[{Christiansen et~al.(2013)Christiansen, Fagerholt, Nygreen, and
  Ronen}]{Christiansen2013-cg_HP}
Christiansen, M., Fagerholt, K., Nygreen, B., Ronen, D., 2013. Ship routing and
  scheduling in the new millennium. European Journal of Operational Research
  228~(3), 467--483.

\bibitem[{Christiansen et~al.(2020)Christiansen, Hellsten, Pisinger,
  Sacramento, and Vilhelmsen}]{Christiansen2021_MH}
Christiansen, M., Hellsten, E., Pisinger, D., Sacramento, D., Vilhelmsen, C.,
  2020. Liner shipping network design. European Journal of Operational Research
  286~(1), 1--20.

\bibitem[{Christofides(1975)}]{C75_SMPT}
Christofides, N., 1975. Graph Theory: An Algorithmic Approach. Academic Press,
  London.

\bibitem[{Christofides et~al.(1979)Christofides, Mingozzi, Toth, and
  Sandi}]{CMTS79_SMPT}
Christofides, N., Mingozzi, A., Toth, P., Sandi, C. (Eds.), 1979. Combinatorial
  Optimization. Wiley, Chichester.

\bibitem[{Christofides and Whitlock(1981)}]{christofides.whitlock:81_BF}
Christofides, N., Whitlock, C., 1981. Network synthesis with connectivity
  constraints --- a survey. In: Brans, J. (Ed.), Operational Research '81.
  North-Holland, pp. 705--723.

\bibitem[{Chu et~al.(2019)Chu, Xu, and Li}]{Chu2019-bw_HL}
Chu, Z., Xu, Z., Li, H., 2019. New heuristics for the {RCPSP} with multiple
  overlapping modes. Computers \& Industrial Engineering 131, 146--156.

\bibitem[{Churchman(1979)}]{Churchman1979_GM}
Churchman, C., 1979. The systems approach, 2nd Edition. Dell, New York.

\bibitem[{Churchman(1970)}]{Churchman1970-hv_GM}
Churchman, C.~W., 1970. Operations research as a profession. Management Science
  17~(2), B37--B53.

\bibitem[{Churchman et~al.(1957)Churchman, Ackoff, and
  Arnoff}]{Churchman1957-ax_MYLW}
Churchman, C.~W., Ackoff, R., Arnoff, E.~L., 1957. Introduction to Operations
  Research. Wiley, New York, NY.

\bibitem[{Chv\'atal(1973)}]{Ch73_ALAL}
Chv\'atal, V., 1973. Edmonds polytopes and a hierarchy of combinatorial
  problems. Discrete Mathematics 4, 305--337.

\bibitem[{Chv{\'a}tal(1983)}]{Chvatal1983-ri_JMB}
Chv{\'a}tal, V., 1983. Linear Programming. W. H. Freeman, New York/San
  Francisco.

\bibitem[{Ciampi et~al.(2011)Ciampi, Dyachenko, Cole, and
  McCusker}]{Ciampi2011-fi_CV}
Ciampi, A., Dyachenko, A., Cole, M., McCusker, J., 2011. Delirium superimposed
  on dementia: defining disease states and course from longitudinal
  measurements of a multivariate index using latent class analysis and hidden
  {M}arkov chains. International Psychogeriatrics 23~(10), 1659--1670.

\bibitem[{Cinelli et~al.(2014)Cinelli, Coles, and Kirwan}]{Cinelli2014-gf_JL}
Cinelli, M., Coles, S.~R., Kirwan, K., 2014. Analysis of the potentials of
  multi criteria decision analysis methods to conduct sustainability
  assessment. Ecological Indicators 46, 138--148.

\bibitem[{Cinelli et~al.(2020)Cinelli, Kadzi{\'n}ski, Gonzalez, and
  S{\l}owi{\'n}ski}]{Cinelli2020-lq_JL}
Cinelli, M., Kadzi{\'n}ski, M., Gonzalez, M., S{\l}owi{\'n}ski, R., 2020. How
  to support the application of multiple criteria decision analysis? {L}et us
  start with a comprehensive taxonomy. Omega 96, 102261.

\bibitem[{Cioppa et~al.(2004)Cioppa, Lucas, and Sanchez}]{Cioppa2004_JLYHK}
Cioppa, T.~M., Lucas, T.~W., Sanchez, S.~M., 2004. Military applications of
  agent-based simulations. In: Ingalis, R.~G., Rossetti, M.~D., Smith, J.~S.,
  Peters, B.~A. (Eds.), Proceedings of the 2004 Winter Simulation Conference,
  2004. Vol.~1. p. 180.

\bibitem[{Claeskens et~al.(2016)Claeskens, Magnus, Vasnev, and
  Wang}]{Claeskens2016-hr}
Claeskens, G., Magnus, J.~R., Vasnev, A.~L., Wang, W., 2016. The forecast
  combination puzzle: A simple theoretical explanation. International Journal
  of Forecasting 32~(3), 754--762.

\bibitem[{Clark and Scarf(1960)}]{clark1960optimal_JSS}
Clark, A.~J., Scarf, H., 1960. Optimal policies for a multi-echelon inventory
  problem. Management Science 6~(4), 475--490.

\bibitem[{Clarke and Wright(1964)}]{ClarkeW1964_CA_MB}
Clarke, G., Wright, J., 1964. Scheduling of vehicles from a central depot to a
  number of delivery points. Operations Research 12, 568--581.

\bibitem[{Cla{\ss}en et~al.(2014)Cla{\ss}en, Koster, Coudert, and
  Nepomuceno}]{Clasen2014-ki_HL}
Cla{\ss}en, G., Koster, A. M. C.~A., Coudert, D., Nepomuceno, N., 2014.
  {Chance-Constrained} optimization of reliable fixed broadband wireless
  networks. INFORMS Journal on Computing 26~(4), 893--909.

\bibitem[{Clautiaux et~al.(2007)Clautiaux, Carlier, and
  Moukrim}]{Clautiaux2007-ti_JB}
Clautiaux, F., Carlier, J., Moukrim, A., 2007. A new exact method for the
  two-dimensional orthogonal packing problem. European Journal of Operational
  Research 183~(3), 1196--1211.

\bibitem[{Clemen(1996)}]{clemen1996making_MESG}
Clemen, R., 1996. Making hard Decisions: An Introduction to Decision Analysis.
  Brooks/Cole, Monterey.

\bibitem[{Cleophas and Ehmke(2014)}]{cleophas2014deliveries_CKTVW}
Cleophas, C., Ehmke, J.~F., 2014. When are deliveries profitable? Business \&
  Information Systems Engineering 6~(3), 153--163.

\bibitem[{Coates and Parshakov(2022)}]{Coates2022-sm_IM}
Coates, D., Parshakov, P., 2022. The wisdom of crowds and transfer market
  values. European Journal of Operational Research 301~(2), 523--534.

\bibitem[{Cobacho et~al.(2010)Cobacho, Caballero, Gonz{\'a}lez, and
  Molina}]{Cobacho2010-mn_JJ}
Cobacho, B., Caballero, R., Gonz{\'a}lez, M., Molina, J., 2010. Planning
  federal public investment in mexico using multiobjective decision making.
  Journal of the Operational Research Society 61~(9), 1328--1339.

\bibitem[{Cochrane(2009)}]{cochrane09_MB}
Cochrane, J., 2009. Asset pricing: Revised edition. Princeton University Press.

\bibitem[{Coelho et~al.(2014)Coelho, Cordeau, and
  Laporte}]{coelho2014heuristics_JLYHK}
Coelho, L.~C., Cordeau, J.-F., Laporte, G., 2014. Heuristics for dynamic and
  stochastic inventory-routing. Computers \& Operations Research 52, 55--67.

\bibitem[{Coffman~Jr. et~al.(1980)Coffman~Jr., Garey, Johnson, and
  Tarjan}]{Coffman1980-vj_JB}
Coffman~Jr., G.~E., Garey, M.~R., Johnson, D.~S., Tarjan, R.~E., 1980.
  Performance bounds for {Level-Oriented} {Two-Dimensional} packing algorithms.
  SIAM Journal on Computing 9~(4), 808--826.

\bibitem[{{COIN-OR Foundation, Inc.}(2022)}]{Coinor2_CTCGE}
{COIN-OR Foundation, Inc.}, 2022. {The Computational Infrastructure for
  Operations Research}.
\newline\urlprefix\url{https://www.coin-or.org/}

\bibitem[{Colapinto et~al.(2020)Colapinto, Jayaraman, Ben~Abdelaziz, and
  La~Torre}]{Colapinto2020-mc_JL}
Colapinto, C., Jayaraman, R., Ben~Abdelaziz, F., La~Torre, D., 2020.
  Environmental sustainability and multifaceted development: multi-criteria
  decision models with applications. Annals of Operations Research 293~(2),
  405--432.

\bibitem[{Comi et~al.(2014)Comi, Bischof, and Eppler~Martin}]{Comi2014-rt_MJE}
Comi, A., Bischof, N., Eppler~Martin, J., 2014. Beyond projection: using
  collaborative visualization to conduct qualitative interviews. Qualitative
  Research in Organizations and Management: An International Journal 9~(2),
  110--133.

\bibitem[{Conforti et~al.(2014)Conforti, Cornu{\'e}jols, and
  Zambelli}]{CCZ14_ALAL}
Conforti, M., Cornu{\'e}jols, G., Zambelli, G., 2014. Integer Programming.
  Springer, Cham, Switzerland.

\bibitem[{Congram et~al.(2002)Congram, Potts, and van~de
  Velde}]{congram2002iterated_COIT}
Congram, R.~K., Potts, C.~N., van~de Velde, S.~L., 2002. An iterated dynasearch
  algorithm for the single-machine total weighted tardiness scheduling problem.
  INFORMS Journal on Computing 14~(1), 52--67.

\bibitem[{Contreras(2020)}]{Contreras2020-rn_SL}
Contreras, I., 2020. A review of the literature on {DEA} models under common
  set of weights. Journal of Modelling in Management 15~(4), 1277--1300.

\bibitem[{Contreras and Fern{\'a}ndez(2012)}]{contreras.fernandez:12_BF}
Contreras, I., Fern{\'a}ndez, E., 2012. General network design: A unified view
  of combined location and network design problems. European Journal of
  Operational Research 219, 680--697.

\bibitem[{Cook(1971)}]{Cook:STOC71_UPCT}
Cook, S.~A., 1971. The complexity of theorem-proving procedures. In:
  Proceedings of the 3rd {ACM} Symposium on the Theory of Computing ({STOC}).
  pp. 151--158.

\bibitem[{Cook(1973)}]{Cook1973-fz_AJG}
Cook, S.~L., 1973. Operational research, social well-being and the zero growth
  concept. Omega 1~(6), 647--667.

\bibitem[{Cook(1984)}]{Cook1984-cu_AJG}
Cook, S.~L., 1984. Operational research, social well-being and the zero growth
  concept. In: Bowen, K., Cook, A., Luck, M. (Eds.), The Writings of Steve
  Cook. Operational Research Society, Birmingham, pp. 34--51.

\bibitem[{Cook et~al.(1998)Cook, Cunningham, Pulleyblank, and
  Schrijver}]{CCPS98_SMPT}
Cook, W., Cunningham, W., Pulleyblank, W., Schrijver, A., 1998. Combinatorial
  Optimization. Wiley, Chichester.

\bibitem[{Cook et~al.(2021)Cook, Held, and Helsgaun}]{Cook-et-al:2022_IL}
Cook, W., Held, S., Helsgaun, K., 2021. Constrained local search for last-mile
  routing. arXiv:2112.15192.

\bibitem[{Cook et~al.(2010)Cook, Liang, and Zhu}]{Cook2010-bk_SL}
Cook, W.~D., Liang, L., Zhu, J., 2010. Measuring performance of two-stage
  network structures by {DEA}: A review and future perspective. Omega 38~(6),
  423--430.

\bibitem[{Cook et~al.(2014)Cook, Tone, and Zhu}]{Cook2014-xo_SL}
Cook, W.~D., Tone, K., Zhu, J., 2014. Data envelopment analysis: Prior to
  choosing a model. Omega 44, 1--4.

\bibitem[{Cook and Zhu(2014)}]{Cook2014-fa_SL}
Cook, W.~D., Zhu, J., 2014. Data Envelopment Analysis: A Handbook of Modeling
  Internal Structure and Network. Springer, New York, NY.

\bibitem[{Cook(2011)}]{Cook:2011_IL}
Cook, W.~J., 2011. In pursuit of the traveling salesman. Princeton University
  Press, Princeton.

\bibitem[{Cooper(1972)}]{cooper_HATI}
Cooper, R.~B., 1972. Introduction to Queueing Theory. Macmillan, New York.

\bibitem[{Cooper et~al.(2007)Cooper, Seiford, and Tone}]{Cooper2007-pc_SL}
Cooper, W.~W., Seiford, L.~M., Tone, K., 2007. Data Envelopment Analysis: A
  Comprehensive Text with Models, Applications, References and {DEA-Solver}
  Software. Springer, New York, NY.

\bibitem[{Cooper et~al.(2011)Cooper, Seiford, and Zhu}]{Cooper2011-sx_SL}
Cooper, W.~W., Seiford, L.~M., Zhu, J., 2011. Handbook on Data Envelopment
  Analysis. Springer, New York, NY.

\bibitem[{Corbett-Davies et~al.(2017)Corbett-Davies, Pierson, Feller, Goel, and
  Huq}]{Corbett-Davies2017-ls_TC}
Corbett-Davies, S., Pierson, E., Feller, A., Goel, S., Huq, A., 2017.
  Algorithmic decision making and the cost of fairness. In: Proceedings of the
  23rd {ACM} {SIGKDD} International Conference on Knowledge Discovery and Data
  Mining. KDD '17. Association for Computing Machinery, New York, NY, USA, pp.
  797--806.

\bibitem[{Cordeau and Laporte(2005)}]{Cordeau2005_CA_MB}
Cordeau, J.-F., Laporte, G., 2005. Tabu search heuristics for the vehicle
  routing problem. In: Sharda, R., Vo{\ss}, S., Rego, C., Alidaee, B. (Eds.),
  Metaheuristic Optimization via Memory and Evolution: Tabu Search and Scatter
  Search. Springer US, Boston, MA, pp. 145--163.

\bibitem[{Cordeau et~al.(2000)Cordeau, Soumis, and Desrosiers}]{Cordeau2000_DC}
Cordeau, J.-F., Soumis, F., Desrosiers, J., 2000. A {B}enders decomposition
  approach for the locomotive and car assignment problem. Transportation
  Science 34~(2), 133--149.

\bibitem[{Cordeau et~al.(2001)Cordeau, Stojkovi{\'c}, Soumis, and
  Desrosiers}]{cordeau2001benders_VLVV}
Cordeau, J.-F., Stojkovi{\'c}, G., Soumis, F., Desrosiers, J., 2001. Benders
  decomposition for simultaneous aircraft routing and crew scheduling.
  Transportation Science 35~(4), 375--388.

\bibitem[{C{\'o}rdoba and Midgley(2006)}]{Cordoba2006-tv_AJG}
C{\'o}rdoba, J.-R., Midgley, G., 2006. Broadening the boundaries: an
  application of critical systems thinking to {IS} planning in {C}olombia.
  Journal of the Operational Research Society 57~(9), 1064--1080.

\bibitem[{Corlu et~al.(2020)Corlu, Akcay, and Xie}]{Corlu2020_CC}
Corlu, C.~G., Akcay, A., Xie, W., 2020. Stochastic simulation under input
  uncertainty: A review. Operations Research Perspectives 7, 100162.

\bibitem[{Cormen et~al.(2022)Cormen, Leiserson, Rivest, and
  Stein}]{cormen01introduction_IL}
Cormen, T.~H., Leiserson, C.~E., Rivest, R.~L., Stein, C., 2022. Introduction
  to Algorithms, 4th Edition. MIT Press.

\bibitem[{Cormen et~al.(2009)Cormen, Stein, Leiserson, and Rivest}]{Cor09_UPCT}
Cormen, T.~H., Stein, C., Leiserson, C.~E., Rivest, R.~L., 2009. Introduction
  to Algorithms, 3rd Edition. MIT Press, Cambridge, MA.

\bibitem[{Cornu{\'e}jols(2008)}]{Co08_ALAL}
Cornu{\'e}jols, G., 2008. Valid inequalities for mixed integer linear programs.
  Mathematical Programming 112, 3--44.

\bibitem[{Cornu{\'e}jols and Dawande(1999)}]{CD99_ALAL}
Cornu{\'e}jols, G., Dawande, M., 1999. A class of hard small 0-1 programs.
  INFORMS Journal on Computing 11, 205--210.

\bibitem[{Cornuejols and T{\"u}t{\"u}nc{\"u}(2006)}]{cornuejols06_MB}
Cornuejols, G., T{\"u}t{\"u}nc{\"u}, R., 2006. Optimization methods in finance.
  Cambridge University Press.

\bibitem[{Cortazar et~al.(1998)Cortazar, Schwartz, and
  Salinas}]{Cortazar1998_EA}
Cortazar, G., Schwartz, E.~S., Salinas, M., 1998. Evaluating environmental
  investments: A real options approach. Management Science 44~(8), 1059--1070.

\bibitem[{Costa(2005)}]{costa2005survey_MH}
Costa, A.~M., 2005. A survey on benders decomposition applied to fixed-charge
  network design problems. Computers \& Operations Research 32~(6), 1429--1450.

\bibitem[{Costa et~al.(2009)Costa, Gomes, and Oliveira}]{Costa2009-il_JB}
Costa, M.~T., Gomes, A.~M., Oliveira, J.~F., 2009. Heuristic approaches to
  large-scale periodic packing of irregular shapes on a rectangular sheet.
  European Journal of Operational Research 192~(1), 29--40.

\bibitem[{{Council for Science and Technology}(2020)}]{GOVUK_2020-yk}
{Council for Science and Technology}, 2020. Achieving net zero carbon emissions
  through a whole systems approach.
  \url{https://www.gov.uk/government/publications/achieving-net-zero-carbon-emissions-through-a-whole-systems-approach},
  accessed on 2023-01-12.

\bibitem[{Cowell(2000)}]{Cow00_JH}
Cowell, F.~A., 2000. Measurement of inequality. In: Atkinson, A.~B.,
  Bourguignon, F. (Eds.), Handbook of Income Distribution. Vol.~1. Elsevier,
  pp. 89--166.

\bibitem[{Cox et~al.(1979)Cox, Ross, and Rubinstein}]{cox79_MB}
Cox, J.~C., Ross, S.~A., Rubinstein, M., 1979. Option pricing: A simplified
  approach. Journal of Financial Economics 7~(3), 229--263.

\bibitem[{Cox(2020)}]{Cox2020-sj_TC}
Cox, Jr, L.~A., 2020. Answerable and unanswerable questions in risk analysis
  with open-world novelty. Risk Analysis 40~(S1), 2144--2177.

\bibitem[{Cox et~al.(2018)Cox, Popken, and Sun}]{Cox2018-ui_TC}
Cox, Jr, L.~A., Popken, D.~A., Sun, R.~X., 2018. Causal Analytics for Applied
  Risk Analysis. Springer.

\bibitem[{Craig and Winchester(2021)}]{Craig2021-fc_IM}
Craig, J.~D., Winchester, N., 2021. Predicting the national football league
  potential of college quarterbacks. European Journal of Operational Research
  295~(2), 733--743.

\bibitem[{Crainic et~al.(1984)Crainic, Ferland, and Rousseau}]{Crainic1984_DC}
Crainic, T.~G., Ferland, J.-A., Rousseau, J.-M., 1984. A tactical planning
  model for rail freigth transportation. Transportation Science 18~(2),
  165--184.

\bibitem[{Crainic et~al.(1990)Crainic, Florian, and L{\'e}al}]{Crainic1990_DC}
Crainic, T.~G., Florian, M., L{\'e}al, J.~E., 1990. A model for the strategic
  planning of national freight transportation by rail. Transportation Science
  24~(1), 1--24.

\bibitem[{Crainic et~al.(2011)Crainic, Fu, Gendreau, Rei, and
  Wallace}]{crainic2011progressive_MH}
Crainic, T.~G., Fu, X., Gendreau, M., Rei, W., Wallace, S.~W., 2011.
  Progressive hedging-based metaheuristics for stochastic network design.
  Networks 58~(2), 114--124.

\bibitem[{Crainic and Gendreau(2021)}]{Crainic2021Heuristic_MH}
Crainic, T.~G., Gendreau, M., 2021. Heuristics and metaheuristics for
  fixed-charge network design. In: Crainic, T.~G., Gendreau, M., Gendron, B.
  (Eds.), Network Design with Applications to Transportation and Logistics.
  Springer, pp. 91--138.

\bibitem[{Crainic et~al.(2021{\natexlab{a}})Crainic, Gendreau, and
  Gendron}]{Crainic2021MIP_MH}
Crainic, T.~G., Gendreau, M., Gendron, B., 2021{\natexlab{a}}. Fixed-charge
  network design problems. In: Crainic, T.~G., Gendreau, M., Gendron, B.
  (Eds.), Network Design with Applications to Transportation and Logistics.
  Springer, Cham, pp. 15--28.

\bibitem[{Crainic et~al.(2021{\natexlab{b}})Crainic, Gendreau, and
  Gendron}]{crainic2021network_MH}
Crainic, T.~G., Gendreau, M., Gendron, B., 2021{\natexlab{b}}. Network Design
  with Applications to Transportation and Logistics. Springer.

\bibitem[{Crainic and Gendron(2021)}]{Crainic2021Exact_MH}
Crainic, T.~G., Gendron, B., 2021. Exact methods for fixed-charge network
  design. In: Crainic, T.~G., Gendreau, M., Gendron, B. (Eds.), Network Design
  with Applications to Transportation and Logistics. Springer, pp. 29--89.

\bibitem[{Crainic et~al.(2021{\natexlab{c}})Crainic, Hewitt, Maggioni, and
  Rei}]{crainic2021partial_MH}
Crainic, T.~G., Hewitt, M., Maggioni, F., Rei, W., 2021{\natexlab{c}}. Partial
  benders decomposition: general methodology and application to stochastic
  network design. Transportation Science 55~(2), 414--435.

\bibitem[{Crainic et~al.(2014)Crainic, Hewitt, and
  Rei}]{crainic2014scenario_MH}
Crainic, T.~G., Hewitt, M., Rei, W., 2014. Scenario grouping in a progressive
  hedging-based meta-heuristic for stochastic network design. Computers \&
  Operations Research 43, 90--99.

\bibitem[{Crainic et~al.(2016)Crainic, Hewitt, Toulouse, and
  Vu}]{HewittVuCrainic2016_MH}
Crainic, T.~G., Hewitt, M., Toulouse, M., Vu, D.~M., 2016. Service network
  design with resource constraints. Transportation Science 50~(4), 1380--1393.

\bibitem[{Crainic et~al.(2018)Crainic, Hewitt, Toulouse, and
  Vu}]{Crainic2018_MH}
Crainic, T.~G., Hewitt, M., Toulouse, M., Vu, D.~M., 2018. Scheduled service
  network design with resource acquisition and management. EURO Journal on
  Transportation and Logistics 7~(3), 277--309.

\bibitem[{Crama et~al.(2022)Crama, Rezaei, Savelsbergh, and
  Van~Woensel}]{Crama2022_JLYHK}
Crama, Y., Rezaei, M., Savelsbergh, M., Van~Woensel, T., 2022. Stochastic
  inventory routing for perishable products. Transportation Science 52~(3),
  526--546.

\bibitem[{{CreateASoft Inc.}(2022)}]{Simcad_CTCGE}
{CreateASoft Inc.}, 2022. {SimCAD}.
\newline\urlprefix\url{https://www.createasoft.com/simcad-pro-healthcare-simulation-software}

\bibitem[{Creemers et~al.(2014)Creemers, Demeulemeester, and Van~de
  Vonder}]{Creemers2014-lb_WH_ED}
Creemers, S., Demeulemeester, E., Van~de Vonder, S., 2014. A new approach for
  quantitative risk analysis. Annals of Operations Research 213~(1), 27--65.

\bibitem[{Cressman(2003)}]{Cr2003_GZ}
Cressman, R., 2003. Evolutionary Dynamics and Extensive Form Games. MIT Press.

\bibitem[{Crouhy et~al.(2000)Crouhy, Galai, and Mark}]{crouhy00_MB}
Crouhy, M., Galai, D., Mark, R., 2000. A comparative analysis of current credit
  risk models. Journal of Banking \& Finance 24~(1-2), 59--117.

\bibitem[{Crowder et~al.(1983)Crowder, Johnson, and Padberg}]{CJP83_ALAL}
Crowder, H., Johnson, E., Padberg, M., 1983. Solving large-scale 0-1 linear
  programming problems. Discrete Mathematics 31, 803--834.

\bibitem[{Crowe and Utley(2022)}]{Crowe2022-bc_CV}
Crowe, S., Utley, M., 2022. Praxis in healthcare {OR}: An empirical behavioural
  {OR} study. Journal of the Operational Research Society 73~(7), 1444--1456.

\bibitem[{Crowe et~al.(2014)Crowe, Vasilakis, Skeen, Storr, Grove, Gallivan,
  and Utley}]{Crowe2014-gr_CV}
Crowe, S., Vasilakis, C., Skeen, A., Storr, P., Grove, P., Gallivan, S., Utley,
  M., 2014. Examining the feasibility of using a modelling tool to assess
  resilience across a health-care system and assist with decisions concerning
  service reconfiguration. Journal of the Operational Research Society 65~(10),
  1522--1532.

\bibitem[{Cui and Wu(2018)}]{Cui2018-xl_AFRH}
Cui, T.~H., Wu, Y., 2018. Incorporating behavioral factors into operations
  theory. In: The Handbook of Behavioral Operations. Wiley, Hoboken, NJ, pp.
  89--119.

\bibitem[{Currie et~al.(2020)Currie, Fowler, Kotiadis, Monks, Onggo, Robertson,
  and Tako}]{Currie2020-tz_CV}
Currie, C. S.~M., Fowler, J.~W., Kotiadis, K., Monks, T., Onggo, B.~S.,
  Robertson, D.~A., Tako, A.~A., 2020. How simulation modelling can help reduce
  the impact of {COVID-19}. Journal of Simulation 14~(2), 83--97.

\bibitem[{Daduna(2001)}]{daduna_HATI}
Daduna, H., 2001. Queueing Networks with Discrete Time Scaling.
  Springer-Verlag, Berlin, Heidelberg.

\bibitem[{Dagkakis and Heavey(2016)}]{Dagkakis2016_CTCGE}
Dagkakis, G., Heavey, C., 2016. A review of open source discrete event
  simulation software for operations research. Journal of Simulation 10~(3),
  193--206.

\bibitem[{Dahl et~al.(1993{\natexlab{a}})Dahl, Meeraus, and
  Zenios}]{dahl93I_MB}
Dahl, H., Meeraus, A., Zenios, S.~A., 1993{\natexlab{a}}. Some financial
  optimization models: {I. Risk Management}. In: Zenios, S.~A. (Ed.), Financial
  Optimization. Cambridge University Press, pp. 3--36.

\bibitem[{Dahl et~al.(1993{\natexlab{b}})Dahl, Meeraus, and
  Zenios}]{dahl93II_MB}
Dahl, H., Meeraus, A., Zenios, S.~A., 1993{\natexlab{b}}. Some financial
  optimization models: {II. Financial Engineering}. In: Zenios, S.~A. (Ed.),
  Financial Optimization. Cambridge University Press, pp. 37--71.

\bibitem[{Dai(1995)}]{dai1995positive_HATI}
Dai, J.~G., 1995. On positive {H}arris recurrence of multiclass queueing
  networks: a unified approach via fluid limit models. The Annals of Applied
  Probability 5~(1), 49--77.

\bibitem[{Dai and Meyn(1995)}]{dai1995stability_HATI}
Dai, J.~G., Meyn, S.~P., 1995. Stability and convergence of moments for
  multiclass queueing networks via fluid limit models. IEEE Transactions on
  Automatic Control 40~(11), 1889--1904.

\bibitem[{Daley and Rolski(1992)}]{daley1992light_HATI}
Daley, D., Rolski, T., 1992. Light traffic approximations in many-server
  queues. Advances in Applied Probability 24~(1), 202--218.

\bibitem[{Dammon et~al.(2001)Dammon, Spatt, and Zhang}]{dammon01_MB}
Dammon, R.~M., Spatt, C.~S., Zhang, H.~H., 2001. Optimal consumption and
  investment with capital gains taxes. The Review of Financial Studies 14~(3),
  583--616.

\bibitem[{Damodaran and Wagner(2020)}]{Damodaran2020-us_KVRH}
Damodaran, S.~K., Wagner, N., 2020. Modeling and simulation to support cyber
  defense. The Journal of Defense Modeling and Simulation 17~(1), 3--4.

\bibitem[{Dando and Bennett(1981)}]{Dando1981-hv_MYLW}
Dando, M.~R., Bennett, P.~G., 1981. {A Kuhnian Crisis in Management Science?}
  Journal of the Operational Research Society 32~(2), 91--103.

\bibitem[{Danna et~al.(2005)Danna, Rothberg, and Le~Pape}]{DRL05_ALAL}
Danna, E., Rothberg, E., Le~Pape, C., 2005. Exploring relaxation induced
  neighborhoods to improve {MIP} solutions. Mathematical Programming 102,
  71--90.

\bibitem[{Dantzig(1960)}]{Da60_ALAL}
Dantzig, G., 1960. On the significance of solving linear programming problems
  with some integer variables. Econometrica 28, 30--44.

\bibitem[{Dantzig et~al.(1954)Dantzig, Fulkerson, and Johnson}]{DFJ:1954_IL}
Dantzig, G., Fulkerson, R., Johnson, S., 1954. Solution of a large-scale
  traveling-salesman problem. Journal of the Operations Research Society of
  America 2~(4), 393--410.

\bibitem[{Dantzig and Ramser(1959)}]{DantzigR1959_CA_MB}
Dantzig, G., Ramser, J., 1959. The truck dispatching problem. Management
  Science 6, 80--91.

\bibitem[{Dantzig(1951)}]{Dantzig1951-qp_GL}
Dantzig, G.~B., 1951. Maximization of a linear function of variables subject to
  linear inequalities. In: Koopmans, T.~C. (Ed.), Activity Analysis of
  Production and Allocation. Wiley, pp. 339--347.

\bibitem[{Dantzig(1955)}]{Dantzig1955-um_HL}
Dantzig, G.~B., 1955. Linear programming under uncertainty. Management Science
  1~(3/4), 197--206.

\bibitem[{Dantzig(1963)}]{Dantzig1963-sy_JMB}
Dantzig, G.~B., 1963. Linear Programming and Extensions. Princeton University
  Press, Princeton.

\bibitem[{Dantzig(1982)}]{Dantzig1982-rx_JMB}
Dantzig, G.~B., 1982. Reminiscences about the origins of linear programming.
  Operations Research Letters 1~(2), 43--48.

\bibitem[{Dantzig(1990)}]{Dantzig1990_EAY}
Dantzig, G.~B., 1990. Origins of the simplex method. In: Nash, S.~G. (Ed.), A
  History of Scientific Computing. Association for Computing Machinery, New
  York, NY, p. 141–151.

\bibitem[{Dantzig(1991)}]{Dantzig1991-tc}
Dantzig, G.~B., 1991. Linear programming. In: Lenstra, J.~K., Rinnooy~Kan,
  A.~H., Schrijver, A. (Eds.), History of Mathematical Programming: A
  Collection of Personal Reminiscences. CWI, North-Holland, Amsterdam.

\bibitem[{Dantzig and Fulkerson(1955)}]{Dantzig-Fulkerson:1955_IL}
Dantzig, G.~B., Fulkerson, D.~R., 1955. On the max flow min cut theorem of
  networks. Tech. rep., RAND Corporation, Santa Monica, CA.

\bibitem[{Dantzig and Infanger(1991)}]{Dantzig1991-yh_HL}
Dantzig, G.~B., Infanger, G., 1991. {Large-Scale} stochastic linear programs:
  Importance sampling and benders decomposition. Tech. Rep. ADA234962, Systems
  Optimization Laboratoty, Stanford University.

\bibitem[{Daraio and Simar(2007)}]{Daraio2007-qr_JJ}
Daraio, C., Simar, L., 2007. Conditional nonparametric frontier models for
  convex and nonconvex technologies: a unifying approach. Journal of
  Productivity Analysis 28~(1), 13--32.

\bibitem[{Daraio et~al.(2018)Daraio, Simar, and Wilson}]{Daraio2018-ut_JJ}
Daraio, C., Simar, L., Wilson, P.~W., 2018. Central limit theorems for
  conditional efficiency measures and tests of the `separability' condition in
  non‐parametric, two‐stage models of production. The Econometrics Journal
  21~(2), 170--191.

\bibitem[{Dasgupta et~al.(2022)Dasgupta, Akhtar, and
  Sen}]{Dasgupta2022-mj_KVRH}
Dasgupta, D., Akhtar, Z., Sen, S., 2022. Machine learning in cybersecurity: a
  comprehensive survey. The Journal of Defense Modeling and Simulation 19~(1),
  57--106.

\bibitem[{Daskin(1995)}]{daskin1995network_SAA}
Daskin, M., 1995. Network and Discrete Location: Models, Algorithms, and
  Applications. Wiley.

\bibitem[{Datta(1995)}]{Datta1995_EA}
Datta, S., 1995. A decision support system for micro-watershed management in
  {I}ndia. Journal of the Operational Research Society 46~(5), 592--603.

\bibitem[{Davenport(2013)}]{Davenport2013-xh_JEB}
Davenport, T.~H., 2013. Analytics 3.0. Harvard Business Review.
\newline\urlprefix\url{https://hbr.org/2013/12/analytics-30}

\bibitem[{Davenport and Harris(2007)}]{Davenport2007-jn_JEB}
Davenport, T.~H., Harris, J.~G., 2007. Competing on Analytics: The New Science
  of Winning. Harvard Business Press.

\bibitem[{Davenport et~al.(2010)Davenport, Harris, and
  Morison}]{Davenport2010-vu_JEB}
Davenport, T.~H., Harris, J.~G., Morison, R., 2010. Analytics at Work: Smarter
  Decisions, Better Results. Harvard Business Press.

\bibitem[{Davies et~al.(2005)Davies, Mabin, and
  Balderstone}]{davies_theory_2005_JHLS}
Davies, J., Mabin, V.~J., Balderstone, S.~J., 2005. The theory of constraints:
  a methodology apart?—a comparison with selected {OR}/{MS} methodologies.
  Omega 33~(6), 506--524.

\bibitem[{Davis et~al.(1993)Davis, Panas, and Zariphopoulou}]{davis93_MB}
Davis, M.~H., Panas, V.~G., Zariphopoulou, T., 1993. European option pricing
  with transaction costs. SIAM Journal on Control and Optimization 31~(2),
  470--493.

\bibitem[{Davis and Bracken(2022)}]{Davis2022-kx_KVRH}
Davis, P.~K., Bracken, P., 2022. Artificial intelligence for wargaming and
  modeling. The Journal of Defense Modeling and Simulation, DOI:
  10.1177/15485129211073126.

\bibitem[{Davtalab-Olyaie et~al.(2023)Davtalab-Olyaie, Mahmudi-Baram, and
  Asgharian}]{Davtalab-Olyaie2022-yb_SL}
Davtalab-Olyaie, M., Mahmudi-Baram, H., Asgharian, M., 2023. Measuring
  individual efficiency and unit influence in centrally managed systems. Annals
  of Operations Research 321, 139--164.

\bibitem[{Davydenko and Fildes(2013)}]{Davydenko2013_FP}
Davydenko, A., Fildes, R., 2013. Measuring forecasting accuracy: The case of
  judgmental adjustments to {SKU}-level demand forecasts. International Journal
  of Forecasting 29~(3), 510--522.

\bibitem[{Dayarian and Savelsbergh(2020)}]{dayarian2020crowdshipping_CKTVW}
Dayarian, I., Savelsbergh, M., 2020. Crowdshipping and same-day delivery:
  Employing in-store customers to deliver online orders. Production and
  Operations Management 29~(9), 2153--2174.

\bibitem[{Dayarian et~al.(2020)Dayarian, Savelsbergh, and
  Clarke}]{dayarian2020same_VLVV}
Dayarian, I., Savelsbergh, M., Clarke, J.-P., 2020. Same-day delivery with
  drone resupply. Transportation Science 54~(1), 229--249.

\bibitem[{De~Baets and Harvey(2020)}]{De_Baets2020-du_FP}
De~Baets, S., Harvey, N., 2020. Using judgment to select and adjust forecasts
  from statistical models. European Journal of Operational Research 284~(3),
  882--895.

\bibitem[{de~Brito and Evers(2016)}]{Madruga_de_Brito2016-ul_JL}
de~Brito, M.~M., Evers, M., 2016. Multi-criteria decision-making for flood risk
  management: a survey of the current state of the art. Natural Hazards and
  Earth System Sciences 16, 1019--1033.

\bibitem[{{De Causmaecker} and {Vanden Berghe}(2011)}]{categorisation_GVBSP}
{De Causmaecker}, P., {Vanden Berghe}, G., 2011. A categorisation of nurse
  rostering problems. Journal of Scheduling 14, 3--16.

\bibitem[{de~Figueiredo and Mayerle(2008)}]{De_Figueiredo2008-rt_AM}
de~Figueiredo, J.~N., Mayerle, S.~F., 2008. Designing minimum-cost recycling
  collection networks with required throughput. Transportation Research Part E:
  Logistics and Transportation Review 44~(5), 731--752.

\bibitem[{de~Finetti(1937)}]{de1937prevision_MESG}
de~Finetti, B., 1937. La pr{\'e}vision: ses lois logiques, ses sources
  subjectives. In: Annales de l'Institut Henri Poincar{\'e}. Vol.~7. pp. 1--68.

\bibitem[{de~Gooyert et~al.(2017)de~Gooyert, Rouwette, van Kranenburg, and
  Freeman}]{De_Gooyert2017-yt_JL}
de~Gooyert, V., Rouwette, E., van Kranenburg, H., Freeman, E., 2017. Reviewing
  the role of stakeholders in operational research: A stakeholder theory
  perspective. European Journal of Operational Research 262~(2), 402--410.

\bibitem[{de~Jomini(1862)}]{Jomini1862_JLYHK}
de~Jomini, A.-H., 1862. The Art of War. J.B. Lippincott, Philadelphia.

\bibitem[{{de Keizer} et~al.(2017){de Keizer}, Akkerman, Grunow, Bloemhof,
  Haijema, and {van der Vorst}}]{DeKeizer2017_JLYHK}
{de Keizer}, M., Akkerman, R., Grunow, M., Bloemhof, J.~M., Haijema, R., {van
  der Vorst}, J. G. A.~J., 2017. Logistics network design for perishable
  products with heterogeneous quality decay. European Journal of Operational
  Research 262~(2), 535--549.

\bibitem[{de~Kok and Graves(2003)}]{de2003supply_JSS}
de~Kok, A.~d., Graves, S.~C., 2003. Supply chain management: Design,
  coordination and operation. Elsevier.

\bibitem[{{de Oliveira} et~al.(2015){de Oliveira}, {de Souza}, and
  Yunes}]{DEOLIVEIRA2015101_GVBSP}
{de Oliveira}, L., {de Souza}, C.~C., Yunes, T., 2015. On the complexity of the
  traveling umpire problem. Theoretical Computer Science 562, 101--111.

\bibitem[{De~Reyck and Herroelen(1998)}]{De_Reyck1998-sh_WH_ED}
De~Reyck, B., Herroelen, W., 1998. A branch-and-bound procedure for the
  resource-constrained project scheduling problem with generalized precedence
  relations. European Journal of Operational Research 111~(1), 152--174.

\bibitem[{De~Reyck and Herroelen(1999)}]{De_Reyck1999-oz_WH_ED}
De~Reyck, B., Herroelen, W., 1999. The multi-mode resource-constrained project
  scheduling problem with generalized precedence relations. European Journal of
  Operational Research 119~(2), 538--556.

\bibitem[{de~Treville and Antonakis(2006)}]{Treville2006-ak_KESXZ}
de~Treville, S., Antonakis, J., 2006. Could lean production job design be
  intrinsically motivating? {C}ontextual, configurational, and
  levels-of-analysis issues. Journal of Operations Management 24~(2), 99--123.

\bibitem[{{de Werra} et~al.(2002){de Werra}, Asratian, and
  Durand}]{DeWerra_GVBSP}
{de Werra}, D., Asratian, A.~S., Durand, S., 2002. Complexity of some special
  types of timetabling problems. Journal of Scheduling 5, 171--183.

\bibitem[{De~Witte and L{\'o}pez-Torres(2017)}]{Witte2017-xg_JJ}
De~Witte, K., L{\'o}pez-Torres, L., 2017. Efficiency in education: a review of
  literature and a way forward. Journal of the Operational Research Society
  68~(4), 339--363.

\bibitem[{Deb(2001)}]{kdeb01_MESG}
Deb, K., 2001. Multi-Objective Optimization using Evolutionary Algorithms. John
  Wiley \& Sons, Hoboken, NJ.

\bibitem[{Debels et~al.(2006)Debels, De~Reyck, Leus, and
  Vanhoucke}]{Debels2006-di_WH_ED}
Debels, D., De~Reyck, B., Leus, R., Vanhoucke, M., 2006. A hybrid scatter
  search/electromagnetism meta-heuristic for project scheduling. European
  Journal of Operational Research 169~(2), 638--653.

\bibitem[{Debo et~al.(2004)Debo, Savaskan, and Van~Wassenhove}]{Debo2004-xh_AM}
Debo, L.~G., Savaskan, R.~C., Van~Wassenhove, L.~N., 2004. Coordination in
  {Closed-Loop} supply chains. In: Dekker, R., Fleischmann, M., Inderfurth, K.,
  Van~Wassenhove, L.~N. (Eds.), Reverse Logistics: Quantitative Models for
  {Closed-Loop} Supply Chains. Springer, Berlin, pp. 295--311.

\bibitem[{Debo et~al.(2005)Debo, Toktay, and Van~Wassenhove}]{Debo2005-vu_AM}
Debo, L.~G., Toktay, L.~B., Van~Wassenhove, L.~N., 2005. Market segmentation
  and product technology selection for remanufacturable products. Management
  Science 51~(8), 1193--1205.

\bibitem[{Dejonckheere et~al.(2003)Dejonckheere, Disney, Lambrecht, and
  Towill}]{Dejonckheere2003_SMD}
Dejonckheere, J., Disney, S.~M., Lambrecht, M.~R., Towill, D.~R., 2003.
  Measuring and avoiding the bullwhip effect: A control theoretic approach.
  European Journal of Operational Research 147~(3), 567--590.

\bibitem[{Dejonckheere et~al.(2004)Dejonckheere, Disney, Lambrecht, and
  Towill}]{Dejonckheere2004_SMD}
Dejonckheere, J., Disney, S.~M., Lambrecht, M.~R., Towill, D.~R., 2004. The
  impact of information enrichment on the bullwhip effect in supply chains: A
  control engineering perspective. European Journal of Operational Research
  153~(3), 727 -- 750.

\bibitem[{Delen and Ram(2018)}]{Delen2018-hf_JEB}
Delen, D., Ram, S., 2018. Research challenges and opportunities in business
  analytics. Journal of Business Analytics 1~(1), 2--12.

\bibitem[{Della~Croce et~al.(2014)Della~Croce, Grosso, and
  Salassa}]{della2014matheuristic_COIT}
Della~Croce, F., Grosso, A., Salassa, F., 2014. A matheuristic approach for the
  two-machine total completion time flow shop problem. Annals of Operations
  Research 213~(1), 67--78.

\bibitem[{Dellnitz(2022)}]{Dellnitz2022-ei_SL}
Dellnitz, A., 2022. Big data efficiency analysis: Improved algorithms for data
  envelopment analysis involving large datasets. Computers \& Operations
  Research 137, 105553.

\bibitem[{Demassey(2008)}]{Demassey2010-pt_WH_ED}
Demassey, S., 2008. Mathematical programming formulations and lower bounds. In:
  Artigues, C., Demassey, S., N{\'e}ron, E. (Eds.), {Resource-Constrained}
  Project Scheduling -- Models, algorithms, extensions and applications. ISTE,
  London, pp. 49--62.

\bibitem[{Demeulemeester(1995)}]{Demeulemeester1995-gy_WH_ED}
Demeulemeester, E., 1995. Minimizing resource availability costs in
  time-limited project networks. Management Science 41~(10), 1590--1598.

\bibitem[{Demeulemeester and Herroelen(1992)}]{Demeulemeester1992-xs_WH_ED}
Demeulemeester, E., Herroelen, W., 1992. A branch-and-bound procedure for the
  multiple resource-constrained project scheduling problem. Management Science
  38~(12), 1803--1818.

\bibitem[{Demeulemeester and Herroelen(2011)}]{Demeulemeester2011-eb_WH_ED}
Demeulemeester, E., Herroelen, W., 2011. Robust project scheduling. Foundations
  and Trends\textregistered{} in Technology, Information and Operations
  Management 3~(3--4), 201--376.

\bibitem[{Demeulemeester and Herroelen(2002)}]{Demeulemeester2002-an_WH_ED}
Demeulemeester, E.~L., Herroelen, W., 2002. Project Scheduling: A Research
  Handbook. Springer Science \& Business Media.

\bibitem[{Demeulemeester and Herroelen(1996)}]{Demeulemeester1996-yl_WH_ED}
Demeulemeester, E.~L., Herroelen, W.~S., 1996. An efficient optimal solution
  procedure for the preemptive resource-constrained project scheduling problem.
  European Journal of Operational Research 90~(2), 334--348.

\bibitem[{Demeulemeester and Herroelen(1997)}]{Demeulemeester1997-bt_WH_ED}
Demeulemeester, E.~L., Herroelen, W.~S., 1997. A {Branch-and-Bound} procedure
  for the generalized resource-constrained project scheduling problem.
  Operations Research 45~(2), 201--212.

\bibitem[{den Boer et~al.(2020)den Boer, Lambrechts, and
  Krikke}]{Boer2020_JLYHK}
den Boer, J., Lambrechts, W., Krikke, H., 2020. Additive manufacturing in
  military and humanitarian missions: Advantages and challenges in the spare
  parts supply chain. Journal of Cleaner Production 257, 120301.

\bibitem[{{Department for Transport}(2022)}]{Department_for_Transport2018-vw}
{Department for Transport}, 2022. Road goods vehicles travelling to europe
  ({RORO}).
  \url{https://www.gov.uk/government/statistical-data-sets/road-goods-vehicles-travelling-to-europe},
  accessed on 2023-01-12.

\bibitem[{Desaulniers et~al.(1997)Desaulniers, Desrosiers, Dumas, Solomon, and
  Soumis}]{desaulniers1997daily_VLVV}
Desaulniers, G., Desrosiers, J., Dumas, Y., Solomon, M.~M., Soumis, F., 1997.
  Daily aircraft routing and scheduling. Management Science 43~(6), 841--855.

\bibitem[{Desaulniers et~al.(2006)Desaulniers, Desrosiers, and
  Solomon}]{DDS06_ALAL}
Desaulniers, G., Desrosiers, J., Solomon, M. (Eds.), 2006. Column Generation.
  Springer, New York.

\bibitem[{Desaulniers et~al.(2002)Desaulniers, Desrosiers, and
  Solomon}]{desaulniers2002accelerating_CA_MB}
Desaulniers, G., Desrosiers, J., Solomon, M.~M., 2002. Accelerating strategies
  in column generation methods for vehicle routing and crew scheduling
  problems. In: Ribeiro, C.~C., Hansen, P. (Eds.), Essays and surveys in
  metaheuristics. Springer, pp. 309--324.

\bibitem[{Desrosiers and L{\"u}bbecke(2005)}]{Desrosiers2005_CA_MB}
Desrosiers, J., L{\"u}bbecke, M., 2005. A primer in column generation. In:
  Desaulniers, G., Desrosiers, J., Solomon, M. (Eds.), Column Generation.
  Springer US, Boston, MA, pp. 1--32.

\bibitem[{Dettmer(1997)}]{dettmer_goldratts_1997_JHLS}
Dettmer, H.~W., 1997. Goldratt's {Theory} of {Constraints}: {A} {Systems}
  {Approach} to {Continuous} {Improvement}. ASQC Quality Press, Milwaukee, WI.

\bibitem[{Deutsch et~al.(2022)Deutsch, Lustfield, and
  Jalali}]{Deutsch2022-ns_AJG}
Deutsch, A.~R., Lustfield, R., Jalali, M.~S., 2022. Community-based system
  dynamics modelling of stigmatized public health issues: Increasing diverse
  representation of individuals with personal experiences. Systems Research and
  Behavioral Science 39~(4), 734--749.

\bibitem[{DeValve et~al.(2023)DeValve, Song, and Wei}]{de2023assemble_JSS}
DeValve, L., Song, J.-S., Wei, Y., 2023. Assemble-to-order systems. In: Song,
  J.-S. (Ed.), Research Handbook on Inventory Management. Edward Elgar
  Publishing.

\bibitem[{DHL(2013)}]{DHL_CKTVW}
DHL, 2013. Logistics trend radar. DHL.
\newline\urlprefix\url{https://www.dhl.com/content/dam/Campaigns/InnovationDay_2013/90310673_HI-RES.PDF}

\bibitem[{Diaz-Balteiro et~al.(2017)Diaz-Balteiro, Gonz{\'a}lez-Pach{\'o}n, and
  Romero}]{Diaz-Balteiro2017-tl_JL}
Diaz-Balteiro, L., Gonz{\'a}lez-Pach{\'o}n, J., Romero, C., 2017. Measuring
  systems sustainability with multi-criteria methods: A critical review.
  European Journal of Operational Research 258~(2), 607--616.

\bibitem[{Dieterich et~al.(2016)Dieterich, Mendoza, and
  Brennan}]{dieterich2016compas_JH}
Dieterich, W., Mendoza, C., Brennan, T., 2016. {COMPAS} risk scales:
  Demonstrating accuracy equity and predictive parity. Tech. rep., Northpointe
  Inc.\ Research Department, 8 July.

\bibitem[{Dietrich et~al.(1993)Dietrich, Escudero, and Chance}]{DEC93_ALAL}
Dietrich, B., Escudero, L., Chance, F., 1993. Efficient reformulation for 0-1
  programs---methods and computational results. Discrete Applied Mathematics
  42, 147--175.

\bibitem[{Dignum(2019)}]{Dignum19_LCAL}
Dignum, V., 2019. Responsible Artificial Intelligence: How to Develop and Use
  AI in a Responsible Way, 1st Edition. Springer, New York NY.

\bibitem[{Dijkstra(1959)}]{Dijkstra:1959_IL}
Dijkstra, E.~W., 1959. A note on two problems in connexion with graphs.
  Numerische Mathematik 1, 269--271.

\bibitem[{Dillenburger et~al.(2019)Dillenburger, Jordan, and
  Cochran}]{Dillenburger2019-ir_KVRH}
Dillenburger, S.~P., Jordan, J.~D., Cochran, J.~K., 2019. Pareto-optimality for
  lethality and collateral risk in the airstrike multi-objective problem.
  Journal of the Operational Research Society 70~(7), 1051--1064.

\bibitem[{DIMACS(2021)}]{DIMACS_CA_MB}
DIMACS, 2021. Implementation challenge: Vehicle routing.
  \url{http://dimacs.rutgers.edu/programs/challenge/vrp/}, accessed on
  2021-09-14.

\bibitem[{Dimopoulou and Miliotis(2001)}]{Dimopoulou2001-rz_JJ}
Dimopoulou, M., Miliotis, P., 2001. Implementation of a university course and
  examination timetabling system. European Journal of Operational Research
  130~(1), 202--213.

\bibitem[{Ding et~al.(2008)Ding, Glazebrook, and
  Kirkbride}]{ding2008allocation_DLLD}
Ding, L., Glazebrook, K.~D., Kirkbride, C., 2008. Allocation models and
  heuristics for the outsourcing of repairs for a dynamic warranty population.
  Management Science 54~(3), 594--607.

\bibitem[{Dinic(1970)}]{Dinic:1970_IL}
Dinic, E., 1970. An algorithm for the solution of the max-flow problem with the
  polynomial estimation. Doklady Akademii Nauk SSSR 194~(4), 1277--1280.

\bibitem[{Dinitz(2006)}]{Dinitz:2006_IL}
Dinitz, Y., 2006. Dinitz' algorithm: The original version and {E}ven's version.
  In: Goldreich, O., Rosenberg, A.~L., Selman, A.~L. (Eds.), Theoretical
  Computer Science, Essays in Memory of Shimon Even. Vol. 3895 of Lecture Notes
  in Computer Science. Springer, pp. 218--240.

\bibitem[{Dixit and Pindyck(2009)}]{dixit09_MB}
Dixit, A.~K., Pindyck, R.~S., 2009. The options approach to capital investment.
  In: Siesfeld, T., Cefola, J., Neef, D. (Eds.), The Economic Impact of
  Knowledge. Routledge, London, pp. 325--340.

\bibitem[{Dixon et~al.(2020)Dixon, Halperin, and Bilokon}]{dixon20_MB}
Dixon, M.~F., Halperin, I., Bilokon, P., 2020. Machine learning in Finance.
  Springer, Cham.

\bibitem[{Dixon and Coles(1997)}]{Dixon1997-cf_IM}
Dixon, M.~J., Coles, S.~G., 1997. Modelling association football scores and
  inefficiencies in the football betting market. Journal of the Royal
  Statistical Society. Series C, Applied statistics 46~(2), 265--280.

\bibitem[{Doganis et~al.(2008)Doganis, Aggelogiannaki, and
  Sarimveis}]{Doganis2008_XW}
Doganis, P., Aggelogiannaki, E., Sarimveis, H., 2008. A combined model
  predictive control and time series forecasting framework for
  production-inventory systems. International Journal of Production Research
  46~(24), 6841--6853.

\bibitem[{Dong et~al.(2016)Dong, Shi, and Zhang}]{Dong2016-aq_KESXZ}
Dong, L., Shi, D., Zhang, F., 2016. {3D} printing vs. traditional flexible
  technology: Implications for manufacturing strategy. Working Paper, St.
  Louis, MO: Washington University.

\bibitem[{Dong and Xu(2002)}]{Dong2002_SMD}
Dong, Y., Xu, K., 2002. A supply chain model of vendor managed inventory.
  Transportation Research Part E: Logistics and Transportation Review 38~(2),
  75--95.

\bibitem[{Donohue et~al.(2020)Donohue, {\"O}zer, and
  Zheng}]{Donohue2020-mw_AFRH}
Donohue, K., {\"O}zer, {\"O}., Zheng, Y., 2020. Behavioral operations: Past,
  present, and future. Manufacturing \& Service Operations Management 22~(1),
  191--202.

\bibitem[{Douglas and Wildavsky(1983)}]{Douglas1983-td_TC}
Douglas, M., Wildavsky, A., 1983. Risk and Culture: An Essay on the Selection
  of Technological and Environmental Dangers. University of California Press.

\bibitem[{Doukas and Nikas(2020)}]{Doukas2020-vl_JL}
Doukas, H., Nikas, A., 2020. Decision support models in climate policy.
  European Journal of Operational Research 280~(1), 1--24.

\bibitem[{Downey and Fellows(1999)}]{DF99_UPCT}
Downey, R.~G., Fellows, M.~R., 1999. Parameterized Complexity. Springer, New
  York, NY.

\bibitem[{Doyle and Green(1994)}]{Doyle1994-of_DEA}
Doyle, J., Green, R., 1994. Efficiency and {Cross-Efficiency} in {DEA}:
  Derivations, meanings and uses. The Journal of the Operational Research
  Society 45~(5), 567--578.

\bibitem[{Drexl and Knust(2007)}]{Drexl2007_GVBSP}
Drexl, A., Knust, S., 2007. {Sports league scheduling: Graph- and
  resource-based models}. Omega 35~(5), 465--471.

\bibitem[{Drexl(2012)}]{drexl2012synchronization_CKTVW}
Drexl, M., 2012. Synchronization in vehicle routing—a survey of vrps with
  multiple synchronization constraints. Transportation Science 46~(3),
  297--316.

\bibitem[{Dreyfus(1969)}]{dreyfus1969appraisal_DLLD}
Dreyfus, S.~E., 1969. An appraisal of some shortest-path algorithms. Operations
  Research 17~(3), 395--412.

\bibitem[{Dreyfus and Law(1977)}]{Dreyfus1977-sl_HL}
Dreyfus, S.~E., Law, A.~M., 1977. The Art and Theory of Dynamic Programming.
  Academic Press.

\bibitem[{Drezner and Hamacher(2004)}]{drezner2004facility_SAA}
Drezner, Z., Hamacher, H. (Eds.), 2004. Facility Location: Applications and
  Theory. Springer, Berlin-Heidelberg-New York.

\bibitem[{Duan and Xiong(2015)}]{Duan2015-bh_JEB}
Duan, L., Xiong, Y., 2015. Big data analytics and business analytics. Journal
  of Management Analytics 2~(1), 1--21.

\bibitem[{Dubitzky et~al.(2019)Dubitzky, Lopes, Davis, and
  Berrar}]{Dubitzky2019-kc_IM}
Dubitzky, W., Lopes, P., Davis, J., Berrar, D., 2019. The open international
  soccer database for machine learning. Machine Learning 108~(1), 9--28.

\bibitem[{Duckworth and Lewis(1998)}]{Duckworth1998-cu_IM}
Duckworth, F.~C., Lewis, A.~J., 1998. A fair method for resetting the target in
  interrupted one-day cricket matches. Journal of the Operational Research
  Society 49~(3), 220--227.

\bibitem[{Dunbar et~al.(2012)Dunbar, Froyland, and Wu}]{dunbar2012robust_VLVV}
Dunbar, M., Froyland, G., Wu, C.-L., 2012. Robust airline schedule planning:
  Minimizing propagated delay in an integrated routing and crewing framework.
  Transportation Science 46~(2), 204--216.

\bibitem[{Dunning et~al.(2017)Dunning, Huchette, and Lubin}]{Jump_CTCGE}
Dunning, I., Huchette, J., Lubin, M., 2017. {JuMP}: A modeling language for
  mathematical optimization. SIAM Review 59~(2), 295--320.

\bibitem[{Dur{\'a}n(2021)}]{duran2021sports_GVBSP}
Dur{\'a}n, G., 2021. Sports scheduling and other topics in sports analytics: a
  survey with special reference to {L}atin {A}merica. Top 29~(1), 125--155.

\bibitem[{Duran et~al.(2011)Duran, Gutierrez, and Keskinocak}]{Duran2011_JLYHK}
Duran, S., Gutierrez, M.~A., Keskinocak, P., 2011. Pre-positioning of emergency
  items for {CARE International}. Interfaces 41~(3), 223--237.

\bibitem[{Dutta(1995)}]{Du1995_GZ}
Dutta, P., 1995. A folk theorem for stochastic games. Journal of Economic
  Theory 66, 1--32.

\bibitem[{Dutton and Walton(1964)}]{Dutton1964-wf_AFRH}
Dutton, J.~M., Walton, R.~E., 1964. Operational research and the behavioural
  sciences. Operations Research Quarterly 15~(3), 207--217.

\bibitem[{Dwork et~al.(2012)Dwork, Hardt, Pitassi, Reingold, and
  Zemel}]{dwork2012fairness_JH}
Dwork, C., Hardt, M., Pitassi, T., Reingold, O., Zemel, R., 2012. Fairness
  through awareness. In: Proceedings of the 3rd Innovations in Theoretical
  Computer Science Conference. pp. 214--226.

\bibitem[{Dyckhoff(1990)}]{Dyckhoff1990-lr_JB}
Dyckhoff, H., 1990. A typology of cutting and packing problems. European
  Journal of Operational Research 44~(2), 145--159.

\bibitem[{Dyer(1990)}]{dyer1990remarks}
Dyer, J.~S., 1990. Remarks on the analytic hierarchy process. Management
  Science 36~(3), 249--258.

\bibitem[{Dyson(2000)}]{Dyson2000-gz_MYLW}
Dyson, R.~G., 2000. Strategy, performance and operational research. Journal of
  the Operational Research Society 51~(1), 5--11.

\bibitem[{Dyson et~al.(2021)Dyson, O'Brien, and Shah}]{Dyson2021-sb_MYLW}
Dyson, R.~G., O'Brien, F.~A., Shah, D.~B., 2021. Soft {OR} and practice: The
  contribution of the founders of operations research. Operations Research
  69~(3), 727--738.

\bibitem[{Dzidonu and Foster(1993)}]{Dzidonu1993_EA}
Dzidonu, C.~K., Foster, F.~G., 1993. Prolegomena to or modelling of the global
  environment-development problem. Journal of the Operational Research Society
  44~(4), 321--331.

\bibitem[{Dönmez et~al.(2021)Dönmez, Kara, Özlem Karsu, and
  da~Gama}]{DONMEZ20211_BYKOK}
Dönmez, Z., Kara, B.~Y., Özlem Karsu, da~Gama, F.~S., 2021. Humanitarian
  facility location under uncertainty: Critical review and future prospects.
  Omega 102, 102393.

\bibitem[{D’Ariano et~al.(2019)D’Ariano, Meng, Centulio, and
  Corman}]{Dariano2019_DC}
D’Ariano, A., Meng, L., Centulio, G., Corman, F., 2019. Integrated stochastic
  optimization approaches for tactical scheduling of trains and railway
  infrastructure maintenance. Computers \& Industrial Engineering 127,
  1315--1335.

\bibitem[{Easton et~al.(2001)Easton, Nemhauser, and Trick}]{Easton2001_GVBSP}
Easton, K., Nemhauser, G., Trick, M., 2001. The traveling tournament problem
  description and benchmarks. In: Walsh, T. (Ed.), Principles and Practice of
  Constraint Programming --- CP 2001: 7th International Conference, CP 2001
  Paphos, Cyprus, 2001 Proceedings. Springer, Berlin, Heidelberg, pp. 580--584.

\bibitem[{Eden(1989)}]{Eden1989-rz_MYLW}
Eden, C., 1989. Using cognitive mapping for strategic options development and
  analysis ({SODA}). In: Rosenhead, J. (Ed.), Rational analysis for a
  problematic world: problem structuring methods for complexity, uncertainty,
  and conflict. Wiley, Chichester, pp. 21--42.

\bibitem[{Eden and Ackermann(2004)}]{Eden2004-ws_MYLW}
Eden, C., Ackermann, F., 2004. Cognitive mapping expert views for policy
  analysis in the public sector. European Journal of Operational Research
  152~(3), 615--630.

\bibitem[{Eden and Ackermann(2006)}]{Eden2006-as_MYLW}
Eden, C., Ackermann, F., 2006. Where next for problem structuring methods.
  Journal of the Operational Research Society 57~(7), 766--768.

\bibitem[{Ederer(2015)}]{Ederer_2015_DSRW}
Ederer, N., 2015. Evaluating capital and operating cost efficiency of offshore
  wind farms: A {DEA} approach. Renewable and Sustainable Energy Reviews
  42~(1), 1034--1046.

\bibitem[{Edgeworth(1888)}]{edgeworth1888mathematical_JSS}
Edgeworth, F.~Y., 1888. The mathematical theory of banking. Journal of the
  Royal Statistical Society 51~(1), 113--127.

\bibitem[{Edmonds(1965{\natexlab{a}})}]{Edmonds1965-xd_JMB}
Edmonds, J., 1965{\natexlab{a}}. Maximum matching and a polyhedron with
  0,1-vertices. Journal of Research of the National Bureau of Standards 69B~(1
  and 2), 125.

\bibitem[{Edmonds(1965{\natexlab{b}})}]{Edmonds1965-zj_JMB}
Edmonds, J., 1965{\natexlab{b}}. Paths, trees, and flowers. Canadian Journal of
  Mathematics. Journal Canadien de Mathematiques 17, 449--467.

\bibitem[{Edmonds and Karp(1972)}]{Edmonds-Karp:1972_IL}
Edmonds, J., Karp, R.~M., 1972. Theoretical improvements in algorithmic
  efficiency for network flow problems. Journal of the {ACM} 19~(2), 248--264.

\bibitem[{Edwards and Barron(1994)}]{edwards1994smarts_MESG}
Edwards, W., Barron, F., 1994. {SMARTS} and {SMARTER}: Improved simple methods
  for multiattribute utility measurement. Organizational Behavior and Human
  Decision Processes 60~(3), 306--325.

\bibitem[{Edwards et~al.(2007)Edwards, Miles, and von
  Winterfeldt}]{edwards2007advances_MESG}
Edwards, W., Miles, R., von Winterfeldt, D., 2007. Advances in decision
  analysis: From foundations to applications. Cambridge University Press,
  Cambridge.

\bibitem[{Egerv\'ary(1931)}]{E31_SMPT}
Egerv\'ary, E., 1931. Matrixok kombinatorius tulajdons\'agairol. Matematikai
  \'es Fizikai Lapok 38, 16--28.

\bibitem[{Ehmke and Campbell(2014)}]{ehmke2014customer_CKTVW}
Ehmke, J.~F., Campbell, A.~M., 2014. Customer acceptance mechanisms for home
  deliveries in metropolitan areas. European Journal of Operational Research
  233~(1), 193--207.

\bibitem[{Ehrgott(2005)}]{ehrgott05_MESG}
Ehrgott, M., 2005. Multicriteria Optimization, 2nd Edition. Springer Verlag,
  Berlin.

\bibitem[{Ehrhardt(1984)}]{ehrhardt1984s_JSS}
Ehrhardt, R., 1984. (s, s) policies for a dynamic inventory model with
  stochastic lead times. Operations Research 32~(1), 121--132.

\bibitem[{Eiselt and Marianov(2011)}]{eiselt2011foundations_SAA}
Eiselt, H., Marianov, V. (Eds.), 2011. Foundations of Location Analysis.
  International Series in Operations Research \& Management Science. Springer,
  Berlin-Heidelberg-New York.

\bibitem[{Eiselt and Marianov(2015)}]{eiselt2015applications_SAA}
Eiselt, H., Marianov, V. (Eds.), 2015. Applications of Location Analysis.
  International Series in Operations Research \& Management Science. Springer,
  Berlin-Heidelberg-New York.

\bibitem[{Eisenbrand et~al.(2010)Eisenbrand, Grandoni, Rothvoß, and
  Schäfer}]{eisenbrand.grandoni.ea:10_BF}
Eisenbrand, F., Grandoni, F., Rothvoß, T., Schäfer, G., 2010. Connected
  facility location via random facility sampling and core detouring. Journal of
  Computer and System Sciences 76, 709 -- 726.

\bibitem[{El-Taha and Stidham(1999)}]{eltaha_HATI}
El-Taha, M., Stidham, Jr, S., 1999. Sample-path analysis of queueing systems.
  Kluwer Academic Publishers, Boston, MA.

\bibitem[{El{\c{c}}i et~al.(2022)El{\c{c}}i, Hooker, and
  Zhang}]{ElcHooZha22_JH}
El{\c{c}}i, {\"{O}}., Hooker, J.~N., Zhang, P., 2022. Structural properties of
  equitable and efficient distributions. Tech. rep., Carnegie Mellon
  University, submitted.

\bibitem[{Eliashberg and Winkler(1981)}]{Eliashberg1981-dt_TC}
Eliashberg, J., Winkler, R.~L., 1981. Risk sharing and group decision making.
  Management Science 27~(11), 1221--1235.

\bibitem[{Eme{\c{c}} et~al.(2016)Eme{\c{c}}, {\c{C}}atay, and
  Bozkaya}]{emec2016adaptive_CKTVW}
Eme{\c{c}}, U., {\c{C}}atay, B., Bozkaya, B., 2016. An adaptive large
  neighborhood search for an e-grocery delivery routing problem. Computers \&
  Operations Research 69, 109--125.

\bibitem[{Emmeche et~al.(1997)Emmeche, K{\o}ppe, and
  Stjernfelt}]{Emmeche1997-cj_GM}
Emmeche, C., K{\o}ppe, S., Stjernfelt, F., 1997. Explaining emergence: Towards
  an ontology of levels. Journal for General Philosophy of Science 28~(1),
  83--117.

\bibitem[{Engwerda(2005)}]{En2005_GZ}
Engwerda, J., 2005. {LQ} dynamic optimization and differential games. John
  Wiley \& Sons.

\bibitem[{Eppen(1979)}]{eppen1979note_JSS}
Eppen, G.~D., 1979. Note—effects of centralization on expected costs in a
  multi-location newsboy problem. Management Science 25~(5), 498--501.

\bibitem[{Eppler and Aeschimann(2009)}]{Eppler2009-qy_MJE}
Eppler, M.~J., Aeschimann, M., 2009. A systematic framework for risk
  visualization in risk management and communication. Risk Management: An
  International Journal 11~(2), 67--89.

\bibitem[{Eppler and Burkhard(2008)}]{Eppler2008-gh_MJE}
Eppler, M.~J., Burkhard, R., 2008. Knowledge visualization. In: Jennex, E.
  (Ed.), Knowledge Management: Concepts, Methodologies, Tools, and
  Applications. IGI Global, pp. 987--999.

\bibitem[{Eppler and Kernbach(2016)}]{Eppler2016-tb_MJE}
Eppler, M.~J., Kernbach, S., 2016. Dynagrams: Enhancing design thinking through
  dynamic diagrams. Design Studies 47, 91--117.

\bibitem[{Eppler and Platts(2009)}]{Eppler2009-qw_MJE}
Eppler, M.~J., Platts, K.~W., 2009. Visual strategizing: The systematic use of
  visualization in the {Strategic-Planning} process. Long Range Planning
  42~(1), 42--74.

\bibitem[{Epure et~al.(2011)Epure, Kerstens, and Prior}]{Epure2011-yw_SL}
Epure, M., Kerstens, K., Prior, D., 2011. Bank productivity and performance
  groups: A decomposition approach based upon the {L}uenberger productivity
  indicator. European Journal of Operational Research 211~(3), 630--641.

\bibitem[{Erdo{\u g}an et~al.(2019)Erdo{\u g}an, Stylianou, and
  Vasilakis}]{Erdogan2019-ms_CV}
Erdo{\u g}an, G., Stylianou, N., Vasilakis, C., 2019. An open source decision
  support system for facility location analysis. Decision Support Systems 125,
  113116.

\bibitem[{Erdo\u{g}an(2017)}]{Erdogan2017_CTCGE}
Erdo\u{g}an, G., 2017. An open source spreadsheet solver for vehicle routing
  problems. Computers \& Operations Research 84, 62--72.

\bibitem[{Erera et~al.(2013)Erera, Hewitt, Savelsbergh, and
  Zhang}]{erera2013improved_MH}
Erera, A., Hewitt, M., Savelsbergh, M., Zhang, Y., 2013. Improved load plan
  design through integer programming based local search. Transportation Science
  47~(3), 412--427.

\bibitem[{Erkut and Neuman(1989)}]{ERKUT1989275_SAA}
Erkut, E., Neuman, S., 1989. Analytical models for locating undesirable
  facilities. European Journal of Operational Research 40~(3), 275--291.

\bibitem[{Erlang(1909)}]{Erlang1090-er_GL}
Erlang, A.~K., 1909. The theory of probabilities and telephone conversations.
  Nyt Tidsskrift for Matematik B 20, 33--39.

\bibitem[{Ernst et~al.(2004)Ernst, Jiang, Krishnamoorthy, and
  Sier}]{Ernst_GVBSP}
Ernst, A.~T., Jiang, H., Krishnamoorthy, M., Sier, D., 2004. Staff scheduling
  and rostering: A review of applications, methods and models. European Journal
  of Operational Research 153, 3--27.

\bibitem[{Esmail and Geneletti(2018)}]{Esmail2018-wp_JL}
Esmail, B.~A., Geneletti, D., 2018. Multi-criteria decision analysis for nature
  conservation: A review of 20 years of applications. Methods in Ecology and
  Evolution 9~(1), 42--53.

\bibitem[{Espejo and Harnden(1989)}]{Espejo1989-fd_GM}
Espejo, R., Harnden, R., 1989. The Viable System Model: Interpretations and
  Applications of Stafford Beer's {VSM}. Wiley.

\bibitem[{Essid et~al.(2010)Essid, Ouellette, and Vigeant}]{Essid2010-at_JJ}
Essid, H., Ouellette, P., Vigeant, S., 2010. Measuring efficiency of {T}unisian
  schools in the presence of quasi-fixed inputs: A bootstrap data envelopment
  analysis approach. Economics of Education Review 29~(4), 589--596.

\bibitem[{Esteve et~al.(2020)Esteve, Aparicio, Rabasa, and
  Rodriguez-Sala}]{Esteve2020-ud_SL}
Esteve, M., Aparicio, J., Rabasa, A., Rodriguez-Sala, J.~J., 2020. Efficiency
  analysis trees: A new methodology for estimating production frontiers through
  decision trees. Expert Systems with Applications 162, 113783.

\bibitem[{{EURO Meets NeurIPS 2022}(2022)}]{EURO_CA_MB}
{EURO Meets NeurIPS 2022}, 2022. Vehicle routing competition.
  \url{http://www.verolog.eu/}, accessed on 2021-09-14.

\bibitem[{Eveborn et~al.(2009)Eveborn, Rönnqvist, Einarsdóttir, Eklund,
  Lidén, and Almroth}]{eveborn_operations_2009_JHLS}
Eveborn, P., Rönnqvist, M., Einarsdóttir, H., Eklund, M., Lidén, K.,
  Almroth, M., 2009. Operations {Research} {Improves} {Quality} and
  {Efficiency} in {Home} {Care}. Interfaces 39~(1), 18--34.

\bibitem[{{EWG PATAT}(1996)}]{EWG_PATAT_GVBSP}
{EWG PATAT}, 1996. {EURO} working group on automated timetabling.
  \url{https://patat.cs.kuleuven.be}, accessed on 2022-10-09.

\bibitem[{Ewing et~al.(2006)Ewing, Tarantino, and Parnell}]{Ewing2006-as_KVRH}
Ewing, P.~L., Tarantino, W., Parnell, G.~S., 2006. Use of decision analysis in
  the army base realignment and closure ({BRAC}) 2005 military value analysis.
  Decision Analysis 3~(1), 33--49.

\bibitem[{{Facsimile Simulation Library}(2021)}]{Facsimile_CTCGE}
{Facsimile Simulation Library}, 2021. {Facsimile}.
\newline\urlprefix\url{https://github.com/facsimile/facsimile}

\bibitem[{Fagerholt et~al.(2010)Fagerholt, Laporte, and
  Norstad}]{Fagerholt2010-lr_HP}
Fagerholt, K., Laporte, G., Norstad, I., 2010. Reducing fuel emissions by
  optimizing speed on shipping routes. Journal of the Operational Research
  Society 61~(3), 523--529.

\bibitem[{Fagerholt and Ronen(2013)}]{Fagerholt2013-on_HP}
Fagerholt, K., Ronen, D., 2013. Bulk ship routing and scheduling: solving
  practical problems may provide better results. Maritime Policy \& Management
  40~(1), 48--64.

\bibitem[{Fairbrother et~al.(2020)Fairbrother, Zografos, and
  Glazebrook}]{fairbrother2020slot_VLVV}
Fairbrother, J., Zografos, K.~G., Glazebrook, K.~D., 2020. A slot-scheduling
  mechanism at congested airports that incorporates efficiency, fairness, and
  airline preferences. Transportation Science 54~(1), 115--138.

\bibitem[{Fairley et~al.(2019)Fairley, Scheinker, and
  Brandeau}]{Fairley2019-eu_CV}
Fairley, M., Scheinker, D., Brandeau, M.~L., 2019. Improving the efficiency of
  the operating room environment with an optimization and machine learning
  model. Health Care Management Science 22~(4), 756--767.

\bibitem[{Fan and Sundaresan(2000)}]{fan00_MB}
Fan, H., Sundaresan, S.~M., 2000. Debt valuation, renegotiation, and optimal
  dividend policy. The Review of Financial Studies 13~(4), 1057--1099.

\bibitem[{Fang and Puthenpura(1993)}]{Fang1993-bg_JMB}
Fang, S.-C., Puthenpura, S., 1993. Linear Optimization and Extensions: Theory
  and Algorithms. Prentice Hall, Englewood Cliffs.

\bibitem[{Farahani et~al.(2013)Farahani, Miandoabchi, Szeto, and
  Rashidi}]{farahani2013review_MH}
Farahani, R.~Z., Miandoabchi, E., Szeto, W.~Y., Rashidi, H., 2013. A review of
  urban transportation network design problems. European Journal of Operational
  Research 229~(2), 281--302.

\bibitem[{Farahani et~al.(2023)Farahani, Ruiz, and
  Van~Wassenhove}]{FARAHANI20231_BKOK}
Farahani, R.~Z., Ruiz, R., Van~Wassenhove, L.~N., 2023. Introduction to the
  special issue on the role of operational research in future epidemics/
  pandemics. European Journal of Operational Research 304~(1), 1--8, the role
  of Operational Research in future epidemics/ pandemics.
\newline\urlprefix\url{https://www.sciencedirect.com/science/article/pii/S0377221722005720}

\bibitem[{F{\"a}re and Grosskopf(2000)}]{Fare2000-wo_JJ}
F{\"a}re, R., Grosskopf, S., 2000. Network {DEA}. Socio-Economic Planning
  Sciences 34~(1), 35--49.

\bibitem[{F{\"a}re et~al.(1992)F{\"a}re, Grosskopf, Lindgren, and
  Roos}]{Fare1992-gf_SL}
F{\"a}re, R., Grosskopf, S., Lindgren, B., Roos, P., 1992. Productivity changes
  in {S}wedish pharamacies 1980--1989: A non-parametric {M}almquist approach.
  Journal of Productivity Analysis 3~(1), 85--101.

\bibitem[{Fasolo and Bana~e Costa(2014)}]{Fasolo2014-mp_AFRH}
Fasolo, B., Bana~e Costa, C.~A., 2014. Tailoring value elicitation to decision
  makers' numeracy and fluency: Expressing value judgments in numbers or words.
  Omega 44, 83--90.

\bibitem[{Fatnassi et~al.(2015)Fatnassi, Chaouachi, and
  Klibi}]{fatnassi2015planning_JLYHK}
Fatnassi, E., Chaouachi, J., Klibi, W., 2015. Planning and operating a shared
  goods and passengers on-demand rapid transit system for sustainable
  city-logistics. Transportation Research Part B: Methodological 81, 440--460.

\bibitem[{Fauske and Hoff(2016)}]{Fauske2016-fi_KVRH}
Fauske, M.~F., Hoff, E.~{\O}., 2016. From {F-16} to {F-35}: Optimizing the
  training of pilots in the {Royal Norwegian Air Force}. Interfaces 46~(4),
  326--333.

\bibitem[{Federgruen(1993)}]{federgruen1993centralized_JSS}
Federgruen, A., 1993. Centralized planning models for multi-echelon inventory
  systems under uncertainty. Handbooks in Operations Research and Management
  Science 4, 133--173.

\bibitem[{Federgruen et~al.(1986)Federgruen, Prastacos, and
  Zipkin}]{Federgruen1986_JLYHK}
Federgruen, A., Prastacos, G., Zipkin, P.~H., 1986. {An Allocation and
  Distribution Model for Perishable Products}. Operations Research 34~(1),
  75--82.

\bibitem[{Federgruen and
  Zipkin(1984{\natexlab{a}})}]{federgruen1984combined_JLYHK}
Federgruen, A., Zipkin, P., 1984{\natexlab{a}}. A combined vehicle routing and
  inventory allocation problem. Operations Research 32~(5), 1019--1037.

\bibitem[{Federgruen and
  Zipkin(1984{\natexlab{b}})}]{federgruen1984computational_JSS}
Federgruen, A., Zipkin, P., 1984{\natexlab{b}}. Computational issues in an
  infinite-horizon, multiechelon inventory model. Operations Research 32~(4),
  818--836.

\bibitem[{Federgruen and
  Zipkin(1986{\natexlab{a}})}]{federgruen1986ainventory_JSS}
Federgruen, A., Zipkin, P., 1986{\natexlab{a}}. An inventory model with limited
  production capacity and uncertain demands {I.} the average-cost criterion.
  Mathematics of Operations Research 11~(2), 193--207.

\bibitem[{Federgruen and
  Zipkin(1986{\natexlab{b}})}]{federgruen1986binventory_JSS}
Federgruen, A., Zipkin, P., 1986{\natexlab{b}}. An inventory model with limited
  production capacity and uncertain demands {II}. the discounted-cost
  criterion. Mathematics of Operations Research 11~(2), 208--215.

\bibitem[{Feitzinger and Lee(1997)}]{Feitzinger1997_SMD}
Feitzinger, E., Lee, H.~L., 1997. Mass customization at {H}ewlett-{P}ackard:
  {T}he power of postponement. Harvard Business Review 75~(1), 116--121.

\bibitem[{Feller(1957{\natexlab{a}})}]{feller1_HATI}
Feller, W., 1957{\natexlab{a}}. An Introduction to Probability Theory and Its
  Applications, Volume I. Wiley, New York.

\bibitem[{Feller(1957{\natexlab{b}})}]{feller2_HATI}
Feller, W., 1957{\natexlab{b}}. An Introduction to Probability Theory and Its
  Applications, Volume II. Wiley, New York.

\bibitem[{Feng and Gallego(1995)}]{feng1995optimal_VLVV}
Feng, Y., Gallego, G., 1995. Optimal starting times for end-of-season sales and
  optimal stopping times for promotional fares. Management Science 41~(8),
  1371--1391.

\bibitem[{Fenger et~al.(1991)Fenger, Halsnaes, Larsen, Schroll, and
  Vidal}]{Fenger1991_EA}
Fenger, J., Halsnaes, K., Larsen, H., Schroll, H., Vidal, V., 1991.
  Environment, Energy and Natural Resources Management in the Baltic Sea
  Region. Nordic Council of Ministers, Copenhagen.

\bibitem[{Ferguson et~al.(2009)Ferguson, Guide, Koca, and
  Souza}]{Ferguson2009-vt_AM}
Ferguson, M., Guide, Jr, V.~D., Koca, E., Souza, G.~C., 2009. The value of
  quality grading in remanufacturing. Production and Operations Management
  18~(3), 300--314.

\bibitem[{Ferguson and Souza(2010)}]{Ferguson2010-pb_AM}
Ferguson, M.~E., Souza, G.~C., 2010. {Closed-Loop} Supply Chains: New
  Developments to Improve the Sustainability of Business Practices. CRC Press,
  Boca Raton, FL.

\bibitem[{Ferguson et~al.(2020)Ferguson, Laydon, Nedjati~Gilani, Ainslie,
  Baguelin, Bhatia, Boonyasiri, Cucunuba~Perez, Cuomo-Dannenburg, Dighe,
  Dorigatti, Fu, Gaythorpe, Green, Hamlet, Hinsley, Okell, Van~Elsland,
  Thompson, Verity, Volz, Wang, Wang, Walker, Walters, Winskill, Whittaker,
  Donnelly, Riley, and Ghani}]{Ferguson2020_CC}
Ferguson, N., Laydon, D., Nedjati~Gilani, G~Imai, N., Ainslie, K., Baguelin,
  M., Bhatia, S., Boonyasiri, A., Cucunuba~Perez, Z., Cuomo-Dannenburg, G.,
  Dighe, A., Dorigatti, I., Fu, H., Gaythorpe, K., Green, W., Hamlet, A.,
  Hinsley, W., Okell, L., Van~Elsland, S., Thompson, H., Verity, R., Volz, E.,
  Wang, H., Wang, Y., Walker, P., Walters, C., Winskill, P., Whittaker, C.,
  Donnelly, C., Riley, S., Ghani, A., 2020. Report 9: Impact of
  non-pharmaceutical interventions ({NPIs}) to reduce {COVID19} mortality and
  healthcare demand. Tech. rep., Imperial College London.

\bibitem[{Fern{\'a}ndez and Landete(2019)}]{fernandez2015fixed_SAA}
Fern{\'a}ndez, E., Landete, M., 2019. Fixed-charge facility location problems.
  In: Laporte, G., Nickel, S., Saldanha~da Gama, F. (Eds.), Location Science.
  Springer, pp. 67--98.

\bibitem[{Fern{\'a}ndez et~al.(2021)Fern{\'a}ndez, Bornn, and
  Cervone}]{Fernandez2021-do_IM}
Fern{\'a}ndez, J., Bornn, L., Cervone, D., 2021. A framework for the
  fine-grained evaluation of the instantaneous expected value of soccer
  possessions. Machine Learning 110~(6), 1389--1427.

\bibitem[{Ferrer and Whybark(2001)}]{Ferrer2001-yn_AM}
Ferrer, G., Whybark, D.~C., 2001. Material planning for a remanufacturing
  facility. Production and Operations Management 10~(2), 112--124.

\bibitem[{Ferrer et~al.(2018)Ferrer, Mart{\'\i}n-Campo, Ortu{\~n}o,
  Pedraza-Mart{\'\i}nez, Tirado, and Vitoriano}]{ferrer2018multi_BYKOK}
Ferrer, J.~M., Mart{\'\i}n-Campo, F.~J., Ortu{\~n}o, M.~T.,
  Pedraza-Mart{\'\i}nez, A.~J., Tirado, G., Vitoriano, B., 2018. Multi-criteria
  optimization for last mile distribution of disaster relief aid: Test cases
  and applications. European Journal of Operational Research 269~(2), 501--515.

\bibitem[{Ferretti and Montibeller(2019)}]{Ferretti2019-jm_JL}
Ferretti, V., Montibeller, G., 2019. An integrated framework for environmental
  {Multi-Impact} spatial risk analysis. Risk Analysis 39~(1), 257--273.

\bibitem[{Ferretti et~al.(2019)Ferretti, Pluchinotta, and
  Tsouki{\`a}s}]{Ferretti2019-qt_MYLW}
Ferretti, V., Pluchinotta, I., Tsouki{\`a}s, A., 2019. Studying the generation
  of alternatives in public policy making processes. European Journal of
  Operational Research 273~(1), 353--363.

\bibitem[{Fiacco and McCormick(1968)}]{FM68_EAY}
Fiacco, A.~V., McCormick, G.~P., 1968. Nonlinear programming : {S}equential
  unconstrained minimization technique. Wiley, New York, NY.

\bibitem[{Figueira et~al.(2013)Figueira, Greco, Roy, and
  S{\l}owi{\'n}ski}]{figueira2013overview_MESG}
Figueira, J., Greco, S., Roy, B., S{\l}owi{\'n}ski, R., 2013. An overview of
  {ELECTRE} methods and their recent extensions. Journal of Multi-Criteria
  Decision Analysis 20~(1-2), 61--85.

\bibitem[{Figueira et~al.(2016)Figueira, Mousseau, and
  Roy}]{figueira2016electre_MESG}
Figueira, J., Mousseau, V., Roy, B., 2016. {ELECTRE} methods. In: Greco, S.,
  Ehrgott, M., Figueira, J.~R. (Eds.), Multiple Criteria Decision Analysis.
  Springer-Verlag, Berlin, pp. 155--185.

\bibitem[{Fildes et~al.(2009)Fildes, Goodwin, Lawrence, and
  Nikolopoulos}]{Fildes2009-tc_FP}
Fildes, R., Goodwin, P., Lawrence, M., Nikolopoulos, K., 2009. Effective
  forecasting and judgmental adjustments: an empirical evaluation and
  strategies for improvement in supply-chain planning. International Journal of
  Forecasting 25~(1), 3--23.

\bibitem[{Fildes et~al.(2019)Fildes, Goodwin, and {\"O}nkal}]{Fildes2019-uj_FP}
Fildes, R., Goodwin, P., {\"O}nkal, D., 2019. Use and misuse of information in
  supply chain forecasting of promotion effects. International Journal of
  Forecasting 35~(1), 144--156.

\bibitem[{Fildes and Ranyard(1997)}]{Fildes1997-bo_JEB}
Fildes, R., Ranyard, J.~C., 1997. Success and survival of operational research
  {Groups-A} review. Journal of the Operational Research Society 48~(4),
  336--360.

\bibitem[{Filsfils et~al.(2015)Filsfils, Kumar~Nainar, Pignataro,
  Camilo~Cardona, and Fran{\c c}ois}]{filsfils.kumar-nainar.ea:15_BF}
Filsfils, C., Kumar~Nainar, N., Pignataro, C., Camilo~Cardona, J., Fran{\c
  c}ois, P., 2015. The segment routing architecture. IEEE Global Communications
  Conference (GLOBECOM), 1--6.

\bibitem[{Fioole et~al.(2006)Fioole, Kroon, Mar{\'o}ti, and
  Schrijver}]{Fioole2006_DC}
Fioole, P.-J., Kroon, L., Mar{\'o}ti, G., Schrijver, A., 2006. A rolling stock
  circulation model for combining and splitting of passenger trains. European
  Journal of Operational Research 174~(2), 1281--1297.

\bibitem[{Firoozi and Caines(2017)}]{Firoozi17_MB}
Firoozi, D., Caines, P.~E., 2017. The execution problem in finance with major
  and minor traders: A mean field game formulation. In: Apaloo, J., Viscolani,
  B. (Eds.), Annals of the International Society of Dynamic Games (ISDG):
  Advances in Dynamic and Mean Field Games. Vol.~15. Birkh\"auser, Basel, pp.
  107--130.

\bibitem[{Fischetti et~al.(2005)Fischetti, Glover, and Lodi}]{FGL05_ALAL}
Fischetti, M., Glover, F., Lodi, A., 2005. The feasibility pump. Mathematical
  Programming 104, 91--104.

\bibitem[{Fischetti et~al.(2017{\natexlab{a}})Fischetti, Leitner, Ljubi{\'{c}},
  Luipersbeck, Monaci, Resch, Salvagnin, and Sinnl}]{Fischetti:2016_IL}
Fischetti, M., Leitner, M., Ljubi{\'{c}}, I., Luipersbeck, M., Monaci, M.,
  Resch, M., Salvagnin, D., Sinnl, M., 2017{\natexlab{a}}. Thinning out
  {{S}teiner} trees: {A} node-based model for uniform edge costs. Mathematical
  Programming Computation 9~(2), 203--229.

\bibitem[{Fischetti et~al.(2017{\natexlab{b}})Fischetti, Ljubi{\'c}, and
  Sinnl}]{fischetti2017redesigning_SAA}
Fischetti, M., Ljubi{\'c}, I., Sinnl, M., 2017{\natexlab{b}}. Redesigning
  benders decomposition for large-scale facility location. Management Science
  63~(7), 2146--2162.

\bibitem[{Fischetti and Lodi(2003)}]{FL03_ALAL}
Fischetti, M., Lodi, A., 2003. Local branching. Mathematical Programming 98,
  23--47.

\bibitem[{Fisher(1997)}]{Fisher1997_SMD}
Fisher, M., 1997. What is the right supply chain for your product? Harvard
  Business Review March-April, 105--116.

\bibitem[{Fisher et~al.(2018)Fisher, Gallino, and
  Li}]{fisherCompetitionbasedDynamicPricing2018_AKSJF}
Fisher, M., Gallino, S., Li, J., 2018. Competition-based dynamic pricing in
  online retailing: A methodology validated with field experiments. Management
  Science 64~(6), 2496--2514.

\bibitem[{Fisher(1981)}]{fisher1981lagrangian_COIT}
Fisher, M.~L., 1981. The {L}agrangian relaxation method for solving integer
  programming problems. Management Science 27~(1), 1--18.

\bibitem[{Fleckenstein et~al.(2022)Fleckenstein, Klein, and
  Steinhardt}]{fleckensteinRecentAdvancesIntegrating2022_AKSJF}
Fleckenstein, D., Klein, R., Steinhardt, C., 2022. Recent advances in
  integrating demand management and vehicle routing: A methodological review.
  European Journal of Operational Research.

\bibitem[{Fleischmann(2003)}]{Fleischmann2003-fx_AM}
Fleischmann, M., 2003. Reverse logistics network structures and design. In:
  Guide, V. D.~R., Van Wassenhove~L, N. (Eds.), Business Aspects of
  {Closed-Loop} Supply Chains. Carnegie Mellon University Press, Pittsburgh,
  PA, USA.

\bibitem[{Fleischmann et~al.(2004)Fleischmann, Bloemhof-Ruwaard, Beullens, and
  Dekker}]{Fleischmann2004-wv_AM}
Fleischmann, M., Bloemhof-Ruwaard, J.~M., Beullens, P., Dekker, R., 2004.
  Reverse logistics network design. In: Dekker, R., Fleischmann, M.,
  Inderfurth, K., Van~Wassenhove, L.~N. (Eds.), Reverse Logistics: Quantitative
  Models for {Closed-Loop} Supply Chains. Springer-Verlag, Berlin, pp. 65--94.

\bibitem[{Flood(1956)}]{Flood1956_CA_MB}
Flood, M., 1956. The traveling-salesman problem. Operations Research 4~(1),
  61--75.

\bibitem[{Flood and Jackson(1991{\natexlab{a}})}]{Flood1991-fb_AJG}
Flood, R.~L., Jackson, M.~C., 1991{\natexlab{a}}. Creative Problem Solving:
  Total Systems Intervention. Wiley, Chichester.

\bibitem[{Flood and Jackson(1991{\natexlab{b}})}]{Flood1991-ej_GM}
Flood, R.~L., Jackson, M.~C. (Eds.), 1991{\natexlab{b}}. Critical Systems
  Thinking: Directed Readings. Wiley, Chichester.

\bibitem[{Flood and Romm(1996)}]{Flood1996-pg_GM}
Flood, R.~L., Romm, N. R.~A., 1996. Critical Systems Thinking: Current Research
  and Practice. Plenum, New York.

\bibitem[{Fone et~al.(2003)Fone, Hollinghurst, Temple, Round, Lester,
  Weightman, Roberts, Coyle, Bevan, and Palmer}]{Fone2003-uq_CV}
Fone, D., Hollinghurst, S., Temple, M., Round, A., Lester, N., Weightman, A.,
  Roberts, K., Coyle, E., Bevan, G., Palmer, S., 2003. Systematic review of the
  use and value of computer simulation modelling in population health and
  health care delivery. Journal of Public Health Medicine 25~(4), 325--335.

\bibitem[{Ford and Fulkerson(1957)}]{Ford-Fulkerson:1957_IL}
Ford, L.~R., Fulkerson, D.~R., 1957. A simple algorithm for finding the maximum
  network flows and an application to the hitchcock problem. Canadian Journal
  of Mathematics 9, 210--218.

\bibitem[{Ford and Fulkerson(1962)}]{Ford-Fulkerson:1962_IL}
Ford, L.~R., Fulkerson, D.~R., 1962. Flows in networks. Princeton University
  Pres, Princeton, NJ.

\bibitem[{Forest et~al.(2022)Forest, Vigerske, Ralphs, Hafer, Gambini,
  Matthew~Saltzman, Kristjansson, and King}]{CLP_CTCGE}
Forest, J., Vigerske, S., Ralphs, T., Hafer, L., Gambini, H., Matthew~Saltzman,
  S., Kristjansson, B., King, A., 2022. {CLP}.
\newline\urlprefix\url{https://projects.coin-or.org/Clp}

\bibitem[{Forrest and Simmons(2002)}]{Forrest2002-rg_IM}
Forrest, D., Simmons, R., 2002. Outcome uncertainty and attendance demand in
  sport: The case of {E}nglish soccer. Journal of the Royal Statistical
  Society. Series D (The Statistician) 51~(2), 229--241.

\bibitem[{Forrest et~al.(2022)Forrest, Ralphs, Santos, Vigerske, Hafer,
  Kristjansson, Straver, Lubin, Brito, Saltzman, Pitrus, and
  Matsushima}]{Cbc_CTCGE}
Forrest, J., Ralphs, T., Santos, H., Vigerske, S., Hafer, L., Kristjansson, B.,
  Straver, E., Lubin, M., Brito, S., Saltzman, M., Pitrus, B., Matsushima, H.,
  2022. {CBC}.
\newline\urlprefix\url{https://projects.coin-or.org/Bcp}

\bibitem[{Forrester(1958)}]{Forrester1958_SMD}
Forrester, J.~W., 1958. Industrial dynamics—a major breakthrough for decision
  makers. Harvard Business Review 36~(4), 37--66.

\bibitem[{Forrester(1961)}]{Forrester1961-px_MCJM}
Forrester, J.~W., 1961. Industrial Dynamics. MIT Press, Cambridge, MA.

\bibitem[{Fortz(2011)}]{fortz:11_BF}
Fortz, B., 2011. Applications of meta-heuristics to traffic engineering in {IP}
  networks. International Transactions in Operational Research 18~(2),
  131--147.

\bibitem[{Fortz(2015)}]{fortz:15_BF}
Fortz, B., 2015. Location problems in telecommunications. In: Laporte, G.,
  Nickel, S., Saldanha~da Gama, F. (Eds.), Location Science. Springer, pp.
  537--554.

\bibitem[{Fortz(2021)}]{fortz:21_BF}
Fortz, B., 2021. Topology-constrained network design. In: Crainic, T.~G.,
  Gendreau, M., Gendron, B. (Eds.), Network Design with Applications to
  Transportation and Logistics. Springer, pp. 187--208.

\bibitem[{Fortz and Labb\'e(2006)}]{fortz.labbe:06_BF}
Fortz, B., Labb\'e, M., 2006. Design of survivable networks. In: Resende, M.,
  Pardalos, P. (Eds.), Handbook of Optimization in Telecommunications.
  Springer, pp. 367--389.

\bibitem[{Fortz et~al.(2000)Fortz, Labb{\'e}, and
  Maffioli}]{fortz.labbe.ea:00_BF}
Fortz, B., Labb{\'e}, M., Maffioli, F., 2000. Solving the two-connected network
  with bounded meshes problem. Operations Research 48~(6), 866--877.

\bibitem[{Fortz and Thorup(2000)}]{fortz.thorup:00_BF}
Fortz, B., Thorup, M., 2000. Internet traffic engineering by optimizing {OSPF}
  weights. In: Proceedings of the 19th IEEE Conference on Computer
  Communications (INFOCOM). pp. 519--528.

\bibitem[{Foss and Korshunov(2012)}]{foss2012large_HATI}
Foss, S., Korshunov, D., 2012. On large delays in multi-server queues with
  heavy tails. Mathematics of Operations Research 37~(2), 201--218.

\bibitem[{Foulds(1983)}]{foulds1983heuristic_COIT}
Foulds, L.~R., 1983. The heuristic problem-solving approach. Journal of the
  Operational Research Society 34~(10), 927--934.

\bibitem[{Fourier(1826{\natexlab{a}})}]{Fourier1826b-hx_GL}
Fourier, J. B.~J., 1826{\natexlab{a}}. {Analyse des travaux de l’Académie
  Royale des Sciences, pendant l’année 1823}. Partie mathématique, Histoire
  de l’Académie Royale des Sciences de l’Institut de France 6, xxix--xli.

\bibitem[{Fourier(1826{\natexlab{b}})}]{Fourier1826a-hx_GL}
Fourier, J. B.~J., 1826{\natexlab{b}}. Solution d’une question particulière
  du calcul des inégalités. Nouveau Bulletin des Sciences par la Société
  Philomatique de Paris, 99--100.

\bibitem[{Fox and Burks(2019)}]{Fox2019-ky_KVRH}
Fox, W.~P., Burks, R., 2019. Applications of Operations Research and Management
  Science for Military Decision Making. Springer, Cham.

\bibitem[{Fragapane et~al.(2021)Fragapane, De~Koster, Sgarbossa, and
  Strandhagen}]{fragapane2021planning_JLYHK}
Fragapane, G., De~Koster, R., Sgarbossa, F., Strandhagen, J.~O., 2021. Planning
  and control of autonomous mobile robots for intralogistics: Literature review
  and research agenda. European Journal of Operational Research 294~(2),
  405--426.

\bibitem[{Francis et~al.(2004)Francis, McGinnis, and
  White}]{francis2004facility_SAA}
Francis, R., McGinnis, L., White, J., 2004. Facility Layout and Location: An
  Analytical Approach. Prentice-Hall International Series in Industrial and
  Systems Engineering. Prentice Hall.

\bibitem[{Franco et~al.(2022)Franco, Herazo-Padilla, and
  Casta{\~n}eda}]{Franco2022-dx_CV}
Franco, C., Herazo-Padilla, N., Casta{\~n}eda, J.~A., 2022. A queueing network
  approach for capacity planning and patient scheduling: A case study for the
  {COVID-19} vaccination process in colombia. Vaccine 40~(49), 7073--7086.

\bibitem[{Franco et~al.(2016{\natexlab{a}})Franco, , and
  H{\"a}m{\"a}l{\"a}inen}]{Franco2016-ns_AFRH}
Franco, L.~A., , H{\"a}m{\"a}l{\"a}inen, R.~P., 2016{\natexlab{a}}. Feature
  cluster: Behavioural operational research. European Journal of Operational
  Research 249~(3), 791--1073.

\bibitem[{Franco(2013)}]{Alberto_Franco2013-yw_MYLW}
Franco, L.~A., 2013. Rethinking soft {OR} interventions: Models as boundary
  objects. European Journal of Operational Research 231~(3), 720--733.

\bibitem[{Franco et~al.(2004)Franco, Cushman, and
  Rosenhead}]{Franco2004-qi_MYLW}
Franco, L.~A., Cushman, M., Rosenhead, J., 2004. Project review and learning in
  the construction industry: Embedding a problem structuring method within a
  partnership context. European Journal of Operational Research 152~(3),
  586--601.

\bibitem[{Franco and Greiffenhagen(2018)}]{Franco2018-nk_AFRH}
Franco, L.~A., Greiffenhagen, C., 2018. Making {OR} practice visible: Using
  ethnomethodology to analyse facilitated modelling workshops. European Journal
  of Operational Research 265~(2), 673--684.

\bibitem[{Franco and H{\"a}m{\"a}l{\"a}inen(2016)}]{Franco2016-qk_MYLW}
Franco, L.~A., H{\"a}m{\"a}l{\"a}inen, R.~P., 2016. Behavioural operational
  research: Returning to the roots of the {OR} profession. European Journal of
  Operational Research 249~(3), 791--795.

\bibitem[{Franco et~al.(2021)Franco, H{\"a}m{\"a}l{\"a}inen, Rouwette, and
  Lepp{\"a}nen}]{Franco2021-sm_JL}
Franco, L.~A., H{\"a}m{\"a}l{\"a}inen, R.~P., Rouwette, E. A. J.~A.,
  Lepp{\"a}nen, I., 2021. Taking stock of behavioural {OR}: A review of
  behavioural studies with an intervention focus. European Journal of
  Operational Research 293~(2), 401--418.

\bibitem[{Franco and Montibeller(2010)}]{Franco2010-pf_MYLW}
Franco, L.~A., Montibeller, G., 2010. Facilitated modelling in operational
  research. European Journal of Operational Research 205~(3), 489--500.

\bibitem[{Franco and Rouwette(2022)}]{Franco2022-mn_AFRH}
Franco, L.~A., Rouwette, E. A. J.~A., 2022. Problem structuring methods: Taking
  stock and looking ahead. In: Salhi, S., Boylan, J. (Eds.), The Palgrave
  Handbook of Operations Research. Springer, Cham, pp. 735--780.

\bibitem[{Franco et~al.(2016{\natexlab{b}})Franco, Rouwette, and
  Korzilius}]{Franco2016-yn_AFRH}
Franco, L.~A., Rouwette, E. A. J.~A., Korzilius, H., 2016{\natexlab{b}}.
  Different paths to consensus? {T}he impact of need for closure on
  model-supported group conflict management. European Journal of Operational
  Research 249~(3), 878--889.

\bibitem[{Frangioni and Gendron(2021)}]{Frangioni2021_MH}
Frangioni, A., Gendron, B., 2021. Piecewise linear cost network design. In:
  Crainic, T.~G., Gendreau, M., Gendron, B. (Eds.), Network Design with
  Applications to Transportation and Logistics. Springer, pp. 167--185.

\bibitem[{Franke(2017)}]{franke2017network_MH}
Franke, M., 2017. Network design strategies. In: Bruce, P.~J., Gao, Y., King,
  J. M.~C. (Eds.), Airline Operations. Routledge, pp. 44--60.

\bibitem[{Franses et~al.(2014)Franses, van Dijk, and
  Opschoor}]{Franses2014-kl_FP}
Franses, P.~H., van Dijk, D., Opschoor, A., 2014. Time Series Models for
  Business and Economic Forecasting. Cambridge University Press.

\bibitem[{Fraunholz et~al.(2021)Fraunholz, Kraft, Keles, and
  Fichtner}]{fra:kra:kel:fic:21_DSRW}
Fraunholz, C., Kraft, E., Keles, D., Fichtner, W., 2021. Advanced price
  forecasting in agent-based electricity market simulation. Applied Energy 290,
  116688.

\bibitem[{Freeman(2003)}]{Fre03_JH}
Freeman, S. (Ed.), 2003. The Cambridge Companion to Rawls. Cambridge University
  Press.

\bibitem[{Fregonara et~al.(2013)Fregonara, Curto, Grosso, Mellano, Rolando, and
  Tulliani}]{Fregonara2013-yd_MYLW}
Fregonara, E., Curto, R., Grosso, M., Mellano, P., Rolando, D., Tulliani,
  J.-M., 2013. Environmental technology, materials science, architectural
  design, and real estate market evaluation: A multidisciplinary approach for
  {Energy-Efficient} buildings. Journal of Urban Technology 20~(4), 57--80.

\bibitem[{French(2022)}]{French2022-kr_JL}
French, S., 2022. From soft to hard elicitation. Journal of the Operational
  Research Society 73~(6), 1181--1197.

\bibitem[{French and Geldermann(2005)}]{French2005-pc_JL}
French, S., Geldermann, J., 2005. The varied contexts of environmental decision
  problems and their implications for decision support. Environmental Science
  \& Policy 8~(4), 378--391.

\bibitem[{Friend(1989)}]{Friend1989-yj_MYLW}
Friend, J., 1989. The strategic choice approach. In: Rosenhead, J. (Ed.),
  Rational analysis for a problematic world: problem structuring methods for
  complexity, uncertainty, and conflict. Wiley, Chichester, pp. 121--158.

\bibitem[{Friesl et~al.(2017)Friesl, Lenten, Libich, and
  Stehl{\'\i}k}]{Friesl2017-ji_IM}
Friesl, M., Lenten, L. J.~A., Libich, J., Stehl{\'\i}k, P., 2017. In search of
  goals: increasing ice hockey's attractiveness by a sides swap. Journal of the
  Operational Research Society 68~(9), 1006--1018.

\bibitem[{Fritzson et~al.(2020)Fritzson, Pop, Abdelhak, Ashgar, Bachmann,
  Braun, Bouskela, Braun, Buffoni, Casella, Castro, Franke, Fritzson,
  Gebremedhin, Heuermann, Lie, Mengist, Mikelsons, Moudgalya, Ochel,
  Palanisamy, Ruge, Schamai, Sjölund, Thiele, Tinnerholm, and
  Östlund}]{Openmodelica_CTCGE}
Fritzson, P., Pop, A., Abdelhak, K., Ashgar, A., Bachmann, B., Braun, W.,
  Bouskela, D., Braun, R., Buffoni, L., Casella, F., Castro, R., Franke, R.,
  Fritzson, D., Gebremedhin, M., Heuermann, A., Lie, B., Mengist, A.,
  Mikelsons, L., Moudgalya, K., Ochel, L., Palanisamy, A., Ruge, V., Schamai,
  W., Sjölund, M., Thiele, B., Tinnerholm, J., Östlund, P., 2020. {The
  OpenModelica Integrated Environment for Modeling, Simulation, and Model-Based
  Development}. Modeling, Identification and Control 41~(4), 241--295.

\bibitem[{Fu et~al.(2005)Fu, Glover, and April}]{Fu2005-qb_HL}
Fu, M.~C., Glover, F.~W., April, J., 2005. Simulation optimization: a review,
  new developments, and applications. In: Kuhl, M.~E., Steiger, N.~M.,
  Armstrong, F.~B., Joines, J.~A. (Eds.), Proceedings of the Winter Simulation
  Conference, 2005. pp. 1--13.

\bibitem[{F{\"u}gener et~al.(2014)F{\"u}gener, Hans, Kolisch, Kortbeek, and
  Vanberkel}]{Fugener2014-jk_CV}
F{\"u}gener, A., Hans, E.~W., Kolisch, R., Kortbeek, N., Vanberkel, P.~T.,
  2014. Master surgery scheduling with consideration of multiple downstream
  units. European Journal of Operational Research 239~(1), 227--236.

\bibitem[{Fukasawa et~al.(2006)Fukasawa, Longo, Lysgaard, and
  Poggi}]{Fukasawa2006_CA_MB}
Fukasawa, R., Longo, H., Lysgaard, J., Poggi, M., 2006. Robust
  branch-and-cut-and-price for the capacitated vehicle routing problem.
  Mathematical Programming 106~(3), 491--511.

\bibitem[{Fukuda(1964)}]{fukuda1964optimal_JSS}
Fukuda, Y., 1964. Optimal policies for the inventory problem with negotiable
  leadtime. Management Science 10~(4), 690--708.

\bibitem[{Fukuyama and Weber(2009)}]{Fukuyama2009-mo_SL}
Fukuyama, H., Weber, W.~L., 2009. A directional slacks-based measure of
  technical inefficiency. Socio-Economic Planning Sciences 43~(4), 274--287.

\bibitem[{Färe et~al.(1983)Färe, Grosskopf, and Logan}]{Fare_1983_DSRW}
Färe, R., Grosskopf, S., Logan, J., 1983. The relative efficiency of
  {I}llinois electric utilities. Resources and Energy 5~(4), 349--367.

\bibitem[{Färe et~al.(1996)Färe, Grosskopf, and Tyteca}]{Fare_1996_DSRW}
Färe, R., Grosskopf, S., Tyteca, D., 1996. An activity analysis model of the
  environmental performance of firms—application to fossil-fuel-fired
  electric utilities. Ecological Economics 18~(2), 161--175.

\bibitem[{Gaalman(2006)}]{Gaalman2006_XW}
Gaalman, G., 2006. Bullwhip reduction for {ARMA} demand: The proportional
  order-up-to policy versus the full-state-feedback policy. Automatica 42~(8),
  1283--1290.

\bibitem[{Gaalman et~al.(2022)Gaalman, Disney, and Wang}]{Gaalman2022_SMD}
Gaalman, G., Disney, S.~M., Wang, X., 2022. When bullwhip increases in the lead
  time: An eigenvalue analysis of {ARMA} demand. International Journal of
  Production Economics 250, 108623.

\bibitem[{G{\'a}cs and Lov{\'a}sz(1981)}]{Gacs1981-kz_JMB}
G{\'a}cs, P., Lov{\'a}sz, L., 1981. Khachiyan's algorithm for linear
  programming. In: K{\"o}nig, H., Korte, B., Ritter, K. (Eds.), Mathematical
  Programming at Oberwolfach. Springer Berlin Heidelberg, Berlin, Heidelberg,
  pp. 61--68.

\bibitem[{Gaillard et~al.(2016)Gaillard, Goude, and
  Nedellec}]{gai:gou:ned:16_DSRW}
Gaillard, P., Goude, Y., Nedellec, R., 2016. Additive models and robust
  aggregation for {GEFCom2014} probabilistic electric load and electricity
  price forecasting. International Journal of Forecasting 32~(3), 1038--1050.

\bibitem[{Gal{\'a}n et~al.(2022)Gal{\'a}n, Carrasco, and
  LaTorre}]{Galan2022-xx_KVRH}
Gal{\'a}n, J.~J., Carrasco, R.~A., LaTorre, A., 2022. Military applications of
  machine learning: A bibliometric perspective. Science in China, Series A:
  Mathematics 10~(9), 1397.

\bibitem[{Galbreth and Blackburn(2006)}]{Galbreth2009-ar_AM}
Galbreth, M.~R., Blackburn, J.~D., 2006. Optimal acquisition and sorting
  policies for remanufacturing. Production and Operations Management 15~(3),
  384--392.

\bibitem[{Galiullina et~al.(2022)Galiullina, Mutlu, Kinable, and
  Van~Woensel}]{Galiullina2022_CKTVW}
Galiullina, A., Mutlu, N., Kinable, J., Van~Woensel, T., 2022. Demand steering
  in a last-mile delivery problem with multiple delivery channels. Working
  paper.

\bibitem[{Gallego and Moon(1993)}]{gallego1993distribution_JSS}
Gallego, G., Moon, I., 1993. The distribution free newsboy problem: review and
  extensions. Journal of the Operational Research Society 44~(8), 825--834.

\bibitem[{Gallego and Topalo{\u
  g}lu(2019)}]{gallegoRevenueManagementPricing2019_AKSJF}
Gallego, G., Topalo{\u g}lu, H., 2019. Revenue Management and Pricing
  Analytics. Vol. 279 of International Series in Operations Research \&
  Management Science. {Springer}, {New York, NY}.

\bibitem[{Gallego and {van
  Ryzin}(1994)}]{gallegoOptimalDynamicPricing1994_AKSJF}
Gallego, G., {van Ryzin}, G., 1994. Optimal {{Dynamic Pricing}} of
  {{Inventories}} with {{Stochastic Demand}} over {{Finite Horizons}}.
  Management Science 40~(8), 999--1020.

\bibitem[{Gallego and Van~Ryzin(1994)}]{gallego1994optimal_VLVV}
Gallego, G., Van~Ryzin, G., 1994. Optimal dynamic pricing of inventories with
  stochastic demand over finite horizons. Management Science 40~(8), 999--1020.

\bibitem[{Gallego and {van
  Ryzin}(1997)}]{gallegoMultiproductDynamicPricing1997_AKSJF}
Gallego, G., {van Ryzin}, G., 1997. A {{Multiproduct Dynamic Pricing Problem}}
  and {{Its Applications}} to {{Network Yield Management}}. Operations Research
  45~(1), 24--41.

\bibitem[{Gamrath et~al.(2017)Gamrath, Koch, Maher, Rehfeldt, and
  Shinano}]{gamrathscip_IL}
Gamrath, G., Koch, T., Maher, S.~J., Rehfeldt, D., Shinano, Y., 2017.
  {SCIP-Jack} --- a solver for {STP} and variants with parallelization
  extensions. Mathematical Programming Computation 9~(2), 231--296.

\bibitem[{Gan and Lee(2022)}]{Gan2021-fz_SL}
Gan, G., Lee, H.-S., 2022. Measuring group performance based on metafrontier.
  Journal of the Operational Research Society 73~(10), 2261--2274.

\bibitem[{Gao et~al.(2019)Gao, Darvishan, Toghani, Mohammadi, Abedinia, and
  Ghadimi}]{gao:etal:19_DSRW}
Gao, W., Darvishan, A., Toghani, M., Mohammadi, M., Abedinia, O., Ghadimi, N.,
  2019. Different states of multi-block based forecast engine for price and
  load prediction. International Journal of Electrical Power and Energy Systems
  104, 423--435.

\bibitem[{Garc{\'\i}a and Mar{\'\i}n(2019)}]{garcia2015covering_SAA}
Garc{\'\i}a, S., Mar{\'\i}n, A., 2019. Covering location problems. In: Laporte,
  G., Nickel, S., {Saldanha da Gama}, F. (Eds.), Location Science. Springer,
  pp. 99--119.

\bibitem[{Gardner(2006)}]{Gardner2006-bv_FP}
Gardner, E.~S., 2006. Exponential smoothing: The state of the art - part {II}.
  International Journal of Forecasting 22~(4), 637--666.

\bibitem[{Gardner and Righter(2020)}]{rhonda_HATI}
Gardner, K., Righter, R., 2020. Product forms for {FCFS} queueing models with
  arbitrary server-job compatibilities: an overview. Queueing Systems 96~(1),
  3--51.

\bibitem[{Garey and Johnson(1979)}]{GJ79_SMPT}
Garey, M., Johnson, D., 1979. Computers and intractability: A guide to the
  theory of {NP}-completeness. W.H. Freeman and Company, New York, NY.

\bibitem[{Garfinkel and Nemhauser(1972)}]{GN72_SMPT}
Garfinkel, R., Nemhauser, G., 1972. Integer programming. John Wiley\&Sons, New
  York.

\bibitem[{Garrett(2014)}]{Garrett2014-on_KESXZ}
Garrett, B., 2014. {3D} printing: New economic paradigms and strategic shifts.
  Global Policy 5~(1), 70--75.

\bibitem[{Gary et~al.(2008)Gary, Kunc, Morecroft, and
  Rockart}]{Gary2008-uu_MCJM}
Gary, M.~S., Kunc, M., Morecroft, J. D.~W., Rockart, S.~F., 2008. System
  dynamics and strategy. System Dynamics Review 24~(4), 407--429.

\bibitem[{Gasse et~al.(2019)Gasse, Chetelat, Ferroni, Charlin, and
  Lodi}]{gasse2019exact_LCAL}
Gasse, M., Chetelat, D., Ferroni, N., Charlin, L., Lodi, A., 2019. Exact
  combinatorial optimization with graph convolutional neural networks. In:
  Wallach, H., Larochelle, H., Beygelzimer, A., d\textquotesingle
  Alch\'{e}-Buc, F., Fox, E., Garnett, R. (Eds.), Advances in Neural
  Information Processing Systems. Vol.~32. Curran Associates, Inc.

\bibitem[{Gaur and Fisher(2004)}]{gaur2004periodic_JLYHK}
Gaur, V., Fisher, M.~L., 2004. A periodic inventory routing problem at a
  supermarket chain. Operations Research 52~(6), 813--822.

\bibitem[{Gautam(2012)}]{gautam_HATI}
Gautam, N., 2012. Analysis of Queues: Methods and Applications. CRC Press.

\bibitem[{Gauthier(1983)}]{Gau87_JH}
Gauthier, D., 1983. Morals by Agreement. Oxford University Press.

\bibitem[{Geis et~al.(2011)Geis, Parnell, Newton, and
  Bresnick}]{Geis2011-zd_KVRH}
Geis, J.~P., Parnell, G.~S., Newton, H., Bresnick, T., 2011. Blue horizons
  study assesses future capabilities and technologies for the {United States
  Air Force}. Interfaces 41~(4), 338--353.

\bibitem[{Geldermann et~al.(2009)Geldermann, Bertsch, Treitz, French,
  Papamichail, and H{\"a}m{\"a}l{\"a}inen}]{Geldermann2009-cz_JL}
Geldermann, J., Bertsch, V., Treitz, M., French, S., Papamichail, K.~N.,
  H{\"a}m{\"a}l{\"a}inen, R.~P., 2009. Multi-criteria decision support and
  evaluation of strategies for nuclear remediation management. Omega 37~(1),
  238--251.

\bibitem[{Gendreau et~al.(2015)Gendreau, Ghiani, and
  Guerriero}]{gendreau2015time_CA_MB}
Gendreau, M., Ghiani, G., Guerriero, E., 2015. Time-dependent routing problems:
  A review. Computers \& Operations Research 64, 189--197.

\bibitem[{Gendreau et~al.(2016)Gendreau, Jabali, and Rei}]{Gendreau2016_CA_MB}
Gendreau, M., Jabali, O., Rei, W., 2016. 50th anniversary invited
  article—future research directions in stochastic vehicle routing.
  Transportation Science 50~(4), 1163--1173.

\bibitem[{Gendreau and Potvin(2010)}]{gendreau2010handbook_CA_MB}
Gendreau, M., Potvin, J.-Y., 2010. Handbook of metaheuristics. Vol.~2.
  Springer.

\bibitem[{Geng and Xie(2019)}]{Geng2019-ph_HL}
Geng, X., Xie, L., 2019. Data-driven decision making in power systems with
  probabilistic guarantees: Theory and applications of chance-constrained
  optimization. Annual Reviews in Control 47, 341--363.

\bibitem[{Gentili et~al.(2022)Gentili, Mirchandani, Agnetis, and
  Ghelichi}]{GENTILI2022108057_BKOK}
Gentili, M., Mirchandani, P.~B., Agnetis, A., Ghelichi, Z., 2022. Locating
  platforms and scheduling a fleet of drones for emergency delivery of
  perishable items. Computers \& Industrial Engineering 168, 108057.
\newline\urlprefix\url{https://www.sciencedirect.com/science/article/pii/S0360835222001279}

\bibitem[{Georgantas et~al.(2021)Georgantas, Doumpos, and
  Zopounidis}]{Georgantas2021-nj_HL}
Georgantas, A., Doumpos, M., Zopounidis, C., 2021. Robust optimization
  approaches for portfolio selection: a comparative analysis. Annals of
  Operations Research.

\bibitem[{Gettinger et~al.(2013)Gettinger, Kiesling, Stummer, and
  Vetschera}]{Gettinger2013-mi_AFRH}
Gettinger, J., Kiesling, E., Stummer, C., Vetschera, R., 2013. A comparison of
  representations for discrete multi-criteria decision problems. Decision
  Support Systems 54~(2), 976--985.

\bibitem[{Ghamlouche et~al.(2003)Ghamlouche, Crainic, and
  Gendreau}]{ghamlouche2003cycle_MH}
Ghamlouche, I., Crainic, T.~G., Gendreau, M., 2003. Cycle-based neighbourhoods
  for fixed-charge capacitated multicommodity network design. Operations
  Research 51~(4), 655--667.

\bibitem[{Ghelichi et~al.(2021)Ghelichi, Gentili, and
  Mirchandani}]{GHELICHI2021105443_BKOK}
Ghelichi, Z., Gentili, M., Mirchandani, P.~B., 2021. Logistics for a fleet of
  drones for medical item delivery: A case study for louisville, ky. Computers
  \& Operations Research 135, 105443.
\newline\urlprefix\url{https://www.sciencedirect.com/science/article/pii/S0305054821001970}

\bibitem[{Giacco et~al.(2014)Giacco, Carillo, D{'}Ariano, Pacciarelli, and
  Mar{\'\i}n}]{Giacco2014_DC}
Giacco, G.~L., Carillo, D., D{'}Ariano, A., Pacciarelli, D., Mar{\'\i}n, A.,
  2014. Short-term rail rolling stock rostering and maintenance scheduling.
  Transportation Research Procedia 3, 651--659.

\bibitem[{Gilbert(1963)}]{Gilbert1963_XW}
Gilbert, E.~G., 1963. Controllability and observability in multivariable
  control systems. Journal of the Society for Industrial and Applied
  Mathematics, Series A: Control 1~(2), 128--151.

\bibitem[{Gilboa et~al.(2017)Gilboa, Minardi, and Samuelson}]{Gilboa2017-hp_TC}
Gilboa, I., Minardi, S., Samuelson, L., 2017. Cases and scenarios in decisions
  under uncertainty. HEC Paris Research Paper No. ECO/SCD-2017-1200, 1--49.

\bibitem[{Gilbreth(1911)}]{Gilbreth1911_SMD}
Gilbreth, F.~B., 1911. Motion study: A method for increasing the efficiency of
  the workman. D. Van Nostrand Company.

\bibitem[{Gillies(1953)}]{Gi1953_GZ}
Gillies, D., 1953. Some theorems on $n$-person games. {Ph.D.} thesis, Princeton
  University.

\bibitem[{Gilmore and Gomory(1965)}]{Gilmore1965-mm_JB}
Gilmore, P.~C., Gomory, R.~E., 1965. Multistage cutting stock problems of two
  and more dimensions. Operations Research 13~(1), 94--120.

\bibitem[{Gini(1912)}]{Gin12_JH}
Gini, C., 1912. Variabilit\`{a} e mutabilit\`{a}. P. Cuppini, reprinted 1955 in
  E.\ Pizetti abd T.\ Salvemini, eds., {\em Memorie di metodologica
  statistica}, Rome: Libreria Eredi Virgilio Veschi.

\bibitem[{Giovannini and Psaraftis(2019)}]{Giovannini2019-fb_HP}
Giovannini, M., Psaraftis, H.~N., 2019. The profit maximizing liner shipping
  problem with flexible frequencies: logistical and environmental
  considerations. Flexible Services and Manufacturing Journal 31~(3), 567--597.

\bibitem[{Glaize et~al.(2019)Glaize, Duenas, Di~Martinelly, and
  Fagnot}]{Glaize2019-jz_CV}
Glaize, A., Duenas, A., Di~Martinelly, C., Fagnot, I., 2019. Healthcare
  decision-making applications using multicriteria decision analysis: A scoping
  review. Journal of Multi-criteria Decision Analysis 26~(1-2), 62--83.

\bibitem[{Glasgow et~al.(2018)Glasgow, Perkins, Tai, Brohi, and
  Vasilakis}]{Glasgow2018-th_CV}
Glasgow, S.~M., Perkins, Z.~B., Tai, N. R.~M., Brohi, K., Vasilakis, C., 2018.
  Development of a discrete event simulation model for evaluating strategies of
  red blood cell provision following mass casualty events. European Journal of
  Operational Research 270~(1), 362--374.

\bibitem[{Glazebrook et~al.(2014)Glazebrook, Hodge, Kirkbride, and
  Minty}]{glazebrook2014stochastic_DLLD}
Glazebrook, K.~D., Hodge, D.~J., Kirkbride, C., Minty, R., 2014. Stochastic
  scheduling: A short history of index policies and new approaches to index
  generation for dynamic resource allocation. Journal of Scheduling 17~(5),
  407--425.

\bibitem[{Gleixner et~al.(2021)Gleixner, Hendel, Gamrath, et~al.}]{Gl17_ALAL}
Gleixner, A., Hendel, G., Gamrath, G., et~al., 2021. {MIPLIB} 2017: data-driven
  compilation of the 6th mixed-integer programming library. Mathematical
  Programming Computation 13, 443--490.

\bibitem[{Glen(1997)}]{Glen1997_EA}
Glen, J.~J., 1997. An infinite horizon mathematical programming model of a
  multicohort single species fishery. Journal of the Operational Research
  Society 48~(11), 1095--1104.

\bibitem[{Glickman(2001)}]{Glickman2001-vc_IM}
Glickman, M.~E., 2001. Dynamic paired comparison models with stochastic
  variances. Journal of Applied Statistics 28~(6), 673--689.

\bibitem[{Glover(1977)}]{glover1977heuristics_COIT}
Glover, F., 1977. Heuristics for integer programming using surrogate
  constraints. Decision Sciences 8~(1), 156--166.

\bibitem[{Glover(1986)}]{glover1986future_COIT}
Glover, F., 1986. Future paths for integer programming and links to artificial
  intelligence. Computers \& Operations Research 13~(5), 533--549.

\bibitem[{Glover(1990)}]{glover1990tabu_COIT}
Glover, F., 1990. Tabu search: {A} tutorial. Interfaces 20~(4), 74--94.

\bibitem[{Glover(1998)}]{Glover1998-hd_HL}
Glover, F., 1998. A template for scatter search and path relinking. In: Hao,
  J.~K., Lutton, E., Ronald, E., Schoenauer, M., Snyers, D. (Eds.), Lecture
  Notes in Computer Science. Vol. 1363 of Lecture Notes in Computer Science.
  Springer, Berlin, Heidelberg, pp. 1--51.

\bibitem[{Glover and Laguna(1997)}]{Glover1997-bi_HL}
Glover, F., Laguna, M., 1997. Tabu Search. Springer US.

\bibitem[{Glover and Kochenberger(2003)}]{Glover2003-cz_HL}
Glover, F.~W., Kochenberger, G.~A., 2003. Handbook of Metaheuristics. Springer
  Science \& Business Media.

\bibitem[{Gneiting and Raftery(2007)}]{Gneiting2007-oj_FP}
Gneiting, T., Raftery, A.~E., 2007. Strictly proper scoring rules, prediction,
  and estimation. Journal of the American Statistical Association 102~(477),
  359--378.

\bibitem[{Goemans(1994)}]{Goemans:1994related_IL}
Goemans, M.~X., 1994. The {{S}teiner} tree polytope and related polyhedra.
  Mathematical Programming 63~(1), 157--182.

\bibitem[{Goemans et~al.(2012)Goemans, Olver, Rothvo{\ss}, and
  Zenklusen}]{Goemans:2012matroids_IL}
Goemans, M.~X., Olver, N., Rothvo{\ss}, T., Zenklusen, R., 2012. Matroids and
  integrality gaps for hypergraphic {S}teiner tree relaxations. In: Karloff,
  H.~J., Pitassi, T. (Eds.), Proceedings of the 44th Symposium on Theory of
  Computing Conference, {STOC} 2012, New York, NY, USA, May 19 - 22, 2012.
  {ACM}, New York, USA, pp. 1161--1176.

\bibitem[{Goldberg et~al.(2015)Goldberg, Hed, Kaplan, Kohli, Tarjan, and
  Werneck}]{Goldberg-et-al:2015_IL}
Goldberg, A.~V., Hed, S., Kaplan, H., Kohli, P., Tarjan, R.~E., Werneck, R.~F.,
  2015. Faster and more dynamic maximum flow by incremental breadth-first
  search. In: Bansal, N., Finocchi, I. (Eds.), Algorithms - {ESA} 2015 - 23rd
  Annual European Symposium, Patras, Greece, September 14-16, 2015,
  Proceedings. Vol. 9294 of Lecture Notes in Computer Science. Springer, pp.
  619--630.

\bibitem[{Goldberg and Tarjan(1988)}]{Goldberg-Tarjan:1988_IL}
Goldberg, A.~V., Tarjan, R.~E., 1988. A new approach to the maximum-flow
  problem. Journal of the {ACM} 35~(4), 921--940.

\bibitem[{Goldberg and Tarjan(2014)}]{Goldberg-Tarjan:2014_IL}
Goldberg, A.~V., Tarjan, R.~E., 2014. Efficient maximum flow algorithms.
  Communications of the {ACM} 57~(8), 82--89.

\bibitem[{Golden et~al.(2008)Golden, Raghavan, and Wasil}]{GoldenRW2008_CA_MB}
Golden, B., Raghavan, S., Wasil, E. (Eds.), 2008. {The Vehicle Routing Problem:
  Latest Advances and New Challenges}. Springer, New York.

\bibitem[{Goldfarb and Todd(1989)}]{Goldfarb1989-bz_JMB}
Goldfarb, D., Todd, M.~J., 1989. Chapter {II} linear programming. In:
  Nemhauser, G.~L., Rinnooy~Kan, A. H.~G., Todd, M.~J. (Eds.), Handbooks in
  Operations Research and Management Science. Vol.~1. Elsevier, North-Holland,
  Amsterdam, pp. 73--170.

\bibitem[{Goldratt(1997)}]{Goldratt1997-tq_WH_ED}
Goldratt, E.~M., 1997. Critical Chain. North River Press.

\bibitem[{Gollowitzer and Ljubi{\'c}(2011)}]{gollowitzer.ljubic:11_BF}
Gollowitzer, S., Ljubi{\'c}, I., 2011. {MIP} models for connected facility
  location: A theoretical and computational study. Computers \& Operations
  Research 38, 435--449.

\bibitem[{Gomes and Sa\'{u}de(2014)}]{Gosa2014_GZ}
Gomes, D., Sa\'{u}de, J., 2014. Mean field games models---{A} brief survey.
  Dynamic Games and Applications 4, 110--154.

\bibitem[{Gomes and Lima(1991)}]{gomes1991todim}
Gomes, L. F. A.~M., Lima, M. M. P.~P., 1991. Todim: Basics and application to
  multicriteria ranking of projects with environmental impact. Foundations of
  Computing and Decision Sciences 16~(3-4), 1--16.

\bibitem[{G{\'o}mez-Rocha and
  Hern{\'a}ndez-Gress(2022)}]{Gomez-Rocha2022-kp_HL}
G{\'o}mez-Rocha, J.~E., Hern{\'a}ndez-Gress, E.~S., 2022. A stochastic
  programming model for {Multi-Product} aggregate production planning using
  valid inequalities. NATO Advanced Science Institutes series E: Applied
  sciences 12~(19), 9903.

\bibitem[{Gomory(1958)}]{Go58_ALAL}
Gomory, R., 1958. Outline of an algorithm for integer solutions to linear
  programs. Bulletin of the American Mathematical Society 64, 275--278.

\bibitem[{Gomory(1960)}]{Go60_ALAL}
Gomory, R., 1960. An algorithm for the mixed integer problem. Tech. Rep.
  RM-2597, RAND Corporation.

\bibitem[{Gomory(1969)}]{Go69_ALAL}
Gomory, R., 1969. Some polyhedra related to combinatorial problems. Linear
  Algebra and its Applications 2, 451--558.

\bibitem[{Gomory and Hu(1961)}]{Gomory-Hu:1961_IL}
Gomory, R.~E., Hu, T.~C., 1961. Multi-terminal network flows. Journal of the
  Society for Industrial and Applied Mathematics 9~(4), 551--570.

\bibitem[{Goodchild and Daganzo(2007)}]{Goodchild2007-rb_HP}
Goodchild, A.~V., Daganzo, C.~F., 2007. Crane double cycling in container
  ports: Planning methods and evaluation. Transportation Research Part B:
  Methodological 41~(8), 875--891.

\bibitem[{Goodfellow et~al.(2016)Goodfellow, Bengio, and
  Courville}]{Goodfellow-et-al-2016_LCAL}
Goodfellow, I., Bengio, Y., Courville, A., 2016. Deep Learning. MIT Press,
  Cambridge MA.

\bibitem[{Goodson et~al.(2013)Goodson, Ohlmann, and Thomas}]{Goodson2013-fd_HL}
Goodson, J.~C., Ohlmann, J.~W., Thomas, B.~W., 2013. Rollout policies for
  dynamic solutions to the multivehicle routing problem with stochastic demand
  and duration limits. Operations Research 61~(1), 138--154.

\bibitem[{{Google}(2022)}]{ORTools_CTCGE}
{Google}, 2022. {OR-Tools}.
\newline\urlprefix\url{https://developers.google.com/optimization}

\bibitem[{Goossens et~al.(2006)Goossens, van Hoesel, and
  Kroon}]{Goossens2006_DC}
Goossens, J.-W., van Hoesel, S., Kroon, L., 2006. On solving multi-type railway
  line planning problems. European Journal of Operational Research 168~(2),
  403--424.

\bibitem[{Gopalan and Talluri(1998)}]{gopalan1998aircraft_VLVV}
Gopalan, R., Talluri, K.~T., 1998. The aircraft maintenance routing problem.
  Operations Research 46~(2), 260--271.

\bibitem[{Gordon and Newell(1967)}]{gordon_HATI}
Gordon, W.~J., Newell, G.~F., 1967. Cyclic queuing systems with restricted
  length queues. Operations Research 15~(2), 266--277.

\bibitem[{Gorissen et~al.(2015)Gorissen, Yan{\i}ko{\u g}lu, and den
  Hertog}]{Gorissen2015-zv_HL}
Gorissen, B.~L., Yan{\i}ko{\u g}lu, {\.I}., den Hertog, D., 2015. A practical
  guide to robust optimization. Omega 53, 124--137.

\bibitem[{Gouveia(1998)}]{gouveia:98_BF}
Gouveia, L., 1998. Using variable redefinition for computing lower bounds for
  minimum spanning and {S}teiner trees with hop constraints. INFORMS Journal on
  Computing 10, 180--188.

\bibitem[{Govindan et~al.(2014)Govindan, Jafarian, Khodaverdi, and
  Devika}]{GOVINDAN2014_JLYHK}
Govindan, K., Jafarian, A., Khodaverdi, R., Devika, K., 2014. Two-echelon
  multiple-vehicle location–routing problem with time windows for
  optimization of sustainable supply chain network of perishable food.
  International Journal of Production Economics 152, 9--28.

\bibitem[{Graefe and Armstrong(2011)}]{Graefe2011-ha_FP}
Graefe, A., Armstrong, J.~S., 2011. Comparing face-to-face meetings, nominal
  groups, {Delphi} and prediction markets on an estimation task. International
  Journal of Forecasting 27~(1), 183--195.

\bibitem[{Grass and Fischer(2016)}]{Grass2016-rg_HL}
Grass, E., Fischer, K., 2016. Two-stage stochastic programming in disaster
  management: A literature survey. Surveys in Operations Research and
  Management Science 21~(2), 85--100.

\bibitem[{Graves et~al.(1993{\natexlab{a}})Graves, McBride, Gershkoff,
  Anderson, and Mahidhara}]{graves1993flight_VLVV}
Graves, G.~W., McBride, R.~D., Gershkoff, I., Anderson, D., Mahidhara, D.,
  1993{\natexlab{a}}. Flight crew scheduling. Management Science 39~(6),
  736--745.

\bibitem[{Graves et~al.(1993{\natexlab{b}})Graves, Kan, and
  Zipkin}]{graves1993logistics_JSS}
Graves, S.~C., Kan, A.~R., Zipkin, P.~H., 1993{\natexlab{b}}. Logistics of
  production and inventory. Vol.~4. Elsevier.

\bibitem[{Greco and Slowinski(2001)}]{greco2001rough_MESG}
Greco, S.and~Matarazzo, B., Slowinski, R., 2001. Rough sets theory for
  multicriteria decision analysis. European Journal of Operational Research
  129~(1), 1--47.

\bibitem[{Greco et~al.(2016)Greco, Ehrgott, and
  Figueira}]{greco2016multiple_MESG}
Greco, S., Ehrgott, M., Figueira, J. (Eds.), 2016. Multiple criteria decision
  analysis: State of the art surveys. Vol. 233 of International Series in
  Operations Research and Management Science. Springer-Verlag, Berlin.

\bibitem[{Green and Armstrong(2007)}]{Green2007-ts_FP}
Green, K.~C., Armstrong, J.~S., 2007. Structured analogies for forecasting.
  International Journal of Forecasting 23~(3), 365--376.

\bibitem[{Green et~al.(2006)Green, Blandford, Church, Roast, and
  Clarke}]{Green2006-cm_MJE}
Green, T. R.~G., Blandford, A.~E., Church, L., Roast, C.~R., Clarke, S., 2006.
  Cognitive dimensions: Achievements, new directions, and open questions.
  Journal of Visual Languages \& Computing 17~(4), 328--365.

\bibitem[{Greenberg et~al.(2020)Greenberg, Cox, Bier, Lambert, Lowrie, North,
  Siegrist, and Wu}]{Greenberg2020-xl_TC}
Greenberg, M., Cox, A., Bier, V., Lambert, J., Lowrie, K., North, W., Siegrist,
  M., Wu, F., 2020. Risk analysis: Celebrating the accomplishments and
  embracing ongoing challenges. Risk Analysis 40~(S1), 2113--2127.

\bibitem[{Greenberg and Cox(2021)}]{Greenberg2021-vu_TC}
Greenberg, M., Cox, Jr, L.~A., 2021. Plutonium disposition: Using and
  explaining complex {Risk-Related} methods. Risk Analysis 41~(12), 2186--2195.

\bibitem[{Gregory and Ronan(2015)}]{Gregory2015-co_AJG}
Gregory, A., Ronan, M., 2015. Insights into the development of strategy from a
  complexity perspective. Journal of the Operational Research Society 66~(4),
  627--636.

\bibitem[{Gregory et~al.(2020)Gregory, Atkins, Midgley, and
  Hodgson}]{Gregory2020-zf_JL}
Gregory, A.~J., Atkins, J.~P., Midgley, G., Hodgson, A.~M., 2020. Stakeholder
  identification and engagement in problem structuring interventions. European
  Journal of Operational Research 283~(1), 321--340.

\bibitem[{Gregory and Jackson(1992{\natexlab{a}})}]{Gregory1992-ab_AJG}
Gregory, A.~J., Jackson, M.~C., 1992{\natexlab{a}}. Evaluating organizations: A
  systems and contingency approach. Systems Practice 5~(1), 37--60.

\bibitem[{Gregory and Jackson(1992{\natexlab{b}})}]{Gregory1992-yb_AJG}
Gregory, A.~J., Jackson, M.~C., 1992{\natexlab{b}}. Evaluation methodologies: A
  system for use. Journal of the Operational Research Society 43~(1), 19--28.

\bibitem[{Gregory et~al.(2012)Gregory, Failing, Harstone, Long, McDaniels, and
  Ohlson}]{Gregory2012-zh_JL}
Gregory, R., Failing, L., Harstone, M., Long, G., McDaniels, T., Ohlson, D.,
  2012. Structured Decision Making: A Practical Guide to Environmental
  Management Choices. Wiley.

\bibitem[{Gregory(2000)}]{Gregory2000-mo_GM}
Gregory, W.~J., 2000. Transforming self and society: A ``critical
  appreciation'' model. Systemic Practice and Action Research 13~(4), 475--501.

\bibitem[{Grieco et~al.(2021)Grieco, Utley, and Crowe}]{Grieco2021-le_CV}
Grieco, L., Utley, M., Crowe, S., 2021. Operational research applied to
  decisions in home health care: A systematic literature review. Journal of the
  Operational Research Society 72~(9), 1960--1991.

\bibitem[{Grinsztajn et~al.(2022)Grinsztajn, Oyallon, and
  Gaël}]{NEURIPS2022_Varoquaux_LCAL}
Grinsztajn, L., Oyallon, E., Gaël, V., 2022. Why do tree-based models still
  outperform deep learning on typical tabular data? In: Advances in Neural
  Information Processing Systems. Vol.~35. Curran Associates, Inc.

\bibitem[{Groop et~al.(2017)Groop, Ketokivi, Gupta, and
  Holmström}]{groop_improving_2017_JHLS}
Groop, J., Ketokivi, M., Gupta, M., Holmström, J., 2017. Improving home care:
  {Knowledge} creation through engagement and design. Journal of Operations
  Management 53-56, 9--22.

\bibitem[{Gross et~al.(2019)Gross, Brunner, and Blobner}]{fairness_GVBSP}
Gross, C.~N., Brunner, J.~O., Blobner, M., 2019. Hospital physicians can’t
  get no long-term satisfaction--an indicator for fairness in preference
  fulfillment on duty schedules. Health Care Management Science 22~(4),
  691--708.

\bibitem[{Gross and Harris(1974)}]{grossharris_HATI}
Gross, D., Harris, C.~M., 1974. Fundamentals of Queueing Theory. Wiley, New
  York.

\bibitem[{Gr{\"o}tschel et~al.(1981)Gr{\"o}tschel, Lov{\'a}sz, and
  Schrijver}]{GLS81_ALAL}
Gr{\"o}tschel, M., Lov{\'a}sz, L., Schrijver, A., 1981. The ellipsoid method
  and its consequences in combinatorial optimization. Combinatorica 1,
  169--197.

\bibitem[{Gr{\"o}tschel et~al.(1988)Gr{\"o}tschel, Lov{\'a}sz, and
  Schrijver}]{Grotschel_1988-ic_JMB}
Gr{\"o}tschel, M., Lov{\'a}sz, L., Schrijver, A., 1988. Geometric Algorithms
  and Combinatorial Optimization. Springer, Berlin.

\bibitem[{Gr{\"{o}}tschel and Monma(1990)}]{grotschel.monma:90_BF}
Gr{\"{o}}tschel, M., Monma, C., 1990. Integer polyhedra arising from certain
  design problems with connectivity constraints. SIAM Journal on Discrete
  Mathematics 3, 502--523.

\bibitem[{Grubbstr{\"o}m(1967)}]{grubbstrom1967_XW}
Grubbstr{\"o}m, R.~W., 1967. On the application of the {L}aplace transform to
  certain economic problems. Management Science 13~(7), 558--567.

\bibitem[{Grus(2019)}]{grus2019data_LCAL}
Grus, J., 2019. Data Science from Scratch: First Principles with Python.
  O'Reilly Media, Sebastopol, CA.

\bibitem[{Gschwind and Irnich(2016)}]{Gschwind2016_CA_MB}
Gschwind, T., Irnich, S., 2016. Dual inequalities for stabilized column
  generation revisited. INFORMS Journal on Computing 28~(1), 175--194.

\bibitem[{Gu et~al.(2020)Gu, Kelly, and Xiu}]{gu20_MB}
Gu, S., Kelly, B., Xiu, D., 2020. Empirical asset pricing via machine learning.
  The Review of Financial Studies 33~(5), 2223--2273.

\bibitem[{Gu et~al.(1998)Gu, Nemhauser, and Savelsbergh}]{GNS98_ALAL}
Gu, Z., Nemhauser, G., Savelsbergh, M., 1998. Lifted cover inequalities for 0-1
  integer programs: computation. {INFORMS} Journal on Computing 10, 427--437.

\bibitem[{Gu et~al.(1999)Gu, Nemhauser, and Savelsbergh}]{gu1999lifted_MH}
Gu, Z., Nemhauser, G.~L., Savelsbergh, M.~W., 1999. Lifted flow cover
  inequalities for mixed 0-1 integer programs. Mathematical Programming 85~(3),
  439--467.

\bibitem[{Guide et~al.(2006)Guide, Souza, Van~Wassenhove, and
  Blackburn}]{Guide2006-eo_AM}
Guide, V. D.~R., Souza, G.~C., Van~Wassenhove, L.~N., Blackburn, J., 2006. Time
  value of commercial product returns. Management Science 52~(8), 1200--1214.

\bibitem[{Guide et~al.(2003)Guide, Teunter, and
  Van~Wassenhove}]{Guide2003-nv_AM}
Guide, V. D.~R., Teunter, R.~H., Van~Wassenhove, L.~N., 2003. Matching demand
  and supply to maximize profits from remanufacturing. Manufacturing \& Service
  Operations Management 5~(4), 303--316.

\bibitem[{Guide and Van~Wassenhove(2003)}]{Guide2003-qd_AM}
Guide, V. D.~R., Van~Wassenhove, L.~N., 2003. Business Aspects of Closed-loop
  Supply Chains. Carnegie Mellon University Press, Pittsburgh, PA, USA.

\bibitem[{Guide and Li(2010)}]{Guide2010-ey_AM}
Guide, Jr, V. D.~R., Li, J., 2010. The potential for cannibalization of new
  products sales by remanufactured products. Decision Sciences 41~(3),
  547--572.

\bibitem[{Guihaire and Hao(2008)}]{Guihaire2008_DC}
Guihaire, V., Hao, J.-K., 2008. Transit network design and scheduling{: A}
  global review. Transportation Research Part A: Policy and Practice 42~(10),
  1251--1273.

\bibitem[{Guikema(2020)}]{Guikema2020-pw_TC}
Guikema, S., 2020. Artificial intelligence for natural hazards risk analysis:
  Potential, challenges, and research needs. Risk Analysis 40~(6), 1117--1123.

\bibitem[{Guo et~al.(2012)Guo, Xiao, and
  Li}]{guoUnconstrainingMethodsRevenue2012_AKSJF}
Guo, P., Xiao, B., Li, J., 2012. Unconstraining methods in revenue management
  systems: Research overview and prospects. Advances in Operations Research
  2012.

\bibitem[{{Gurobi}(2022)}]{Gurobi_CTCGE}
{Gurobi}, 2022. {Gurobi}.
\newline\urlprefix\url{https://www.gurobi.com/}

\bibitem[{Gusfield(1990)}]{Gusfield:1990_IL}
Gusfield, D., 1990. Very simple methods for all pairs network flow analysis.
  {SIAM} Journal on Computing 19~(1), 143--155.

\bibitem[{Guti\'{e}rrez-Jarpa et~al.(2013)Guti\'{e}rrez-Jarpa, Obreque,
  Laporte, and Marianov}]{GutierrezJarpa2013_DC}
Guti\'{e}rrez-Jarpa, G., Obreque, C., Laporte, G., Marianov, V., 2013. Rapid
  transit network design for optimal cost and origin-destination demand
  capture. Computers \& {O}perations {R}esearch 40, 3000--3009.

\bibitem[{Gutin and Punnen(2006)}]{GP06_SMPT}
Gutin, G., Punnen, A. (Eds.), 2006. The traveling salesman problem and its
  variations. Kluwer, Dordrecht.

\bibitem[{Haase and M{\"u}ller(2013)}]{haase2013management_SAA}
Haase, K., M{\"u}ller, S., 2013. Management of school locations allowing for
  free school choice. Omega 41~(5), 847--855.

\bibitem[{Habermas(1972)}]{Habermas1972-yr_AJG}
Habermas, J., 1972. Knowledge and Human Interests. Heinemann, London,
  translated by J.J. Shapiro.

\bibitem[{Hadley and Within(1963)}]{hadley1963within_JSS}
Hadley, G., Within, T., 1963. Analysis of Inventory Systems. Prentice-Hall.

\bibitem[{Haeringer(2018)}]{Haeringer18_BC}
Haeringer, G., 2018. Market Design: Auctions and Matching. MIT Press,
  Cambridge, MA.

\bibitem[{Hahler and Fleischmann(2017)}]{Hahler2017-qt_AM}
Hahler, S., Fleischmann, M., 2017. Strategic grading in the product acquisition
  process of a reverse supply chain. Production and Operations Management
  26~(8), 1498--1511.

\bibitem[{Hakes and Sauer(2006)}]{Hakes2006-tt_IM}
Hakes, J.~K., Sauer, R.~D., 2006. An economic evaluation of the moneyball
  hypothesis. The Journal of Economic Perspectives 20~(3), 173--186.

\bibitem[{Hakimi(1965)}]{hakimi1965optimum_SAA}
Hakimi, S.~L., 1965. Optimum distribution of switching centers in a
  communication network and some related graph theoretic problems. Operations
  Research 13~(3), 462--475.

\bibitem[{Halfin and Whitt(1981)}]{halfin1981heavy_HATI}
Halfin, S., Whitt, W., 1981. Heavy-traffic limits for queues with many
  exponential servers. Operations Research 29~(3), 567--588.

\bibitem[{Halkos et~al.(2014)Halkos, Tzeremes, and
  Kourtzidis}]{Halkos2014-de_SL}
Halkos, G.~E., Tzeremes, N.~G., Kourtzidis, S.~A., 2014. A unified
  classification of two-stage {DEA} models. Surveys in Operations Research and
  Management Science 19~(1), 1--16.

\bibitem[{Hall(1962)}]{Hall1962-by_GM}
Hall, A.~D., 1962. A Methodology for Systems Engineering. Van Nostrand.

\bibitem[{Halvorsen-Weare and Fagerholt(2011)}]{Halvorsen-Weare2011-fm_HP}
Halvorsen-Weare, E.~E., Fagerholt, K., 2011. Robust supply vessel planning. In:
  Network Optimization. Springer, Berlin, Heidelberg, pp. 559--573.

\bibitem[{H{\"a}m{\"a}l{\"a}inen(2015)}]{Hamalainen2015-fv_JL}
H{\"a}m{\"a}l{\"a}inen, R.~P., 2015. Behavioural issues in environmental
  modelling -- the missing perspective. Environmental Modelling \& Software 73,
  244--253.

\bibitem[{H{\"a}m{\"a}l{\"a}inen and Lahtinen(2016)}]{Hamalainen2016-gr_AFRH}
H{\"a}m{\"a}l{\"a}inen, R.~P., Lahtinen, T.~J., 2016. Path dependence in
  {O}perational {Research}--{How} the modeling process can influence the
  results. Operations Research Perspectives 3, 14--20.

\bibitem[{H{\"a}m{\"a}l{\"a}inen et~al.(2013)H{\"a}m{\"a}l{\"a}inen, Luoma, and
  Saarinen}]{Hamalainen2013-zr_AFRH}
H{\"a}m{\"a}l{\"a}inen, R.~P., Luoma, J., Saarinen, E., 2013. On the importance
  of behavioral operational research: The case of understanding and
  communicating about dynamic systems. European Journal of Operational Research
  228~(3), 623--634.

\bibitem[{Hammant et~al.(1999)Hammant, Disney, Childerhouse, and
  Naim}]{Hammant1999_SMD}
Hammant, J., Disney, S.~M., Childerhouse, P., Naim, M.~M., 1999. Modelling the
  consequences of a strategic supply chain initiative of an automotive
  aftermarket operation. International Journal of Physical Distribution \&
  Logistics Management 29~(9), 535--550.

\bibitem[{Hammond(1992)}]{Hammond1992-gv_TC}
Hammond, P.~J., 1992. Harsanyi's utilitarian theorem: A simpler proof and some
  ethical connotations. In: Selten, R. (Ed.), Rational Interaction: Essays in
  Honor of John C. Harsanyi. Springer, Berlin, Heidelberg, pp. 305--319.

\bibitem[{Hane et~al.(1995)Hane, Barnhart, Johnson, Marsten, Nemhauser, and
  Sigismondi}]{hane1995fleet_VLVV}
Hane, C.~A., Barnhart, C., Johnson, E.~L., Marsten, R.~E., Nemhauser, G.~L.,
  Sigismondi, G., 1995. The fleet assignment problem: Solving a large-scale
  integer program. Mathematical Programming 70~(1), 211--232.

\bibitem[{Hanea et~al.(2022)Hanea, Christophersen, and Alday}]{Hanea2022-ls_TC}
Hanea, A.~M., Christophersen, A., Alday, S., 2022. Bayesian networks for risk
  analysis and decision support. Risk Analysis 42~(6), 1149--1154.

\bibitem[{Hansen and Mladenovi{\'c}(1999)}]{hansen1999introduction_COIT}
Hansen, P., Mladenovi{\'c}, N., 1999. An introduction to variable neighborhood
  search. In: Vo{\ss}, S., Martello, S., Osman, I.~H., Roucairol, C. (Eds.),
  {Meta-Heuristics}: Advances and Trends in Local Search Paradigms for
  Optimization. Springer, pp. 433--458,
  \url{https://doi.org/10.1007/978-1-4615-5775-3_30}.

\bibitem[{Hao and Orlin(1994)}]{Hao-Orlin:1994_IL}
Hao, J., Orlin, J.~B., 1994. A faster algorithm for finding the minimum cut in
  a directed graph. Journal of Algorithms 17~(3), 424--446.

\bibitem[{Harchol-Balter(2013)}]{mor_HATI}
Harchol-Balter, M., 2013. Performance Modeling and Design of Computer Systems:
  Queueing Theory in Action. Cambridge University Press, New York.

\bibitem[{Hardt et~al.(2016)Hardt, Price, and Srebro}]{HarPriSre16_JH}
Hardt, M., Price, E., Srebro, N., 2016. Equality of opportunity in supervised
  learning. In: Proceedings of 30th International Conference on Neural
  Information Processing. pp. 3323--3331.

\bibitem[{Hargreaves et~al.(2022)Hargreaves, Langan, Oswald, Halliday,
  Sturgess, Phelan, Nguipdop-Djomo, Ford, Allen, Sundaram, Ireland, Poh, Ijaz,
  Diamond, Rourke, Dawe, Judd, Warren-Gash, Clark, Glynn, Edmunds, Bonell,
  Mangtani, Ladhani, and {COVID-19 Schools Infection Survey Study
  Group}}]{Hargreaves2022-yt}
Hargreaves, J.~R., Langan, S.~M., Oswald, W.~E., Halliday, K.~E., Sturgess, J.,
  Phelan, J., Nguipdop-Djomo, P., Ford, B., Allen, E., Sundaram, N., Ireland,
  G., Poh, J., Ijaz, S., Diamond, I., Rourke, E., Dawe, F., Judd, A.,
  Warren-Gash, C., Clark, T.~G., Glynn, J.~R., Edmunds, W.~J., Bonell, C.,
  Mangtani, P., Ladhani, S.~N., {COVID-19 Schools Infection Survey Study
  Group}, 2022. Epidemiology of {SARS-CoV-2} infection among staff and students
  in a cohort of english primary and secondary schools during 2020-2021. The
  Lancet regional health. Europe 21, 100471.

\bibitem[{Harju et~al.(2019)Harju, Liesi{\"o}, and
  Virtanen}]{Harju2019-mb_KVRH}
Harju, M., Liesi{\"o}, J., Virtanen, K., 2019. Spatial multi-attribute decision
  analysis: Axiomatic foundations and incomplete preference information.
  European Journal of Operational Research 275~(1), 167--181.

\bibitem[{Harper et~al.(2021)Harper, Pitt, De~Prez, Dumlu, Vasilakis, Forte,
  and Wood}]{Harper2021-qk_CV}
Harper, A., Pitt, M., De~Prez, M., Dumlu, Z.~{\"O}., Vasilakis, C., Forte, P.,
  Wood, R., 2021. A demand and capacity model for {Home-Based} intermediate
  care: Optimizing the `step down' pathway. In: 2021 Winter Simulation
  Conference ({WSC}). pp. 1--12.

\bibitem[{Harris(1913)}]{harris1913many_JSS}
Harris, F.~W., 1913. How many parts to make at once. Factory, The Magazine of
  Management 10~(2), 135--136, 152.

\bibitem[{Harris(1990)}]{harris1990many_JSS}
Harris, F.~W., 1990. How many parts to make at once. Operations Research
  38~(6), 947--950.

\bibitem[{Harrison(1985)}]{harrison1985brownian_HATI}
Harrison, J.~M., 1985. Brownian motion and stochastic flow systems. Wiley New
  York.

\bibitem[{Harsanyi(1967)}]{Ha1967_GZ_1}
Harsanyi, J., 1967. Games with incomplete information played by
  {``{B}ayesian''} players, {Part I}. Management Science 14~(3), 159--182.

\bibitem[{Harsanyi(1968{\natexlab{a}})}]{Ha1967_GZ_2}
Harsanyi, J., 1968{\natexlab{a}}. Games with incomplete information played by
  {``{B}ayesian''} players, {Part II}. Management Science 14~(5), 320--334.

\bibitem[{Harsanyi(1968{\natexlab{b}})}]{Ha1967_GZ_3}
Harsanyi, J., 1968{\natexlab{b}}. Games with incomplete information played by
  {``{B}ayesian''} players, {Part III}. Management Science 14~(7), 486--502.

\bibitem[{Harsanyi(1977)}]{Har77_JH}
Harsanyi, J.~C., 1977. Rational Behavior and Bargaining Equilibrium in Games
  and Social Situations. Cambridge University Press.

\bibitem[{Hartert et~al.(2015)Hartert, Vissicchio, Schaus, Bonaventure,
  Filsfils, Telkamp, and Francois}]{hartert.vissicchio.ea:15_BF}
Hartert, R., Vissicchio, S., Schaus, P., Bonaventure, O., Filsfils, C.,
  Telkamp, T., Francois, P., 2015. A declarative and expressive approach to
  control forwarding paths in carrier-grade networks. ACM SIGCOMM Computer
  Communication Review 45~(4), 15--28.

\bibitem[{Hartmann(2001)}]{Hartmann2001-bu_WH_ED}
Hartmann, S., 2001. Project scheduling with multiple modes: A genetic
  algorithm. Annals of Operations Research 102~(1), 111--135.

\bibitem[{Hartmann(2002)}]{Hartmann2002-vg_WH_ED}
Hartmann, S., 2002. A self‐adapting genetic algorithm for project scheduling
  under resource constraints. Naval Research Logistics 49~(5), 433--448.

\bibitem[{Hartmann and Briskorn(2010)}]{Hartmann2010-dc_WH_ED}
Hartmann, S., Briskorn, D., 2010. A survey of variants and extensions of the
  resource-constrained project scheduling problem. European Journal of
  Operational Research 207~(1), 1--14.

\bibitem[{Hartmann and Drexl(1998)}]{Hartmann1998-ny_WH_ED}
Hartmann, S., Drexl, A., 1998. Project scheduling with multiple modes: A
  comparison of exact algorithms. Networks. An International Journal 32~(4),
  283--297.

\bibitem[{Harwood(2019)}]{Harwood2019-kp_MYLW}
Harwood, S., 2019. Whither is problem structuring methods ({PSMs)}? Journal of
  the Operational Research Society 70~(8), 1391--1392.

\bibitem[{Hasbrouck(1995)}]{hasbrouck95_FP}
Hasbrouck, J., 1995. One security, many markets: Determining the contributions
  to price discovery. Journal of Finance 50~(4), 1175--1199.

\bibitem[{Haspeslagh et~al.(2014)Haspeslagh, {De Causmaecker}, Schaerf, and
  Stølevik}]{INRCI_GVBSP}
Haspeslagh, S., {De Causmaecker}, P., Schaerf, A., Stølevik, M., 2014. The
  first international nurse rostering competition 2010. Annals of Operations
  Research 218, 221--236.

\bibitem[{Hastie et~al.(2015)Hastie, Tibshirani, and
  Wainwright}]{has:tib:wai:15_DSRW}
Hastie, T., Tibshirani, R., Wainwright, M., 2015. Statistical Learning with
  Sparsity: The Lasso and Generalizations. CRC Press, Boca Raton, FL.

\bibitem[{Hattori et~al.(2005)Hattori, Jamasb, and Pollitt}]{Hattori_2005_DSRW}
Hattori, T., Jamasb, T., Pollitt, M., 2005. Electricity distribution in the
  {UK} and {J}apan: {A} comparative efficiency analysis 1985--1998. The Energy
  Journal 26~(2), 23--47.

\bibitem[{Haug and Blackburn(2017)}]{Haug2017-bv_JJ}
Haug, A.~A., Blackburn, V.~C., 2017. Government secondary school finances in
  {New South Wales}: accounting for students' prior achievements in a two-stage
  {DEA} at the school level. Journal of Productivity Analysis 48~(1), 69--83.

\bibitem[{Haurie and Pohjola(1987)}]{Hapo1987_GZ}
Haurie, A., Pohjola, M., 1987. Efficient equilibria in a differential game of
  capitalism. Journal of Economic Dynamics and Control 11~(1), 65--78.

\bibitem[{Haurie and Tolwinski(1985)}]{Hato1985_GZ}
Haurie, A., Tolwinski, B., 1985. Definition and properties of cooperative
  equilibria in a two-player game of infinite duration. Journal of Optimization
  Theory and Applications 46, 525--534.

\bibitem[{Hausken(2018)}]{Hausken2018-gv_KVRH}
Hausken, K., 2018. A cost--benefit analysis of terrorist attacks. Defence and
  Peace Economics 29~(2), 111--129.

\bibitem[{Hax and Meal(1973)}]{hax1973hierarchical_COIT}
Hax, A.~C., Meal, H.~C., 1973. Hierarchical integration of production planning
  and scheduling. Massachusetts Institute of Technology, Sloan School of
  Management 656-73.

\bibitem[{He et~al.(2018)He, Yang, and Li}]{he2018vehicle_DLLD}
He, F., Yang, J., Li, M., 2018. Vehicle scheduling under stochastic trip times:
  {A}n approximate dynamic programming approach. Transportation Research Part
  C: Emerging Technologies 96, 144--159.

\bibitem[{Helsgaun(2000)}]{Helsgaun:2000_IL}
Helsgaun, K., 2000. An effective implementation of the {L}in-{K}ernighan
  traveling salesman heuristic. European Journal of Operational Research
  126~(1), 106--130.

\bibitem[{Hemmelmayr et~al.(2009)Hemmelmayr, Doerner, Hartl, and
  Savelsbergh}]{Hemmel2009_JLYHK}
Hemmelmayr, V., Doerner, K.~F., Hartl, R.~F., Savelsbergh, M. W.~P., 2009.
  Delivery strategies for blood product supplies. OR Spectrum 31, 707--725.

\bibitem[{Hermans and Thissen(2009)}]{Hermans2009-rp_JL}
Hermans, L.~M., Thissen, W. A.~H., 2009. Actor analysis methods and their use
  for public policy analysts. European Journal of Operational Research 196~(2),
  808--818.

\bibitem[{Herrmann(2012)}]{Herrmann2012-cg_KVRH}
Herrmann, J., 2012. Handbook of Operations Research for Homeland Security.
  Springer, New York, NY.

\bibitem[{Herroelen(2007)}]{Herroelen2007-hh_WH_ED}
Herroelen, W., 2007. Generating robust project baseline schedules. In:
  Klastorin, T. (Ed.), {OR} Tools and Applications: Glimpses of Future
  Technologies. INFORMS TutORials in Operations Research. INFORMS, pp.
  124--144.

\bibitem[{Herroelen and Leus(2001)}]{Herroelen2001-vl_WH_ED}
Herroelen, W., Leus, R., 2001. On the merits and pitfalls of critical chain
  scheduling. Journal of Operations Management 19~(5), 559--577.

\bibitem[{Herroelen and Leus(2005)}]{Herroelen2005-me_WH_ED}
Herroelen, W., Leus, R., 2005. Project scheduling under uncertainty: Survey and
  research potentials. European Journal of Operational Research 165~(2),
  289--306.

\bibitem[{Herron and Mendiwelso-Bendek(2018)}]{Herron2018-ff_AJG}
Herron, R., Mendiwelso-Bendek, Z., 2018. Supporting self-organised community
  research through informal learning. European Journal of Operational Research
  268~(3), 825--835.

\bibitem[{Hewitt(2019)}]{EnhancedDDD_MH}
Hewitt, M., 2019. Enhanced dynamic discretization discovery for the continuous
  time load plan design problem. Transportation Science 53~(6), 1731--1750.

\bibitem[{Hewitt et~al.(2019)Hewitt, Crainic, Nowak, and
  Rei}]{HEWITT2019324_MH}
Hewitt, M., Crainic, T.~G., Nowak, M., Rei, W., 2019. Scheduled service network
  design with resource acquisition and management under uncertainty.
  Transportation Research Part B: Methodological 128, 324--343.

\bibitem[{Hewitt et~al.(2010)Hewitt, Nemhauser, and
  Savelsbergh}]{hewitt2010combining_MH}
Hewitt, M., Nemhauser, G.~L., Savelsbergh, M.~W., 2010. Combining exact and
  heuristic approaches for the capacitated fixed-charge network flow problem.
  INFORMS Journal on Computing 22~(2), 314--325.

\bibitem[{Hewitt et~al.(2021)Hewitt, Rei, and Wallace}]{Hewitt2021SND_MH}
Hewitt, M., Rei, W., Wallace, S.~W., 2021. Stochastic network design. In:
  Crainic, T.~G., Gendreau, M., Gendron, B. (Eds.), Network Design with
  Applications to Transportation and Logistics. Springer, pp. 283--315.

\bibitem[{Heydenreich et~al.(2009)Heydenreich, Müller, Uetz, and
  Vohra}]{Heydenreich09_BC}
Heydenreich, B., Müller, R., Uetz, M., Vohra, R.~V., 2009. Characterization of
  revenue equivalence. Econometrica 77~(1), 307--316.

\bibitem[{Hickman(2010)}]{Hickman08_BC}
Hickman, B.~R., 2010. On the pricing rule in electronic auctions. International
  Journal of Industrial Organization 28~(5), 423--433.

\bibitem[{Hifi and M'Hallah(2009)}]{Hifi2009-ce_JB}
Hifi, M., M'Hallah, R., 2009. A literature review on circle and sphere packing
  problems: Models and methodologies. Advances in Operations Research 2009,
  150624.

\bibitem[{Higgins et~al.(2010)Higgins, Miller, Archer, Ton, Fletcher, and
  McAllister}]{Higgins2010_EA}
Higgins, A.~J., Miller, C.~J., Archer, A.~A., Ton, T., Fletcher, C.~S.,
  McAllister, R.~R., 2010. Challenges of operations research practice in
  agricultural value chains. Journal of the Operational Research Society
  61~(6), 964--973.

\bibitem[{Higle and Sen(1991)}]{Higle1991-lz_HL}
Higle, J.~L., Sen, S., 1991. Stochastic decomposition: An algorithm for
  {Two-Stage} linear programs with recourse. Mathematics of Operations Research
  16~(3), 650--669.

\bibitem[{Hill(2020)}]{Hill2020-sv_KVRH}
Hill, R.~R., 2020. Modern data analytics for the military operational
  researcher. In: Scala, N., Howard, J. (Eds.), Handbook of military and
  defense operations research. CRC Press, Boca Raton, FL, pp. 3--18.

\bibitem[{Hindle and Vidgen(2018)}]{Hindle2018-zo_JEB}
Hindle, G.~A., Vidgen, R., 2018. Developing a business analytics methodology: A
  case study in the foodbank sector. European Journal of Operational Research
  268~(3), 836--851.

\bibitem[{Hines et~al.(2004)Hines, Holweg, and Rich}]{Hines2004-fo_KESXZ}
Hines, P., Holweg, M., Rich, N., 2004. Learning to evolve: A review of
  contemporary lean thinking. International Journal of Operations \& Production
  Management 24~(10), 994--1011.

\bibitem[{Hines and Rich(1997)}]{Hines1997_SMD}
Hines, P., Rich, N., 1997. The seven value stream mapping tools. International
  Journal of Operations \& Production Management 17, 46--64.

\bibitem[{Hitchcock(1941)}]{Hitchcock1941-mc_JMB}
Hitchcock, F.~L., 1941. The distribution of a product from several sources to
  numerous localities. Journal of Mathematics and Physics 20~(1-4), 224--230.

\bibitem[{Hjorts{\o}(2004)}]{Hjortso2004-aa_MYLW}
Hjorts{\o}, C.~N., 2004. Enhancing public participation in natural resource
  management using soft {OR----an} application of strategic option development
  and analysis in tactical forest planning. European Journal of Operational
  Research 152~(3), 667--683.

\bibitem[{Hochbaum(1997)}]{Hoch96_UPCT}
Hochbaum, D.~S. (Ed.), 1997. Approximation Algorithms for {NP}-Hard Problems.
  PWS Publishing Co., Boston, MA.

\bibitem[{Hochbaum(2008)}]{Hochbaum:2008_IL}
Hochbaum, D.~S., 2008. The pseudoflow algorithm: {A} new algorithm for the
  maximum-flow problem. Operations Research 56~(4), 992--1009.

\bibitem[{Hochbaum and Orlin(2013)}]{Hochbaum-Orlin:2013_IL}
Hochbaum, D.~S., Orlin, J.~B., 2013. Simplifications and speedups of the
  pseudoflow algorithm. Networks 61~(1), 40--57.

\bibitem[{Hodson and Hill(2014)}]{Hodson2014-ft_KVRH}
Hodson, D.~D., Hill, R.~R., 2014. The art and science of live, virtual, and
  constructive simulation for test and analysis. The Journal of Defense
  Modeling and Simulation 11~(2), 77--89.

\bibitem[{Hofbauer and Sigmund(1998)}]{Hosi1998_GZ}
Hofbauer, J., Sigmund, K., 1998. Evolutionary games and population dynamics.
  Cambridge University Press, Cambridge.

\bibitem[{Hoffman and Padberg(1991)}]{HP91_ALAL}
Hoffman, K., Padberg, M., 1991. Improving {LP}-representations of zero-one
  linear programs for branch-and-cut. {ORSA} Journal on Computing 3, 121--134.

\bibitem[{Hoffman et~al.(2013)Hoffman, Blei, Wang, and
  Paisley}]{hoffman13_LCAL}
Hoffman, M.~D., Blei, D.~M., Wang, C., Paisley, J., 2013. Stochastic
  variational inference. Journal of Machine Learning Research 14~(1),
  1303–1347.

\bibitem[{Holgu{\'\i}n-Veras et~al.(2013)Holgu{\'\i}n-Veras, P{\'e}rez, Jaller,
  Van~Wassenhove, and Aros-Vera}]{holguin2013appropriate_BYKOK}
Holgu{\'\i}n-Veras, J., P{\'e}rez, N., Jaller, M., Van~Wassenhove, L.~N.,
  Aros-Vera, F., 2013. On the appropriate objective function for post-disaster
  humanitarian logistics models. Journal of Operations Management 31~(5),
  262--280.

\bibitem[{Holland(1975)}]{Holland1975-qj_HL}
Holland, J.~H., 1975. Adaptation in Natural and Artificial Systems: An
  Introductory Analysis with Applications to Biology, Control, and Artificial
  Intelligence. University of Michigan Press, Ann Arbor, MI.

\bibitem[{Hollyman et~al.(2021)Hollyman, Petropoulos, and
  Tipping}]{Hollyman2021-wj_FP}
Hollyman, R., Petropoulos, F., Tipping, M.~E., 2021. Understanding forecast
  reconciliation. European Journal of Operational Research 294~(1), 149--160.

\bibitem[{Holmström et~al.(2010)Holmström, Brax, and
  Ala‐Risku}]{holmstrom_comparing_2010_JHLS}
Holmström, J., Brax, S., Ala‐Risku, T., 2010. Comparing provider‐customer
  constellations of visibility‐based service. Journal of Service Management
  21~(5), 675--692.

\bibitem[{Holmström et~al.(2009)Holmström, Ketokivi, and
  Hameri}]{holmstrom_bridging_2009_JHLS}
Holmström, J., Ketokivi, M., Hameri, A.-P., 2009. Bridging {Practice} and
  {Theory}: {A} {Design} {Science} {Approach}. Decision Sciences 40~(1),
  65--87.

\bibitem[{Holsapple et~al.(2014)Holsapple, Lee-Post, and
  Pakath}]{Holsapple2014-gf_JEB}
Holsapple, C., Lee-Post, A., Pakath, R., 2014. A unified foundation for
  business analytics. Decision Support Systems 64, 130--141.

\bibitem[{Holt(2004)}]{Holt20045_FP}
Holt, C.~C., 2004. Forecasting seasonals and trends by exponentially weighted
  moving averages. International Journal of Forecasting 20~(1), 5--10.

\bibitem[{Holweg(2007)}]{Holweg2007_SMD}
Holweg, M., 2007. The genealogy of lean production. Journal of Operations
  Management 25~(2), 420--437.

\bibitem[{Holweg et~al.(2005)Holweg, Disney, Holmström, and
  Småros}]{Holweg2005_SMD}
Holweg, M., Disney, S.~M., Holmström, J., Småros, J., 2005. Supply chain
  collaboration: Making sense of the strategy continuum. European Management
  Journal 23~(2), 170--181.

\bibitem[{Holzmann and Smith(2021)}]{Holzmann2021-ib_KVRH}
Holzmann, T., Smith, J.~C., 2021. The shortest path interdiction problem with
  randomized interdiction strategies: Complexity and algorithms. Operations
  Research 69~(1), 82--99.

\bibitem[{Hong and Nelson(2009)}]{HongNelsonWSC2009_CC}
Hong, L.~J., Nelson, B.~L., 2009. A brief introduction to optimization via
  simulation. In: Rossetti, M.~D., Hill, R.~R., Johansson, B., Dunkin, A.,
  Ingalls, R.~G. (Eds.), Proceedings of the 2009 Winter Simulation Conference.
  IEEE Piscataway, pp. 75--85.

\bibitem[{Hong(2014)}]{hon:14_DSRW}
Hong, T., 2014. Energy forecasting: Past, present, and future. Foresight 32,
  43--48.

\bibitem[{Hong and Fan(2016)}]{hon:fan:16_DSRW}
Hong, T., Fan, S., 2016. Probabilistic electric load forecasting: {A} tutorial
  review. International Journal of Forecasting 32, 914--938.

\bibitem[{Hong et~al.(2020)Hong, Pinson, Wang, Weron, Yang, and
  Zareipour}]{hon:pin:etal:20_DSRW}
Hong, T., Pinson, P., Wang, Y., Weron, R., Yang, D., Zareipour, H., 2020.
  Energy forecasting: A review and outlook. IEEE Open Access Journal of Power
  and Energy 7, 376--388.

\bibitem[{Hooghiemstra et~al.(1999)Hooghiemstra, Kroon, Odijk, Salomon, and
  Zwaneveld}]{hooghiemstra1999decision_MH}
Hooghiemstra, J.~S., Kroon, L.~G., Odijk, M.~A., Salomon, M., Zwaneveld, P.~J.,
  1999. Decision support systems support the search for win-win solutions in
  railway network design. Interfaces 29~(2), 15--32.

\bibitem[{Hooker and Williams(2012)}]{HooWil12_JH}
Hooker, J.~N., Williams, H.~P., 2012. Combining equity and utilitarianism in a
  mathematical programming model. Management Science 58, 1682--1693.

\bibitem[{Hoos(1972)}]{Hoos1972-cq_GM}
Hoos, I.~R., 1972. Systems Analysis in Public Policy: A Critique. University of
  California Press.

\bibitem[{Hoot et~al.(2008)Hoot, LeBlanc, Jones, Levin, Zhou, Gadd, and
  Aronsky}]{Hoot2008_CC}
Hoot, N.~R., LeBlanc, L.~J., Jones, I., Levin, S.~R., Zhou, C., Gadd, C.~S.,
  Aronsky, D., 2008. Forecasting emergency department crowding: A discrete
  event simulation. Annals of Emergency Medicine 52~(2), 116--125.

\bibitem[{Hoover(1936)}]{Hoov36_JH}
Hoover, E.~M., 1936. The measurement of industrial localization. Review of
  Economics and Statistics 18, 162--171.

\bibitem[{Hopper and Turton(2001)}]{Hopper2001-dz_JB}
Hopper, E., Turton, B. C.~H., 2001. An empirical investigation of
  meta-heuristic and heuristic algorithms for a {2D} packing problem. European
  Journal of Operational Research 128~(1), 34--57.

\bibitem[{Howard(1966)}]{howarddecision_MESG}
Howard, R., 1966. Decision analysis: Applied decision theory. In: Hertz,
  D.B.and~Melese, J. (Ed.), Proceedings of the Fourth International Conference
  on Operational Research. John Wiley \& Sons, Hoboken, NJ, pp. 55--71.

\bibitem[{Howick et~al.(2017)Howick, Ackermann, Walls, Quigley, and
  Houghton}]{Howick2017-of_MYLW}
Howick, S., Ackermann, F., Walls, L., Quigley, J., Houghton, T., 2017. Learning
  from mixed {OR} method practice: The {NINES} case study. Omega 69, 70--81.

\bibitem[{Hu et~al.(2018)Hu, Toriello, and Dessouky}]{Hu2018_JLYHK}
Hu, W., Toriello, A., Dessouky, M., 2018. Integrated inventory routing and
  freight consolidation for perishable goods. European Journal of Operational
  Research 271~(2), 548--560.

\bibitem[{Huang et~al.(2011)Huang, Keisler, and Linkov}]{Huang2011-hy_JL}
Huang, I.~B., Keisler, J., Linkov, I., 2011. Multi-criteria decision analysis
  in environmental sciences: Ten years of applications and trends. The Science
  of The Total Environment 409~(19), 3578--3594.

\bibitem[{Huang et~al.(2003)Huang, Caines, and Malham\'{e}}]{Huetal2003_GZ}
Huang, M., Caines, P., Malham\'{e}, R., 2003. Individual and mass behaviour in
  large population stochastic wireless power control problems: {C}entralized
  and {N}ash equilibrium solutions. In: Proceedings of the 42nd {IEEE}
  Conference on Decision and Control, Maui, Hawaii. pp. 98--103.

\bibitem[{Huang et~al.(2007)Huang, Caines, and Malham\'{e}}]{Huetal2007_GZ}
Huang, M., Caines, P., Malham\'{e}, R., 2007. Large-population cost coupled
  {LQG} problems with nonuniform agents: {I}ndividual-mass behavior and
  decentralized epsilon-{N}ash equilibria. {IEEE} Transactions on Automatic
  Control 52, 1560--1571.

\bibitem[{Huang et~al.(2006)Huang, Malham\'{e}, and Caines}]{Huetal2006_GZ}
Huang, M., Malham\'{e}, R., Caines, P., 2006. Large population stochastic
  dynamic games: {C}losed-loop {McKean-Vlasov} systems and the {Nash} certainty
  equivalence principle. Communications in Information \& Systems 6, 221--252.

\bibitem[{Huang et~al.(2012)Huang, Smilowitz, and Balcik}]{Huang2012_JLYHK}
Huang, M., Smilowitz, K., Balcik, B., 2012. Models for relief routing: Equity,
  efficiency and efficacy. Transportation Research Part E: Logistics and
  Transportation Review 48~(1), 2--18.

\bibitem[{Huang et~al.(2019)Huang, Jaimungal, and Nourian}]{Mojtaba19_MB}
Huang, X., Jaimungal, S., Nourian, M., 2019. Mean-field game strategies for
  optimal execution. Applied Mathematical Finance 26~(2), 153--185.

\bibitem[{Huangfu and Hall(2018)}]{Huangfu2018_CTCGE}
Huangfu, Q., Hall, J. A.~J., 2018. Parallelizing the dual revised simplex
  method. Mathematical Programming Computation 10, 119--142.

\bibitem[{Hubicka et~al.(2019)Hubicka, Marcjasz, and
  Weron}]{hub:mar:wer:19_DSRW}
Hubicka, K., Marcjasz, G., Weron, R., 2019. A note on averaging day-ahead
  electricity price forecasts across calibration windows. IEEE Transactions on
  Sustainable Energy 10~(1), 321--323.

\bibitem[{Hughes et~al.(2019)Hughes, Dwivedi, Misra, Rana, Raghavan, and
  Akella}]{Hughes2019-iv_KESXZ}
Hughes, L., Dwivedi, Y.~K., Misra, S.~K., Rana, N.~P., Raghavan, V., Akella,
  V., 2019. Blockchain research, practice and policy: Applications, benefits,
  limitations, emerging research themes and research agenda. International
  Journal of Information Management 49, 114--129.

\bibitem[{Hughes et~al.(2022)Hughes, Dwivedi, Rana, Williams, and
  Raghavan}]{Hughes2022-jg_KESXZ}
Hughes, L., Dwivedi, Y.~K., Rana, N.~P., Williams, M.~D., Raghavan, V., 2022.
  Perspectives on the future of manufacturing within the industry 4.0 era.
  Production Planning \& Control 33~(2-3), 138--158.

\bibitem[{Hughes(1995)}]{Hughes1995-za_KVRH}
Hughes, W.~P., 1995. A salvo model of warships in missile combat used to
  evaluate their staying power. Naval Research Logistics 42~(2), 267--289.

\bibitem[{Huisman et~al.(2005)Huisman, Kroon, Lentink, and
  Vromans}]{Huisman2005_DC}
Huisman, D., Kroon, L.~G., Lentink, R.~M., Vromans, M. J. C.~M., 2005.
  Operations research in passenger railway transportation. Statistica
  Neerlandica 59~(4), 467--497.

\bibitem[{Hunt et~al.(2016)Hunt, Jacobi, Gannon, Zorn, Moore, and
  Lonsdorf}]{Hunt2016_EA}
Hunt, V.~M., Jacobi, S.~K., Gannon, J.~J., Zorn, J.~E., Moore, C.~T., Lonsdorf,
  E.~V., 2016. A decision support tool for adaptive management of native
  prairie ecosystems. Interfaces 46~(4), 334--344.

\bibitem[{Hunter(2007)}]{Hunter:2007_LCAL}
Hunter, J.~D., 2007. Matplotlib: A {2D} graphics environment. Computing in
  Science \& Engineering 9~(3), 90--95.

\bibitem[{Hunter et~al.(2019)Hunter, Applegate, Arora, Chong, Cooper,
  Rinc\'{o}n-Guevara, and Vivas-Valencia}]{Hunter2019_CC}
Hunter, S.~R., Applegate, E.~A., Arora, V., Chong, B., Cooper, K.,
  Rinc\'{o}n-Guevara, O., Vivas-Valencia, C., 2019. An introduction to
  multiobjective simulation optimization. ACM Transactions on Modeling and
  Computer Simulation 29~(1).

\bibitem[{Hussain et~al.(2019)Hussain, Mohd~Salleh, Cheng, and
  Shi}]{hussain2019metaheuristic_MH}
Hussain, K., Mohd~Salleh, M.~N., Cheng, S., Shi, Y., 2019. Metaheuristic
  research: a comprehensive survey. Artificial Intelligence Review 52~(4),
  2191--2233.

\bibitem[{Hussein et~al.(2017)Hussein, Gaber, Elyan, and
  Jayne}]{Hussein17_LCAL}
Hussein, A., Gaber, M.~M., Elyan, E., Jayne, C., 2017. Imitation learning: A
  survey of learning methods. ACM Computing Surveys 50~(2), 1--35.

\bibitem[{Hvattum and Arntzen(2010)}]{Hvattum2010-nn_IM}
Hvattum, L.~M., Arntzen, H., 2010. Using {ELO} ratings for match result
  prediction in association football. International Journal of Forecasting
  26~(3), 460--470.

\bibitem[{Hwang and Yoon(1981)}]{Hwang1981multiple_MESG}
Hwang, C., Yoon, K., 1981. Multiple Attribute Decision Making: Methods and
  Applications. Springer-Verlag, Berlin.

\bibitem[{Hwang et~al.(1992)Hwang, Richards, and Winter}]{HRW92_SMPT}
Hwang, F., Richards, D., Winter, P., 1992. The Steiner Tree Problem. No.~53 in
  Annals of Discrete Mathematics. North-Holland, Amsterdam.

\bibitem[{Hyndman et~al.(2022)Hyndman, Athanasopoulos, Bergmeir, Caceres,
  Chhay, Kuroptevand, Petropoulos, Razbash, Wang, Yasmeen, Reid, Shaub, Carza,
  Team, Ihaka, Wang, Tang, and Zhou}]{Hyndman2022Rforecast_FP}
Hyndman, R., Athanasopoulos, G., Bergmeir, C., Caceres, G., Chhay, L.,
  Kuroptevand, K. O'Hara-Wild, M., Petropoulos, F., Razbash, S., Wang, E.,
  Yasmeen, F., Reid, D., Shaub, D., Carza, F., Team, R.~C., Ihaka, R., Wang,
  X., Tang, Y., Zhou, Z., 2022. forecast: Forecasting functions for time series
  and linear models. R package version 8.17.

\bibitem[{Hyndman and Athanasopoulos(2021)}]{HyndmanAthanasopoulos2021_FP}
Hyndman, R.~J., Athanasopoulos, G., 2021. Forecasting: principles and practice,
  3rd Edition. OTexts, Melbourne, Australia.

\bibitem[{Hyndman and Khandakar(2008)}]{Hyndman2008_FP}
Hyndman, R.~J., Khandakar, Y., 2008. Automatic time series forecasting: The
  forecast package for \textsc{R}. Journal of Statistical Software 27~(3),
  1--22.

\bibitem[{Hyndman and Koehler(2006)}]{HYNDMAN2006_FP}
Hyndman, R.~J., Koehler, A.~B., 2006. Another look at measures of forecast
  accuracy. International Journal of Forecasting 22~(4), 679--688.

\bibitem[{Hyndman et~al.(2008)Hyndman, Koehler, Ord, and
  Snyder}]{Hyndman2008-iu_FP}
Hyndman, R.~J., Koehler, A.~B., Ord, J.~K., Snyder, R.~D., 2008. Forecasting
  with Exponential Smoothing: The State Space Approach. Springer Verlag,
  Berlin.

\bibitem[{Hyndman et~al.(2002)Hyndman, Koehler, Snyder, and
  Grose}]{Hyndman2002-jp_FP}
Hyndman, R.~J., Koehler, A.~B., Snyder, R., Grose, S., 2002. A state space
  framework for automatic forecasting using exponential smoothing methods.
  International Journal of Forecasting 18~(3), 439--454.

\bibitem[{Ibarra-Rojas et~al.(2015)Ibarra-Rojas, Delgado, Giesen, and
  Mu\~noz}]{IbarraRojas2015_DC}
Ibarra-Rojas, O.~J., Delgado, F., Giesen, R., Mu\~noz, J.~C., 2015. Planning,
  operation, and control of bus transport systems{: A} literature review.
  Transportation Research Part B: Methodological 77, 38--75.

\bibitem[{{IBM}(2022)}]{Cplex_CTCGE}
{IBM}, 2022. {ILOG CPLEX Optimization Studio}.
\newline\urlprefix\url{https://www.ibm.com/uk-en/products/ilog-cplex-optimization-studio}

\bibitem[{Iglehart and Karlin(1960)}]{iglehart1960optimal_JSS}
Iglehart, D., Karlin, S., 1960. Optimal policy for dynamic inventory process
  with non-stationary stochastic demands. Tech. rep., Stanford University,
  Applied Mathematics and Statistics Labs.

\bibitem[{Iglehart(1965)}]{iglehart1965limiting_HATI}
Iglehart, D.~L., 1965. Limiting diffusion approximations for the many server
  queue and the repairman problem. Journal of Applied Probability 2~(2),
  429--441.

\bibitem[{{IMO}(2018)}]{Imo2018-zj_HP}
{IMO}, 2018. {Resolution {MEPC.304(72)}, Initial {IMO} Strategy on reduction of
  {GHG} emissions from ships}. Tech. Rep. MEPC 72/17/Add.1, Annex 11,
  International Maritime Organization.

\bibitem[{Inderfurth et~al.(2001)Inderfurth, de~Kok, and
  Flapper}]{Inderfurth2001-ye_AM}
Inderfurth, K., de~Kok, A.~G., Flapper, S. D.~P., 2001. Product recovery in
  stochastic remanufacturing systems with multiple reuse options. European
  Journal of Operational Research 133~(1), 130--152.

\bibitem[{Inderfurth et~al.(2004)Inderfurth, Flapper, Lambert, Pappis, and
  Voutsinas}]{Inderfurth2004-qp_AM}
Inderfurth, K., Flapper, S. D.~P., Lambert, A. J.~D., Pappis, C.~P., Voutsinas,
  T.~G., 2004. Production planning for product recovery management. In: Dekker,
  R., Fleischmann, M., Inderfurth, K., Van~Wassenhove, L.~N. (Eds.), Reverse
  Logistics: Quantitative Models for {Closed-Loop} Supply Chains. Springer,
  Berlin, pp. 249--274.

\bibitem[{Inderfurth and Jensen(1999)}]{Inderfurth1999-kv_AM}
Inderfurth, K., Jensen, T., 1999. Analysis of {MRP} policies with recovery
  options. In: Leopold-Wildburger, U., Feichtinger, G., Kistner, K.-P. (Eds.),
  Modelling and Decisions in Economics: Essays in Honor of Franz Ferschl.
  Physica-Verlag HD, Heidelberg, pp. 189--228.

\bibitem[{{INFORMS}(2022)}]{INFORMS2022}
{INFORMS}, 2022. Certified analytics professional ({CAP}).
  \url{www.certifiedanalytics.org}, accessed on 2022-11-01.

\bibitem[{{InsightMaker}(2016)}]{Insightmaker_CTCGE}
{InsightMaker}, 2016. {InsightMaker}.
\newline\urlprefix\url{https://github.com/scottfr/simulation}

\bibitem[{{Institute for Apprenticeships \& Technical
  Education}(2021)}]{Institute_for_Apprenticeships_Technical_Education_undated-av_MYLW}
{Institute for Apprenticeships \& Technical Education}, 2021. Operational
  research specialist.
  \url{https://www.instituteforapprenticeships.org/apprenticeship-standards/operational-research-specialist-v1-0},
  accessed on 2021-12-9.

\bibitem[{Iori et~al.(2021)Iori, de~Lima, Martello, Miyazawa, and
  Monaci}]{Iori2021-hs_JB}
Iori, M., de~Lima, V.~L., Martello, S., Miyazawa, F.~K., Monaci, M., 2021.
  Exact solution techniques for two-dimensional cutting and packing. European
  Journal of Operational Research 289~(2), 399--415.

\bibitem[{Irnich and Desaulniers(2005)}]{Irnich2005_CA_MB}
Irnich, S., Desaulniers, G., 2005. Shortest path problems with resource
  constraints. In: Desaulniers, G., Desrosiers, J., Solomon, M. (Eds.), Column
  Generation. Springer US, Boston, MA, pp. 33--65.

\bibitem[{Irohara et~al.(2013)Irohara, Kuo, and Leung}]{Irohara2013_JLYHK}
Irohara, T., Kuo, Y.-H., Leung, J. M.~Y., 2013. From preparedness to recovery:
  A tri-level programming model for disaster relief planning. In: Pacino, D.,
  Vo{\ss}, S., Jensen, R.~M. (Eds.), Computational Logistics. Springer, Berlin,
  Heidelberg, pp. 213--228.

\bibitem[{Isaacs(1975)}]{Is1975_GZ}
Isaacs, R., 1975. Differential games, 2nd Edition. Kruger, New York.

\bibitem[{Ismail and Pham(2019)}]{ismail19_MB}
Ismail, A., Pham, H., 2019. Robust {M}arkowitz mean-variance portfolio
  selection under ambiguous covariance matrix. Mathematical Finance 29~(1),
  174--207.

\bibitem[{Israr et~al.(2022)Israr, Ali, Alkhammash, and
  Jussila}]{Israr2022-woKAI}
Israr, A., Ali, Z.~A., Alkhammash, E.~H., Jussila, J.~J., 2022. Optimization
  methods applied to motion planning of unmanned aerial vehicles: A review.
  Drones 6~(5), 126.

\bibitem[{Ivan and Yin(2017)}]{Ivan2017-ux_KESXZ}
Ivan, S., Yin, Y., 2017. The benefits of {3D} printing for supply chains in the
  context of industry 4.0: Cases from automotive industry. Working Paper,
  Kyoto: Doshisha University.

\bibitem[{Ivanov et~al.(2017)Ivanov, Dolgui, and Sokolov}]{Ivanov2017-rv_KESXZ}
Ivanov, D., Dolgui, A., Sokolov, B., 2017. A dynamic approach to multi-stage
  job shop scheduling in an industry 4.0-based flexible assembly system. In:
  Lödding, H., Riedel, R., Thoben, K.-D., von Cieminski, G., Kiritsis, D.
  (Eds.), Advances in Production Management Systems. The Path to Intelligent,
  Collaborative and Sustainable Manufacturing. Springer, pp. 475--482.

\bibitem[{Ivanov et~al.(2016)Ivanov, Dolgui, Sokolov, Werner, and
  Ivanova}]{Ivanov2016-qo_KESXZ}
Ivanov, D., Dolgui, A., Sokolov, B., Werner, F., Ivanova, M., 2016. A dynamic
  model and an algorithm for short-term supply chain scheduling in the smart
  factory industry 4.0. International Journal of Production Research 54~(2),
  386--402.

\bibitem[{{JaamSim Development Team}(2016)}]{Jaamsim_CTCGE}
{JaamSim Development Team}, 2016. {JaamSim: Discrete-Event Simulation
  Software}.
\newline\urlprefix\url{http://jaamsim.com}

\bibitem[{Jackson(1957)}]{jackson1_HATI}
Jackson, J.~R., 1957. Networks of waiting lines. Operations Research 5~(4),
  518--521.

\bibitem[{Jackson(1963)}]{jackson2_HATI}
Jackson, J.~R., 1963. Jobshop-like queueing systems. Management Science 10~(1),
  131--142.

\bibitem[{Jackson(1982)}]{Jackson1982-rr_GM}
Jackson, M.~C., 1982. {The nature of soft systems thinking: The work of
  Churchman, Ackoff and Checkland}. Journal of Applied Systems Analysis 9,
  17--29.

\bibitem[{Jackson(1985)}]{Jackson1985-ws_GM}
Jackson, M.~C., 1985. Social systems theory and practice: The nees for a
  critical approach. International Journal of General Systems 10~(2-3),
  135--151.

\bibitem[{Jackson(1991)}]{Jackson1991-pi_GM}
Jackson, M.~C., 1991. Systems Methodology for the Management Sciences. Plenum,
  New York.

\bibitem[{Jackson(2000)}]{Jackson2000-ei_AJG}
Jackson, M.~C., 2000. Systems Approaches to Management. Kluwer/Plenum, London.

\bibitem[{Jackson(2003)}]{Jackson2003-ln_MYLW}
Jackson, M.~C., 2003. Systems Thinking: Creative Holism for Managers. Wiley,
  Chichester.

\bibitem[{Jackson(2004)}]{Jackson2004-qf_AJG}
Jackson, M.~C., 2004. Community operational research: purposes, theory and
  practice. In: Midgley, G., Ochoa-Arias, A.~E. (Eds.), Community Operational
  Research: {OR} and Systems Thinking for Community Development. Springer, New
  York, NY, pp. 57--74.

\bibitem[{Jackson(2006)}]{Jackson2006-cb_AJG}
Jackson, M.~C., 2006. Beyond problem structuring methods: reinventing the
  future of {OR/MS}. Journal of the Operational Research Society 57~(7),
  868--878.

\bibitem[{Jackson(2019)}]{Jackson2019-eb_GM}
Jackson, M.~C., 2019. Critical Systems Thinking and the Management of
  Complexity. John Wiley \& Sons.

\bibitem[{Jackson and Keys(1984)}]{Jackson1984-zt_GM}
Jackson, M.~C., Keys, P., 1984. Towards a system of systems methodologies.
  Journal of the Operational Research Society 35~(6), 473--486.

\bibitem[{Jackson et~al.(1989)Jackson, Keys, and Cropper}]{Jackson1989-fk_AFRH}
Jackson, M.~C., Keys, P., Cropper, S.~A. (Eds.), 1989. {OR} and the social
  sciences. Plenum Press, New York.

\bibitem[{Jacquet-Lagreze and Siskos(1982)}]{jacquet1982assessing_MESG}
Jacquet-Lagreze, E., Siskos, J., 1982. Assessing a set of additive utility
  functions for multicriteria decision-making, the {UTA} method. European
  Journal of Operational Research 10~(2), 151--164.

\bibitem[{Jacquillat and Odoni(2015)}]{jacquillat2015integrated_VLVV}
Jacquillat, A., Odoni, A.~R., 2015. An integrated scheduling and operations
  approach to airport congestion mitigation. Operations Research 63~(6),
  1390--1410.

\bibitem[{Jadin et~al.(2019)Jadin, Aubry, Schaus, and
  Bonaventure}]{jadin.aubry.ea:19_BF}
Jadin, M., Aubry, F., Schaus, P., Bonaventure, O., 2019. {CG4SR}: Near optimal
  traffic engineering for segment routing with column generation. IEEE INFOCOM
  2019 - IEEE Conference on Computer Communications, 1333--1341.

\bibitem[{Jaehn and Michaelis(2016)}]{Jaehn2016_DC}
Jaehn, F., Michaelis, S., 2016. Shunting of trains in succeeding yards.
  Computers \& Industrial Engineering 102, 1--9.

\bibitem[{Jaehn et~al.(2015)Jaehn, Rieder, and Wiehl}]{Jaehn2015_DC}
Jaehn, F., Rieder, J., Wiehl, A., 2015. Minimizing delays in a shunting yard.
  OR Spectrum 37, 407--429.

\bibitem[{Jain et~al.(1984)Jain, Chiu, and Hawe}]{JaiChiHaw84_JH}
Jain, R., Chiu, D.~M., Hawe, W., 1984. A quantitative measure of fairness and
  discrimination for resource allocation in shared computer systems. Tech. Rep.
  TR--301, Eastern Research Laboratory, DEC, Hudson, MA.

\bibitem[{Jaiswal(2012)}]{Jaiswal2012-md_KVRH}
Jaiswal, N.~K., 2012. Military Operations Research: Quantitative Decision
  Making. Springer, New York, NY.

\bibitem[{Janczura and Wójcik(2022)}]{jan:woj:22_DSRW}
Janczura, J., Wójcik, E., 2022. Dynamic short-term risk management strategies
  for the choice of electricity market based on probabilistic forecasts of
  profit and risk measures. {The} {German} and the {Polish} market case study.
  Energy Economics 110, 106015.

\bibitem[{Jarn\'{\i}k(1930)}]{J30_SMPT}
Jarn\'{\i}k, V., 1930. O jist\'em probl\'emu minim\'aln\'{\i}m (z dopisupanu
  {O.} {B}or{{\r{u}}}vkovi). Pr\'ace Moravsk'e
  P\v{r}\'{\i}rodov\v{e}\-deck\'{e} Spole\v{c}nosti VI~(4), 57--63.

\bibitem[{Ja\'{s}kiewicz and Nowak(2018{\natexlab{a}})}]{Jano2018b_GZ}
Ja\'{s}kiewicz, A., Nowak, A., 2018{\natexlab{a}}. Nonzero-sum stochastic
  games. In: Ba\c{s}ar, T., Zaccour, G. (Eds.), Handbook of Dynamic Game
  Theory. Springer, Cham, pp. 281--344.

\bibitem[{Ja\'{s}kiewicz and Nowak(2018{\natexlab{b}})}]{Jano2018a_GZ}
Ja\'{s}kiewicz, A., Nowak, A., 2018{\natexlab{b}}. Zero-sum stochastic games.
  In: Ba\c{s}ar, T., Zaccour, G. (Eds.), Handbook of Dynamic Game Theory.
  Springer, Cham, pp. 1--65.

\bibitem[{Jauch and Glueck(1975)}]{Jauch1975-hd_JJ}
Jauch, L.~R., Glueck, W.~F., 1975. Evaluation of university professors'
  research performance. Management Science 22~(1), 66--75.

\bibitem[{Jenkins(1969)}]{Jenkins1969-vq_GM}
Jenkins, G., 1969. The systems approach. Journal of Systems Engineering 1,
  3--49.

\bibitem[{Jenkins et~al.(2021)Jenkins, Robbins, and
  Lunday}]{Jenkins2021-nx_KVRH}
Jenkins, P.~R., Robbins, M.~J., Lunday, B.~J., 2021. Approximate dynamic
  programming for military medical evacuation dispatching policies. INFORMS
  Journal on Computing 33~(1), 2--26.

\bibitem[{Jenkins and {Van Kerm}(2011)}]{JenVanKer11_JH}
Jenkins, S.~P., {Van Kerm}, P., 2011. The measurement of economic inequality.
  In: Nolan, B., Salverda, W., Smeeding, T.~M. (Eds.), The Oxford Handbook of
  Economic Inequality. Oxford University Press, pp. 40--68.

\bibitem[{Jepsen et~al.(2008)Jepsen, Petersen, Spoorendonk, and
  Pisinger}]{JepsenPSP2008_CA_MB}
Jepsen, M., Petersen, B., Spoorendonk, S., Pisinger, D., 2008. Subset-row
  inequalities applied to the vehicle-routing problem with time windows.
  Operations Research 56~(2), 497--511.

\bibitem[{Jia et~al.(2016)Jia, Xu, and Guide}]{Jia2016-vy_AM}
Jia, J., Xu, S.~H., Guide, Jr, V. D.~R., 2016. Addressing {Supply--Demand}
  imbalance: Designing efficient remanufacturing strategies. Production and
  Operations Management 25~(11), 1958--1967.

\bibitem[{Jnitova et~al.(2017)Jnitova, Elsawah, and Ryan}]{Jnitova2017-vp_KVRH}
Jnitova, V., Elsawah, S., Ryan, M., 2017. Review of simulation models in
  military workforce planning and management context. The Journal of Defense
  Modeling and Simulation 14~(4), 447--463.

\bibitem[{John et~al.(1994)John, Naim, and Towill}]{John1994_SMD}
John, S., Naim, M.~M., Towill, D.~R., 1994. Dynamic analysis of a {WIP}
  compensated decision support system. International Journal of Manufacturing
  Systems Design 1, 283--297.

\bibitem[{Johnes and Johnes(2016)}]{Johnes2016-di_JJ}
Johnes, G., Johnes, J., 2016. Costs, efficiency, and economies of scale and
  scope in the english higher education sector. Oxford Review of Economic
  Policy 32~(4), 596--614.

\bibitem[{Johnes et~al.(2005)Johnes, Johnes, Thanassoulis, Lenton, and
  Emrouznejad}]{Johnes2005-ue_JJ}
Johnes, G., Johnes, J., Thanassoulis, E., Lenton, P., Emrouznejad, A., 2005. An
  exploratory analysis of the cost structure of higher education in {E}ngland.
  Tech. Rep. 641, Department for Education and Skills, London.

\bibitem[{Johnes and Schwarzenberger(2011)}]{Johnes2011-hd_JJ}
Johnes, G., Schwarzenberger, A., 2011. Differences in cost structure and the
  evaluation of efficiency: the case of {G}erman universities. Education
  Economics 19~(5), 487--499.

\bibitem[{Johnes(1996)}]{Johnes1996-pk_JJ}
Johnes, J., 1996. Performance assessment in higher education in {B}ritain.
  European Journal of Operational Research 89~(1), 18--33.

\bibitem[{Johnes(2006)}]{Johnes2006-id_JJ}
Johnes, J., 2006. Data envelopment analysis and its application to the
  measurement of efficiency in higher education. Economics of Education Review
  25~(3), 273--288.

\bibitem[{Johnes(2014)}]{Johnes2014-jm_JJ}
Johnes, J., 2014. Efficiency and mergers in english higher education 1996/97 to
  2008/9: Parametric and non-parametric estimation of the multi-input
  multi-output distance function. Manchester School 82~(4), 465--487.

\bibitem[{Johnes(2022)}]{Johnes2022-jv_JJ}
Johnes, J., 2022. Applications of production economics in education. In: Ray,
  S.~C., Chambers, R., Kumbhakar, S.~C. (Eds.), Handbook of Production
  Economics: Survey of Applications. Springer, New York, pp. 1193--1239.

\bibitem[{Johnes and Johnes(2013)}]{Johnes2013-eg_JJ}
Johnes, J., Johnes, G., 2013. Efficiency in the higher education sector: A
  technical exploration. Department for Business Innovation and Skills - BIS
  Research Paper 113.

\bibitem[{Johnes and Taylor(1990)}]{Johnes1990-zt_JJ}
Johnes, J., Taylor, J., 1990. Performance Indicators in Higher Education: {UK}
  Universities. Society for Research into Higher Education, Buckingham.

\bibitem[{Johnson et~al.(1978)Johnson, Lenstra, and
  Rinnooy~Kan}]{johnson1978complexity_MH}
Johnson, D.~S., Lenstra, J.~K., Rinnooy~Kan, A. H.~G., 1978. The complexity of
  the network design problem. Networks 8~(4), 279--285.

\bibitem[{Johnson(2012{\natexlab{a}})}]{Johnson2012-jm_AJG}
Johnson, M.~P. (Ed.), 2012{\natexlab{a}}. {Community-Based} Operations
  Research: Decision Modeling for Local Impact and Diverse Populations. Vol.
  167 of International Series in Operations Research \& Management Science.
  Springer, New York, NY.

\bibitem[{Johnson(2012{\natexlab{b}})}]{Johnson2012-iq_AJG}
Johnson, M.~P., 2012{\natexlab{b}}. Community-based operations research:
  introduction, theory, and applications. In: Johnson, M.~P. (Ed.),
  {Community-Based} Operations Research: Decision Modeling for Local Impact and
  Diverse Populations. Vol. 167 of International Series in Operations Research
  \& Management Science. Springer, New York, NY, pp. 3--36.

\bibitem[{Johnson et~al.(2018)Johnson, Midgley, Wright, and
  Chichirau}]{Johnson2018-sj}
Johnson, M.~P., Midgley, G., Wright, J., Chichirau, G., 2018. Community
  operational research: Innovations, internationalization and agenda-setting
  applications. European Journal of Operational Research 268~(3), 761--770.

\bibitem[{Johnson and Smilowitz(2012)}]{Johnson2012-vr_AJG}
Johnson, M.~P., Smilowitz, K., 2012. What is community-based {OR}? In: Johnson,
  M.~P. (Ed.), {Community-Based} Operations Research: Decision Modeling for
  Local Impact and Diverse Populations. Vol. 167 of International Series in
  Operations Research \& Management Science. Springer, New York, NY, pp.
  37--65.

\bibitem[{Jondrow et~al.(1982)Jondrow, Knox~Lovell, Materov, and
  Schmidt}]{Jondrow1982-fq_JJ}
Jondrow, J., Knox~Lovell, C.~A., Materov, I.~S., Schmidt, P., 1982. On the
  estimation of technical inefficiency in the stochastic frontier production
  function model. Journal of Econometrics 19~(2), 233--238.

\bibitem[{Jones and Eden(1981)}]{Jones1981-xt_MYLW}
Jones, S., Eden, C., 1981. {O.R}. in the community. Journal of the Operational
  Research Society 32~(5), 335--345.

\bibitem[{Jotrao and Batta(2021)}]{Jotrao2021_JLYHK}
Jotrao, S., Batta, R., 2021. Time-constrained {UAV} path planning in {3D}
  network for maximim information gain. Military Operations Research 26~(3),
  5--25.

\bibitem[{Jun et~al.(1999)Jun, Jacobson, and Swisher}]{Jun1999-oz_CV}
Jun, J.~B., Jacobson, S.~H., Swisher, J.~R., 1999. Application of
  discrete-event simulation in health care clinics: A survey. Journal of the
  Operational Research Society 50~(2), 109--123.

\bibitem[{Jung et~al.(2019)Jung, Pinedo, Sriskandarajah, and
  Tiwari}]{Jung2019-ls_CV}
Jung, K.~S., Pinedo, M., Sriskandarajah, C., Tiwari, V., 2019. Scheduling
  elective surgeries with emergency patients at shared operating rooms.
  Production and Operations Management 28~(6), 1407--1430.

\bibitem[{J{\"u}nger et~al.(1995)J{\"u}nger, Reinelt, and
  Rinaldi}]{JungerReineltRinaldi:1995_IL}
J{\"u}nger, M., Reinelt, G., Rinaldi, G., 1995. The traveling salesman problem.
  In: Ball, M.~O., Magnanti, T.~L., Monma, C.~L., Nemhauser, G.~L. (Eds.),
  Handbooks in Operations Research and Management Science. Vol.~7. Elsevier,
  pp. 225--330.

\bibitem[{J{\"{u}}nger et~al.(2000)J{\"{u}}nger, Rinaldi, and
  Thienel}]{Junger-et-al:2000_IL}
J{\"{u}}nger, M., Rinaldi, G., Thienel, S., 2000. Practical performance of
  efficient minimum cut algorithms. Algorithmica 26~(1), 172--195.

\bibitem[{Jury and Paynter(1975)}]{jury1975_XW}
Jury, E.~I., Paynter, H., 1975. Inners and stability of dynamic systems. The
  American Society of Mechanical Engineers (ASME).

\bibitem[{Jünger and Thienel(2000)}]{Abacus_CTCGE}
Jünger, M., Thienel, S., 2000. The {ABACUS} system for
  branch-and-cut-and-price algorithms in integer programming and combinatorial
  optimization. Software: Practice and Experience 30~(11), 1325--1352.

\bibitem[{Kaelbling et~al.(1998)Kaelbling, Littman, and
  Cassandra}]{KAELBLING199899_LCAL}
Kaelbling, L.~P., Littman, M.~L., Cassandra, A.~R., 1998. Planning and acting
  in partially observable stochastic domains. Artificial Intelligence 101~(1),
  99--134.

\bibitem[{Kagel(2020)}]{Kagel20_BC}
Kagel, J.~H., 2020. Auctions: A survey of experimental research. In: Kagel,
  J.~H., Roth, A.~E. (Eds.), The Handbook of Experimental Economics. Princeton
  University Press, Princeton, NJ, pp. 501--586.

\bibitem[{Kahneman(2011)}]{Kahneman2011-vb_TC}
Kahneman, D., 2011. Thinking, Fast and Slow. Penguin.

\bibitem[{Kaipia et~al.(2017)Kaipia, Holmström, Småros, and
  Rajala}]{Kaipia2017_SMD}
Kaipia, R., Holmström, J., Småros, J., Rajala, R., 2017. Information sharing
  for sales and operations planning: Contextualized solutions and mechanisms.
  Journal of Operations Management 52, 15--29.

\bibitem[{K{\"a}ki et~al.(2019)K{\"a}ki, Kemppainen, and
  Liesi{\"o}}]{Kaki2019-bo_AFRH}
K{\"a}ki, A., Kemppainen, K., Liesi{\"o}, J., 2019. What to do when
  decision-makers deviate from model recommendations? empirical evidence from
  hydropower industry. European Journal of Operational Research 278~(3),
  869--882.

\bibitem[{Kalai and Smorodinsky(1975)}]{KalSmo75_JH}
Kalai, E., Smorodinsky, M., 1975. Other solutions to {Nash's} bargaining
  problem. Econometrica 43, 513--518.

\bibitem[{Kalman(1960)}]{kalman1960_XW}
Kalman, R.~E., 1960. A new approach to linear filtering and prediction
  problems. Journal of Basic Engineering 82, 35--45.

\bibitem[{Kandakoglu et~al.(2019)Kandakoglu, Frini, and
  Ben~Amor}]{Kandakoglu2019-pm_JL}
Kandakoglu, A., Frini, A., Ben~Amor, S., 2019. Multicriteria decision making
  for sustainable development: A systematic review. Journal of Multi-Criteria
  Decision Analysis 26~(5-6), 202--251.

\bibitem[{Kantorovich(1960)}]{Kantorovich1960-ky_JMB}
Kantorovich, L.~V., 1960. Mathematical methods of organizing and planning
  production. Management Science 6~(4), 366--422.

\bibitem[{Kantorovitch(1958)}]{Kantorovitch1958-xh_JMB}
Kantorovitch, L., 1958. On the translocation of masses. Management Science
  5~(1), 1--4.

\bibitem[{Kao(2014)}]{Kao2014-jt_SL}
Kao, C., 2014. Network data envelopment analysis: A review. European Journal of
  Operational Research 239~(1), 1--16.

\bibitem[{Kao(2016)}]{Kao2016-sl_SL}
Kao, C., 2016. Efficiency decomposition and aggregation in network data
  envelopment analysis. European Journal of Operational Research 255~(3),
  778--786.

\bibitem[{Kao and Liu(2014)}]{Kao2014-dj_SL}
Kao, C., Liu, S.-T., 2014. Multi-period efficiency measurement in data
  envelopment analysis: The case of {T}aiwanese commercial banks. Omega 47,
  90--98.

\bibitem[{Kapelko et~al.(2015)Kapelko, Horta, Camanho, and
  Oude~Lansink}]{Kapelko2015-wf_SL}
Kapelko, M., Horta, I.~M., Camanho, A.~S., Oude~Lansink, A., 2015. Measurement
  of input-specific productivity growth with an application to the construction
  industry in {S}pain and {P}ortugal. International Journal of Production
  Economics 166, 64--71.

\bibitem[{Kaplan(2010)}]{Kaplan2010-dt_KVRH}
Kaplan, E.~H., 2010. Terror queues. Operations Research 58~(4), 773--784.

\bibitem[{Kaplan(1970)}]{kaplan1970dynamic_JSS}
Kaplan, R.~S., 1970. A dynamic inventory model with stochastic lead times.
  Management Science 16~(7), 491--507.

\bibitem[{Kaplan and Garrick(1981)}]{Kaplan1981-vg_TC}
Kaplan, S., Garrick, B.~J., 1981. On the quantitative definition of risk. Risk
  Analysis 1~(1), 11--27.

\bibitem[{Kapu{\'s}ci{\'n}ski and
  Parker(2023)}]{kapuscinski2023capacitated_JSS}
Kapu{\'s}ci{\'n}ski, R., Parker, R.~P., 2023. Capacitated inventory systems.
  In: Song, J.-S. (Ed.), Research Handbook on Inventory Management. Edward
  Elgar Publishing.

\bibitem[{Kara and Sava{\c{s}}er(2017)}]{kara2017humanitarian_BYKOK}
Kara, B.~Y., Sava{\c{s}}er, S., 2017. Humanitarian logistics. In: Batta, R.,
  Peng, J. (Eds.), Leading developments from INFORMS communities. INFORMS
  TutORials in Operations Research. INFORMS, pp. 263--303.

\bibitem[{Karelahti et~al.(2007)Karelahti, Virtanen, and
  Raivio}]{Karelahti2007-jb_KVRH}
Karelahti, J., Virtanen, K., Raivio, T., 2007. {Near-Optimal} missile avoidance
  trajectories via receding horizon control. Journal of Guidance, Control, and
  Dynamics 30~(5), 1287--1298.

\bibitem[{Karger(2000)}]{Karger:2000_IL}
Karger, D.~R., 2000. Minimum cuts in near-linear time. Journal of the {ACM}
  47~(1), 46--76.

\bibitem[{Karger and Stein(1996)}]{Karger-Stein:1996_IL}
Karger, D.~R., Stein, C., 1996. A new approach to the minimum cut problem.
  Journal of the {ACM} 43~(4), 601--640.

\bibitem[{Karmarkar(1984)}]{Karmarkar84_EAY}
Karmarkar, N., 1984. A new polynomial-time algorithm for linear programming.
  In: Proceedings of the Sixteenth Annual ACM Symposium on Theory of Computing.
  STOC '84. Association for Computing Machinery, New York, NY, USA, p.
  302–311.

\bibitem[{Karp(1972{\natexlab{a}})}]{Karp:complexity_UPCT}
Karp, R.~M., 1972{\natexlab{a}}. Reducibility among combinatorial problems. In:
  Miller, R.~E., Thatcher, J.~W., Bohlinger, J.~D. (Eds.), Complexity of
  Computer Computations: Proceedings of a Symposium on the Complexity of
  Computer Computations. The IBM Research Symposia Series. Plenum Press, New
  York, pp. 85--103.

\bibitem[{Karp(1972{\natexlab{b}})}]{Karp:1972_IL}
Karp, R.~M., 1972{\natexlab{b}}. Reducibility among combinatorial problems. In:
  Complexity of Computer Computations. Springer {US}, {Boston, MA}, pp.
  85--103.

\bibitem[{Karras et~al.(2021)Karras, Laine, and Aila}]{Karas19_LCAL}
Karras, T., Laine, S., Aila, T., 2021. A style-based generator architecture for
  generative adversarial networks. {IEEE} Transactions on Pattern Analysis and
  Machine Intelligence 43~(12), 4217--4228.

\bibitem[{Karsu et~al.(2019)Karsu, Kara, and Selvi}]{karsu2019refugee_BYKOK}
Karsu, O., Kara, B.~Y., Selvi, B., 2019. The refugee camp management: a general
  framework and a unifying decision-making model. Journal of Humanitarian
  Logistics and Supply Chain Management 9~(2), 131--150.

\bibitem[{Karsu and Morton(2015)}]{KARSU2015_BYKOK}
Karsu, {\"O}., Morton, A., 2015. Inequity averse optimization in operational
  research. European Journal of Operational Research 245~(2), 343--359.

\bibitem[{Karvetski et~al.(2011)Karvetski, Lambert, Keisler, and
  Linkov}]{Karvetski2011-sk_JL}
Karvetski, C.~W., Lambert, J.~H., Keisler, J.~M., Linkov, I., 2011. Integration
  of decision analysis and scenario planning for coastal engineering and
  climate change. IEEE Transactions on Systems, Man, and Cybernetics - Part A:
  Systems and Humans 41~(1), 63--73.

\bibitem[{Karzanov(1974)}]{Karzanov:1974_IL}
Karzanov, A.~V., 1974. Determining the maximal flow in a network by the method
  of preflows. Soviet Mathematics Doklady 15, 434--437.

\bibitem[{Kasperson et~al.(2022)Kasperson, Webler, Ram, and
  Sutton}]{Kasperson2022-dv_TC}
Kasperson, R.~E., Webler, T., Ram, B., Sutton, J., 2022. The social
  amplification of risk framework: New perspectives. Risk Analysis 42~(7),
  1367--1380.

\bibitem[{Katsaliaki et~al.(2010)Katsaliaki, Mustafee, Dwivedi, Williams, and
  Wilson}]{Katsaliaki2010-ae_CV}
Katsaliaki, K., Mustafee, N., Dwivedi, Y.~K., Williams, T., Wilson, J.~M.,
  2010. A profile of {OR} research and practice published in {Journal of the
  Operational Research Society}. Journal of the Operational Research Society
  61~(1), 82--94.

\bibitem[{Kavadias and Loch(2004)}]{Kavadias2012-mr_WH_ED}
Kavadias, S., Loch, C.~H., 2004. Project Selection Under Uncertainty:
  Dynamically Allocating Resources to Maximize Value. Springer Science \&
  Business Media.

\bibitem[{Kazemi~Matin and Kuosmanen(2009)}]{Kazemi_Matin2009-ha_SL}
Kazemi~Matin, R., Kuosmanen, T., 2009. Theory of integer-valued data
  envelopment analysis under alternative returns to scale axioms. Omega 37~(5),
  988--995.

\bibitem[{Keenan and Jankowski(2019)}]{Keenan2019-mp_JL}
Keenan, P.~B., Jankowski, P., 2019. Spatial decision support systems: Three
  decades on. Decision Support Systems 116, 64--76.

\bibitem[{Keeney(1982)}]{keeney1982decision_MESG}
Keeney, R., 1982. Decision analysis: An overview. Operations Research 30~(5),
  803--838.

\bibitem[{Keeney(1996{\natexlab{a}})}]{keeney1996value_MESG}
Keeney, R., 1996{\natexlab{a}}. Value-focused Thinking: A Path to Creative
  Decisionmaking. Harvard University Press, Cambridge.

\bibitem[{Keeney(1996{\natexlab{b}})}]{Keeney1996-yo}
Keeney, R., 1996{\natexlab{b}}. Value-focused thinking: Identifying decision
  opportunities and creating alternatives. European Journal of Operational
  Research 92~(3), 537--549.

\bibitem[{Keeney and Raiffa(1976)}]{KeeneyRaiffa1976_MESG}
Keeney, R., Raiffa, H., 1976. Decisions with Multiple Objectives: Preferences
  and Value Tradeoffs. John Wiley \& Sons, Hoboken, NJ.

\bibitem[{Kellenbrink and Helber(2015)}]{Kellenbrink2015-av_WH_ED}
Kellenbrink, C., Helber, S., 2015. Scheduling resource-constrained projects
  with a flexible project structure. European Journal of Operational Research
  246~(2), 379--391.

\bibitem[{Keller and Simon(2019)}]{Keller2019-uq_JL}
Keller, L.~R., Simon, J., 2019. Preference functions for spatial risk analysis.
  Risk Analysis 39~(1), 244--256.

\bibitem[{Kellerer et~al.(2004)Kellerer, Pferschy, and Pisinger}]{KPP04_SMPT}
Kellerer, H., Pferschy, U., Pisinger, D., 2004. Knapsack Problems. Springer,
  Berlin.

\bibitem[{Kelley(1961)}]{Kelley1961-ko_WH_ED}
Kelley, J.~E., 1961. {Critical-Path} planning and scheduling: Mathematical
  basis. Operations Research 9~(3), 296--320.

\bibitem[{Kelly(1979)}]{kelly_HATI}
Kelly, F.~P., 1979. Reversibility and Stochastic Networks. Wiley, New York.

\bibitem[{Kelly et~al.(1998)Kelly, Maulloo, and Tan}]{Kelly1998_JH}
Kelly, F.~P., Maulloo, A.~K., Tan, D. K.~H., 1998. Rate control for
  communication networks: {Shadow} prices, proportional fairness and stability.
  Journal of the Operational Research Society 49~(3), 237--252.

\bibitem[{Kemball-Cook and Vaughan(1983)}]{Kemball-Cook1983-yd_CV}
Kemball-Cook, D., Vaughan, J.~P., 1983. Operational research in primary health
  care planning: a theoretical model for estimating the coverage achieved by
  different distributions of staff and facilities. Bulletin of the World Health
  Organization 61~(2), 361--369.

\bibitem[{Kemmer et~al.(2012)Kemmer, Strauss, and
  Winter}]{kemmerDynamicSimultaneousFare2012_AKSJF}
Kemmer, P., Strauss, A.~K., Winter, T., 2012. Dynamic simultaneous fare
  proration for large-scale network revenue management. Journal of the
  Operational Research Society 63~(10), 1336--1350.

\bibitem[{Ker{\"a}nen(2018)}]{Keranen2018-ok_KVRH}
Ker{\"a}nen, J.-P., 2018. Capability assessment finds a multirole fighter
  suitable for {F}inland's defence. Tech. rep., Finnish Air Force.

\bibitem[{Kerivin and Mahjoub(2005)}]{kerivin.mahjoub:05_BF}
Kerivin, H., Mahjoub, A.~R., 2005. Design of survivable networks: A survey.
  Networks 46~(1), 1--21.

\bibitem[{Kerkhove et~al.(2017)Kerkhove, Vanhoucke, and
  Maenhout}]{Kerkhove2017-et_WH_ED}
Kerkhove, L.-P., Vanhoucke, M., Maenhout, B., 2017. On the resource renting
  problem with overtime. Computers \& Industrial Engineering 111, 303--319.

\bibitem[{Kersten et~al.(2017)Kersten, Schr{\"o}der, and
  Indorf}]{Kersten2017-zq_KESXZ}
Kersten, W., Schr{\"o}der, M., Indorf, M., 2017. Potentials of digitalization
  for supply chain risk management: An empirical analysis.
  Betriebswirtschaftliche Aspekte Von Industrie 4.0 71~(17), 47--74.

\bibitem[{Ketzenberg et~al.(2003)Ketzenberg, Souza, and
  Guide}]{Ketzenberg2003-ek_AM}
Ketzenberg, M.~E., Souza, G.~C., Guide, Jr, V. D.~R., 2003. Mixed assembly and
  disassembly operations for remanufacturing. Production and Operations
  Management 12~(3), 320--335.

\bibitem[{Keys(1998)}]{Keys1998-bo_MYLW}
Keys, P., 1998. {OR} as technology revisited. Journal of the Operational
  Research Society 49~(2), 99--108.

\bibitem[{Kharrat et~al.(2020)Kharrat, McHale, and
  Pe{\~n}a}]{Kharrat2020-gg_IM}
Kharrat, T., McHale, I.~G., Pe{\~n}a, J.~L., 2020. Plus--minus player ratings
  for soccer. European Journal of Operational Research 283~(2), 726--736.

\bibitem[{Kiesel and Kusterman(2016)}]{kie:kus:16_DSRW}
Kiesel, R., Kusterman, M., 2016. Structural models for coupled electricity
  markets. Journal of Commodity Markets 3~(1), 16--38.

\bibitem[{K{\i}nay et~al.(2021)K{\i}nay, Gzara, and Alumur}]{kinay2021full_SAA}
K{\i}nay, {\"O}.~B., Gzara, F., Alumur, S.~A., 2021. Full cover charging
  station location problem with routing. Transportation Research Part B:
  Methodological 144, 1--22.

\bibitem[{King(2022)}]{Smi_CTCGE}
King, A., 2022. {SMI}.
\newline\urlprefix\url{https://projects.coin-or.org/Smi}

\bibitem[{King et~al.(1994)King, Rao, and Tarjan}]{King-et-al:1994_IL}
King, V., Rao, S., Tarjan, R., 1994. A faster deterministic maximum flow
  algorithm. Journal of Algorithms 17~(3), 447--474.

\bibitem[{Kingman(1961)}]{kingman1961single_HATI}
Kingman, J. F.~C., 1961. The single server queue in heavy traffic. In: Green,
  B.~J. (Ed.), Mathematical Proceedings of the Cambridge Philosophical Society.
  Vol.~57. Cambridge University Press, pp. 902--904.

\bibitem[{Kingman(1962)}]{kingman1962queues_HATI}
Kingman, J. F.~C., 1962. On queues in heavy traffic. Journal of the Royal
  Statistical Society: Series B (Methodological) 24~(2), 383--392.

\bibitem[{Kingman(1965)}]{kingman1965heavy_HATI}
Kingman, J. F.~C., 1965. The heavy traffic approximation in the theory of
  queues. In: Smith, W.~L., Wilkinson, W.~E. (Eds.), Proceedings of the
  Symposium on Congestion Theory. Vol.~2. University of North Carolina Press,
  Chapel Hill, NC.

\bibitem[{Kingston et~al.(2018{\natexlab{a}})Kingston, Comas-Herrera, and
  Jagger}]{Kingston2018-fi_CV}
Kingston, A., Comas-Herrera, A., Jagger, C., 2018{\natexlab{a}}. Forecasting
  the care needs of the older population in {E}ngland over the next 20 years:
  estimates from the population ageing and care simulation ({PACSim}) modelling
  study. The Lancet Public Health 3~(9), e447--e455.

\bibitem[{Kingston et~al.(2018{\natexlab{b}})Kingston, Post, and {Vanden
  Berghe}}]{kingston2018unifed_GVBSP}
Kingston, J.~H., Post, G., {Vanden Berghe}, G., 2018{\natexlab{b}}. A unifed
  nurse rostering model based on {XHSTT}. In: Proceedings of the 12th
  International Conference of the Practice and Theory of Automated Timetabling.
  Vol.~12. pp. 81--96.

\bibitem[{Kirby(2003)}]{Kirby2003_JLYHK}
Kirby, M.~W., 2003. Operational Research In War And Peace: The British
  Experience from 1930s to 1970. Imperial College Press and The Operational
  Research Society.

\bibitem[{Kirjavainen(2012)}]{Kirjavainen2012-cr_JJ}
Kirjavainen, T., 2012. Efficiency of {F}innish general upper secondary schools:
  an application of stochastic frontier analysis with panel data. Education
  Economics 20~(4), 343--364.

\bibitem[{Kirkpatrick et~al.(1983)Kirkpatrick, Gelatt~Jr, and
  Vecchi}]{kirkpatrick1983optimization_COIT}
Kirkpatrick, S., Gelatt~Jr, C.~D., Vecchi, M.~P., 1983. Optimization by
  simulated annealing. Science 220~(4598), 671--680.

\bibitem[{Klapp et~al.(2020)Klapp, Erera, and
  Toriello}]{klapp2020request_CKTVW}
Klapp, M.~A., Erera, A.~L., Toriello, A., 2020. Request acceptance in same-day
  delivery. Transportation Research Part E: Logistics and Transportation Review
  143, 102083.

\bibitem[{Klee and Minty(1972)}]{Klee1972-oq_JMB}
Klee, V., Minty, G.~J., 1972. How good is the simplex algorithm? In: Shisha, O.
  (Ed.), Inequalities {III}. Academic Press, New York-London, pp. 159--175.

\bibitem[{Klein et~al.(2007)Klein, Connell, and Meyer}]{Klein2007-dj_AJG}
Klein, J.~H., Connell, N. A.~D., Meyer, E., 2007. Operational research practice
  as storytelling. Journal of the Operational Research Society 58~(12),
  1535--1542.

\bibitem[{Klein et~al.(2020)Klein, Koch, Steinhardt, and
  Strauss}]{kleinReviewRevenueManagement2020_AKSJF}
Klein, R., Koch, S., Steinhardt, C., Strauss, A.~K., 2020. A review of revenue
  management: Recent generalizations and advances in industry applications.
  European Journal of Operational Research 284~(2), 397--412.

\bibitem[{Klein et~al.(2018)Klein, Mackert, Neugebauer, and
  Steinhardt}]{klein2018model_CKTVW}
Klein, R., Mackert, J., Neugebauer, M., Steinhardt, C., 2018. A model-based
  approximation of opportunity cost for dynamic pricing in attended home
  delivery. OR Spectrum 40~(4), 969--996.

\bibitem[{Klein et~al.(2019)Klein, Neugebauer, Ratkovitch, and
  Steinhardt}]{klein2019differentiated_CKTVW}
Klein, R., Neugebauer, M., Ratkovitch, D., Steinhardt, C., 2019. Differentiated
  time slot pricing under routing considerations in attended home delivery.
  Transportation Science 53~(1), 236--255.

\bibitem[{Kleinberg and Tardos(2006)}]{Kleinberg-Tardos:2006_IL}
Kleinberg, J., Tardos, E., 2006. Algorithm Design. Addison Wesley.

\bibitem[{Kleinrock(1975{\natexlab{a}})}]{kleinrock_HATI}
Kleinrock, L., 1975{\natexlab{a}}. Queueing Systems, Volume 1: Theory. Wiley
  Intersciences, New York.

\bibitem[{Kleinrock(1975{\natexlab{b}})}]{kleinrock_HATI2}
Kleinrock, L., 1975{\natexlab{b}}. Queueing Systems, Volume 2: Computer
  Applications. Wiley Intersciences, New York.

\bibitem[{Kleywegt et~al.(2002)Kleywegt, Shapiro, and Homem-de
  Mello}]{Kleywegt2002-mi_HL}
Kleywegt, A.~J., Shapiro, A., Homem-de Mello, T., 2002. The sample average
  approximation method for stochastic discrete optimization. SIAM Journal on
  Optimization 12~(2), 479--502.

\bibitem[{Kline et~al.(2019)Kline, Ahner, and Hill}]{Kline2019-fx_KVRH}
Kline, A., Ahner, D., Hill, R., 2019. The weapon-target assignment problem.
  Computers \& Operations Research 105, 226--236.

\bibitem[{Klinke and Renn(2021)}]{Klinke2021-zn_TC}
Klinke, A., Renn, O., 2021. The coming of age of risk governance. Risk Analysis
  41~(3), 544--557.

\bibitem[{Ko et~al.(2018)Ko, Lee, and Ryu}]{Ko2018-pi_KESXZ}
Ko, T., Lee, J., Ryu, D., 2018. Blockchain technology and manufacturing
  industry: Real-time transparency and cost savings. Sustainability 10~(11),
  4274.

\bibitem[{Koch and Klein(2020)}]{Koch2020-bf_HL}
Koch, S., Klein, R., 2020. Route-based approximate dynamic programming for
  dynamic pricing in attended home delivery. European Journal of Operational
  Research 287~(2), 633--652.

\bibitem[{Kohavi et~al.(2002)Kohavi, Rothleder, and
  Simoudis}]{Kohavi2002-rl_JEB}
Kohavi, R., Rothleder, N.~J., Simoudis, E., 2002. Emerging trends in business
  analytics. Communications of the ACM 45~(8), 45--48.

\bibitem[{K{\"o}hler et~al.(2022)K{\"o}hler, Campbell, and
  Ehmke}]{kohler2022data_CKTVW}
K{\"o}hler, C., Campbell, A.~M., Ehmke, J.~F., 2022. Data-driven customer
  acceptance for attended home delivery. Working paper.

\bibitem[{K{\"o}hler et~al.(2020)K{\"o}hler, Ehmke, and
  Campbell}]{kohler2020flexible_CKTVW}
K{\"o}hler, C., Ehmke, J.~F., Campbell, A.~M., 2020. Flexible time window
  management for attended home deliveries. Omega 91, 102023.

\bibitem[{K{\"o}hler and Haferkamp(2019)}]{kohler2019evaluation_CKTVW}
K{\"o}hler, C., Haferkamp, J., 2019. Evaluation of delivery cost approximation
  for attended home deliveries. Transportation Research Procedia 37, 67--74.

\bibitem[{Kohtamäki et~al.(2018)Kohtamäki, Baines, Rabetino, and
  Bigdeli}]{kohtamaki_practices_2018_JHLS}
Kohtamäki, M., Baines, T., Rabetino, R., Bigdeli, A.~Z., 2018. Practices in
  {Servitization}. In: Kohtamäki, M., Baines, T., Rabetino, R., Bigdeli, A.~Z.
  (Eds.), Practices and {Tools} for {Servitization}: {Managing} {Service}
  {Transition}. Springer, Cham, pp. 1--21.

\bibitem[{Kolassa(2011)}]{kolassa2011combining_FP}
Kolassa, S., 2011. Combining exponential smoothing forecasts using {A}kaike
  weights. International Journal of Forecasting 27~(2), 238--251.

\bibitem[{Kolassa(2020)}]{Kolassa2020-jh}
Kolassa, S., 2020. Why the ``best'' point forecast depends on the error or
  accuracy measure. International Journal of Forecasting 36~(1), 208--211.

\bibitem[{K{\"o}llerstr{\"o}m(1974)}]{kollerstrom1974heavy_HATI}
K{\"o}llerstr{\"o}m, J., 1974. Heavy traffic theory for queues with several
  servers. {I}. Journal of Applied Probability 11~(3), 544--552.

\bibitem[{Komarudin et~al.(2020)Komarudin, {De Feyter}, Guerry, and {Vanden
  Berghe}}]{Komarudin_GVBSP}
Komarudin, {De Feyter}, T., Guerry, M.-A., {Vanden Berghe}, G., 2020. The
  extended roster quality staffing problem: addressing roster quality variation
  within a staffing planning period. Journal of Scheduling 23, 253--–264.

\bibitem[{K\"onig(1916)}]{K16_SMPT}
K\"onig, D., 1916. {\"Uber graphen und ihre anwendungen (in {German})}.
  Mathematische Annalen 77, 453--465.

\bibitem[{Konings(2003)}]{konings2003network_MH}
Konings, R., 2003. Network design for intermodal barge transport.
  Transportation Research Record 1820~(1), 17--25.

\bibitem[{Koopmans(1949)}]{Koopmans1949-qx_JMB}
Koopmans, T.~C., 1949. Optimum utilization of the transportation system.
  Econometrica 17, 136--146.

\bibitem[{Kordzadeh and Ghasemaghaei(2022)}]{Kordzadeh2022-yp_JEB}
Kordzadeh, N., Ghasemaghaei, M., 2022. Algorithmic bias: review, synthesis, and
  future research directions. European Journal of Information Systems 31~(3),
  388--409.

\bibitem[{Korte and Vygen(2008)}]{Korte-Vygen:2008_IL}
Korte, B., Vygen, J., 2008. Combinatorial optimization: Theory and algorithms,
  3rd Edition. Springer.

\bibitem[{Koster and Schmidt(2021)}]{Koster2021_MH}
Koster, A. M. C.~A., Schmidt, D.~R., 2021. Robust network design. In: Crainic,
  T.~G., Gendreau, M., Gendron, B. (Eds.), Network Design with Applications to
  Transportation and Logistics. Springer, pp. 317--343.

\bibitem[{Kotary et~al.(2021)Kotary, Fioretto, Hentenryck, and
  Wilder}]{fioretto2021_LCAL}
Kotary, J., Fioretto, F., Hentenryck, P.~V., Wilder, B., 2021. End-to-end
  constrained optimization learning: {A} survey. In: Zhou, Z. (Ed.),
  Proceedings of the Thirtieth International Joint Conference on Artificial
  Intelligence, {IJCAI} 2021, Virtual Event / Montreal, Canada, 19-27 August
  2021. ijcai.org, pp. 4475--4482.

\bibitem[{Kotiadis and Mingers(2006)}]{Kotiadis2006-gv_MYLW}
Kotiadis, K., Mingers, J., 2006. Combining {PSMs} with hard {OR} methods: the
  philosophical and practical challenges. Journal of the Operational Research
  Society 57~(7), 856--867.

\bibitem[{Koutsandreas et~al.(2022)Koutsandreas, Spiliotis, Petropoulos, and
  Assimakopoulos}]{Koutsandreas2022_FP}
Koutsandreas, D., Spiliotis, E., Petropoulos, F., Assimakopoulos, V., 2022. On
  the selection of forecasting accuracy measures. Journal of the Operational
  Research Society 73~(5), 937--954.

\bibitem[{Kov{\'a}cs and Falagara~Sigala(2021)}]{kovacs2021lessons_BYKOK}
Kov{\'a}cs, G., Falagara~Sigala, I., 2021. Lessons learned from humanitarian
  logistics to manage supply chain disruptions. Journal of Supply Chain
  Management 57~(1), 41--49.

\bibitem[{Kovács(2015)}]{Kovacs:2015_IL}
Kovács, P., 2015. Minimum-cost flow algorithms: an experimental evaluation.
  Optimization Methods and Software 30~(1), 94--127.

\bibitem[{Krawczyk and Petkov(2018)}]{Krpe2018_GZ}
Krawczyk, J., Petkov, V., 2018. Multistage games. In: Ba\c{s}ar, T., Zaccour,
  G. (Eds.), Handbook of Dynamic Game Theory. Springer, Cham, pp. 157--213.

\bibitem[{Kress and Van~Leeuwen(2006)}]{Kress2006-ug_MJE}
Kress, G.~R., Van~Leeuwen, T., 2006. Reading Images: The Grammar of Visual
  Design. Routledge.

\bibitem[{Krijkamp et~al.(2018)Krijkamp, Alarid-Escudero, Enns, Jalal, Hunink,
  and Pechlivanoglou}]{Krijkamp2018-iz_CV}
Krijkamp, E.~M., Alarid-Escudero, F., Enns, E.~A., Jalal, H.~J., Hunink, M.
  G.~M., Pechlivanoglou, P., 2018. Microsimulation modeling for health decision
  sciences using {R}: A tutorial. Medical Decision Making 38~(3), 400--422.

\bibitem[{Krishna(2010)}]{Krishna10_BC}
Krishna, V., 2010. Auction Theory. Elsevier, Amsterdam.

\bibitem[{Kristiansen et~al.(2022)Kristiansen, Sandberg, Hansen, Jensen,
  Friederich, and Lazarova-Molnar}]{Kristiansen2022_CTCGE}
Kristiansen, O.~S., Sandberg, U., Hansen, C., Jensen, M.~S., Friederich, J.,
  Lazarova-Molnar, S., 2022. Experimental comparison of open source
  discrete-event simulation frameworks. In: 315–330 (Ed.), Simulation Tools
  and Techniques. SIMUtools 2021. Lecture Notes of the Institute for Computer
  Sciences, Social Informatics and Telecommunications Engineering. Vol. 424.
  Springer, Cham, pp. 315--330.

\bibitem[{Kr{\'o}l and Zdonek(2020)}]{Krol2020-yp_JEB}
Kr{\'o}l, K., Zdonek, D., 2020. Analytics maturity models: An overview.
  Information. An International Interdisciplinary Journal 11~(3), 142.

\bibitem[{Kronqvist et~al.(2019)Kronqvist, Bernal, Lundell, and
  Grossmann}]{Kr19_ALAL}
Kronqvist, J., Bernal, D., Lundell, A., Grossmann, I., 2019. A review and
  comparison of solvers for convex {MINLP}. Optimization and Engineering 20,
  397--455.

\bibitem[{Kroon et~al.(2009)Kroon, Huisman, Abbink, Fioole, Fischetti,
  Mar\'{o}ti, Schrijver, Steenbeek, and Ybema}]{Kroon2009_DC}
Kroon, L., Huisman, D., Abbink, E., Fioole, P.~J., Fischetti, M., Mar\'{o}ti,
  G., Schrijver, A., Steenbeek, A., Ybema, R., 2009. The new {D}utch
  timetable{:} {T}he {OR} revolution. Interfaces 39~(1), 6--17.

\bibitem[{Kroon et~al.(2008)Kroon, Lentink, and Schrijver}]{Kroon2008_DC}
Kroon, L.~G., Lentink, R.~M., Schrijver, A., 2008. Shunting of passenger train
  units{: A}n integrated approach. Transportation Science 42, 436--449.

\bibitem[{Kruskal(1957)}]{K57_SMPT}
Kruskal, J., 1957. On the shortest spanning subtree of a graph and the
  traveling salesman problem. Proceedings of the American Mathematical Society
  7, 48--50.

\bibitem[{K{\"u}{\c{c}}{\"u}kyavuz and Sen(2017)}]{KS17_ALAL}
K{\"u}{\c{c}}{\"u}kyavuz, S., Sen, S., 2017. An introduction to two-stage
  stochastic mixed-integer programming. In: Batta, R., Peng, J. (Eds.), Leading
  Developments from INFORMS Communities. INFORMS TutORials in Operations
  Research. INFORMS, Catonsville, MD, pp. 1--27.

\bibitem[{Kuhn(1955)}]{K55_SMPT}
Kuhn, H., 1955. The {Hungarian} method for the assignment problem. Naval
  Research Logistics Quarterly 2, 83--97.

\bibitem[{Kuhn(1956)}]{K56_SMPT}
Kuhn, H., 1956. Variants of the {Hungarian} method for the assignment problem.
  Naval Research Logistics Quarterly 3, 253--258.

\bibitem[{Kumar and Swaminathan(2003)}]{Kumar2003_XW}
Kumar, S., Swaminathan, J.~M., 2003. Diffusion of innovations under supply
  constraints. Operations Research 51~(6), 866--879.

\bibitem[{Kunc(2017{\natexlab{a}})}]{Kunc2017-bq_MCJM}
Kunc, M., 2017{\natexlab{a}}. System dynamics: A soft and hard approach to
  modelling. In: Chan, W. K.~V., D'Ambrogio, A., Zacharewicz, G., Mustafee, N.,
  Wainer, G., Page, E. (Eds.), Proceedings of the 2017 Winter Simulation
  Conference, Las Vegas, NV. pp. 597--606.

\bibitem[{Kunc(2017{\natexlab{b}})}]{Kunc2017-pw_MCJM}
Kunc, M., 2017{\natexlab{b}}. System Dynamics: Soft and Hard Operational
  Research. Springer, London.

\bibitem[{Kunc et~al.(2020)Kunc, Harper, and Katsikopoulos}]{Kunc2020-fn_MCJM}
Kunc, M., Harper, P., Katsikopoulos, K., 2020. A review of implementation of
  behavioural aspects in the application of {OR} in healthcare. Journal of the
  Operational Research Society 71~(7), 1055--1072.

\bibitem[{Kunc et~al.(2016)Kunc, Malpass, and White}]{Kunc2016-gy_AFRH}
Kunc, M., Malpass, J., White, L., 2016. Behavioral Operational Research:
  Theory, Methodology and Practice. Springer.

\bibitem[{Kunc et~al.(2018)Kunc, Mortenson, and Vidgen}]{Kunc2018-uj_MCJM}
Kunc, M., Mortenson, M.~J., Vidgen, R., 2018. A computational literature review
  of the field of system dynamics from 1974 to 2017. Journal of Simulation
  12~(2), 115--127.

\bibitem[{K{\"u}nnen and
  Strauss(2022)}]{kunnenValueFlexibleFlighttoroute2022_AKSJF}
K{\"u}nnen, J.~R., Strauss, A.~K., 2022. The value of flexible flight-to-route
  assignments in pre-tactical air traffic management. Transportation Research
  Part B: Methodological 160, 76--96.

\bibitem[{Kunnumkal and Topaloglu(2009)}]{kunnumkal2009stochastic_VLVV}
Kunnumkal, S., Topaloglu, H., 2009. A stochastic approximation method for the
  single-leg revenue management problem with discrete demand distributions.
  Mathematical Methods of Operations Research 70~(3), 477--504.

\bibitem[{Kunnumkal and Topaloglu(2010)}]{kunnumkal2010new_DLLD}
Kunnumkal, S., Topaloglu, H., 2010. A new dynamic programming decomposition
  method for the network revenue management problem with customer choice
  behavior. Production and Operations Management 19~(5), 575--590.

\bibitem[{Kunz et~al.(2017)Kunz, Van~Wassenhove, Besiou, Hambye, and
  Kovacs}]{kunz2017relevance_BYKOK}
Kunz, N., Van~Wassenhove, L.~N., Besiou, M., Hambye, C., Kovacs, G., 2017.
  Relevance of humanitarian logistics research: best practices and way forward.
  International Journal of Operations \& Production Management 37~(11),
  1585--1599.

\bibitem[{Kuosmanen(2005)}]{Kuosmanen2005-og_SL}
Kuosmanen, T., 2005. Weak disposability in nonparametric production analysis
  with undesirable outputs. American Journal of Agricultural Economics 87~(4),
  1077--1082.

\bibitem[{Kusner et~al.(2017)Kusner, Loftus, Russell, and
  Silva}]{KusLofRusSil17_JH}
Kusner, M.~J., Loftus, J., Russell, C., Silva, R., 2017. Counterfactual
  fairness. In: Guyon, I., Von~Luxburg, U., Bengio, S., Wallach, H., Fergus,
  R., Vishwanathan, S., Garnett, R. (Eds.), Proceedings of Advances in Neural
  Information Processing Systems 30 (NIPS 2017). Curran Associates, Inc.

\bibitem[{Labb\'e and Violin(2016)}]{Lab16_UPCT}
Labb\'e, M., Violin, A., 2016. Bilevel programming and price setting problems.
  Annals of Operations Research 240, 141--169.

\bibitem[{Laborie and Nuijten(2008)}]{Laborie2010-li_WH_ED}
Laborie, P., Nuijten, W., 2008. Constraint programming formulations and
  propagation algorithms. In: Artigues, C., Demassey, S., N{\'e}ron, E. (Eds.),
  {Resource-Constrained} Project Scheduling. ISTE, London, pp. 63--72.

\bibitem[{Lad`anyi(2004)}]{Bcp_CTCGE}
Lad`anyi, L., 2004. {BCP}.
\newline\urlprefix\url{https://projects.coin-or.org/Bcp}

\bibitem[{Lago et~al.(2018)Lago, {De Ridder}, and {De
  Schutter}}]{lag:rid:sch:18_DSRW}
Lago, J., {De Ridder}, F., {De Schutter}, B., 2018. Forecasting spot
  electricity prices: deep learning approaches and empirical comparison of
  traditional algorithms. Applied Energy 221, 386--405.

\bibitem[{Lago et~al.(2021)Lago, Marcjasz, {De Schutter}, and
  Weron}]{lag:mar:sch:wer:21_DSRW}
Lago, J., Marcjasz, G., {De Schutter}, B., Weron, R., 2021. Forecasting
  day-ahead electricity prices: A review of state-of-the-art algorithms, best
  practices and an open-access benchmark. Applied Energy 293, 116983.

\bibitem[{Laguna and Mart{\'i}(2013)}]{laguna2013heuristics_COIT}
Laguna, M., Mart{\'i}, R., 2013. Heuristics. In: Gass, S.~I., Fu, M.~C. (Eds.),
  Encyclopedia of Operations Research and Management Science. Springer, Boston,
  MA, pp. 695--703.

\bibitem[{Lahtinen et~al.(2017)Lahtinen, H{\"a}m{\"a}l{\"a}inen, and
  Liesi{\"o}}]{Lahtinen2017-sh_JL}
Lahtinen, T.~J., H{\"a}m{\"a}l{\"a}inen, R.~P., Liesi{\"o}, J., 2017. Portfolio
  decision analysis methods in environmental decision making. Environmental
  Modelling \& Software 94, 73--86.

\bibitem[{Lai et~al.(2022)Lai, Wu, Wang, and Wang}]{Lai2022-nz_HP}
Lai, X., Wu, L., Wang, K., Wang, F., 2022. Robust ship fleet deployment with
  shipping revenue management. Transportation Research Part B: Methodological
  161, 169--196.

\bibitem[{Lai et~al.(2015)Lai, Fan, and Huang}]{Lai2015_DC}
Lai, Y.-C., Fan, D.-C., Huang, K.-L., 2015. Optimizing rolling stock assignment
  and maintenance plan for passenger railway operations. Computers \&
  Industrial Engineering 85, 284--295.

\bibitem[{Lamas-Fernandez et~al.(2022)Lamas-Fernandez, Bennell, and
  Martinez-Sykora}]{Lamas-Fernandez2022-fe_JB}
Lamas-Fernandez, C., Bennell, J.~A., Martinez-Sykora, A., 2022. {Voxel-Based}
  solution approaches to the {Three-Dimensional} irregular packing problem.
  Operations Research, DOI: 10.1287/opre.2022.2260.

\bibitem[{Lambert(2003)}]{Lambert2003-au_AM}
Lambert, A. J.~D., 2003. Disassembly sequencing: A survey. International
  Journal of Production Research 41~(16), 3721--3759.

\bibitem[{Lambrechts et~al.(2008)Lambrechts, Demeulemeester, and
  Herroelen}]{Lambrechts2008-fp_WH_ED}
Lambrechts, O., Demeulemeester, E., Herroelen, W., 2008. A tabu search
  procedure for developing robust predictive project schedules. International
  Journal of Production Economics 111~(2), 493--508.

\bibitem[{Lan et~al.(2006)Lan, Clarke, and Barnhart}]{lan2006planning_VLVV}
Lan, S., Clarke, J.-P., Barnhart, C., 2006. Planning for robust airline
  operations: Optimizing aircraft routings and flight departure times to
  minimize passenger disruptions. Transportation Science 40~(1), 15--28.

\bibitem[{Lan and Chiang(2011)}]{LanChi11_JH}
Lan, T., Chiang, M., 2011. An axiomatic theory of fairness in resource
  allocation. Tech. rep., Princeton University.

\bibitem[{{Lan} et~al.(2010){Lan}, {Kao}, and {Chiang}}]{Lan2010_JH}
{Lan}, T., {Kao}, D., {Chiang}, M., 2010. An axiomatic theory of fairness in
  network resource allocation. In: 2010 Proceedings of the 30th IEEE
  International Conference on Computer Communications (INFOCOM). pp. 1--9.

\bibitem[{Land and Doig(1960)}]{LD60_ALAL}
Land, A., Doig, A., 1960. An automatic method of solving discrete programming
  problems. Econometrica 28, 497--520.

\bibitem[{Lane et~al.(2000)Lane, Monefeldt, and Rosenhead}]{Lane2000-qz_CV}
Lane, D.~C., Monefeldt, C., Rosenhead, J.~V., 2000. Looking in the wrong place
  for healthcare improvements: A system dynamics study of an accident and
  emergency department. Journal of the Operational Research Society 51~(5),
  518--531.

\bibitem[{Lane and Oliva(1998)}]{Lane1998-ff_MYLW}
Lane, D.~C., Oliva, R., 1998. The greater whole: Towards a synthesis of system
  dynamics and soft systems methodology. European Journal of Operational
  Research 107~(1), 214--235.

\bibitem[{Lane and Rouwette(2023)}]{Lane2022-md_MCJM}
Lane, D.~C., Rouwette, E. A. J.~A., 2023. Towards a behavioural system
  dynamics: Exploring its scope and delineating its promise. European Journal
  of Operational Research 306~(2), 777--794.

\bibitem[{Lang et~al.(2021{\natexlab{a}})Lang, Cleophas, and
  Ehmke}]{lang2021anticipative_CKTVW}
Lang, M.~A., Cleophas, C., Ehmke, J.~F., 2021{\natexlab{a}}. Anticipative
  dynamic slotting for attended home deliveries. In: Operations Research Forum.
  Vol.~2. Springer, pp. 1--39.

\bibitem[{Lang et~al.(2021{\natexlab{b}})Lang, Reggelin, Müller, and
  Nahhas}]{Lang2021_CTCGE}
Lang, S., Reggelin, T., Müller, M., Nahhas, A., 2021{\natexlab{b}}.
  Open-source discrete-event simulation software for applications in production
  and logistics: {A}n alternative to commercial tools? Procedia Computer
  Science 180, 978--987.

\bibitem[{Langville and Meyer(2013)}]{Langville2013-cf_IM}
Langville, A.~N., Meyer, C.~D., 2013. Who's \#1?: The Science of Rating and
  Ranking. Princeton University Press.

\bibitem[{Laporte(1992)}]{laporte1992vehicle_CA_MB}
Laporte, G., 1992. The vehicle routing problem: An overview of exact and
  approximate algorithms. European Journal of Operational Research 59~(3),
  345--358.

\bibitem[{Laporte(2009)}]{laporte2009fifty_CA_MB}
Laporte, G., 2009. Fifty years of vehicle routing. Transportation Science
  43~(4), 408--416.

\bibitem[{Laporte et~al.(2007)Laporte, Mar\'{i}n, Mesa, and
  Ortega}]{Laporte2007_DC}
Laporte, G., Mar\'{i}n, A., Mesa, J.~A., Ortega, F.~A., 2007. An integrated
  methodology for the rapid transit network design problem. In: Geraets, F.,
  Kroon, L., Sch\"{o}ebel, A., Wagner, D., Zaroliagis, C. (Eds.), Algorithmic
  Methods for Railway Optimization. Vol. 4359 of Lecture Notes in Computer
  Science. Springer, Berlin, Heidelberg, pp. 187--199.

\bibitem[{Laporte et~al.(2015)Laporte, Nickel, and {Saldanha da
  Gama}}]{LNS15_SMPT}
Laporte, G., Nickel, S., {Saldanha da Gama}, F. (Eds.), 2015. Location Science.
  Springer, Berlin.

\bibitem[{Laporte and Pascoal(2015)}]{Laporte2015_DC}
Laporte, G., Pascoal, M. M.~B., 2015. Path based algorithms for metro network
  design. Computers \& Operations Research 62, 78--94.

\bibitem[{Laporte et~al.(2022)Laporte, Rancourt, Rodr{\'\i}guez-Pereira, and
  Silvestri}]{LAPORTE2022_BYKOK}
Laporte, G., Rancourt, M.-{\`E}., Rodr{\'\i}guez-Pereira, J., Silvestri, S.,
  2022. Optimizing access to drinking water in remote areas. {A}pplication to
  {N}epal. Computers \& Operations Research 140, 105669.

\bibitem[{Lariviere(2016)}]{Lariviere2016_SMD}
Lariviere, M.~A., 2016. {O}{M} {F}orum—{S}upply chain contracting:
  {D}oughnuts to bubbles. Manufacturing \& Service Operations Management
  18~(3), 309--313.

\bibitem[{Larkin and Simon(1987)}]{Larkin1987-se_MJE}
Larkin, J.~H., Simon, H.~A., 1987. Why a diagram is (sometimes) worth ten
  thousand words. Cognitive Science 11~(1), 65--100.

\bibitem[{Lasry and Lions(2006{\natexlab{a}})}]{Lali2006a_GZ}
Lasry, J., Lions, P., 2006{\natexlab{a}}. Jeux à champ moyen. {I}--le cas
  stationnaire. Comptes Rendus Math\'{e}matique 343, 619--625.

\bibitem[{Lasry and Lions(2006{\natexlab{b}})}]{Lali2006b_GZ}
Lasry, J., Lions, P., 2006{\natexlab{b}}. Jeux à champ moyen. {II}--horizon
  fini et contr\^{o}le optimal. Comptes Rendus Math\'{e}matique 343, 679--684.

\bibitem[{Lasry and Lions(2007)}]{Lali2007_GZ}
Lasry, J., Lions, P., 2007. Mean field games. Japanese Journal of Mathematics
  2, 229--260.

\bibitem[{Latouche and Ramaswami(1999)}]{latouche_HATI}
Latouche, G., Ramaswami, V., 1999. Introduction to Matrix Analytic Methods in
  Stochastic Modeling. Society for Industrial and Applied Mathematics,
  Philadelphia.

\bibitem[{Lauwens(2021)}]{SimJulia_CTCGE}
Lauwens, B., 2021. {SimJulia}.
\newline\urlprefix\url{https://github.com/BenLauwens/SimJulia.jl}

\bibitem[{Lavoie et~al.(1988)Lavoie, Minoux, and Odier}]{lavoie1988new_VLVV}
Lavoie, S., Minoux, M., Odier, E., 1988. A new approach for crew pairing
  problems by column generation with an application to air transportation.
  European Journal of Operational Research 35~(1), 45--58.

\bibitem[{{Law}(2003)}]{law_how_2003_JHLS}
{Law}, 2003. How to conduct a successful simulation study. In: Proceedings of
  the 2003 Winter Simulation Conference. Vol.~1. pp. 66--70.

\bibitem[{Law(2015)}]{Law2015_CC}
Law, A.~M., 2015. {Simulation Modeling and Analysis}, 5th Edition. McGraw-Hill
  Education, New York.

\bibitem[{Lawler(1976)}]{L76_SMPT}
Lawler, E., 1976. Combinatorial optimization: networks and matroids. Holt,
  Rinehart and Winston, New York.

\bibitem[{Lawler et~al.(1985)Lawler, Lenstra, Rinnooy~Kan, and
  Shmoys}]{LLRS85_SMPT}
Lawler, E.~L., Lenstra, J.~K., Rinnooy~Kan, A. H.~G., Shmoys, D.~B. (Eds.),
  1985. The Traveling Salesman Problem (A guided tour of combinatorial
  optimization). Wiley, Chichester.

\bibitem[{Lawrence et~al.(2006)Lawrence, Goodwin, O'Connor, and
  {\"O}nkal}]{Lawrence2006-ey_FP}
Lawrence, M., Goodwin, P., O'Connor, M., {\"O}nkal, D., 2006. Judgmental
  forecasting: A review of progress over the last 25 years. International
  Journal of Forecasting 22~(3), 493--518.

\bibitem[{Le~Menestrel and Van~Wassenhove(2004)}]{Le_Menestrel2004-wy}
Le~Menestrel, M., Van~Wassenhove, L.~N., 2004. Ethics outside, within, or
  beyond {OR} models? European Journal of Operational Research 153~(2),
  477--484.

\bibitem[{Lebedev et~al.(2021)Lebedev, Goulart, and
  Margellos}]{lebedev2021dynamic_CKTVW}
Lebedev, D., Goulart, P., Margellos, K., 2021. A dynamic programming framework
  for optimal delivery time slot pricing. European Journal of Operational
  Research 292~(2), 456--468.

\bibitem[{Lecun(1989)}]{LeCun89_LCAL}
Lecun, Y., 1989. Generalization and network design strategies. In: Pfeifer, R.,
  Schreter, Z., Fogelman, F., Steels, L. (Eds.), Connectionism in perspective.
  Elsevier, Amsterdam.

\bibitem[{Lee and Johnes(2022)}]{Lee2022-an_JJ}
Lee, B.~L., Johnes, J., 2022. Using network {DEA} to inform policy: The case of
  the teaching quality of higher education in {E}ngland. Higher Education
  Quarterly 76~(2), 399--421.

\bibitem[{Lee(1973)}]{Lee1973-us_GM}
Lee, D.~B., 1973. Requiem for {Large-Scale} models. Journal of the American
  Institute of Planners 39~(3), 163--178.

\bibitem[{Lee et~al.(1997)Lee, Padmanabhan, and Whang}]{Lee1997_XW}
Lee, H.~L., Padmanabhan, V., Whang, S., 1997. Information distortion in a
  supply chain: The bullwhip effect. Management Science 43~(4), 546--558.

\bibitem[{Lee et~al.(2000)Lee, So, and Tang}]{Lee2000_SMD}
Lee, H.~L., So, K.~C., Tang, C.~S., 2000. The value of information sharing in a
  two-level supply chain. Management Science 46~(5), 626--643.

\bibitem[{Lee and Leyffer(2012)}]{LL12_ALAL}
Lee, J., Leyffer, S. (Eds.), 2012. Mixed Integer Nonlinear Programming.
  Springer, New York.

\bibitem[{Lee et~al.(2011)Lee, Mogi, Li, Hui, Lee, Hui, Park, Ha, and
  Kim}]{Lee_2011_DSRW}
Lee, S.~K., Mogi, G., Li, Z., Hui, K.~S., Lee, S.~K., Hui, K.~N., Park, S.~Y.,
  Ha, Y.~J., Kim, J.~W., 2011. Measuring the relative efficiency of hydrogen
  energy technologies for implementing the hydrogen economy: An integrated
  fuzzy {AHP/DEA} approach. International Journal of Hydrogen Energy 36~(20),
  12655--12663.

\bibitem[{Lehtonen and Virtanen(2022)}]{Lehtonen2022-yc_KVRH}
Lehtonen, J.-M., Virtanen, K., 2022. Choosing the most economically
  advantageous tender using a multi-criteria decision analysis approach.
  Journal of Public Procurement 22~(2), 164--179.

\bibitem[{Lempert and Turner(2021)}]{Lempert2021-os_TC}
Lempert, R.~J., Turner, S., 2021. Engaging multiple worldviews with
  quantitative decision support: A robust decision-making demonstration using
  the lake model. Risk Analysis 41~(6), 845--865.

\bibitem[{Lenstra(1983)}]{Le83_ALAL}
Lenstra, H., 1983. Integer programming with a fixed number of variables.
  Mathematics of Operations Research 8, 538--548.

\bibitem[{Leontief(1936)}]{Leontief1936-gg_JMB}
Leontief, W.~W., 1936. Quantitative input and output relations in the economic
  systems of the {U}nited {S}tates. Review of Economics and Statistics 18~(3),
  105--125.

\bibitem[{Lepenioti et~al.(2020)Lepenioti, Bousdekis, Apostolou, and
  Mentzas}]{Lepenioti2020-mq_JEB}
Lepenioti, K., Bousdekis, A., Apostolou, D., Mentzas, G., 2020. Prescriptive
  analytics: Literature review and research challenges. International Journal
  of Information Management 50, 57--70.

\bibitem[{Lettovsk{\`y} et~al.(2000)Lettovsk{\`y}, Johnson, and
  Nemhauser}]{lettovsky2000airline_VLVV}
Lettovsk{\`y}, L., Johnson, E.~L., Nemhauser, G.~L., 2000. Airline crew
  recovery. Transportation Science 34~(4), 337--348.

\bibitem[{Leus and Herroelen(2004)}]{Leus2004-kp_WH_ED}
Leus, R., Herroelen, W., 2004. Stability and resource allocation in project
  planning. IIE Transactions 36~(7), 667--682.

\bibitem[{Lev et~al.(1987)Lev, Bloom, Gleit, Murphy, and
  Shoemaker}]{Lev1987_EA}
Lev, B., Bloom, J., Gleit, A., Murphy, F., Shoemaker, C., 1987. Strategic
  planning in energy and natural resources. North-Holland, New York, NY.

\bibitem[{Levin(1994)}]{Levin1994-yg_AJG}
Levin, M., 1994. Action research and critical systems thinking: Two icons
  carved out of the same log? Systems Practice 7~(1), 25--41.

\bibitem[{Levine(2005)}]{Levine2007-bi_WH_ED}
Levine, H.~A., 2005. Project Portfolio Management: A Practical Guide to
  Selecting Projects, Managing Portfolios, and Maximizing Benefits. Wiley India
  Pvt. Limited.

\bibitem[{Levitt(1972)}]{levitt_production-line_1972_JHLS}
Levitt, T., 1972. Production-line approach to service. Harvard Business Review
  50~(5), 41--52.

\bibitem[{Lewis(2003)}]{Lewis2003-hl_IM}
Lewis, M., 2003. Moneyball: The Art Of Winning An Unfair Game. WW Norton.

\bibitem[{Li and Grossmann(2021)}]{Li2021-qy_HL}
Li, C., Grossmann, I.~E., 2021. A review of stochastic programming methods for
  optimization of process systems under uncertainty. Frontiers of Chemical
  Engineering in China 2.

\bibitem[{Li et~al.(2020)Li, Ding, and Connor}]{li2020switch_DLLD}
Li, D., Ding, L., Connor, S., 2020. When to switch? {I}ndex policies for
  resource scheduling in emergency response. Production and Operations
  Management 29~(2), 241--262.

\bibitem[{Li and Ng(2000)}]{li00_MB}
Li, D., Ng, W.-L., 2000. Optimal dynamic portfolio selection: Multiperiod
  mean-variance formulation. Mathematical Finance 10~(3), 387--406.

\bibitem[{Li and Pang(2017)}]{li2017dynamic_DLLD}
Li, D., Pang, Z., 2017. Dynamic booking control for car rental revenue
  management: A decomposition approach. European Journal of Operational
  Research 256~(3), 850--867.

\bibitem[{Li et~al.(2022{\natexlab{a}})Li, Pang, and Qian}]{li2022bid_DLLD}
Li, D., Pang, Z., Qian, L., 2022{\natexlab{a}}. Bid price controls for car
  rental network revenue management. Production and Operations Management.

\bibitem[{Li et~al.(2021)Li, Wang, Emrouznejad, Zhu, and Kou}]{Li2021-ia_SL}
Li, F., Wang, Y., Emrouznejad, A., Zhu, Q., Kou, G., 2021. Allocating a fixed
  cost across decision-making units with undesirable outputs: A bargaining game
  approach. Journal of the Operational Research Society 73~(10), 2309--2325.

\bibitem[{Li and Womer(2015)}]{Li2015-ey_HL}
Li, H., Womer, N.~K., 2015. Solving stochastic resource-constrained project
  scheduling problems by closed-loop approximate dynamic programming. European
  Journal of Operational Research 246~(1), 20--33.

\bibitem[{Li(2021)}]{Li:2021_IL}
Li, J., 2021. Deterministic mincut in almost-linear time. In: Proceedings of
  the 53rd Annual ACM SIGACT Symposium on Theory of Computing. STOC 2021.
  Association for Computing Machinery, New York, NY, USA, pp. 384--395.

\bibitem[{Li et~al.(2022{\natexlab{b}})Li, Huan, and Kunc}]{Li2022_MC}
Li, L., Huan, Y., Kunc, M., 2022{\natexlab{b}}. The impact of forum content on
  data science open innovation performance: a system dynamics-based causal
  machine learning approach. Working paper, Southampton Business School.

\bibitem[{Li and Yu(2023)}]{li2023perishable_JSS}
Li, Q., Yu, P., 2023. Perishable inventory systems. In: Song, J.-S. (Ed.),
  Research Handbook on Inventory Management. Edward Elgar Publishing.

\bibitem[{Li et~al.(2023)Li, Mutha, Ryan, and Sun}]{Li2022-ad_AM}
Li, Z., Mutha, A., Ryan, J., Sun, D., 2023. Mechanism design under asymmetric
  information regarding demand for remanufactured products. Working paper,
  University of Vermont, Burlington, Vermont, USA.

\bibitem[{Liao et~al.(2017)Liao, Deschamps, Loures, and
  Ramos}]{Liao2017-qi_KESXZ}
Liao, Y., Deschamps, F., Loures, E. d. F.~R., Ramos, L. F.~P., 2017. Past,
  present and future of industry 4.0 - a systematic literature review and
  research agenda proposal. International Journal of Production Research
  55~(12), 3609--3629.

\bibitem[{Liebchen and M\"ohring(2002)}]{Liebchen2002_DC}
Liebchen, C., M\"ohring, R., 2002. A case study in periodic timetabling.
  Electronic Notes in Theoretical Computer Science 66~(6), 1--14.

\bibitem[{Liesi{\"o} et~al.(2021)Liesi{\"o}, Salo, Keisler, and
  Morton}]{Liesio2021-zv_JL}
Liesi{\"o}, J., Salo, A., Keisler, J.~M., Morton, A., 2021. Portfolio decision
  analysis: Recent developments and future prospects. European Journal of
  Operational Research 293~(3), 811--825.

\bibitem[{Lilienfeld(1978)}]{Lilienfeld1978-kd_GM}
Lilienfeld, R., 1978. The Rise of Systems Theory: An Ideological Analysis.
  Wiley.

\bibitem[{Lilley et~al.(2022)Lilley, Whitehead, and Midgley}]{Lilley2022-ud_GM}
Lilley, R., Whitehead, M., Midgley, G., 2022. Mindfulness and behavioural
  insights: Reflections on the meditative brain, systems theory and
  organizational change. Journal of Awareness-Based System Change 2~(2),
  29--57.

\bibitem[{Lin and Liu(2011)}]{Lin2011-yp_HP}
Lin, D.-Y., Liu, H.-Y., 2011. Combined ship allocation, routing and freight
  assignment in tramp shipping. Transportation Research Part E: Logistics and
  Transportation Review 47~(4), 414--431.

\bibitem[{Lin and Kernighan(1973)}]{lin1973effective_COIT}
Lin, S., Kernighan, B.~W., 1973. An effective heuristic algorithm for the
  traveling-salesman problem. Operations Research 21~(2), 498--516.

\bibitem[{Lin et~al.(2020{\natexlab{a}})Lin, Leung, Zhang, and
  Gu}]{Lin2020_JLYHK}
Lin, Y., Leung, J. M.~Y., Zhang, L., Gu, J.-W., 2020{\natexlab{a}}. Single-item
  repairable inventory system with stochastic new and warranty demands.
  Transportation Research Part E: Logistics and Transportation Review 142,
  102035.

\bibitem[{Lin et~al.(2020{\natexlab{b}})Lin, Wang, He, and
  Lee}]{lin2020last_CKTVW}
Lin, Y.~H., Wang, Y., He, D., Lee, L.~H., 2020{\natexlab{b}}. Last-mile
  delivery: Optimal locker location under multinomial logit choice model.
  Transportation Research Part E: Logistics and Transportation Review 142,
  102059.

\bibitem[{Linderoth and Ralphs(2005)}]{Linderoth2005_CTCGE}
Linderoth, J., Ralphs, T., 2005. Noncommercial software for mixed-integer
  linear programming. In: Karlof, J. (Ed.), Integer Programming, 1st Edition.
  CRC Press, Boca Raton, Ch.~10, pp. 1--52.

\bibitem[{Linderoth et~al.(2006)Linderoth, Shapiro, and
  Wright}]{Linderoth2006-ff_HL}
Linderoth, J., Shapiro, A., Wright, S., 2006. The empirical behavior of
  sampling methods for stochastic programming. Annals of Operations Research
  142~(1), 215--241.

\bibitem[{Lindner et~al.(2022)Lindner, Reiner, and Keil}]{Lindner2022-lk_MJE}
Lindner, F., Reiner, G., Keil, S., 2022. Review, trends, and opportunities of
  visualizations in manufacturing and production management: a behavioural
  operations perspective. In: {Machuca et al.} (Ed.), Proceedings of the 6th
  World Conference on Production and Operations Management. pp. 546--555.

\bibitem[{Lismont et~al.(2017)Lismont, Vanthienen, Baesens, and
  Lemahieu}]{Lismont2017-qz_JEB}
Lismont, J., Vanthienen, J., Baesens, B., Lemahieu, W., 2017. Defining
  analytics maturity indicators: A survey approach. International Journal of
  Information Management 37~(3), 114--124.

\bibitem[{Littlewood(1972)}]{littlewood1972forecasting_VLVV}
Littlewood, K., 1972. Forecasting and control of passenger bookings. Airline
  Group International Federation of Operational Research Societies Proceedings,
  1972 12, 95--117.

\bibitem[{Littlewood(2005)}]{littlewoodSpecialIssuePapers2005_AKSJF}
Littlewood, K., 2005. Special issue papers: Forecasting and control of
  passenger bookings. Journal of Revenue and Pricing Management 4~(2),
  111--123.

\bibitem[{Liu et~al.(2014)Liu, Stecke, Lian, and Yin}]{Liu2014-mk_KESXZ}
Liu, C., Stecke, K.~E., Lian, J., Yin, Y., 2014. An implementation framework
  for \emph{seru} production. International Transactions in Operational
  Research 21~(1), 1--19.

\bibitem[{Liu et~al.(2020)Liu, Luo, Schulte, and Kharrat}]{Liu2020-ps_IM}
Liu, G., Luo, Y., Schulte, O., Kharrat, T., 2020. Deep soccer analytics:
  learning an action-value function for evaluating soccer players. Data Mining
  and Knowledge Discovery 34~(5), 1531--1559.

\bibitem[{Liu and {van
  Ryzin}(2008)}]{liuChoicebasedLinearProgramming2008_AKSJF}
Liu, Q., {van Ryzin}, G., 2008. On the choice-based linear programming model
  for network revenue management. Manufacturing \& Service Operations
  Management 10~(2), 288--310.

\bibitem[{Ljubi{\'c}(2007)}]{ljubic:07_BF}
Ljubi{\'c}, I., 2007. A hybrid {VNS} for connected facility location. Lecture
  Notes in Computer Science 4771, 157--169.

\bibitem[{Ljubi\'c(2021)}]{Ljubic:2021_IL}
Ljubi\'c, I., 2021. Solving {S}teiner trees: {R}ecent advances, challenges, and
  perspectives. Networks 77~(2), 177--204.

\bibitem[{Lo and McCord(1998)}]{Lo1998-ix_HP}
Lo, H.~K., McCord, M.~R., 1998. Adaptive ship routing through stochastic ocean
  currents: general formulations and empirical results. Transportation Research
  Part A: Policy and Practice 32~(7), 547--561.

\bibitem[{Loch and Wu(2007)}]{Loch2007-ij_AFRH}
Loch, C.~H., Wu, Y., 2007. Behavioral Operations Management. Now Publishers
  Inc.

\bibitem[{Locoro et~al.(2021)Locoro, Fisher, and Mari}]{Locoro2021-dc_MJE}
Locoro, A., Fisher, W.~P., Mari, L., 2021. Visual information literacy:
  Definition, construct modeling and assessment. IEEE Access 9, 71053--71071.

\bibitem[{Lodi et~al.(2014)Lodi, Martello, Monaci, and Vigo}]{Lodi2014-jq_JB}
Lodi, A., Martello, S., Monaci, M., Vigo, D., 2014. Two-dimensional bin packing
  problems. In: Paschos, V.~T. (Ed.), Paradigms of Combinatorial Optimization.
  John Wiley \& Sons, Hoboken, NJ, pp. 107--129.

\bibitem[{Lodi and Zarpellon(2017)}]{lodi2017learning_LCAL}
Lodi, A., Zarpellon, G., 2017. On learning and branching: a survey. Top 25~(2),
  207--236.

\bibitem[{Lohatepanont and Barnhart(2004)}]{lohatepanont2004airline_VLVV}
Lohatepanont, M., Barnhart, C., 2004. Airline schedule planning: Integrated
  models and algorithms for schedule design and fleet assignment.
  Transportation Science 38~(1), 19--32.

\bibitem[{Long and Meadows(2018)}]{Long2018-ji_CV}
Long, K.~M., Meadows, G.~N., 2018. Simulation modelling in mental health: A
  systematic review. Journal of Simulation 12~(1), 76--85.

\bibitem[{Long et~al.(2017)Long, Pan, Zhang, and Hao}]{Long2017-ie_KESXZ}
Long, Y., Pan, J., Zhang, Q., Hao, Y., 2017. {3D} printing technology and its
  impact on chinese manufacturing. International Journal of Production Research
  55~(5), 1488--1497.

\bibitem[{Longstaff and Schwartz(2001)}]{longstaff01_MB}
Longstaff, F.~A., Schwartz, E.~S., 2001. Valuing {A}merican options by
  simulation: A simple least-squares approach. The Review of Financial Studies
  14~(1), 113--147.

\bibitem[{L{\'o}pez-Torres et~al.(2021)L{\'o}pez-Torres, Johnes, Elliott, and
  Polo}]{Lopez-Torres2021-yl_JJ}
L{\'o}pez-Torres, L., Johnes, J., Elliott, C., Polo, C., 2021. The effects of
  competition and collaboration on efficiency in the {UK} independent school
  sector. Economic Modelling 96, 40--53.

\bibitem[{Lov{\'a}sz and Schrijver(1991)}]{LS91_ALAL}
Lov{\'a}sz, L., Schrijver, A., 1991. Cones of matrices and set-functions and
  0--1 optimization. SIAM Journal on Optimization 1, 166--190.

\bibitem[{Lowe et~al.(2020)Lowe, Espinosa, and Yearworth}]{Lowe2020-ai_MYLW}
Lowe, D., Espinosa, A., Yearworth, M., 2020. Constitutive rules for guiding the
  use of the viable system model: Reflections on practice. European Journal of
  Operational Research 287~(3), 1014--1035.

\bibitem[{Lowe and Yearworth(2019)}]{Lowe2019-gm_MYLW}
Lowe, D., Yearworth, M., 2019. Response to viewpoint: Whither problem
  structuring methods ({PSMs)}? Journal of the Operational Research Society
  70~(8), 1393--1395.

\bibitem[{Lozano(2014)}]{Lozano2014-bb_SL}
Lozano, S., 2014. Company-wide production planning using a multiple technology
  {DEA} approach. Journal of the Operational Research Society 65~(5), 723--734.

\bibitem[{Lozano and Villa(2004)}]{Lozano2004-mz_SL}
Lozano, S., Villa, G., 2004. Centralized resource allocation using data
  envelopment analysis. Journal of Productivity Analysis 22~(1), 143--161.

\bibitem[{Lozano and Villa(2005)}]{Lozano2005-hz_SL}
Lozano, S., Villa, G., 2005. Determining a sequence of targets in {DEA}.
  Journal of the Operational Research Society 56~(12), 1439--1447.

\bibitem[{Lu et~al.(2001)Lu, Yu, and Lu}]{Lu2001-ow_AFRH}
Lu, H.-P., Yu, H.-J., Lu, S. S.~K., 2001. The effects of cognitive style and
  model type on {DSS} acceptance: An empirical study. European Journal of
  Operational Research 131~(3), 649--663.

\bibitem[{Lu et~al.(2017)Lu, Chen, Ma, Ke, Li, Zhang, and
  Maciejewski}]{Lu2017-oe_JEB}
Lu, J., Chen, W., Ma, Y., Ke, J., Li, Z., Zhang, F., Maciejewski, R., 2017.
  Recent progress and trends in predictive visual analytics. Frontiers of
  Computer Science 11~(2), 192--207.

\bibitem[{Luck(1984)}]{Luck1984-sj_AJG}
Luck, M., 1984. Working with inner city community organizations. In: Bowen, K.,
  Cook, A., Luck, M. (Eds.), The Writings of Steve Cook. Operational Research
  Society, Birmingham, pp. 77--78.

\bibitem[{Luenberger and Ye(2016)}]{LY2016_EAY}
Luenberger, D.~G., Ye, Y., 2016. Linear and Nonlinear Programming, 4th Edition.
  International Series in Operations Research \& Management Science. Springer,
  Cham.

\bibitem[{Luis et~al.(2012)Luis, Dolinskaya, and
  Smilowitz}]{luis2012disaster_BYKOK}
Luis, E., Dolinskaya, I.~S., Smilowitz, K.~R., 2012. Disaster relief routing:
  Integrating research and practice. Socio-economic Planning Sciences 46~(1),
  88--97.

\bibitem[{Lund and Oksendal(1991)}]{Lund1991_EA}
Lund, D., Oksendal, B., 1991. Stochastic Models and Option Values: Applications
  to Resources, Environment and Investment Problems. Emerald, Bingley.

\bibitem[{Lundell et~al.(2022)Lundell, Kronqvist, and Westerlund}]{Shot_CTCGE}
Lundell, A., Kronqvist, J., Westerlund, T., 2022. The supporting hyperplane
  optimization toolkit for convex {MINLP}. Journal of Global Optimization 84,
  1--41.

\bibitem[{Lustig et~al.(2010)Lustig, Dietrich, Johnson, and
  Dziekan}]{Lustig2010-ty_JEB}
Lustig, I., Dietrich, B., Johnson, C., Dziekan, C., 2010. The analytics
  journey. Analytics Nov/Dec 2010, 11--18.

\bibitem[{Lysgaard et~al.(2004)Lysgaard, Letchford, and
  Eglese}]{LysgaardLE2004_CA_MB}
Lysgaard, J., Letchford, A., Eglese, R., 2004. A new branch-and-cut algorithm
  for the capacitated vehicle routing problem. Mathematical Programming 100,
  423--445.

\bibitem[{L’Ecuyer(2009)}]{lEcuyer09_MB}
L’Ecuyer, P., 2009. Quasi-{Monte Carlo} methods with applications in finance.
  Finance and Stochastics 13~(3), 307--349.

\bibitem[{Ma et~al.(2022)Ma, Zhang, and Branke}]{Ma2022_EA}
Ma, Y., Zhang, W., Branke, J., 2022. Multi-objective optimisation of
  multifaceted maintenance strategies for wind farms. Journal of the
  Operational Research Society, 1--16, DOI: 10.1080/01605682.2022.2085066.

\bibitem[{Maas et~al.(2022)Maas, Lahr, Uyttenboogaart, Buskens, van~der Zee,
  and {CONTRAST investigators}}]{Maas2022-ty_CV}
Maas, W.~J., Lahr, M. M.~H., Uyttenboogaart, M., Buskens, E., van~der Zee,
  D.-J., {CONTRAST investigators}, 2022. Expediting workflow in the acute
  stroke pathway for endovascular thrombectomy in the northern {N}etherlands: a
  simulation model. BMJ Open 12~(4), e056415.

\bibitem[{Macal(2016)}]{Macal2016_CC}
Macal, C.~M., 2016. Everything you need to know about agent-based modelling and
  simulation. Journal of Simulation 10~(2), 144--156.

\bibitem[{Macdonald(2012)}]{Macdonald2012-xc_IM}
Macdonald, B., 2012. Adjusted {Plus-Minus} for {NHL} players using ridge
  regression with goals, shots, fenwick, and corsi. Journal of Quantitative
  Analysis in Sports 8~(3), 1--24.

\bibitem[{Maciejowska et~al.(2021)Maciejowska, Nitka, and
  Weron}]{mac:nit:wer:21_DSRW}
Maciejowska, K., Nitka, W., Weron, T., 2021. Enhancing load, wind and solar
  generation for day-ahead forecasting of electricity prices. Energy Economics
  99, 105273.

\bibitem[{Maciejowska and Nowotarski(2016)}]{mac:now:16_DSRW}
Maciejowska, K., Nowotarski, J., 2016. A hybrid model for {GEFCom2014}
  probabilistic electricity price forecasting. International Journal of
  Forecasting 32~(3), 1051--1056.

\bibitem[{Mackert(2019)}]{mackert2019choice_CKTVW}
Mackert, J., 2019. Choice-based dynamic time slot management in attended home
  delivery. Computers \& Industrial Engineering 129, 333--345.

\bibitem[{MacQueen(1967)}]{MacQueen1967_LCAL}
MacQueen, J.~B., 1967. Some methods for classification and analysis of
  multivariate observations. In: Cam, L. M.~L., Neyman, J. (Eds.), Proc. of the
  fifth Berkeley Symposium on Mathematical Statistics and Probability. Vol.~1.
  University of California Press, pp. 281--297.

\bibitem[{Macrina et~al.(2019)Macrina, Laporte, Guerriero, and
  Di~Puglia~Pugliese}]{macrina2019energy_COIT}
Macrina, G., Laporte, G., Guerriero, F., Di~Puglia~Pugliese, L., 2019. An
  energy-efficient green-vehicle routing problem with mixed vehicle fleet,
  partial battery recharging and time windows. European Journal of Operational
  Research 276~(3), 971--982.

\bibitem[{Maculan(1987)}]{Maculan:1987_IL}
Maculan, N., 1987. The {{S}teiner} problem in graphs. Annals of Discrete
  Mathematics 31, 185--211.

\bibitem[{Magirou et~al.(2015)Magirou, Psaraftis, and
  Bouritas}]{Magirou2015-bq_HP}
Magirou, E.~F., Psaraftis, H.~N., Bouritas, T., 2015. The economic speed of an
  oceangoing vessel in a dynamic setting. Transportation Research Part B:
  Methodological 76, 48--67.

\bibitem[{Magnanti and Wolsey(1995)}]{magnanti.wolsey:95_BF}
Magnanti, T.~L., Wolsey, L.~A., 1995. Optimal trees. Handbooks in Operations
  Research and Management Science 7, 503--615.

\bibitem[{Magnanti and Wong(1981)}]{magnanti1981accelerating_MH}
Magnanti, T.~L., Wong, R.~T., 1981. Accelerating benders decomposition:
  Algorithmic enhancement and model selection criteria. Operations Research
  29~(3), 464--484.

\bibitem[{Magnanti and Wong(1984)}]{magnanti1984network_MH}
Magnanti, T.~L., Wong, R.~T., 1984. Network design and transportation planning:
  Models and algorithms. Transportation Science 18~(1), 1--55.

\bibitem[{Mahar and Wright(2017)}]{mahar2017store_TVWCK}
Mahar, S., Wright, P.~D., 2017. In-store pickup and returns for a dual channel
  retailer. IEEE Transactions on Engineering Management 64~(4), 491--504.

\bibitem[{Maharjan et~al.(2020)Maharjan, Shrestha, Rakhal, Suman, Hulst, and
  Hanaoka}]{Maharjan2020_JLYHK}
Maharjan, R., Shrestha, Y., Rakhal, B., Suman, S., Hulst, J., Hanaoka, S.,
  2020. Mobile logistics hubs prepositioning for emergency preparedness and
  response in {Nepal}. Journal of Humanitarian Logistics and Supply Chain
  Management 10~(4), 555--572.

\bibitem[{Maher(2016)}]{maher2016solving_VLVV}
Maher, S.~J., 2016. Solving the integrated airline recovery problem using
  column-and-row generation. Transportation Science 50~(1), 216--239.

\bibitem[{Mahler et~al.(2022)Mahler, Girard, and
  Kariniotakis}]{mah:gir:kar:22_DSRW}
Mahler, V., Girard, R., Kariniotakis, G., 2022. Data-driven structural modeling
  of electricity price dynamics. Energy Economics 107, 105811.

\bibitem[{Mahmoudi et~al.(2022)Mahmoudi, Shirzad, and
  Verter}]{mahmoudi2022decision_BYKOK}
Mahmoudi, M., Shirzad, K., Verter, V., 2022. Decision support models for
  managing food aid supply chains: A systematic literature review.
  Socio-Economic Planning Sciences, 101255.

\bibitem[{Maister(1976)}]{Maister1976_SMD}
Maister, D.~H., 1976. Centralisation of inventories and the ``square root
  law''. International Journal of Physical Distribution 6~(3), 124--134.

\bibitem[{Makhorin(2020)}]{Glpk_CTCGE}
Makhorin, A., 2020. {GNU Linear Programming Kit (GLPK)}.
\newline\urlprefix\url{https://www.gnu.org/software/glpk}

\bibitem[{Makridakis et~al.(2020{\natexlab{a}})Makridakis, Hyndman, and
  Petropoulos}]{Makridakis2020-eu_FP}
Makridakis, S., Hyndman, R.~J., Petropoulos, F., 2020{\natexlab{a}}.
  Forecasting in social settings: The state of the art. International Journal
  of Forecasting 36~(1), 15--28.

\bibitem[{Makridakis et~al.(2020{\natexlab{b}})Makridakis, Spiliotis, and
  Assimakopoulos}]{Makridakis2019-oy_FP}
Makridakis, S., Spiliotis, E., Assimakopoulos, V., 2020{\natexlab{b}}. The {M4}
  competition: 100,000 time series and 61 forecasting methods. International
  Journal of Forecasting 36~(1), 54--74.

\bibitem[{Malcolm et~al.(1959)Malcolm, Roseboom, Clark, and
  Fazar}]{Malcolm1959-yn_WH_ED}
Malcolm, D.~G., Roseboom, J.~H., Clark, C.~E., Fazar, W., 1959. Application of
  a technique for research and development program evaluation. Operations
  Research 7~(5), 646--669.

\bibitem[{Malczewski and Jankowski(2020)}]{Malczewski2020-gs_JL}
Malczewski, J., Jankowski, P., 2020. Emerging trends and research frontiers in
  spatial multicriteria analysis. International Journal of Geographical
  Information Science 34~(7), 1257--1282.

\bibitem[{Malecki et~al.(2020)Malecki, Keating, and Safdar}]{Malecki2020-cm_TC}
Malecki, K. M.~C., Keating, J.~A., Safdar, N., 2020. Crisis communication and
  public perception of {COVID-19} risk in the era of social media. Clinical
  Infectious Diseases 72~(4), 697--702.

\bibitem[{Malmquist(1953)}]{Malmquist1953-ve_JJ}
Malmquist, S., 1953. Index numbers and indifference surfaces. Trabajos de
  Estadistica 4~(2), 209--242.

\bibitem[{Mandelbaum et~al.(1999)Mandelbaum, Massey, Reiman, and
  Rider}]{mandelbaum1999time_HATI}
Mandelbaum, A., Massey, W.~A., Reiman, M.~I., Rider, B., 1999. Time varying
  multiserver queues with abandonment and retrials. In: Proceedings of the 16th
  International Teletraffic Conference. Vol.~4. pp. 4--7.

\bibitem[{Mandelbaum and Stolyar(2004)}]{mandelbaum2004scheduling_HATI}
Mandelbaum, A., Stolyar, A.~L., 2004. Scheduling flexible servers with convex
  delay costs: Heavy-traffic optimality of the generalized c$\mu$-rule.
  Operations Research 52~(6), 836--855.

\bibitem[{Manerba and Mansini(2014)}]{manerba2014effective_COIT}
Manerba, D., Mansini, R., 2014. An effective matheuristic for the capacitated
  total quantity discount problem. Computers \& Operations Research 41, 1--11.

\bibitem[{Mangasarian(1994)}]{Man94_EAY}
Mangasarian, O.~L., 1994. Nonlinear Programming. SIAM, Philadelphia, PA.

\bibitem[{Maniezzo et~al.(2021)Maniezzo, St{\"u}tzle, and
  Vo{\ss}}]{maniezzo2021matheuristics_MH}
Maniezzo, V., St{\"u}tzle, T., Vo{\ss}, S., 2021. Matheuristics. Springer.

\bibitem[{Mansikka et~al.(2019)Mansikka, Virtanen, and
  Harris}]{Mansikka2019-ol_KVRH}
Mansikka, H., Virtanen, K., Harris, D., 2019. Dissociation between mental
  workload, performance, and task awareness in pilots of high performance
  aircraft. IEEE Transactions on Human-Machine Systems 49~(1), 1--9.

\bibitem[{Mansikka et~al.(2021{\natexlab{a}})Mansikka, Virtanen, Harris, and
  Jalava}]{Mansikka2021-rw_KVRH}
Mansikka, H., Virtanen, K., Harris, D., Jalava, M., 2021{\natexlab{a}}.
  Measurement of team performance in air combat -- have we been
  underperforming? Theoretical Issues in Ergonomics Science 22~(3), 338--359.

\bibitem[{Mansikka et~al.(2021{\natexlab{b}})Mansikka, Virtanen, Harris, and
  Salom{\"a}ki}]{Mansikka2021-is_KVRH}
Mansikka, H., Virtanen, K., Harris, D., Salom{\"a}ki, J., 2021{\natexlab{b}}.
  Live--virtual--constructive simulation for testing and evaluation of air
  combat tactics, techniques, and procedures, {P}art 1: assessment framework.
  The Journal of Defense Modeling and Simulation 18~(4), 285--293.

\bibitem[{Maoz(2013)}]{Maoz2013-zz_JEB}
Maoz, M., 2013. How {IT} should deepen big data analysis to support
  customer-centricity. Tech. Rep. G00248980, Gartner.

\bibitem[{Mar-Molinero and Mingers(2007)}]{Mar-Molinero2007-pi_JJ}
Mar-Molinero, C., Mingers, J., 2007. An evaluation of the limitations of, and
  alternatives to, the {Co-Plot} methodology. Journal of the Operational
  Research Society 58~(7), 874--886.

\bibitem[{Marcjasz et~al.(2022)Marcjasz, Narajewski, Weron, and
  Ziel}]{mar:nar:wer:zie:22_DSRW}
Marcjasz, G., Narajewski, M., Weron, R., Ziel, F., 2022. Distributional neural
  networks for electricity price forecasting. arXiv:2207.02832.

\bibitem[{Mardani et~al.(2017)Mardani, Zavadskas, Streimikiene, Jusoh, and
  Khoshnoudi}]{Mardani_2017_DSRW}
Mardani, A., Zavadskas, E.~K., Streimikiene, D., Jusoh, A., Khoshnoudi, M.,
  2017. A comprehensive review of data envelopment analysis {(DEA)} approach in
  energy efficiency. Renewable and Sustainable Energy Reviews 70~(1),
  1298--1322.

\bibitem[{Mar{\'\i}n and Pelegr{\'\i}n(2019)}]{marin2019p_SAA}
Mar{\'\i}n, A., Pelegr{\'\i}n, M., 2019. p-median problems. In: Laporte, G.,
  Nickel, S., Saldanha~da Gama, F. (Eds.), Location Science. Springer, pp.
  25--50.

\bibitem[{Markowitz(1952)}]{Markowitz52_MB}
Markowitz, H., 1952. Portfolio selection. Journal of Finance 7~(1), 77--91.

\bibitem[{Markowitz and Manne(1957)}]{MM57_ALAL}
Markowitz, H., Manne, A., 1957. On the solution of discrete programming
  problems. Econometrica 25, 84--110.

\bibitem[{Markowitz and Todd(2000)}]{markowitz00_MB}
Markowitz, H.~M., Todd, G.~P., 2000. Mean-variance analysis in portfolio choice
  and capital markets. Wiley, Hoboken, NJ.

\bibitem[{Marla et~al.(2017)Marla, Vaaben, and
  Barnhart}]{marla2017integrated_VLVV}
Marla, L., Vaaben, B., Barnhart, C., 2017. Integrated disruption management and
  flight planning to trade off delays and fuel burn. Transportation Science
  51~(1), 88--111.

\bibitem[{Marler and Arora(2004)}]{marler2004survey}
Marler, R.~T., Arora, J.~S., 2004. Survey of multi-objective optimization
  methods for engineering. Structural and Multidisciplinary Optimization 26,
  369--395.

\bibitem[{Mar{\'o}ti and Kroon(2005)}]{Maroti2005_DC}
Mar{\'o}ti, G., Kroon, L., 2005. Maintenance routing for train units: The
  transition model. Transportation Science 39~(4), 518--525.

\bibitem[{Mar{\'o}ti and Kroon(2007)}]{Maroti2007_DC}
Mar{\'o}ti, G., Kroon, L., 2007. Maintenance routing for train units: The
  interchange model. Computers \& Operations Research 34~(4), 1121--1140.

\bibitem[{Martello(2010)}]{M10_SMPT}
Martello, S., 2010. {J}eno {E}gerv\'ary: from the origins of the {H}ungarian
  algorithm to satellite communication. Central European Journal of Operations
  Research 18, 47--58.

\bibitem[{Martello et~al.(1987)Martello, Laporte, Minoux, and
  Ribeiro}]{MLMR87_SMPT}
Martello, S., Laporte, G., Minoux, M., Ribeiro, C. (Eds.), 1987. Surveys in
  Combinatorial Optimization. Vol.~31 of Annals of Discrete Mathematics.
  North-Holland, Amsterdam.

\bibitem[{Martello et~al.(2003)Martello, Monaci, and Vigo}]{Martello2003-mo_JB}
Martello, S., Monaci, M., Vigo, D., 2003. An exact approach to the
  {Strip-Packing} problem. INFORMS Journal on Computing 15~(3), 310--319.

\bibitem[{Martello and Toth(1990)}]{MT90_SMPT}
Martello, S., Toth, P., 1990. Knapsack Problems: Algorithms and Computer
  Implementations. Wiley, Chichester.

\bibitem[{Mart{\'\i} et~al.(2018)Mart{\'\i}, Pardalos, and
  Resende}]{marte2018handbook_COIT}
Mart{\'\i}, R., Pardalos, P.~M., Resende, M. G.~C., 2018. Handbook of
  Heuristics. Springer.

\bibitem[{Martin-Martinez et~al.(2022)Martin-Martinez, Samsó, Houghton, and
  Solé}]{Pysd_CTCGE}
Martin-Martinez, E., Samsó, R., Houghton, J., Solé, J., 2022. {PySD: S}ystem
  dynamics modeling in python. Journal of Open Source Software 7~(78), 4329.

\bibitem[{Martinovic et~al.(2018)Martinovic, Scheithauer, and Val{\'e}rio~de
  Carvalho}]{Martinovic2018-by_JB}
Martinovic, J., Scheithauer, G., Val{\'e}rio~de Carvalho, J.~M., 2018. A
  comparative study of the arcflow model and the one-cut model for
  one-dimensional cutting stock problems. European Journal of Operational
  Research 266~(2), 458--471.

\bibitem[{Marttunen et~al.(2018)Marttunen, Belton, and
  Lienert}]{Marttunen2018-qs_JL}
Marttunen, M., Belton, V., Lienert, J., 2018. Are objectives hierarchy related
  biases observed in practice? {A} meta-analysis of environmental and energy
  applications of {Multi-Criteria} decision analysis. European Journal of
  Operational Research 265~(1), 178--194.

\bibitem[{Marttunen et~al.(2019)Marttunen, Haag, Belton, Mustajoki, and
  Lienert}]{Marttunen2019-nl_JL}
Marttunen, M., Haag, F., Belton, V., Mustajoki, J., Lienert, J., 2019. Methods
  to inform the development of concise objectives hierarchies in multi-criteria
  decision analysis. European Journal of Operational Research 277~(2),
  604--620.

\bibitem[{Marttunen et~al.(2017)Marttunen, Lienert, and
  Belton}]{Marttunen2017-yk_JL}
Marttunen, M., Lienert, J., Belton, V., 2017. Structuring problems for
  {Multi-Criteria} decision analysis in practice: A literature review of method
  combinations. European Journal of Operational Research 263~(1), 1--17.

\bibitem[{Martzoukos(2009)}]{Martzoukos2009_EA}
Martzoukos, S.~H., 2009. Real {R}\&{D} options and optimal activation of
  two-dimensional random controls. Journal of the Operational Research Society
  60~(6), 843--858.

\bibitem[{Mashlakov et~al.(2021)Mashlakov, Kuronen, Lensu, Kaarna, and
  Honkapuro}]{mas:etal:21_DSRW}
Mashlakov, A., Kuronen, T., Lensu, L., Kaarna, A., Honkapuro, S., 2021.
  Assessing the performance of deep learning models for multivariate
  probabilistic energy forecasting. Applied Energy 285, 116405.

\bibitem[{Masmoudi and Abdelaziz(2018)}]{Masmoudi2018-ya_HL}
Masmoudi, M., Abdelaziz, F.~B., 2018. Portfolio selection problem: a review of
  deterministic and stochastic multiple objective programming models. Annals of
  Operations Research 267~(1), 335--352.

\bibitem[{Masmoudi et~al.(2022)Masmoudi, Mancini, Baldacci, and
  Kuo}]{masmoudi2022vehicle_JLYHK}
Masmoudi, M.~A., Mancini, S., Baldacci, R., Kuo, Y.-H., 2022. Vehicle routing
  problems with drones equipped with multi-package payload compartments.
  Transportation Research Part E: Logistics and Transportation Review 164,
  102757.

\bibitem[{Mason and Mitroff(1981)}]{Mason1981-wa_GM}
Mason, R.~O., Mitroff, I.~I., 1981. Challenging Strategic Planning Assumptions:
  Theory, Cases, and Techniques. Wiley.

\bibitem[{Massey(1981)}]{massey_HATI}
Massey, W.~A., 1981. Non-stationary queues. Ph.D. thesis, Stanford University.

\bibitem[{Massey and Whitt(1998)}]{massey4_HATI}
Massey, W.~A., Whitt, W., 1998. Uniform acceleration expansions for {M}arkov
  chains with time-varying rates. Annals of Applied Probability 8~(4),
  1130--1155.

\bibitem[{Mattila and Virtanen(2014)}]{Mattila2014-mm_KVRH}
Mattila, V., Virtanen, K., 2014. Maintenance scheduling of a fleet of fighter
  aircraft through multi-objective simulation-optimization. Simulation 90~(9),
  1023--1040.

\bibitem[{Mauttone et~al.(2021)Mauttone, Cancela, and
  Urquhart}]{Mauttone2021_MH}
Mauttone, A., Cancela, H., Urquhart, M.~E., 2021. Public transportation. In:
  Crainic, T.~G., Gendreau, M., Gendron, B. (Eds.), Network Design with
  Applications to Transportation and Logistics. Springer, pp. 539--565.

\bibitem[{M\'{a}ximo and Nascimento(2021)}]{MAXIMO20211108_CA_MB}
M\'{a}ximo, V., Nascimento, M., 2021. A hybrid adaptive iterated local search
  with diversification control to the capacitated vehicle routing problem.
  European Journal of Operational Research 294~(3), 1108--1119.

\bibitem[{Mayer and Tr\"uck(2018)}]{may:tru:18_DSRW}
Mayer, K., Tr\"uck, S., 2018. Electricity markets around the world. Journal of
  Commodity Markets 9, 77--100.

\bibitem[{Mazumdar et~al.(1991)Mazumdar, Mason, and
  Douligeris}]{Mazumdar1991_JH}
Mazumdar, R., Mason, L., Douligeris, C., 1991. {Fairness in network optimal
  flow control: Optimality of product forms}. IEEE Transactions on
  Communications 39~(5), 775--782.

\bibitem[{McAfee and McMillan(1987)}]{McAfee87_BC}
McAfee, R.~P., McMillan, J., 1987. Auctions and bidding. Journal of Economic
  Literature 25~(2), 699--738.

\bibitem[{McCloskey(1987)}]{McCloskey1987-gw_GL}
McCloskey, J.~F., 1987. British operational research in {World War} {II}.
  Operations Research 35~(3), 453--470.

\bibitem[{McCollum(2002)}]{ITC-2002_GVBSP}
McCollum, B., 2002. Integrating human abilities and automated systems for
  timetabling: {A} competition using {STARK} and {HuSSH} representations. In:
  Proceedings of the International Conference of the Practice and Theory of
  Automated Timetabling. pp. 265--273.

\bibitem[{McCollum et~al.(2007)McCollum, McMullan, Burke, Parkes, and
  Qu}]{McCollumEtal2010_GVBSP}
McCollum, B., McMullan, P., Burke, E.~K., Parkes, A.~J., Qu, R., 2007. The
  second international timetabling competition: Examination timetabling track.
  Tech. rep., Technical Report QUB/IEEE/Tech/ITC2007/-Exam/v4. 0/17, Queen’s
  University.

\bibitem[{McConnell et~al.(2021)McConnell, Hodgson, Kay, King, Liu, Parlier,
  Thoney-Barletta, and Wilson}]{McConnell2021_JLYHK}
McConnell, B.~M., Hodgson, T.~J., Kay, M.~G., King, R.~E., Liu, Y., Parlier,
  G.~H., Thoney-Barletta, K., Wilson, J.~R., 2021. Assessing uncertainty and
  risk in an expeditionary military logistics network. Journal of Defense
  Modeling and Simulation: Applications, Methodology, Technology 18~(2),
  135--156.

\bibitem[{{McElfresh} and {Dickerson}(2018)}]{McElfresh2018_JH}
{McElfresh}, C., {Dickerson}, J., 2018. Balancing lexicographic fairness and a
  utilitarian objective with application to kidney exchange. In: 32nd AAAI
  Conference on Artificial Intelligence. pp. 1161--1168.

\bibitem[{McEvoy et~al.(2017)McEvoy, Gilbertz, Anderson, Ormerod, and
  Bergmann}]{McEvoy2017-zp_TC}
McEvoy, J., Gilbertz, S.~J., Anderson, M.~B., Ormerod, K.~J., Bergmann, N.~T.,
  2017. Cultural theory of risk as a heuristic for understanding perceptions of
  oil and gas development in {Eastern Montana}, {USA}. The Extractive
  Industries and Society 4~(4), 852--859.

\bibitem[{McHale and Morton(2011)}]{McHale2011-qs_IM}
McHale, I., Morton, A., 2011. A {Bradley-Terry} type model for forecasting
  tennis match results. International Journal of Forecasting 27~(2), 619--630.

\bibitem[{McHale and Holmes(2022)}]{McHale2022-vx_IM}
McHale, I.~G., Holmes, B., 2022. Estimating transfer fees of professional
  footballers using advanced performance metrics and machine learning. European
  Journal of Operational Research 306~(1), 389--399.

\bibitem[{{M}c{K}inney(2010)}]{mckinney-proc-scipy-2010_LCAL}
{M}c{K}inney, W., 2010. {D}ata {S}tructures for {S}tatistical {C}omputing in
  {P}ython. In: {S}t\'efan van~der {W}alt, {J}arrod {M}illman (Eds.),
  {P}roceedings of the 9th {P}ython in {S}cience {C}onference. pp. 56 -- 61.

\bibitem[{Meeusen and van~den Broeck(1977)}]{Meeusen1977-gq_JJ}
Meeusen, W., van~den Broeck, J., 1977. Efficiency estimation from
  {Cobb-Douglas} production functions with composed error. International
  Economic Review 18~(2), 435--444.

\bibitem[{Mehrabi et~al.(2022)Mehrabi, Morstatter, Saxena, Lerman, and
  Galstyan}]{MehMorSaxLerGal22_JH}
Mehrabi, N., Morstatter, F., Saxena, N., Lerman, K., Galstyan, A., 2022. A
  survey on bias and fairness in machine learning. arXiv:1908-09635.

\bibitem[{Melese et~al.(2015)Melese, Richter, and Solomon}]{Melese2015-zf_KVRH}
Melese, F., Richter, A., Solomon, B., 2015. Military {Cost-Benefit} Analysis:
  Theory and Practice. Routledge.

\bibitem[{Melo et~al.(2009)Melo, Nickel, and
  Saldanha-Da-Gama}]{melo2009facility_SAA}
Melo, M.~T., Nickel, S., Saldanha-Da-Gama, F., 2009. Facility location and
  supply chain management--a review. European Journal of Operational Research
  196~(2), 401--412.

\bibitem[{Meng et~al.(2015)Meng, Qi, Zhang, Ang, Chu, and Sim}]{Meng2015-dg_HL}
Meng, F., Qi, J., Zhang, M., Ang, J., Chu, S., Sim, M., 2015. A robust
  optimization model for managing elective admission in a public hospital.
  Operations Research 63~(6), 1452--1467.

\bibitem[{Meng and Wang(2011)}]{Meng2011-of_HP}
Meng, Q., Wang, S., 2011. Optimal operating strategy for a long-haul liner
  service route. European Journal of Operational Research 215~(1), 105--114.

\bibitem[{Menger(1927)}]{M27_SMPT}
Menger, K., 1927. Zur allgemeinen {K}urventheorie. Fundamenta Mathematicae 10,
  96--115.

\bibitem[{Meredith and Mantel(2003)}]{Meredith2011-uu_WH_ED}
Meredith, J.~R., Mantel, Jr, S.~J., 2003. Project Management: A Managerial
  Approach. John Wiley \& Sons.

\bibitem[{Meretoja et~al.(2014)Meretoja, Keshtkaran, Saver, Tatlisumak,
  Parsons, Kaste, Davis, Donnan, and Churilov}]{Meretoja2014-he_CV}
Meretoja, A., Keshtkaran, M., Saver, J.~L., Tatlisumak, T., Parsons, M.~W.,
  Kaste, M., Davis, S.~M., Donnan, G.~A., Churilov, L., 2014. Stroke
  thrombolysis: save a minute, save a day. Stroke 45~(4), 1053--1058.

\bibitem[{Merkle et~al.(2002)Merkle, Middendorf, and
  Schmeck}]{Merkle2002-xi_WH_ED}
Merkle, D., Middendorf, M., Schmeck, H., 2002. Ant colony optimization for
  resource-constrained project scheduling. IEEE Transactions on Evolutionary
  Computation 6~(4), 333--346.

\bibitem[{Mertens and Neyman(1981)}]{Mene1981_GZ}
Mertens, J.-F., Neyman, A., 1981. Stochastic games. International Journal of
  Game Theory 10, 53--56.

\bibitem[{Mertens et~al.(2015)Mertens, Sorin, and Zamir}]{Metal2015_GZ}
Mertens, J.-F., Sorin, S., Zamir, S., 2015. Repeated games. Cambridge
  University Press.

\bibitem[{Messner and Pinson(2019)}]{mes:pin:19_DSRW}
Messner, J., Pinson, P., 2019. Online adaptive lasso estimation in vector
  autoregressive models for high dimensional wind power forecasting.
  International Journal of Forecasting 35~(4), 1485--1498.

\bibitem[{Mete and Zabinsky(2010)}]{Mete2010_JLYHK}
Mete, H.~O., Zabinsky, Z.~B., 2010. Stochastic optimization of medical supply
  location and distribution in disaster management. International Journal of
  Production Economics 126~(1), 76--84.

\bibitem[{Meyer et~al.(2018)Meyer, Jancsary, H{\"o}llerer, and
  Boxenbaum}]{Meyer2018-ui_MJE}
Meyer, R.~E., Jancsary, D., H{\"o}llerer, M.~A., Boxenbaum, E., 2018. The role
  of verbal and visual text in the process of institutionalization. Academy of
  Management Review 43~(3), 392--418.

\bibitem[{Michaud(1989)}]{michaud89_MB}
Michaud, R.~O., 1989. The {M}arkowitz optimization enigma: Is ‘optimized’
  optimal? Financial Analysts Journal 45~(1), 31--42.

\bibitem[{Michna et~al.(2020)Michna, Disney, and Nielsen}]{Michna2020_SMD}
Michna, Z., Disney, S.~M., Nielsen, P., 2020. The impact of stochastic lead
  times on the bullwhip effect under correlated demand and moving average
  forecasts. Omega 93, 102033.

\bibitem[{Midgley(1992)}]{Midgley1992-cr_GM}
Midgley, G., 1992. The sacred and profane in critical systems thinking. Systems
  Practice 5~(1), 5--16.

\bibitem[{Midgley(1994)}]{Midgley1994-cj_GM}
Midgley, G., 1994. Ecology and the poverty of humanism: A critical systems
  perspective. Systems Research 11~(4), 67--76.

\bibitem[{Midgley(2000)}]{Midgley2000-gr_GM}
Midgley, G., 2000. Systemic Intervention: Philosophy, Methodology, and
  Practice. Kluwer/Plenum, New York.

\bibitem[{Midgley et~al.(2013)Midgley, Cavana, Brocklesby, Foote, Wood, and
  Ahuriri-Driscoll}]{Midgley2013-ec_MYLW}
Midgley, G., Cavana, R.~Y., Brocklesby, J., Foote, J.~L., Wood, D. R.~R.,
  Ahuriri-Driscoll, A., 2013. Towards a new framework for evaluating systemic
  problem structuring methods. European Journal of Operational Research
  229~(1), 143--154.

\bibitem[{Midgley et~al.(2018)Midgley, Johnson, and
  Chichirau}]{Midgley2018-jv_AJG}
Midgley, G., Johnson, M.~P., Chichirau, G., 2018. What is {Community
  Operational Research}? European Journal of Operational Research 268~(3),
  771--783.

\bibitem[{Midgley et~al.(1998)Midgley, Munlo, and Brown}]{Midgley1998-rw_GM}
Midgley, G., Munlo, I., Brown, M., 1998. The theory and practice of boundary
  critique: Developing housing services for older people. Journal of the
  Operational Research Society 49~(5), 467--478.

\bibitem[{Midgley and Ochoa-Arias(2004)}]{Midgley2004-xx_AJG}
Midgley, G., Ochoa-Arias, A.~E., 2004. Visions of community for community {OR}.
  In: Midgley, G., Ochoa-Arias, A.~E. (Eds.), Community Operational Research:
  OR and Systems Thinking for Community Development. Kluwer Academic/Plenum
  Publishers, New York, NY, pp. 75--105.

\bibitem[{Midgley and Pinz{\'o}n(2011)}]{Midgley2011-gm_GM}
Midgley, G., Pinz{\'o}n, L.~A., 2011. Boundary critique and its implications
  for conflict prevention. Journal of the Operational Research Society 62~(8),
  1543--1554.

\bibitem[{Midgley and Rajagopalan(2021)}]{Midgley2021-hy_GM}
Midgley, G., Rajagopalan, R., 2021. Critical systems thinking, systemic
  intervention, and beyond. In: Metcalf, G.~S., Kijima, K., Deguchi, H. (Eds.),
  Handbook of Systems Sciences. Springer Singapore, Singapore, pp. 107--157.

\bibitem[{Miettinen(1999)}]{miett99_MESG}
Miettinen, K., 1999. Nonlinear Multiobjective Optimization. Vol.~12 of
  International Series in Operations Research and Management Science. Kluwer
  Academic Publishers, Dordrecht.

\bibitem[{Milgrom(1985)}]{Milgrom85_BC}
Milgrom, P.~R., 1985. The economics of competitive bidding: A selective survey.
  In: Hurwicz, L., Schmeidler, D., Sonnenschein, H. (Eds.), Social Goals and
  Social Organization: Essays in Memory of Elisha Pazner. Cambridge University
  Press, New York, NY, pp. 261--289.

\bibitem[{Milgrom(1987)}]{Milgrom87_BC}
Milgrom, P.~R., 1987. Auction theory. In: Hart, O., Holmstrom, B., Bewley, T.
  (Eds.), Advances in Economic Theory: Fifth World Congress. Econometric
  Society Monographs. Cambridge University Press, Cambridge, pp. 1--32.

\bibitem[{Milgrom(2004)}]{Milgrom04_BC}
Milgrom, P.~R., 2004. Putting Auction Theory to Work. Churchill Lectures in
  Economics. Cambridge University Press, Cambridge.

\bibitem[{Min(2010)}]{Min2010-tg_MCJM}
Min, H., 2010. Artificial intelligence in supply chain management: theory and
  applications. International Journal of Logistics Research and Applications
  13~(1), 13--39.

\bibitem[{Mingers(2000)}]{Mingers2000-oa_MYLW}
Mingers, J., 2000. The contribution of critical realism as an underpinning
  philosophy for {OR/MS} and systems. Journal of the Operational Research
  Society 51~(11), 1256--1270.

\bibitem[{Mingers(2011{\natexlab{a}})}]{Mingers2011-id}
Mingers, J., 2011{\natexlab{a}}. Ethics and {OR}: Operationalising discourse
  ethics. European Journal of Operational Research 210~(1), 114--124.

\bibitem[{Mingers(2011{\natexlab{b}})}]{mingers_soft_2011_JHLS}
Mingers, J., 2011{\natexlab{b}}. Soft {OR} comes of age—but not everywhere!
  Omega 39~(6), 729--741.

\bibitem[{Mingers(2015)}]{mingers_helping_2015_JHLS}
Mingers, J., 2015. Helping business schools engage with real problems: {The}
  contribution of critical realism and systems thinking. European Journal of
  Operational Research 242~(1), 316--331.

\bibitem[{Mingers and Brocklesby(1997)}]{Mingers1997-nd_MYLW}
Mingers, J., Brocklesby, J., 1997. Multimethodology: Towards a framework for
  mixing methodologies. Omega 25~(5), 489--509.

\bibitem[{Mingers and Gill(1997)}]{Mingers1997-vb_GM}
Mingers, J., Gill, A., 1997. Multimethodology: Towards Theory and Practice and
  Mixing and Matching Methodologies. Wiley.

\bibitem[{Mingers and Rosenhead(2004)}]{Mingers2004-qf_MYLW}
Mingers, J., Rosenhead, J., 2004. Problem structuring methods in action.
  European Journal of Operational Research 152~(3), 530--554.

\bibitem[{Mingers and White(2010)}]{Mingers2010-ei_MYLW}
Mingers, J., White, L., 2010. A review of the recent contribution of systems
  thinking to operational research and management science. European Journal of
  Operational Research 207~(3), 1147--1161.

\bibitem[{Mingers(1980)}]{Mingers1980-ya_GM}
Mingers, J.~C., 1980. Towards an appropriate social theory for applied systems
  thinking: Critical theory and soft systems methodology. Journal of Applied
  Systems Analysis 7, 41--50.

\bibitem[{Mingers(1984)}]{Mingers1984-yt_GM}
Mingers, J.~C., 1984. Subjectivism and soft systems methodology -- a critique.
  Journal of Applied Systems Analysis 11, 85--103.

\bibitem[{Mirchandani and Francis(1990)}]{MPL90_SMPT}
Mirchandani, P., Francis, R. (Eds.), 1990. Discrete Location Theory. Wiley,
  Chichester.

\bibitem[{Mirzaei and Seifi(2015)}]{Mirzaei2015_JLYHK}
Mirzaei, S., Seifi, A., 2015. Considering lost sale in inventory routing
  problems for perishable goods. Computers and Industrial Engineering 87,
  213--227.

\bibitem[{Miser and Quade(1985)}]{Miser1985-yx_GM}
Miser, H.~J., Quade, E.~S., 1985. Handbook of Systems Analysis: Overview of
  Uses, Procedures, Applications, and Practice. North-Holland.

\bibitem[{Miser and Quade(1988)}]{Miser1988-ev_GM}
Miser, H.~J., Quade, E.~S. (Eds.), 1988. Handbook of Systems Analysis: Craft
  Issues and Procedural Choices. Wiley, New York.

\bibitem[{Mishra and Parkes(2009)}]{Mishra09_BC}
Mishra, D., Parkes, D.~C., 2009. Multi-item {Vickrey–Dutch} auctions. Games
  and Economic Behavior 66~(1), 326--347.

\bibitem[{Mitchell et~al.(2022)Mitchell, Kean, Mason, O'Sullivan, Phillips, and
  Peschiera}]{Pulp_CTCGE}
Mitchell, S., Kean, A., Mason, A., O'Sullivan, M., Phillips, A., Peschiera, F.,
  2022. {PuLP}.
\newline\urlprefix\url{https://projects.coin-or.org/pulp}

\bibitem[{Mitchell(1997)}]{mitchell1997machine_LCAL}
Mitchell, T.~M., 1997. Machine learning. McGraw-hill, New York NY.

\bibitem[{Mnih et~al.(2015)Mnih, Kavukcuoglu, Silver, Rusu, Veness, Bellemare,
  Graves, Riedmiller, Fidjeland, Ostrovski, Petersen, Beattie, Sadik,
  Antonoglou, King, Kumaran, Wierstra, Legg, and Hassabis}]{Mnih_DRL_2015_LCAL}
Mnih, V., Kavukcuoglu, K., Silver, D., Rusu, A.~A., Veness, J., Bellemare,
  M.~G., Graves, A., Riedmiller, M., Fidjeland, A.~K., Ostrovski, G., Petersen,
  S., Beattie, C., Sadik, A., Antonoglou, I., King, H., Kumaran, D., Wierstra,
  D., Legg, S., Hassabis, D., 2015. Human-level control through deep
  reinforcement learning. Nature 518~(7540), 529--533.

\bibitem[{Mo and Walrand(2000)}]{MoWal00_JH}
Mo, J., Walrand, J., 2000. Fair end-to-end window-based congestion control.
  IEEE/ACM Transactions on Networking 8, 556--567.

\bibitem[{Moccia et~al.(2006)Moccia, Cordeau, Gaudioso, and
  Laporte}]{Moccia2006-cl_HP}
Moccia, L., Cordeau, J.-F., Gaudioso, M., Laporte, G., 2006. A branch-and-cut
  algorithm for the quay crane scheduling problem in a container terminal.
  Naval Research Logistics 53~(1), 45--59.

\bibitem[{Moeller(2010)}]{moeller_characteristics_2010_JHLS}
Moeller, S., 2010. Characteristics of services: a customer integration
  perspective uncovers their value. Journal of Services Marketing 24~(5),
  359--368.

\bibitem[{Moghaddam and DePuy(2011)}]{Moghaddam2011-yh_HL}
Moghaddam, K.~S., DePuy, G.~W., 2011. Farm management optimization using chance
  constrained programming method. Computers and Electronics in Agriculture
  77~(2), 229--237.

\bibitem[{Mohiuddin et~al.(2017)Mohiuddin, Busby, Savovi{\'c}, Richards,
  Northstone, Hollingworth, Donovan, and Vasilakis}]{Mohiuddin2017-zh_CV}
Mohiuddin, S., Busby, J., Savovi{\'c}, J., Richards, A., Northstone, K.,
  Hollingworth, W., Donovan, J.~L., Vasilakis, C., 2017. Patient flow within
  {UK} emergency departments: a systematic review of the use of computer
  simulation modelling methods. BMJ Open 7~(5), e015007.

\bibitem[{Monks et~al.(2019)Monks, Currie, Onggo, Robinson, Kunc, and
  Taylor}]{Monks2019-wt_CV}
Monks, T., Currie, C. S.~M., Onggo, B.~S., Robinson, S., Kunc, M., Taylor, S.
  J.~E., 2019. Strengthening the reporting of empirical simulation studies:
  Introducing the {STRESS} guidelines. Journal of Simulation 13~(1), 55--67.

\bibitem[{Monks et~al.(2015)Monks, Pearson, Pitt, Stein, and
  James}]{Monks2015-yy_CV}
Monks, T., Pearson, M., Pitt, M., Stein, K., James, M.~A., 2015. Evaluating the
  impact of a simulation study in emergency stroke care. Operations Research
  for Health Care 6, 40--49.

\bibitem[{Monks et~al.(2012)Monks, Pitt, Stein, and James}]{Monks2012-ec_CV}
Monks, T., Pitt, M., Stein, K., James, M., 2012. Maximizing the population
  benefit from thrombolysis in acute ischemic stroke. Stroke 43~(10),
  2706--2711.

\bibitem[{Monma and Shallcross(1989)}]{monma.shallcross:89_BF}
Monma, C., Shallcross, D., 1989. Methods for designing communications networks
  with certain two-connected survivability constraints. Operations Research
  37~(4), 531--541.

\bibitem[{Montero-Manso et~al.(2020)Montero-Manso, Athanasopoulos, Hyndman, and
  Talagala}]{Montero-Manso2020_FP}
Montero-Manso, P., Athanasopoulos, G., Hyndman, R.~J., Talagala, T.~S., 2020.
  {FFORMA}: Feature-based forecast model averaging. International Journal of
  Forecasting 36~(1), 86--92.

\bibitem[{Mor and Speranza(2020)}]{mor2022vehicle_CA_MB}
Mor, A., Speranza, M., 2020. Vehicle routing problems over time: a survey. 4OR
  - A Quarterly Journal of Operations Research 18, 129--149.

\bibitem[{Morales et~al.(2014)Morales, Conejo, Madsen, Pinson, and
  Zugno}]{mor:con:mad:pin:zug:14_DSRW}
Morales, J.~M., Conejo, A.~J., Madsen, H., Pinson, P., Zugno, M., 2014.
  Integrating Renewables in Electricity Markets: {O}perational Problems.
  Springer, New York, NY.

\bibitem[{Morecroft(2010)}]{Morecroft2010-ui_MCJM}
Morecroft, J. D.~W., 2010. {Romeo and Juliet in Brazil}: Use of metaphorical
  models for feedback systems thinking. In: Richmond, J., Stuntz, L., Richmond,
  K., {Egner}, {Joanne} (Eds.), Tracing Connections: Voices of Systems
  Thinkers. ISEE Systems, Inc. and The Creative Learning Exchange, pp. 95--119.

\bibitem[{Morecroft(2012)}]{Morecroft2012-zf_MCJM}
Morecroft, J. D.~W., 2012. Metaphorical models for limits to growth and
  industrialization. Systems Research and Behavioral Science 29~(6), 645--666.

\bibitem[{Morecroft(2015)}]{Morecroft2015-vr_MCJM}
Morecroft, J. D.~W., 2015. Strategic Modelling and Business Dynamics: A
  feedback systems approach. John Wiley \& Sons.

\bibitem[{Mortenson et~al.(2015)Mortenson, Doherty, and
  Robinson}]{Mortenson2015-ak_JEB}
Mortenson, M.~J., Doherty, N.~F., Robinson, S., 2015. Operational research from
  {Taylorism to Terabytes: A} research agenda for the analytics age. European
  Journal of Operational Research 241~(3), 583--595.

\bibitem[{Morton and Pentico(1993)}]{morton1993heuristic_COIT}
Morton, T., Pentico, D.~W., 1993. Heuristic Scheduling Systems: With
  Applications to Production Systems and Project Management. Vol.~3.
  Wiley-Interscience, New York.

\bibitem[{Moss et~al.(2022)Moss, Vasilakis, and Wood}]{Moss2022-ao_CV}
Moss, S.~J., Vasilakis, C., Wood, R.~M., 2022. Exploring financially
  sustainable initiatives to address out-of-area placements in psychiatric
  {ICUs}: a computer simulation study. Journal of Mental Health, 1--9.

\bibitem[{Mostajabdaveh et~al.(2019)Mostajabdaveh, Gutjahr, and
  Sibel~Salman}]{Sibel2019inequity_JH}
Mostajabdaveh, M., Gutjahr, W.~J., Sibel~Salman, F., 2019. Inequity-averse
  shelter location for disaster preparedness. IISE Transactions 51~(8),
  809--829.

\bibitem[{Moulin(1988)}]{Mo1988_GZ}
Moulin, H., 1988. Axioms of Cooperative Decision Making, 1st Edition. Cambridge
  University Press.

\bibitem[{Muckstadt and Roundy(1993)}]{muckstadt1993analysis_JSS}
Muckstadt, J.~A., Roundy, R.~O., 1993. Analysis of multistage production
  systems. Handbooks in Operations Research and Management Science 4, 59--131.

\bibitem[{Muehlheusser et~al.(2018)Muehlheusser, Schneemann, Sliwka, and
  Wallmeier}]{Muehlheusser2018-no_IM}
Muehlheusser, G., Schneemann, S., Sliwka, D., Wallmeier, N., 2018. {The
  Contribution of Managers to Organizational Success: Evidence from {G}erman
  Soccer}. Journal of Sports Economics 19~(6), 786--819.

\bibitem[{M\"uller et~al.(2018)M\"uller, Rudov\'a, and
  M\"ullerov\'a}]{muller_GVBSP}
M\"uller, T., Rudov\'a, H., M\"ullerov\'a, Z., 2018. University course
  timetabling and international timetabling competition 2019. In: Proceedings
  of the 12th International Conference of the Practice and Theory of Automated
  Timetabling. pp. 5--31.

\bibitem[{M{\"u}ller-Merbach(1981)}]{muller1981heuristics_COIT}
M{\"u}ller-Merbach, H., 1981. Heuristics and their design: a survey. European
  Journal of Operational Research 8~(1), 1--23.

\bibitem[{Munien and Ezugwu(2021)}]{Munien2021-sh_JB}
Munien, C., Ezugwu, A.~E., 2021. Metaheuristic algorithms for one-dimensional
  bin-packing problems: A survey of recent advances and applications. Journal
  of Intelligent Systems 30~(1), 636--663.

\bibitem[{Murphy(2022)}]{pml1Book_LCAL}
Murphy, K.~P., 2022. Probabilistic Machine Learning: An introduction. MIT
  Press, Cambridge MA.

\bibitem[{Murphy(2023)}]{pml2Book_LCAL}
Murphy, K.~P., 2023. Probabilistic Machine Learning: Advanced Topics. MIT
  Press, Cambridge MA.

\bibitem[{Murty and Kabadi(1987)}]{MK87_EAY}
Murty, K.~G., Kabadi, S.~N., 1987. Some {NP}-complete problems in quadratic and
  nonlinear programming. Mathematical Programming 39~(2), 117--129.

\bibitem[{Musliu(2006)}]{musliu_GVBSP}
Musliu, N., 2006. Heuristic methods for automatic rotating workforce
  scheduling. International Journal of Computational Intelligence Research 2,
  309--326.

\bibitem[{Mustajoki and Marttunen(2017)}]{Mustajoki2017-eq_JL}
Mustajoki, J., Marttunen, M., 2017. Comparison of multi-criteria decision
  analytical software for supporting environmental planning processes.
  Environmental Modelling \& Software 93, 78--91.

\bibitem[{Mutha and Bansal(2023)}]{Mutha2022-ry_AM}
Mutha, A., Bansal, S., 2023. Determining assortments of used products for {B2B}
  transactions in reverse supply chain. IISE Transactions.

\bibitem[{Mutha et~al.(2016)Mutha, Bansal, and Guide}]{Mutha2016-dc_AM}
Mutha, A., Bansal, S., Guide, V. D.~R., 2016. Managing demand uncertainty
  through core acquisition in remanufacturing. Production and Operations
  Management 25~(8), 1449--1464.

\bibitem[{Mutha et~al.(2019)Mutha, Bansal, and Guide}]{Mutha2019-nv_AM}
Mutha, A., Bansal, S., Guide, V. D.~R., 2019. Selling assortments of used
  products to third‐party remanufacturers. Production and Operations
  Management 28~(7), 1792--1817.

\bibitem[{Nagamochi et~al.(1994)Nagamochi, Ono, and
  Ibaraki}]{Nagamochi-et-al:1994_IL}
Nagamochi, H., Ono, T., Ibaraki, T., 1994. Implementing an efficient minimum
  capacity cut algorithm. Mathematical Programming 67, 325--341.

\bibitem[{Nahmias(1979)}]{nahmias1979simple_JSS}
Nahmias, S., 1979. Simple approximations for a variety of dynamic leadtime
  lost-sales inventory models. Operations Research 27~(5), 904--924.

\bibitem[{Nahmias(2011)}]{nahmias2011perishable_JSS}
Nahmias, S., 2011. Perishable inventory systems. Vol. 160 of International
  Series in Operations Research \& Management Science. Springer Science \&
  Business Media.

\bibitem[{Narasimhan et~al.(2006)Narasimhan, Swink, and
  Kim}]{Narasimhan2006-qt_KESXZ}
Narasimhan, R., Swink, M., Kim, S.~W., 2006. Disentangling leanness and
  agility: An empirical investigation. Journal of Operations Management 24~(5),
  440--457.

\bibitem[{Nash(1950{\natexlab{a}})}]{Nas50_JH}
Nash, J., 1950{\natexlab{a}}. The bargaining problem. Econometrica 18,
  155--162.

\bibitem[{Nash(1950{\natexlab{b}})}]{Na1950_GZ}
Nash, J., 1950{\natexlab{b}}. Equilibrium points in $n$-person games.
  Proceedings of the National Academy of Sciences USA 36~(1), 48--49.

\bibitem[{Nash(1951)}]{Na1951_GZ}
Nash, J., 1951. Non-cooperative games. Annals of Mathematics 54~(2), 286--295.

\bibitem[{{National Grid ESO}(2022)}]{NG22_BC}
{National Grid ESO}, 2022. Short term operating reserve.
  https://www.nationalgrideso.com/industry-information/balancing-services/%
  reserve-services/short-term-operating-reserve, accessed on 2022-11-01.

\bibitem[{Naylor et~al.(1999)Naylor, Naim, and Berry}]{Naylor1999_SMD}
Naylor, J.~B., Naim, M.~M., Berry, D., 1999. Leagility: Integrating the lean
  and agile manufacturing paradigms in the total supply chain. International
  Journal of Production Economics 62~(1), 107--118.

\bibitem[{Nehring(2007)}]{Nehring2007-zh_TC}
Nehring, K., 2007. The impossibility of a {P}aretian rational: A {B}ayesian
  perspective. Economics Letters 96~(1), 45--50.

\bibitem[{Nemhauser and Wolsey(1988)}]{NW88_SMPT}
Nemhauser, G., Wolsey, L., 1988. Integer and combinatorial optimization. John
  Wiley \& Sons, Chichester.

\bibitem[{Nesterov and Nemirovskii(1994)}]{NN94_EAY}
Nesterov, Y.~E., Nemirovskii, A., 1994. Interior-Point Polynomial Algorithms in
  Convex Programming. Vol.~13 of {SIAM} Studies in Applied Mathematics. SIAM,
  Philadelphia, PA.

\bibitem[{Neumann et~al.(2003)Neumann, Schwindt, and
  Zimmermann}]{Neumann2002-ve_WH_ED}
Neumann, K., Schwindt, C., Zimmermann, J., 2003. Project Scheduling with Time
  Windows and Scarce Resources: Temporal and Resource-constrained Project
  Scheduling with Regular and Nonregular Objective Functions. Springer.

\bibitem[{Neumann and Zimmermann(2000)}]{Neumann2000-yl_WH_ED}
Neumann, K., Zimmermann, J., 2000. Procedures for resource leveling and net
  present value problems in project scheduling with general temporal and
  resource constraints. European Journal of Operational Research 127~(2),
  425--443.

\bibitem[{Neuts(1981)}]{neuts_HATI}
Neuts, M.~F., 1981. Matrix-Geometric Solutions in Stochastic Models: An
  Algorithmic Approach. Vol.~2 of Johns Hopkins Series in the Mathematical
  Sciences. Johns Hopkins University Press.

\bibitem[{Neuts(1989)}]{neuts2_HATI}
Neuts, M.~F., 1989. Structured Stochastic Matrices of {M/G/1} Type and Their
  Applications. Probability: Pure and Applied 5. Marcel Dekker, New York.

\bibitem[{Newbold(1998)}]{Newbold1998-xg_WH_ED}
Newbold, R.~C., 1998. Project Management in the Fast Lane: Applying the Theory
  of Constraints. CRC Press.

\bibitem[{Newman et~al.(2010)Newman, Rubio, Caro, Weintraub, and
  Eurek}]{Newman2010_EA}
Newman, A.~M., Rubio, E., Caro, R., Weintraub, A., Eurek, K., 2010. A review of
  operations research in mine planning. Interfaces 40~(3), 222--245.

\bibitem[{Nicholls(2009)}]{Nicholls2009-vc_JJ}
Nicholls, M.~G., 2009. The use of markov models as an aid to the evaluation,
  planning and benchmarking of doctoral programs. Journal of the Operational
  Research Society 60~(9), 1183--1190.

\bibitem[{Niedermeier(2006)}]{Nied06_UPCT}
Niedermeier, R., 2006. Invitation to Fixed-Parameter Algorithms. Oxford
  University Press, Oxford.

\bibitem[{Nifakos et~al.(2021)Nifakos, Chandramouli, Nikolaou, Papachristou,
  Koch, Panaousis, and Bonacina}]{Nifakos2021-jk_TC}
Nifakos, S., Chandramouli, K., Nikolaou, C.~K., Papachristou, P., Koch, S.,
  Panaousis, E., Bonacina, S., 2021. Influence of human factors on cyber
  security within healthcare organisations: A systematic review. Sensors
  21~(15).

\bibitem[{Nikolopoulos et~al.(2015)Nikolopoulos, Litsa, Petropoulos,
  Bougioukos, and Khammash}]{Nikolopoulos2015-mz_FP}
Nikolopoulos, K., Litsa, A., Petropoulos, F., Bougioukos, V., Khammash, M.,
  2015. Relative performance of methods for forecasting special events. Journal
  of Business Research 68~(8), 1785--1791.

\bibitem[{Nisan(2007)}]{Nisan07_BC}
Nisan, N., 2007. Introduction to mechanism design. In: Nisan, N., Roughgarden,
  T., Tardos, E., Vazirani, V.~V. (Eds.), Algorithmic Game Theory. Cambridge
  University Press, Cambridge, pp. 209--242.

\bibitem[{Nishihara(2012)}]{Nishihara2012_EA}
Nishihara, M., 2012. Real options with synergies: Static versus dynamic
  policies. Journal of the Operational Research Society 63~(1), 107--121.

\bibitem[{Niu and Zhou(2013)}]{Niu2013_DC}
Niu, H., Zhou, X., 2013. Optimizing urban rail timetable under time-dependent
  demand and oversaturated conditions. Transportation Research Part C: Emerging
  Technologies 36, 212--230.

\bibitem[{Niu et~al.(2015)Niu, Zhou, and Gao}]{Niu2015_DC}
Niu, H., Zhou, X., Gao, R., 2015. Train scheduling for minimizing passenger
  waiting time with time-dependent demand and skip-stop patterns: Nonlinear
  integer programming models with linear constraints. Transportation Research
  Part B: Methodological 76, 117--135.

\bibitem[{Nocedal and Wright(2006)}]{JW06_EAY}
Nocedal, J., Wright, S.~J., 2006. Numerical optimization, 2nd Edition. Springer
  series in Operations Research and Financial Engineering. Springer, New York,
  NY.

\bibitem[{Norlund and Gribkovskaia(2013)}]{Norlund2013-qg_HP}
Norlund, E.~K., Gribkovskaia, I., 2013. Reducing emissions through speed
  optimization in supply vessel operations. Transportation Research Part D:
  Transport and Environment 23, 105--113.

\bibitem[{Nowotarski et~al.(2016)Nowotarski, Liu, Weron, and
  Hong}]{now:liu:wer:hon:16_DSRW}
Nowotarski, J., Liu, B., Weron, R., Hong, T., 2016. Improving short term load
  forecast accuracy via combining sister forecasts. Energy 98, 40--49.

\bibitem[{Nowotarski and Weron(2018)}]{now:wer:18_DSRW}
Nowotarski, J., Weron, R., 2018. Recent advances in electricity price
  forecasting: {A} review of probabilistic forecasting. Renewable and
  Sustainable Energy Reviews 81, 1548--1568.

\bibitem[{N{\"u}bel(2001)}]{Nubel2001-ou_WH_ED}
N{\"u}bel, H., 2001. The resource renting problem subject to temporal
  constraints. OR Spectrum 23~(3), 359--381.

\bibitem[{Numrich and Picucci(2012)}]{Numrich2012-ov_KVRH}
Numrich, S.~K., Picucci, P.~M., 2012. New challenges: Human, social, cultural,
  and behavioral modeling. In: Tolk, A. (Ed.), Engineering Principles of Combat
  Modeling and Distributed Simulation. Wiley, Hoboken, NJ, pp. 641--667.

\bibitem[{Oesterreich et~al.(2022)Oesterreich, Anton, and
  Teuteberg}]{Oesterreich2022-ez_JEB}
Oesterreich, T.~D., Anton, E., Teuteberg, F., 2022. What translates big data
  into business value? {A} meta-analysis of the impacts of business analytics
  on firm performance. Information \& Management 59~(6), 103685.

\bibitem[{{Office for National
  Statistics}(2022{\natexlab{a}})}]{Office_for_National_Statistics2022-ey}
{Office for National Statistics}, 2022{\natexlab{a}}. Coronavirus ({COVID-19})
  infection survey: methods and further information. Tech. rep., Office for
  National Statistics.

\bibitem[{{Office for National
  Statistics}(2022{\natexlab{b}})}]{Office_for_National_Statistics2022-rx}
{Office for National Statistics}, 2022{\natexlab{b}}. Coronavirus ({COVID-19})
  latest insights. Tech. rep., Office for National Statistics.

\bibitem[{{Office for National
  Statistics}(2022{\natexlab{c}})}]{Office_for_National_Statistics2022-gy}
{Office for National Statistics}, 2022{\natexlab{c}}. {COVID-19} schools
  infection survey, england statistical bulletins. Tech. rep., Office for
  National Statistics.

\bibitem[{Ogata et~al.(2010)}]{Ogata2010_XW}
Ogata, K., et~al., 2010. Modern control engineering, 5th Edition. Prentice
  Hall, Upper Saddle River, NJ.

\bibitem[{Ogryczak et~al.(2014)Ogryczak, Luss, {Pi\'{o}ro}, Nace, and
  Tomaszewski}]{OgrLusPioNacTom14_JH}
Ogryczak, W., Luss, H., {Pi\'{o}ro}, M., Nace, D., Tomaszewski, A., 2014. Fair
  optimization and networks: {A} survey. Journal of Applied Mathematics 2014,
  1--25.

\bibitem[{Ogryczak and {\'S}liwi{\'n}ski(2003)}]{ogryczak2003solving_JH}
Ogryczak, W., {\'S}liwi{\'n}ski, T., 2003. On solving linear programs with the
  ordered weighted averaging objective. European Journal of Operational
  Research 148~(1), 80--91.

\bibitem[{{\"O}hman et~al.(2021){\"O}hman, Hiltunen, Virtanen, and
  Holmström}]{ohman_frontlog_2021_JHLS}
{\"O}hman, M., Hiltunen, M., Virtanen, K., Holmström, J., 2021. Frontlog
  scheduling in aircraft line maintenance: {From} explorative solution design
  to theoretical insight into buffer management. Journal of Operations
  Management 67~(2), 120--151.

\bibitem[{Ohno(1988)}]{Ohno1988_SMD}
Ohno, T., 1988. Toyota {P}roduction {S}ystem: Beyond Large-Scale Production.
  Productivity Press, Portland, OR.

\bibitem[{Olesen et~al.(2022)Olesen, Petersen, and
  Podinovski}]{Olesen2022-is_SL}
Olesen, O.~B., Petersen, N.~C., Podinovski, V.~V., 2022. The structure of
  production technologies with ratio inputs and outputs. Journal of
  Productivity Analysis 57~(3), 255--267.

\bibitem[{Olhager(2010)}]{Olhager2010_SMD}
Olhager, J., 2010. The role of the customer order decoupling point in
  production and supply chain management. Computers in Industry 61~(9),
  863--868.

\bibitem[{Olivares et~al.(2023)Olivares, Challu, Marcjasz, Weron, and
  Dubrawski}]{oli:cha:mar:wer:dub:22_DSRW}
Olivares, K., Challu, C., Marcjasz, G., Weron, R., Dubrawski, A., 2023. Neural
  basis expansion analysis with exogenous variables: Forecasting electricity
  prices with {NBEATSx}. International Journal of Forecasting 39~(2), 884--900.

\bibitem[{Oliveira and Ferreira(1990)}]{Oliveira1990-cd_JB}
Oliveira, J., Ferreira, J., 1990. An improved version of {W}ang's algorithm for
  two-dimensional cutting problems. European Journal of Operational Research
  44~(2), 256--266.

\bibitem[{O'Neil(2016)}]{ONeil2016-rs_JEB}
O'Neil, C., 2016. Weapons of Math Destruction: How Big Data Increases
  Inequality and Threatens Democracy. Crown Publishing Group, New York.

\bibitem[{Onggo(2019)}]{Onggo_chapter_2019_CC}
Onggo, B.~S., 2019. {Symbiotic Simulation System (S3) for Industry 4.0}. In:
  Gunal, M. (Ed.), Simulation for Industry 4.0: Past, Present and Future.
  Springer, Switzerland, pp. 153--165.

\bibitem[{Opricovic and Tzeng(2004)}]{opricovic2004compromise}
Opricovic, S., Tzeng, G.-H., 2004. Compromise solution by mcdm methods: A
  comparative analysis of vikor and topsis. European Journal of Operational
  Research 156~(2), 445--455.

\bibitem[{Ord et~al.(2017)Ord, Fildes, and Kourentzes}]{FildesOrd2017_FP}
Ord, K., Fildes, R., Kourentzes, N., 2017. Principles of Business Forecasting,
  2nd Edition. Wessex Press Inc.

\bibitem[{Oreshkin et~al.(2021)Oreshkin, Dudek, Pełka, and
  Turkina}]{ore:dud:pel:tur:21_DSRW}
Oreshkin, B., Dudek, G., Pełka, P., Turkina, E., 2021. {N-BEATS} neural
  network for mid-term electricity load forecasting. Applied Energy 293,
  116918.

\bibitem[{Orlin(1993)}]{Orlin:1993_IL}
Orlin, J.~B., 1993. A faster strongly polynomial minimum cost flow algorithm.
  Operations Research 41~(2), 338--350.

\bibitem[{Orlin(2013)}]{Orlin:2013_IL}
Orlin, J.~B., 2013. Max flows in {$O(nm)$} time, or better. In: Boneh, D.,
  Roughgarden, T., Feigenbaum, J. (Eds.), Symposium on Theory of Computing
  Conference, STOC'13, Palo Alto, CA, USA, June 1-4, 2013. pp. 765--774.

\bibitem[{Ormerod et~al.(2023{\natexlab{a}})Ormerod, Yearworth, and
  White}]{Ormerod2023-qp}
Ormerod, R., Yearworth, M., White, L., 2023{\natexlab{a}}. Understanding
  participant actions in {OR} interventions using practice theories: A research
  agenda. European Journal of Operational Research 306~(2), 810--827.

\bibitem[{Ormerod(2014{\natexlab{a}})}]{Ormerod2014-ax_MYLW}
Ormerod, R.~J., 2014{\natexlab{a}}. The mangle of {OR} practice: towards more
  informative case studies of `technical' projects. Journal of the Operational
  Research Society 65~(8), 1245--1260.

\bibitem[{Ormerod(2014{\natexlab{b}})}]{Ormerod2014-kb_MYLW}
Ormerod, R.~J., 2014{\natexlab{b}}. {OR} competences: the demands of problem
  structuring methods. EURO Journal on Decision Processes 2~(3), 313--340.

\bibitem[{Ormerod and Ulrich(2013)}]{Ormerod2013-km_JH}
Ormerod, R.~J., Ulrich, W., 2013. Operational research and ethics: A literature
  review. European Journal of Operational Research 228~(2), 291--307.

\bibitem[{Ormerod et~al.(2023{\natexlab{b}})Ormerod, Yearworth, and
  White}]{Ormerod2022-qp_MYLW}
Ormerod, R.~J., Yearworth, M., White, L., 2023{\natexlab{b}}. Understanding
  participant actions in {OR} interventions using practice theories: A research
  agenda. European Journal of Operational Research 306~(2), 810--827.

\bibitem[{Osborne and Rubinstein(1994)}]{Osru1994_GZ}
Osborne, M., Rubinstein, A., 1994. A Course in Game Theory. {MIT} Press,
  Cambridge, MA.

\bibitem[{Otto et~al.(2018)Otto, Agatz, Campbell, Golden, and
  Pesch}]{OttoACGP2018_CA_MB}
Otto, A., Agatz, N., Campbell, J., Golden, B., Pesch, E., 2018. Optimization
  approaches for civil applications of unmanned aerial vehicles ({UAV}s) or
  aerial drones: A survey. Networks 72~(4), 411--458.

\bibitem[{Oude~Vrielink et~al.(2019)Oude~Vrielink, Jansen, Hans, and van
  Hillegersberg}]{Oude_Vrielink2019-qh_JJ}
Oude~Vrielink, R.~A., Jansen, E.~A., Hans, E.~W., van Hillegersberg, J., 2019.
  Practices in timetabling in higher education institutions: a systematic
  review. Annals of Operations Research 275~(1), 145--160.

\bibitem[{Ovchinnikov(2011)}]{Ovchinnikov2011-qd_AM}
Ovchinnikov, A., 2011. Revenue and cost management for remanufactured products.
  Production and Operations Management 20~(6), 824--840.

\bibitem[{Owen(1973)}]{Ow73_ALAL}
Owen, G., 1973. Cutting planes for programs with disjunctive constraints.
  Journal of Optimization Theory and Applications 11, 49--55.

\bibitem[{Owen(1995)}]{Ow1995_GZ}
Owen, G., 1995. Game theory, 3rd Edition. Academic Press, San Diego.

\bibitem[{{\"O}zdemir-Aky{\i}ld{\i}r{\i}m
  et~al.(2014){\"O}zdemir-Aky{\i}ld{\i}r{\i}m, Denizel, and
  Ferguson}]{Ozdemir-Akyildirim2014-by_AM}
{\"O}zdemir-Aky{\i}ld{\i}r{\i}m, {\"O}., Denizel, M., Ferguson, M., 2014.
  Allocation of returned products among different recovery options through an
  opportunity cost-based dynamic approach. Decision Sciences 45~(6),
  1083--1116.

\bibitem[{O’Hanley and Church(2011)}]{OHanley2011_JLYHK}
O’Hanley, J.~R., Church, R.~L., 2011. Designing robust coverage networks to
  hedge against worst-case facility losses. European Journal of Operational
  Research 209~(1), 23--36.

\bibitem[{Padberg and Rinaldi(1987)}]{PR87_ALAL}
Padberg, M., Rinaldi, G., 1987. Optimization of a 532 city symmetric traveling
  salesman problem by branch-and-cut. Operations Research Letters 6, 1--7.

\bibitem[{Padberg and Rinaldi(1990{\natexlab{a}})}]{Padberg-Rinaldi:1990_IL}
Padberg, M., Rinaldi, G., 1990{\natexlab{a}}. An efficient algorithm for the
  minimum capacity cut problem. Mathematical Programming 47, 19--36.

\bibitem[{Padberg and Rinaldi(1990{\natexlab{b}})}]{Padberg-Rinaldi:1990a_IL}
Padberg, M., Rinaldi, G., 1990{\natexlab{b}}. Facet identification for the
  symmetric traveling salesman polytope. Mathematical Programming 47, 219--257.

\bibitem[{Padberg and Rinaldi(1991)}]{Padberg-Rinaldi:1991_IL}
Padberg, M., Rinaldi, G., 1991. A branch-and-cut algorithm for the resolution
  of large-scale symmetric traveling salesman problems. {SIAM} Review 33~(1),
  60--100.

\bibitem[{Padberg and Sung(1991)}]{Padberg-Sung:1991_IL}
Padberg, M., Sung, T., 1991. An analytical comparison of different formulations
  of the travelling salesman problem. Mathematical Programming 52, 315--357.

\bibitem[{Padberg et~al.(1984)Padberg, Van~Roy, and Wolsey}]{PVW84_ALAL}
Padberg, M., Van~Roy, T., Wolsey, L., 1984. Valid inequalities for fixed charge
  problems. Operations Research 32, 842--861.

\bibitem[{Pagel and Yates(2022)}]{Pagel2022-vg_CV}
Pagel, C., Yates, C.~A., 2022. Role of mathematical modelling in future
  pandemic response policy. BMJ 378.

\bibitem[{Pahl-Wostl(2007)}]{Pahl-Wostl2007-br_MYLW}
Pahl-Wostl, C., 2007. The implications of complexity for integrated resources
  management. Environmental Modelling \& Software 22~(5), 561--569.

\bibitem[{Paivio(1990)}]{Paivio1990-mw_MJE}
Paivio, A., 1990. Mental Representations: A dual coding approach. Oxford
  University Press.

\bibitem[{Pajor et~al.(2018)Pajor, Uchoa, and Werneck}]{Pajor:2018_IL}
Pajor, T., Uchoa, E., Werneck, R.~F., 2018. A robust and scalable algorithm for
  the {S}teiner problem in graphs. Mathematical Programming Computation 10~(1),
  69--118.

\bibitem[{Pang and Whitt(2010)}]{pang2010two_HATI}
Pang, G., Whitt, W., 2010. Two-parameter heavy-traffic limits for
  infinite-server queues. Queueing Systems 65~(4), 325--364.

\bibitem[{Pantuso et~al.(2014)Pantuso, Fagerholt, and
  Hvattum}]{Pantuso2014-bk_HP}
Pantuso, G., Fagerholt, K., Hvattum, L.~M., 2014. A survey on maritime fleet
  size and mix problems. European Journal of Operational Research 235~(2),
  341--349.

\bibitem[{Papadimitriou and Steiglitz(1982)}]{PS82_SMPT}
Papadimitriou, C., Steiglitz, K., 1982. Combinatorial Optimization: Algorithms
  and Complexity. Prentice Hall, Englewood Cliffs, NJ.

\bibitem[{Papadimitriou and Johnes(2019)}]{Papadimitriou2019-ap_JJ}
Papadimitriou, M., Johnes, J., 2019. Does merging improve efficiency? {A} study
  of {E}nglish universities. Studies in Higher Education 44~(8), 1454--1474.

\bibitem[{Pardalos and Vavasis(1991)}]{PV91_EAY}
Pardalos, P.~M., Vavasis, S.~A., 1991. Quadratic programming with one negative
  eigenvalue is {NP}-hard. Journal of Global Optimization 1~(1), 15--22.

\bibitem[{Parilina and Zaccour(2015)}]{Paza2015_GZ}
Parilina, E., Zaccour, G., 2015. Approximated cooperative equilibria for games
  played over event trees. Operations Research Letters 43, 507--513.

\bibitem[{Parry and Mingers(1991)}]{Parry1991-px_MYLW}
Parry, R., Mingers, J., 1991. Community operational research: Its context and
  its future. Omega 19~(6), 577--586.

\bibitem[{Partyka and Hall(2014)}]{PartykaH2014_CA_MB}
Partyka, J., Hall, R., 2014. Vehicle routing software survey: Rising customer
  expectations drive software innovation, integration. OR/MS Today 41.

\bibitem[{Pastor and Lovell(2005)}]{Pastor2005-rg_SL}
Pastor, J.~T., Lovell, C. A.~K., 2005. A global {M}almquist productivity index.
  Economics Letters 88~(2), 266--271.

\bibitem[{Paszke et~al.(2019)Paszke, Gross, Massa, Lerer, Bradbury, Chanan,
  Killeen, Lin, Gimelshein, Antiga, Desmaison, Kopf, Yang, DeVito, Raison,
  Tejani, Chilamkurthy, Steiner, Fang, Bai, and
  Chintala}]{PyToch_NEURIPS2019_9015_LCAL}
Paszke, A., Gross, S., Massa, F., Lerer, A., Bradbury, J., Chanan, G., Killeen,
  T., Lin, Z., Gimelshein, N., Antiga, L., Desmaison, A., Kopf, A., Yang, E.,
  DeVito, Z., Raison, M., Tejani, A., Chilamkurthy, S., Steiner, B., Fang, L.,
  Bai, J., Chintala, S., 2019. Pytorch: An imperative style, high-performance
  deep learning library. In: Advances in Neural Information Processing Systems
  32 (NeurIPS 2019). Curran Associates, Inc., pp. 8024--8035.

\bibitem[{Paucar-Caceres et~al.(2022)Paucar-Caceres, Cavalcanti-Bandos,
  Quispe-Prieto, Huerta-Tantalean, and
  Werner-Masters}]{Paucar-Caceres2022-dy_AJG}
Paucar-Caceres, A., Cavalcanti-Bandos, M.~F., Quispe-Prieto, S.~C.,
  Huerta-Tantalean, L.~N., Werner-Masters, K., 2022. Using soft systems
  methodology to align community projects with sustainability development in
  higher education stakeholders' networks in a {B}razilian university. Systems
  Research and Behavioral Science 39~(4), 750--764.

\bibitem[{Paucar-Caceres and Espinosa(2011)}]{paucar2011_EA}
Paucar-Caceres, A., Espinosa, A., 2011. Management science methodologies in
  environmental management and sustainability: Discourses and applications.
  Journal of the Operational Research Society 62~(9), 1601--1620.

\bibitem[{Paucar-Caceres and Thorpe(2005)}]{Paucar-Caceres2005-kz_JJ}
Paucar-Caceres, A., Thorpe, R., 2005. Mapping the structure of {MBA}
  programmes: a comparative study of the structure of accredited {AMBA}
  programmes in the {U}nited {K}ingdom. Journal of the Operational Research
  Society 56~(1), 25--38.

\bibitem[{Paul and Knust(2015)}]{paul2015classification_GVBSP}
Paul, M., Knust, S., 2015. A classification scheme for integrated staff
  rostering and scheduling problems. RAIRO-Operations Research 49~(2),
  393--412.

\bibitem[{Pearl(1988)}]{pearl1988probabilistic_MESG}
Pearl, J., 1988. Probabilistic Reasoning in Intelligent Systems: Networks of
  Plausible Inference. Morgan Kaufmann, Burligton.

\bibitem[{Pecin et~al.(2017)Pecin, Pessoa, Poggi, and Uchoa}]{Pecin2017_CA_MB}
Pecin, D., Pessoa, A., Poggi, M., Uchoa, E., 2017. Improved
  branch-cut-and-price for capacitated vehicle routing. Mathematical
  Programming Computation 9~(1), 61--100.

\bibitem[{Pedraza-Martinez and
  Van~Wassenhove(2016)}]{pedraza2016empirically_BYKOK}
Pedraza-Martinez, A.~J., Van~Wassenhove, L.~N., 2016. Empirically grounded
  research in humanitarian operations management: The way forward. Journal of
  Operations Management 45~(1), 1--10.

\bibitem[{Pedregosa et~al.(2011)Pedregosa, Varoquaux, Gramfort, Michel,
  Thirion, Grisel, Blondel, Prettenhofer, Weiss, Dubourg, Vanderplas, Passos,
  Cournapeau, Brucher, Perrot, and Duchesnay}]{scikit-learn_LCAL}
Pedregosa, F., Varoquaux, G., Gramfort, A., Michel, V., Thirion, B., Grisel,
  O., Blondel, M., Prettenhofer, P., Weiss, R., Dubourg, V., Vanderplas, J.,
  Passos, A., Cournapeau, D., Brucher, M., Perrot, M., Duchesnay, E., 2011.
  Scikit-learn: Machine learning in {P}ython. Journal of Machine Learning
  Research 12, 2825--2830.

\bibitem[{Peeters and Kroon(2008)}]{Peeters2008_DC}
Peeters, M., Kroon, L., 2008. Circulation of railway rolling stock: a
  branch-and-price approach. Computers \& Operations Research 35~(2), 538--556.

\bibitem[{Peeters et~al.(2020)Peeters, Salaga, and
  Juravich}]{Peeters2020-yb_IM}
Peeters, T. L. P.~R., Salaga, S., Juravich, M., 2020. Matching and winning?
  {T}he impact of upper and middle managers on firm performance in major league
  baseball. Management Science 66~(6), 2735--2751.

\bibitem[{Pelletier et~al.(2016)Pelletier, Jabali, and
  Laporte}]{PelletierJL2016_CA_MB}
Pelletier, S., Jabali, O., Laporte, G., 2016. 50th anniversary invited
  article—goods distribution with electric vehicles: Review and research
  perspectives. Transportation Science 50~(1), 3--22.

\bibitem[{Pellizzoni and Ungaro(2000)}]{Pellizzoni2000-fl_TC}
Pellizzoni, L., Ungaro, D., 2000. Technological risk, participation and
  deliberation. some results from three {I}talian case studies. Journal of
  Hazardous Materials 78~(1-3), 261--280.

\bibitem[{Perakis and Papadakis(1989)}]{Perakis1989-zi_HP}
Perakis, A.~N., Papadakis, N.~A., 1989. Minimal time vessel routing in a
  {Time-Dependent} environment. Transportation Science 23~(4), 266--276.

\bibitem[{Perera et~al.(2016)Perera, Davis, and Swartz}]{Perera2016-iw_IM}
Perera, H., Davis, J., Swartz, T.~B., 2016. Optimal lineups in {Twenty20}
  cricket. Journal of Statistical Computation and Simulation 86~(14),
  2888--2900.

\bibitem[{Perera and Sethi(2022{\natexlab{a}})}]{perera2022bsurvey_JSS}
Perera, S.~C., Sethi, S.~P., 2022{\natexlab{a}}. A survey of stochastic
  inventory models with fixed costs: Optimality of (s, {S}) and (s, {S})-type
  policies--{C}ontinuous-time case. Production and Operations Management, DOI:
  10.1111/poms.13819.

\bibitem[{Perera and Sethi(2022{\natexlab{b}})}]{perera2022asurvey_JSS}
Perera, S.~C., Sethi, S.~P., 2022{\natexlab{b}}. A survey of stochastic
  inventory models with fixed costs: Optimality of (s, {S}) and (s, {S})-type
  policies--{D}iscrete-time case. Production and Operations Management, DOI:
  10.1111/poms.13820.

\bibitem[{Pessoa et~al.(2018)Pessoa, Sadykov, Uchoa, and
  Vanderbeck}]{Pessoa2018_CA_MB}
Pessoa, A., Sadykov, R., Uchoa, E., Vanderbeck, F., 2018. Automation and
  combination of linear-programming based stabilization techniques in column
  generation. {INFORMS} Journal on Computing 30~(2), 339--360.

\bibitem[{Pessoa et~al.(2020)Pessoa, Sadykov, Uchoa, and
  Vanderbeck}]{PessoaSUV2020_CA_MB}
Pessoa, A., Sadykov, R., Uchoa, E., Vanderbeck, F., 2020. A generic exact
  solver for vehicle routing and related problems. Mathematical Programming
  183~(1), 483--523.

\bibitem[{Petersen et~al.(2012)Petersen, S{\"o}lveling, Clarke, Johnson, and
  Shebalov}]{petersen2012optimization_VLVV}
Petersen, J.~D., S{\"o}lveling, G., Clarke, J.-P., Johnson, E.~L., Shebalov,
  S., 2012. An optimization approach to airline integrated recovery.
  Transportation Science 46~(4), 482--500.

\bibitem[{Petropoulos et~al.(2022)Petropoulos, Apiletti, Assimakopoulos, Babai,
  Barrow, {Ben Taieb}, Bergmeir, Bessa, Bijak, Boylan, Browell, Carnevale,
  Castle, Cirillo, Clements, Cordeiro, {Cyrino Oliveira}, {De Baets},
  Dokumentov, Ellison, Fiszeder, Franses, Frazier, Gilliland, Gönül, Goodwin,
  Grossi, Grushka-Cockayne, Guidolin, Guidolin, Gunter, Guo, Guseo, Harvey,
  Hendry, Hollyman, Januschowski, Jeon, Jose, Kang, Koehler, Kolassa,
  Kourentzes, Leva, Li, Litsiou, Makridakis, Martin, Martinez, Meeran, Modis,
  Nikolopoulos, Önkal, Paccagnini, Panagiotelis, Panapakidis, Pavía, Pedio,
  Pedregal, Pinson, Ramos, Rapach, Reade, Rostami-Tabar, Rubaszek, Sermpinis,
  Shang, Spiliotis, Syntetos, Talagala, Talagala, Tashman, Thomakos,
  Thorarinsdottir, Todini, {Trapero Arenas}, Wang, Winkler, Yusupova, and
  Ziel}]{PETROPOULOS2022705_FP}
Petropoulos, F., Apiletti, D., Assimakopoulos, V., Babai, M.~Z., Barrow, D.~K.,
  {Ben Taieb}, S., Bergmeir, C., Bessa, R.~J., Bijak, J., Boylan, J.~E.,
  Browell, J., Carnevale, C., Castle, J.~L., Cirillo, P., Clements, M.~P.,
  Cordeiro, C., {Cyrino Oliveira}, F.~L., {De Baets}, S., Dokumentov, A.,
  Ellison, J., Fiszeder, P., Franses, P.~H., Frazier, D.~T., Gilliland, M.,
  Gönül, M.~S., Goodwin, P., Grossi, L., Grushka-Cockayne, Y., Guidolin, M.,
  Guidolin, M., Gunter, U., Guo, X., Guseo, R., Harvey, N., Hendry, D.~F.,
  Hollyman, R., Januschowski, T., Jeon, J., Jose, V. R.~R., Kang, Y., Koehler,
  A.~B., Kolassa, S., Kourentzes, N., Leva, S., Li, F., Litsiou, K.,
  Makridakis, S., Martin, G.~M., Martinez, A.~B., Meeran, S., Modis, T.,
  Nikolopoulos, K., Önkal, D., Paccagnini, A., Panagiotelis, A., Panapakidis,
  I., Pavía, J.~M., Pedio, M., Pedregal, D.~J., Pinson, P., Ramos, P., Rapach,
  D.~E., Reade, J.~J., Rostami-Tabar, B., Rubaszek, M., Sermpinis, G., Shang,
  H.~L., Spiliotis, E., Syntetos, A.~A., Talagala, P.~D., Talagala, T.~S.,
  Tashman, L., Thomakos, D., Thorarinsdottir, T., Todini, E., {Trapero Arenas},
  J.~R., Wang, X., Winkler, R.~L., Yusupova, A., Ziel, F., 2022. Forecasting:
  theory and practice. International Journal of Forecasting 38~(3), 705--871.

\bibitem[{Petropoulos et~al.(2016)Petropoulos, Fildes, and
  Goodwin}]{Petropoulos2016-uk_FP}
Petropoulos, F., Fildes, R., Goodwin, P., 2016. Do `big losses' in judgmental
  adjustments to statistical forecasts affect experts' behaviour? European
  Journal of Operational Research 249~(3), 842--852.

\bibitem[{Petropoulos et~al.(2018)Petropoulos, Kourentzes, Nikolopoulos, and
  Siemsen}]{Petropoulos2018-mt_FP}
Petropoulos, F., Kourentzes, N., Nikolopoulos, K., Siemsen, E., 2018.
  Judgmental selection of forecasting models. Journal of Operations Management
  60, 34--46.

\bibitem[{Petropoulos and Makridakis(2020)}]{Petropoulos2020-qm_CV}
Petropoulos, F., Makridakis, S., 2020. Forecasting the novel coronavirus
  {COVID-19}. PLoS One 15~(3), e0231236.

\bibitem[{Petropoulos and Siemsen(2022)}]{Petropoulos2022REP_FP}
Petropoulos, F., Siemsen, E., 2022. Forecast selection and representativeness.
  Management Science, DOI: 10.1287/mnsc.2022.4485.

\bibitem[{Petrosjan(1977)}]{Pe1977_GZ}
Petrosjan, L., 1977. Stable solutions of differential games with many
  participants. Viestnik of Leningrad University 19, 46--52.

\bibitem[{Petrosyan and Zaccour(2018)}]{Peza2018_GZ}
Petrosyan, L., Zaccour, G., 2018. Cooperative differential games with
  transferable payoffs. In: Ba\c{s}ar, T., Zaccour, G. (Eds.), Handbook of
  Dynamic Game Theory. Springer, Cham, pp. 1--38.

\bibitem[{Petrovic and Burke(2004)}]{Petrovic2004-ly_JJ}
Petrovic, S., Burke, E., 2004. University timetabling. In: Leung, J. Y.-T.
  (Ed.), Handbook of Scheduling: Algorithms, Models, and Performance Analysis.
  Chapman and Hall/CRC, Boca Raton, FL, Ch.~45.

\bibitem[{Petrovic et~al.(2020)Petrovic, Parkin, and
  Wrigley}]{well-being_GVBSP}
Petrovic, S., Parkin, J., Wrigley, D., 2020. Personnel scheduling considering
  employee well-being: insights from case studies. In: Proceedings of the 13th
  International Conference on the Practice and Theory of Automated Timetabling
  - PATAT 2021: Volume I. pp. 10--23.

\bibitem[{Petrovic et~al.(2007)Petrovic, Yang, and Dror}]{Petrovic2007-ls_JJ}
Petrovic, S., Yang, Y., Dror, M., 2007. Case-based selection of initialisation
  heuristics for metaheuristic examination timetabling. Expert Systems with
  Applications 33~(3), 772--785.

\bibitem[{Petruzzi and Dada(1999)}]{petruzzi1999pricing_JSS}
Petruzzi, N.~C., Dada, M., 1999. Pricing and the newsvendor problem: A review
  with extensions. Operations Research 47~(2), 183--194.

\bibitem[{Peykani et~al.(2021)Peykani, Farzipoor~Saen, Seyed~Esmaeili, and
  Gheidar-Kheljani}]{Peykani2021-qb_SL}
Peykani, P., Farzipoor~Saen, R., Seyed~Esmaeili, F.~S., Gheidar-Kheljani, J.,
  2021. Window data envelopment analysis approach: A review and bibliometric
  analysis. Expert Systems 38~(7), e12721.

\bibitem[{Peykani et~al.(2020)Peykani, Mohammadi, Saen, Sadjadi, and
  Rostamy-Malkhalifeh}]{Peykani2020-zz_SL}
Peykani, P., Mohammadi, E., Saen, R.~F., Sadjadi, S.~J., Rostamy-Malkhalifeh,
  M., 2020. Data envelopment analysis and robust optimization: A review. Expert
  Systems 37~(4), e12534.

\bibitem[{Phelps and K{\"o}ksalan(2003)}]{phelps2003interactive_MESG}
Phelps, S., K{\"o}ksalan, M., 2003. An interactive evolutionary metaheuristic
  for multiobjective combinatorial optimization. Management Science 49~(12),
  1726--1738.

\bibitem[{Phillips et~al.(2021)Phillips, Hahn, Fontana, Yates, Greene,
  Broniatowski, and Przybocki}]{XAI_NIST_21}
Phillips, P.~J., Hahn, C., Fontana, P., Yates, A., Greene, K.~K., Broniatowski,
  D., Przybocki, M.~A., 2021-09-29 04:09:00 2021. Four principles of
  explainable artificial intelligence.
\newline\urlprefix\url{https://tsapps.nist.gov/publication/get_pdf.cfm?pub_id=933399}

\bibitem[{Pidd(2009)}]{Pidd2009-ny_JEB}
Pidd, M., 2009. Tools for Thinking: Modelling in Management Science. Wiley.

\bibitem[{Pillac et~al.(2013)Pillac, Guéret, and Medaglia}]{Pillac2013_COIT}
Pillac, V., Guéret, C., Medaglia, A.~L., 2013. A parallel matheuristic for the
  technician routing and scheduling problem. Optimization Letters 7,
  1525--1535.

\bibitem[{Pillay(2016)}]{Pillay2016-pp_JJ}
Pillay, N., 2016. A review of hyper-heuristics for educational timetabling.
  Annals of Operations Research 239~(1), 3--38.

\bibitem[{Pin{\c c}e et~al.(2016)Pin{\c c}e, Ferguson, and
  Toktay}]{Pince2016-ol_AM}
Pin{\c c}e, {\c C}., Ferguson, M., Toktay, B., 2016. Extracting maximum value
  from consumer returns: Allocating between remarketing and refurbishing for
  warranty claims. Manufacturing \& Service Operations Management 18~(4),
  475--492.

\bibitem[{Pinedo(2012)}]{P12_SMPT}
Pinedo, M., 2012. Scheduling. Springer, Berlin.

\bibitem[{Pinzon-Salcedo and Torres-Cuello(2022)}]{Pinzon-Salcedo2022-im_AJG}
Pinzon-Salcedo, L.~A., Torres-Cuello, M.~A., 2022. Systems thinking concepts
  within a collaborative programme evaluation methodology: {The Hermes
  Programme} evaluation. Systems Research and Behavioral Science 39~(4),
  708--722.

\bibitem[{Pi\'{o}ro and Medhi(2004)}]{pioro.medhi:04_BF}
Pi\'{o}ro, M., Medhi, D., 2004. Routing, Flow, and Capacity Design in
  Communication and Computer Networks. Morgan Kaufman.

\bibitem[{Pi{\'o}ro et~al.(2000)Pi{\'o}ro, Szentesi, Harmatos, and
  J{\"u}ttner}]{pioro.szentesi.ea:00_BF}
Pi{\'o}ro, M., Szentesi, A., Harmatos, J., J{\"u}ttner, A., 2000. On {OSPF}
  related network optimization problems. In: {8th IFIP Workshop on Performance
  Modelling and Evaluation of ATM \& IP Networks}. Ilkley, UK, pp. 70/1--70/14.

\bibitem[{Pirab{\'a}n et~al.(2019)Pirab{\'a}n, Guerrero, and
  Labadie}]{Piraban2019_JLYHK}
Pirab{\'a}n, A., Guerrero, W.~J., Labadie, N., 2019. Survey on blood supply
  chain management: Models and methods. Computers and Operations Research 112,
  104756.

\bibitem[{Pirkul(1987)}]{pirkul:87_BF}
Pirkul, H., 1987. Efficient algorithms for the capacitated concentrator
  location problem. Computers \& Operations Research 14, 197 -- 208.

\bibitem[{Pishvaee et~al.(2011)Pishvaee, Rabbani, and
  Torabi}]{Pishvaee2011-hu_HL}
Pishvaee, M.~S., Rabbani, M., Torabi, S.~A., 2011. A robust optimization
  approach to closed-loop supply chain network design under uncertainty.
  Applied Mathematical Modelling 35~(2), 637--649.

\bibitem[{Pisinger and Ropke(2007)}]{PISINGER20072403_CA_MB}
Pisinger, D., Ropke, S., 2007. A general heuristic for vehicle routing
  problems. Computers \& Operations Research 34~(8), 2403--2435.

\bibitem[{Pisinger and Sigurd(2007)}]{Pisinger2007-rx_JB}
Pisinger, D., Sigurd, M., 2007. Using decomposition techniques and constraint
  programming for solving the {Two-Dimensional} {Bin-Packing} problem. INFORMS
  Journal on Computing 19~(1), 36--51.

\bibitem[{Plastria(2002)}]{plastria2002continuous_SAA}
Plastria, F., 2002. Continuous covering location problems. In: Drezner, Z.,
  Hamacher, H. (Eds.), Facility Location: Applications and Theory. Springer,
  Berlin, pp. 37--79.

\bibitem[{Plà et~al.(2014)Plà, Sandars, and Higgins}]{pla2014_EA}
Plà, L.~M., Sandars, D.~L., Higgins, A.~J., 2014. A perspective on operational
  research prospects for agriculture. Journal of the Operational Research
  Society 65~(7), 1078--1089.

\bibitem[{Polzin(2003)}]{Polzin:2003_IL}
Polzin, T., 2003. Algorithms for the {{S}teiner} problem in networks. Ph.D.
  thesis, Saarland University, Saarbr\"ucken, Germany.

\bibitem[{Polzin and Vahdati~Daneshmand(2009)}]{Polzin:2009approaches_IL}
Polzin, T., Vahdati~Daneshmand, S., 2009. Approaches to the {{S}teiner} problem
  in networks. In: Lerner, J., Wagner, D., Zweig, K.~A. (Eds.), Algorithmics of
  Large and Complex Networks. {Springer}, {Boston, MA}, pp. 81--103.

\bibitem[{Poole(2004)}]{Poole2004-cq_AFRH}
Poole, M.~S., 2004. Central issues in the study of change and innovation. In:
  Poole, M.~S., Van~de Ven, A.~H. (Eds.), Handbook of organizational change and
  innovation. Oxford University Press, New York, pp. 3--31.

\bibitem[{Poole(2007)}]{Poole2007-cg_AFRH}
Poole, M.~S., 2007. Generalization in process theories of communication.
  Communication Methods and Measures 1~(3), 181--190.

\bibitem[{Poropudas and Virtanen(2010)}]{Poropudas2010-jx_KVRH}
Poropudas, J., Virtanen, K., 2010. {Game-Theoretic} validation and analysis of
  air combat simulation models. IEEE Transactions on Systems, Man, and
  Cybernetics - Part A: Systems and Humans 40~(5), 1057--1070.

\bibitem[{Poropudas and Virtanen(2011)}]{Poropudas2011-ns_KVRH}
Poropudas, J., Virtanen, K., 2011. Simulation metamodeling with dynamic
  {B}ayesian networks. European Journal of Operational Research 214~(3),
  644--655.

\bibitem[{Portela and Thanassoulis(2010)}]{Portela2010-zl_SL}
Portela, M. C. A.~S., Thanassoulis, E., 2010. Malmquist-type indices in the
  presence of negative data: An application to bank branches. Journal of
  Banking \& Finance 34~(7), 1472--1483.

\bibitem[{Portela et~al.(2012)Portela, Camanho, and Borges}]{Portela2012-bw_JJ}
Portela, M. C.~S., Camanho, A.~S., Borges, D., 2012. Performance assessment of
  secondary schools: the snapshot of a country taken by {DEA}. Journal of the
  Operational Research Society 63~(8), 1098--1115.

\bibitem[{Porteus(1990)}]{porteus1990stochastic_JSS}
Porteus, E.~L., 1990. Stochastic inventory theory. In: Heyman, D.~P., Sobel,
  M.~J. (Eds.), Stochastic Models. Vol.~2 of Handbooks in Operations Research
  and Management Science. Elsevier, pp. 605--652.

\bibitem[{Porteus(2002)}]{porteus2002foundations_JSS}
Porteus, E.~L., 2002. Foundations of stochastic inventory theory. Stanford
  University Press.

\bibitem[{Post et~al.(2012)Post, Ahmadi, Daskalaki, Kingston, Kyngas, Nurmi,
  and Ranson}]{post2012_GVBSP}
Post, G., Ahmadi, S., Daskalaki, S., Kingston, J.~H., Kyngas, J., Nurmi, C.,
  Ranson, D., 2012. An {XML} format for benchmarks in high school timetabling.
  Annals of Operations Research 194, 385--397.

\bibitem[{Post et~al.(2016)Post, {Di Gaspero}, Kingston, McCollum, and
  Schaerf}]{post2016_GVBSP}
Post, G., {Di Gaspero}, L., Kingston, J.~H., McCollum, B., Schaerf, A., 2016.
  The third international timetabling competition. Annals of Operations
  Research 239, 69--–75.

\bibitem[{Post et~al.(2014)Post, Kingston, Ahmadi, Daskalaki, Gogos, Kyngas,
  Nurmi, Musliu, Pillay, Santos, and Schaerf}]{post2014_GVBSP}
Post, G., Kingston, J.~H., Ahmadi, S., Daskalaki, S., Gogos, C., Kyngas, J.,
  Nurmi, C., Musliu, N., Pillay, N., Santos, H., Schaerf, A., 2014. {XHSTT}: an
  {XML} archive for high school timetabling problems in different countries.
  Annals of Operations Research 218, 295--301.

\bibitem[{Posta et~al.(2014)Posta, Ferland, and Michelon}]{posta2014exact_SAA}
Posta, M., Ferland, J.~A., Michelon, P., 2014. An exact cooperative method for
  the uncapacitated facility location problem. Mathematical Programming
  Computation 6~(3), 199--231.

\bibitem[{Pouwels et~al.(2021)Pouwels, House, Pritchard, Robotham, Birrell,
  Gelman, Vihta, Bowers, Boreham, Thomas, Lewis, Bell, Bell, Newton, Farrar,
  Diamond, Benton, Walker, and {COVID-19 Infection Survey
  Team}}]{Pouwels2021-xu}
Pouwels, K.~B., House, T., Pritchard, E., Robotham, J.~V., Birrell, P.~J.,
  Gelman, A., Vihta, K.-D., Bowers, N., Boreham, I., Thomas, H., Lewis, J.,
  Bell, I., Bell, J.~I., Newton, J.~N., Farrar, J., Diamond, I., Benton, P.,
  Walker, A.~S., {COVID-19 Infection Survey Team}, 2021. Community prevalence
  of {SARS-CoV-2} in {E}ngland from {April to November}, 2020: results from the
  {ONS} coronavirus infection survey. The Lancet. Public Health 6~(1),
  e30--e38.

\bibitem[{Powell(2011)}]{Powell2011-km_HL}
Powell, W.~B., 2011. Approximate Dynamic Programming: Solving the Curses of
  Dimensionality. Wiley.

\bibitem[{Powell and Sheffi(1983)}]{POWELL1983471_MH}
Powell, W.~B., Sheffi, Y., 1983. The load planning problem of motor carriers:
  Problem description and a proposed solution approach. Transportation Research
  Part A: General 17~(6), 471--480.

\bibitem[{Powell and Sheffi(1989)}]{powell/sheffi:1989_MH}
Powell, W.~B., Sheffi, Y., 1989. {Design and Implementation of an Interactive
  Optimization System for Network Design in the Motor Carrier Industry}.
  Operations Research 37, 12--29.

\bibitem[{Power et~al.(2018)Power, Heavin, McDermott, and
  Daly}]{Power2018-kx_JEB}
Power, D.~J., Heavin, C., McDermott, J., Daly, M., 2018. Defining business
  analytics: an empirical approach. Journal of Business Analytics 1~(1),
  40--53.

\bibitem[{P{\"o}yh{\"o}nen et~al.(2001)P{\"o}yh{\"o}nen, Vrolijk, and
  H{\"a}m{\"a}l{\"a}inen}]{Poyhonen2001-dw_AFRH}
P{\"o}yh{\"o}nen, M., Vrolijk, H., H{\"a}m{\"a}l{\"a}inen, R.~P., 2001.
  Behavioral and procedural consequences of structural variation in value
  trees. European Journal of Operational Research 134~(1), 216--227.

\bibitem[{Pr{\'e}kopa(2013)}]{Prekopa2013-du_HL}
Pr{\'e}kopa, A., 2013. Stochastic Programming. Springer Science \& Business
  Media.

\bibitem[{Prim(1957)}]{P57_SMPT}
Prim, R., 1957. Shortest connection networks and some generalizations. Bell
  System Technical Journal 36, 1389--1401.

\bibitem[{Pr{\"o}mel and Steger(2012)}]{PS12_SMPT}
Pr{\"o}mel, H., Steger, A., 2012. The Steiner tree problem: {A} tour through
  graphs, algorithms, and complexity. Springer, Berlin.

\bibitem[{Pruhs and Woeginger(2007)}]{PW07_UPCT}
Pruhs, K., Woeginger, G.~J., 2007. Approximation schemes for a class of subset
  selection problems. Theoretical Computer Science 382~(2), 151--156.

\bibitem[{Psaraftis and Kontovas(2010)}]{Psaraftis2010-rt_HP}
Psaraftis, H.~N., Kontovas, C.~A., 2010. Balancing the economic and
  environmental performance of maritime transportation. Transportation Research
  Part D: Transport and Environment 15~(8), 458--462.

\bibitem[{Psaraftis and Kontovas(2013)}]{Psaraftis2013-zv_HP}
Psaraftis, H.~N., Kontovas, C.~A., 2013. Speed models for energy-efficient
  maritime transportation: A taxonomy and survey. Transportation Research Part
  C: Emerging Technologies 26, 331--351.

\bibitem[{Psaraftis and Kontovas(2014)}]{Psaraftis2014-lu_HP}
Psaraftis, H.~N., Kontovas, C.~A., 2014. Ship speed optimization: Concepts,
  models and combined speed-routing scenarios. Transportation Research Part C:
  Emerging Technologies 44, 52--69.

\bibitem[{Pulleyblank et~al.(2000)Pulleyblank, Dietrich, Forrest, and
  Lougee-Heimer}]{Coinor_CTCGE}
Pulleyblank, W., Dietrich, B., Forrest, J., Lougee-Heimer, R., 2000. Open
  source for optimization software. In: International Symposium for
  Mathematical Programming {(ISMP)}.

\bibitem[{Puterman(2014)}]{Puterman2005-lr_HL}
Puterman, M.~L., 2014. Markov Decision Processes: Discrete Stochastic Dynamic
  Programming. Wiley.

\bibitem[{Pyrgiotis et~al.(2013)Pyrgiotis, Malone, and
  Odoni}]{pyrgiotis2013modelling_VLVV}
Pyrgiotis, N., Malone, K.~M., Odoni, A., 2013. Modelling delay propagation
  within an airport network. Transportation Research Part C: Emerging
  Technologies 27, 60--75.

\bibitem[{Qi(2015)}]{Qi2015-rp_HP}
Qi, X., 2015. Disruption management for liner shipping. In: Lee, C.-Y., Meng,
  Q. (Eds.), Handbook of Ocean Container Transport Logistics: Making Global
  Supply Chains Effective. Springer, Cham, pp. 231--249.

\bibitem[{Qi et~al.(2022)Qi, Harrod, Psaraftis, and Lang}]{Qi2022-le_HP}
Qi, Y., Harrod, S., Psaraftis, H.~N., Lang, M., 2022. Transport service
  selection and routing with carbon emissions and inventory costs consideration
  in the context of the belt and road initiative. Transportation Research Part
  E: Logistics and Transportation Review 159, 102630.

\bibitem[{Qin et~al.(2011)Qin, Wang, Vakharia, Chen, and
  Seref}]{qin2011newsvendor_JSS}
Qin, Y., Wang, R., Vakharia, A.~J., Chen, Y., Seref, M.~M., 2011. The
  newsvendor problem: Review and directions for future research. European
  Journal of Operational Research 213~(2), 361--374.

\bibitem[{Qu et~al.(2015)Qu, Pham, Bai, and Kendall}]{Qu2015-oi_JJ}
Qu, R., Pham, N., Bai, R., Kendall, G., 2015. Hybridising heuristics within an
  estimation distribution algorithm for examination timetabling. Applied
  Intelligence 42~(4), 679--693.

\bibitem[{Queiroga et~al.(2021)Queiroga, Sadykov, and
  Uchoa}]{QUEIROGA2021105475_CA_MB}
Queiroga, E., Sadykov, R., Uchoa, E., 2021. A popmusic matheuristic for the
  capacitated vehicle routing problem. Computers \& Operations Research 136,
  105475.

\bibitem[{Queiroga et~al.(2022)Queiroga, Sadykov, Uchoa, and
  Vidal}]{queiroga2022_CA_MB}
Queiroga, E., Sadykov, R., Uchoa, E., Vidal, T., 2022. 10,000 optimal {CVRP}
  solutions for testing machine learning based heuristics. In: AAAI-22 Workshop
  on Machine Learning for Operations Research (ML4OR). pp. 1--2.

\bibitem[{Raack et~al.(2011)Raack, Koster, Orlowski, and
  Wess{\"a}ly}]{raack2011cut_MH}
Raack, C., Koster, A.~M., Orlowski, S., Wess{\"a}ly, R., 2011. On cut-based
  inequalities for capacitated network design polyhedra. Networks 57~(2),
  141--156.

\bibitem[{Rad and Roy(2021)}]{Rad2021-yr_TC}
Rad, S.~R., Roy, O., 2021. Deliberation, {Single-Peakedness}, and coherent
  aggregation. The American Political Science Review 115~(2), 629--648.

\bibitem[{Rahmaniani et~al.(2017)Rahmaniani, Crainic, Gendreau, and
  Rei}]{rahmaniani2017benders_COIT}
Rahmaniani, R., Crainic, T.~G., Gendreau, M., Rei, W., 2017. The {B}enders
  decomposition algorithm: {A} literature review. European Journal of
  Operational Research 259~(3), 801--817.

\bibitem[{Rajagopalan(2020)}]{Rajagopalan2021-cv_GM}
Rajagopalan, R., 2020. Immersive Systemic Knowing: Advancing Systems Thinking
  Beyond Rational Analysis. Springer.

\bibitem[{Ralphs et~al.(2022)Ralphs, Bulut, and Vigerske}]{Gimpy_CTCGE}
Ralphs, T., Bulut, A., Vigerske, S., 2022. {GiMPy}.
\newline\urlprefix\url{https://projects.coin-or.org/Gimpy}

\bibitem[{Ralphs and Guzelsoy(2005)}]{Symphony_CTCGE}
Ralphs, T.~K., Guzelsoy, M., 2005. The {Symphony} callable library for mixed
  integer programming. In: Golden, B., Raghavan, S., Wasil, E. (Eds.), The Next
  Wave in Computing, Optimization, and Decision Technologies. Springer US,
  Boston, MA, pp. 61--76.

\bibitem[{Ram Mohan~Rao et~al.(2018)Ram Mohan~Rao, Murali~Krishna, and
  Siva~Kumar}]{Ram_Mohan_Rao2018-jb_JEB}
Ram Mohan~Rao, P., Murali~Krishna, S., Siva~Kumar, A.~P., 2018. Privacy
  preservation techniques in big data analytics: a survey. Journal of Big Data
  5~(1), 1--12.

\bibitem[{Ramsey(1931)}]{RePEc:hay:hetcha:ramsey1926_MESG}
Ramsey, F.~P., 1931. Truth and probability. In: Braithwaite, R.~B. (Ed.), The
  Foundations of Mathematics and other Logical Essays. Routledge and Kegan
  Paul, London, England, Ch.~7, pp. 156--198.

\bibitem[{Rancourt et~al.(2015)Rancourt, Cordeau, Laporte, and
  Watkins}]{foodaid2015_BYKOK}
Rancourt, M., Cordeau, J., Laporte, G., Watkins, B., 2015. Tactical network
  planning for food aid distribution in kenya. Computers \& Operations Research
  56, 68--83.

\bibitem[{Ranyard et~al.(2021)Ranyard, Hopes, and Murray}]{Ranyard2021-kw_MYLW}
Ranyard, J., Hopes, J., Murray, E., 2021. Reinvigorating soft {OR} for
  practitioners: Report to {HORAF}. 63rd Conference of the UK Operational
  Research Society (OR63).

\bibitem[{Ranyard et~al.(2015)Ranyard, Fildes, and Hu}]{Ranyard2015-xx_MYLW}
Ranyard, J.~C., Fildes, R., Hu, T.-I., 2015. Reassessing the scope of {OR}
  practice: The influences of problem structuring methods and the analytics
  movement. European Journal of Operational Research 245~(1), 1--13.

\bibitem[{Rasmussen and Williams(2005)}]{Rasmussen05_LCAL}
Rasmussen, C.~E., Williams, C. K.~I., 2005. {Gaussian Processes for Machine
  Learning}. The MIT Press, Cambridge MA.

\bibitem[{Rasmussen and Trick(2008)}]{Rasmussen2008_GVBSP}
Rasmussen, R.~V., Trick, M.~A., 2008. Round robin scheduling -- a survey.
  European Journal of Operational Research 188~(3), 617--636.

\bibitem[{Rawls and Turnquist(2010)}]{rawls2010pre_BYKOK}
Rawls, C.~G., Turnquist, M.~A., 2010. Pre-positioning of emergency supplies for
  disaster response. Transportation Research Part B: Methodological 44~(4),
  521--534.

\bibitem[{Rawls and Turnquist(2011)}]{Rawls2011_JLYHK}
Rawls, C.~G., Turnquist, M.~A., 2011. Pre-positioning planning for emergency
  response with service quality constraints. OR Spectrum 33~(3), 481--498.

\bibitem[{Rawls(1971)}]{Raw99_JH}
Rawls, J., 1971. A Theory of Justice. Harvard University Press.

\bibitem[{Rea et~al.(2021)Rea, Froehle, Masterson, Stettler, Fermann, and
  Pancioli}]{ReaFroMasSteFerPan21_JH}
Rea, D., Froehle, C., Masterson, S., Stettler, B., Fermann, G., Pancioli, A.,
  2021. Unequal but fair: {Incorporating} distributive justice in operational
  allocation models. Production and Operations Management 30, 2304--2320.

\bibitem[{{Red Hat}(2022)}]{Optaplanner_CTCGE}
{Red Hat}, 2022. Optaplanner.
\newline\urlprefix\url{https://www.optaplanner.org/}

\bibitem[{Reddie et~al.(2018)Reddie, Goldblum, Lakkaraju, Reinhardt, Nacht, and
  Epifanovskaya}]{Reddie2018-zs_KVRH}
Reddie, A.~W., Goldblum, B.~L., Lakkaraju, K., Reinhardt, J., Nacht, M.,
  Epifanovskaya, L., 2018. Next-generation wargames. Science 362~(6421),
  1362--1364.

\bibitem[{Reed(2009)}]{reed2009g_HATI}
Reed, J., 2009. The {G/GI/N} queue in the halfin--whitt regime. The Annals of
  Applied Probability 19~(6), 2211--2269.

\bibitem[{Reed et~al.(2022)Reed, Campbell, and Thomas}]{reed2022impact_JLYHK}
Reed, S., Campbell, A.~M., Thomas, B.~W., 2022. Impact of autonomous vehicle
  assisted last-mile delivery in urban to rural settings. Transportation
  Science 56~(6), 1530--1548.

\bibitem[{Rehfeldt(2021)}]{Rehfeldt:2021_IL}
Rehfeldt, D., 2021. Faster algorithms for steiner tree and related problems:
  From theory to practice. Ph.D. thesis, Technische Universit\"{a}t Berlin.

\bibitem[{Rehfeldt and Koch(2021)}]{Rehfeldt-Koch:2021_IL}
Rehfeldt, D., Koch, T., 2021. Implications, conflicts, and reductions for
  {S}teiner trees. In: Singh, M., Williamson, D.~P. (Eds.), Integer Programming
  and Combinatorial Optimization - 22nd International Conference, {IPCO} 2021,
  Atlanta, GA, USA, May 19-21, 2021, Proceedings. Vol. 12707 of Lecture Notes
  in Computer Science. Springer, pp. 473--487.

\bibitem[{Reichert et~al.(2015)Reichert, Langhans, Lienert, and
  Schuwirth}]{Reichert2015-jd_JL}
Reichert, P., Langhans, S.~D., Lienert, J., Schuwirth, N., 2015. The conceptual
  foundation of environmental decision support. Journal of Environmental
  Management 154, 316--332.

\bibitem[{Reinhardt and Pisinger(2012)}]{Reinhardt2012-do_HP}
Reinhardt, L.~B., Pisinger, D., 2012. A branch and cut algorithm for the
  container shipping network design problem. Flexible Services and
  Manufacturing Journal 24~(3), 349--374.

\bibitem[{Ren and Huang(2022)}]{Ren2022_SMD}
Ren, H., Huang, T., 2022. Opaque selling and inventory management in vertically
  differentiated markets. Manufacturing \& Service Operations Management
  24~(5), 2543--2557.

\bibitem[{Ren et~al.(2021)Ren, Jiang, Khoveyni, Guan, and Yang}]{Ren2021-lh_SL}
Ren, X., Jiang, C., Khoveyni, M., Guan, Z., Yang, G., 2021. A review of {DEA}
  methods to identify and measure congestion. Journal of Management Science and
  Engineering 6~(4), 345--362.

\bibitem[{Resnick(2007)}]{resnick2007heavy_HATI}
Resnick, S.~I., 2007. Heavy-Tail Phenomena: Probabilistic and Statistical
  Modeling. Springer Science \& Business Media.

\bibitem[{Ribeiro(2012)}]{ribeiro_GVBSP}
Ribeiro, C.~C., 2012. Sports scheduling: Problems and applications.
  International Transactions in Operational Research 19~(2), 201--226.

\bibitem[{Richardson and Weithman(1999)}]{RicWei99_JH}
Richardson, H.~S., Weithman, P.~J. (Eds.), 1999. The Philosophy of Rawls {\rm
  (5 volumes)}. Garland.

\bibitem[{Righini and Salani(2008)}]{RighiniS2008_CA_MB}
Righini, G., Salani, M., 2008. New dynamic programming algorithms for the
  resource constrained elementary shortest path problem. Networks 51~(3),
  155--170.

\bibitem[{Rintam{\"a}ki et~al.(2021)Rintam{\"a}ki, Spence, Saarij{\"a}rvi,
  Joensuu, and Yrj{\"o}l{\"a}}]{rintamaki2021customers_TVWCK}
Rintam{\"a}ki, T., Spence, M.~T., Saarij{\"a}rvi, H., Joensuu, J.,
  Yrj{\"o}l{\"a}, M., 2021. Customers' perceptions of returning items purchased
  online: planned versus unplanned product returners. International Journal of
  Physical Distribution \& Logistics Management 51~(4), 403--422.

\bibitem[{Rios~Insua et~al.(2021)Rios~Insua, Couce-Vieira, Rubio, Pieters,
  Labunets, and G~Rasines}]{Rios_Insua2021-yx_KVRH}
Rios~Insua, D., Couce-Vieira, A., Rubio, J.~A., Pieters, W., Labunets, K.,
  G~Rasines, D., 2021. An adversarial risk analysis framework for
  cybersecurity. Risk Analysis 41~(1), 16--36.

\bibitem[{Ritchie et~al.(1994)Ritchie, Taket, and Bryant}]{Ritchie1994-rw_AJG}
Ritchie, C., Taket, A., Bryant, J. (Eds.), 1994. Community Works: 26 Case
  studies Showing Community Operational Research in Action. Pavic Press,
  Sheffield.

\bibitem[{Rittel and Webber(1973)}]{Rittel1973-wn_JL}
Rittel, H. W.~J., Webber, M.~M., 1973. Dilemmas in a general theory of
  planning. Policy Sciences 4~(2), 155--169.

\bibitem[{Robenek et~al.(2017)Robenek, Azadeh, Maknoona, and
  Bierlaire}]{Robenek2017_DC}
Robenek, T., Azadeh, S.~S., Maknoona, Y., Bierlaire, M., 2017. Hybrid
  cyclicity{: C}ombining the benefits of cyclic and non-cyclic timetables.
  Transportation Research Part {C: E}merging {T}echnologies 75, 228--253.

\bibitem[{Roberti and Ruthmair(2021)}]{roberti2021exact_VLVV}
Roberti, R., Ruthmair, M., 2021. Exact methods for the traveling salesman
  problem with drone. Transportation Science 55~(2), 315--335.

\bibitem[{Robinson and Robinson(1994)}]{Robinson1994_SMD}
Robinson, A., Robinson, M., 1994. On the tabletop improvement experiments of
  {J}apan. Production and Operations Management 3~(3), 201--217.

\bibitem[{Robinson(2014)}]{Robinson2014-bd_AFRH}
Robinson, S., 2014. Simulation: The Practice of Model Development and Use.
  Palgrave Macmillan, London.

\bibitem[{Robinson(2020)}]{Robinson2020_CC}
Robinson, S., 2020. Conceptual modelling for simulation: Progress and grand
  challenges. Journal of Simulation 14~(1), 1--20.

\bibitem[{Robinson(2008)}]{Robinson2008_CC}
Robinson, S.~L., 2008. Conceptual modelling for simulation {Part I}: definition
  and requirements. Journal of the Operational Research Society 59~(3),
  278--290.

\bibitem[{Rocha et~al.(2013)Rocha, Oliveira, and Carravilla}]{cyclic_GVBSP}
Rocha, M., Oliveira, J.~F., Carravilla, M.~A., 2013. Cyclic staff scheduling:
  optimization models for some real-life problems. Journal of Scheduling
  16~(2), 231--242.

\bibitem[{Rockafellar and Uryasev(2002)}]{rockafellar02_MB}
Rockafellar, R.~T., Uryasev, S., 2002. Conditional value-at-risk for general
  loss distributions. Journal of Banking \& Finance 26~(7), 1443--1471.

\bibitem[{Rockafellar et~al.(2000)Rockafellar, Uryasev,
  et~al.}]{rockafellar00_MB}
Rockafellar, R.~T., Uryasev, S., et~al., 2000. Optimization of conditional
  value-at-risk. Journal of Risk 2, 21--42.

\bibitem[{Rockafellar and Wets(1991)}]{RockWetsPH_MH}
Rockafellar, R.~T., Wets, R. J.-B., 1991. Scenarios and policy aggregation in
  optimization under uncertainty. Mathematics of Operations Research 16~(1),
  119--147.

\bibitem[{{Rockwell Automation}(2022)}]{Arena_CTCGE}
{Rockwell Automation}, 2022. {Arena Simulation}.
\newline\urlprefix\url{https://www.rockwellautomation.com/en-us/products/software/arena-simulation.html}

\bibitem[{Rogers and Tibben-Lembke(2001)}]{Rogers2001_JLYHK}
Rogers, D.~S., Tibben-Lembke, R., 2001. An examination of reverse logistics
  practices. Journal of Business Logistics 22~(2), 129--148.

\bibitem[{Romero and Abad(2022)}]{Romero2022-nw_JEB}
Romero, J.~A., Abad, C., 2022. Cloud-based big data analytics integration with
  {ERP} platforms. Management Decision 60~(12), 3416--3437.

\bibitem[{Rose(2016)}]{Rose2016-lq_JEB}
Rose, B., 2016. Defining analytics: a conceptual framework. ORMS Today 43~(3),
  36--40.

\bibitem[{Rosenbaum(2004)}]{Rosenbaum2004-ut_IM}
Rosenbaum, D.~T., 2004. Measuring how {NBA} players help their teams win.
  \url{http://www.82games.com/comm30.htm}, accessed on 2022-11-1.

\bibitem[{Rosenberger et~al.(2003)Rosenberger, Johnson, and
  Nemhauser}]{rosenberger2003rerouting_VLVV}
Rosenberger, J.~M., Johnson, E.~L., Nemhauser, G.~L., 2003. Rerouting aircraft
  for airline recovery. Transportation Science 37~(4), 408--421.

\bibitem[{Rosenhead(1986)}]{Rosenhead1986-ti_MYLW}
Rosenhead, J., 1986. Custom and practice. Journal of the Operational Research
  Society 37~(4), 335--343.

\bibitem[{Rosenhead(1989)}]{Rosenhead1989-nu_MYLW}
Rosenhead, J. (Ed.), 1989. Rational Analysis for a Problematic World: Problem
  Structuring Methods for Complexity, Uncertainty and Conflict. Wiley,
  Chichester.

\bibitem[{Rosenhead(1993)}]{Rosenhead1993-dj_AJG}
Rosenhead, J., 1993. Enabling analysis: Across the developmental divide.
  Systems Practice 6~(2), 117--138.

\bibitem[{Rosenhead(1996)}]{Rosenhead1996-xa_MYLW}
Rosenhead, J., 1996. What's the problem? {A}n introduction to problem
  structuring methods. Interfaces 26~(6), 117--131.

\bibitem[{Rosenhead(2006)}]{Rosenhead2006-ih_MYLW}
Rosenhead, J., 2006. Past, present and future of problem structuring methods.
  Journal of the Operational Research Society 57~(7), 759--765.

\bibitem[{Rosenhead(2013)}]{Rosenhead2013-zp_AJG}
Rosenhead, J., 2013. Book review of {Johnson, Michael P., ed. 2012}.
  {Community-Based} operations research: Decision modeling for local impact and
  diverse populations. Interfaces 43, 609--610.

\bibitem[{Rosenhead and Mingers(2001)}]{Rosenhead2001-ai_JL}
Rosenhead, J., Mingers, J. (Eds.), 2001. Rational Analysis for a Problematic
  World Revisited: Problem Structuring Methods for Complexity, Uncertainty and
  Conflict, 2nd Edition. Wiley.

\bibitem[{Rosenhead and White(1996)}]{Rosenhead1996-zi_AJG}
Rosenhead, J., White, L., 1996. Nuclear fusion: Some linked case studies in
  community operational research. Journal of the Operational Research Society
  47~(4), 479--489.

\bibitem[{Rosling(1989)}]{rosling1989optimal_JSS}
Rosling, K., 1989. Optimal inventory policies for assembly systems under random
  demands. Operations Research 37~(4), 565--579.

\bibitem[{Ross(1983)}]{Ross1983-yc_HL}
Ross, S.~M., 1983. Introduction to Stochastic Dynamic Programming. Academic
  Press.

\bibitem[{Roth et~al.(2016)Roth, Singhal, Singhal, and
  Tang}]{Roth2016-pg_KESXZ}
Roth, A., Singhal, J., Singhal, K., Tang, C.~S., 2016. Knowledge creation and
  dissemination in operations and supply chain management. International
  Journal of Operations \& Production Management 25~(9), 1473--1488.

\bibitem[{Rother and Shook(1999)}]{Rother1999_SMD}
Rother, M., Shook, J., 1999. Learning to See: Value Stream Mapping to Add Value
  and Eliminate Muda. The Lean Enterprise Institute, Cambridge, MA.

\bibitem[{Rothstein(1971)}]{rothstein1971airline_VLVV}
Rothstein, M., 1971. An airline overbooking model. Transportation Science
  5~(2), 180--192.

\bibitem[{Roundy(1985)}]{roundy198598_JSS}
Roundy, R., 1985. 98\%-effective integer-ratio lot-sizing for one-warehouse
  multi-retailer systems. Management Science 31~(11), 1416--1430.

\bibitem[{Roundy(1986)}]{roundy198698_JSS}
Roundy, R., 1986. A 98\%-effective lot-sizing rule for a multi-product,
  multi-stage production/inventory system. Mathematics of Operations Research
  11~(4), 699--727.

\bibitem[{Rouwette(2016)}]{Rouwette2016-su_MCJM}
Rouwette, E. A. J.~A., 2016. The impact of group model building on behavior.
  In: Kunc, M., Malpass, J., White, L. (Eds.), Behavioral Operational Research:
  Theory, Methodology and Practice. Palgrave Macmillan, London, pp. 213--241.

\bibitem[{Rouwette et~al.(2010)Rouwette, Korzilius, Vennix, and
  Jacobs}]{Rouwette2010-sc_MCJM}
Rouwette, E. A. J.~A., Korzilius, H., Vennix, J. A.~M., Jacobs, E., 2010.
  Modeling as persuasion: the impact of group model building on attitudes and
  behavior. System Dynamics Review 27~(1), 1--21.

\bibitem[{Rowe and Wright(1999)}]{Rowe1999-vs_FP}
Rowe, G., Wright, G., 1999. The delphi technique as a forecasting tool: issues
  and analysis. International Journal of Forecasting 15~(4), 353--375.

\bibitem[{Roy(1993)}]{roy1993decision_MESG}
Roy, B., 1993. Decision science or decision-aid science? European Journal of
  Operational Research 66~(2), 184--203.

\bibitem[{Roy(1996)}]{roy1996multicriteria_MESG}
Roy, B., 1996. Multicriteria Methodology for Decision Aiding. Vol.~12 of
  Nonconvex Optimization and its Applications. Springer-Verlag, Berlin.

\bibitem[{Roy(2005)}]{roy2005paradigms_MESG}
Roy, B., 2005. Paradigms and challenges. In: Greco, S., Ehrgott, M., Figueira,
  J.~R. (Eds.), Multiple Criteria Decision Analysis: State of the Art Surveys.
  Springer-Verlag, Berlin, pp. 3--24.

\bibitem[{Roy and S{\l}owi{\'n}ski(2013)}]{Roy2013-ii_JL}
Roy, B., S{\l}owi{\'n}ski, R., 2013. Questions guiding the choice of a
  multicriteria decision aiding method. EURO Journal on Decision Processes
  1~(1), 69--97.

\bibitem[{Royset et~al.(2009)Royset, Carlyle, and Wood}]{Royset2009-ps_KVRH}
Royset, J.~O., Carlyle, W.~M., Wood, R.~K., 2009. Routing military aircraft
  with a constrained {Shortest-Path} algorithm. Military Operations Research
  14~(3), 31--52.

\bibitem[{Royston(2009)}]{Royston2009-ae_CV}
Royston, G., 2009. One hundred years of operational research in {Health---UK}
  1948--2048. Journal of the Operational Research Society 60~(Supplement 1),
  s169--s179.

\bibitem[{RSAS(2020)}]{RSAS20_BC}
RSAS, 2020. The prize in economic sciences 2020. Press release, October 12,
  Royal Swedish Academy of Sciences.

\bibitem[{Rubinstein(1982)}]{Rub82_JH}
Rubinstein, A., 1982. Perfect equilibrium in a bargaining model. Econometrica
  50, 97--109.

\bibitem[{Ruf and Cordeau(2021)}]{Ruf2021_DC}
Ruf, M., Cordeau, J.~F., 2021. Adaptive large neighborhood search for
  integrated planning in railroad classification yards. Transportation Research
  Part B: Methodological 150, 26--51.

\bibitem[{Ruibal and Mazumdar(2008)}]{rui:maz:08_DSRW}
Ruibal, C., Mazumdar, M., 2008. Forecasting the mean and the variance of
  electricity prices in deregulated markets. IEEE Transactions on Power Systems
  23~(1), 25--32.

\bibitem[{Ruivo et~al.(2020)Ruivo, Johansson, Sarker, and
  Oliveira}]{Ruivo2020-kq_JEB}
Ruivo, P., Johansson, B., Sarker, S., Oliveira, T., 2020. The relationship
  between {ERP} capabilities, use, and value. Computers in Industry 117,
  103209.

\bibitem[{Ruiz-Tagle et~al.(2022)Ruiz-Tagle, Lopez~Droguett, and
  Groth}]{Ruiz-Tagle2022-ts_TC}
Ruiz-Tagle, A., Lopez~Droguett, E., Groth, K.~M., 2022. Exploiting the
  capabilities of {B}ayesian networks for engineering risk assessment: Causal
  reasoning through interventions. Risk Analysis 42~(6), 1306--1324.

\bibitem[{Rumelhart et~al.(1986)Rumelhart, Hinton, and
  Williams}]{Rumelhart86_LCAL}
Rumelhart, D.~E., Hinton, G.~E., Williams, R.~J., 1986. Learning
  representations by back-propagating errors. Nature 323~(6088), 533--536.

\bibitem[{Russell et~al.(2017)Russell, Kusner, Loftus, and
  Silva}]{RusKusLofSil17_JH}
Russell, C., Kusner, M.~J., Loftus, J.~R., Silva, R., 2017. When worlds
  collide: {Integrating} different counterfactual assumptions in fairness. In:
  Proceedings of 31st International Conference on Neural Information Processing
  Systems. pp. 6417--6426.

\bibitem[{Russell(1998)}]{Russell98_LCAL}
Russell, S., 1998. Learning agents for uncertain environments (extended
  abstract). In: Bartlett, P.~L., Mansour, Y. (Eds.), Proceedings of the
  Eleventh Annual Conference on Computational Learning Theory, {COLT} 1998,
  Madison, Wisconsin, USA, July 24-26, 1998. {ACM}, pp. 101--103.

\bibitem[{Saaty(1977)}]{saaty1977scaling_MESG}
Saaty, T., 1977. A scaling method for priorities in hierarchical structures.
  Journal of Mathematical Psychology 15~(3), 234--281.

\bibitem[{Sadykov(2022)}]{Sadykov2022_CA_MB}
Sadykov, R., 2022. Academic webpage.
  \url{https://www.math.u-bordeaux.fr/~rsadykov/}, accessed on 2022-09-16.

\bibitem[{Sadykov et~al.(2021)Sadykov, Uchoa, and Pessoa}]{Sadykov2021_CA_MB}
Sadykov, R., Uchoa, E., Pessoa, A., 2021. A bucket graph–based labeling
  algorithm with application to vehicle routing. Transportation Science 55~(1),
  4--28.

\bibitem[{Sadykov et~al.(2019)Sadykov, Vanderbeck, Pessoa, Tahiri, and
  Uchoa}]{Sadykov2019_CA_MB}
Sadykov, R., Vanderbeck, F., Pessoa, A., Tahiri, I., Uchoa, E., 2019. Primal
  heuristics for branch and price: The assets of diving methods. {INFORMS}
  Journal on Computing 31~(2), 251--267.

\bibitem[{{\c{S}}ahin and S{\"u}ral(2007)}]{csahin2007review_SAA}
{\c{S}}ahin, G., S{\"u}ral, H., 2007. A review of hierarchical facility
  location models. Computers \& Operations Research 34~(8), 2310--2331.

\bibitem[{Sahinidis(1996)}]{Baron_CTCGE}
Sahinidis, N., 1996. {BARON}: A general purpose global optimization software
  package. Journal of Global Optimization 8, 201--205.

\bibitem[{Salahi et~al.(2021)Salahi, Toloo, and Torabi}]{Salahi2021-eh_SL}
Salahi, M., Toloo, M., Torabi, N., 2021. A new robust optimization approach to
  common weights formulation in {DEA}. Journal of the Operational Research
  Society 72~(6), 1390--1402.

\bibitem[{Salo et~al.(2011)Salo, Keisler, and Morton}]{salo2011portfolio_MESG}
Salo, A., Keisler, J., Morton, A., 2011. Portfolio Decision Analysis: Improved
  Methods for Resource Allocation. Vol. 162 of International Series in
  Operations Research and Management Science. Springer-Verlag, Berlin.

\bibitem[{Salo and Punkka(2011)}]{Salo2011-io_DEA}
Salo, A., Punkka, A., 2011. Ranking intervals and dominance relations for
  {Ratio-Based} efficiency analysis. Management Science 57~(1), 200--214.

\bibitem[{Salo et~al.(2022)Salo, Tosoni, Roponen, and Bunn}]{Salo2022-yz_JL}
Salo, A., Tosoni, E., Roponen, J., Bunn, D.~W., 2022. Using cross‐impact
  analysis for probabilistic risk assessment. Futures \& Foresight Science
  4~(2), e2103.

\bibitem[{Samorani et~al.(2019)Samorani, Alptekino{\u g}lu, and
  Messinger}]{Samorani2019-rw_AM}
Samorani, M., Alptekino{\u g}lu, A., Messinger, P.~R., 2019. Product return
  episodes in retailing. Service Science 11~(4), 263--278.

\bibitem[{Sampaio et~al.(2020)Sampaio, Savelsbergh, Veelenturf, and
  Van~Woensel}]{sampaio2020delivery_JLYHK}
Sampaio, A., Savelsbergh, M., Veelenturf, L.~P., Van~Woensel, T., 2020.
  Delivery systems with crowd-sourced drivers: A pickup and delivery problem
  with transfers. Networks 76~(2), 232--255.

\bibitem[{Samuel(1959)}]{Samuel_IBM_5392560_LCAL}
Samuel, A.~L., 1959. Some studies in machine learning using the game of
  checkers. IBM Journal of Research and Development 3~(3), 210--229.

\bibitem[{Samuelson(2014)}]{Samuelson14_BC}
Samuelson, W., 2014. Auctions: Advances in theory and practice. In: Chatterjee,
  K., Samuelson, W. (Eds.), Game Theory and Business Applications. Springer,
  Boston, MA, pp. 323--366.

\bibitem[{Sandholm(2010)}]{Sa2010_GZ}
Sandholm, W., 2010. Population games and evolutionary dynamics. {MIT} Press.

\bibitem[{Santibanez~Gonzalez et~al.(2019)Santibanez~Gonzalez, Koh, and
  Leung}]{SG2019_JLYHK}
Santibanez~Gonzalez, E.~D., Koh, L., Leung, J., 2019. Towards a circular
  economy production system: trends and challenges for operations management.
  International Journal of Production Research 57~(23), 7209--7218.

\bibitem[{Sarimveis et~al.(2008)Sarimveis, Patrinos, Tarantilis, and
  Kiranoudis}]{Sarimveis2008_XW}
Sarimveis, H., Patrinos, P., Tarantilis, C.~D., Kiranoudis, C.~T., 2008.
  Dynamic modeling and control of supply chain systems: A review. Computers \&
  Operations Research 35~(11), 3530--3561.

\bibitem[{Sartori et~al.(2021)Sartori, Gandra, {\c{C}}al{\i}k, and
  Smet}]{altachem_GVBSP}
Sartori, C.~S., Gandra, V., {\c{C}}al{\i}k, H., Smet, P., 2021. Production
  scheduling with stock- and staff-related restrictions. In: Mes, M.,
  Lalla-Ruiz, E., Vo{\ss}, S. (Eds.), Computational Logistics. Springer, Cham,
  pp. 142--162.

\bibitem[{Savage(1954)}]{savage1954foundations_MESG}
Savage, L., 1954. The Foundations of Statistics. John Wiley \& Sons, Hoboken,
  NJ.

\bibitem[{Savaskan and Van~Wassenhove(2006)}]{Savaskan2006-ou_AM}
Savaskan, R.~C., Van~Wassenhove, L.~N., 2006. Reverse channel design: The case
  of competing retailers. Management Science 52~(1), 1--14.

\bibitem[{Savelsbergh(1994)}]{Sa94_ALAL}
Savelsbergh, M., 1994. Preprocessing and probing techniques for mixed integer
  programming problems. ORSA Journal on Computing 6, 445--454.

\bibitem[{Savelsbergh and Van~Woensel(2016)}]{savelsbergh201650th_CKTVW}
Savelsbergh, M., Van~Woensel, T., 2016. 50th anniversary invited article—city
  logistics: Challenges and opportunities. Transportation Science 50~(2),
  579--590.

\bibitem[{Sawaragi et~al.(1985)Sawaragi, Nakayama, and Tanino}]{sawnata85_MESG}
Sawaragi, Y., Nakayama, H., Tanino, T., 1985. Theory of Multiobjective
  Optimization. Academic Press, Orlando.

\bibitem[{Saxena et~al.(2020)Saxena, Huang, DeFilippis, Radanovic, Parkes, and
  Liu}]{SaxHuaDefRadParLiu20_JH}
Saxena, N.~A., Huang, K., DeFilippis, E., Radanovic, G., Parkes, D.~C., Liu,
  Y., 2020. How do fairness definitions fare? {Testing} public attitudes
  towards three algorithmic definitions of fairness in loan allocations.
  Artificial Intelligence 283, 103238.

\bibitem[{Scala and {Howard II}(2020)}]{Scala2020_JLYHK}
Scala, N.~M., {Howard II}, J.~P. (Eds.), 2020. Handbook of Military and Defense
  Operations Research. CRC Press, Boca Raton, FL.

\bibitem[{Scarf(1959)}]{scarf1959bayes_JSS}
Scarf, H., 1959. Bayes solutions of the statistical inventory problem. The
  Annals of Mathematical Statistics 30~(2), 490--508.

\bibitem[{Scarf(1960{\natexlab{a}})}]{scarf1960optimality_JSS}
Scarf, H., 1960{\natexlab{a}}. The optimality of ({S}, s) policies in the
  dynamic inventory problem. In: Mathematical Methods in the Social Sciences.
  Stanford University Press, Stanford, pp. 196--202.

\bibitem[{Scarf(1960{\natexlab{b}})}]{scarf1960some_JSS}
Scarf, H.~E., 1960{\natexlab{b}}. Some remarks on {B}ayes solutions to the
  inventory problem. Naval Research Logistics Quarterly 7~(4), 591--596.

\bibitem[{Scarf et~al.(2019)Scarf, Parma, and McHale}]{Scarf2019-mk_IM}
Scarf, P., Parma, R., McHale, I., 2019. On outcome uncertainty and scoring
  rates in sport: The case of international rugby union. European Journal of
  Operational Research 273~(2), 721--730.

\bibitem[{Scarf et~al.(2009)Scarf, Yusof, and Bilbao}]{Scarf2009-lc_IM}
Scarf, P., Yusof, M.~M., Bilbao, M., 2009. A numerical study of designs for
  sporting contests. European Journal of Operational Research 198~(1),
  190--198.

\bibitem[{Scarselli et~al.(2009)Scarselli, Gori, Tsoi, Hagenbuchner, and
  Monfardini}]{Scarselli09_LCAL}
Scarselli, F., Gori, M., Tsoi, A.~C., Hagenbuchner, M., Monfardini, G., 2009.
  The graph neural network model. IEEE Transactions on Neural Networks 20~(1),
  61--80.

\bibitem[{Schaefer et~al.(2005)Schaefer, Johnson, Kleywegt, and
  Nemhauser}]{schaefer2005airline_VLVV}
Schaefer, A.~J., Johnson, E.~L., Kleywegt, A.~J., Nemhauser, G.~L., 2005.
  Airline crew scheduling under uncertainty. Transportation Science 39~(3),
  340--348.

\bibitem[{Scheithauer(2017)}]{Scheithauer2018-xm_JB}
Scheithauer, G., 2017. Introduction to Cutting and Packing Optimization:
  Problems, Modeling Approaches, Solution Methods. International Series in
  Operations Research \& Management Science. Springer.
\newline\urlprefix\url{https://books.google.co.uk/books?id=FlM7DwAAQBAJ}

\bibitem[{Scherfke(2021)}]{Simpy_CTCGE}
Scherfke, S., 2021. {SimPy}.
\newline\urlprefix\url{https://github.com/simpx/simpy}

\bibitem[{Schiffer et~al.(2022)Schiffer, Boysen, Klein, Laporte, and
  Pavone}]{Schiffer2022_CA_MB}
Schiffer, M., Boysen, N., Klein, P., Laporte, G., Pavone, M., 2022. Optimal
  picking policies in e-commerce warehouses. Management Science 68~(10),
  7497--7517.

\bibitem[{Schlkopf et~al.(2018)Schlkopf, Smola, and Bach}]{Schlkopf18_LCAL}
Schlkopf, B., Smola, A.~J., Bach, F., 2018. Learning with Kernels: Support
  Vector Machines, Regularization, Optimization, and Beyond. The MIT Press,
  Cambridge MA.

\bibitem[{Schmenner(2004)}]{schmenner_service_2004_JHLS}
Schmenner, R.~W., 2004. Service {Businesses} and {Productivity}. Decision
  Sciences 35~(3), 333--347.

\bibitem[{Sch\"obel(2012)}]{Shoebel2012_DC}
Sch\"obel, A., 2012. Line planning in public transportation: models and
  methods. OR Spectrum 34, 491--510.

\bibitem[{Scholten et~al.(2015)Scholten, Schuwirth, Reichert, and
  Lienert}]{Scholten2015-ja_JL}
Scholten, L., Schuwirth, N., Reichert, P., Lienert, J., 2015. Tackling
  uncertainty in multi-criteria decision analysis -- an application to water
  supply infrastructure planning. European Journal of Operational Research
  242~(1), 243--260.

\bibitem[{Schrijver(1986)}]{Sc86_ALAL}
Schrijver, A., 1986. Theory of Linear and Integer Programming. Wiley, New York.

\bibitem[{Schrijver(1998)}]{S98_SMPT}
Schrijver, A., 1998. Theory of linear and integer programming. John Wiley \&
  Sons, Chichester.

\bibitem[{Schrijver(2002)}]{schrijver2002history_JLYHK}
Schrijver, A., 2002. On the history of the transportation and maximum flow
  problems. Mathematical Programming 91~(3), 437--445.

\bibitem[{Schrijver(2003)}]{Schrijver:2003_IL}
Schrijver, A., 2003. Combinatorial Optimization - Polyhedra and Efficiency.
  Springer.

\bibitem[{Schrijver(2005)}]{S05_SMPT}
Schrijver, A., 2005. On the history of combinatorial optimization (till 1960).
  In: Aardal, K., Nemhauser, G.~L., Weismantel, R. (Eds.), Handbooks in
  Operations Research and Management Science. Vol.~12. Elsevier, pp. 1--68.

\bibitem[{Schummer and Vohra(2003)}]{Schummer03_BC}
Schummer, J., Vohra, R.~V., 2003. Auctions for procuring options. Operations
  Research 51~(1), 41--51.

\bibitem[{Sch{\"u}tz and Stanley-Lockman(2017)}]{Schutz2017_JLYHK}
Sch{\"u}tz, T., Stanley-Lockman, Z., 2017. {Smart logistics for future armed
  forces}. Brief~30, European Union Institute for Security Studies.

\bibitem[{Schwartz and Trigeorgis(2004)}]{schwartz04_MB}
Schwartz, E.~S., Trigeorgis, L., 2004. Real options and investment under
  uncertainty: Classical readings and recent contributions. MIT press,
  Cambridge, MA.

\bibitem[{Schwetschenau et~al.(2019)Schwetschenau, VanBriesen, and
  Cohon}]{Schwetschenau_S_E2019-ws_HL}
Schwetschenau, S.~E., VanBriesen, J.~M., Cohon, J.~L., 2019. Integrated
  simulation and optimization models for treatment plant placement in drinking
  water systems. Journal of Water Resources Planning and Management 145~(11),
  04019047.

\bibitem[{Scott et~al.(2013)Scott, Cavana, and Cameron}]{Scott2013-fg_AFRH}
Scott, R.~J., Cavana, R.~Y., Cameron, D., 2013. Evaluating immediate and
  long-term impacts of qualitative group model building workshops on
  participants' mental models. System Dynamics Review 29~(4), 216--236.

\bibitem[{{SDG}(2022)}]{SDG_BYKOK}
{SDG}, 2022. Sustainable development goalsu. \url{https://sdgs.un.org/goals},
  accessed on 2022-10-24.

\bibitem[{Secomandi(2001)}]{Secomandi2001-lw_HL}
Secomandi, N., 2001. A rollout policy for the vehicle routing problem with
  stochastic demands. Operations Research 49~(5), 796--802.

\bibitem[{See et~al.(2021)See, Md~Hamzah, and Yu}]{See2021-yy_SL}
See, K.~F., Md~Hamzah, N., Yu, M.-M., 2021. Metafrontier efficiency analysis
  for hospital pharmacy services using dynamic network {DEA} framework.
  Socio-Economic Planning Sciences 78, 101044.

\bibitem[{Seidl et~al.(2016)Seidl, Kaplan, Caulkins, Wrzaczek, and
  Feichtinger}]{Seidl2016-ga_KVRH}
Seidl, A., Kaplan, E.~H., Caulkins, J.~P., Wrzaczek, S., Feichtinger, G., 2016.
  Optimal control of a terror queue. European Journal of Operational Research
  248~(1), 246--256.

\bibitem[{Sein et~al.(2011)Sein, Henfridsson, Purao, Rossi, and
  Lindgren}]{sein_action_2011_JHLS}
Sein, M.~K., Henfridsson, O., Purao, S., Rossi, M., Lindgren, R., 2011. Action
  {Design} {Research}. MIS Quarterly 35~(1), 37--56.

\bibitem[{Selten(1975)}]{Selten1975-rl_GZ}
Selten, R., 1975. Reexamination of the perfectness concept for equilibrium
  points in extensive games. International Journal of Game Theory 4~(1),
  25--55.

\bibitem[{Serfozo(1999)}]{dick_HATI}
Serfozo, R., 1999. Introduction to Stochastic Networks. Springer-Verlag, New
  York.

\bibitem[{Servranckx and
  Vanhoucke(2019{\natexlab{a}})}]{Servranckx2019-zs_WH_ED}
Servranckx, T., Vanhoucke, M., 2019{\natexlab{a}}. Strategies for project
  scheduling with alternative subgraphs under uncertainty: similar and
  dissimilar sets of schedules. European Journal of Operational Research
  279~(1), 38--53.

\bibitem[{Servranckx and
  Vanhoucke(2019{\natexlab{b}})}]{Servranckx2019-vw_WH_ED}
Servranckx, T., Vanhoucke, M., 2019{\natexlab{b}}. A tabu search procedure for
  the resource-constrained project scheduling problem with alternative
  subgraphs. European Journal of Operational Research 273~(3), 841--860.

\bibitem[{Sethi and Cheng(1997)}]{sethi1997optimality_JSS}
Sethi, S.~P., Cheng, F., 1997. Optimality of (s, s) policies in inventory
  models with markovian demand. Operations Research 45~(6), 931--939.

\bibitem[{Sethi and Thompson(1970)}]{Sethi1970_XW}
Sethi, S.~P., Thompson, G.~L., 1970. Applications of mathematical control
  theory to finance: Modeling simple dynamic cash balance problems. Journal of
  Financial and Quantitative Analysis 5~(4-5), 381--394.

\bibitem[{Sethi and Thompson(2009)}]{Sethi2009_XW}
Sethi, S.~P., Thompson, G.~L., 2009. Optimal Control Theory: Applications to
  Management Science and Economics. Springer.

\bibitem[{Shaabani(2022)}]{Shaabani2022_JLYHK}
Shaabani, H., 2022. A literature review of the perishable inventory routing
  problem. The Asian Journal of Shipping and Logistics 38~(3), 143--161.

\bibitem[{Shah et~al.(2022)Shah, Hao, Yan, Yasmeen, Padda, and
  Ullah}]{Shah_2022_DSRW}
Shah, W. U.~H., Hao, G., Yan, H., Yasmeen, R., Padda, I. U.~H., Ullah, A.,
  2022. The impact of trade, financial development and government integrity on
  energy efficiency: An analysis from {G}7-{C}ountries. Energy 255~(1), 124507.

\bibitem[{Shamir(1987)}]{Shamir1987-pu_JMB}
Shamir, R., 1987. The efficiency of the simplex method: A survey. Management
  Science 33~(3), 301--334.

\bibitem[{Shang and Song(2003)}]{shang2003newsvendor_JSS}
Shang, K.~H., Song, J.-S., 2003. Newsvendor bounds and heuristic for optimal
  policies in serial supply chains. Management Science 49~(5), 618--638.

\bibitem[{Shang et~al.(2023)Shang, Song, and Zhou}]{shang2023single_JSS}
Shang, K.~H., Song, J.-S., Zhou, S.~X., 2023. Single-stage approximations of
  multiechelon inventory models. In: Song, J.-S. (Ed.), Research Handbook on
  Inventory Management. Edward Elgar Publishing.

\bibitem[{Shapiro(2003)}]{Shapiro2003-sm_HL}
Shapiro, A., 2003. Monte {C}arlo sampling methods. In: Ruszczynski, A.,
  Shapiro, A. (Eds.), Handbooks in Operations Research and Management Science.
  Horth-Holland, Amsterdam, pp. 353--425.

\bibitem[{Shapiro et~al.(2021)Shapiro, Dentcheva, and
  Ruszczy{\'n}ski}]{Shapiro2021-wq_HL}
Shapiro, A., Dentcheva, D., Ruszczy{\'n}ski, A., 2021. Lectures on Stochastic
  Programming: Modeling and Theory, 3rd Edition. SIAM.

\bibitem[{Shapley(1953)}]{Sh1953_GZ}
Shapley, L., 1953. A value for $n$-person games. In: Kuhn, H., Tucker, A.
  (Eds.), Contributions to the Theory of Games {II}. Princeton University
  Press, Princeton, NJ, pp. 307--317.

\bibitem[{Sharp et~al.(2007)Sharp, Meng, and Liu}]{Sharp2007-xj_SL}
Sharp, J.~A., Meng, W., Liu, W., 2007. A modified slacks-based measure model
  for data envelopment analysis with `natural' negative outputs and inputs.
  Journal of the Operational Research Society 58~(12), 1672--1677.

\bibitem[{Shaw et~al.(2004)Shaw, Westcombe, Hodgkin, and
  Montibeller}]{Shaw2004-gw_MYLW}
Shaw, D., Westcombe, M., Hodgkin, J., Montibeller, G., 2004. Problem
  structuring methods for large group interventions. Journal of the Operational
  Research Society 55~(5), 453--463.

\bibitem[{Shaw et~al.(2016)Shaw, Irfan, Shankar, and Yadav}]{Shaw2016-sh_HL}
Shaw, K., Irfan, M., Shankar, R., Yadav, S.~S., 2016. Low carbon chance
  constrained supply chain network design problem: a benders decomposition
  based approach. Computers \& Industrial Engineering 98, 483--497.

\bibitem[{Sherali(1984)}]{Sh1984_GZ}
Sherali, H., 1984. A multiple leader {S}tackelberg model and analysis.
  Operations Research 32~(2), 229--476.

\bibitem[{Sherali and Adams(1990)}]{SA90_ALAL}
Sherali, H., Adams, W., 1990. A hierarchy of relaxations between the continuous
  and convex hull representations for zero-one programming problems. SIAM
  Journal on Discrete Mathematics 3, 411--430.

\bibitem[{Sherali et~al.(1991)Sherali, Carter, and Hobeika}]{Sherali1991_JLYHK}
Sherali, H.~D., Carter, T.~B., Hobeika, A.~G., 1991. A location-allocation
  model and algorithm for evacuation planning under hurricane/flood conditions.
  Transportation Research Part B: Methodological 25~(6), 439--452.

\bibitem[{Sheu(2014)}]{Sheu2014_JLYHK}
Sheu, J.-B., 2014. Post-disaster relief–service centralized logistics
  distribution with survivor resilience maximization. Transportation Research
  Part B: Methodological 68, 288--314.

\bibitem[{Shi(2023)}]{shi2023approximation_JSS}
Shi, C., 2023. Approximation algorithms for stochastic inventory systems. In:
  Song, J.-S. (Ed.), Research Handbook on Inventory Management. Edward Elgar
  Publishing.

\bibitem[{Shi and Wang(2018)}]{Shi2018-eb_JEB}
Shi, Z., Wang, G., 2018. Integration of big-data {ERP} and business analytics
  ({BA}). The Journal of High Technology Management Research 29~(2), 141--150.

\bibitem[{Shipley et~al.(2018)Shipley, Coy, and Shipley}]{Shipley2018_EA}
Shipley, M.~F., Coy, S.~P., Shipley, J.~B., 2018. Utilizing statistical
  significance in fuzzy interval-valued evidence sets for assessing artificial
  reef structure impact. Journal of the Operational Research Society 69~(6),
  905--918.

\bibitem[{Shneiderman(1996)}]{Shneiderman1996-an_MJE}
Shneiderman, B., 1996. The eyes have it: a task by data type taxonomy for
  information visualizations. In: Proceedings 1996 {IEEE} Symposium on Visual
  Languages. pp. 336--343.

\bibitem[{Shor(1977)}]{Shor77_EAY}
Shor, N.~Z., 1977. Cut-off method with space extension in convex programming
  problems. Cybernetics 13~(1), 94--96.

\bibitem[{Shtub et~al.(2004)Shtub, Bard, and Globerson}]{Shtub2004-ob_WH_ED}
Shtub, A., Bard, J.~F., Globerson, S., 2004. Project Management: Engineering,
  Technology and Implementation, 15th Edition. Pearson Custom Publishing.

\bibitem[{Si et~al.(2004)Si, Barto, Powell, and Wunsch}]{Si2004-mp_HL}
Si, J., Barto, A.~G., Powell, W.~B., Wunsch, D., 2004. Handbook of Learning and
  Approximate Dynamic Programming. John Wiley \& Sons.

\bibitem[{Siegrist and {\'A}rvai(2020)}]{Siegrist2020-ul_TC}
Siegrist, M., {\'A}rvai, J., 2020. Risk perception: Reflections on 40 years of
  research. Risk Analysis 40~(S1), 2191--2206.

\bibitem[{Silva et~al.(2016)Silva, Oliveira, and W{\"a}scher}]{Silva2016-jj_JB}
Silva, E., Oliveira, J.~F., W{\"a}scher, G., 2016. The pallet loading problem:
  a review of solution methods and computational experiments. International
  Transactions in Operational Research 23~(1-2), 147--172.

\bibitem[{Silver(2004)}]{silver2004overview_COIT}
Silver, E.~A., 2004. An overview of heuristic solution methods. Journal of the
  Operational Research Society 55~(9), 936--956.

\bibitem[{Silver and Meal(1973)}]{silver1973heuristic_JSS}
Silver, E.~A., Meal, C., 1973. A heuristic for selecting lot size quantities
  for the case of a deterministic time-varying demand rate and discrete
  opportunities for replenishment. Production and Inventory Management 2,
  64--74.

\bibitem[{Silver et~al.(1988)Silver, Pyke, and
  Petterson}]{silver1988inventory_JSS}
Silver, E.~A., Pyke, D.~F., Petterson, R., 1988. Inventory Management and
  Production Planning and Scheduling. John Wiley \& Sons.

\bibitem[{Silver et~al.(1980)Silver, Victor, Vidal, and
  de~Werra}]{silver1980tutorial_COIT}
Silver, E.~A., Victor, R., Vidal, V., de~Werra, D., 1980. A tutorial on
  heuristic methods. European Journal of Operational Research 5~(3), 153--162.

\bibitem[{{Simantics System Dynamics}(2017)}]{Sysdyn_CTCGE}
{Simantics System Dynamics}, 2017. {SysDyn}.
\newline\urlprefix\url{http://sysdyn.simantics.org/}

\bibitem[{Sim{\~a}o et~al.(2009)Sim{\~a}o, Day, George, Gifford, Nienow, and
  Powell}]{Simao2009-ag_HL}
Sim{\~a}o, H.~P., Day, J., George, A.~P., Gifford, T., Nienow, J., Powell,
  W.~B., 2009. An approximate dynamic programming algorithm for {Large-Scale}
  fleet management: A case application. Transportation Science 43~(2),
  178--197.

\bibitem[{Simar and Wilson(2011)}]{Simar2011-oz_JJ}
Simar, L., Wilson, P.~W., 2011. Two-stage {DEA}: \emph{caveat emptor}. Journal
  of Productivity Analysis 36~(2), 205--218.

\bibitem[{Simchi-Levi et~al.(2014)Simchi-Levi, Chen, and
  Bramel}]{simchi2014logic_JSS}
Simchi-Levi, D., Chen, X., Bramel, J., 2014. The logic of logistics: Theory,
  algorithms, and applications for logistics and supply chain management.
  Springer.

\bibitem[{Simchi-Levi et~al.(2019)Simchi-Levi, Trichakis, and
  Zhang}]{Simchi-Levi2019-pv_KVRH}
Simchi-Levi, D., Trichakis, N., Zhang, P.~Y., 2019. Designing response supply
  chain against bioattacks. Operations Research 67~(5), 1246--1268.

\bibitem[{Simon(1992)}]{simon1992simple_HATI}
Simon, B., 1992. A simple relationship between light and heavy traffic limits.
  Operations Research 40~(3-supplement-2), S342--S345.

\bibitem[{Simon(1952)}]{Simon1952_XW}
Simon, H.~A., 1952. On the application of servomechanism theory in the study of
  production control. Econometrica 20~(2), 247--268.

\bibitem[{Simon(1968)}]{simonAlmostPracticalSolution1968_AKSJF}
Simon, J.~L., 1968. An {{Almost Practical Solution}} to {{Airline
  Overbooking}}. Journal of Transport Economics and Policy 2~(2), 201--202.

\bibitem[{Simoni et~al.(2020)Simoni, Kutanoglu, and
  Claudel}]{simoni2020optimization_CKTVW}
Simoni, M.~D., Kutanoglu, E., Claudel, C.~G., 2020. Optimization and analysis
  of a robot-assisted last mile delivery system. Transportation Research Part
  E: Logistics and Transportation Review 142, 102049.

\bibitem[{Simpson(1978)}]{Simpson1978_JLYHK}
Simpson, V.~P., 1978. Optimum solution structure for a repairable inventory
  problem. Operations Research 26~(2), 270--281.

\bibitem[{Sims(1980)}]{sims1980_FP}
Sims, C., 1980. {Macroeconomics and Reality}. Econometrica 48~(1), 1--48.

\bibitem[{Sims and Smithin(1982)}]{Sims1982-sn_AJG}
Sims, D., Smithin, T., 1982. Voluntary operational research. Journal of the
  Operational Research Society 33~(1), 21--28.

\bibitem[{{Simul8 Corporation}(2022)}]{Simul8_CTCGE}
{Simul8 Corporation}, 2022. {Simul8}.
\newline\urlprefix\url{https://www.simul8.com/}

\bibitem[{Skagerlund et~al.(2020)Skagerlund, Forsblad, Slovic, and
  V{\"a}stfj{\"a}ll}]{Skagerlund2020-cs_TC}
Skagerlund, K., Forsblad, M., Slovic, P., V{\"a}stfj{\"a}ll, D., 2020. The
  affect heuristic and risk perception - stability across elicitation methods
  and individual cognitive abilities. Frontiers in Psychology 11, 970.

\bibitem[{Skorin-Kapov et~al.(2006)Skorin-Kapov, Skorin-Kapov, and
  Boljunčic}]{skorin-kapov.skorin-kapov.ea:06_BF}
Skorin-Kapov, D., Skorin-Kapov, J., Boljunčic, V., 2006. Location problems in
  telecommunications. In: Resende, M., Pardalos, P. (Eds.), Handbook of
  Optimization in Telecommunications. Springer, New-York, pp. 517--544.

\bibitem[{Slotine et~al.(1991)Slotine, Li, et~al.}]{Slotine1991_XW}
Slotine, J.-J.~E., Li, W., et~al., 1991. Applied nonlinear control. Vol. 199.
  Prentice Hall, Englewood Cliffs, NJ.

\bibitem[{Slovic(2000)}]{Slovic2000-zr_TC}
Slovic, P., 2000. The Perception of Risk. Earthscan Publications.

\bibitem[{Sluijk et~al.(2022{\natexlab{a}})Sluijk, Florio, Kinable, Dellaert,
  and Van~Woensel}]{sluijk2022chance_CKTVW}
Sluijk, N., Florio, A.~M., Kinable, J., Dellaert, N., Van~Woensel, T.,
  2022{\natexlab{a}}. A chance-constrained two-echelon vehicle routing problem
  with stochastic demands. Transportation Science, DOI: 10.1287/trsc.2022.1162.

\bibitem[{Sluijk et~al.(2022{\natexlab{b}})Sluijk, Florio, Kinable, Dellaert,
  and Van~Woensel}]{sluijk2022two_CKTVW}
Sluijk, N., Florio, A.~M., Kinable, J., Dellaert, N., Van~Woensel, T.,
  2022{\natexlab{b}}. Two-echelon vehicle routing problems: A literature
  review. European Journal of Operational Research 304~(3), 865--886.

\bibitem[{Smet et~al.(2016)Smet, Brucker, De~Causmaecker, and
  Vanden~Berghe}]{smet2016polynomially_GVBSP}
Smet, P., Brucker, P., De~Causmaecker, P., Vanden~Berghe, G., 2016.
  Polynomially solvable personnel rostering problems. European Journal of
  Operational Research 249~(1), 67--75.

\bibitem[{Smidts(1997)}]{Smidts1997-ma_TC}
Smidts, A., 1997. The relationship between risk attitude and strength of
  preference: A test of intrinsic risk attitude. Management Science 43~(3),
  357--370.

\bibitem[{Smith and Shaw(2019)}]{Smith2019-pw_JL}
Smith, C.~M., Shaw, D., 2019. The characteristics of problem structuring
  methods: A literature review. European Journal of Operational Research
  274~(2), 403--416.

\bibitem[{Smits(2010)}]{Smits2010-dq_CV}
Smits, M., 2010. Impact of policy and process design on the performance of
  intake and treatment processes in mental health care: a system dynamics case
  study. Journal of the Operational Research Society 61~(10), 1437--1445.

\bibitem[{Sniedovich and Vi{\ss}(2006)}]{sniedovich2006corridor_COIT}
Sniedovich, M., Vi{\ss}, S., 2006. The corridor method: a dynamic programming
  inspired metaheuristic. Control and Cybernetics 35~(3), 551--578.

\bibitem[{Snoeck et~al.(2020)Snoeck, Merch{\'a}n, and
  Winkenbach}]{snoeck2020revenue_CKTVW}
Snoeck, A., Merch{\'a}n, D., Winkenbach, M., 2020. Revenue management in
  last-mile delivery: state-of-the-art and future research directions.
  Transportation Research Procedia 46, 109--116.

\bibitem[{Snyder and Shen(2019)}]{snyder2019fundamentals_JSS}
Snyder, L.~V., Shen, Z.-J.~M., 2019. Fundamentals of supply chain theory. John
  Wiley \& Sons.

\bibitem[{Sobolev et~al.(2011)Sobolev, Sanchez, and
  Vasilakis}]{Sobolev2011-xv_CV}
Sobolev, B.~G., Sanchez, V., Vasilakis, C., 2011. Systematic review of the use
  of computer simulation modeling of patient flow in surgical care. Journal of
  Medical Systems 35~(1), 1--16.

\bibitem[{Soeffker et~al.(2022)Soeffker, Ulmer, and
  Mattfeld}]{soeffker2021stochastic_CA_MB}
Soeffker, N., Ulmer, M., Mattfeld, D., 2022. Stochastic dynamic vehicle routing
  in the light of prescriptive analytics: A review. European Journal of
  Operational Research 298~(3), 801--820.

\bibitem[{Solomon et~al.(2019)Solomon, Li, Womer, and
  Santos}]{Solomon2019-cu_HL}
Solomon, S., Li, H., Womer, K., Santos, C., 2019. Multiperiod stochastic
  resource planning in professional services organizations. Decision Sciences
  50~(6), 1281--1318.

\bibitem[{Solomon et~al.(2022)Solomon, Pannirselvam, and
  Li}]{Solomon2022-ou_HL}
Solomon, S., Pannirselvam, G.~P., Li, H., 2022. Using integrated
  simulation-optimisation to optimise staffing decisions in a service supply
  chain. International Journal of Integrated Supply Management 15~(1), 1--26.

\bibitem[{Soltani and Lozano(2020)}]{Soltani2020-xv_SL}
Soltani, N., Lozano, S., 2020. Interactive multiobjective {DEA} target setting
  using lexicographic {DDF}. RAIRO - Operations Research 54~(6), 1703--1722.

\bibitem[{Solyal{\i} and S{\"u}ral(2022)}]{solyali2022effective_COIT}
Solyal{\i}, O., S{\"u}ral, H., 2022. An effective matheuristic for the
  multivehicle inventory routing problem. Transportation Science 56,
  1044--1057.

\bibitem[{Song et~al.(2014)Song, Nelson, and Pegden}]{SongWSC2014_CC}
Song, E., Nelson, B.~L., Pegden, C.~D., 2014. Advanced tutorial: Input
  uncertainty quantification. In: Tolk, A., Diallo, S.~Y., Ryzhov, I.~O.,
  Yilmaz, L., Buckley, S., Miller, J.~A. (Eds.), Proceedings of the 2014 Winter
  Simulation Conference. IEEE Piscataway, pp. 162--176.

\bibitem[{Song(1994)}]{song1994effect_JSS}
Song, J.-S., 1994. The effect of leadtime uncertainty in a simple stochastic
  inventory model. Management Science 40~(5), 603--613.

\bibitem[{Song(2023)}]{song2023handbook_JSS}
Song, J.-S. (Ed.), 2023. Research Handbook on Inventory Management. Edward
  Elgar Publishing.

\bibitem[{Song et~al.(2017)Song, Xiao, Zhang, and Zipkin}]{song2017optimal_JSS}
Song, J.-S., Xiao, L., Zhang, H., Zipkin, P., 2017. Optimal policies for a
  dual-sourcing inventory problem with endogenous stochastic lead times.
  Operations Research 65~(2), 379--395.

\bibitem[{Song and Zhang(2020)}]{Song2020-ti_KESXZ}
Song, J.-S., Zhang, Y., 2020. Stock or print? {I}mpact of {3-D} printing on
  spare parts logistics. Management Science 66~(9), 3860--3878.

\bibitem[{Song and Zipkin(1993)}]{song1993inventory_JSS}
Song, J.-S., Zipkin, P., 1993. Inventory control in a fluctuating demand
  environment. Operations Research 41~(2), 351--370.

\bibitem[{Song and Zipkin(2003)}]{song2003supply_JSS}
Song, J.-S., Zipkin, P., 2003. Supply chain operations: {Assemble-to-Order}
  systems. In: Graves, S.~C., de~Kok, A.~G. (Eds.), Handbooks in Operations
  Research and Management Science. Vol.~11. Elsevier, pp. 561--596.

\bibitem[{Song and Zipkin(1996)}]{song1996inventory_JSS}
Song, J.-S., Zipkin, P.~H., 1996. Inventory control with information about
  supply conditions. Management Science 42~(10), 1409--1419.

\bibitem[{Sonnessa et~al.(2017)Sonnessa, T{\`a}nfani, and
  Testi}]{Sonnessa2017-la_AFRH}
Sonnessa, M., T{\`a}nfani, E., Testi, A., 2017. An agent-based simulation model
  to evaluate alternative co-payment scenarios for contributing to health
  systems financing. Journal of the Operational Research Society 68~(5),
  591--604.

\bibitem[{Soorapanth et~al.(2022)Soorapanth, Eldabi, and
  Young}]{Soorapanth2022-jt_CV}
Soorapanth, S., Eldabi, T., Young, T., 2022. Towards a framework for evaluating
  the costs and benefits of simulation modelling in healthcare. Journal of the
  Operational Research Society, DOI: 10.1080/01605682.2022.2064780.

\bibitem[{Soyster(1973)}]{Soyster1973-hb_HL}
Soyster, A.~L., 1973. Convex programming with {Set-Inclusive} constraints and
  applications to inexact linear programming. Operations Research 21~(5),
  1154--1157.

\bibitem[{Speranza and Ukovich(1994)}]{speranza1994minimizing_JLYHK}
Speranza, M.~G., Ukovich, W., 1994. Minimizing transportation and inventory
  costs for several products on a single link. Operations Research 42~(5),
  879--894.

\bibitem[{Spliet et~al.(2018)Spliet, Dabia, and
  Van~Woensel}]{spliet2018time_CKTVW}
Spliet, R., Dabia, S., Van~Woensel, T., 2018. The time window assignment
  vehicle routing problem with time-dependent travel times. Transportation
  Science 52~(2), 261--276.

\bibitem[{Srinivasan et~al.(2018)Srinivasan, Satyajit, Behera, and
  Panigrahi}]{Srinivasan2018_SMD}
Srinivasan, K., Satyajit, S., Behera, B.~K., Panigrahi, P.~K., 2018. Efficient
  quantum algorithm for solving travelling salesman problem: An {IBM} quantum
  experience. arXiv:1805.10928.

\bibitem[{Stadtler(2003)}]{stadtler2003multilevel_COIT}
Stadtler, H., 2003. Multilevel lot sizing with setup times and multiple
  constrained resources: Internally rolling schedules with lot-sizing windows.
  Operations Research 51~(3), 487--502.

\bibitem[{Starita et~al.(2020)Starita, Strauss, Fei, Jovanovi{\'c}, Ivanov,
  Pavlovi{\'c}, and Fichert}]{starita2020air_VLVV}
Starita, S., Strauss, A.~K., Fei, X., Jovanovi{\'c}, R., Ivanov, N.,
  Pavlovi{\'c}, G., Fichert, F., 2020. Air traffic control capacity planning
  under demand and capacity provision uncertainty. Transportation Science
  54~(4), 882--896.

\bibitem[{Stecke(1983)}]{Stecke1983-lc_KESXZ}
Stecke, K.~E., 1983. Formulation and solution of nonlinear integer production
  planning problems for flexible manufacturing systems. Management Science
  29~(3), 273--288.

\bibitem[{Stecke and Solberg(1981)}]{Stecke1981-eb_KESXZ}
Stecke, K.~E., Solberg, J.~J., 1981. Loading and control policies for a
  flexible manufacturing system. International Journal of Production Research
  19~(5), 481--490.

\bibitem[{Stecke et~al.(2012)Stecke, Yin, Kaku, and
  Murase}]{Stecke2012-rg_KESXZ}
Stecke, K.~E., Yin, Y., Kaku, I., Murase, Y., 2012. Seru: The organizational
  extension of {JIT} for a {Super-Talent} factory. International Journal of
  Strategic Decision Sciences 3~(1), 106--119.

\bibitem[{Steinbach(2001)}]{Steinbach01_MB}
Steinbach, M.~C., 2001. Markowitz revisited: Mean-variance models in financial
  portfolio analysis. SIAM Review 43~(1), 31--85.

\bibitem[{Stephens et~al.(2018)Stephens, Lewis, and
  Reddy}]{Stephens2018-ov_AJG}
Stephens, A., Lewis, E.~D., Reddy, S.~M., 2018. Inclusive Systemic Evaluation
  ({ISE4GEMs)}: A New Approach for the {SDG} Era. UN Women, New York, NY.

\bibitem[{Sterman(1989)}]{Sterman1989_XW}
Sterman, J.~D., 1989. Modeling managerial behavior: Misperceptions of feedback
  in a dynamic decision making experiment. Management Science 35~(3), 321--339.

\bibitem[{Sterman(2000)}]{Sterman2000_SMD}
Sterman, J.~D., 2000. Business Dynamics: Systems Thinking and Modeling for a
  Complex World. Irwin/McGraw-Hill, Boston, MA.

\bibitem[{Sterman(2014)}]{Sterman2014-du_MCJM}
Sterman, J.~D., 2014. Interactive web-based simulations for strategy and
  sustainability: The {MIT} sloan {LearningEdge} management flight simulators,
  {P}art {I}. System Dynamics Review 30~(1-2), 89--121.

\bibitem[{Sterman(2018)}]{Sterman2018-cg_MCJM}
Sterman, J.~D., 2018. System dynamics at sixty: the path forward. System
  Dynamics Review 34~(1-2), 5--47.

\bibitem[{Steuer(1985)}]{steuer_MESG}
Steuer, R., 1985. Multiple Criteria Optimization: Theory, Computation and
  Application. John Wiley \& Sons, New York.

\bibitem[{Stewart et~al.(2013)Stewart, French, and Rios}]{Stewart2013-ez_JL}
Stewart, T.~J., French, S., Rios, J., 2013. Integrating multicriteria decision
  analysis and scenario planning---review and extension. Omega 41~(4),
  679--688.

\bibitem[{Stidham(2002)}]{stidham_HATI}
Stidham, Jr, S., 2002. Analysis, design, and control of queueing systems.
  Operations Research 50~(1), 197--216.

\bibitem[{Stolyar(2004)}]{stolyar2004maxweight_HATI}
Stolyar, A.~L., 2004. Maxweight scheduling in a generalized switch: State space
  collapse and workload minimization in heavy traffic. The Annals of Applied
  Probability 14~(1), 1--53.

\bibitem[{Stoyan et~al.(2016)Stoyan, Pankratov, and
  Romanova}]{Stoyan2016-fk_JB}
Stoyan, Y., Pankratov, A., Romanova, T., 2016. Cutting and packing problems for
  irregular objects with continuous rotations: mathematical modelling and
  non-linear optimization. Journal of the Operational Research Society 67~(5),
  786--800.

\bibitem[{Strang(1987)}]{Strang1987-zh_JMB}
Strang, G., 1987. Karmarkar's algorithm and its place in applied mathematics.
  The Mathematical Intelligencer 9~(2), 4--10.

\bibitem[{Strauss et~al.(2021)Strauss, G{\"u}lp{\i}nar, and
  Zheng}]{strauss2021dynamic_CKTVW}
Strauss, A., G{\"u}lp{\i}nar, N., Zheng, Y., 2021. Dynamic pricing of flexible
  time slots for attended home delivery. European Journal of Operational
  Research 294~(3), 1022--1041.

\bibitem[{Strauss et~al.(2018)Strauss, Klein, and
  Steinhardt}]{straussReviewChoicebasedRevenue2018_AKSJF}
Strauss, A.~K., Klein, R., Steinhardt, C., 2018. A review of choice-based
  revenue management: Theory and methods. European Journal of Operational
  Research 271~(2), 375--387.

\bibitem[{Street et~al.(2016)Street, Di~Mauro, Humphrey, Johns, Boyd,
  Crawford-Brown, Evans, Kitchen, Hunt, Knox, Low, McCall, Watkiss, and
  Wilby}]{report_BYKOK}
Street, R., Di~Mauro, M., Humphrey, K., Johns, D., Boyd, E., Crawford-Brown,
  D., Evans, J., Kitchen, J., Hunt, A., Knox, K., Low, R., McCall, R., Watkiss,
  P., Wilby, R., 2016. {UK} climate change risk assessment evidence report:
  Chapter 8, {Cross-cutting Issues}. Report prepared for the Adaptation
  Sub-Committee of the Committee on Climate Change, London.

\bibitem[{Strogatz(2018)}]{strogatz2018_XW}
Strogatz, S.~H., 2018. Nonlinear dynamics and chaos: with applications to
  physics, biology, chemistry, and engineering. CRC press.

\bibitem[{Strozzi et~al.(2017)Strozzi, Colicchia, Creazza, and
  No{\`e}}]{Strozzi2017-ff_KESXZ}
Strozzi, F., Colicchia, C., Creazza, A., No{\`e}, C., 2017. Literature review
  on the `smart factory' concept using bibliometric tools. International
  Journal of Production Research 55~(22), 6572--6591.

\bibitem[{STT(2021)}]{UGENT_STT_competition_GVBSP}
STT, 2021. {ITC} 2021: International timetabling competition on sports
  timetabling. \url{https://www.sportscheduling.ugent.be/ITC2021/}, accessed on
  2022-10-09.

\bibitem[{Su(2007)}]{suIntertemporalPricingStrategic2007_AKSJF}
Su, X., 2007. Intertemporal pricing with strategic customer behavior.
  Management Science 53~(5), 726--741.

\bibitem[{Subramanian et~al.(2013)Subramanian, Uchoa, and
  Ochi}]{SUBRAMANIAN20132519_CA_MB}
Subramanian, A., Uchoa, E., Ochi, L., 2013. A hybrid algorithm for a class of
  vehicle routing problems. Computers \& Operations Research 40~(10),
  2519--2531.

\bibitem[{Subramanian and Subramanyam(2012)}]{Subramanian2012-om_AM}
Subramanian, R., Subramanyam, R., 2012. Key factors in the market for
  remanufactured products. Manufacturing \& Service Operations Management
  14~(2), 315--326.

\bibitem[{Sueyoshi et~al.(2017)Sueyoshi, Yuan, and Goto}]{Sueyoshi_2017_DSRW}
Sueyoshi, T., Yuan, Y., Goto, M., 2017. A literature study for {DEA} applied to
  energy and environment. Energy Economics 62~(1), 104--124.

\bibitem[{Sun and Conejo(2021)}]{Sun2021-pm_HL}
Sun, X.~A., Conejo, A.~J., 2021. Robust Optimization in Electric Energy
  Systems. Springer International Publishing.

\bibitem[{Sun and Hassanlou(2019)}]{Sun2019-lq_HL}
Sun, Y., Hassanlou, K., 2019. Equity trading server allocation using chance
  constrained programming. Journal of Supply Chain and Operations Management
  17~(1), 1--13.

\bibitem[{Sutton and Barto(2018)}]{sutton2018reinforcement_LCAL}
Sutton, R.~S., Barto, A.~G., 2018. Reinforcement learning: An introduction. MIT
  press, Cambridge MA.

\bibitem[{Svetunkov(2022)}]{Svetunkov2022smooth_FP}
Svetunkov, I., 2022. smooth: Forecasting Using State Space Models. R package
  version 3.1.6.

\bibitem[{Swamy and Kumar(2004)}]{swamy.kumar:04_BF}
Swamy, C., Kumar, A., 2004. Primal-dual algorithms for connected facility
  location problems. Algorithmica 40, 245--269.

\bibitem[{Swan et~al.(2022)Swan, Adriaensen, Brownlee, Hammond, Johnson,
  Kheiri, Krawiec, Merelo, Minku, Özcan, Pappa, García-Sánchez, Sörensen,
  Voß, Wagner, and White}]{swan2022metaheuristics_COIT}
Swan, J., Adriaensen, S., Brownlee, A.~E., Hammond, K., Johnson, C.~G., Kheiri,
  A., Krawiec, F., Merelo, J., Minku, L.~L., Özcan, E., Pappa, G.~L.,
  García-Sánchez, P., Sörensen, K., Voß, S., Wagner, M., White, D.~R.,
  2022. Metaheuristics “in the large”. European Journal of Operational
  Research 297~(2), 393--406.

\bibitem[{Swaroop et~al.(2012)Swaroop, Zou, Ball, and
  Hansen}]{swaroop2012more_VLVV}
Swaroop, P., Zou, B., Ball, M.~O., Hansen, M., 2012. Do more us airports need
  slot controls? {A} welfare based approach to determine slot levels.
  Transportation Research Part B: Methodological 46~(9), 1239--1259.

\bibitem[{Sweeney et~al.(2020)Sweeney, Bessa, Browell, and
  Pinson}]{swe:bes:bro:pin:20_DSRW}
Sweeney, C., Bessa, R., Browell, J., Pinson, P., 2020. The future of
  forecasting for renewable energy. Wiley Interdisciplinary Reviews: Energy and
  Environment 9~(2), e365.

\bibitem[{Sydelko et~al.(2023)Sydelko, Espinosa, and
  Midgley}]{Sydelko2023-pa_GM}
Sydelko, P., Espinosa, A., Midgley, G., 2023. Designing interagency responses
  to wicked problems: A viable system model board game. European Journal of
  Operational Research, DOI: 10.1016/j.ejor.2023.06.040.

\bibitem[{Sydelko et~al.(2021)Sydelko, Midgley, and
  Espinosa}]{Sydelko2021-qp_GM}
Sydelko, P., Midgley, G., Espinosa, A., 2021. Designing interagency responses
  to wicked problems: Creating a common, cross-agency understanding. European
  Journal of Operational Research 294~(1), 250--263.

\bibitem[{Szeider(2003)}]{Szei03}
Szeider, S., 2003. On fixed-parameter tractable parameterizations of {SAT}. In:
  Proceedings of the 6th International Conference on Theory and Applications of
  Satisfiability Testing ({SAT}). Vol. 2919 of LNCS. pp. 188--202.

\bibitem[{Szeto and Jiang(2014)}]{Szeto2014_DC}
Szeto, W.~Y., Jiang, Y., 2014. Transit route and frequency design{: B}i-level
  modeling and hybrid artificial bee colony algorithm approach. Transportation
  Research Part B: Methodological 67, 235--263.

\bibitem[{Taillandier et~al.(2019)Taillandier, Gaudou, Grignard, Huynh,
  Marilleau, Caillou, Philippon, and Drogoul}]{Gama_CTCGE}
Taillandier, P., Gaudou, B., Grignard, A., Huynh, Q., Marilleau, N., Caillou,
  P., Philippon, D., Drogoul, A., 2019. Building, composing and experimenting
  complex spatial models with the {GAMA} platform. Geoinformatica 23, 299--322.

\bibitem[{Tak\'{a}cs(1962)}]{takas_HATI}
Tak\'{a}cs, L., 1962. Introduction to the Theory of Queues. University Texts in
  the Mathematical Sciences, Oxford University Press, New York.

\bibitem[{Tako(2015)}]{Tako2015-dk_AFRH}
Tako, A.~A., 2015. Exploring the model development process in discrete-event
  simulation: insights from six expert modellers. Journal of the Operational
  Research Society 66~(5), 747--760.

\bibitem[{Tako and Kotiadis(2015)}]{Tako2015-wi_MYLW}
Tako, A.~A., Kotiadis, K., 2015. {PartiSim}: A multi-methodology framework to
  support facilitated simulation modelling in healthcare. European Journal of
  Operational Research 244~(2), 555--564.

\bibitem[{Tako and Robinson(2010)}]{Tako2010-ua_AFRH}
Tako, A.~A., Robinson, S., 2010. Model development in discrete-event simulation
  and system dynamics: An empirical study of expert modellers. European Journal
  of Operational Research 207~(2), 784--794.

\bibitem[{Talluri and Van~Ryzin(1999)}]{talluri1999randomized_VLVV}
Talluri, K., Van~Ryzin, G., 1999. A randomized linear programming method for
  computing network bid prices. Transportation Science 33~(2), 207--216.

\bibitem[{Talluri and
  Van~Ryzin(2004{\natexlab{a}})}]{talluriRevenueManagementGeneral2004_AKSJF}
Talluri, K., Van~Ryzin, G., 2004{\natexlab{a}}. Revenue management under a
  general discrete choice model of consumer behavior. Management Science
  50~(1), 15--33.

\bibitem[{Talluri and
  Van~Ryzin(2004{\natexlab{b}})}]{talluriTheoryPracticeRevenue2004a_AKSJF}
Talluri, K., Van~Ryzin, G., 2004{\natexlab{b}}. The Theory and Practice of
  Revenue Management. International Series in Operations Research \& Management
  Science. {Springer}, {New York, NY}.

\bibitem[{Talmar et~al.(2020)Talmar, Walrave, Podoynitsyna, Holmström, and
  Romme}]{talmar_mapping_2020_JHLS}
Talmar, M., Walrave, B., Podoynitsyna, K.~S., Holmström, J., Romme, A. G.~L.,
  2020. Mapping, analyzing and designing innovation ecosystems: The ecosystem
  pie model. Long Range Planning 53~(4), 101850.

\bibitem[{Tan et~al.(2021)Tan, Goh, Kendall, and Sabar}]{Tan2021-ba_JJ}
Tan, J.~S., Goh, S.~L., Kendall, G., Sabar, N.~R., 2021. A survey of the
  state-of-the-art of optimisation methodologies in school timetabling
  problems. Expert Systems with Applications 165, 113943.

\bibitem[{Tang and Veelenturf(2019)}]{tang2019strategic_JLYHK}
Tang, C.~S., Veelenturf, L.~P., 2019. The strategic role of logistics in the
  industry 4.0 era. Transportation Research Part E: Logistics and
  Transportation Review 129, 1--11.

\bibitem[{Tarhan and O{\u{g}}uz(2022)}]{tarhan2022matheuristic_COIT}
Tarhan, {\.I}., O{\u{g}}uz, C., 2022. A matheuristic for the generalized order
  acceptance and scheduling problem. European Journal of Operational Research
  299~(1), 87--103.

\bibitem[{Tashman(2000)}]{tashman2000out_FP}
Tashman, L.~J., 2000. Out-of-sample tests of forecasting accuracy: an analysis
  and review. International Journal of Forecasting 16~(4), 437--450.

\bibitem[{Tavella et~al.(2021)Tavella, Papadopoulos, and
  Paroutis}]{Tavella2021-mj_AFRH}
Tavella, E., Papadopoulos, T., Paroutis, S., 2021. Artefact appropriation in
  facilitated modelling: An adaptive structuration theory approach. Journal of
  the Operational Research Society 72~(11), 2381--2395.

\bibitem[{Tena and Forrest(2007)}]{De_Dios_Tena2007-dp_IM}
Tena, J. d.~D., Forrest, D., 2007. Within-season dismissal of football coaches:
  Statistical analysis of causes and consequences. European Journal of
  Operational Researchh 181~(1), 362--373.

\bibitem[{Teng et~al.(2022)Teng, Zhang, and Sun}]{Teng2022-ra_JJ}
Teng, Y., Zhang, J., Sun, T., 2022. Data‐driven decision‐making model based
  on artificial intelligence in higher education system of colleges and
  universities. Expert Systems, e12820, DOI: 10.1111/exsy.12820.

\bibitem[{Terrab and Odoni(1993)}]{terrab1993strategic_VLVV}
Terrab, M., Odoni, A.~R., 1993. Strategic flow management for air traffic
  control. Operations Research 41~(1), 138--152.

\bibitem[{Teunter et~al.(2004)Teunter, van~der Laan, and
  Vlachos}]{R_Teunter2004-cp_AM}
Teunter, R., van~der Laan, E., Vlachos, D., 2004. Inventory strategies for
  systems with fast remanufacturing. Journal of the Operational Research
  Society 55~(5), 475--484.

\bibitem[{Thanassoulis et~al.(2016)Thanassoulis, De~Witte, Johnes, Johnes,
  Karagiannis, and Portela}]{Thanassoulis2016-du_JJ}
Thanassoulis, E., De~Witte, K., Johnes, J., Johnes, G., Karagiannis, G.,
  Portela, C.~S., 2016. Applications of data envelopment analysis in education.
  In: Zhu, J. (Ed.), Data Envelopment Analysis: A Handbook of Empirical Studies
  and Applications. Springer US, Boston, MA, pp. 367--438.

\bibitem[{Thanassoulis et~al.(2011)Thanassoulis, Kortelainen, Johnes, and
  Johnes}]{Thanassoulis2011-yg_JJ}
Thanassoulis, E., Kortelainen, M., Johnes, G., Johnes, J., 2011. Costs and
  efficiency of higher education institutions in {E}ngland: a {DEA} analysis.
  Journal of the Operational Research Society 62~(7), 1282--1297.

\bibitem[{{The AnyLogic Company}(2022)}]{Anylogistix_CTCGE}
{The AnyLogic Company}, 2022. {AnyLogistix Supply Chain Software}.
\newline\urlprefix\url{https://www.anylogistix.com/}

\bibitem[{Theorin et~al.(2017)Theorin, Bengtsson, Provost, Lieder, Johnsson,
  Lundholm, and Lennartson}]{Theorin2017-ae_KESXZ}
Theorin, A., Bengtsson, K., Provost, J., Lieder, M., Johnsson, C., Lundholm,
  T., Lennartson, B., 2017. An event-driven manufacturing information system
  architecture for industry 4.0. International Journal of Production Research
  55~(5), 1297--1311.

\bibitem[{Thies et~al.(2019)Thies, Kieckh{\"a}fer, Spengler, and
  Sodhi}]{Thies2019-by_JL}
Thies, C., Kieckh{\"a}fer, K., Spengler, T.~S., Sodhi, M.~S., 2019. Operations
  research for sustainability assessment of products: A review. European
  Journal of Operational Research 274~(1), 1--21.

\bibitem[{Thoben et~al.(2017)Thoben, Wiesner, and Wuest}]{Thoben2017-il_KESXZ}
Thoben, K.-D., Wiesner, S., Wuest, T., 2017. ``{I}ndustrie 4.0'' and {Smart
  Manufacturing - A Review of Research Issues and Application Examples}.
  International Journal of Automation Technology 11~(1), 4--16.

\bibitem[{Thompson(1994)}]{Tho94_JH}
Thompson, W., 1994. Cooperative models of bargaining. In: Aumann, R.~J., Hart,
  S. (Eds.), Handbook of Game Theory. Vol.~2. North-Holland, pp. 1237--1284.

\bibitem[{Thunhurst(1992)}]{Thunhurst1992-vg_AJG}
Thunhurst, C., 1992. Operational research: A role in strengthening community
  participation? Journal of Management in Medicine 6~(4), 56--71.

\bibitem[{Tijms(1994)}]{tijms1994stochastic_DLLD}
Tijms, H.~C., 1994. Stochastic models: An algorithmic approach. Vol. 303.
  Wiley, Chichester.

\bibitem[{Tine(2005)}]{Tine2005_JLYHK}
Tine, G.~C., 2005. Berlin airlift: Logistics, humanitarian aid, and strategic
  success. Army Logistician 37~(5), 39--41.

\bibitem[{Tippong et~al.(2022)Tippong, Petrovic, and
  Akbari}]{Tippong2022-ec_CV}
Tippong, D., Petrovic, S., Akbari, V., 2022. A review of applications of
  operational research in healthcare coordination in disaster management.
  European Journal of Operational Research 301~(1), 1--17.

\bibitem[{Todd(1988)}]{Todd1988-uj_JMB}
Todd, M.~J., 1988. Improved bounds and containing ellipsoids in karmarkar's
  linear programming algorithm. Mathematics of Operations Research 13~(4),
  650--659.

\bibitem[{Toffolo et~al.(2015)Toffolo, Wauters, and Trick}]{TUP_Website_GVBSP}
Toffolo, T. A.~M., Wauters, T., Trick, M., 2015. An automated benchmark website
  for the traveling umpire problem.
\newline\urlprefix\url{http://gent.cs.kuleuven.be/tup}

\bibitem[{Tohidi et~al.(2012)Tohidi, Razavyan, and Tohidnia}]{Tohidi2012-nr_SL}
Tohidi, G., Razavyan, S., Tohidnia, S., 2012. A global cost {M}almquist
  productivity index using data envelopment analysis. Journal of the
  Operational Research Society 63~(1), 72--78.

\bibitem[{Toktay et~al.(2000)Toktay, Wein, and Zenios}]{Toktay2000-iy_AM}
Toktay, L.~B., Wein, L.~M., Zenios, S.~A., 2000. Inventory management of
  remanufacturable products. Management Science 46~(11), 1412--1426.

\bibitem[{Toledo et~al.(2013)Toledo, Carravilla, Ribeiro, Oliveira, and
  Gomes}]{Toledo2013-xy_JB}
Toledo, F. M.~B., Carravilla, M.~A., Ribeiro, C., Oliveira, J.~F., Gomes,
  A.~M., 2013. The {Dotted-Board} model: A new {MIP} model for nesting
  irregular shapes. International Journal of Production Economics 145~(2),
  478--487.

\bibitem[{Tolk(2012)}]{Tolk2012-fx_KVRH}
Tolk, A., 2012. Engineering Principles of Combat Modeling and Distributed
  Simulation. Wiley, Hoboken, NJ.

\bibitem[{Tolsto{\i}(1930)}]{tolstoi1930metody_JLYHK}
Tolsto{\i}, A.~N., 1930. Metody nakhozhdeniya naimen’shego summovogo
  kilometrazha pri planirovanii perevozok v prostranstve. Planirovanie
  Perevozok, Sbornik pervy{\i}, Transpechat’NKPS, Moscow, 23--55.

\bibitem[{Tolwinski et~al.(1986)Tolwinski, A., and Leitmann}]{Toetal1986_GZ}
Tolwinski, B., A., H., Leitmann, G., 1986. Cooperative equilibria in
  differential games. Journal of Mathematical Analysis and Applications 119,
  182--202.

\bibitem[{Tompkins et~al.(2010)Tompkins, White, Bozer, and
  Tanchoco}]{tompkins2010facilities_SAA}
Tompkins, J., White, J., Bozer, Y., Tanchoco, J., 2010. Facilities Planning.
  Wiley.

\bibitem[{Tone(2002)}]{Tone2002-ha_SL}
Tone, K., 2002. A slacks-based measure of super-efficiency in data envelopment
  analysis. European Journal of Operational Research 143~(1), 32--41.

\bibitem[{Tone and Tsutsui(2009)}]{Tone2009-ph_SL}
Tone, K., Tsutsui, M., 2009. Network {DEA}: A slacks-based measure approach.
  European Journal of Operational Research 197~(1), 243--252.

\bibitem[{Tone and Tsutsui(2010)}]{Tone2010-gh_SL}
Tone, K., Tsutsui, M., 2010. Dynamic {DEA}: A slacks-based measure approach.
  Omega 38~(3), 145--156.

\bibitem[{Tone and Tsutsui(2014)}]{Tone2014-hp_SL}
Tone, K., Tsutsui, M., 2014. Dynamic {DEA} with network structure: A
  slacks-based measure approach. Omega 42~(1), 124--131.

\bibitem[{Topaloglu and Powell(2006)}]{Topaloglu2006-bn_HL}
Topaloglu, H., Powell, W.~B., 2006. Dynamic-programming approximations for
  stochastic time-staged integer multicommodity-flow problems. INFORMS Journal
  on Computing 18~(1), 31--42.

\bibitem[{Torres et~al.(2017)Torres, Kunc, and O'Brien}]{Torres2017-vm_MCJM}
Torres, J.~P., Kunc, M., O'Brien, F., 2017. Supporting strategy using system
  dynamics. European Journal of Operational Research 260~(3), 1081--1094.

\bibitem[{Toth and Vigo(2002)}]{TothV2002_CA_MB}
Toth, P., Vigo, D. (Eds.), 2002. The Vehicle Routing Problem. SIAM Monographs
  on Discrete Mathematics and Applications, Philadelphia.

\bibitem[{Toth and Vigo(2014)}]{TV14_SMPT}
Toth, P., Vigo, D. (Eds.), 2014. Vehicle Routing: Problems, Methods, and
  Applications. SIAM, Philadelphia, PA.

\bibitem[{Towill(1970)}]{towill1970_XW}
Towill, D.~R., 1970. Transfer function techniques for control engineers.
  Iliffe.

\bibitem[{Towill(1982)}]{Towill1982_SMD}
Towill, D.~R., 1982. Dynamic analysis of an inventory and order based
  production control system. International Journal of Production Research
  20~(6), 671--687.

\bibitem[{Traub(2020)}]{Traub:2020_IL}
Traub, V., 2020. Approximation algorithms for traveling salesman problems.
  Ph.D. thesis, Universit{\"a}ts-und Landesbibliothek Bonn.

\bibitem[{Trick(2001)}]{TTP_GVBSP}
Trick, M., 2001. Challenge traveling tournament instances.
\newline\urlprefix\url{https://mat.tepper.cmu.edu/TOURN/}

\bibitem[{Trick et~al.(2012)Trick, Yildiz, and Yunes}]{TUP_GVBSP}
Trick, M.~A., Yildiz, H., Yunes, T., 2012. Scheduling major league baseball
  umpires and the traveling umpire problem. Interfaces 42~(3), 232--244.

\bibitem[{Trigeorgis(1996)}]{trigeorgis96_MB}
Trigeorgis, L., 1996. Real options: Managerial flexibility and strategy in
  resource allocation. MIT press, Cambridge, MA.

\bibitem[{Trigeorgis and Tsekrekos(2018)}]{Trigeorgis2018_EA}
Trigeorgis, L., Tsekrekos, A.~E., 2018. Real options in operations research: A
  review. European Journal of Operational Research 270~(1), 1--24.

\bibitem[{Trkman et~al.(2010)Trkman, McCormack, de~Oliveira, and
  Ladeira}]{Trkman2010-nx_JEB}
Trkman, P., McCormack, K., de~Oliveira, M. P.~V., Ladeira, M.~B., 2010. The
  impact of business analytics on supply chain performance. Decision Support
  Systems 49~(3), 318--327.

\bibitem[{Trudeau and Dror(1992)}]{trudeau1992stochastic_JLYHK}
Trudeau, P., Dror, M., 1992. Stochastic inventory routing: Route design with
  stockouts and route failures. Transportation Science 26~(3), 171--184.

\bibitem[{Tsai and D.~Gemmill(1998)}]{Tsai1998-tz_HL}
Tsai, Y.-W., D.~Gemmill, D., 1998. Using tabu search to schedule activities of
  stochastic resource-constrained projects. European Journal of Operational
  Research 111~(1), 129--141.

\bibitem[{Tsoukias(2021)}]{Tsoukias2021-pz_JH}
Tsoukias, A., 2021. Social responsibility of algorithms: An overview. In:
  Papathanasiou, J., Zarat{\'e}, P., Freire~de Sousa, J. (Eds.), {EURO} Working
  Group on {DSS}: A Tour of the {DSS} Developments Over the Last 30 Years.
  Springer, Cham, pp. 153--166.

\bibitem[{Tufte(2001)}]{Tufte2001-bh_MJE}
Tufte, E.~R., 2001. The Visual Display of Quantitative Information. Graphics
  Press.

\bibitem[{Tully et~al.(2019)Tully, White, and Yearworth}]{Tully2019-dd_MYLW}
Tully, P., White, L., Yearworth, M., 2019. The value paradox of problem
  structuring methods. Systems Research and Behavioral Science 36~(4),
  424--444.

\bibitem[{Turnitsa et~al.(2021)Turnitsa, Blais, and
  Tolk}]{Turnitsa2021-dz_KVRH}
Turnitsa, C., Blais, C., Tolk, A., 2021. Simulation and Wargaming. Wiley,
  Hoboken, NJ.

\bibitem[{Tustin(1953)}]{tustin1953_XW}
Tustin, A., 1953. The Mechanism of Economic Systems: An approach to the problem
  of economic stabilization from the point of view of control-system
  engineering. MA: Harvard University.

\bibitem[{Tversky and Kahneman(1974)}]{tversky1974judgment_MESG}
Tversky, A., Kahneman, D., 1974. Judgment under uncertainty: Heuristics and
  biases: Biases in judgments reveal some heuristics of thinking under
  uncertainty. Science 185~(4157), 1124--1131.

\bibitem[{Uchoa et~al.(2017)Uchoa, Pecin, Pessoa, Poggi, Vidal, and
  Subramanian}]{UCHOA2017845_CA_MB}
Uchoa, E., Pecin, D., Pessoa, A., Poggi, M., Vidal, T., Subramanian, A., 2017.
  New benchmark instances for the capacitated vehicle routing problem. European
  Journal of Operational Research 257~(3), 845--858.

\bibitem[{Ulmer(2020)}]{ulmerDynamicPricingRouting2020_AKSJF}
Ulmer, M.~W., 2020. Dynamic pricing and routing for same-day delivery.
  Transportation Science 54~(4), 1016--1033.

\bibitem[{Ulmer and Thomas(2018)}]{ulmer2018same_CKTVW}
Ulmer, M.~W., Thomas, B.~W., 2018. Same-day delivery with heterogeneous fleets
  of drones and vehicles. Networks 72~(4), 475--505.

\bibitem[{Ulrich et~al.(2021)Ulrich, Jahnke, Langrock, Pesch, and
  Senge}]{ulrich2021distributional_TVWCK}
Ulrich, M., Jahnke, H., Langrock, R., Pesch, R., Senge, R., 2021.
  Distributional regression for demand forecasting in e-grocery. European
  Journal of Operational Research 294~(3), 831--842.

\bibitem[{Ulrich(1983)}]{Ulrich1983-ki_AJG}
Ulrich, W., 1983. Critical Heuristics of Social Planning. Wiley, Chichester.

\bibitem[{Ulrich(1987)}]{Ulrich1987-hh_GM}
Ulrich, W., 1987. Critical heuristics of social systems design. European
  Journal of Operational Research 31~(3), 276--283.

\bibitem[{Ulrich(1994)}]{Ulrich1994-fw_GM}
Ulrich, W., 1994. Critical Heuristics of Social Planning: A New Approach to
  Practical Philosophy. Wiley.

\bibitem[{{UNCTAD}(2022)}]{Unctad2022-sd_HP}
{UNCTAD}, 2022. Review of maritime transport 2022. Tech. rep., United Nations
  Conference on Trade And Development.

\bibitem[{Uniejewski et~al.(2016)Uniejewski, Nowotarski, and
  Weron}]{uni:now:wer:16_DSRW}
Uniejewski, B., Nowotarski, J., Weron, R., 2016. Automated variable selection
  and shrinkage for day-ahead electricity price forecasting. Energies 9~(8),
  621.

\bibitem[{Utley et~al.(2022)Utley, Crowe, and
  Pagel}]{utley_crowe_pagel_2022_CV}
Utley, M., Crowe, S., Pagel, C., 2022. Operational Research Approaches.
  Elements of Improving Quality and Safety in Healthcare. Cambridge University
  Press.

\bibitem[{Utomo et~al.(2018)Utomo, Onggo, and Eldridge}]{Utomo2018-uy_AFRH}
Utomo, D.~S., Onggo, B.~S., Eldridge, S., 2018. Applications of agent-based
  modelling and simulation in the agri-food supply chains. European Journal of
  Operational Research 269~(3), 794--805.

\bibitem[{Vahdati~Daneshmand(2003)}]{Daneshmand:2003_IL}
Vahdati~Daneshmand, S., 2003. Algorithmic approaches to the {{S}teiner} problem
  in networks. Ph.D. thesis, University of Mannheim, Germany.

\bibitem[{Valero-Carreras et~al.(2022)Valero-Carreras, Aparicio, and
  Guerrero}]{Valero-Carreras2022-gd_SL}
Valero-Carreras, D., Aparicio, J., Guerrero, N.~M., 2022. Multi-output support
  vector frontiers. Computers \& Operations Research 143, 105765.

\bibitem[{Valladares et~al.(2022)Valladares, Nino, Mart{\'\i}nez, Sobek,
  Claudio, and Moyce}]{Valladares2022-yk_CV}
Valladares, L., Nino, V., Mart{\'\i}nez, K., Sobek, D., Claudio, D., Moyce, S.,
  2022. Optimizing patient flow, capacity, and performance of {COVID-19}
  vaccination clinics. IISE Transactions on Healthcare Systems Engineering,
  1--13.

\bibitem[{{Van Bulck} et~al.(2021){Van Bulck}, Goossens, Beli\"en, and
  Davari}]{sports_competition_GVBSP}
{Van Bulck}, D., Goossens, D., Beli\"en, J., Davari, M., 2021. {The Fifth
  International Timetabling Competition} ({ITC} 2021): Sports timetabling. In:
  Proceedings of MathSport International. University of Reading, pp. 117--122.

\bibitem[{Van~de Ven and Delbeco(1971)}]{Van_de_Ven1971-vd_FP}
Van~de Ven, A., Delbeco, A.~L., 1971. Nominal versus interacting group
  processes for committee {Decision-Making} effectiveness. Academy of
  Management Journal 14~(2), 203--212.

\bibitem[{Van~de Ven and Poole(2005)}]{Van_de_Ven2005-ug_AFRH}
Van~de Ven, A.~H., Poole, M.~S., 2005. Alternative approaches for studying
  organizational change. Organization Studies 26~(9), 1377--1404.

\bibitem[{Van~de Vonder et~al.(2008)Van~de Vonder, Demeulemeester, and
  Herroelen}]{Van_de_Vonder2008-hv_WH_ED}
Van~de Vonder, S., Demeulemeester, E., Herroelen, W., 2008. Proactive heuristic
  procedures for robust project scheduling: An experimental analysis. European
  Journal of Operational Research 189~(3), 723--733.

\bibitem[{{van der Hagen} et~al.(2022){van der Hagen}, Agatz, Spliet, Visser,
  and Kok}]{vanderhagen2022machine_CKTVW}
{van der Hagen}, L., Agatz, N., Spliet, R., Visser, T.~R., Kok, A.~L., 2022.
  Machine learning-based feasability checks for dynamic time slot management.
  SSRN Electronic Journal 4011237.

\bibitem[{van~der Laan et~al.(1999)van~der Laan, Salomon, Dekker, and
  Van~Wassenhove}]{Van_der_Laan1999-kp_AM}
van~der Laan, E., Salomon, M., Dekker, R., Van~Wassenhove, L., 1999. Inventory
  control in hybrid systems with remanufacturing. Management Science 45~(5),
  733--747.

\bibitem[{Van~Engeland et~al.(2020)Van~Engeland, Beli{\"e}n, De~Boeck, and
  De~Jaeger}]{Van_Engeland2020-qq_JL}
Van~Engeland, J., Beli{\"e}n, J., De~Boeck, L., De~Jaeger, S., 2020. Literature
  review: Strategic network optimization models in waste reverse supply chains.
  Omega 91, 102012.

\bibitem[{van Leeuwaarden et~al.(2019)van Leeuwaarden, Mathijsen, and
  Zwart}]{van2019economies_HATI}
van Leeuwaarden, J.~S., Mathijsen, B.~W., Zwart, B., 2019. {Economies-of-Scale
  in Many-Server Queueing Systems: Tutorial and Partial Review of the QED
  Halfin--Whitt Heavy-Traffic Regime}. SIAM Review 61~(3), 403--440.

\bibitem[{Van~Roy and Wolsey(1986)}]{VW86_ALAL}
Van~Roy, T., Wolsey, L., 1986. Valid inequalities for mixed 0–1 programs.
  Discrete Applied Mathematics 14, 199--213.

\bibitem[{Van~Roy and Wolsey(1987)}]{VW87_ALAL}
Van~Roy, T., Wolsey, L., 1987. Solving mixed integer programming problems using
  automatic reformulation. Operations Research 35, 45--57.

\bibitem[{Van~Ryzin and McGill(2000)}]{van2000revenue_VLVV}
Van~Ryzin, G., McGill, J., 2000. Revenue management without forecasting or
  optimization: An adaptive algorithm for determining airline seat protection
  levels. Management Science 46~(6), 760--775.

\bibitem[{Van~Slyke and Wets(1969)}]{Van_Slyke1969-qz_HL}
Van~Slyke, R.~M., Wets, R., 1969. {L-Shaped} linear programs with applications
  to optimal control and stochastic programming. SIAM Journal on Applied
  Mathematics 17~(4), 638--663.

\bibitem[{Vanhoucke(2018)}]{Vanhoucke2018-vv_WH_ED}
Vanhoucke, M., 2018. Planning projects with scarce resources: Yesterday, today
  and tomorrow's research challenges. Frontiers of Engineering Management
  5~(2), 133--149.

\bibitem[{Vanhoucke et~al.(2001{\natexlab{a}})Vanhoucke, Demeulemeester, and
  Herroelen}]{Vanhoucke2001-fb_WH_ED}
Vanhoucke, M., Demeulemeester, E., Herroelen, W., 2001{\natexlab{a}}. An exact
  procedure for the resource-constrained weighted earliness--tardiness project
  scheduling problem. Annals of Operations Research 102~(1), 179--196.

\bibitem[{Vanhoucke et~al.(2001{\natexlab{b}})Vanhoucke, Demeulemeester, and
  Herroelen}]{Vanhoucke2001-kp_WH_ED}
Vanhoucke, M., Demeulemeester, E., Herroelen, W., 2001{\natexlab{b}}.
  Maximizing the net present value of a project with linear time-dependent cash
  flows. International Journal of Production Research 39~(14), 3159--3181.

\bibitem[{Vanhoucke et~al.(2001{\natexlab{c}})Vanhoucke, Demeulemeester, and
  Herroelen}]{Vanhoucke2001-nb_WH_ED}
Vanhoucke, M., Demeulemeester, E., Herroelen, W., 2001{\natexlab{c}}. On
  maximizing the net present value of a project under renewable resource
  constraints. Management Science 47~(8), 1113--1121.

\bibitem[{Vasilakis et~al.(2013)Vasilakis, Pagel, Gallivan, Richards, Weaver,
  and Utley}]{Vasilakis2013-sn_CV}
Vasilakis, C., Pagel, C., Gallivan, S., Richards, D., Weaver, A., Utley, M.,
  2013. Modelling toolkit to assist with introducing a stepped care system
  design in mental health care. Journal of the Operational Research Society
  64~(7), 1049--1059.

\bibitem[{Vassian(1955)}]{Vassian1955_XW}
Vassian, H.~J., 1955. Application of discrete variable servo theory to
  inventory control. Journal of the Operations Research Society of America
  3~(3), 272--282.

\bibitem[{Vaswani et~al.(2017)Vaswani, Shazeer, Parmar, Uszkoreit, Jones,
  Gomez, Kaiser, and Polosukhin}]{NIPS2017_3f5ee243_LCAL}
Vaswani, A., Shazeer, N., Parmar, N., Uszkoreit, J., Jones, L., Gomez, A.~N.,
  Kaiser, L.~u., Polosukhin, I., 2017. Attention is all you need. In: Guyon,
  I., Luxburg, U.~V., Bengio, S., Wallach, H., Fergus, R., Vishwanathan, S.,
  Garnett, R. (Eds.), Advances in Neural Information Processing Systems.
  Vol.~30. Curran Associates, Inc., pp. 1--11.

\bibitem[{Vazirani(2001)}]{Vazi10_UPCT}
Vazirani, V.~V., 2001. Approximation Algorithms. Springer, Berlin.

\bibitem[{Vedantam and Iyer(2021)}]{Vedantam2021-zj_AM}
Vedantam, A., Iyer, A., 2021. Revenue‐sharing contracts under quality
  uncertainty in remanufacturing. Production and Operations Management 30~(7),
  2008--2026.

\bibitem[{Veinott~Jr(1965)}]{veinott1965optimal_JSS}
Veinott~Jr, A.~F., 1965. Optimal policy for a multi-product, dynamic,
  nonstationary inventory problem. Management Science 12~(3), 206--222.

\bibitem[{Veinott~Jr(1966)}]{veinott1966status_JSS}
Veinott~Jr, A.~F., 1966. The status of mathematical inventory theory.
  Management Science 12~(11), 745--777.

\bibitem[{Velez-Castiblanco et~al.(2016)Velez-Castiblanco, Brocklesby, and
  Midgley}]{Velez-Castiblanco2016-kh_AFRH}
Velez-Castiblanco, J., Brocklesby, J., Midgley, G., 2016. Boundary games: How
  teams of {OR} practitioners explore the boundaries of intervention. European
  Journal of Operational Research 249~(3), 968--982.

\bibitem[{Vennix and Vennix(1996)}]{Vennix1996-ei_GM}
Vennix, J. A.~M., Vennix, J., 1996. Group Model Building: Facilitating Team
  Learning Using System Dynamics. Wiley.

\bibitem[{{Ventana Systems Inc.}(2022)}]{Vensim_CTCGE}
{Ventana Systems Inc.}, 2022. {Vensim}.
\newline\urlprefix\url{https://vensim.com/}

\bibitem[{Ventosa et~al.(2005)Ventosa, Baillo, Ramos, and
  Rivier}]{ven:bai:ram:riv:05_DSRW}
Ventosa, M., Baillo, A., Ramos, A., Rivier, M., 2005. Electricity market
  modeling trends. Energy Policy 33~(7), 897--913.

\bibitem[{Verloop et~al.(2010)Verloop, Ayesta, and Borst}]{VerAyeBor10_JH}
Verloop, I.~M., Ayesta, U., Borst, S., 2010. Monotonicity properties for
  multi-class queueing systems. Discrete Event Dynamic Systems 20, 473--509.

\bibitem[{Verma and Rubin(2018)}]{VerRub18_JH}
Verma, S., Rubin, J., 2018. Fairness definitions explained. In: Proceedings of
  the International Workshop on Software Fairness (FairWare). pp. 1--7.

\bibitem[{Verstegen(1996)}]{Ver96_JH}
Verstegen, D.~A., 1996. Concepts and measures of fiscal inequality: {A} new
  approach and effects for five states. Journal of Education Finance 22,
  145--160.

\bibitem[{Vidal(2022{\natexlab{a}})}]{Vidal2022_CA_MB}
Vidal, T., 2022{\natexlab{a}}. Github webpage.
  \url{https://github.com/vidalt/}, accessed on 2022-09-16.

\bibitem[{Vidal(2022{\natexlab{b}})}]{VIDAL2022105643_CA_MB}
Vidal, T., 2022{\natexlab{b}}. Hybrid genetic search for the {CVRP}:
  Open-source implementation and {SWAP}* neighborhood. Computers \& Operations
  Research 140, 105643.

\bibitem[{Vidal et~al.(2013)Vidal, Crainic, Gendreau, and
  Prins}]{VIDAL20131_CA_MB}
Vidal, T., Crainic, T., Gendreau, M., Prins, C., 2013. Heuristics for
  multi-attribute vehicle routing problems: A survey and synthesis. European
  Journal of Operational Research 231~(1), 1--21.

\bibitem[{Vidal et~al.(2020)Vidal, Laporte, and Matl}]{VIDAL2020401_CA_MB}
Vidal, T., Laporte, G., Matl, P., 2020. A concise guide to existing and
  emerging vehicle routing problem variants. European Journal of Operational
  Research 286, 401--416.

\bibitem[{Vidgen et~al.(2020)Vidgen, Hindle, and Randolph}]{Vidgen2020-vy_JEB}
Vidgen, R., Hindle, G., Randolph, I., 2020. Exploring the ethical implications
  of business analytics with a business ethics canvas. European Journal of
  Operational Research 281~(3), 491--501.

\bibitem[{Vidgen et~al.(2017)Vidgen, Shaw, and Grant}]{Vidgen2017-vp_JEB}
Vidgen, R., Shaw, S., Grant, D.~B., 2017. Management challenges in creating
  value from business analytics. European Journal of Operational Research
  261~(2), 626--639.

\bibitem[{Virtanen et~al.(2022)Virtanen, Mansikka, Kontio, and
  Harris}]{Virtanen2021-ft_KVRH}
Virtanen, K., Mansikka, H., Kontio, H., Harris, D., 2022. Weight watchers:
  {NASA-TLX} weights revisited. Theoretical Issues in Ergonomics Science
  23~(6), 725--748.

\bibitem[{Virtanen et~al.(2004)Virtanen, Raivio, and
  H{\"a}m{\"a}l{\"a}inen}]{Virtanen2004-su_KVRH}
Virtanen, K., Raivio, T., H{\"a}m{\"a}l{\"a}inen, R.~P., 2004. Modeling pilot's
  sequential maneuvering decisions by a multistage influence diagram. Journal
  of Guidance, Control, and Dynamics 27~(4), 665--677.

\bibitem[{von Bertalanffy(1968)}]{Von_Bertalanffy1968-yb_GM}
von Bertalanffy, L., 1968. General System Theory: Foundations, Development,
  Applications. G. Braziller.

\bibitem[{von Neumann(1945)}]{Neumann1945-wq_JMB}
von Neumann, J., 1945. A model of general economic equilibrium. The Review of
  Economic Studies 13~(1), 1--9.

\bibitem[{von Neumann and Morgenstern(1944)}]{Vomo1944_GZ}
von Neumann, J., Morgenstern, O., 1944. Theory of games and economic behavior,
  2nd Edition. Princeton University Press, Princeton.

\bibitem[{von Nitzsch and Weber(1993)}]{Von_Nitzsch1993-hh_AFRH}
von Nitzsch, R., Weber, M., 1993. The effect of attribute ranges on weights in
  multiattribute utility measurements. Management Science 39~(8), 937--943.

\bibitem[{von Stackelberg(1934)}]{vo1934_GZ}
von Stackelberg, H., 1934. Marktform und gleichgewicht. Springer Verlag,
  Vienna.

\bibitem[{Vossen and Ball(2006)}]{vossen2006slot_VLVV}
Vossen, T.~W., Ball, M.~O., 2006. Slot trading opportunities in collaborative
  ground delay programs. Transportation Science 40~(1), 29--43.

\bibitem[{Vranas et~al.(1994)Vranas, Bertsimas, and
  Odoni}]{vranas1994multi_VLVV}
Vranas, P.~B., Bertsimas, D.~J., Odoni, A.~R., 1994. The multi-airport
  ground-holding problem in air traffic control. Operations Research 42~(2),
  249--261.

\bibitem[{Vygen(2002)}]{Vygen:2002_IL}
Vygen, J., 2002. On dual minimum cost flow algorithms. Mathematical Methods of
  Operations Research 56~(1), 101--126.

\bibitem[{W{\"a}chter and Biegler(2006)}]{Ipopt_CTCGE}
W{\"a}chter, A., Biegler, L., 2006. On the implementation of a primal-dual
  interior point filter line search algorithm for large-scale nonlinear
  programming. Mathematical Programming 1, 25--57.

\bibitem[{Wagner and Whitin(1958)}]{wagner1958dynamic_JSS}
Wagner, H.~M., Whitin, T.~M., 1958. Dynamic version of the economic lot size
  model. Management Science 5~(1), 89--96.

\bibitem[{Wagner and Radovilsky(2012)}]{Wagner2012-rt_KVRH}
Wagner, M.~R., Radovilsky, Z., 2012. Optimizing boat resources at the {U.S}.
  coast guard: Deterministic and stochastic models. Operations Research 60~(5),
  1035--1049.

\bibitem[{Waisel et~al.(2008)Waisel, Wallace, and
  Willemain}]{Waisel2008-uw_AFRH}
Waisel, L.~B., Wallace, W.~A., Willemain, T.~R., 2008. Visualization and model
  formulation: an analysis of the sketches of expert modellers. Journal of the
  Operational Research Society 59~(3), 353--361.

\bibitem[{Wallace and Ziemba(2005)}]{Wallace2005-lf_HL}
Wallace, S.~W., Ziemba, W.~T., 2005. Applications of Stochastic Programming.
  SIAM.

\bibitem[{Walling and Vaneeckhaute(2020)}]{Walling2020-kn_JL}
Walling, E., Vaneeckhaute, C., 2020. Developing successful environmental
  decision support systems: Challenges and best practices. Journal of
  Environmental Management 264, 110513.

\bibitem[{Wang et~al.(2022{\natexlab{a}})Wang, Dang, Nguyen, and
  Wang}]{Wang_a_2022_DSRW}
Wang, C.~N., Dang, T.~T., Nguyen, N. A.~T., Wang, J.~W., 2022{\natexlab{a}}. A
  combined {D}ata {E}nvelopment {A}nalysis {(DEA)} and {G}rey {B}ased
  {M}ultiple {C}riteria {D}ecision {M}aking {(G-MCDM)} for solar {PV} power
  plants site selection: {A} case study in {V}ietnam. Energy Reports 8~(1),
  1124--1142.

\bibitem[{Wang et~al.(2017)Wang, Li, Wang, Peng, Jiang, and
  Liu}]{wan:etal:17_DSRW}
Wang, H.-Z., Li, G.-Q., Wang, G.-B., Peng, J.-C., Jiang, H., Liu, Y.-T., 2017.
  Deep learning based ensemble approach for probabilistic wind power
  forecasting. Applied Energy 188, 56--70.

\bibitem[{Wang et~al.(2021)Wang, Zhao, and Huchzermeier}]{wang21_MB}
Wang, J., Zhao, L., Huchzermeier, A., 2021. Operations-finance interface in
  risk management: Research evolution and opportunities. Production and
  Operations Management 30~(2), 355--389.

\bibitem[{Wang and Jacquillat(2020)}]{wang2020stochastic_VLVV}
Wang, K., Jacquillat, A., 2020. A stochastic integer programming approach to
  air traffic scheduling and operations. Operations Research 68~(5),
  1375--1402.

\bibitem[{Wang et~al.(2022{\natexlab{b}})Wang, Jacquillat, and
  Vaze}]{wang2022vertiport_VLVV}
Wang, K., Jacquillat, A., Vaze, V., 2022{\natexlab{b}}. Vertiport planning for
  urban aerial mobility: An adaptive discretization approach. Manufacturing \&
  Service Operations Management 24~(6), 3215--3235.

\bibitem[{Wang et~al.(2019{\natexlab{a}})Wang, Xian, Lee, Wei, and
  Huang}]{Wang2019-sv_SL}
Wang, K., Xian, Y., Lee, C.-Y., Wei, Y.-M., Huang, Z., 2019{\natexlab{a}}. On
  selecting directions for directional distance functions in a non-parametric
  framework: a review. Annals of Operations Research 278~(1), 43--76.

\bibitem[{Wang et~al.(2019{\natexlab{b}})Wang, Wu, and Chen}]{Wang2019-qo_JJ}
Wang, Q., Wu, Z., Chen, X., 2019{\natexlab{b}}. Decomposition weights and
  overall efficiency in a two-stage {DEA} model with shared resources.
  Computers \& Industrial Engineering 136, 135--148.

\bibitem[{Wang and Disney(2016)}]{Wang2016_SMD}
Wang, X., Disney, S.~M., 2016. The bullwhip effect: Progress, trends and
  directions. European Journal of Operational Research 250~(3), 691--701.

\bibitem[{Wang et~al.(2022{\natexlab{c}})Wang, Hyndman, Li, and
  Kang}]{Wang2022comb_FP}
Wang, X., Hyndman, R.~J., Li, F., Kang, Y., 2022{\natexlab{c}}. Forecast
  combinations: an over 50-year review. arXiv:2205.04216.

\bibitem[{Wang et~al.(2022{\natexlab{d}})Wang, Kang, Petropoulos, and
  Li}]{Wang2022_FP}
Wang, X., Kang, Y., Petropoulos, F., Li, F., 2022{\natexlab{d}}. The
  uncertainty estimation of feature-based forecast combinations. Journal of the
  Operational Research Society 73~(5), 979--993.

\bibitem[{Wang and Curry(2012)}]{Wang2012-qi_HL}
Wang, X.~j., Curry, D.~J., 2012. A robust approach to the share-of-choice
  product design problem. Omega 40~(6), 818--826.

\bibitem[{Wang(2021)}]{Wang2021-rm_JJ}
Wang, Y., 2021. An improved machine learning and artificial intelligence
  algorithm for classroom management of {E}nglish distance education. Journal
  of Intelligent \& Fuzzy Systems 40~(2), 3477--3488.

\bibitem[{Wang et~al.(2020)Wang, Bi, Lai, and Chen}]{wang2020locating_CKTVW}
Wang, Y., Bi, M., Lai, J., Chen, Y., 2020. Locating movable parcel lockers
  under stochastic demands. Symmetry 12~(12), 2033.

\bibitem[{Wang et~al.(2022{\natexlab{e}})Wang, Kuo, Huang, Gu, and
  Hu}]{wang2022dynamic_JLYHK}
Wang, Y.-J., Kuo, Y.-H., Huang, G.~Q., Gu, W., Hu, Y., 2022{\natexlab{e}}.
  Dynamic demand-driven bike station clustering. Transportation Research Part
  E: Logistics and Transportation Review 160, 102656.

\bibitem[{Ware(2020)}]{Ware2020-sj_MJE}
Ware, C., 2020. Information Visualization: Perception for Design (Interactive
  Technologies), 4th Edition. Morgan Kaufmann.

\bibitem[{Warfield(1994)}]{Warfield1994-fa_GM}
Warfield, J.~N., 1994. A Science of Generic Design: Managing Complexity Through
  Systems Design. Iowa State University Press.

\bibitem[{W{\"a}scher et~al.(2007)W{\"a}scher, Hau{\ss}ner, and
  Schumann}]{Wascher2007-bg_JB}
W{\"a}scher, G., Hau{\ss}ner, H., Schumann, H., 2007. An improved typology of
  cutting and packing problems. European Journal of Operational Research
  183~(3), 1109--1130.

\bibitem[{Wa{\ss}muth et~al.(2022)Wa{\ss}muth, K{\"o}hler, Agatz, and
  Fleischmann}]{wassmuth2022demand_CKTVW}
Wa{\ss}muth, K., K{\"o}hler, C., Agatz, N., Fleischmann, M., 2022. Demand
  management for attended home delivery--a literature review. ERIM Report
  Series ERS-2022-002-LIS.

\bibitem[{Watson et~al.(2021)Watson, Hendricks, Stewart, and
  Durbach}]{Watson2021-od_IM}
Watson, N., Hendricks, S., Stewart, T., Durbach, I., 2021. Integrating machine
  learning and decision support in tactical decision-making in rugby union.
  Journal of the Operational Research Society 72~(10), 2274--2285.

\bibitem[{Wei and Vaze(2018)}]{wei2018modeling_VLVV}
Wei, K., Vaze, V., 2018. Modeling crew itineraries and delays in the national
  air transportation system. Transportation Science 52~(5), 1276--1296.

\bibitem[{Wei et~al.(2020)Wei, Vaze, and Jacquillat}]{wei2020airline_VLVV}
Wei, K., Vaze, V., Jacquillat, A., 2020. Airline timetable development and
  fleet assignment incorporating passenger choice. Transportation Science
  54~(1), 139--163.

\bibitem[{Weintraub et~al.(2007)Weintraub, Romero, Miranda, Epstein, and
  Bjørndal}]{Weintraub2007_EA}
Weintraub, A., Romero, C., Miranda, J.~P., Epstein, R., Bjørndal, T., 2007.
  Handbook of Operations Research in Natural Resources. Springer, New York, NY.

\bibitem[{Wen et~al.(2017)Wen, Pacino, Kontovas, and Psaraftis}]{Wen2017-ku_HP}
Wen, M., Pacino, D., Kontovas, C.~A., Psaraftis, H.~N., 2017. A multiple ship
  routing and speed optimization problem under time, cost and environmental
  objectives. Transportation Research Part D: Transport and Environment 52,
  303--321.

\bibitem[{Weron(2014)}]{wer:14_DSRW}
Weron, R., 2014. Electricity price forecasting: {A} review of the
  state-of-the-art with a look into the future. International Journal of
  Forecasting 30~(4), 1030--1081.

\bibitem[{Westcombe et~al.(2006)Westcombe, Franco, and
  Shaw}]{Westcombe2006-vm_MYLW}
Westcombe, M., Franco, L.~A., Shaw, D., 2006. Where next for {PSMs---A}
  grassroots revolution? Journal of the Operational Research Society 57~(7),
  776--778.

\bibitem[{White(2002)}]{White2002-ck_MYLW}
White, L., 2002. Size matters: Large group methods and the process of
  operational research. Journal of the Operational Research Society 53~(2),
  149--160.

\bibitem[{White(2006)}]{White2006-wo_MYLW}
White, L., 2006. Evaluating problem-structuring methods: developing an approach
  to show the value and effectiveness of {PSMs}. Journal of the Operational
  Research Society 57~(7), 842--855.

\bibitem[{White(2009)}]{White2009-di_MYLW}
White, L., 2009. Understanding problem structuring methods interventions.
  European Journal of Operational Research 199~(3), 823--833.

\bibitem[{White(2016)}]{White2016-jg_MYLW}
White, L., 2016. Behavioural operational research: Towards a framework for
  understanding behaviour in {OR} interventions. European Journal of
  Operational Research 249~(3), 827--841.

\bibitem[{White(2018)}]{White2018-rh_AJG}
White, L., 2018. A {C}ook's tour: Towards a framework for measuring the social
  impact of social purpose organisations. European Journal of Operational
  Research 268~(3), 784--797.

\bibitem[{White et~al.(2016)White, Burger, and Yearworth}]{White2016-sv_MYLW}
White, L., Burger, K., Yearworth, M., 2016. Understanding behaviour in problem
  structuring methods interventions with activity theory. European Journal of
  Operational Research 249~(3), 983--1004.

\bibitem[{White et~al.(2020)White, Kunc, Burger, and
  Malpass}]{White2020-fw_AFRH}
White, L., Kunc, M., Burger, K., Malpass, J., 2020. Behavioral Operational
  Research: A Capabilities Approach. Palgrave Macmillan, London.

\bibitem[{White and Lee(2009)}]{White2009-dd_JL}
White, L., Lee, G.~J., 2009. Operational research and sustainable development:
  Tackling the social dimension. European Journal of Operational Research
  193~(3), 683--692.

\bibitem[{Whitt(1982)}]{whitt1982heavy_HATI}
Whitt, W., 1982. On the heavy-traffic limit theorem for {GI/G/$\infty$} queues.
  Advances in Applied Probability 14~(1), 171--190.

\bibitem[{Whitt(1991)}]{whitt_HATI}
Whitt, W., 1991. A review of {$L= \lambda W$} and extensions. Queueing Systems
  9~(3), 235--268.

\bibitem[{Whitt(2002)}]{whitt2002stochastic_HATI}
Whitt, W., 2002. Stochastic-process limits: an introduction to
  stochastic-process limits and their application to queues. Space 500,
  391--426.

\bibitem[{Whitt(2018)}]{whitt2018time_HATI}
Whitt, W., 2018. Time-varying queues. Queueing Models and Service Management
  1~(2), 79--164.

\bibitem[{Whittle(1988)}]{whittle1988restless_DLLD}
Whittle, P., 1988. Restless bandits: Activity allocation in a changing world.
  Journal of Applied Probability 25~(A), 287--298.

\bibitem[{Wickramasuriya et~al.(2019)Wickramasuriya, Athanasopoulos, and
  Hyndman}]{Wickramasuriya2019_FP}
Wickramasuriya, S.~L., Athanasopoulos, G., Hyndman, R.~J., 2019. Optimal
  forecast reconciliation for hierarchical and grouped time series through
  trace minimization. Journal of the American Statistical Association
  114~(526), 804--19.

\bibitem[{Wilensky(1999)}]{NetLogo_CTCGE}
Wilensky, U., 1999. {NetLogo}.
\newline\urlprefix\url{https://github.com/NetLogo/NetLogo}

\bibitem[{Willemain(1995)}]{Willemain1995-qz_AFRH}
Willemain, T.~R., 1995. Model formulation: What experts think about and when.
  Operations Research 43~(6), 916--932.

\bibitem[{Willemain and Powell(2007)}]{Willemain2007-hj_AFRH}
Willemain, T.~R., Powell, S.~G., 2007. How novices formulate models. {P}art
  {II}: a quantitative description of behaviour. Journal of the Operational
  Research Society 58~(10), 1271--1283.

\bibitem[{{Williams} and {Cookson}(2000)}]{WilliamsCookson2000_JH}
{Williams}, A., {Cookson}, R., 2000. {Equity in Health}. In: Culyer, A.,
  Newhouse, J. (Eds.), Handbook of Health Economics. Elsevier, pp. 1863--1910.

\bibitem[{Williamson(2019)}]{Williamson:2019_IL}
Williamson, D.~P., 2019. Network Flow Algorithms. Cambridge University Press.

\bibitem[{Williamson and Shmoys(2011)}]{SW11_UPCT}
Williamson, D.~P., Shmoys, D.~B., 2011. The Design of Approximation Algorithms.
  Cambridge University Press, Cambridge.

\bibitem[{Wilson(1934)}]{wilson1934scientific_JSS}
Wilson, R., 1934. A scientific routine for stock control. Harvard Business
  Review 13, 116--128.

\bibitem[{Wilson(1998)}]{Wilson98_BC}
Wilson, R., 1998. Sequential equilibria of asymmetric ascending auctions: The
  case of log-normal distributions. Economic Theory 12, 433--440.

\bibitem[{Winkenbach et~al.(2016)Winkenbach, Roset, and
  Spinler}]{winkenbach2016strategic_MH}
Winkenbach, M., Roset, A., Spinler, S., 2016. Strategic redesign of urban mail
  and parcel networks at la poste. Interfaces 46~(5), 445--458.

\bibitem[{Winston and Venkataramanan(2003)}]{Winston2003-ls_JMB}
Winston, W.~L., Venkataramanan, M., 2003. Introduction to Mathematical
  Programming. Operations Research: Volume One, 4th Edition. Brooks/Cole --
  Thomson Learning, Pacific Grove.

\bibitem[{Winter(1987)}]{Winter:1987_IL}
Winter, P., 1987. {{S}teiner} problem in networks: {A} survey. Networks 17~(2),
  129--167.

\bibitem[{Winters(1960)}]{Winters1960_FP}
Winters, P.~R., 1960. Forecasting sales by exponentially weighted moving
  averages. Management Science 6, 324--342.

\bibitem[{Woeginger(2000)}]{Woe00_UPCT}
Woeginger, G.~J., 2000. When does a dynamic programming formulation guarantee
  the existence of a fully polynomial time approximation scheme ({FPTAS})?
  INFORMS Journal on Computing 12~(1), 57--74.

\bibitem[{Woeginger(2021)}]{Woe21_UPCT}
Woeginger, G.~J., 2021. The trouble with the second quantifier. 4OR - A
  Quarterly Journal of Operations Research 19, 157--181.

\bibitem[{Wojtczak(2018)}]{Woi18_UPCT}
Wojtczak, D., 2018. On strong {NP}-completeness of rational problems. In:
  Fomin, F.~V., Podolskii, V.~V. (Eds.), Computer Science – Theory and
  Applications. Vol. 10846 of Lecture Notes in Computer Science. pp. 308--320.

\bibitem[{Wolfers and Zitzewitz(2004)}]{Wolfers2004-mq_FP}
Wolfers, J., Zitzewitz, E., 2004. Prediction markets. The Journal of Economic
  Perspectives 18~(2), 107--126.

\bibitem[{Wollmer(1992)}]{wollmerAirlineSeatManagement1992_AKSJF}
Wollmer, R.~D., 1992. An {{Airline Seat Management Model}} for a {{Single Leg
  Route When Lower Fare Classes Book First}}. Operations Research 40~(1),
  26--37.

\bibitem[{Wolsey(1975)}]{Wo75_ALAL}
Wolsey, L., 1975. Faces for a linear inequality in 0--1 variables. Mathematical
  Programming 8, 165--178.

\bibitem[{Wolstenholme(1990)}]{Wolstenholme1990-df_MCJM}
Wolstenholme, E.~F., 1990. System Enquiry: A System Dynamics Approach. Wiley,
  Chichester.

\bibitem[{Wolszczak-Derlacz(2018)}]{Wolszczak-Derlacz2018-xo_JJ}
Wolszczak-Derlacz, J., 2018. Assessment of {TFP} in {E}uropean and {A}merican
  higher education institutions -- application of {M}almquist indices.
  Technological and Economic Development of Economy 24~(2), 467--488.

\bibitem[{Womack et~al.(1990)Womack, Jones, and Roos}]{Womack1990_SMD}
Womack, J., Jones, D., Roos, D., 1990. The Machine that Changed the World.
  Mandarin Books, London.

\bibitem[{Wong and Hiew(2020)}]{Wong2020-so_AJG}
Wong, D., Hiew, Y., 2020. Community operational research ({OR}) and design
  thinking for the health and social services: A comparative analysis. In: 2020
  {IEEE} International Conference on Industrial Engineering and Engineering
  Management ({IEEM}). Singapore, pp. 1042--1047.

\bibitem[{Wong and Mingers(1994)}]{Wong1994-lb_AJG}
Wong, N., Mingers, J., 1994. The nature of community {OR}. Journal of the
  Operational Research Society 45~(3), 245--254.

\bibitem[{Wood et~al.(2020)Wood, McWilliams, Thomas, Bourdeaux, and
  Vasilakis}]{Wood2020-wx_CV}
Wood, R.~M., McWilliams, C.~J., Thomas, M.~J., Bourdeaux, C.~P., Vasilakis, C.,
  2020. {COVID-19} scenario modelling for the mitigation of capacity-dependent
  deaths in intensive care. Health Care Management Science 23~(3), 315--324.

\bibitem[{Wood et~al.(2022)Wood, Moss, Murch, Vasilakis, and
  Clatworthy}]{Wood2022-ot_CV}
Wood, R.~M., Moss, S.~J., Murch, B.~J., Vasilakis, C., Clatworthy, P.~L., 2022.
  Optimising acute stroke pathways through flexible use of bed capacity: a
  computer modelling study. BMC Health Services Research 22~(1), 1068.

\bibitem[{Wood and Murch(2020)}]{Wood2020-jf_CV}
Wood, R.~M., Murch, B.~J., 2020. Modelling capacity along a patient pathway
  with delays to transfer and discharge. Journal of the Operational Research
  Society 71~(10), 1530--1544.

\bibitem[{Wood et~al.(2021{\natexlab{a}})Wood, Murch, Moss, Tyler, Thompson,
  and Vasilakis}]{Wood2021-jt_CV}
Wood, R.~M., Murch, B.~J., Moss, S.~J., Tyler, J. M.~B., Thompson, A.~L.,
  Vasilakis, C., 2021{\natexlab{a}}. Operational research for the safe and
  effective design of {COVID-19} mass vaccination centres. Vaccine 39~(27),
  3537--3540.

\bibitem[{Wood et~al.(2021{\natexlab{b}})Wood, Pratt, Kenward, McWilliams,
  Booton, Thomas, Bourdeaux, and Vasilakis}]{Wood2021-jv_CV}
Wood, R.~M., Pratt, A.~C., Kenward, C., McWilliams, C.~J., Booton, R.~D.,
  Thomas, M.~J., Bourdeaux, C.~P., Vasilakis, C., 2021{\natexlab{b}}. The value
  of triage during periods of intense {COVID-19} demand: Simulation modeling
  study. Medical Decision Making 41~(4), 393--407.

\bibitem[{Woodhouse and Goldstein(1988)}]{Woodhouse1988-fr_JJ}
Woodhouse, G., Goldstein, H., 1988. Educational performance indicators and
  {LEA} league tables. Oxford Review of Education 14~(3), 301--320.

\bibitem[{Woolley and Pidd(1981)}]{Woolley1981-go_MYLW}
Woolley, R.~N., Pidd, M., 1981. Problem structuring --- a literature review.
  Journal of the Operational Research Society 32~(3), 197--206.

\bibitem[{Wright(1983)}]{Wright1983_EA}
Wright, D.~J., 1983. Catastrophe theory in management forecasting and decision
  making. Journal of the Operational Research Society 34~(10), 935--942.

\bibitem[{Wright et~al.(2019)Wright, Cairns, O'Brien, and
  Goodwin}]{Wright2019-av_JL}
Wright, G., Cairns, G., O'Brien, F.~A., Goodwin, P., 2019. Scenario analysis to
  support decision making in addressing wicked problems: Pitfalls and
  potential. European Journal of Operational Research 278~(1), 3--19.

\bibitem[{Wright et~al.(2006)Wright, Liberatore, and
  Nydick}]{Wright2006-ug_KVRH}
Wright, P.~D., Liberatore, M.~J., Nydick, R.~L., 2006. A survey of operations
  research models and applications in homeland security. Interfaces 36~(6),
  514--529.

\bibitem[{Wright(1997)}]{Wright97_EAY}
Wright, S.~J., 1997. Primal-Dual Interior-Point Methods. SIAM, Philadelphia,
  PA.

\bibitem[{Xiang et~al.(2015)Xiang, Yin, and Lim}]{Xiang2015-qe_CV}
Xiang, W., Yin, J., Lim, G., 2015. A short-term operating room surgery
  scheduling problem integrating multiple nurses roster constraints. Artificial
  Intelligence in Medicine 63~(2), 91--106.

\bibitem[{Xin and Van~Mieghem(2023)}]{xin2023dual_JSS}
Xin, L., Van~Mieghem, J.~A., 2023. Dual-sourcing, dual-mode dynamic stochastic
  inventory models. In: Song, J.-S. (Ed.), Research Handbook on Inventory
  Management. Edward Elgar Publishing.

\bibitem[{Xu et~al.(2016)Xu, Huang, Hsieh, Lee, Jia, and Chen}]{Xu2016_CC}
Xu, J., Huang, E., Hsieh, L., Lee, L.~H., Jia, Q.-S., Chen, C.-H., 2016.
  Simulation optimization in the era of industrial 4.0 and the industrial
  internet. Journal of Simulation 10~(4), 310--320.

\bibitem[{{Yager}(1997)}]{Yager1997_JH}
{Yager}, R., 1997. On the analytic representation of the leximin ordering and
  its application to flexible constraint propagation. European Journal of
  Operational Research 102~(1), 176 -- 192.

\bibitem[{Yaman(2005)}]{yaman:05*1_BF}
Yaman, H., 2005. {Concentrator Location in Telecommunications Networks}.
  Vol.~16 of Combinatorial Optimization. Springer, New-York.

\bibitem[{Yan et~al.(2022{\natexlab{a}})Yan, Barnhart, and
  Vaze}]{yan2022choice_VLVV}
Yan, C., Barnhart, C., Vaze, V., 2022{\natexlab{a}}. Choice-based airline
  schedule design and fleet assignment: A decomposition approach.
  Transportation Science, DOI: 10.1287/trsc.2022.1141.

\bibitem[{Yan and Kung(2018)}]{yan2018robust_VLVV}
Yan, C., Kung, J., 2018. Robust aircraft routing. Transportation Science
  52~(1), 118--133.

\bibitem[{Yan et~al.(2022{\natexlab{b}})Yan, Archibald, Han, and
  Bian}]{yan2022whether_TVWCK}
Yan, S., Archibald, T.~W., Han, X., Bian, Y., 2022{\natexlab{b}}. Whether to
  adopt "buy online and return to store" strategy in a competitive market?
  European Journal of Operational Research 301~(3), 974--986.

\bibitem[{Yan and Wang(2021)}]{Yan2021_JJ}
Yan, W., Wang, G., 2021. Research on the development trend of foreign education
  based on machine learning and artificial intelligence simulation analysis.
  Journal of Intelligent and Fuzzy Systems, 1--10.

\bibitem[{Yan et~al.(2022{\natexlab{c}})Yan, Chow, Ho, Kuo, Wu, and
  Ying}]{yan2022reinforcement_JLYHK}
Yan, Y., Chow, A.~H., Ho, C.~P., Kuo, Y.-H., Wu, Q., Ying, C.,
  2022{\natexlab{c}}. Reinforcement learning for logistics and supply chain
  management: Methodologies, state of the art, and future opportunities.
  Transportation Research Part E: Logistics and Transportation Review 162,
  102712.

\bibitem[{Yang(2018)}]{yan:18_DSRW}
Yang, D., 2018. Ultra-fast preselection in lasso-type spatio-temporal solar
  forecasting problems. Solar Energy 176, 788--796.

\bibitem[{Yang et~al.(2022)Yang, Wang, Gueymard, Hong, Kleissl, Huang, Perez,
  Perez, Bright, Xia, van~der Meer, and Peters}]{yan:etal:22_DSRW}
Yang, D., Wang, W., Gueymard, C., Hong, T., Kleissl, J., Huang, J., Perez, M.,
  Perez, R., Bright, J., Xia, X., van~der Meer, D., Peters, I., 2022. A review
  of solar forecasting, its dependence on atmospheric sciences and implications
  for grid integration: Towards carbon neutrality. Renewable and Sustainable
  Energy Reviews 161, 112348.

\bibitem[{Yang and Strauss(2017)}]{yang2017approximate_CKTVW}
Yang, X., Strauss, A.~K., 2017. An approximate dynamic programming approach to
  attended home delivery management. European Journal of Operational Research
  263~(3), 935--945.

\bibitem[{Yang et~al.(2016)Yang, Strauss, Currie, and
  Eglese}]{yang2016choice_CKTVW}
Yang, X., Strauss, A.~K., Currie, C.~S., Eglese, R., 2016. Choice-based demand
  management and vehicle routing in e-fulfillment. Transportation Science
  50~(2), 473--488.

\bibitem[{Yardley and Petropoulos(2021)}]{Yardley2021_FP}
Yardley, E., Petropoulos, F., 2021. Beyond error measures to the utility and
  cost of the forecasts. Foresight: The International Journal of Applied
  Forecasting 63, 36--45.

\bibitem[{Ye(1987)}]{Ye1987-kd_JMB}
Ye, Y., 1987. Karmarkar's algorithm and the ellipsoid method. Operations
  Research Letters 6~(4), 177--182.

\bibitem[{Ye(1997)}]{Ye1997_EAY}
Ye, Y., 1997. Interior Point Algorithms: Theory and Analysis. Wiley, New York,
  NY.

\bibitem[{Yearworth and White(2014)}]{Yearworth2014-xf_MYLW}
Yearworth, M., White, L., 2014. The non-codified use of problem structuring
  methods and the need for a generic constitutive definition. European Journal
  of Operational Research 237~(3), 932--945.

\bibitem[{Yearworth and White(2019)}]{Yearworth2019-qi_MYLW}
Yearworth, M., White, L., 2019. Group support systems: Experiments with an
  online system and implications for {Same-Time/Different-Places} working. In:
  Kilgour, D.~M., Eden, C. (Eds.), Handbook of Group Decision and Negotiation.
  Springer, Cham, pp. 681--706.

\bibitem[{Yen and Birge(2006)}]{yen2006stochastic_VLVV}
Yen, J.~W., Birge, J.~R., 2006. A stochastic programming approach to the
  airline crew scheduling problem. Transportation Science 40~(1), 3--14.

\bibitem[{Yeung and Petrosyan(2018)}]{Yepe2018_GZ}
Yeung, D., Petrosyan, L., 2018. Nontransferable utility cooperative dynamic
  games. In: Ba\c{s}ar, T., Zaccour, G. (Eds.), Handbook of Dynamic Game
  Theory. Springer, Cham, pp. 633--670.

\bibitem[{Y{\i}ld{\i}z and Kara{\c{s}}an(2017)}]{yildiz2017regenerator_SAA}
Y{\i}ld{\i}z, B., Kara{\c{s}}an, O.~E., 2017. Regenerator location problem in
  flexible optical networks. Operations Research 65~(3), 595--620.

\bibitem[{Yin et~al.(2021)Yin, Perchet, and Soup{\'e}}]{yin21_MBl}
Yin, C., Perchet, R., Soup{\'e}, F., 2021. A practical guide to robust
  portfolio optimization. Quantitative Finance 21~(6), 911--928.

\bibitem[{Yin et~al.(2008)Yin, Kaku, and Stecke}]{Yin2008-jw_KESXZ}
Yin, Y., Kaku, I., Stecke, K.~E., 2008. The evolution of seru production
  systems throughout canon. Operations Management Education Review 2, 27--40.

\bibitem[{Yin et~al.(2017)Yin, Stecke, Swink, and Kaku}]{Yin2017-gh_KESXZ}
Yin, Y., Stecke, K.~E., Swink, M., Kaku, I., 2017. Lessons from seru production
  on manufacturing competitively in a high cost environment. Journal of
  Operations Management 49-51, 67--76.

\bibitem[{Yin and Yasuda(2006)}]{Yin2006-vm_KESXZ}
Yin, Y., Yasuda, K., 2006. Similarity coefficient methods applied to the cell
  formation problem: A taxonomy and review. International Journal of Production
  Economics 101~(2), 329--352.

\bibitem[{Yitzhaki and Schechtman(2013)}]{yitzhaki2013more_JH}
Yitzhaki, S., Schechtman, E., 2013. More than a dozen alternative ways of
  spelling {Gini}. In: Yitzhaki, S., Schechtman, E. (Eds.), The Gini
  Methodology. Springer, pp. 11--31.

\bibitem[{You(1999)}]{you1999dynamic_VLVV}
You, P.-S., 1999. Dynamic pricing in airline seat management for flights with
  multiple flight legs. Transportation Science 33~(2), 192--206.

\bibitem[{Young et~al.(2021)Young, Eyre, Kendrick, White, Smith, Beveridge,
  Nonnenmacher, Ichofu, Hillier, Oakley, Diamond, Rourke, Dawe, Day, Davies,
  Staite, Lacey, McCrae, Jones, Kelly, Bankiewicz, Tunkel, Ovens, Chapman,
  Bhalla, Marks, Hicks, Fowler, Hopkins, Yardley, and Peto}]{Young2021-ws}
Young, B.~C., Eyre, D.~W., Kendrick, S., White, C., Smith, S., Beveridge, G.,
  Nonnenmacher, T., Ichofu, F., Hillier, J., Oakley, S., Diamond, I., Rourke,
  E., Dawe, F., Day, I., Davies, L., Staite, P., Lacey, A., McCrae, J., Jones,
  F., Kelly, J., Bankiewicz, U., Tunkel, S., Ovens, R., Chapman, D., Bhalla,
  V., Marks, P., Hicks, N., Fowler, T., Hopkins, S., Yardley, L., Peto, T.
  E.~A., 2021. Daily testing for contacts of individuals with {SARS-CoV-2}
  infection and attendance and {SARS-CoV-2} transmission in english secondary
  schools and colleges: an open-label, cluster-randomised trial. The Lancet
  398~(10307), 1217--1229.

\bibitem[{Yu and He(2020)}]{Yu_2020_DSRW}
Yu, D., He, X., 2020. A bibliometric study for {DEA} applied to energy
  efficiency: Trends and future challenges. Applied Energy 268~(1), 115048.

\bibitem[{Yuan et~al.(2022)Yuan, Gao, Li, Liu, and Yang}]{Yuan2022_DC}
Yuan, J., Gao, Y., Li, S., Liu, P., Yang, L., 2022. Integrated optimization of
  train timetable, rolling stock assignment and short-turning strategy for a
  metro line. European Journal of Operational Research 301~(3), 855--874.

\bibitem[{Yudin and Nemirovskii(1976)}]{YN76_EAY}
Yudin, D., Nemirovskii, A.~S., 1976. Informational complexity and efficient
  methods for the solution of convex extremal problems. Ekonomika i
  Matematicheskie Metody 12, 128--142.

\bibitem[{Zanakis et~al.(1989)Zanakis, Evans, and
  Vazacopoulos}]{zanakis1989heuristic_COIT}
Zanakis, S.~H., Evans, J.~R., Vazacopoulos, A.~A., 1989. Heuristic methods and
  applications: A categorized survey. European Journal of Operational Research
  43~(1), 88--110.

\bibitem[{Zawacki-Richter et~al.(2019)Zawacki-Richter, Mar{\'\i}n, Bond, and
  Gouverneur}]{Zawacki-Richter2019-ki_JJ}
Zawacki-Richter, O., Mar{\'\i}n, V.~I., Bond, M., Gouverneur, F., 2019.
  Systematic review of research on artificial intelligence applications in
  higher education – where are the educators? International Journal of
  Educational Technology in Higher Education 16~(1), 1--27.

\bibitem[{Zawadzki and {\.Z}ywicki(2016)}]{Zawadzki2016-ek_KESXZ}
Zawadzki, P., {\.Z}ywicki, K., 2016. Smart product design and production
  control for effective mass customization in the industry 4.0 concept.
  Management and Production Engineering Review 7~(3), 105--112.

\bibitem[{Zeng et~al.(2007)Zeng, , and Yang}]{Zeng2007-xf_HP}
Zeng, Q., , Yang, Z., 2007. Model integrating fleet design and ship routing
  problems for coal shipping. In: Shi, Y., van Albada, G.~D., Dongarra, J.,
  Sloot, P. M.~A. (Eds.), Computational {Science--ICCS}. Springer, Berlin,
  Heidelberg, pp. 1000--1003.

\bibitem[{Zhang and Dimitrakopoulos(2018)}]{Zhang2018_EA}
Zhang, J., Dimitrakopoulos, R.~G., 2018. Stochastic optimization for a mineral
  value chain with nonlinear recovery and forward contracts. Journal of the
  Operational Research Society 69~(6), 864--875.

\bibitem[{Zhang and Alipour(2022)}]{Zhang2022-su_HL}
Zhang, N., Alipour, A., 2022. A stochastic programming approach to enhance the
  resilience of infrastructure under weather‐related risk. Computer-aided
  Civil and Infrastructure Engineering.

\bibitem[{Zhang et~al.(2021)Zhang, Zhang, Lim, and Sim}]{zhang2021robust_JLYHK}
Zhang, Y., Zhang, Z., Lim, A., Sim, M., 2021. Robust data-driven vehicle
  routing with time windows. Operations Research 69~(2), 469--485.

\bibitem[{Zhang(2023)}]{Zhang2023-xu_HL}
Zhang, Z.~G., 2023. Fundamentals of Stochastic Models. CRC Press.

\bibitem[{Zhao et~al.(2017)Zhao, Jin, and Lee}]{Zhao2017-jt_HL}
Zhao, K., Jin, J.~G., Lee, D.-H., 2017. {Two-Stage} stochastic programming
  model for robust personal rapid transit network design. Transportation
  Research Record 2650~(1), 152--162.

\bibitem[{Zhao et~al.(2016)Zhao, Bennell, Bekta{\c s}, and
  Dowsland}]{Zhao2016-nt_JB}
Zhao, X., Bennell, J.~A., Bekta{\c s}, T., Dowsland, K., 2016. A comparative
  review of {3D} container loading algorithms. International Transactions in
  Operational Research 23~(1-2), 287--320.

\bibitem[{Zhen(2015)}]{Zhen2015-mh_HP}
Zhen, L., 2015. Tactical berth allocation under uncertainty. European Journal
  of Operational Research 247~(3), 928--944.

\bibitem[{Zheng(1992)}]{zheng1992properties_JSS}
Zheng, Y.-S., 1992. On properties of stochastic inventory systems. Management
  Science 38~(1), 87--103.

\bibitem[{Zhou et~al.(2021)Zhou, Ma, Cao, Lee, and Chew}]{Zhou2021-vs_HL}
Zhou, C., Ma, N., Cao, X., Lee, L.~H., Chew, E.~P., 2021. Classification and
  literature review on the integration of simulation and optimization in
  maritime logistics studies. IISE Transactions 53~(10), 1157--1176.

\bibitem[{Zhou et~al.(2022)Zhou, Qi, Yang, Shi, Pan, and Gao}]{Zhou2022_DC}
Zhou, H., Qi, J., Yang, L., Shi, J., Pan, H., Gao, Y., 2022. Joint optimization
  of train timetabling and rolling stock circulation planning{: A} novel
  flexible train composition mode. Transportation Research Part B:
  Methodological 162, 352--385.

\bibitem[{Zhou et~al.(2018)Zhou, Yang, Chen, and Zhu}]{Zhou2018-et_JL}
Zhou, H., Yang, Y., Chen, Y., Zhu, J., 2018. Data envelopment analysis
  application in sustainability: The origins, development and future
  directions. European Journal of Operational Research 264~(1), 1--16.

\bibitem[{Zhou and Doyle(1998)}]{Zhou1998_XW}
Zhou, K., Doyle, J.~C., 1998. Essentials of robust control. Prentice Hall,
  Upper Saddle River, NJ.

\bibitem[{Zhou et~al.(2011)Zhou, Tao, and Chao}]{Zhou2011-zd_AM}
Zhou, S.~X., Tao, Z., Chao, X., 2011. Optimal control of inventory systems with
  multiple types of remanufacturable products. Manufacturing \& Service
  Operations Management 13~(1), 20--34.

\bibitem[{Zhu et~al.(2014)Zhu, Crainic, and Gendreau}]{Zhu2014_DC}
Zhu, E., Crainic, T.~G., Gendreau, M., 2014. Scheduled service network design
  for freight rail transportation. Operations Research 62~(2), 383--400.

\bibitem[{Zhu(2015)}]{Zhu2015-co_SL}
Zhu, J., 2015. Data Envelopment Analysis: A Handbook of Models and Methods.
  Springer, New York, NY.

\bibitem[{Zhu et~al.(2017)Zhu, Wu, Song, and Liang}]{Zhu2017-ol_SL}
Zhu, Q., Wu, J., Song, M., Liang, L., 2017. A unique equilibrium efficient
  frontier with fixed-sum outputs in data envelopment analysis. Journal of the
  Operational Research Society 68~(12), 1483--1490.

\bibitem[{Zhuang et~al.(2021)Zhuang, Qi, Duan, Xi, Zhu, Zhu, Xiong, and
  He}]{Zhuang_TransferLearning21_LCAL}
Zhuang, F., Qi, Z., Duan, K., Xi, D., Zhu, Y., Zhu, H., Xiong, H., He, Q.,
  2021. A comprehensive survey on transfer learning. Proceedings of the IEEE
  109~(1), 43--76.

\bibitem[{Ziel(2016)}]{zie:16_DSRW}
Ziel, F., 2016. Forecasting electricity spot prices using lasso: On capturing
  the autoregressive intraday structure. IEEE Transactions on Power Systems
  31~(6), 4977--4987.

\bibitem[{Ziel and Liu(2016)}]{zie:liu:16_DSRW}
Ziel, F., Liu, B., 2016. Lasso estimation for {GEFCom2014} probabilistic
  electric load forecasting. International Journal of Forecasting 32~(3),
  1029--1037.

\bibitem[{Ziel and Steinert(2018)}]{zie:ste:18_DSRW}
Ziel, F., Steinert, R., 2018. Probabilistic mid- and long-term electricity
  price forecasting. Renewable and Sustainable Energy Reviews 94, 251--266.

\bibitem[{Ziel and Weron(2018)}]{zie:wer:18_DSRW}
Ziel, F., Weron, R., 2018. Day-ahead electricity price forecasting with
  high-dimensional structures: Univariate vs. multivariate modeling frameworks.
  Energy Economics 70, 396--420.

\bibitem[{Zikopoulos and Tagaras(2008)}]{Zikopoulos2008-sq_AM}
Zikopoulos, C., Tagaras, G., 2008. On the attractiveness of sorting before
  disassembly in remanufacturing. IIE Transactions 40~(3), 313--323.

\bibitem[{Zipkin(2000)}]{Zipkin2000_JSS}
Zipkin, P., 2000. Foundations of inventory management. McGraw-Hill.

\bibitem[{Zis and Psaraftis(2017)}]{Zis2017-ka_HP}
Zis, T., Psaraftis, H.~N., 2017. The implications of the new sulphur limits on
  the european {Ro-Ro} sector. Transportation Research Part D: Transport and
  Environment 52, 185--201.

\bibitem[{Zis and Psaraftis(2019)}]{Zis2019-ne_HP}
Zis, T., Psaraftis, H.~N., 2019. Operational measures to mitigate and reverse
  the potential modal shifts due to environmental legislation. Maritime Policy
  \& Management 46~(1), 117--132.

\bibitem[{Zis et~al.(2020)Zis, Psaraftis, and Ding}]{Zis2020-uk_HP}
Zis, T. P.~V., Psaraftis, H.~N., Ding, L., 2020. Ship weather routing: A
  taxonomy and survey. Ocean Engineering 213, 107697.

\bibitem[{Zografos et~al.(2012)Zografos, Salouras, and
  Madas}]{zografos2012dealing_VLVV}
Zografos, K.~G., Salouras, Y., Madas, M.~A., 2012. Dealing with the efficient
  allocation of scarce resources at congested airports. Transportation Research
  Part C: Emerging Technologies 21~(1), 244--256.

\end{thebibliography}

\end{document}